\documentclass[9pt]{report}
\usepackage[english]{babel}
\usepackage{amssymb}
\usepackage{amsmath}
\usepackage{amsthm}
\usepackage[toc]{appendix}
\usepackage[dvips]{graphicx}

\input pictex.tex
\usepackage{graphicx}
\usepackage{pict2e}

\newcommand{\lh}{\ell ength}
\newcommand{\lip}{\mathrm{Diff}^{1 + \mathrm{Lip}}_{+}(\mathrm{S}^1)}
\newcommand{\vl}{\mathrm{Diff}^{1 + \mathrm{bv}}_{+}(\mathrm{S}^1)}

\newcommand{\clo}{\mathrm{S}^1}
\newcommand{\R}{\mathbb{R}}
\newcommand{\dos}{\mathrm{Diff}_{+}^2(\mathrm{S}^1)}
\newcommand{\ger}{\mathcal{G}^r_{+}(\mathbb{R},0)}

\newcommand{\ce}{\mathrm{C}}
\newcommand{\V}{\mathrm S}

\newcommand{\C}{\mathrm C}
\newcommand{\ner}{\mathcal{N}^p}
\newcommand{\db}{\partial \mathrm{ist}}
\newcommand{\tepe}{\mathcal{T}^p}

\newcommand{\sgn}{\vartriangleleft}
\newcommand{\sol}{\mathrm{sol}}
\newcommand{\precede}{\preceq}

\newcommand{\esp}{\hspace{0.1cm}}

\newcommand{\vf}{\varphi}

\theoremstyle{definition}

\newtheorem{thm}{Theorem}[section]

\newtheorem{prop}[thm]{Proposition}

\newtheorem{lem}[thm]{Lemma}

\newtheorem{obs}[thm]{Remark}

\newtheorem{defn}[thm]{Definition}

\newtheorem{cor}[thm]{Corollary}

\newtheorem{ejem}[thm]{Example}

\newtheorem{ejer}[thm]{Exercise}

\usepackage{makeidx}
\makeindex

\input epsf

\author{\Large \sc Andr\'es Navas \\ \\ 
\Large University of Santiago, Chile}
\date{}
\title{\LARGE{Groups of Circle Diffeomorphisms}}
\begin{document}

\maketitle


\thispagestyle{empty}

\thispagestyle{empty}

\pagenumbering{arabic}

$\mbox{}$

\newpage







\pagenumbering{roman}
\tableofcontents


\newpage
${}$
\thispagestyle{empty}
\newpage

\pagenumbering{arabic}


\chapter*{INTRODUCTION}

\addcontentsline{toc}{chapter}{INTRODUCTION}


\hspace{0.45cm} The theory of dynamical systems concerns the 
quantitative and qualitative study of the orbits of a map or 
the flow associated to a vector field. In the case the map is 
invertible or the vector field is regular, these systems may be 
thought of as actions of the group $\mathbb{Z}$ or $\mathbb{R}$, 
respectively. From this point of view, classical dynamics may be 
considered as a particular subject in the general theory of group 
actions.

From an algebraic point of view, the latter theory may be considered as a 
``non-linear'' version of that of group representations. In this direction, 
the internal structure of the groups is intended to be revealed by looking 
at their actions on nice spaces, preferably differentiable manifolds. Several 
algebraic notions become natural, and many topological and/or analytical 
aspects of the underlying spaces turn out to be relevant. In recent years, 
this approach has been revealed very fruitful, leading (sometimes quite 
unexpectedly) to dynamical proofs of some results of pure algebraic nature.

Without any doubt, it would be too ambitious (and perhaps impossible) to provide a complete 
treatment of the general theory of group actions on manifolds. It is then natural to 
restrict the study to certain groups or manifolds. In this book, we follow the latter 
direction by studying group actions on the simplest closed manifold, namely the 
circle. In spite of the apparent simplicity of the subject, it turns out that 
it is highly nontrivial and quite extensive. Its relevance lies mainly on its 
connections with many other branches like low dimensional geometry and topology 
(including Foliation Theory), mathematical logic (through the theory of orderable 
groups), mathematical physics (through the cohomological aspects of the theory), 
etc. In this spirit, this book has been conceived as a text where people who 
have encountered some of the topics contained here in their research may  
see them integrated in an independent theory. 

\vspace{0.1cm}

There already exists a very nice and complete survey by Ghys \cite{Gh1} on the 
theory of groups of circle {\em homeomorphisms}. Unlike \cite{Gh1}, here we  
mainly focus on the theory of groups of circle {\em diffeomorphisms}, which 
is essentially different in many aspects. Nevertheless, even in the 
context of diffeomorphisms, this text is still incomplete. 
We would have liked to add at least one section on the 
theory of small denominators (including recent solutions of both 
Moser's problem on the simultaneous conjugacy to rotations of commuting 
diffeomorphisms by Fayad and Khanin \cite{bassam}, and Frank's problem 
on distortion properties for irrational Euclidean rotations by Avila 
\cite{avila}), extend the treatment of Sacksteder's theorem by providing a 
discussion of the so-called Level Theory \cite{CC1,level,hector}, include two 
small sections on groups of real-analytic diffeomorphisms and piecewise affine 
homeomorphisms respectively, develop the notion of topological entropy for 
group actions on compact metric spaces by focusing on the case of one-dimensional 
manifolds \cite{GLW,hurder2,jorquera,walzac}, 
\index{topological entropy}
confront the dynamical and 
cohomological aspects involved in the study of the so-called Godbillon-Vey class 
(or Bott-Virasoro-Thurston cocycle) providing a proof of 
Duminy's ``third'' theorem \cite{ghys-bourbaki,Hur,hurder-survey,HL,taba}, 
and explore the representations of fundamental groups of surfaces 
\cite{ghys-ihes,Gh4}. 
\index{Duminy!third theorem}
Moreover, since we mainly focus on actions of discrete groups, we decided not to  
include some relevant topics as for example 
simplicity properties and diffusion processes 
for the groups of interval and circle diffeomorphisms (see \cite{ban} 
and \cite{mall,mal}, respectively). 

\vspace{0.25cm}

The book begins with a brief section where we establish most of the 
notation and recall some general definitions. More elaborate concepts, 
as for instance group amenability, are discussed in the Appendix.

Chapter 1 studies some simple but relevant examples of groups which do act on 
the circle. After recalling some fundamental properties of the rotation group, 
the affine group, and the M\"obius group, it treats the general case of Lie 
group actions on the circle, and concludes by discussing the very important 
Thompson groups.

In Chapter 2, we study some fundamental results about the dynamics 
of groups of interval and circle homeomorphisms. In the first part, 
we discuss some of their 
combinatorial aspects, as for example Poincar\'e's theory of rotation number 
and its relation with the invariant probability measures. We next provide necessary 
and sufficient conditions for the existence of faithful actions 
on the interval, and free actions on the interval and the circle. 
As an application, we give a dynamical 
proof of a recent and beautiful result due to Witte-Morris  
asserting that left-orderable, amenable groups are locally indicable. 
At this point, it would have been natural to treat Ghys' characterization 
of faithful actions on the circle in terms of bounded cohomology. 
\index{Ghys!theorem on bounded cohomology} 
However, we decided not to treat this topic, mainly because it is 
well developed in \cite{Gh1}, and furthermore because, though it is very important 
for describing continuous actions on the circle, its relevance for the smooth 
theory is much smaller. 
In the second part of Chapter 2, we essentially treat a result due to Margulis, 
which may be thought of as a weak form of the so-called Tits alternative 
for groups of circle homeomorphisms. Although we do not give Margulis' original 
proof of this result, we develop an alternative one due to Ghys that is more 
suitable for a probabilistic interpretation in the context of random walks 
on groups. 

In Chapter 3, we collect several results of dynamical nature which require some degree of 
smoothness, in general $\mathrm{C}^2$. First, we extensively study the most important 
one, namely Denjoy's theorem. Next we  present some closely related results, as for 
example Sacksteder's and Duminy's theorem. Then we discuss two of the major open 
problems of the theory, namely the zero Lebesgue measure conjecture for minimal 
invariant Cantor sets, and the ergodicity conjecture for minimal actions. 
Finally, we treat the smoothness properties of 
conjugacies between group actions. At this point, it should have been 
natural to treat some particular properties in smoothness higher than $\ce^2$. 
Unfortunately, though there exist several interesting and promising works in 
this direction (see for example \cite{wilk,CC-sm,Ts-no}), it still seems 
to be impossible to collect them in a systematic and coherent way.

Chapter 4 corresponds to a tentative description of the dynamics of groups of 
diffeomorphisms of one-dimensional manifolds based on some relevant algebraic 
information. It begins with the case of Abelian and nilpotent groups, where Denjoy's  
theorem and the well-known Kopell's lemma become relevant. After a digression 
concerning ``growth'' for groups of diffeomorphisms, it continues with the 
case of polycyclic and solvable groups, concluding by showing the 
difficulties encountered in the case of amenable groups. 

Chapter 5 concerns obstructions to smooth actions for groups satisfying particular 
cohomological properties. It begins with the already classical Thurston's stability 
theorem for groups of $\ce^1$ diffeomorphisms of the interval. It then passes to 
a rigidity theorem for groups having the so-called Kazhdan Property (T) (roughly, 
Kazhdan group actions on the circle have finite image). The chapter concludes 
with a closely related super-rigidity result for actions of irreducible higher 
rank lattices. These last two theorems (due to the author) may be viewed as natural 
(but still non-definitive) generalizations of a series of results concerning obstructions 
for actions on one-dimensional spaces of higher-rank simple and semisimple Lie groups, 
in the spirit of the seminal works by Margulis and the quite inspiring  
Zimmer's program \cite{zimmer-chic,Zi}.

\vspace{0.2cm}

We have made an effort to make this book mostly self-contained. Although most of 
the results treated here are very recent, the techniques involved are, in general, 
elementary. We have also included a large list of complementary exercises  
where we pursue a little bit on some topics or briefly explain some related 
results. However, we must alert the reader that these ``exercises" may vary 
drastically in level of difficulty. Indeed, in many cases the small results which 
are presented in them do not appear in the literature. This is also the case of 
certain sections of the book. With no doubt, the most relevant case is that of \S 
\ref{super-du} where we give the original proof of a theorem proved by Duminy (more 
than 30 years ago) on the existence of infinitely many ends for semi-exceptional 
leaves of transversely $\ce^2$ codimension-one foliations. The pressing need 
to publishing Duminy's brilliant proof of this remarkable result (for 
which an alternative reference is \cite{CC2}) was an extra 
motivation for the author for writing this book.


\newpage
\thispagestyle{empty}
${}$
\newpage

\section*{Acknowledgments}

\addcontentsline{toc}{section}{Acknowledgments}


\vspace{0.15cm} 

\hspace{0.45cm} The original version (in Spanish) of this text was prepared for 
a mini-course in Antofagasta, Chile. Subsequently, enlarged and revised versions 
were published in the series {\em Monograf\'{\i}as del IMCA} (Per\'u) and 
{\em Ensaios Matem\'aticos} (Brasil). This translation arose from the necessity 
of making this text accessible to a larger audience. I would like to thank Juan 
Rivera-Letelier for his invitation to the II Workshop on Dynamical Systems (2001), 
for which the original version was prepared, and Roger Metzger for his invitation 
to IMCA (2006), where part of this material was presented. I would also like 
to thank both \'Etienne Ghys and Maria Eul\'alia Vares for motivating me 
and allowing me to publish the revised Spanish version in Brasil, as 
well as all the people who strongly encouraged me to conclude  
this English version. 

This work was partially supported by CONICYT and PBCT via the Research Network on Low Dimensional 
Dynamics. I would also to acknowledge the support of both the UMPA Department of the \'Ecole 
Normale Sup\'erieure de Lyon, where the idea of writing this text was born during my PhD 
thesis, and the Institut des Hautes \'Etudes Scientifiques, where the original notes 
started taking its definite form while I benefited from a one-year postdoctoral position.

\vspace{0.1cm}

This work owes a lot to many of my colleges. Several remarks spread throughout the text  
and the content of some examples and exercises were born during fruitful 
and quite stimulating discussions. It is then a pleasure to thank  
Sylvain Crovisier (Proposition \ref{facil}), 
Albert Fathi (Exercise \ref{fat}), 
Tsachik Gelander (Exercise \ref{ejer-gelander}), 
Adolfo Guillot (Remark \ref{adolfo}), 
Carlos Moreira (Exercise \ref{gugu}), 
Pierre Py (Remark \ref{pierre}), 
Takashi Tsuboi (Exercise \ref{no-anal-jap}), 
Dave Witte-Morris (Proposition \ref{unpuntofijo}), 
and Jean-Christophe Yoccoz (Exercises \ref{conjcorestable} and \ref{rec-frappend}). 
It is also a pleasure to thank Levon Beklaryan, 
Rostislav Grigorchuk, Leslie Jim\'enez, Daniel Jorquera, 
Yoshifumi Matsuda, Daniel Pons, Juan Rivera-Letelier, Pierre Paul Romagnoli, 
Eugenio Trucco, and Dave Witte-Morris, 
for their many corrections to this and earlier versions 
of the text, as well as Marina Flores for her help with the English translation. 

Finally, it is a very big pleasure to thank \'Etienne Ghys for his 
encouragements, and for numerous and fascinating discussions on the 
subject of this book and many other aspects of mathematics. 

\vspace{0.85cm}

\hfill Santiago, May 28, 2009

\newpage
\thispagestyle{empty}
${}$
\newpage


\chapter*{NOTATION AND GENERAL DEFINITIONS}

\addcontentsline{toc}{chapter}{NOTATION AND GENERAL DEFINITIONS}

\hspace{0.45cm} We 
will commonly denote the circle by $\clo$. As usual, we will consider the counterclockwise 
orientation on it. We will denote by $]a,b[$ the open interval from $a$ to $b$ according 
to this orientation. Notice that if $b$ belongs to $]a,c[$, then $c \! \in ]b,a[$ and 
$a \! \in ]c,b[$. We will sometimes write $a<b<c<a$ for these relations. One defines 
similarly the intervals $[a,b]$, $[a,b[$, and $]a,b]$. The distance between $a$ and 
$b$ is the shortest among the lengths of the intervals $]a,b[$ and $]b,a[$. We will denote 
this distance by $dist(a,b)$ or $|a-b|$. We will also use the notation $|I|$ for the length 
of an interval $I$ (on the circle or the real line). The Lebesgue measure of a measurable 
subset $A$ of either $\mathrm{S}^1$ or $\mathbb{R}$ will be denoted $Leb(A)$.

Although the contrary is said, along this text we will only deal with orientation preserving 
maps. The group of (orientation preserving) circle homeomorphisms will be denoted by 
$\mathrm{Homeo}_{+}(\mathrm{S}^1)$, and for $k \!\in\! \mathbb{N} \cup  \{ \infty \}$ we 
will denote the subgroup of $\ce^k$ diffeomorphisms by $\mathrm{Diff}^{k}_{+}(\mathrm{S}^1)$. 
For $\tau \! \in ]0,1[$ we will also deal with the group 
$\mathrm{Diff}_{+}^{1 + \tau}(\mathrm{S}^1)$ of circle 
diffeomorphisms whose derivatives are H\"older 
\index{H\"older!derivative} 
continuous of exponent $\tau$ (or just $\tau$-H\"older continuous), 
that is, for which there exists a constant $C > 0$ such that 
\esp $\big| f'(x)-f'(y) \big| \leq C |x-y|^{\tau}$ \esp 
for all $x,y$. 

The group of circle diffeomorphisms having Lipschitz derivative 
\index{Lipschitz!derivative}
will be denoted by $\mathrm{Diff}_{+}^{1 +\mathrm{Lip}}(\mathrm{S}^1)$. 
Finally, the notation $\mathrm{Diff}_{+}(\mathrm{S}^1)$ will be employed when 
the involved regularity is clear from the context or it is irrelevant.

We will sometimes see the real line as the universal covering of the 
circle through the maps $x \mapsto e^{ix}$ or $x \mapsto e^{2\pi ix}$, 
depending if we parameterize $\clo$ either by $[0,2\pi]$ or $[0,1]$. In general,  
we will consider the first of these parameterizations. By a slight abuse of notation, 
for each orientation preserving circle homeomorphism $f$ we will commonly denote 
also by $f$ (sometimes by $F$) each of its lifts to the real line. In this way, 
$f: \mathbb{R} \rightarrow \mathbb{R}$ will be an increasing continuous function 
such that, for every $x \!\in\! \R$, either \esp $f(x+2\pi) \!=\! f(x)+2\pi$ \esp  
or \esp $f(x+1) \!=\! f(x)+1$, \esp depending on the chosen parameterization. We will denote 
by $\widetilde{\mathrm{Homeo}_+}(\clo)$ the group of homeomorphisms of the line obtained 
as lifts of circle homeomorphisms. The circle rotation of angle $\theta$ will be 
\index{rotation}
denoted by $R_{\theta}$. Once again, notice that $\theta$ is an angle either in 
$[0,2\pi]$ or $[0,1]$, depending on the parameterization.\\

In many cases we will deal directly with subgroups of 
$\mathrm{Homeo}_{+}(\mathrm{S}^1)$ or $\mathrm{Diff}_{+}(\mathrm{S}^1)$. 
However, we will also consider representations of a group  
$\Gamma$ in the group of circle homeomorphisms or diffeomorphisms. In 
other words, we will deal with homomorphisms $\Phi$ from $\Gamma$ into 
$\mathrm{Homeo}_{+}(\mathrm{S}^1)$ or $\mathrm{Diff}_{+}(\mathrm{S}^1)$. 
These representations, which we may also thought of as actions, will be commonly denoted 
by $\Phi$. For simplicity, we will generally identify the element $g \in \Gamma$ 
with the map $\Phi(g)$. To avoid any confusion, we will denote the 
identity transformation by $Id$, and the neutral element in the underlying 
group by $id$.

\index{action!faithful} \index{action!free} 
Let us recall that, in general, an action $\Phi$ of a group $\Gamma$ on a space M is 
{\bf \em{faithful}} if for every $g \neq id$ the map $\Phi(g)$ is not the identity 
of $\mathrm{M}$. The action is {\bf \textit{free}} if for every $g \neq id$ one has 
$\Phi(g)(x) \neq x$ for all $x \in \mathrm{M}.$ The {\bf \em{orbit}} of a point 
$x \in \mathrm{M}$ is the set $\{\Phi(g)(x) \!: g \in \Gamma\}$. More generally, 
the orbit of a subset $A \subset \mathrm{M}$ is $\{\Phi(g)(x)\!: x \in A, \esp 
g \in \Gamma\}$. If M is endowed 
with a probability measure $\mu$, then $\Phi$ is {\bf {\em ergodic}} with 
\index{action!ergodic} 
respect to $\mu$ if every measurable set $A$ which is {\bf {\em invariant}} 
(that is, $\Phi(g) (A) \!=\! A$ for every $g \!\in\! \Gamma$) has either 
null or total $\mu$-measure. This is equivalent to the fact that every 
measurable function $\phi \!: \mathrm{M} \rightarrow \mathbb{R}$ satisfying 
$\phi \circ \Phi(g) = \phi$ for every $g \!\in\! \Gamma$ is $\mu$-almost surely 
constant. In the case the measure 
$\mu$ is invariant (that is, $\mu(\Phi(g)(A)) = \mu(A)$ 
for every measurable set $A$ and every $g \!\in\! \Gamma$), this is also equivalent 
to the fact that $\mu$ cannot be written as a nontrivial convex combination of two 
invariant, probability measures.

Two elements $f,g$ in $\mathrm{Homeo}_+(\clo)$ are {\bf \textit{topologically conjugate}} 
\index{conjugacy}
if there exists a circle homeomorphism $h$ such that $h \circ f = g \circ h$ (in this 
case, we say that $hfh^{-1}$ is the {\bf \textit{conjugate}} of $f$ by $h$). Similarly, 
two actions $\Phi_1$ and $\Phi_2$ of a group $\Gamma$ by circle homeomorphisms 
are topologically conjugate if there exists $h \!\in\! \mathrm{Homeo}_+(\clo)$ such that 
$h \circ \Phi_1(g) = \Phi_2(g) \circ h$ for every $g \in \Gamma$. 
In many cases, we will avoid the 
use of the symbol $\circ$ when composing maps. The non-Abelian free group on $n$ 
generators will be denoted by $\mathbb{F}_n$. 
\index{group!free} 
To avoid confusions, we will denote the one-dimensional torus by $\mathbb{T}^1$, 
thus emphasizing the group structure on it (which identifies with that of 
$(\mathbb{R} \hspace{0.15cm} \mathrm{mod} \hspace{0.15cm} 1,+)$, or 
$\mathrm{SO}(2,\mathbb{R})$). Similarly, we will sometimes write $(\mathbb{R},+)$ 
\esp (resp. $(\mathbb{Z},+)$) \esp instead of $\mathbb{R}$ (resp. $\mathbb{Z}$) 
to emphasize the corresponding additive group structure. 


\chapter{EXAMPLES OF GROUP ACTIONS ON THE CIRCLE}

\section{The Group of Rotations}

\hspace{0.45cm} The rotation group 
\index{rotation!group|(}
$\mathrm{SO}(2,\mathbb{R})$ is the simplest 
one acting transitively by circle homeomorphisms. Up to topological conjugacy, 
it may be characterized as the group of the homeomorphisms of $\clo$ that preserve 
a probability measure having properties similar to those of Lebesgue measure.  
Recall that the {\bf{\em support}} of a measure is the complement of 
the largest open set with null measure. 
\index{support!of a measure} 
Hence, the measure has {\bf \textit{total support}} if 
the measure of every non-empty open set is positive.

\vspace{0.1cm}

\begin{prop} {\em If $\Gamma$ is a subgroup of $\mathrm{Homeo}_{+}(\mathrm{S}^1)$ which   
preserves a probability measure having total support and no atoms, then $\Gamma$ is 
topologically conjugate to a subgroup of $\mathrm{SO}(2,\mathbb{R})$.}
\label{gatomedible}
\end{prop}

\noindent{\bf Proof.} The measure $\mu$ on $\clo$ given by the hypothesis 
induces in a natural way a $\sigma$-finite measure $\tilde{\mu}$ on the real 
line satisfying $\tilde{\mu}([x,x+2\pi]) \!=\! 1$ for all $x \in \R$. Let 
us define $\varphi: \mathbb{R} \rightarrow \mathbb{R}$ by letting 
$\varphi(x) = 2\pi\tilde{\mu} \big( [0,x] \big)$ if $x \geq 0$, and 
$\varphi(x) = -2 \pi \tilde{\mu} \big( [x,0] \big)$ if $x<0$. 
It is easy to check that $\varphi$ is a homeomorphism. For 
an arbitrary element $g \!\in\! \Gamma$, let us fix a lift  
$\tilde{g}$ such that $\tilde{g}(0)>0$. For $y > 0$ the 
point $\varphi \tilde{g} \varphi^{-1}(y)$ coincides with  
\begin{small}
$$2\pi \tilde{\mu} \big( [0,\tilde{g}\varphi^{-1}(y)] \big) = 
2\pi\tilde{\mu} \big( [0,\tilde{g}(0)] \big) + 
2\pi \tilde{\mu} \big([\tilde{g}(0),\tilde{g}\varphi^{-1}(y)] \big) 
= 2\pi \tilde{\mu} \big( [0,\tilde{g}(0)] \big) + 2\pi \tilde{\mu}
\big( [0,\varphi^{-1}(y)] \big).$$
\end{small}Therefore, 
$$\varphi \tilde{g}\varphi^{-1}(y)
= 2\pi\tilde{\mu} \big( [0,\tilde{g}(0)] \big) + y,$$ 
and a similar argument shows that the same holds for $y \leq 0$. 
Thus, $\varphi \tilde{g} \varphi^{-1}$ is the translation by 
$2 \pi \tilde{\mu} \big( [0,\tilde{g}(0)] \big)$. The map $\vf$ is 
$2\pi$-periodic, and hence induces a homeomorphism of $\clo$. The above 
computation then shows that, 
after conjugating by this homeomorphism, each $g \in \Gamma$ becomes the 
rotation by angle \hspace{0.1cm} $2 \pi \tilde{\mu} \big( [0,g(0)] \big) 
\hspace{0.2cm} \mathrm{mod} \hspace{0.2cm} 2\pi$. $\hfill\square$

\vspace{0.35cm}

Compact topological groups do satisfy the hypothesis of the preceding proposition. 
Indeed, if $\Gamma$ is a compact group, we may consider the (normalized) Haar 
measure $dg$ on it. If $\Gamma$ acts by circle homeomorphisms, we define a 
probability measure $\mu$ on the Borel sets of $\mathrm{S}^1$ by letting 
$$\mu(A) = \int_{\Gamma} Leb(gA) \hspace{0.1cm} dg.$$
This measure $\mu$ is invariant by $\Gamma$, and has total 
support and no atoms. We thus conclude the following. \index{group!compact}

\begin{prop} {\em Every compact subgroup of $\mathrm{Homeo}_{+}(\mathrm{S}^1)$ 
\index{rotation!group|)} 
is topologically conjugate to a subgroup of $\mathrm{SO}(2,\mathbb{R}).$}
\label{gatorot}
\end{prop}


\section{The Group of Translations and the Affine Group}

\hspace{0.45cm} Group actions on the circle with a global fixed point may be thought 
of as actions on the real line. Due to this, to understand general actions on $\clo$, 
we first need to understand actions on $\mathbb{R}$. 

An example of a nice group of homeomorphisms of the real line is the affine group 
\index{affine group|(}
$\mathrm{Aff}_{+}(\mathbb{R})$. To each element 
\esp $g(x) = ax + b$ \esp ($a > 0$) \esp in this group, we associate the matrix 
$$ \left(
\begin{array}
{cc}
a & b  \\
0 & 1  \\
\end{array}
\right) \in \mathrm{GL}_{+}(2,\mathbb{R}).$$
Via this correspondence, the group $\mathrm{Aff}_{+}(\mathbb{R})$ 
identifies with a subgroup of $\mathrm{GL}_+(2,\mathbb{R})$.

Recall that a {\bf \textit{Radon measure}} is a (nontrivial) \index{measure!Radon}
measure defined on the Borel sets of a topological space which 
is finite on compact sets. An example of a Radon measure is the Lebesgue measure. 
Since this measure is preserved up to a multiplicative constant by each element 
in \esp $\mathrm{Aff}_{+}(\mathbb{R})$, \esp the following definition becomes natural.

\vspace{0.04cm}

\begin{defn} Let $\upsilon$ be a Radon measure on the real line and $\Gamma$ 
a subgroup of $\mathrm{Homeo}_{+}(\mathbb{R})$. We say that $\upsilon$ is 
\index{measure!quasi-invariant}
{\bf \textit{quasi-invariant}} by $\Gamma$ if for every $g \in \Gamma$ 
there exists a positive real number $\kappa(g)$ such that 
$g^{*}(\upsilon) = \kappa(g) \cdot \upsilon$ \esp (that 
is, for each Borel set $A \subset \mathbb{R}$ one has 
$\upsilon \big( g(A) \big) = \kappa(g) \cdot \upsilon(A)$).
\end{defn}

\vspace{0.04cm}

As an analogue to Proposition \ref{gatomedible}, we have 
the following characterization of the affine group.

\vspace{0.05cm}

\begin{prop} {\em Let $\Gamma$ be a subgroup of $\mathrm{Homeo}_{+}(\mathbb{R})$. If 
there exists a Radon measure having total support and no atoms which is quasi-invariant 
by $\Gamma$, then $\Gamma$ is conjugate to a subgroup of the affine group.}
\label{conjugaralafin}
\end{prop}

\noindent{\bf Proof.} Let us define $\varphi \!: \mathbb{R} \rightarrow \mathbb{R}$ 
by letting $\varphi(x) \!=\! \upsilon \big( [0,x] \big)$ if $x \geq 0$, and  
$\varphi(x) \!=\! -\upsilon \big( [x,0] \big)$ if $x < 0$. If \esp $g \!\in\! \Gamma$ 
\esp and \esp $x \!\geq\! 0$ \esp are such that \esp $g(x) \!\geq\! 0$ 
\esp and \esp $g(0) \geq 0$, \esp then 
\begin{eqnarray*}
\varphi \big( g(x) \big) 
\!\!\!&=&\!\!\! \upsilon \big( [0,g(x)] \big) = \kappa(g) \upsilon \big( [g^{-1}(0),x] \big)\\
&=&\!\! \kappa(g) \upsilon \big( [0,x] \big) + \kappa(g) \upsilon \big( [g^{-1}(0),0] \big) 
\esp = \esp \kappa(g) \varphi(x) - \kappa(g) \varphi \big( g^{-1}(0) \big),
\end{eqnarray*}
and therefore,
$$\varphi g \varphi^{-1}(x) \!=\! 
\kappa(g) x - \kappa(g) \varphi \big( g^{-1}(0) \big).$$ 
Looking at all other cases, it is not difficult to check that the 
last equality actually holds for every $x\!\in\!\mathbb{R}$ and 
every $g\!\in\!\Gamma$, which concludes the proof. $\hfill\square$

\vspace{0.35cm}

The affine group contains the group of translations of the real line, and a similar 
argument to that of the proof of Proposition \ref{gatomedible} shows the following.

\vspace{0.05cm}

\begin{prop} {\em Let $\Gamma$ be a subgroup of $\mathrm{Homeo}_{+}(\mathbb{R})$. If 
\esp $\Gamma$ preserves a Radon measure having total support and no atoms, then \esp 
$\Gamma$ is topologically conjugate to a subgroup of the group of translations.} 
\label{conjugaratraslac}
\end{prop}

\begin{small} \begin{ejer} Define {\bf \textit{logarithmic derivative}} 
(of the derivative) of a $\ce^2$ diffeomorphism \esp 
\index{derivative!logarithmic}
$f\!\!: I \!\subset\! \mathbb{R} \!\rightarrow \! J \!\subset\!\mathbb{R}$ 
\esp as $LD ( f ) (x) \!=\!  ( \log(f') )' (x).$  Prove 
that $LD(f) \equiv 0$ if and only if $f$ is the restriction of an 
element of $\mathrm{Aff}_{+}(\mathbb{R})$. From the equality
\begin{equation}
\log ( (f \circ g)' ) (x) \esp = \esp \log(g')(x) + \log(f')(g(x)),
\label{cociclolog}
\end{equation}
deduce the cocycle relation
$$LD ( f \circ g ) (x) \esp = \esp LD(g)(x) + g'(x) \!\cdot\! LD(f) ( g(x) ).$$
\end{ejer} \end{small}
\index{affine group|)}\index{cocycle}


\section{The Group $\mathbf{\mathrm{PSL}(2,\mathbb{R})}$}

\subsection{$\mathbf{\mathrm{PSL}(2,\mathbb{R})}$ as the M\"obius group}

\hspace{0.45cm} We denote by \index{M\"obius group|(}
$\mathbb{D}$ the {\bf \textit{Poincar\'e disk}}, 
\index{Poincar\'e!disk}
that is, the unit disk endowed with the hyperbolic metric
$$\frac{4 du}{(1 - |u|^2)^2}, \quad u \in \mathbb{D}.$$
The group of non-necessarily orientation preserving diffeomorphisms of 
$\mathbb{D}$ which preserve this metric coincide with the group of conformal 
diffeomorphisms of $\mathbb{D}$, and it contains the M\"obius group as an index 
2 subgroup. This means that the only (orientation preserving) diffeomorphisms 
$g\!:\mathbb{D} \rightarrow \mathbb{D}$ for which the equality
$$\frac{2 \esp \! \| Dg(u) (\zeta) \|}{1-|g(u)|^2} \esp 
= \esp \frac{2 \esp \! \| \zeta \|}{1-|u|^2},$$
holds for every \esp $u\!\in\! \mathbb{D}$ \esp and every vector 
\esp $\zeta \!\in\! T_u (\mathbb{D}) \!\sim\! \mathbb{R}^2$ \esp 
are those which in complex notation may be written in the form 
$$g(z) = e^{i\theta} \cdot \frac{z-a}{1-\bar{a}z}, \quad \esp \theta \in [0,2\pi],
\esp \esp a \in \mathbb{C}, \esp \esp |a| < 1, \esp \esp z \in \mathbb{D}.$$
The M\"obius group is denoted by $\mathcal{M}$. Each 
one of its elements induce a real-analytic diffeomorphism of the circle 
(identified with the boundary $\partial \mathbb{D}$).

Let us now consider the map $\varphi (z) \!=\!
(z+i)/(1+iz).$ Notice that $\varphi(\mathrm{S}^1) 
\!=\! \mathbb{R} \cup \{ \infty \}$; moreover,  
the image of $\mathbb{D}$ by $\varphi$ is the upper half-plane 
$\mathbb{R}^2_+$, which endowed with the induced metric 
corresponds to the {\bf \textit{hyperbolic plane}} 
\index{hyperbolic!plane} $\mathbb{H}^2$. In complex notation, the 
action of each element of $\mathcal{M}$ on $\mathbb{H}^2$ is 
of the form $z \mapsto (a_1 z+a_2)/(a_3 z+a_4)$ for some 
$a_1,a_2,a_3,a_4$ in $\mathbb{R}$ such that $a_1 a_4-a_2 a_3=1$.
For $a_3 = 0$ we obtain an affine transformation, thus showing  
that $\mathrm{Aff}_{+}(\mathbb{R})$ is a subgroup of $\mathcal{M}$.

To each element $z \mapsto (a_1 z+a_2)/(a_3 z+a_4)$ 
in $\mathcal{M}$ we associate the matrix
\begin{equation}
\left(
\begin{array}
{cc}
a_1 & a_2  \\
a_3 & a_4  \\
\end{array}
\right) \in \mathrm{SL}(2,\mathbb{R}).
\label{matriz-psl2}
\end{equation}
An easy computation shows that the matrix associated to the composition 
of two elements in $\mathcal{M}$ corresponds to the product of the matrices 
associated to these elements. Moreover, two matrices $M_1,M_2$ in 
$\mathrm{SL}(2,\mathbb{R})$ induce the same element of $\mathcal{M}$ 
if and only if they coincide or $M_1 = - M_2$. In this way, the 
M\"obius group naturally identifies with the projective 
\index{projective!group|see{M\"obius group}} 
group $\mathrm{PSL}(2,\mathbb{R})$.

The action of $\mathrm{PSL}(2,\mathbb{R})$ satisfies a remarkable property of 
transitivity and rigidity: given two triples of cyclically ordered points on 
$\clo$, say $(a,b,c)$ and $(a',b',c')$, there exists a unique element 
$g \in \mathrm{PSL}(2,\mathbb{R})$ sending $a$, $b$, and $c$, into 
$a'$, $b'$, and $c'$, respectively. In particular, if $g$ fixes
three points, then $g = Id$.\\

The elements of $\mathcal{M} \sim \mathrm{PSL}(2,\mathbb{R})$ 
may be classified according to their  fixed points on 
$\mathrm{S^1} \subset \overline{\mathbb{D}}$. Remark 
that, to find these points in the upper half-plane 
model, we need to solve the equation
\begin{equation}
\frac{a_1 z+a_2}{a_3 z+a_4} = z.
\label{ecuaciondeclasificacion}
\end{equation}
A simple analysis shows that there are three cases:\\

\vspace{0.08cm}

\noindent{(i) $|a_1+a_4| < 2.$ In this case, the solutions to 
(\ref{ecuaciondeclasificacion}) are different (conjugate) points of the complex 
plane. Therefore, in the Poincar\'e disk model, the map $g$ has no fixed point 
on the unit circle. Actually, the map $g$ is conjugate to a rotation.}

\vspace{0.08cm}

\noindent{(ii) $|a_1+a_4|=2.$ In this case, the solutions of 
(\ref{ecuaciondeclasificacion}) coincide and are situated on the real line. 
Hence, in the Poincar\'e disk model, $g$ fixes a unique point on the circle.}

\vspace{0.08cm}

\noindent{(iii) $|a_1+a_4|>2.$ In this case, there exist two distinct solutions to 
(\ref{ecuaciondeclasificacion}), which are also on the real line. Therefore, the map 
$g$ fixes two points on the circle, one of them attracting and the other one repelling.}

\vspace{0.08cm}

Notice that $|a_1+a_4|$ corresponds to the absolute value of the trace of the 
corresponding matrix. (Though the function $M \mapsto a_1 \!+\! a_4$ is not 
well-defined on $\mathrm{PSL}(2,\mathbb{R})$, there is no ambiguity for the 
definition of its absolute value.) Figures 1, 2, and 3 below illustrate the 
cases (i), (ii), and (iii), respectively. In case (i) we say that the element 
is {\bf \textit{elliptic}}, in case (ii) it is {\bf \textit{parabolic}} if 
it does not coincide with the identity, and in case (iii) the map is 
{\bf \textit{hyperbolic}}.\\

\vspace{0.2cm}


\beginpicture

\setcoordinatesystem units <0.8cm,0.8cm>


\put{$a$} at 3.2  0
\put{$\bullet$} at 3.5 0

\put{$g$} at 0.08  0

\put{$g$} at 3.9  0

\put{Figure 1} at 2 -2.4


\plot
3.6522712    0.4
3.6881943    0.2
3.7          0 /

\plot
0.3477288    0.4
0.3118057    0.2
0.3          0 /

\plot
3.6522712    -0.4
3.6881943    -0.2
3.7          0 /

\plot
0.3477288    -0.4
0.3118057    -0.2
0.3          0 /

\plot
0.347   -0.4
0.4   -0.22 /

\plot
0.347    -0.4
0.22    -0.235 /

\plot
3.6522712    0.4
3.625        0.21 /

\plot
3.6522712    0.4
3.75        0.22 /


\plot
2             1.5
2.1           1.496663
2.3           1.4696938
2.5           1.4142136
2.7           1.3266499
2.9           1.2
3.0198039     1.1
3.1           1.0198039
3.2           0.9
3.3           0.74833148
3.3266499     0.7
3.4           0.53851648
3.4142136     0.5
3.4696938     0.3
3.496663      0.1
3.5           0 /

\plot
2             1.5
1.9           1.496663
1.7           1.4696938
1.5           1.4142136
1.3           1.3266499
1.1           1.2
0.9801961     1.1
0.9           1.0198039
0.8           0.9
0.7           0.74833148
0.6733501     0.7
0.6           0.53851648
0.5857864     0.5
0.5303062     0.3
0.503337      0.1
0.5           0 /

\plot
2             -1.5
2.1           -1.496663
2.3           -1.4696938
2.5           -1.4142136
2.7           -1.3266499
2.9           -1.2
3.0198039     -1.1
3.1           -1.0198039
3.2           -0.9
3.3           -0.74833148
3.3266499     -0.7
3.4           -0.53851648
3.4142136     -0.5
3.4696938     -0.3
3.496663      -0.1
3.5           0 /

\plot
2             -1.5
1.9           -1.496663
1.7           -1.4696938
1.5           -1.4142136
1.3           -1.3266499
1.1           -1.2
0.9801961     -1.1
0.9           -1.0198039
0.8           -0.9
0.7           -0.74833148
0.6733501     -0.7
0.6           -0.53851648
0.5857864     -0.5
0.5303062     -0.3
0.503337      -0.1
0.5           0 /




\put{Figure 2} at 8.3 -2.4


\plot
8.3            1.7
8.5          1.6881943
8.7          1.6522712
8.9          1.5905974 /

\plot
8.3            1.7
8.1          1.6881943
7.9          1.6522712 /

\plot
7.88     1.6522
8.07    1.8  /

\plot
7.88     1.6522
8.09     1.6 /

\plot
9.9522712    0.4
9.9881943    0.2
10          0 /

\plot
6.6477288    0.4
6.6118057    0.2
6.6          0 /

\plot
9.9522712    -0.4
9.9881943    -0.2
10          0 /

\plot
6.6477288    -0.4
6.6118057    -0.2
6.6          0 /

\plot
9.952   0.4
10.07    0.23 /

\plot
9.952   0.4
9.9     0.21 /

\plot
6.647    -0.4
6.7      -0.21 /

\plot
6.647     -0.4
6.5       -0.23 /


\plot
8.3             1.5
8.4           1.496663
8.6           1.4696938
8.8           1.4142136
9           1.3266499
9.2           1.2
9.3198039     1.1
9.4           1.0198039
9.5          0.9
9.6           0.74833148
9.6266499     0.7
9.7           0.53851648
9.7142136     0.5
9.7696938     0.3
9.796663      0.1
9.8           0 /

\plot
8.3             1.5
8.2           1.496663
8             1.4696938
7.8           1.4142136
7.6           1.3266499
7.4           1.2
7.2801961     1.1
7.2           1.0198039
7.1           0.9
7           0.74833148
6.9733501     0.7
6.9           0.53851648
6.8857864     0.5
6.8303062     0.3
6.803337      0.1
6.8           0 /

\plot
8.3           -1.5
8.4           -1.496663
8.6           -1.4696938
8.8           -1.4142136
9             -1.3266499
9.2           -1.2
9.3198039     -1.1
9.4           -1.0198039
9.5           -0.9
9.6           -0.74833148
9.6266499     -0.7
9.7           -0.53851648
9.7142136     -0.5
9.7696938     -0.3
9.796663      -0.1
9.8           0 /

\plot
8.3           -1.5
8.2           -1.496663
8             -1.4696938
7.8           -1.4142136
7.6           -1.3266499
7.4           -1.2
7.2801961     -1.1
7.2           -1.0198039
7.1           -0.9
7             -0.74833148
6.9733501     -0.7
6.9           -0.53851648
6.8857864     -0.5
6.8303062     -0.3
6.803337      -0.1
6.8           0 /

\put{$a$} at 8.1  -1.8
\put{$\bullet$} at 8.3 -1.5

\put{$g$} at 8.6 2

\put{$a$} at 14.25  -1.7

\put{$b$} at 14.38  1.87

\put{$g$} at 16.6 0

\put{$g$} at 12.6 0


\put{Figure 3} at 14.53 -2.4


\plot
16.2          0.6
16.2522712    0.4
16.2881943    0.2
16.3          0 /

\plot
13.0094026    0.6
12.9477288    0.4
12.9118057    0.2
12.9          0 /

\plot
16.2522712    -0.4
16.2881943    -0.2
16.3          0 /

\plot
12.9477288    -0.4
12.9118057    -0.2
12.9          0 /

\plot
16.2   0.6
16.13   0.42 /

\plot
16.2   0.6
16.34   0.42 /

\plot
13    0.6
13.07  0.42 /

\plot
13    0.6
12.86  0.42 /


\plot
14.6             1.5
14.7           1.496663
14.9           1.4696938
15.1           1.4142136
15.3           1.3266499
15.5           1.2
15.6198039     1.1
15.7           1.0198039
15.8           0.9
15.9           0.74833148
15.9266499     0.7
16             0.53851648
16.0142136     0.5
16.0696938     0.3
16.096663      0.1
16.1           0 /

\plot
14.6            1.5
14.5           1.496663
14.3           1.4696938
14.1           1.4142136
13.9          1.3266499
13.7           1.2
13.5801961     1.1
13.5           1.0198039
13.4           0.9
13.3           0.74833148
13.2733501     0.7
13.2           0.53851648
13.1857864     0.5
13.1303062     0.3
13.103337      0.1
13.1           0 /

\plot
14.6           -1.5
14.7           -1.496663
14.9           -1.4696938
15.1           -1.4142136
15.3           -1.3266499
15.5           -1.2
15.6198039     -1.1
15.7           -1.0198039
15.8           -0.9
15.9           -0.74833148
15.9266499     -0.7
16             -0.53851648
16.0142136     -0.5
16.0696938     -0.3
16.096663      -0.1
16.1           0 /

\plot
14.6           -1.5
14.5           -1.496663
14.3           -1.4696938
14.1           -1.4142136
13.9           -1.3266499
13.7           -1.2
13.5801961     -1.1
13.5           -1.0198039
13.4           -0.9
13.3           -0.74833148
13.2733501     -0.7
13.2           -0.53851648
13.1857864     -0.5
13.1303062     -0.3
13.103337      -0.1
13.1           0 /


\put{$\bullet$} at 14.6    -1.5

\put{$\bullet$} at 14.6 1.5


\endpicture

\vspace{0.3cm}

\begin{small} \begin{ejer} Prove that any two hyperbolic elements of 
$\mathrm{PSL}(2,\mathbb{R})$ are topologically conjugate. Show that 
the same holds for parabolic elements, but not for the elliptic ones. 
\end{ejer} \end{small}

There is another way to view the action of $\mathrm{PSL}(2,\mathrm{R})$ 
on the circle. For this, recall that the {\bf \textit{projective space}} 
\index{projective!space} $\mathbb{PR}^1$ 
is the set of lines through the origin in the plane. This space naturally identifies 
with the circle parameterized by the interval $[0,\pi]$, since every such a line 
is uniquely determined by the angle 
$\alpha \!\in \![0,\pi[$ which makes with the $x$-axis. Since a 
linear map sends lines into lines and fixes the origin, 
it induces a map from $\mathbb{PR}^1$ into itself; 
moreover, the map induced by two matrices in $\mathrm{GL}_+ (2,\mathbb{R})$ coincide 
if and only if these matrices represent the same element in $\mathrm{PSL}(2,\mathbb{R})$, 
{\em i.e.}, if any of them is a scalar multiple of the other one. 
Actually, if we consider $s \!=\! \mathrm{ctg} (\alpha)$ as a 
parameter, then the action of (\ref{matriz-psl2})
on $\mathbb{PR}^1$ is given by \esp $s \mapsto (a_1 s + a_2)/(a_3 s + a_4)$.
Although this view of the action of $\mathbb{RP}^1$ is perhaps simpler than the one 
arising from extensions of isometries of the Poincar\'e disk, the latter one is more 
suitable as a source of motivation: in a certain sense, most groups of circle 
diffeomorphisms tend to have a ``negative curvature'' behavior.

Despite the fact that the dynamics of each element in $\mathrm{PSL}(2,\mathbb{R})$
is very simple, the global structure of its subgroups is not uniquely determined 
by the individual dynamics. More precisely, for a subgroup $\Gamma$ of 
$\mathrm{Homeo}_{+}(\clo)$ the condition that each of its elements is topologically 
conjugate to an element of $\mathrm{PSL}(2,\mathbb{R})$ does not imply that 
$\Gamma$ itself is conjugate to a subgroup of $\mathrm{PSL}(2,\mathbb{R})$, 
even when the orbits are dense \cite{Kov} and the elements are real-analytic 
diffeomorphisms \cite{Na-kova}.\\

\vspace{0.16cm}

\begin{small} \begin{ejer} The 
{\bf \textit{Schwarzian derivative}}\index{Schwarzian derivative} 
$S(f)$ of a $\mathrm{C}^3$ diffeomorphism \esp $
f\!\!: I \subset \mathbb{R} \rightarrow J \subset \mathbb{R}$ 
\esp is defined by 
$$S(f) = \frac{f'''}{f'} - \frac{3}{2} \Big( \frac{f''}{f'} \Big)^2.$$
Show that $S(f)$ is identically zero if and only if $f$ is the restriction 
of a M\"obius transformation, {\em i.e.}, a map of the form \esp 
$x \mapsto (ax+b)/(cx+d)$. \esp Prove also the cocycle relation
\begin{equation} \index{cocycle}
S ( f \circ g ) (x) = S (g)(x) + ( g'(x) )^2 \cdot S (f) ( g(x) ).
\index{derivative!Schwarzian}
\label{cocicloschwarz}
\end{equation}
\end{ejer} \end{small}

\vspace{0.05cm}
 
\begin{small} \begin{ejer} Show the following formulae for 
the Schwarzian derivative of $\ce^3$ diffeomorphisms between 
intervals in the line:
\begin{equation}
S (g)(y) = 6 \lim\limits_{x \rightarrow y}
\left[ \frac{g'(x)g'(y)}{(g(x) - g(y))^2} - \frac{1}{(x-y)^2} \right] =
6 \lim\limits_{x \rightarrow y}
\frac{\partial^2}{\partial y \partial x} \log \left( \frac{g(x)-g(y)}{x-y} \right),
\label{cochina}
\end{equation}
\begin{equation}
-\frac{1}{2\sqrt{d g / d x}} \hspace{0.02cm} S (g) =
\frac{d^2}{dx^2} \Big( \frac{1}{\sqrt{dg/dx}} \Big).
\label{doscochina}
\end{equation}
\end{ejer} \end{small}

\vspace{0.05cm}

\begin{small}
\begin{ejer} A {\bf \em{projective structure}} on the circle is given by a system 
of local coordinates $\varphi_i \!: I_i \rightarrow \clo$ so that the change of 
coordinates $\varphi_j^{-1} \circ \varphi_i$ are restrictions of M\"obius 
transformations.
\index{projective!structure}

\vspace{-0.2cm}

\noindent (i) Check that the coordinates $\phi_1$ and 
$\phi_2$ with inverses $\alpha \mapsto \tan (\alpha)$ 
and $\alpha \mapsto \mathrm{ctg}(\alpha)$, respectively, 
define a projective structure on $\clo$ (to be referred 
to as the {\bf \em{canonical}} projective structure). 

\vspace{0.1cm}

\noindent (ii) Show that, given a projective structure on $\clo$, the Schwarzian 
derivative of a diffeomorphism $f \!: \clo \rightarrow \clo$ is well-defined as 
a {\bf \em{quadratic differential}}. In other words, given systems of coordinates 
$\varphi_1$, $\bar{\varphi}_1$, $\varphi_2$, and $\bar{\varphi}_2$, which are compatible 
with the prescribed projective structure, for every $x$ in the domain of $\varphi_1$ one has
$$S \big( \bar{\varphi}_2^{-1} \circ f \circ \bar{\varphi}_1 \big) (x) = 
\big((\bar{\varphi}_1^{-1} \circ \varphi_1)'\big)^2 \cdot 
S \big( \varphi_2^{-1} \circ f \circ \varphi_1 \big) 
(\bar{\varphi}_1^{-1} \circ \varphi_1 (x)).$$
\label{estructura-proyectiva}
\index{quadratic differential}
\end{ejer}
\end{small}


\subsection{$\mathbf{\mathrm{PSL}(2,\mathbb{R})}$ and the Liouville 
geodesic current}
\label{liouv}

\hspace{0.45cm} Recall that each geodesic in the Poincar\'e disk is uniquely 
determined by its endpoints on the circle, which are necessarily different. 
Hence, the space of geodesics may be naturally identified with the quotient space  
\esp $\clo \! \times \clo \setminus \Delta$ \esp under the equivalence relation 
which identifies the pairs $(s,t)$ and $(t,s)$, where $s \neq t$. A 
{\bf \textit{geodesic current}} \index{geodesic current|(} 
is a Radon \index{measure!Radon}
measure defined on this space of geodesics. It may be though as a measure 
$L$ defined on the Borel subsets of \esp $\clo \! \times \clo$ which are 
disjoint from the diagonal $\Delta$, which is finite on compact subsets 
of \esp $\mathrm{S}^1 \! \times \mathrm{S}^1 \setminus \Delta$, \esp 
and which satisfies the symmetry condition 
\begin{equation}
L \big( [a,b]\times [c,d] \big) =
L \big( [c,d]\times [a,b] \big), \qquad a<b<c<d<a.
\label{gcuno}
\end{equation}

\vspace{0.1cm}

\begin{prop} {\em The diagonal action of $\mathrm{PSL}(2,\mathbb{R})$ on 
\esp $\mathrm{S}^1 \!\times \mathrm{S}^1 \setminus \Delta$ \esp 
preserves the geodesic current}
$$Lv= \frac{ds \hspace{0.1cm} dt}{4 \hspace{0.02cm} \mathrm{sin}^2
\left( \frac{s-t}{2} \right) }.$$
\end{prop}

\noindent{\bf Proof.} Let us first recall that the cross-ratio of four points 
$e^{ia},e^{ib},e^{ic},e^{id}$ in $\clo$ is defined as 
$$[e^{ia},e^{ib},e^{ic},e^{id}] = 
\frac{(e^{ia}-e^{ic})(e^{ib}-e^{id})}{(e^{ia}-e^{id})(e^{ic}-e^{ib})}.$$
One easily checks that cross-ratios are invariant under M\"obius transformations. 
Conversely, if a circle homeomorphism preserves cross-ratios, then it belongs 
to the M\"obius group. Now notice that, for $a<b<c<d<a$, the measure 
$Lv([a,b]\times[c,d])$ equals
\begin{eqnarray*}
\int_{c}^{d}\!\!\!\int_{a}^{b} \frac{ds \hspace{0.15cm}
dt}{4\mathrm{sin}^2(\frac{s-t}{2})}
&=& \int_{c}^{d} \left[ -
\frac{\mathrm{cos}(\frac{s-t}{2})}{2\mathrm{sin}(\frac{s-t}{2})}
\right]_{s=a}^{s=b}dt\\
&=& \int_{c}^{d} \frac{1}{2} \left[
\mathrm{cot} \left( \frac{a-t}{2} \right) - \mathrm{cot} \left(
\frac{b-t}{2}
\right) \right]dt\\
&=& \log \left( \Big|
\frac{\mathrm{sin}(\frac{b-d}{2})
\mathrm{sin}(\frac{a-c}{2})}{\mathrm{sin}(\frac{b-c}{2})\mathrm{sin}(\frac{a-d}{2})}
\Big| \right).
\end{eqnarray*}
Since
$\left| \mathrm{sin} (\frac{x-y}{2}) \right| = \frac{|e^{ix}-e^{iy}|}{2},$
this yields
$$Lv \big( [a,b]\times[c,d] \big)
= \log \big( \big| [e^{ia},e^{ib},e^{ic},e^{id}] \big| \big) 
= \log \big( [e^{ia},e^{ib},e^{ic},e^{id}] \big),$$
where the last inequality follows from the fact that the cross-ratio of 
cyclically ordered points on the circle is a {\em positive} real number. 
Since the action of $\mathrm{PSL}(2,\mathbb{R})$ on $\mathrm{S}^1$ 
preserves cross-ratios, it also preserves the measure $Lv$. 
$\hfill\square$

\vspace{0.35cm}

The measure $Lv$, called {\bf \textit{Liouville measure}}, satisfies the equality
\begin{equation}
e^{-Lv([a,b] \times [c,d])} + e^{-Lv([b,c]
\times [d,a])} = 1 \label{gcdos}
\end{equation}
for all $a < b < c < d < a.$ Indeed, 
\begin{eqnarray*}
e^{-Lv([a,b] \times [c,d])} + e^{-Lv([b,c] \times [d,a])} 
\! &=& \! 
\frac{1}{[e^{ia},e^{ib},e^{ic},e^{id}]} + \frac{1}{[e^{ib},e^{ic},e^{id},e^{ia}]}\\
\! &=& \! \frac{(e^{ia}-e^{id})(e^{ib}-e^{ic})}{(e^{ia}-e^{ic})(e^{ib}-e^{id})} +
          \frac{(e^{ib}-e^{ia})(e^{ic}-e^{id})}{(e^{ib}-e^{id})(e^{ic}-e^{ia})}\\
\! &=& \! \frac{(e^{ia}-e^{id})(e^{ib}-e^{ic}) -
  (e^{ib}-e^{ia})(e^{ic}-e^{id})}{(e^{ia}-e^{ic})(e^{ib}-e^{id})} 
\esp \esp = \esp \esp 1.
\end{eqnarray*}
As we will next see, property (\ref{gcdos}) 
characterizes Liouville measure. \index{measure!Liouville} 

Given a geodesic current $L$, we denote the group of circle homeomorphisms 
preserving $L$ by $\Gamma_{L}$. For instance, it is easy to check that 
$\Gamma_{Lv} = \mathrm{PSL}(2,\mathbb{R})$. In general, the group $\Gamma_{L}$ 
is very small (it is ``generically'' trivial). However, there exists a very 
simple condition which ensures that it is topologically conjugate to 
$\mathrm{PSL}(2,\mathbb{R})$. The following result may be considered as an 
analogous for the M\"obius group of Propositions \ref{gatomedible}, 
\ref{conjugaralafin}, or \ref{conjugaratraslac}. We refer to \cite{Bo} 
for its proof (see also Exercise \ref{conjcorestable}).

\vspace{0.18cm}

\begin{prop} {\em If $L$ is a geodesic current verifying property 
{\em (\ref{gcdos})}, then $\Gamma_L$ is conjugate to 
$\mathrm{PSL}(2,\mathbb{R})$ by a homeomorphism 
sending $L$ into $Lv$.}
\end{prop}

\vspace{0.02cm}

\begin{small} \begin{ejer} Using the invariance of the Liouville geodesic current 
under M\"obius transformations, show that if $f\!: I \rightarrow \mathbb{R}$ is a 
$\C^1$ local diffeomorphism satisfying  
$$f'(x)f'(y) = \frac{(f(x)-f(y))^2}{(x-y)^2}$$
for all $x \neq y$ in $I$, then $f$ is of the form \esp 
$x \mapsto (ax+b)/(cx+d)$ \esp (compare (\ref{cochina})).
\end{ejer} \end{small}
\index{geodesic current|)}


\subsection{$\mathbf{\mathrm{PSL}(2,\mathbb{R})}$ and the convergence property}

\hspace{0.45cm} A sequence $(g_n)$ of circle homeomorphisms has the {\bf \textit{convergence 
property}}\index{convergence property|(} 
if it contains a subsequence $(g_{n_k})$ satisfying one of the following properties:

\vspace{0.08cm}

\noindent{(i) there exist $a,b$ in $\mathrm{S}^1$ (not necessarily different) such 
that $g_{n_k}$ converges punctually to $b$ on $\mathrm{S}^1 - \{ a \}$, \esp and 
\esp $g_{n_k}^{-1}$ converges punctually to $a$ on $\clo \setminus \{ b \}$;}

\vspace{0.08cm}

\noindent{(ii) there exists $g \in \mathrm{Homeo}_{+}(\mathrm{S}^1)$ such that $g_{n_k}$ 
converges to $g$ and $g_{n_k}^{-1}$ converges to $g^{-1}$ on the circle.}

\vspace{0.08cm}

\noindent{A subgroup $\Gamma$ of $\mathrm{Homeo}_{+}(\mathrm{S}^1)$ has the convergence 
property if every sequence of elements in $\Gamma$ satisfies the property above. 
Notice that the convergence property is invariant under topological conjugacy. 

One easily checks that every subgroup of $\mathrm{PSL}(2,\mathbb{R})$ has the convergence  
property (see Exercises \ref{quasisim} and \ref{conjcorestable}). Conversely, a difficult 
theorem due to Casson-Jungreis, Gabai, Hinkkanen, and Tukia \cite{CJ,Ga,GM,Hi1,Tu} 
asserts that this property characterizes (up to topological conjugacy) the 
subgroups of $\mathrm{PSL}(2,\mathbb{R})$.

\vspace{0.1cm}

\begin{thm} {\em A group of circle homeomorphisms is topologically 
conjugate to a subgroup of $\mathrm{PSL}(2,\mathbb{R})$ 
if and only if it satisfies the convergence property.}
\end{thm}

It is not difficult to show that, 
for discrete subgroups of $\mathrm{Homeo}_{+}(\mathrm{S}^1)$, 
the convergence property is equivalent to that the action on the 
space of ordered triples of points in $\clo$ is free 
\index{action!free}
and properly discontinuous \cite{Tu}.

\vspace{0.1cm}

\begin{small} \begin{ejer} Prove directly from the definition that if $\Gamma$ 
satisfies the convergence property and $g \!\in\! \Gamma$ fixes three points on 
$\clo$, then $g \!=\! Id$.
\end{ejer}

\begin{ejer} A circle homeomorphism $g$ is {\bf \textit{$C$-quasisymmetric}} 
\index{quasisymmetric homeomorphism} 
if for all \esp $a\!<\!b\!<\!c\!<\!d\!<\!a$ \esp for which \esp $[a,b,c,d]\!=\!2$ 
\esp one has \esp $1/C \!\leq\! [g(a),g(b),g(c),g(d)] \!\leq\! C.$ \esp Prove that 
if $\Gamma$ is a {\bf \em{uniformly quasisymmetric}} 
subgroup of $\mathrm{Homeo}_{+}(\clo)$ 
(that is, all of its elements are $C$-quasisymmetric with respect to the same 
constant $C$), then $\Gamma$ satisfies the convergence property. Conclude that 
$\Gamma$ is topologically conjugate to a subgroup of $\mathrm{PSL}(2,\mathbb{R})$.

\noindent{\underbar{Remark.}} According to a difficult result due to Markovic \cite{mark}, 
under the preceding hypothesis the group $\Gamma$ is conjugate to a subgroup of 
$\mathrm{PSL}(2,\mathbb{R})$ by a {\em quasisymmetric} homeomorphism.
\label{quasisim}
\end{ejer}

\begin{ejer} Let $L$ be a geodesic current satisfying  
$L ([a,a] \times [b,c]) = 0$ for all $a < b \leq c < a$ and $L
\left( [a,b[ \times ]b,c] \right) = \infty$ for all $a < b < c <
a$ (notice that Liouville measure \index{geodesic current}
satisfies these properties). Prove that $\Gamma_{L}$ has the convergence 
property (see also \cite{Na-cociclo}). Conclude that $\Gamma_{L}$ is 
topologically conjugate to a subgroup of $\mathrm{PSL}(2,\mathbb{R})$. 
Remark, however, that this conjugacy does not necessarily send $L$ to 
the Liouville geodesic current, since property (\ref{gcdos}) is 
invariant under conjugacy, and is not necessarily satisfied 
\index{M\"obius group|)} by the elements of $\Gamma_{L}$. 
\label{conjcorestable}
\end{ejer} \end{small} 
\index{convergence property|)}


\section{Actions of Lie Groups}
\label{Lie}

\hspace{0.45cm} The 
\index{Lie group} \index{group!Lie|see{Lie group}}object 
of this section is to show that, among locally compact groups 
which do act by circle homeomorphisms, those which may provide new phenomena 
are the discrete ones (that we thought of as zero-dimensional Lie groups). This is 
the reason why we will mostly consider actions on the circle of {\em discrete} groups. 

Recall that a deep (and already classical) result by Montgomery and Zippin \cite{MZ} 
asserts that a locally compact topological group is a Lie group if and only if it 
has no ``small compact subgroups'', that is, there exists a neighborhood of the 
identity which does contain no nontrivial compact subgroups.
\index{Montgomery-Zippin's theorem}

\vspace{0.1cm}

\begin{prop} {\em Every locally compact subgroup of 
$\mathrm{Homeo}_{+}(\mathrm{S}^1)$ is a Lie group.}
\label{localmente compacto}
\end{prop}

\noindent{\bf Proof.} Let $\Gamma$ be a locally compact subgroup of 
$\mathrm{Homeo}_{+}(\mathrm{S}^1)$. The set 
\index{group!compact} 
$$V= \left\{ g \in \Gamma: \hspace{0.1cm} dist(x,g(x)) < 2\pi/3 \right\}$$  
is a neighborhood of the identity in $\Gamma$. We will show that 
$V$ contains no nontrivial compact subgroups, which due to 
Montgomery-Zippin's theorem implies the proposition.

Let $\Gamma_0$ be a compact subgroup of $\mathrm{Homeo}_{+}(\clo)$ contained in $V$. 
For each $f \in \Gamma$ let $\tilde{f} \in \widetilde{\mathrm{Homeo}_+}(\mathrm{S}^1)$ 
be the (unique) lift of $f$ such that $dist(\tilde{f}(x), x) < 2\pi / 3$ for all 
$x \in \mathbb{R}$. One readily checks that $\widetilde{gh} = \tilde{g} \tilde{h}$ 
for all $g,h$ in $\Gamma_0$. Therefore, $\Gamma_0$ embeds into the group 
$\widetilde{\mathrm{Homeo}_{+}}(\mathrm{S}^1).$ Now 
$\widetilde{\mathrm{Homeo}_{+}}(\mathrm{S}^1)$ is torsion-free, though 
Proposition \ref{gatorot} implies that every nontrivial compact subgroup 
of $\mathrm{Homeo}_{+}(\mathrm{S}^1)$ has torsion elements. This implies 
that $\Gamma_0$ must be trivial. $\hfill\square$

\vspace{0.15cm}

\begin{small} \begin{ejer} In order to avoid referring to Proposition \ref{gatorot} in 
the preceding proof, prove the following lemma due to Newman \cite{new} (see also \cite{kolev}): 
If $f$ is a  nontrivial finite order homeomorphism of the sphere $\mathrm{S}^n$ (normalized so that 
its diameter is $1$), then there exists $i \!\in\! \mathbb{N}$ such that \esp $dist(f^i,Id) > 1/2$.

\noindent{\underbar{Hint.}} If the opposite inequality holds for all $i$, then each orbit 
is contained in a hemisphere. We may then define a continuous map \esp $bar\!: \mathrm{S}^n 
\rightarrow \mathrm{S}^n$ \esp by associating to $x$ the ``barycenter'' $bar(x)$ of its orbit 
inside the corresponding hemisphere. This map satisfies $bar(f(x)) = bar(x)$ for all $x$, 
from where one easily deduces that the order of $f$ divides its topological degree. 
Nevertheless, since $f$ is homotopic to the identity, its topological degree is 
equal to $1$.
\index{degree!of a map}

\noindent{\underbar{Remark.}} From the above argument one easily 
deduces that \esp $dist(f,Id) > 1 / 2k$, \esp where $k$ is the order of $f$.
\end{ejer} \end{small} 

\vspace{0.12cm}

It is not very difficult to obtain the classification of transitive actions 
of {\em connected} \esp Lie groups on one-dimensional manifolds: see \cite{Gh1}. 
Up to topological conjugacy, the complete list consists of:

\vspace{0.1cm}

\noindent{(i) the action of $(\mathbb{R},+)$ by translations of the line;}

\vspace{0.1cm}

\noindent{(ii) the action of the rotation group $\mathrm{SO}(2,\mathbb{R})$ 
\index{rotation!group} on the circle;}

\vspace{0.1cm}

\noindent{(iii) the action of the affine group \index{affine group}
$\mathrm{Aff}_{+}(\mathbb{R})$ on the line;}

\vspace{0.1cm}

\noindent{(iv) the action of the group $\mathrm{PSL}_k(2,\mathbb{R})$ whose elements 
are the lifts of the elements of $\mathrm{PSL}(2,\R)$ to the $k$-fold covering 
of $\clo$ (which is topologically a circle);}

\vspace{0.1cm}

\noindent{(v) the action of the group $\widetilde{\mathrm{PSL}}(2,\mathbb{R})$ whose 
elements are the lifts to the real line of the elements of $\mathrm{PSL}(2,\mathbb{R})$.}

\vspace{0.1cm}

In a certain sense, this classification says that there exist only three distinct 
types of geometries on one-dimensional manifolds: Euclidean, affine, and 
projective (compare \cite{mark-nuevo}). The classification of (faithful) 
non transitive actions of connected Lie groups follows from the 
preceding one. Indeed, the orbits of such an action correspond to points 
or whole intervals. Therefore, denoting by $\mathrm{Fix}(\Gamma)$ the set of 
global fixed points of the action, on each connected component of the complement 
of $\mathrm{Fix}(\Gamma)$ we obtain an action given by a surjective homomorphism 
from $\Gamma$ into $(\mathbb{R},+)$, $\mathrm{SO}(2,\mathbb{R})$,
$\mathrm{Aff}_{+}(\mathbb{R})$, $\mathrm{PSL}_k(2,\mathbb{R})$, 
or $\widetilde{\mathrm{PSL}}(2,\mathbb{R})$. \index{M\"obius group} 


\section{Thompson's Groups\index{Thompson's groups|(}}

\hspace{0.45cm} For simplicity, in this section we will use the parameterization 
of the circle by the interval $[0,1]$ via the map $x \mapsto e^{2\pi ix}$. 
Let us consider the group of the homeomorphisms \index{piecewise affine homeomorphism}
$\tilde{f}: \mathbb{R} \rightarrow \mathbb{R}$ satisfying:\\

\vspace{0.1cm}

\noindent{(i) $\tilde{f}(0) = 0$,}

\vspace{0.1cm}

\noindent{(ii) there exists a sequence $\ldots x_{-1} < x_0 < x_1 < \ldots$ 
(diverging in both directions) of dyadic rational numbers such that 
\index{dyadic!rational number}
each of the restrictions $\tilde{f}|_{[x_i,x_{i+1}]}$ is affine with 
derivative an integer power of $2$,}

\vspace{0.1cm}

\noindent{(iii) $\tilde{f}(x+1) = \tilde{f}(x) + 1$ for all $x \in
\mathbb{R}$.}

\vspace{0.1cm}

Each such an $\tilde{f}$ induces a homeomorphism $f$ of 
$[0,1]$ by letting $f(0) \!=\! 0$, $f(1) \!=\! 1$, and 
$f(s) \!=\! \tilde{f}(s)$ $\mathrm{mod}$ $1$ \esp for 
$s \in ]0,1[$. We obtain in this way a group of homeomorphisms of $[0,1]$. This 
group was first introduced by Thompson and is commonly denoted by $\mathrm{F}$.\\

\vspace{0.64cm}


\beginpicture

\setcoordinatesystem units <0.6cm,0.6cm>

\putrule from 0 0 to 6 0

\putrule from 0 0  to 0 6

\putrule from 6 0 to 6 6

\putrule from 0 6 to 6 6

\put{$f$} at 4 2
\put{$g$} at 14.1 3.1

\put{Figure 4} at 2.85 -0.6

\putrule from 10 0 to 16 0
\putrule from 10 6 to 16 6
\putrule from 10 0 to 10 6
\putrule from 16 0 to 16 6

\plot
0  0
3  1.5 /

\plot
3   1.5
4.5 3 /

\plot
4.5 3
6   6 /

\plot
10 0
13 3 /

\put {Figure 5} at 12.98 -0.6

\plot
13   3
14.5 3.75 /

\plot
14.5  3.75
15.25 4.5 /

\plot
15.25 4.5
16  6 /

\setdots

\putrule from 3 0 to 3 6
\putrule from 4.5 0  to 4.5 6
\putrule from 0 1.5  to 6 1.5
\putrule from 0 3 to 6 3

\putrule from 13 0 to 13 6
\putrule from 14.5 0 to 14.5 6
\putrule from 15.25 0 to 15.25 6
\putrule from 10 3 to 16 3
\putrule from 10 3.75 to 16 3.75
\putrule from 10 4.5 to 16 4.5

\put{} at -3.3 0

\endpicture


\vspace{0.65cm}

Thompson's group F has several remarkable properties, which are not always 
easy to prove. First of all, $\mathrm{F}$ admits the finite presentation 
$$\mathrm{F} = \big\langle f,g \hspace{0.05cm}: \hspace{0.2cm}
[fg^{-1},f^{-1}gf] = [fg^{-1},f^{-2}gf^2] = id \big\rangle,$$
where $[\cdot,\cdot]$ denotes the commutator between two elements, and $f,g$ 
are the homeomorphisms whose graphs are depicted in Figures 4 and 5, respectively.

Since every nontrivial homeomorphism of the interval has infinite order, $\mathrm{F}$ 
is torsion-free. \index{group!torsion-free} From the proof of Theorem \ref{arefad} 
later on, it will appear as evident that $\mathrm{F}$ has free Abelian 
subgroups of infinite index. Moreover, the abelianized quotient 
$\mathrm{F}/[\mathrm{F},\mathrm{F}]$ is isomorphic to 
$\mathbb{Z} \times \mathbb{Z}$. To show this it suffices to take, for each 
$[h] \in \mathrm{F}/[\mathrm{F},\mathrm{F}]$, the value of the derivative 
of $h$ at the endpoints (notice that these value do not depend on the 
representative of $h$, since $[f,g]'(0) \!=\! [f,g]'(1) \!=\! 1$ for all $f,g$ in 
$\mathrm{F}$). Then taking the logarithm in base $2$ of these values, we obtain an 
homomorphism from $\mathrm{F}/[\mathrm{F},\mathrm{F}]$ into $\mathbb{Z} \times \mathbb{Z}$. 
Actually, it can be easily checked that this homomorphism is an isomorphism. 
(One may also use the nontrivial fact that the derived group 
$[\mathrm{F},\mathrm{F}]$ is simple \cite{CFP}.) 
\index{group!simple}

Further information concerning Thompson's group $\mathrm{F}$ may be found for instance in 
\cite{brin-ubi,CFP,GS}. In particular, in the former reference one may find a discussion 
of the relevant problem of the 
{\em amenability} of $\mathrm{F}$ (see Appendix B for the notion 
of amenability for groups). According to Exercise \ref{librenomoy}, one of the main 
difficulties for this is the fact that $\mathrm{F}$ does not contain free subgroups 
on two generators (see however Exercise \ref{no-ley}). 
This is a corollary of a much more general and nice result due to Brin and Squier 
\cite{BS} which we reproduce below. We remark that if $\mathrm{F}$ is non-amenable, 
\index{group!amenable} 
then this would lead to the first example of a finitely presented, {\em torsion-free}, 
non-amenable group which does not contain $\mathbb{F}_2$. (We point out that a 
finitely presented, non-amenable group not containing $\mathbb{F}_2$ but having 
torsion has been constructed by Olshanski and Sapir in \cite{OS}.)

\vspace{0.15cm}

\begin{thm} {\em The group $\mathrm{PAff}_+\big( [0,1] \big)$ of piecewise 
affine homeomorphisms of $[0,1]$ does not contain free subgroups on two generators.}
\label{arefad}
\end{thm}

\noindent{\bf Proof.} Suppose for a contradiction that $f,g$ in 
$\mathrm{PAff}_+\big( [0,1] \big)$ generate a free group. 
For each $h \!\in\! \mathrm{PAff}_+\big( [0,1] \big)$ let us 
denote by $supp_0(h)$ the {\em open support} 
\index{support!of a map} 
of $h$, that is, the set of points in $[0,1]$ which are not fixed by $h$. 
The set $I\!=\! supp_0(f) \!\cup\! supp_0(g)$ may be written as the union 
of finitely many open intervals $I_1, \ldots, I_n$. Notice that the 
closure of the set $supp_0([f,g])$ is contained in $I$, since in a 
neighborhood of each endpoint $x_0$ of each $I_i$ the maps $f$ and 
$g$ are of the form \esp $x \!\mapsto\! \lambda (x\!-\!x_0) \!+\! x_0$, 
\esp and hence commute.

Among the nontrivial elements $h \!\in\! \langle f,g \rangle$ such that 
$\overline{supp_0(h)}$ is contained in $I$, let us choose one, say $h_0$, 
such that the number of connected components of $I$ which do intersect 
$\overline{supp_0(h_0)}$ is minimal. Let $]a,b[$ be one of the connected 
components of the intersection, and let $[c,d]$ be an interval contained in 
the interior of $]a,b[$ and which contains $supp_0(h_0) \esp \cap \esp ]a,b[$. 
If $x$ belongs to $]a,b[$, then the orbit of $x$ by the group generated by $f$ 
and $g$ is contained in $]a,b[$; moreover, 
the supremum of this orbit is a point which 
is fixed by $f$ and $g$, and hence 
coincides with $b$. One then deduces the existence 
of an element $\bar{h} \! \in \!\langle f,g \rangle$ sending the interval $[c,d]$ 
to the right of $d$. In particular, the restrictions of $h_0$ and 
$\bar{h} h_0 \bar{h}^{-1}$ to $[a,b]$ commute, and they generate 
a subgroup isomorphic to $\mathbb{Z} \times \mathbb{Z}$. On the 
other hand, $h_0$ and $\bar{h} h_0 \bar{h}^{-1}$ do not commute in 
$\langle f,g \rangle \!\sim \mathbb{F}_2$, since otherwise they 
would generate a subgroup isomorphic to $\mathbb{Z}$. The commutator 
between $h_0$ and $\bar{h} h_0 \bar{h}^{-1}$ is then a nontrivial element 
whose open support does not intersect $]a,b[$, and so it intersects fewer 
components of $I$ than the open support of $h_0$. However, this contradicts 
the choice of $h_0$. $\hfill\square$

\vspace{0.35cm}

If we consider the homeomorphisms $\tilde{f}$ of the real line which satisfy 
only the properties (ii) and (iii) corresponding to the lifts of the 
elements in $\mathrm{F}$, and such that $\tilde{f}(0)$ is a dyadic rational 
number, then after passing to the quotient we obtain a group G of circle 
homeomorphisms. This group is infinite, has a finite presentation, 
and is simple. \index{group!simple} In fact, $\mathrm{G}$ was the first 
example of a group satisfying these three properties simultaneously. 
(This was one of the original motivations which leaded Thompson to 
introduce these groups.)


\subsection{Thurston's piecewise projective realization}
\label{sec-th}

\hspace{0.45cm} In order to understand Thompson's groups better, we will 
give two alternative definitions in this section. One is based on Thompson's 
original work, and the other follows an idea due to Thurston. 
\index{Thurston!realization of Thompson's groups}
We begin with some definitions.

\vspace{0.1cm}

A (non-degenerate) {\bf \textit{dyadic tree}} 
\index{tree} \index{dyadic!tree}
$\mathcal{T}$ is a finite union of closed 
{\bf \textit{edges}} (that is, homeomorphic copies of the unite interval, 
including its endpoints or {\bf \textit{vertexes}}) such that:

\vspace{0.1cm}

\noindent{(i) there exists a  marked vertex, called the {\bf \textit{root}} 
of the tree and denoted by $\sigma$;}

\vspace{0.1cm}

\noindent{(ii) each vertex different from the root is the final point of either 
one or three edges, while the root is the final point of either two edges or no edge 
(the latter case only appears when the tree is degenerate: it contains no edge  
and is reduced to a single point, namely the root);}

\vspace{0.1cm}

\noindent{(iii) $\mathcal{T}$ is connected.}

\vspace{0.1cm}

If a vertex is the final point of three edges, then they may be labeled $\Upsilon_{d}$, 
$\Upsilon_{l}$, and $\Upsilon_{r}$, according as they point down, to the left, or to the 
right, respectively. The edges starting from $\sigma$ are labeled $\Upsilon_l$ and 
$\Upsilon_r$. A {\bf \textit{leaf}} of the tree is a vertex $v$ which is the final 
point of a single edge. The root $\sigma$ will also be considered as a leaf of the 
degenerate tree. The set of all the leaves of a dyadic tree $\mathcal{T}$ will be 
denoted by $lv (\mathcal{T})$. Notice that there exists a natural cyclic order for 
the leaves of a dyadic tree, and the notion of {\bf \textit{first leaf}} may 
also be defined naturally.

Given a dyadic tree $\mathcal{T}$ and a leaf $p \!\in\! lv(\mathcal{T})$, we will say 
that a dyadic tree $\mathcal{T}'$ ``germinates'' from $\mathcal{T}$ at the leaf $p$ if 
$\mathcal{T}'$ is the union of $\mathcal{T}$ and two edges starting from $p$. 
Notice that the number of leaves of the new tree equals that of the original one plus $1$.

Let us now consider two trees $\mathcal{T}_1$ and $\mathcal{T}_2$ having the same number 
of leaves. We will say that a map from $lv(\mathcal{T}_1)$ to $lv(\mathcal{T}_2)$ is 
G-admissible if it preserves the cyclic order between the leaves, and we will say 
that it is F-admissible if it moreover sends the first leaf of $\mathcal{T}_1$ into the 
first leaf of $\mathcal{T}_2$. We will define an equivalence relation for G-admissible 
maps, leaving the task of adapting the definition to the case of F-admissible maps.

Let us consider a G-admissible map $g : lv (\mathcal{T}_1) \rightarrow lv (\mathcal{T}_2)$
and a leaf $p$ of $\mathcal{T}_1$. Let $\mathcal{T}_1'$ and $\mathcal{T}_2'$ be 
dyadic trees germinating from $\mathcal{T}_1$ and $\mathcal{T}_2$ 
at $p$ and $g(p)$, respectively. Let us define the map 
\esp $g': lv (\mathcal{T}_1') \rightarrow lv (\mathcal{T}_2')$ \esp 
by letting $g'(q) = g(q)$ if $q$ is a leaf of $\mathcal{T}_1$ different from $p$, 
and by $g'(p_1) = p_1'$ and $g'(p_2) = p_2'$, where $p_1 \neq p$ and $p_2 \neq p$ 
are, respectively, the vertexes of the edges $\Upsilon_l$ and $\Upsilon_r$ starting 
from $p$ (and analogously for $p_1'$ and $p_2'$ in relation to $g(p)$). The map 
$g'$ will be called a {\bf \textit{germination}} of $g$.

In general, given two G-admissible maps $g \!:
lv (\mathcal{R}_1) \rightarrow lv (\mathcal{S}_1)$  and  $h \!:
lv (\mathcal{R}_2) \rightarrow lv (\mathcal{S}_2)$, we will say that $g$ 
is G-equivalent to $h$ if there exists a finite sequence 
$g_0 \! = \! g, g_1, \ldots, g_n \! = \! h$ of G-admissible maps such that, 
for each $k \!\in\! \{1,\ldots,n\}$, either $g_k$ is a germination of $g_{k-1}$,  
or $g_{k-1}$ is a germination of $g_{k}$. Let us denote by $G$ the set of 
G-admissible maps modulo this equivalence relation.

\vspace{0.5cm}


\beginpicture

\setcoordinatesystem units <0.9cm,0.9cm>

\plot
2 6.5
2.5 5.5 /

\plot
3 6.5
2.5 5.5 /

\plot
2 6.5
2.25 7.4 /

\plot
2 6.5
1.75 7.4 /


\plot
5 6.5
5.5 5.5 /

\plot
6 6.5
5.5 5.5 /

\plot
6 6.5
6.25 7.4 /

\plot
6 6.5
5.75 7.4 /


\plot
10 6.5
10.5 5.5 /

\plot
11 6.5
10.5 5.5 /

\plot
10 6.5
10.25 7.4 /

\plot
10 6.5
9.75 7.4 /

\plot
11 6.5
11.25 7.4 /

\plot
11 6.5
10.75 7.4 /


\plot
13 6.5
13.5 5.5 /

\plot
14 6.5
13.5 5.5 /

\plot
14 6.5
14.25 7.4 /

\plot
14 6.5
13.75 7.4 /

\plot
14.25 7.4
14.4 8.3 /

\plot
14.25 7.4
14.1 8.3 /


\plot
3.4 6
4.4 6 /

\plot
4.4 6
4.2 5.9 /

\plot
4.4 6
4.2 6.1 /


\plot
11.4 6
12.4 6 /

\plot
12.4 6
12.2 5.9 /

\plot
12.4 6
12.2 6.1 /

\plot
4.7 1.5
5.2 0.5 /

\plot
5.7 1.5
5.2 0.5 /

\plot
4.7 1.5
4.95 2.4 /

\plot
4.7 1.5
4.45 2.4 /


\plot
7.7 1.5
8.2 0.5 /

\plot
8.7 1.5
8.2 0.5 /

\plot
8.7 1.5
8.95 2.4 /

\plot
8.7 1.5
8.45 2.4 /


\plot
10.7 1.5
11.2 0.5 /

\plot
11.7 1.5
11.2 0.5 /

\plot
11.7 1.5
11.95 2.4 /

\plot
11.7 1.5
11.45 2.4 /

\plot
11.95 2.4
12.1 3.3 /

\plot
11.95 2.4
11.8 3.3 /


\plot
6.1 1
7.1 1 /

\plot
7.1 1
6.9 0.9 /

\plot
7.1 1
6.9 1.1 /


\plot
9.1 1
10.1 1 /

\plot
10.1 1
9.9 0.9 /

\plot
10.1 1
9.9 1.1 /

\setdots

\plot
4.45 2.4
4.6 3.3 /

\plot
4.45 2.4
4.3 3.3 /

\plot
7.7 1.5
7.95 2.4 /

\plot
7.7 1.5
7.45 2.4 /

\put{$\sigma$} at 2.5 5.2
\put{$\sigma$} at 5.5 5.2
\put{$\sigma$} at 10.5 5.2
\put{$\sigma$} at 13.5 5.2
\put{$g$} at 3.9 5.6
\put{$f$} at 11.9 5.6
\put{$p_1$} at 1.7 7.6
\put{$p_2$} at 2.3 7.6
\put{$p_3$} at 3.25 6.75
\put{$g(p_2)$} at 5.35 7.6
\put{$g(p_3)$} at 6.7 7.6
\put{$g(p_1)$} at 4.75 6.75
\put{$q_1$} at 9.6 7.6
\put{$q_2$} at 10.2 7.6
\put{$q_3$} at 10.8 7.6
\put{$q_4$} at 11.4 7.6
\put{$f(q_4)$} at 12.7 6.75
\put{$f(q_1)$} at 13.55  7.7
\put{$f(q_2)$} at 13.68 8.6
\put{$f(q_3)$} at 14.82 8.6

\put{$\sigma$} at 11.2 0.2
\put{$\sigma$} at 8.2 0.2
\put{$\sigma$} at 5.2 0.2
\put{$g$} at 6.6 0.6
\put{$f$} at 9.6 0.6
\put{$r_1$} at 4.2 3.5
\put{$r_2$} at 4.7 3.5
\put{$r_3$} at 4.95 2.6
\put{$r_4$} at 5.7 1.75
\put{$g(r_1)$} at 6.8 2.75
\put{$g(r_2)$} at 7.75 2.75
\put{$g(r_3)$} at 8.65 2.75
\put{$g(r_4)$} at 9.6 2.75
\put{$fg(r_2)$} at 11.4 3.6
\put{$fg(r_3)$} at 12.55 3.6
\put{$fg(r_4)$} at 10.35 1.8
\put{$fg(r_1)$} at 11.2 2.7


\put{} at 0.8 0

\put {Figure 6} at 8.2 -0.5

\endpicture


\vspace{0.5cm}

We will now proceed to define a group structure on $F$ and $G$. Again, 
we will give the definition only for $G$, since that for $F$ is analogous. 
Fix two elements $f,g$ in $G$. It is not difficult to verify that there exist 
dyadic trees $\mathcal{R}$, $\mathcal{S}$, $\mathcal{T}$ such that in the 
class of $f$ and in that of $g$ there exist maps --which we will still denote by $f$ 
and $g$, respectively-- so that $g \! : lv(\mathcal{R}) \rightarrow lv(\mathcal{S})$ 
and $f \! : lv(\mathcal{S}) \rightarrow lv(\mathcal{T})$. We then define the element 
$fg \in G$ as the class of the map 
$$fg \!: lv (\mathcal{R}) \rightarrow lv (\mathcal{T}).$$
The reader can easily check that this definition does not depend on the chosen 
representatives and that, endowed with this product, $G$ becomes a group. Figure 
6 illustrates the composition of two elements in $G$. Notice that the 
neutral element is the class of the map sending the root (viewed as 
the unique leaf of the degenerate tree) into itself.

We now explain the relationship between the groups $G$ and $F$ defined above, and 
those acting on the circle and the interval, respectively. For this, to each 
vertex of a dyadic tree we associate a {\bf {\em dyadic interval}} 
\index{dyadic!interval} 
in $[0,1]$ ({\em i.e.}, an interval of type $\big[ i/2^n,(i+1)/2^n \big]$, 
with $i \in \{1,\ldots,2^n\}$) as follows:

\vspace{0.08cm}

\noindent{(i) to the root we associate the interval $[0,1]$;}

\vspace{0.08cm}

\noindent{(ii) if to the vertex $p$ which is not a leaf we have associated the 
interval $[a,b]$, then to $p_1$ and $p_2$ we associate the intervals 
$\big[ a,(a+b)/2 \big]$ and $\big[ (a+b)/2,b \big]$, respectively, 
where $p_1 \neq p$ and 
$p_2 \neq p$ are the final vertexes of the edges $\Upsilon_l$ and $\Upsilon_r$ 
starting from $p$.}

\vspace{0.04cm}

For each $g \!\in\! G$ we choose a representative  
$g\!: lv (\mathcal{T}_1) \rightarrow lv (\mathcal{T}_2)$, and 
we associate to it the circle homeomorphism sending affinely the 
interval corresponding to each leaf $p$ of $\mathcal{T}_1$ into the 
interval corresponding to the leaf $g(p)$. It is not difficult to check  
that this definition does not depend on the chosen representative of $g$. 
Thus we obtain a group homomorphism from the recently defined 
group $G$ into the group G acting on the circle, and one is easily 
convinced that this homomorphism is, in fact, an isomorphism.

\vspace{0.25cm}

To conclude this section, we show that $\mathrm{G}$ embeds into 
$\mathrm{Diff}_{+}^{1+\mathrm{Lip}}(\clo)$. The proof that 
we give is based on an idea due to Thurston. Following a construction due 
to Ghys and Sergiescu, we will show in the next section that a stronger 
result holds: $\mathrm{G}$ is topologically conjugate to a subgroup of 
the group of $\ce^{\infty}$ circle diffeomorphisms.
\index{Farey sequence} 

Thurston's idea uses finite partitions of $\clo$ given by the Farey sequence 
instead of partitions into dyadic intervals. In other words, to each vertex 
of a dyadic tree we associate a subinterval of $[0,1]$ in the following way:

\vspace{0.08cm}

\noindent{(i) to the root we associate the interval $[0,1]$;}

\vspace{0.08cm}

\noindent{(ii) if to the vertex $p$ we have associated the interval $[a/b,c/d]$ 
and $p$ is not a leaf, then to $p_1$ and $p_2$ we associate the intervals 
$\big[ a/b,(a+b)/(c+d) \big]$ and $\big[(a+b)/(c+d),c/d \big]$, 
respectively. Here, $p_1 \neq p$ and $p_2\neq p$ are the 
vertexes different from $p$  of the edges $\Upsilon_l$ 
and $\Upsilon_r$ starting from $p$.}

\vspace{0.08cm}

Then for each $g \!\in\! G$ we choose a representative $g \!: lv (\mathcal{T}_1)
\rightarrow lv (\mathcal{T}_2)$, to which we associate the circle homeomorphism 
sending the interval associated to each leaf $p$ of $\mathcal{T}_1$ into 
the interval associated to $g(p)$ by a (uniquely determined) map in 
$\mathrm{PSL}(2,\mathbb{Z})$. As in the previous case, everything is 
well-defined up to the equivalence relation defining the group structure on $G$. 
Actually, one can explicit the elements in $\mathrm{PSL}(2,\mathbb{Z})$ used in 
this definition. Indeed, it is not difficult to verify by induction 
that if the interval associated to some vertex is $[a/b,c/d]$, 
then $bc-ad=1$. Therefore, the unique map in $\mathrm{PSL}(2,\mathbb{Z})$ 
sending $I \!=\! [a/b,c/d]$ into $J \!=\! [a'/b',c'/d']$ is \esp 
$\gamma_{I,J} = \gamma_{J} \circ \gamma_I^{-1}$, \esp where 
$$\gamma_I(x) = \frac{(c-a)x+a}{(d-b)x+b}, \qquad
\gamma_J(x) = \frac{(c'-a')x+a'}{(d'-b')x+b'}.$$
Notice that \esp $\gamma_{I}'(x) = 1/((d-b)x+b)^2$, \esp and therefore 
$$\gamma_{I,J}' \left( \frac{a}{b} \right) = \left( \frac{b}{b'} \right)^2, \qquad
\gamma_{I,J}' \left( \frac{a'}{b'} \right) = \left( \frac{d}{d'} \right)^2.$$
These equalities show that, for each $g \!\in\! G$, the associated 
piecewise $\mathrm{PSL}(2,\mathbb{Z})$ circle 
homeomorphism is actually of class 
$\mathrm{C}^{1+\mathrm{Lip}}$; indeed, the values of the derivatives 
to the left and to the right of the ``break points'' coincide.\\

We have thus constructed a new action of $\mathrm{G}$ on the circle, this time 
by $\mathrm{C}^{1+\mathrm{Lip}}$ diffeomorphisms. It is not very difficult to 
show that the resulting group of piecewise $\mathrm{PSL}(2,\mathbb{Z})$ maps 
and that of piecewise dyadically affine maps are topologically conjugate.


\subsection{Ghys-Sergiescu's smooth realization}
\label{ghys-sergi}\index{Ghys!Ghys-Sergiescu's realization of Thompson's groups}

\hspace{0.45cm} A remarkable (and at first glance surprising) property of Thompson's group G 
\index{Thompson's groups!G} 
is the fact that it can be realized as a group of $\ce^{\infty}$ circle diffeomorphisms. 
We mention, however, that the general problem of knowing what are the subgroups of 
$\mathrm{PAff}_+([0,1])$ and $\mathrm{PAff}_+(\clo)$ sharing this property is wide 
open. This problem is particularly interesting (both from the dynamical and algebraic 
viewpoints) for the groups studied by Stein in \cite{stein} and their natural 
analogues acting on the circle. 

Following (part of) \cite{GS}, we will associate 
a representation of G in $\mathrm{Homeo}_+(\clo)$ to each homeomorphism 
$H\!: \mathbb{R} \rightarrow \mathbb{R}$ satisfying the following properties:

\vspace{0.15cm}

\noindent (i) for each $x \in \mathbb{R}$ one has $H(x+1) = H(x) + 2$;

\vspace{0.15cm}

\noindent (ii) $H(0) = 0.$

\vspace{0.15cm}

Notice that the function \esp $H(x) \!=\! 2x$ \esp satisfies these 
two properties: the associated representation will correspond to the 
canonical action of G by piecewise affine circle homeomorphisms.

For the construction let us first introduce some notation. \index{dyadic!rational number} 
First of all, let us denote by $\mathbb{Q}_2 (\mathbb{R})$ the 
group of dyadic rational numbers (which may be thought of as a subgroup 
of the translation group). By $\mathrm{Aff}_+(\mathbb{Q}_2,\mathbb{R})$ we 
will denote the group of affine transformations of the real line which preserve the set 
of the dyadic rationals, and by $\mathrm{PAff}_+(\mathbb{Q}_2,\mathbb{R})$ we will denote 
the group of homeomorphisms which are piecewise $\mathrm{Aff}_+(\mathbb{Q}_2,\mathbb{R})$. 
Similarly, $\mathbb{Q}_2(\clo)$ will denote the group of dyadic rational 
rotations of the circle, 
\index{dyadic!rotation} \index{rotation!dyadic|see{dyadic, rotation}} 
and $\mathrm{PAff}_+ (\mathbb{Q}_2,\clo)$ will denote the 
group of piecewise dyadically affine circle homeomorphisms. Recall finally that, 
for each $a \!\in\! \mathbb{R}$, the translation by $a$ is denoted $T_a$.

\vspace{0.15cm}

\begin{lem} {\em The correspondence 
$\Phi_H \! : \mathbb{Q}_2(\mathbb{R}) \rightarrow \mathrm{Homeo}_+(\mathbb{R})$ 
sending $p/2^q$ into $H^{-q} T_p H^q$ is well-defined and is a group homomorphism.}
\end{lem}

\noindent{\bf Proof.} To show that the definition does not lead to a 
contradiction we need to check that, for all integers $q\geq 0$ and $p$, 
one has \esp $\Phi_H \big( p / 2^q\big) = \Phi_H( 2p / 2^{q+1})$, \esp that is,
$$H^{-q} T_p H^q = H^{-(q+1)} T_{2p} H^{q+1}.$$
To do this,  notice that this equality is equivalent to $H(x+p) \!=\! H(x) + 2p$, 
which follows directly from property (i). To show that $\Phi_H$ is a 
group homomorphism, just remark that 
$$\hspace{0.6cm} \Phi_H \Big( \frac{p}{2^q} + \frac{p'}{2^{q}} \Big) = H^{-q}
T_{p+p'} H^q = H^{-q} T_p H^{q} H^{-q} T_{p'} H^q = \Phi_H \Big(
\frac{p}{2^q} \Big) \Phi_H \Big( \frac{p'}{2^q} \Big). \esp \esp
\esp \mbox{        } \hspace{0.6cm} \square$$

\vspace{0.15cm}

\begin{lem} {\em The homomorphism $\Phi_H$ 
of the preceding lemma extends to a homomorphism from 
$\mathrm{Aff}_+(\mathbb{Q}_2,\mathbb{R})$ into $\mathrm{Homeo}_+(\mathbb{R})$ by} 
$$\left(
\begin{array}
{cc}
2^n & p/2^q  \\
0 & 1  \\
\end{array}
\right) \longmapsto \Phi_H \Big( \frac{p}{2^q} \Big) \circ H^n.$$
\end{lem}

\noindent{\bf Proof.} The claim follows directly from the equality 
$$H \circ \Phi_H \Big( \frac{p}{2^{q+1}} \Big) = \Phi_H \Big( \frac{p}{2^q} \Big)
\circ H.$$

\vspace{-0.8cm}

$\hfill\square$

\vspace{0.65cm}

The extension of the homomorphism $\Phi_H$ 
to $\mathrm{Aff}_+(\mathbb{Q}_2,\mathbb{R})$ 
will be also denoted by $\Phi_H$. The definition for each 
$g \!\in\! \mathrm{PAff}_+(\mathbb{Q}_2,\mathbb{R})$ is a little more subtle. 
Let us fix a strictly increasing sequence $(a_n)_{n \in \mathbb{Z}}$ of 
dyadic rational numbers without accumulation points, as well as a sequence  
of elements $h_n \!\in\! \mathrm{Aff}_+(\mathbb{Q}_2,\mathbb{R})$, in 
such a way that for all $n \in \mathbb{Z}$ one has
$$g \big|_{[a_n,a_{n+1}]} = h_n \big|_{[a_n,a_{n+1}]}.$$
If we define $b_n = \Phi_H (a_n) (0)$, then it is easy to see 
that the sequence $(b_n)_{n\in\mathbb{Z}}$ is also strictly 
increasing and does not have accumulation points. 

\vspace{0.2cm}

\begin{prop} {\em If to each $g \in \mathrm{PAff}_+(\mathbb{Q}_2,\mathbb{R})$ we 
associate the map that on each interval $[b_n,b_{n+1}[$ coincides with $\Phi_H (h_n)$, 
then we obtain a homomorphism from $\mathrm{PAff}_+(\mathbb{Q}_2,\mathbb{R})$ into 
$\mathrm{Homeo}_+(\clo)$ which extends $\Phi_H$.}
\end{prop}

\noindent{\bf Proof.} The fact that the map associated to each $g$ is 
well-defined (\textit{i.e.,} it does not depend on the choice of 
the $a_n$'s) follows readily from the definition, as well as 
the fact that the map associated to each 
$g \!\in\! \mathrm{Aff}_+(\mathbb{Q}_2,\mathbb{R})$ coincides 
with $\Phi_H (g)$. To prove that the map associated to each 
$g \in \mathrm{PAff}_+(\mathbb{Q}_2,\mathbb{R})$ is a homeomorphism, we need to 
check the continuity on each point $b_n$, which reduces to show that
$$\Phi_H (h_n)(b_n) = \Phi_H (h_{n-1})(b_n).$$
Notice that the above equality is equivalent to 
$$\Phi_H (h_n T_{a_n})(0) = \Phi_H (h_{n-1} T_{a_n})(0),$$
that is,
\begin{equation}
\Phi_H (T_{-a_n} h_{n-1}^{-1} h_n T_{a_n})(0) = 0.
\label{averi}
\end{equation}
Now since $g$ is continuous, we have $h_n (a_n) = h_{n-1}(a_n)$. Therefore 
$T_{-a_n} h_{n-1}^{-1} h_n T_{a_n}$ is an element of $\mathrm{Aff}_+(\mathbb{Q}_2,\mathbb{R})$ 
which fixes the origin, that is, a map $f$ of the form $x \mapsto 2^k x$. The desired 
equality (\ref{averi}) then follows from \esp $\Phi_H (f) = H^k$ \esp and from the 
fact that, by property (ii), $H$ fixes the origin. $\hfill\square$

\vspace{0.4cm}

Notice that from property (i) it follows that $\Phi_H (p) \!=\! T_p$ for every 
integer $p$. Therefore, $\Phi_H$ induces an injective homomorphism (which we will 
still denote by $\Phi_H$) from $\mathrm{G}$ into $\mathrm{Homeo}_+(\clo)$.

\vspace{0.2cm}

\begin{prop} {\em Assume that for some positive integer $r$ or for $r = \infty$, 
the map $H$ is a $\ce^r$ diffeomorphism satisfying the following condition:

\vspace{0.1cm}

\noindent{$(\mathrm{iii})_r$ \esp $H'(0)=1$ \esp and \esp 
$H^{(i)}(0) = 0 $ for all $i \in \{2,\ldots,r\}$.}

\vspace{0.1cm}

\noindent Then the image $\Phi_H (\mathrm{G})$ is 
contained in the group of $\ce^r$ circle diffeomorphisms.}

\end{prop}

\noindent{\bf Proof.} Using the notation of the preceding proposition, 
we need to show that for each $i \in \{1,\ldots,r\}$ one has 
$$\Phi_H (h_n)^{(i)}(b_n) = \Phi_H (h_{n-1})^{(i)} (b_n).$$
But this follows from the fact that the Taylor series expansion of 
order $r$ of the map $\Phi_H (T_{-a_n} h_{n-1}^{-1} h_n T_{a_n}) = H^k$
\esp coincides with that of the identity. $\hfill\square$

\vspace{0.5cm}

Notice that, when $r \!=\! \infty$, property (iii)$_r$ cannot be satisfied by any 
real-analytic diffeomorphism. In fact, the groups F and G cannot act faithfuly 
\index{action!faithful} on $\clo$ 
by real-analytic diffeomorphisms. (Since G is a simple group, this 
implies that every action of G by real-analytic diffeomorphisms of the 
circle is trivial.) \index{group!simple} This may be shown in many 
distinct ways, but it will appear as evident after \S \ref{solvito}: 
Thompson's group F contains \index{group!solvable} 
solvable subgroups of arbitrary length of solvability, while every solvable 
group of real-analytic diffeomorphisms of the (closed) interval is metabelian.
 
In \S \ref{sec-minimal} we will deal again with some dynamical aspects 
of the preceding realization. To conclude, let us point out that 
the dyadic feature of the preceding arguments is not essential: for 
each integer $m \!\geq\! 2$, an analogous construction starting with a 
map $H$ satisfying $H(x+1) = H(x) + m$ for each $x \!\in\! \mathbb{R}$ 
leads to an {\bf \textit{$m$-adic Thompson's group}}. From an algebraic point 
of view, the case $m\!=\!2$ is special in relation to the automorphism 
group: we refer the reader to \cite{brin-IHES,brin-gen} for further information 
on this very interesting topic. For a complete survey on the recent progress 
(especially on cohomological aspects) on Thompson's groups, see \cite{vlad}.
\index{Thompson's groups|)}
\index{Thompson's groups!$m$-adic}


\chapter{DYNAMICS OF GROUPS OF HOMEOMORPHISMS}


\section{Minimal Invariant Sets}

\subsection{The case of the circle}
\label{sec-minimal}

\hspace{0.45cm} Recall that a subset $\Lambda$ of $\mathrm{S}^1$ is {\bf \textit{invariant}} 
by a subgroup $\Gamma$ of $\mathrm{Homeo}_{+}(\mathrm{S}^1)$ if $g(x) \!\in\!
\Lambda$ for every $x \!\in\! \Lambda$ and every $g \!\in\! \Gamma$. A 
compact invariant set $\Lambda$ is {\bf \textit{minimal}} if its only 
closed invariant subsets are the empty set and $\Lambda$ itself.

\begin{thm} {\em If $\Gamma$ is a subgroup of $\mathrm{Homeo}_{+}(\mathrm{S}^1)$, 
then one (and only one) of the following possibilities occurs: \\

\vspace{0.05cm}

\noindent{{\em (i)} there exists a finite orbit;}

\vspace{0.05cm}

\noindent{{\em (ii)} all the orbits are dense;}

\vspace{0.05cm}

\noindent{{\em (iii)} there exists a unique minimal invariant 
compact set which is homeomorphic to the Cantor set (and which 
is contained in the set of accumulation points of every orbit).}}
\label{invariante}
\end{thm}

\noindent {\bf Proof.} The family of non-empty closed invariant subsets of $\mathrm{S}^1$ 
is ordered by inclusion. Since the intersection of nested compact sets is (compact and) 
non-empty, the Zorn Lemma allows 
\index{Zorn Lemma}
us to conclude the existence of a minimal non-empty closed 
invariant set $\Lambda$. The boundary $\partial \Lambda$ and the set 
$\Lambda'$ of the accumulation points of $\Lambda$ are closed invariant sets contained 
in $\Lambda$. By the minimality of $\Lambda$, one of the following possibilities occurs:\\

\vspace{0.08cm}

\noindent{(i) $\Lambda'$ is empty: in this case $\Lambda$ is a finite orbit;}

\vspace{0.08cm}

\noindent{(ii) $\partial \Lambda$ is empty: in this case $\Lambda = \clo$, 
and therefore all the orbits are dense;}\\

\vspace{0.08cm}

\noindent{(iii) $\Lambda = \Lambda' = \partial \Lambda$: in this 
case $\Lambda$ is a closed set with empty interior and having no isolated 
point; in other words, $\Lambda$ is homeomorphic to the Cantor set.}

\vspace{0.08cm}

We will show that, in the last case, $\Lambda$ is contained in the set of accumulation 
points of every orbit, which clearly implies its uniqueness. Let $x$ and $y$ be arbitrary 
points in $\clo$ and $\Lambda$, respectively. We need to show that there exists a 
sequence $(g_n)$ of elements in $\Gamma$ such that $g_n(x)$ converges to $y$. For $x \in
\Lambda$, this follows from the minimality of $\Lambda$. If $x \in \clo \setminus \Lambda$, 
let us consider the interval $I\!= ]a,b[$ contained in $\clo \setminus \Lambda$ such that 
both $a,b$ belong to $\Lambda$ and $x \!\in\! I$. Since the orbit of $a$ is dense in 
$\Lambda$, and since $\Lambda$ 
does not have isolated points, there must exist a sequence $(g_n)$ in $\Gamma$ 
for which $g_n(a)$ tends to $y$ in such a way that the intervals $g_n(I)$ are 
two-by-two disjoint. The length of $g_n(I)$ must go to zero, and thus 
$g_n(x)$ converges to $y$. This concludes the proof. $\hfill\square$

\vspace{0.01cm}

\begin{small}\begin{ejer} Prove that if a group of circle homeomorphisms 
has finite orbits, then all of them have the same cardinality. 
\end{ejer}
\end{small}

In the case where all the orbits are dense, the action is said to be {\bf \textit{minimal}}. If 
\index{action!minimal} 
there exists a minimal invariant Cantor set, this set is called an {\bf \textit{exceptional minimal 
set}}.\index{exceptional minimal set} 
To understand this case better, it is useful to introduce the following terminology.

\begin{defn} A circle homeomorphism $f$ is {\bf {\em semiconjugate}} to $g$ if 
\index{semiconjugacy} there exists 
a continuous degree-1 map $\varphi\!: \mathrm{S}^1 \!\rightarrow\! \mathrm{S}^1$ 
whose lifts to $\mathbb{R}$ are non-decreasing functions and such that 
$\varphi f = g \esp \! \varphi$. Similarly, a group action $\Phi_1$ by circle 
homeomorphisms is semiconjugate to an action $\Phi_2$ if there exists $\varphi$ 
satisfying the preceding properties and such that \esp 
$\varphi \esp \Phi_1(g) = \Phi_2(g) \esp \varphi$ \esp  
for every element $g$ in the acting  group.
\end{defn}

The map $\varphi$ may be non injective; in the case it is injective, the 
semiconjugacy is in fact a conjugacy. Notice that if a group $\Gamma$ acts on 
$\clo$ with an exceptional minimal set $\Lambda$, then replacing the closure 
of each connected component of $\clo \setminus \Lambda$ by a point we obtain a 
topological circle $\clo_{\Lambda}$ upon which $\Gamma$ acts in a natural way 
by homeomorphisms. The original action is semiconjugate to the induced 
minimal action on $\clo_{\Lambda}$.

\begin{small}
\begin{obs} The relation of semiconjugacy is not an equivalence relation. The equivalence 
relation which generates is sometimes called {\bf {\em monotonic equivalence}}. More precisely, 
two actions $\Phi_1$ and $\Phi_2$ are monotonically equivalent if there exists an action $\Phi$ 
which is semiconjugate to both $\Phi_1$ and $\Phi_2$. The interested reader may find more on 
this notion in \cite{cal-inv}.
\end{obs} \end{small}

\vspace{0.65cm}

\begin{center}

\beginpicture
\setcoordinatesystem units <0.8cm,0.8cm>

\circulararc 360 degrees from 4 0
center at 0 0

\circulararc 106 degrees from  2.4  -3.2
center at 0 -5

\circulararc 106 degrees from -3.971264  -0.478448
center at -4.3301270     2.5

\circulararc -106 degrees from  3.971264  -0.478448
center at 4.3301270     2.5

\circulararc -145 degrees from  0.35  -3.96
center at 1.12  -4

\circulararc 145 degrees from  -0.35  -3.96
center at -1.12  -4

\circulararc -145 degrees from  -3.604  1.6769
center at -4.04 1

\circulararc 145 degrees from  -3.2543  2.2831
center at -2.9  2.98

\circulararc 145 degrees from  3.604  1.6769
center at 4.04  1

\circulararc -145 degrees from  3.2543  2.2831
center at 2.9  2.98

\circulararc 110 degrees from 0.5 -0.22
center at 0 0

\circulararc 160 degrees from 0.5 -1.92
center at 0 -2

\plot
-0.53 -1.92
-0.53 -1.72 /

\plot
-0.53 -1.92
-0.36 -1.77 /

\put{$g$} at 0 -1.2

\plot 0.05 0.55 0.2 0.6 /

\plot 0.05 0.55 0.2 0.4 /

\plot 0 0
0 4 /

\plot
0 0
3.4641016  -2 /

\plot
0 0
-3.4641016   -2 /

\put{$\bullet$} at 0 0

\put{$\bullet$} at 0 -2

\put{$O$} at -0.34 0.2

\put{$P$} at -0.1 -2.28

\put{$A$} at 2 0.2

\put{$B$} at -2 0.2

\put{$g(A)$} at -1 -2.9

\put{$g(B)$} at 1 -2.9

\put{$fg(A)$} at 2.6 1.15

\put{$fg(B)$} at 2.08 2

\put{$f^2g(A)$} at -2.1 2

\put{$f^2g(B)$} at  -2.5 1.35

\put{$f$} at 0.8 0.3

\put{} at -8.4 3.4

\put{Figure 7} at 0 -4.8

\endpicture


\end{center}

\vspace{-0.15cm}

It is not difficult to construct homeomorphisms of $\clo$ which do admit exceptional 
minimal sets. Indeed, every subset of $\mathrm{S}^1$ homeomorphic to the Cantor 
set may be exhibited as the exceptional minimal set of a circle homeomorphism 
(see \S \ref{teoria-poinc}). We now show an example of a group 
of real-analytic circle diffeomorphisms having an exceptional minimal set. This 
group is generated by two M\"obius transformations 
\index{M\"obius group} 
$f$ and $g$. The diffeomorphism $f$ is the rotation $R_{2\pi/3}$ 
(centered at the origin $O = (0,0)$). The diffeomorphism $g$  
is also elliptic and corresponds to a rotation of angle $\pi$ 
centered at a point $P \in \mathbb{D}$ situated 
at an Euclidean distance larger than $2-\sqrt{3}$ from $O$ (see Figure 7). Equivalently, 
$g$ is the ``hyperbolic reflexion'' with respect to the geodesic passing through $P$ and 
perpendicular to the geodesic joining this point to $O$. 

The acting group in this example coincides with the so-called 
{\bf {\em modular group}}, which admits the finite presentation \esp 
$\Gamma \!=\! \langle f,g \!: \hspace{0.02cm} g^2 \!=\! f^3 \!=\! id \rangle.$ 
\esp The standard injection of $\Gamma$ into $\mathrm{PSL}(2,\mathbb{R})$ 
is obtained by identifying $f$ with $R_{2\pi/3}$ and $g$ with the 
rotation by $\pi$ centered at the point $(0,\sqrt{3}-2) \!\in\! \mathbb{D}$. 
Via this injection, the modular group identifies with $\mathrm{PSL}(2,\mathbb{Z})$, 
and the corresponding action is minimal.
\index{group!modular}

A small perturbation of a ``piecewise version" of the preceding example allows us 
to produce an action by $\ce^{\infty}$ diffeomorphisms having the triadic Cantor 
set as an exceptional minimal set. For this, it suffices to consider a circle 
diffeomorphism $h$ satisfying:

\vspace{0.09cm}

\noindent -- its restriction to $[0,\pi/6]$ (resp. $[\pi,3\pi/2]$)  
is affine, with derivative $3$ (resp. $1/3$); 

\vspace{0.09cm}

\noindent -- $h(x) = x + 2\pi/3$ \esp for all \esp $x \in [\pi/3,\pi/2]$; 

\vspace{0.09cm}

\noindent -- $h$ is defined in a coherent way on the remaining intervals so that $(hg)^3 = Id$, 
where $g$ denotes the rotation by an angle $\pi$ (see Figure 8). 

\vspace{0.09cm}

\noindent Letting $f = hg$, we obtain the desired action.

\vspace{0.1cm}

It is very interesting to notice that, although the example illustrated by Figure 7  
already appears in the works by Klein and Poincar\'e, for many years this example 
was ``forgotten'', and actually the ``first example'' of a group of $\ce^{\infty}$ 
circle diffeomorphisms having an exceptional minimal set is commonly attributed to 
Sacksteder \cite{sac-four}. \index{Sacksteder!example}
His example corresponds to (a slight modification of) 
the one illustrated below...
\index{Klein}\index{Poincar\'e}

\vspace{0.5cm}


\beginpicture

\setcoordinatesystem units <0.52cm,0.52cm>

\putrule from 2 0 to 14 0

\putrule from 2 0  to 2 12

\putrule from 14 0 to 14 12

\putrule from 2 12 to 14 12

\put{$0$} at 2 -0.4

\put{$2\pi$} at 14 -0.4

\put{$h$} at 7.07 7.35

\put{$g$} at 5.4  9.8

\put{$g$} at 10.7  2.3

\put{} at -4.75 -0.2

\put{Figure 8} at 8 -0.65

\plot
2  6
8  12 /

\plot
8   0
14  6 /

\plot
2    0
3    3 /

\plot
8  8
11  9 /

\plot
4  6
5  7 /

\plot
3    3
3.1  3.318
3.2  3.664
3.3  4.026
3.4  4.392
3.5  4.75
3.6  5.088
3.7  5.394
3.8  5.656
3.9  5.862
4    6 /

\plot
8     8
7.682 7.9
7.336 7.8
6.974 7.7
6.608 7.6
6.25  7.5
5.912 7.4
5.606 7.3
5.344 7.2
5.138 7.1
5     7  /

\plot
11    9
11.1  9.03126
11.2  9.05896
11.5  9.12963
11.8  9.200296
12   9.25926
12.3 9.38326
12.7 9.652296
13   9.96296
13.3 10.39363
13.5 10.75926
13.7 11.196
13.9 11.710963
14   12 /

\setdots

\putrule from 5 0 to 5 12

\putrule from 8 0 to 8 12

\putrule from 11 0 to 11 12

\putrule from 2 3 to 14 3

\putrule from 2 6 to 14 6

\putrule from 2 9 to 14 9

\plot
2 0
14 12 /

\plot
3 0
3 9 /

\plot
4 0
4 9 /

\plot
9 0
9 9 /

\plot
10 0
10 9 /

\plot
2 1
11 1 /

\plot
2 2
11 2 /

\plot
2 7
11 7 /

\plot
2 8
11 8 /

\endpicture


\vspace{0.46cm}

In the case of {\bf {\em Fuchsian groups}} 
\index{group!Fuchsian}
(that is, discrete subgroups 
of $\mathrm{PSL}(2,\mathbb{R})$), a particular terminology for Theorem 
\ref{invariante} is sometimes used. If there exists a finite orbit, then 
the group is called {\bf \textit{elementary}}. In the case the orbits are dense, 
the group is said to be of {\bf \textit{first kind}}. Finally, in the case of an exceptional 
minimal set, the group is of {\bf \textit{second kind}} \cite{Ka}. An interesting family of 
examples of groups of second kind are the so-called {\bf \textit{Schottky groups}}. 
\index{group!Schottky}
These correspond to groups generated by hyperbolic elements $g_0$ and $g_1$ 
in $\mathrm{PSL}(2,\mathbb{R})$ for which there exist disjoint intervals 
$I_0,I_1,J_0,J_1$ in $\clo$ such that \esp $g_i(I_i \cup
J_i \cup J_{1-i}) \subset I_i$ \esp and \esp $g_i^{-1}(I_i \cup
I_{1-i} \cup J_{1-i}) \subset J_i$ \esp for $i \in \{0,1\}$. Figure 
9 illustrates the corresponding combinatorics of the dynamics. 
It is easy to check that $\langle g_0,g_1 \rangle$ acts on $\clo$
admitting an exceptional minimal set, namely 
$$\Lambda =  \bigcap_{g \in \langle g_0,g_1 \rangle} g(\overline{I_0}
\cup \overline{I_1} \cup \overline{J_0} \cup \overline{J_1}).$$
Moreover, the group generated by $g_0$ and $g_1$ is free 
(see \S \ref{sec-margulis}).

\vspace{0.55cm}

\beginpicture

\setcoordinatesystem units <0.7cm,0.7cm>

\circulararc 360 degrees from 4 0 center at 0 0

\circulararc 106 degrees from  2.4  -3.2 center at 0 -5

\circulararc -106 degrees from 2.4  3.2 center at 0  5

\circulararc 106 degrees from  3.2  2.4 center at 5  0

\circulararc -106 degrees from  -3.2  2.4 center at -5  0

\circulararc 146 degrees from  0.65  -3,95 center at 0 -4.125

\circulararc -146 degrees from  -3,95  0.65 center at -4.125 0

\circulararc -146 degrees from  0.65  3,95 center at 0 4.125

\circulararc 146 degrees from  3,95  0.65 center at 4.125 0
\circulararc 146 degrees from  2  -3.45 center at 1.53 -3.85

\circulararc -146 degrees from  -2  -3.45 center at -1.53 -3.85

\circulararc 146 degrees from  3.45  2 center at 3.85 1.53 

\circulararc -146 degrees from  3.45  -2 center at 3.85 -1.53 

\circulararc 146 degrees from  -2  3.45 center at -1.53 3.85

\circulararc -146 degrees from  2  3.45 center at 1.53 3.85

\circulararc 146 degrees from -3.45 -2 center at -3.85 -1.53

\circulararc -146 degrees from  -3.45 2 center at -3.85 1.53
\put{$g_0(I_0)$} at 4.8 0 \put{$g_1(I_1)$} at 0 4.4
\put{$g_1^{-1}(J_1)$} at 0 -4.4 \put{$g_0^{-1}(J_0)$} at -4.9 0

\put{$g_0(I_1)$} at 4.5 1.5 \put{$g_0(J_1)$} at 4.5 -1.5
\put{$g_0^{-1}(I_1)$} at -4.6 1.5 \put{$g_0^{-1}(J_1)$} at -4.6 -1.5

\put{$g_1(J_0)$} at -2 4.05 \put{$g_1(I_0)$} at 2 4.08
\put{$g_1^{-1}(J_0)$} at -2 -4.1 \put{$g_1^{-1}(I_0)$} at 2 -4

\plot -1 0 1 0 /

\plot 0 -1 0 1 /

\plot 1 0 0.8 -0.1 /

\plot 1 0 0.8 0.1 /

\put{$g_0$} at 0.6 -0.3

\plot 0 1 -0.1 0.8 /

\plot 0 1 0.1 0.8 /

\put{$g_1$} at 0.35 0.6

\put{$I_0$} at 1.6 0
\put{$J_1$} at 0 -1.6
\put{$I_1$} at 0 1.6
\put{$J_0$} at -1.6 0

\put{Figure 9} at 0 -5.3 

\put{} at -9.52 4.3

\endpicture


\vspace{0.3cm}

\begin{small} \begin{ejer}
Give examples of finitely generated subgroups of $\mathrm{Diff}_+^{\infty}(\clo)$ 
admitting an exceptional minimal set $\Lambda$ so that the orbits of all points 
in $\clo \setminus \Lambda$ are dense.

\noindent{\underbar{Remark.}} Using a result of Hector 
\index{Hector} 
it is possible to show that, if $\Gamma$ is a (non-necessarily finitely generated) 
group of real-analytic circle diffeomorphisms having an exceptional minimal 
set, then none of its orbits is dense (see for instance \cite{Na-kova}).
\end{ejer}

\begin{ejer} Give an example of a (non finitely generated) group of real-analytic 
circle diffeomorphisms whose action is minimal and such that all of its finitely generated 
subgroups do admit an exceptional minimal set. Analogously, give an example of a group of 
real-analytic circle diffeomorphisms acting minimally and such that all of its finitely 
generated subgroups do admit finite orbits (see Example \ref{ejemplohirsh} in case of 
problems with this).
\end{ejer} \end{small}

We close this section with another important example of a group of $\ce^{\infty}$ circle 
diffeomorphisms having an exceptional minimal set, 
\index{Thompson's groups!G} namely Thompson's group G. Indeed, in \S \ref{ghys-sergi}, to 
each homeomorphism $H\!: \mathbb{R} \rightarrow \mathbb{R}$ satisfying certain properties 
(i) and (ii), we associated 
\index{Ghys!Ghys-Sergiescu's realization of Thompson's groups}
a homomorphism $\Phi_H \!
:  \mathrm{G} \rightarrow \mathrm{Homeo}_+(\clo)$, which takes values in 
$\mathrm{Diff}_+^{r}(\clo)$ when $H$ is a $\ce^r$ diffeomorphism satisfying 
condition (iii)$_r$. It is fairly clear that the groups $\Phi_H(\mathrm{G})$ 
are topologically conjugate to (the canonical inclusion of) 
G (in $\mathrm{Homeo}_+(\clo)$). However, according to the 
following proposition, some of them are not conjugate to G.

\vspace{0.1cm}

\begin{prop} {\em If the homeomorphism 
$H \!: \mathbb{R} \rightarrow \mathbb{R}$ satisfies the 
properties} (i) {\em and } (ii) {\em from} \S \ref{ghys-sergi} 
{\em and has at least two fixed points, then the group 
$\Phi_H (\mathrm{G})$ admits an exceptional minimal set.}
\end{prop}

\noindent{\bf Proof.} If $a,b$ are fixed points of $H$, then property 
(i) implies that they belong to the same ``fundamental domain''. In other words, 
the open interval in the real line whose endpoints are $a$ and $b$ projects injectively 
into an open interval $I$ of the circle satisfying $\bar{H}^n(I)=I$ for every $n \geq 0$ 
(where $\bar{H}$ denotes the degree-2 map of $\clo$ induced by $H$). On the other hand, 
for each $n \geq 0$ the set $\bar{H}^{-n}(I)$ is the union of a family consisting of 
$2^n$ disjoint open intervals, and $\bar{H}^{-n}(I) \subset \bar{H}^{-m}(I)$ for all 
$m \geq n$. Therefore, the union $\cup_{n \geq 0} \bar{H}^{-n}(I)$ is open and invariant by 
$\bar{H}$, and its complementary set is non-empty. It is then not very difficult to conclude 
that this set is invariant by $\Phi_H (\mathrm{G})$ (see Exercise \ref{orbitas-thompson}). 
Therefore, not every orbit of $\Phi_H(\mathrm{G})$ is dense, and since $\Phi_H(\mathrm{G})$ 
has no finite orbit, it must preserve an exceptional minimal set. $\hfill\square$

\vspace{0.32cm}

Using the techniques of Chapter \ref{en-dos}, it is possible to show that every 
faithful action of G by $\ce^2$ circle diffeomorphisms is semiconjugate (by a 
non-necessarily orientation preserving map) to the standard piecewise affine 
action \cite{GS}. However, we will not pursue on this point, since it seems 
that the regularity hypothesis is superfluous for this  claim 
\cite{isa-liousse,Gh1} (compare Exercise \ref{no-F}).

\vspace{0.35cm}


\beginpicture

\setcoordinatesystem units <0.7cm,0.7cm>


\putrule from 0 0 to 6 0
\putrule from 0 0 to 0 6
\putrule from 6 0 to 6 6
\putrule from 0 6 to 6 6

\plot
3 0
6 6 /

\plot
1.5 1.53
3 6 /

\plot
0 0
0.1 0.119
0.2 0.267
0.4 0.593
0.6 0.891
0.75 1.066
0.9 1.191
1.1 1.293
1.3 1.367
1.4 1.42
1.5 1.523 /

\put{$\bar{H}$} at 1.5 3 
\put{$\bar{H}$} at 4.05 3 
\put{} at -6.5 0
\put{Figure 10} at 3 -1

\put{$\bullet$} at 0 0
\put{$\bullet$} at 1.5 0
\put{$\bullet$} at 6 0
\put{$b$} at  1.5 -0.32
\put{$a$} at  0 -0.32
\put{$I$} at  0.75 -0.8

\putrule from 0 -0.8 to 0.6 -0.8
\putrule from 0.9 -0.8 to 1.5 -0.8

\setdots

\plot
0 0
6 6 /

\putrule from 3 0 to 3 6

\endpicture


\vspace{0.25cm}

\begin{small}\begin{ejer}  Using Theorems \ref{arefad} and \ref{margulis}, show 
that every action of F by circle homeomorphisms admits a global fixed point. 
\end{ejer}

\begin{ejer} Prove that if the map $H$ is {\bf {\em expanding}}, that is, if for 
every pair of distinct points $x,y$ in the real line one has \esp 
$dist(H(x),H(y)) > dist(x,y)$, \esp 
then all the orbits of the group $\Phi_H (\mathrm{G})$ are dense.
\index{expanding!map}
\end{ejer}

\begin{ejer} Prove that G admits actions by piecewise affine circle homeomorphisms 
\index{piecewise affine homeomorphism} having an exceptional minimal set.
\end{ejer} 

\begin{ejer} Given a homeomorphism $H$ satisfying the properties (i) and (ii) from \S 
\ref{ghys-sergi}, consider the equivalence relation of the induced map $\bar{H}$ on the 
circle. Prove that the equivalence classes of this relation coincide with the orbits 
of the group $\Phi_H (\mathrm{G})$.
\label{orbitas-thompson} 
\end{ejer} \end{small}


\subsection{The case of the real line}

\hspace{0.45cm} There is no analogue of Theorem \ref{invariante} for general 
groups of homeomorphisms of the real line. In fact, it is not difficult to 
construct examples of subgroups of $\mathrm{Homeo}_+ (\mathbb{R})$ whose action 
is not minimaland so that the only closed invariant subsets are the empty set 
and the real line itself. Nevertheless, a weak version (which in many cases is 
enough for applications) persists for {\em finitely generated} subgroups.

\vspace{0.1cm}

\begin{prop} {\em Every finitely generated subgroup of $\mathrm{Homeo}_+(\mathbb{R})$ 
admits a non-empty minimal invariant closed set.}
\label{para-la-recta}
\end{prop}

\noindent{\bf Proof.} Let $\mathcal{G} \!=\! \{ f_1,\ldots,f_k \}$ be a finite and symmetric 
system of generators for a subgroup $\Gamma$ of $\mathrm{Homeo}_+(\mathbb{R})$ (where 
{\bf{\em symmetric}} means that $f^{-1}$ belongs to $\mathcal{G}$ for every 
$f \in \mathcal{G}$). If $\Gamma$ admits a global fixed point, then the 
claim is obvious. 
If not, fix any point $x_0 \!\in \mathbb{R}$, and let $x_1$ be the maximum 
among the points $f_{i} (x_0)$. We claim that every orbit of $\Gamma$ must intersect 
the interval $[x_0,x_1]$. Indeed, for every $x \!\in\! \mathbb{R}$ the supremum and the 
infimum of its orbit are global fixed points, and therefore coincide with $+ \infty$ 
and $- \infty$, respectively. Thus, we can take $x_0',x_1'$ in this orbit so that 
$x_0' < x_0 < x_1 < x_1'$. Let $f = f_{i_n} \ldots f_{i_1} \in \Gamma$ be 
such that $f(x_0') = x_1'$, where each $f_{i_j}$ belongs to $\mathcal{G}$. 
Let $m \!\in\! \{1,\ldots,n-1\}$ be the largest index 
for which $f_{i_m} \ldots f_{i_1} (x_0') < x_0$. Then  
$f_{i_{m+1}} f_{i_m} \ldots f_{i_1} (x_0')$ is in the orbit of $x$ and 
is greater than or equal to $x_0$; by the definition, 
it is smaller than or equal to $x_1$.

Now let $I= [x_0,x_1]$, and on the family $\mathcal{F}$ 
of non-empty closed invariant subsets of $\mathbb{R}$ let us 
consider the order relation $\precede$ given by \esp 
$\Lambda_1 \succeq \Lambda_2$ if \esp 
$\Lambda_1 \cap I \subset \Lambda_2 \cap I$. \esp 
By the discussion above, every orbit by $\Gamma$ must 
intersect the interval $I$, and so 
$\Lambda \cap I$ is a non-empty compact set for all 
$\Lambda \in \mathcal{F}$. Therefore, Zorn's Lemma 
\index{Zorn Lemma} 
provides us with a maximal element for the order $\precede$, which 
corresponds to the intersection with $I$ of a non-empty minimal 
$\Gamma$-invariant closed subset of $\mathbb{R}$. $\hfill\square$

\vspace{0.3cm}

As in the proof of Theorem \ref{invariante}, one 
can analyze the distinct possibilities for the non-empty 
minimal invariant closed set $\Lambda$ obtained above. For this, notice that the  
boundary $\partial \Lambda$ and the set of accumulation points $\Lambda'$ are also 
closed sets invariant by $\Gamma$. Because of the minimality of $\Lambda$, there 
are three possibilities:

\vspace{0.25cm}

\noindent (i) $\Lambda' = \emptyset$: in this case $\Lambda$ is discrete. If it is 
finite, then it is made up of global fixed points. If it is infinite, then it 
corresponds to a sequence $(y_n)_{n \in \mathbb{Z}}$ satisfying $y_{n} \!<\! y_{n+1}$ 
for all $n$ and without accumulation points in $\mathbb{R}$.

\vspace{0.25cm}

\noindent (ii) $\partial \Lambda = \emptyset$: in this case $\Lambda$ 
coincides with the whole line, and hence the action is minimal.  

\vspace{0.25cm}

\noindent (iii) $\partial \Lambda = \Lambda' = \Lambda$: in this case $\Lambda$ 
is ``locally'' a Cantor set. Therefore, collapsing to a point the closure of each 
connected component of the complement of $\Lambda$, we obtain a topological line 
on which the original action induces (by semi-conjugacy) a minimal action of $\Gamma$. 
However, the reader may easily construct examples showing that, unlike the case of the 
circle, in this case the ``exceptional'' minimal invariant set is not necessarily unique.


\section{Some Combinatorial Results}

\subsection{Poincar\'e's theory}
\label{teoria-poinc}

\index{Poincar\'e!rotation number|see{rotation number}}

\hspace{0.45cm} In this section, we revisit the most important dynamical 
invariant of circle homeomorphisms, namely the {\em rotation number}. 
\index{rotation number|(}
We begin with the following well-known lemma about almost subadditive sequences.

\begin{lem} {\em Let $(a_n)_{n \in \mathbb{Z}}$ be a sequence of real 
numbers. Assume that there exists a constant $C \in \mathbb{R}$ such 
that, for all $m,n$ in $\mathbb{Z}$,  
\begin{equation}
|a_{m+n} - a_m - a_n| \leq C. 
\label{subaditividad}
\end{equation} 
Then there exists a unique $\rho \in \mathbb{R}$ 
such that the sequence $(|a_n - n\rho|)_{n \in \mathbb{Z}}$ is bounded. This 
number $\rho$ is equal to the limit of the sequence $(a_n / n)$ as $n$ goes to 
\esp $\pm \infty$ \esp (in particular, this limit exists).}
\label{subaditivo}
\end{lem}

\noindent{\bf Proof.} For each $n \!\in\! \mathbb{N}$ let us consider the 
interval $I_n = \big[ (a_n - C)/n,(a_n + C)/n \big]$. We claim that 
$I_{mn}$ is contained in $I_n$ for every $m,n$ in $\mathbb{N}$. 
Indeed, by (\ref{subaditividad}) we have \esp 
$|a_{mn} - ma_n| \leq (m-1)C$, \esp from where one concludes that
$a_{mn} + C \leq ma_n + mC$, \esp and therefore
$$\frac{a_{mn} + C}{mn} \leq \frac{a_n + C}{n}.$$
Analogously, 
$$\frac{a_{mn} - C}{mn} \geq \frac{a_n - C}{n},$$
and these two inequalities together imply that $I_{mn} \subset I_n$.\\

Due to the claim above, a direct application of the Finite Intersection Property shows 
that the set $I = \cap_{n \in \mathbb{N}} I_n$ is non-empty. If $\rho$ belongs to $I$, 
then $\rho$ is contained in each of the intervals $I_n$. This allows 
to conclude that, for every $n \in \mathbb{N}$,
\begin{equation}
|a_n -  n \rho| \leq C,
\label{mejorar}
\end{equation}
thus showing that $\rho$ satisfies the claim of the lemma. If $\rho' \neq \rho$ then 
$$|a_n - n \rho'| = |(a_n - n\rho) + n(\rho - \rho')| \geq n|\rho - \rho'| - C,$$
and therefore $|a_n - n\rho'|$ goes to infinity. Finally, dividing by $n$ the 
expressions in both sides of (\ref{mejorar}), and then passing to the limit, 
one concludes that $\rho = \lim_{n \rightarrow \infty} (a_n/n)$. The case 
where $n<0$ is similar and we leave it to the reader. $\hfill\square$

\vspace{0.35cm}

We now consider the parameterization of $\clo$ by $[0,1]$. Given a circle homeomorphism $f$, 
we denote by $F \!: \mathbb{R} \rightarrow \mathbb{R}$ any lift of $f$ to the real line.

\vspace{0.1cm}

\begin{prop} {\em For each $x \in \mathbb{R}$ there exists the limit 
$$\lim\limits_{n \rightarrow \pm \infty} \frac{1}{n}[ F^n(x)-x ],$$
and this limit does not depend on $x$.}
\label{limitar}
\end{prop}

\noindent{\bf Proof.} First notice that, for every 
$x,y$ in $\mathbb{R}$ and every $m \in \mathbb{Z}$, 
\begin{equation}
\big| |F^m(x)-x|-|F^m(y)-y| \big| \leq 2. \label{desig}
\end{equation}
Indeed, since $F^m(x+n)=F^m (x)+n$ for every $n \!\in\! \mathbb{Z}$, 
to show (\ref{desig}) we may assume that $x$ and $y$ belong to $[0,1[$, 
and in this case inequality (\ref{desig}) is clear. Let us now 
fix $x \in \mathbb{R}$, and let us denote $a_n=F^{n}(x)-x$. From 
$$a_{m+n} = F^{m+n}(x)-x = \left[
F^{m}(F^{n}(x))-F^{n}(x) \right] + \left[ F^{n}(x)-x \right]$$ 
one concludes that
\begin{eqnarray*}
| a_{m+n} - a_m - a_n | &=& \big| F^{m+n}(x)-x - (F^m(x)-x) -
(F^n(x)-x) \big|\\
&=& \big| F^m(F^n(x)) - F^n(x) - (F^m(x) - x) \big| \\
&\leq& 2.
\end{eqnarray*}
By the preceding lemma, the expression $[F^n(x)-x ] / n$ converges as $n$ 
goes to infinity, and inequality (\ref{desig}) shows that the 
corresponding limit does not depend on $x$. $\hfill\square$

\vspace{0.35cm}

If we consider two different lifts of $f$ to the real line, then  
the limits given by the proposition above coincide up to an integer 
number. We then define the {\bf \textit{rotation number}} of $f$ by
$$\rho(f) = \lim\limits_{n \rightarrow \pm \infty} \frac{1}{n} [F^n(x)-x]
\hspace{0.5cm} \mathrm{mod} \hspace{0.2cm} 1.$$
As a matter of example, it is easy to check that, for the Euclidean rotation 
of angle $\theta \in [0,1[$, one has $\rho(R_{\theta}) = \theta$. 
\index{rotation}

Notice that for every circle homeomorphism $f$, 
every $x \in \mathrm{S}^1$, and every $m \in \mathbb{Z}$, one has
$$\rho (f^m) = \lim\limits_{n \rightarrow \pm \infty} \frac{1}{n}[ F^{mn}(x)-x ] =
 m \cdot \lim\limits_{n \rightarrow \pm \infty} \frac{1}{mn}[
F^{mn}(x)-x ],$$
from where one concludes that $\rho(f^m) = m \rho(f).$

If $f$ has a {\bf \em{periodic point}}, then $\rho(f)$ is a rational number. Indeed, 
if $f^p(x) = x$ then for each lift $F$ of $f$ we have $F^p(x) = x + q$ for some 
$q \in \mathbb{Z}$ (here, $x$ denotes either a point in $\mathrm{S}^1$ or one 
of its lifts in $\mathbb{R}$). We then obtain 
$$\lim\limits_{n \rightarrow \pm \infty} \frac{1}{pn}[ F^{pn}(x)-x ] =
\lim\limits_{n \rightarrow \pm \infty} \frac{1}{pn}[ (x + nq) -x ] = \frac{q}{p},$$
from where one concludes that 
$$\rho(f) = \frac{q}{p} \hspace{0.5cm} \mathrm{mod} \hspace {0.2cm} 1.$$
Conversely, the following result holds.

\vspace{0.1cm}

\begin{prop} {\em If $f \!\in\! \mathrm{Homeo}_{+}(\mathrm{S}^1)$ 
has rational rotation number, then $f$ has a periodic point.}
\end{prop}

\noindent{\bf Proof.} Due to the equality $\rho(g^m) = m \rho(g)$, it suffices 
to show that every circle homeomorphism with zero rotation number has fixed 
points. To show this, assume that $f \in \mathrm{Homeo}_{+}(\clo)$ has no 
fixed point, and let $F\!: \mathbb{R} \rightarrow \mathbb{R}$ be a lift of $f$ 
such that $F(0) \! \in ]0,1[$. Notice that the function $x \mapsto F(x)-x$ 
has no zero on the line. Therefore, by continuity and periodicity, there exists 
a constant $\delta \! \in ]0,1[$ such that, for all $x \in \mathbb{R}$,
$$\delta \leq F(x) - x \leq 1 - \delta.$$
Letting $x = F^i(0)$ in this inequality, taking the sum 
from $i=0$ up to $i = n-1$, and dividing by $n$, we obtain
$$\delta \leq \frac{F^n(0)}{n} \leq 1 - \delta.$$
Taking the limit as $n$ goes to infinity, the last inequality gives 
\esp $\delta \!\leq\! \rho(f) \!\leq\! 1-\delta.$ \esp Therefore, 
if $f \!\in\! \mathrm{Homeo}_{+}(\mathrm{S}^1)$ has no 
fixed point, then $\rho(f) \neq 0$. $\hfill\square$

\vspace{0.35cm}

We leave to the reader the task of showing that if the rotation number of a circle 
homeomorphism is rational, then all of its periodic points have the same period.

\vspace{0.1cm}

\begin{prop} {\em If two circle homeomorphims are topologically 
conjugate, then their rotation numbers coincide.}
\end{prop}

\noindent{\bf Proof.} We need to show that $\rho(f)=\rho(gfg^{-1})$ 
for every $f,g$ in $\mathrm{Homeo}_+(\clo)$. Let $F$ and $G$ be the 
lifts to the real line of $f$ and $g$, respectively, such that 
$F(0)$ and $G(0)$ belong to $[0,1[$. It is clear that $G^{-1}$ 
is a lift of $g^{-1}$. To show the proposition, we need 
to estimate the value of the expression 
$$\big| (G F G^{-1})^n(x) - F^n(x) \big| =
\big| G F^n G^{-1} (x) - F^n(x) \big|.$$ It is not very difficult 
to check that $|G(x)-x|<2$ \hspace{0.07cm} and \hspace{0.03cm} 
$|G^{-1}(x)-x|<2$ for every $x \in \mathbb{R}$. Moreover, if 
$|x-y|<2$ then $|F^n(x)-F^n(y)|<3$. We thus conclude that 
$$\big| G F^n G^{-1} (x) - F^n(x) \big| \leq \big| G F^n G^{-1} (x) -
F^n G^{-1}(x) \big| + \big| F^n G^{-1}(x) - F^n(x) \big| < 5.$$
Dividing by $n$ and passing to the limit, this clearly yields 
$\rho(f) = \rho(g f g^{-1})$. $\hfill\square$

\vspace{0.25cm}

\begin{obs}
The preceding proposition may be extended to homeomorphisms which 
are semiconjugate. We leave the proof of this fact to the reader.
\end{obs}

As a consequence of what precedes, the dynamics of a circle homeomorphism $f$ 
having rational rotation number $p/q$ (where $p$ and $q$ are relatively prime 
integers) is completely determined by this number, the topology of the set 
$\mathrm{Per}(f)$ of periodic points of $f$, and the ``direction" of the 
dynamics of $f^q$ on each of the connected components of the complement 
of $\mathrm{Per}(f)$.

The case of irrational rotation number is much more interesting. The fact 
that all the orbits of the rotation by an angle $\theta \notin \mathbb{Q}$ 
are dense might suggest the existence of  a unique model (depending on 
$\theta$) for this case. Nevertheless, it is not very difficult to construct 
circle homeomorphisms having irrational rotation number and which are not 
topologically conjugate to the corresponding rotation. To do this, 
choose a point $x \in \mathrm{S}^1$, and  
replace each point $R_{\theta}^i(x)$, $i \!\in\! \mathbb{Z}$, by an interval of length 
$1/2^{|i|}$. We then obtain a ``larger'' circle $\mathrm{S}^1_x$, and the rotation 
$R_{\theta}$ induces a homeomorphism $R_{\theta,x}$ of $\mathrm{S}^1_x$ by extending 
$R_{\theta}$ affinely to each interval that we added to the original circle. The map 
$R_{\theta,x}$ is semiconjugate to $R_{\theta}$, and therefore its rotation number 
equals $\theta$. However, none of the orbits of $R_{\theta,x}$ is dense; in 
particular, $R_{\theta,x}$ is not topologically conjugate to $R_{\theta}$.

Nevertheless, the model is unique up to topological semiconjugacy.

\vspace{0.1cm}

\begin{thm} {\em If the rotation number $\rho(f)$ of $f \!\in\! \mathrm{Homeo}_{+}(\clo)$ 
is irrational, then $f$ is semiconjugate to the rotation of angle $\rho(f)$. The 
semiconjugacy is a conjugacy if and only if all the orbits of $f$ are dense.}
\label{po-ref}
\end{thm}

\noindent{\bf Proof.} Let $F \!: \mathbb{R} \rightarrow \mathbb{R}$ be a lift of 
$f$ to the real line such that $F(0) \in [0,1[$. By Lemma \ref{subaditivo} and 
the proof of Proposition \ref{limitar}, for each $x \in \mathbb{R}$ the value of 
$$\varphi(x) = \sup\limits_{n \in \mathbb{Z}} \big( F^n(x) - n \rho(F) \big)$$
is finite. Moreover, the map $\varphi: \mathbb{R} \rightarrow \mathbb{R}$ 
satisfies the following properties:

\vspace{0.06cm}

\noindent{(i) $\varphi$ is non-decreasing,}

\vspace{0.06cm}

\noindent{(ii) $\varphi(x+1) = \varphi(x) + 1$ \hspace{0.05cm}
for all \hspace{0.05cm} $x \in \mathbb{R}$,}

\vspace{0.06cm}

\noindent{(iii) $\varphi(F(x)) = \varphi(x) + \rho(F)$ \hspace{0.05cm} 
for all \hspace{0.05cm} $x \in \mathbb{R}$.}

\vspace{0.06cm}

From these properties we see that, in order to prove that $\varphi$ is a 
semiconjugacy, we need to show that $\varphi$ is continuous. To do this first notice 
that, for each $x \!\in\! \varphi(\mathbb{R})$, the set $\varphi^{-1}(x)$ is either 
a point or a non-degenerate interval. Let us denote by $\widetilde{\mathrm{Plan}}(F)$ 
the union of the interior of these intervals, and by $\widetilde{\mathrm{Salt}}(F)$ 
the interior of the complement of $\varphi(\mathbb{R})$. 
The sets $\widetilde{\mathrm{Plan}}(F)$ and $\widetilde{\mathrm{Salt}}(F)$ are invariant 
by the integer translations on the line, and therefore they project into 
subsets $\mathrm{Plan}(f)$ and $\mathrm{Salt}(f)$ of the 
circle, respectively. It is easy to 
see that $\mathrm{Salt}(f)$ is invariant under the rotation of angle $\rho(f)$. 
Since $\rho(f)$ is irrational, $\mathrm{Salt}(f)$ must be the empty set, which 
implies that $\varphi$ is continuous, thus inducing a semiconjugacy from $f$ 
to $R_{\rho(f)}$. Finally, notice that $\mathrm{Plan}(f)$ is invariant by $f$. 
Therefore, if the orbits by $f$ are dense, then $\mathrm{Plan}(f)$ is empty. 
In this case $\varphi$ is injective, and thus it induces a conjugacy 
between $f$ and $R_{\rho(f)}$. $\hfill\square$

\vspace{0.35cm}

By the preceding theorem, the combinatorics 
of the dynamics of a homeomorphism of irrational rotation number reduces 
to that of the corresponding rotation. To understand this combinatorics better, 
for $\theta \in [0,1] \setminus \mathbb{Q}$ we define inductively the integers
$$q_1 = 1, \qquad
q_{n+1} = \min \big\{q > q_n \!: dist(q\theta,\mathbb{N}) 
< dist(q_n\theta,\mathbb{N}) \big\}.$$
It is well-known (and easy to check) that the 
sequence $(q_n)_{n \in \mathbb{N}}$ satisfies
\begin{equation}
dist(q_n\theta,\mathbb{N}) = \{q_n\theta\} \esp\esp \mbox{ if and only if }
\esp\esp dist(q_{n+1}\theta,\mathbb{N}) = 1 - \{q_{n+1}\theta\},
\label{alternando}
\end{equation}
where $\{a\} = a - [a]$ denotes the fractional part of $a$ (see for instance \cite{HW}). 
If we project these points on the circle and we keep the same notation, 
property (\ref{alternando}) implies that either 
$$-q_n\theta < q_{n+1}\theta < 0 < -q_{n+1}\theta < q_n\theta < -q_n\theta,$$
or 
$$q_n\theta < -q_{n+1}\theta < 0 < q_{n+1}\theta < -q_n\theta < q_n\theta.$$We 
denote by $I_n$ the interval in $\clo$ of endpoints $0$ and $q_n\theta$ which is  
closed at $0$ and open at $q_n\theta$. The open interval with endpoints 
$-q_n\theta$ and $q_n\theta$ will be denoted by $J_n$. Notice that 
$(q_n+q_{n+1})\theta$ belongs to $J_n$ (see Figures 11 and 12).

Now we claim that the intervals $R_{j\theta}(I_n)$, $j \in \{0,1,\ldots,q_{n+1}-1\}$, 
are disjoint. Indeed, if $R_{j\theta}(I_n) \cap R_{k\theta}(I_n) \neq \emptyset$ 
for some $0 \leq j < k <q_{n+1}$, then $R_{(k-j)\theta}(I_n) \cap I_n \neq \emptyset$, 
which implies that \esp $dist\big( (k-j)\theta,0 \big) < dist(q_n\theta,0)$. \esp 
However, since \esp $k\!-\!j < q_{n+1}$, \esp this 
contradicts the definition of $q_{n+1}$.

From what precedes we deduce that the intervals 
$R_{j\theta}(J_n)$, $j \!\in\! \{ 0,1,\ldots,q_{n+1}-1\}$,
cover the circle, and each point $x \in \clo$ is contained in at most 
two of them. Replacing $\theta$ by $-\theta$, we obtain two sequences 
of intervals $I_{-n}$ and $J_{-n}$ such that $J_n = I_n \cup I_{-n} = J_{-n}$. 
Thus we see that each point of the circle is contained in at most two 
intervals of the form $R_{j\theta}(J_n)$, where 
$j \in \{ -(q_{n+1}-1),-(q_{n+1}-2),\ldots,0 \}$.

\vspace{0.1cm}

The preceding notation is standard and 
will often be used in the next sections.

\vspace{0.3cm}


\beginpicture
\setcoordinatesystem units <0.88cm,0.88cm>

\circulararc 360 degrees from 2.5 0
center at 0 0

\circulararc 40 degrees from 2.8 0
center at 0 0

\circulararc 40 degrees from 2.2 0
center at 0 0

\circulararc -40 degrees from 2.2 0
center at 0 0

\put{$0$} at 2.8 -0.3

\put{$\bullet$} at 2.5 0

\put{$\bullet$} at 2.3 -1

\put{$\bullet$} at 2.4 0.75

\put{$\bullet$} at 1.93 1.6

\put{$\bullet$} at 1.93 -1.6

\put{$I_n$} at 3 1.2

\put{$J_n$} at 1.8 0

\put{$q_n \theta$} at 2.1 2

\put{$-q_n \theta$} at 2.1 -2

\put{$q_{n+1} \theta$} at 2.9 -1.2

\put{$(q_n+q_{n+1})\theta$} at 0.75 0.8


\circulararc 360 degrees from 9.5 0
center at 7 0

\circulararc -40 degrees from 9.8 0
center at 7 0

\circulararc 40 degrees from 9.2 0
center at 7 0

\circulararc -40 degrees from 9.2 0
center at 7 0

\put{$0$} at 9.8 0.3

\put{$\bullet$} at 9.5 0

\put{$\bullet$} at 9.3 1

\put{$\bullet$} at 9.4 -0.75

\put{$\bullet$} at 8.93 -1.6

\put{$\bullet$} at 8.93 1.6

\put{$I_n$} at 10 -1.2

\put{$J_n$} at 8.8 0

\put{$q_n \theta$} at 9.1 -2

\put{$-q_n \theta$} at 9.1 2

\put{$q_{n+1} \theta$} at 9.9 1.2

\put{$(q_n+q_{n+1})\theta$} at 7.8 -0.8

\put{} at -3.8 3

\put{Figure 11} at 0 -3

\put{Figure 12} at 7 -3

\endpicture


\vspace{0.5cm}

\begin{small} \begin{ejer}
Prove that for every $f \in \mathrm{Homeo}_{+}(\clo)$ there exists 
an angle $\theta \in [0,1[$ such that the rotation number of 
\esp $R_{\theta} \circ f$ \esp is non-zero.
\end{ejer}

\begin{ejer}
Given two positive parameters $\lambda_1 > 1$ and $\lambda_2 < 1$, 
let us define 
$a = a(\lambda_1,\lambda_2)$ as the only real number such that 
$\lambda_1 a + \lambda_2 (1-a) = 1$. Consider the (unique) piecewise 
affine \index{piecewise affine homeomorphism} homeomorphism 
$f = f_{\lambda_1,\lambda_2}\!: \clo \rightarrow \clo$ satisfying $f(a)=0$ and 
whose derivative equals $\lambda_1$ on $]0,a[$ and $\lambda_2$ on $]a,1[$  
(this example is due to Boshernitzan \cite{bosher}).

\noindent{(i) For $\sigma = \lambda_1 / \lambda_2$, let 
$h_{\sigma}$ be the homeomorphism  of $[0,1]$ defined by 
$h_{\sigma}(x) = (\sigma^x - 1)/(\sigma - 1)$. Show that 
$h_{\sigma}^{-1} \circ f_{\lambda_1,\lambda_2} \circ h_{\sigma}$
coincides with the rotation $R_{\rho}$, where $\rho$ satisfies 
the equality $\sigma^{\rho} = \lambda_1$.}

\noindent{(ii) Conclude that $\rho (f_{\lambda_1,\lambda_2}) 
= \log(\lambda_1) / ( \log(\lambda_1)-\log(\lambda_2) )$.}
\end{ejer}

\begin{ejer} Let us now fix a positive real number $\sigma \neq 1$, and for 
each $\rho \!\in \esp \! ]0,1[$ let us consider the circle homeomorphism 
$g_{\sigma,\rho} \!=\! h_{\sigma} \circ R_{\rho} \circ h_{\sigma}^{-1}$.

\noindent{(i) Check that $g_{\sigma,\rho}$ 
coincides with $f_{\lambda_1,\lambda_2}$, where
$\lambda_1 = \sigma^{\rho}$ and $\lambda_2 = \sigma^{\rho - 1}$.}

\noindent{(ii) Conclude that inside the group of piecewise affine 
homeomorphisms of the circle, there exist continuous embeddings of 
the group of rotations which are not conjugate to the natural embedding 
by any piecewise affine homeomorphism (see \cite{Min1,Min2} for more on this).}
\end{ejer} \end{small}
\index{rotation number|)}


\subsection{Rotation numbers and invariant measures}
\label{medinv}

\hspace{0.23cm} According to Bogolioubov-Krylov's theorem, 
\index{Bogolioubov-Krylov's theorem}
every circle homeomorphism $f$ preserves a probability measure $\mu$ 
(see Appendix B). Notice that the value of $\mu([x,f(x)[)$ is 
independent of \esp $x \!\in\! \clo$. For instance, if $y \in \clo$ 
is such that \esp $y < f(x) < f(y)$, \esp then 
$$\mu ( [y,f(y)[ ) \!=\! \mu ( [y,f(x)[ ) + \mu ([f(x),f(y)[)
\!=\! \mu ( [y,f(x)[ ) + \mu ( [x,y[ ) \!=\! \mu ([x,f(x)[).$$
This common value will be denoted by $\rho_{\mu}(f)$.

\vspace{0.15cm}

\begin{thm} {\em The value of $\rho_{\mu}(f)$ coincides 
{(\em mod 1)} with the rotation number of $f$.}
\end{thm}

\noindent{\bf Proof.} The measure $\mu$ lifts to a $\sigma$-finite 
measure $\tilde{\mu}$ on $\mathbb{R}$. Let us fix a point $x \in \clo$, 
and let us consider one of its preimages in $\mathbb{R}$, 
which we will still denote by $x$. Notice that every lift $F$ of $f$ 
preserves $\tilde{\mu}$. Moreover, $\tilde{\mu} \big( [x,x+k[ \big) = k$
for all $k \in \mathbb{N}$. Therefore, if $F^n(x) \in [x+k,x+k+1[$, then
$$F^n(x) -x-1 \esp \leq \esp k \esp \leq \esp 
\tilde{\mu} \big( [x,F^n(x)[ \big)
\esp \leq \esp k+1 \esp \leq \esp F^n(x)-x+1.$$
We thus conclude that
$$\lim\limits_{n \rightarrow \infty} \frac{F^n (x)-x}{n} \esp = \esp 
\lim\limits_{n \rightarrow \infty} \frac{\tilde{\mu} \big(
[x,F^{n}(x)[ \big) }{n} \esp = \esp \lim\limits_{n \rightarrow \infty}
\frac{1}{n} \sum_{i=0}^{n-1} \tilde{\mu} \big( [F^{i}
(x),F^{i+1}(x)[ \big) ,$$ 
and since for every $i \in \mathbb{N}$ one has \esp 
$\tilde{\mu} \big( [F^{i} (x),F^{i+1}(x)[ \big) =
\tilde{\mu}\big( [x,F(x)[ \big),$ 
\esp this yields 
$$\lim\limits_{n \rightarrow \infty} \frac{F^n (x)-x}{n} 
\esp = \esp \tilde{\mu} \big( [x,F(x)[ \big).$$
Therefore, (mod 1) the equality 
$\rho(f) = \mu \big( [x,f(x)[ \big) = \rho_{\mu}(f)$ holds. $\hfill\square$

\vspace{0.25cm}

It is very useful to describe the support $supp(\mu)$ of a probability measure 
$\mu$ which is invariant by a circle homeomorphism $f$. If $\rho(f)$ 
is rational, then $f$ has periodic points, and $\mu$ is supported on these 
points. If $\rho(f)$ is irrational, then two cases may occur: if $f$ admits 
an exceptional minimal set\index{exceptional minimal set} 
$\Lambda$, then $supp (\mu) =\Lambda$ and $\mu$ has no atom, 
and if the orbits by $f$ are dense, then the support of 
$\mu$ is the whole circle and $\mu$ has no atom either.

Given a probability measure $\mu$ on the circle, we will denote the group of 
homeomorphisms that preserve $\mu$ by $\Gamma_{\mu}$. Notice that the rotation number 
function restricted to $\Gamma_{\mu}$ is a group homomorphism into $\mathbb{T}^1$. 
Indeed, to check for instance that for every $f,g$ in $\Gamma_{\mu}$ one has \esp 
$\rho_{\mu}(fg) = \rho_{\mu}(f) + \rho_{\mu}(g)$, it suffices to notice that 
$$\rho_{\mu}(fg) \esp = \esp 
\big( [x,fg(x)[ \big) \esp = \esp 
\mu \big( [x,g(x)[ \big) + \mu \big( [g(x),fg(x)[ \big) 
\esp = \esp \rho_{\mu}(f) + \rho_{\mu}(g),$$ 
where the second equality holds (mod 1).

If $\Gamma$ is an amenable \index{group!amenable} 
subgroup of $\mathrm{Homeo}_+(\clo)$, then there exists 
a probability measure on $\clo$ which is invariant by $\Gamma$ (see Appendix B). 
As a consequence, we obtain the following proposition (an alternative proof using 
bounded cohomology appears in \cite{Gh1}).

\vspace{0.08cm}

\begin{prop} The restriction of the rotation number function to every amenable 
subgroup of $\mathrm{Homeo}_+(\clo)$ is a group homomorphism into $\mathbb{T}^1.$
\label{rothomo}
\end{prop}

\vspace{0.001cm}

\begin{small} \begin{ejer} Prove that an irrational rotation is {\bf \textit{uniquely 
ergodic}}, that is, it has a {\em unique} invariant 
probability measure (namely the Lebesgue measure). 
Conclude that if $\Gamma$ is a subgroup of $\mathrm{Homeo}_+(\clo)$ all of whose elements 
do commute with a prescribed minimal circle homeomorphism, then $\Gamma$ is topologically 
conjugate to a group of rotations.\index{rotation!group}
\label{uni-inv}
\end{ejer}

\begin{ejer} Give an example of a subgroup $\Gamma$ of $\mathrm{Homeo}_+(\clo)$ 
which does not preserve any probability measure on the circle, but such that the 
restriction of the rotation number function to $\Gamma$ is a group homomorphism 
into $\mathbb{T}^1$.

\vspace{0.06cm}

\noindent{\underbar{Remark.}} After reading \S \ref{Tits-prob}, the reader 
should be able to prove that, for every group $\Gamma$ satisfying the  
properties above, the image $\rho(\Gamma)$ is finite. 
\end{ejer}\end{small}


\subsection{Faithful actions on the line}
\label{sec-witte}
\index{action!faithful|(}

\hspace{0.23cm} In this section, we will show that the existence of faithful actions on 
the line is closely related to the possibility of endowing the corresponding group with 
a total order relation which is invariant by left-multiplication.

\vspace{0.1cm}

\begin{defn} An order relation $\preceq$ on a group $\Gamma$ is 
{\bf \textit{left-invariant}} (resp. {\bf \textit{right-invariant}}) if for 
every $g,h$ in $\Gamma$ such that $g \preceq h$ one has $fg \preceq fh$ (resp. 
$gf \preceq hf$) for all $f \!\in\! \Gamma.$ The relation is {\bf \textit{bi-invariant}} 
if it is invariant under multiplication on both 
the left and the right, simultaneously. To simplify, we will use 
the term {\bf \textit{ordering}} for a left-invariant total order relation on a group. 
A group is said to be {\bf \textit{orderable}} \index{group!orderable}
(resp. {\bf \textit{bi-orderable}}) \index{group!bi-orderable}
if it admits an ordering (resp. a bi-invariant ordering). 
\end{defn}

\vspace{0.05cm}

For an ordering $\preceq$ on a group $\Gamma$, we will say that 
$f \!\in\! \Gamma$ is {\bf \textit{positive}} if $f \succ id$. 
Notice that the set of these elements forms a semigroup, 
\index{semigroup} 
which is called the {\bf {\em positive cone}} of the ordering. 
An element $f$ is {\bf {\em between g and h}} \esp if \esp 
either $g \prec  f \prec h$ \esp or \esp $h \prec f \prec g$. 

\begin{small}
\begin{ejer} 
Show that every orderable group is torsion-free. \index{group!torsion-free}
\end{ejer}

\begin{ejer}
Prove that a group $\Gamma$ is orderable if and only if it contains 
a subsemigroup $\Gamma_{+}$ such that $\Gamma \setminus \{id\}$ 
is the disjoint union of $\Gamma_{+}$ and the semigroup 
$\Gamma_{-} = \{g\!: g^{-1} \in \Gamma_{+}\}$.
\label{cono-positivo}
\end{ejer}

\begin{ejer} Show that $\{(m,n) \! : m > 0, \esp \mbox{ or } \esp 
m=0 \mbox{ and } n > 0\}$ corresponds to the positive cone of an 
\index{lexicographic ordering}
ordering (called the {\bf \em{lexicographic ordering}}) on $\mathbb{Z}^2$.
\end{ejer}

\begin{ejer} Following the indications below, prove that the free group 
$\mathbb{F}_2$ is bi-orderable.

\vspace{0.08cm}

\noindent (i) Consider the (non-Abelian) ring $\mathbb{A} \! = \! \mathbb{Z} 
\langle X,Y \rangle$ formed by the formal power series with integer coefficients 
in two independent variables $X$, $Y$. Denoting by $o(k)$ the subset of $\mathbb{A}$ 
formed by the elements all of whose terms have degree at least $k$, 
\index{degree!of a power series}
show that \esp $L = 1 + o(1) = \{1 + S \! : \esp S \in o(1)\}$ is a subgroup (under 
multiplication) of $\mathbb{A}$.

\vspace{0.08cm}

\noindent (ii) If $f,g$ are the generators of $\mathbb{F}_2$, prove that the map  
$\phi$ sending $f$ (resp. $g$) to the element $1\! +\! X$ (resp. $1\! +\! Y$) 
in $\mathbb{A}$ extends in a unique way into an injective homomorphism 
$\phi \!: \mathbb{F}_2 \rightarrow L$.  

\vspace{0.08cm}

\noindent (iii) Define a lexicographic type order relation on $L$ which is bi-invariant 
under multiplication by elements in $L$ (notice that this order will be not 
invariant under multiplication by elements in $\mathbb{A}$). Using this order 
and the homomorphism $\phi$, endow $\mathbb{F}_2$ with a bi-invariant ordering.

\vspace{0.06cm}

\noindent{\underbar{Remark.}} The above technique, due to Magnus, 
\index{Magnus} allows easily showing that $\mathbb{F}_2$ is 
{\em residually nilpotent} (see Appendix A). 
\index{group!nilpotent} 
Indeed, it is easy to check that \esp $\Phi (\Gamma_i^{\mathrm{nil}})$ \esp is 
contained in $1 + o(i+1)$ for every $i \geq 0$ (compare Exercise \ref{res-nilp}).
\label{en-el-libre}
\end{ejer} \end{small}

The following theorem gives a dynamical characterization of group orderability.

\vspace{0.1cm}

\begin{thm} {\em For a countable group $\Gamma$, the following are equivalent:

\vspace{0.05cm}

\noindent{{\em (i)} $\Gamma$ acts faithfuly on the line by orientation preserving 
homeomorphisms,}

\vspace{0.05cm}

\noindent{{\em (ii)} $\Gamma$ admits an ordering.}}
\label{thorden}
\end{thm}

\noindent{\bf Proof.} Suppose that $\Gamma$ acts faithfuly on the line by 
orientation preserving homeomorphisms, and let us consider a dense sequence 
$(x_n)$ in $\mathbb{R}$. Let us define \esp $g \prec h$ \esp if the 
smallest index $n$ for which $g(x_n) \neq h(x_n)$ is such that 
$g(x_n) < h(x_n)$. It is not difficult to check 
that $\preceq$ is an ordering. 

Now suppose that $\Gamma$ admits an ordering $\preceq$. Choose a 
numbering $(g_i)$ for $\Gamma$, put $t(g_0) \! = \! 0$, and assume that 
$t(g_0),\ldots,t(g_i)$ have been already defined. If $g_{i+1}$ is greater 
(resp. smaller) than $g_0,\ldots,g_i$, then put $t(g_{i+1}) = \mathrm{max} \{
t(g_0),\ldots,t(g_i) \} + 1$ (resp. $\mathrm{min} \{
t(g_0),\ldots,t(g_i) \} - 1$). Finally, if $g_m \prec g_{i+1}
\prec g_n$ for some $m,n$ in $\{0,\ldots,i\}$, and if $g_j$ 
is not between $g_m$ and $g_n$ for any $0 \leq j \leq i$, 
then define $t(g_{i+1}) = (t(g_m) + t(g_n))/2.$ 
The group $\Gamma$ acts naturally on $t(\Gamma)$ by letting 
$g (t(g_i)) = t(gg_i)$, and this action continuously extends to the 
closure of $t(\Gamma)$. Finally, this action may be extended to the 
whole line in such a way that the map $g$ is affine on each interval 
of the complement of the closure of $t(\Gamma)$. We leave the 
details to the reader. $\hfill\square$

\vspace{0.05cm}

\begin{obs} Notice that the first part of the proof does not use the countability 
assumption. Although this hypothesis is necessary for the second part, many properties 
of orderable groups involve only finitely many elements. To treat such a property, 
one may still use dynamical methods by considering the preceding construction for 
the finitely generated subgroups of the underlying group. 
\end{obs}

If we fix an ordering $\preceq$ on a countable group $\Gamma$, as well as a 
numbering $(g_i)$ of it, then we will call the {\bf \textit{dynamical realization}} 
the action constructed in the proof of Theorem \ref{thorden}. It is easy to 
see that, if $\Gamma$ is nontrivial, then this realization has no global 
fixed point. Another important (and also easy to check) property is the 
fact that, if $f$ is an element of $\Gamma$ whose dynamical realization has 
two fixed points $a\!<\!b$ (which may be equal to $- \infty$ and/or $+ \infty$, 
respectively) such that $]a,b[$ contains no fixed point of $f$, then there 
must exist some point of the form $t(g)$ inside $]a,b[$.

\begin{small} 
\begin{ejer} The construction above is very interesting in the case  
the ordering is bi-invariant, as is shown below.

\noindent (i) Prove that for every element $f$ in the dynamical realization 
of a bi-invariant ordering on a countable group, either $f(x) \leq x$ for 
all $x \in \mathbb{R}$, or $f(x) \geq x$ for every $x \in \mathbb{R}$. 
Conversely, show that every group of homeomorphisms of the real line 
all of whose elements satisfy this property is bi-orderable. 

\noindent (ii) On $\mathrm{PAff}_+ ([0,1])$ define an order relation $\preceq$ by 
letting $f \succ id$ if the right derivative of $f$ at the rightmost point $x_f$ 
such that $f$ coincides with the identity on $[0,x_f]$ is bigger than $1$. Show 
that $\preceq$ is a bi-invariant ordering.

\noindent (iii) From (i) and (ii), conclude that Thompson's group F admits actions 
on the interval which are not semiconjugate to its standard piecewise affine 
\index{piecewise affine homeomorphism} action.

\noindent{\underbar{Remark.}} Bi-orderings on F 
were completely classified in \cite{navas-rivas}.
\label{no-F}
\end{ejer}

\begin{ejer} Give explicit examples of free groups of homeomorphisms 
of the real line to conclude that each $\mathbb{F}_n$ is orderable (compare 
Exercise \ref{en-el-libre} and \cite[Proposition 4.5]{Gh1}).
\end{ejer}\end{small}

For further information concerning orderable groups, we recommend 
\cite{remt,koko,ordering2}. In an opposite direction, the problem of showing that 
some particular classes of groups are non-orderable is also very interesting. An 
important result in this direction, due to Witte-Morris \cite{Wi}, establishes 
that finite index subgroups of $\mathrm{SL}(n,\mathbb{Z})$ are non-orderable 
for $n \geq 3$. (Notice that most of these groups are torsion-free.) 
\index{Witte-Morris!theorem of non orderability of lattices}

\vspace{0.15cm}

\begin{thm} {\em If $n \geq 3$ and $\Gamma$ is a finite index subgroup 
of $\mathrm{SL}(n,\mathbb{Z})$, then $\Gamma$ is not orderable.}
\label{witte}
\end{thm}

\noindent{\bf Proof.} Since $\mathrm{SL}(3,\mathbb{Z})$ injects 
into $\mathrm{SL}(n,\mathbb{Z})$ for every $n \geq 3$, it suffices 
to consider the case $n = 3$. Assume for a contradiction that 
$\preceq$ is an ordering on a finite index subgroup $\Gamma$ of 
$\mathrm{SL}(n,\mathbb{Z})$. Notice that for $k \!\in \!\mathbb{N}$ 
large enough, the following elements must belong to $\Gamma$:
$$\begin{array}
{ccccccccccccccccccccccccccccccccc}
\hspace{1.2cm}
g_1  = \left(
\begin{array}
{ccc}
1 & k & 0 \\
0 & 1 & 0 \\
0 & 0 & 1 \\
\end{array}
\right) ,& &
g_2 = \left(
\begin{array}
{ccc}
1 & 0 & k \\
0 & 1 & 0 \\
0 & 0 & 1 \\
\end{array}
\right) ,& &
g_3  = \left(
\begin{array}
{ccc}
1 & 0 & 0 \\
0 & 1 & k \\
0 & 0 & 1 \\
\end{array}
\right), \\

& & & & & & & & & & & & & & & & & & & & & &   \\

\hspace{1.2cm}
g_4 =
\left(
\begin{array}
{ccc}
1 & 0 & 0 \\
k & 1 & 0 \\
0 & 0 & 1 \\
\end{array}
\right) ,& &
g_5 = \left(
\begin{array}
{ccc}
1 & 0 & 0 \\
0 & 1 & 0 \\
k & 0 & 1 \\
\end{array}
\right) ,& &
g_6 = \left(
\begin{array}
{ccc}
1 & 0 & 0 \\
0 & 1 & 0 \\
0 & k & 1 \\
\end{array}
\right) . \\
\end{array}$$
\noindent It is easy to check that for each $i \in \mathbb{Z}/ 6 \mathbb{Z}$  
the following relations hold:\\ 
$$g_i g_{i+1} = g_{i+1} g_i,
     \qquad [g_{i-1},g_{i+1}] = g_i^k.$$
For $g \!\in\! \Gamma$ we define $|g|\!=\!g$ if $g \!\succeq\! id$, 
and we let $|g|\!=\!g^{-1}$ in the other case. We also write 
$g \gg h$ if $g  \succ  h^n$ for every $n \! \geq \! 1$.

Let us now fix an index $i$, and let us consider the ordering \esp $\preceq$ \esp 
restricted to the subgroup of $\Gamma$ generated by $g_{i-1}$, $g_i$, and $g_{i+1}$. 
One can easily check that it is possible to 
choose three elements $a$, $b$, and $c$, which 
are positive with respect to $\preceq$, such that either 
\hspace{0.025cm} $a = g_{i-1}^{\pm 1}, b = g_{i+1}^{\pm 1}, c = g_{i}^{\pm k}$,
\hspace{0.025cm} or \hspace{0.025cm} $a = g_{i+1}^{\pm 1}$,
$b = g_{i-1}^{\pm 1}, c = g_{i}^{\pm k}$, \esp and such that 
$$a c = c a , \quad b c = c b ,
\quad a b a^{-1} b^{-1} = c^{-1}.$$
We claim that either \esp $a \gg c$ \esp or \esp $b \gg c$. \esp To show this, 
assume that for some $n \geq 1$ one has $c^{n} \succ a \hspace{0.1cm}$ and 
$\hspace{0.1cm} c^{n} \succ b$. Let $d_m = a^{m} b^{m} (a^{-1} c^{n})^{m} 
(b^{-1} c^{n})^{m}.$ Since $d_m$ is a product of positive elements, $d_m$ is 
positive. On the other hand, it is not difficult to check that $d_m = c^{-m^2 + 2mn}$, 
and therefore \esp $d_{m} \prec id$ \esp for $m$ large enough, which is a contradiction.

The claim above allows us to conclude that either $|g_i| \ll |g_{i-1}|$ or 
$|g_i| \ll |g_{i+1}|$. If we assume that $|g_1| \ll |g_2|$, then we obtain 
$|g_1| \ll |g_2| \ll |g_3| \ll |g_4| \ll |g_5| \ll |g_6| \ll |g_1|,$ which is 
a contradiction. The case where $|g_1| \gg |g_2|$ is analogous. $\hfill\square$

\vspace{0.3cm}

It follows from an important theorem due to Margulis that, for $n \geq 3$, 
\index{Margulis!normal subgroup theorem}
every normal subgroup of a finite index subgroup of $\mathrm{SL}(n,\mathbb{Z})$ 
is either finite or of finite index (see \cite{Ma2}). As a corollary, we obtain 
the following strong version of Witte-Morris' theorem.

\vspace{0.1cm}

\begin{thm} For $n \geq 3$, every action of a finite index subgroup of 
$\mathrm{SL}(n,\mathbb{Z})$ by homeomorphisms of the line is trivial.
\label{wittefuerte}
\end{thm}

\begin{small} \begin{ejer} Prove directly that every torsion-free nilpotent 
group is bi-orderable. More generally, prove that the same is true for every 
residually nilpotent group $\Gamma$ for which the (Abelian) quotients \esp 
$\Gamma^{\mathrm{nil}}_{i} / \Gamma^{\mathrm{nil}}_{i+1}$ \esp are 
torsion-free (see Appendix A).

\noindent{\underbar{Hint.}} Use Exercise \ref{se-uso}.
\label{res-nilp}
\end{ejer}

\begin{ejer} Let $\mathrm{N}_n$ be the group of $n \times n$ upper triangular 
matrices with integer entries and such that each entry in the diagonal equals 
1. Prove that $\mathrm{N}_n$ is orderable by using its natural action on 
$\mathbb{Z}^n$ and the lexicographic order on $\mathbb{Z}^n$.\\

\noindent{\underbar{Remark.}} It is easy to check that $\mathrm{N}_n$ is nilpotent 
and torsion-free. On the other hand, according to a classical result due to Malcev 
\cite{Ra}, every nilpotent, finitely generated, torsion-free group embeds into 
$\mathrm{N}_n$ for some $n$. This allows to reobtain indirectly the first 
claim of the preceding exercise.
\label{columnas}
\end{ejer} \end{small}
\index{action!faithful|)}


\subsection{Free actions and H\"older's theorem\index{H\"older!theorem|(}}
\index{action!free|(}

\hspace{0.23cm} The main results of this section are classical and essentially due 
to H\"older. Roughly, they state that free actions on the line exist only for groups 
admitting an order relation satisfying an Archimedean 
\index{Archimedean property} type property. Moreover, 
these groups are necessarily isomorphic to subgroups of $(\mathbb{R},+)$, and 
the corresponding actions are semiconjugate to actions by translations.

\vspace{0.035cm}

\begin{defn} An ordering $\precede$ on a group $\Gamma$ is said to be {\bf 
\textit{Archimedean}} if for all $g, h$ in $\Gamma$ such that $g \!\neq\! id$ 
there exists $n \!\in\! \mathbb{Z}$ satisfying $g^n \!\succ\! h.$\\
\end{defn}

\vspace{0.02cm}

\begin{prop} \textit{If \esp $\Gamma$ is a group acting freely by homeomorphisms 
of the real line, then $\Gamma$ admits a bi-invariant Archimedean ordering.} 
\label{existeordenarq}
\end{prop}

\noindent{\bf Proof.} Let $\precede$ be the left-invariant order relation on $\Gamma$ 
defined by $g \prec h$ if $g(x) < h(x)$ for some (equivalently, for all) $x \!\in\! 
\mathbb{R}$. This order relation is total, and using the fact that the action is 
free, one easily checks that it is also right-invariant and Archimedean. $\hfill\square$

\vspace{0.45cm}

The converse to the proposition above is a direct consequence of the following one. 
As we will see in Exercises \ref{conrid1} and \ref{conrid2}, the hypothesis of 
bi-invariance can be weakened, and left-invariance is sufficient.

\vspace{0.15cm}

\begin{prop} \textit{Every group admitting a bi-invariant Archimedean 
ordering is isomorphic to a subgroup of \esp $(\mathbb{R},+)$.} 
\label{completa}
\end{prop}

\noindent{\bf Proof.} Assume that a nontrivial group $\Gamma$ admits a bi-invariant Archimedean 
ordering $\precede$, and let us fix a positive element $f \in \Gamma$. For each $g \in \Gamma$ 
and each $p \in \mathbb{N}$, let us consider the unique integer $q=q(p)$ such that 
$\hspace{0.15cm} f^q \precede g^p \prec f^{q+1}$.

\vspace{0.25cm}

\noindent{\underbar{Claim (i).}} The sequence \esp $q(p) / p$ \esp 
converges to a real number as \esp $p$ \esp goes to infinity.

\vspace{0.15cm}

Indeed, if $\hspace{0.15cm} f^{q(p_1)} \precede g^{p_1} \prec
f^{q(p_1)+1} \hspace{0.15cm}$ and $\hspace{0.15cm} f^{q(p_2)}
\precede g^{p_2} \prec f^{q(p_2)+1}, \hspace{0.15cm}$ then
$$f^{q(p_1)+q(p_2)} \precede g^{p_1+p_2} \prec f^{q(p_1)+q(p_2)+2}.$$
Therefore, 
$\hspace{0.15cm} q(p_1) + q(p_2) \leq q(p_1+p_2) \leq
q(p_1)+q(p_2)+1. \hspace{0.15cm}$ The convergence of the sequence 
$(q(p) / p)$ to some point in $[-\infty,\infty[$ then follows 
from Lemma \ref{subaditivo}. Moreover, if we denote by 
$\phi(g)$ the limit of \esp $q(p)/p$, \esp then for the integer 
$n \!\in\! \mathbb{Z}$ satisfying \esp $f^n \precede g \prec f^{n+1}$ 
\esp one has \esp $ f^{np} \precede g^p \prec f^{(n+1)p}$, \esp and therefore
$$n = \lim\limits_{p \rightarrow \infty} \frac{np}{p} \leq \phi(g)
\leq \lim\limits_{p \rightarrow \infty} \frac{(n+1)p - 1}{p} = n+1.$$

\vspace{0.25cm}

\noindent{\underbar{Claim (ii).}} The map $\phi: \Gamma
\rightarrow (\mathbb{R},+)$ is a group homomorphism.\\

\vspace{0.15cm}

Indeed, let $g_1,g_2$ be arbitrary elements in $\Gamma$. Let us suppose that 
$g_1g_2 \precede g_2g_1$ (the case where $g_2g_1 \precede g_1g_2$ is analogous). 
Since $\preceq$ is bi-invariant, if \esp $f^{q_1} \precede g_1^p \prec f^{q_1 +1}$ 
\esp and \esp $f^{q_2} \precede g_2^p \prec f^{q_2 +1}$ \esp then
$$\hspace{0.15cm} f^{q_1+q_2} \precede g_1^p g_2^p \precede
(g_1g_2)^p \precede g_2^p g_1^p \prec f^{q_1 + q_2 + 2} \hspace{0.15cm}.$$
From this one concludes that 
$$\phi(g_1) + \phi(g_2) = \lim\limits_{p \rightarrow \infty}
\frac{q_1+q_2}{p} \leq \phi(g_1g_2) \leq \lim\limits_{p
\rightarrow \infty} \frac{q_1+q_2 +1}{p} = \phi(g_1)+\phi(g_2),$$ 
and therefore $\phi(g_1g_2) = \phi(g_1) + \phi(g_2).$

\vspace{0.25cm}

\noindent{\underbar{Claim (iii).}} The homomorphism $\phi$ is one to one.\\

\vspace{0.15cm}

Notice that $\phi$ is order preserving, in the sense that if $g_1 \precede g_2$ 
then $\phi(g_1) \leq \phi(g_2)$. Moreover, $\phi(f)=1.$ Let $h$ be an element 
in $\Gamma$ such that $\phi(h)=0$. Assume that $h \neq id$. Then there exists 
$n \!\in\! \mathbb{Z}$ such that $h^n \succeq f$. From this one concludes that 
$0 = n \phi(h) = \phi(h^n) \geq \phi(f) = 1,$ which is absurd. Therefore, 
if $\phi(h) = 0$ then $h = id$, and this concludes the proof. $\hfill\square$

\vspace{0.65cm}

If $\Gamma$ is an infinite group acting freely on the line, then we can endow it with the order 
relation introduced in the proof of Proposition \ref{existeordenarq}. This order allows us 
to construct an embedding $\phi$ from $\Gamma$ into $(\mathbb{R},+)$. If $\phi(\Gamma)$ 
is isomorphic to $(\mathbb{Z},+)$, then the action of $\Gamma$ is conjugate to the  
action by integer translations. In the other case, the group $\phi(\Gamma)$ is 
dense in $(\mathrm{R},+)$. For each point $x$ in the line we then define
$$\varphi(x) = \sup \{ \phi(h) \in \mathbb{R} \!: h(0) \leq x \}.$$
It is easy to see that $\varphi \!\! : \mathbb{R} \rightarrow \mathbb{R}$ 
is a non-decreasing map. Moreover, it satisfies the equality \esp 
$\varphi(h(x)) = \varphi(x) + \phi(h)$ \esp for all $x \!\in\! \mathbb{R}$ 
and all $h \in \Gamma$. Finally, $\varphi$ is continuous, since 
otherwise the set $\mathbb{R} \setminus \varphi(\mathbb{R})$ would 
be a non-empty open set invariant by the translations of $\phi(\Gamma)$, 
which is impossible.

\vspace{0.15cm}

To summarize, if $\Gamma$ is a group acting freely on the line,  
then its action is semiconjugate to an action by translations.

\vspace{0.18cm}

\begin{small} \begin{ejer} If $\Gamma$ admits an Archimedean ordering, then 
this ordering is necessarily bi-invariant. To show this claim (first 
remarked by Conrad \cite{conrad}) 
\index{Conrad} 
just follow the following steps.

\vspace{0.08cm}

\noindent (i) Prove that an ordering $\preceq$ is bi-invariant if and only if its positive 
cone is a normal subsemigroup, that is, 
$hgh^{-1}$ belongs to $\Gamma_+$ for every $g \in \Gamma_+$ and every 
$h \in \Gamma$ (see Exercise \ref{cono-positivo}).

\vspace{0.08cm}

\noindent (ii) Let $\preceq$ be an Archimedean ordering on a group $\Gamma$. Suppose 
that $g \in \Gamma_+$ and $h \in \Gamma_{-}$ are such that $hgh^{-1} \notin \Gamma_+$,  
and consider the smallest positive integer $n$ for which $h^{-1} \prec g^n$. 
Using the relation \esp $hgh^{-1} \prec id$, \esp show that 
$h^{-1} \prec g^{-1} h^{-1} \prec g^{n-1}$, thus contradicting 
de definition of $n$. Conclude that $\Gamma_+$ is stable under 
conjugacy by elements in $\Gamma_{-}$.

\vspace{0.08cm}

\noindent (iii) Assume now that $g \in \Gamma_+$ and $h \in \Gamma_+$ verify \esp 
$hgh^{-1} \notin \Gamma_+$. \esp In this case, $hg^{-1}h^{-1} \succ id$, and since 
$h^{-1} \in \Gamma_{-}$, by (ii) one has \esp $h^{-1}(hg^{-1}h^{-1})h \in \Gamma_+$, 
\esp that is, $g^{-1} \in \Gamma_+$, which is absurd.
\label{conrid1}
\end{ejer} 


\begin{ejer} As an alternative argument to that of the preceding exercise, 
show that if a countable group is endowed with an Archimedean ordering, 
then the action of the corresponding dynamical realization is free.
\label{conrid2}
\end{ejer} \end{small}

\vspace{0.1cm}

Let us now consider a group $\Gamma$ acting freely by circle homeomorphisms. In this case, 
the preimage $\tilde{\Gamma}$ of $\Gamma$ in $\widetilde{\mathrm{Homeo}_{+}}(\mathrm{S}^1)$ 
acts freely on the line. If we repeat the arguments of the proof of Proposition \ref{completa} 
by considering the translation $x \mapsto x+1$ as being the positive element $f$, then 
one obtains that $\tilde{\Gamma}$ is isomorphic to a subgroup of $(\mathbb{R},+)$, 
and this isomorphism projects to an isomorphism between $\Gamma$ and a subgroup 
of $\mathrm{SO}(2,\mathbb{R})$. We record this fact as a theorem.

\vspace{0.15cm}

\begin{thm} {\em If $\Gamma$ is a group acting freely by circle homeomorphisms, 
then $\Gamma$ is isomorphic to a subgroup of \esp $\mathrm{SO}(2,\mathbb{R})$.}
\end{thm}

\vspace{0.15cm}

As in the case of the real line, under the above hypothesis the 
action of $\Gamma$ is semiconjugate to that of the corresponding group of 
\index{rotation!group} rotations.\\

\vspace{0.1cm}

\begin{small} \begin{ejer} Prove that if $\Gamma$ is a finitely generated subgroup of  
$\mathrm{Homeo}_{+}(\mathrm{S}^1)$ all of whose elements are torsion, then $\Gamma$ 
is finite.

\noindent{\underbar{Remark.}} It is unknown whether the same is true 
for torsion subgroups of $\mathrm{Homeo}_{+}(\mathrm{S}^2)$.
\end{ejer}

\begin{ejer} Let $\Gamma$ be a subgroup of $\mathrm{PSL}(2,\mathbb{R})$ all of 
whose elements are elliptic. Prove that $\Gamma$ is conjugate to a group of rotations.
\label{asi-no-mas}
\end{ejer}

\begin{ejer} Give an alternative proof for Proposition \ref{gatorot} using 
H\"older's theorem.
\end{ejer}\end{small}

\vspace{0.15cm}

The preceding results show that free actions on the circle or the real line are 
topologically semiconjugate to the actions of groups of rotations or translations, 
respectively. In a similar direction, remark that nontrivial elements of the affine 
group fix at most one point. In what follows, we will see that, up to topological 
semiconjugacy and discarding a degenerate case, this property characterizes the 
affine group. 
\index{affine group}\index{Solodov}The next result is due to Solodov.

\vspace{0.1cm}

\begin{thm} {\em Let $\Gamma$ be a subgroup of $\mathrm{Homeo}_{+}(\mathbb{R})$ 
such that every nontrivial element in $\Gamma$ fixes at most one point. Suppose 
that there is no global fixed point for the action. Then $\Gamma$ is 
semiconjugate to a subgroup of the affine group.} 
\label{solodov}
\end{thm}

Notice that the case we are not considering, namely when there exists a 
point $x_0$ which is fixed by every element, the actions of $\Gamma$ on 
$]\!-\!\infty,x_0[$ and $]x_0,\infty[$ are free, and therefore semiconjugate 
to actions by translations. However, the action of $\Gamma$ on the line is 
not necessarily conjugate to the action of a subgroup of the stabilizer of 
some point in the affine group, since $\Gamma$ may contain elements 
for which $x_0$ is a parabolic fixed point...

\vspace{0.5cm}

\noindent{\bf Proof of Theorem \ref{solodov}.} If the action is free, then 
the claim of the theorem follows from H\"older's theorem. We will assume 
throughout that the action is not free. Under this assumption, $\Gamma$ 
cannot be Abelian. Indeed, if $\Gamma$ were Abelian, then the orbit 
of a fixed point $x_0$ of a nontrivial element $g \in \Gamma$ would be 
contained in the set of fixed points of $g$. Thus, $x_0$ would be 
a global fixed point of $\Gamma$, which is a contradiction.

\vspace{0.21cm}

\noindent{\underbar{First step.}} We claim that if $\Gamma_0$ is a normal 
subgroup of $\Gamma$ containing a nontrivial element with a fixed point, 
then $\Gamma_0$ has an element with an attracting fixed point.

\vspace{0.1cm}

Indeed, let $h_0 \!\in\! \Gamma_0$ be a nontrivial element such that $h_0(x_0) \!=\! x_0$ 
for some point $x_0 \!\in\! \mathbb{R}$. Suppose that $x_0$ is a parabolic fixed point of 
$h_0$. Replacing $h_0$ by its inverse if necessary, we may assume that $h_0(y) > y$ for 
$y \neq x_0$. By hypothesis, there exists an element $g \in \Gamma$ such that 
$x_1 = g(x_0) \neq x_0$. Changing $g$ by $g^{-1}$ if necessary, we may assume that 
$x_1 \!>\! x_0$. Let us consider the elements $h_1 = gh_0g^{-1} \in \Gamma_0$ and 
$h = h_0h_1^{-1} \in \Gamma_0$. It is easy to see that $h(x_0) < x_0$ and $h(x_1) > x_1$. 
Thus, $h$ has a repelling fixed point in $]x_0,x_1[$, which proves the claim.

\vspace{0.21cm}

\noindent{\underbar{Second step.}} The definition of an order relation.\\

\vspace{0.1cm}

Given $g,h$ in $\Gamma$, we write $g \preceq h$ if there exists $x \!\in\! \mathbb{R}$
such that $g(y) \leq h(y)$ for every $y \geq x$. It is easy to check that this defines 
a total and bi-invariant order relation. We claim that this ordering 
satisfies the following weak form of the \index{Archimedean property}
Archimedean property: if $f \in \Gamma$ has a repelling fixed point and $g \in \Gamma$, 
then there exists $n \!\in\! \mathbb{N}$ such that $g \preceq f^n$. Indeed, letting $x_0$ 
be the fixed point of $f$, let $x_{-},x_+$ be such that $x_{-} < x_0 < x_{+}$. For 
$n \in \mathbb{N}$ large enough we have $f^n(x_{-}) < g(x_{-})$ and $f^n(x_{+}) > g(x_{+})$, 
and therefore $g^{-1}f^n$ has a fixed point in the interval $]x_{-},x_{+}[$. Since 
$g^{-1}f^n(x_{+}) > x_{+}$, this implies that $f^n(x) > g(x)$ for all 
$x \geq x_{+}$, and thus $g \preceq f^n$.

\vspace{0.21cm}

\noindent{\underbar{Third step.}} A homomorphism into the reals.\\

\vspace{0.1cm}

Let us fix an element $f \!\in\! \Gamma$ with a repelling fixed point. As in the proof 
of H\"older's theorem, for $g \!\in\! \Gamma$ such that $g \succeq id$ we define
$$\phi (g) = \lim\limits_{p \rightarrow \infty}  \left\{  \frac{q}{p}:
f^q \preceq g^p \prec f^{q+1} \right\},$$
and for $g \preceq id$ we let $\phi (g^{-1}) = - \phi (g)$. The map 
$\phi\!: \Gamma \rightarrow (\mathbb{R},+)$ is a group homomorphism 
satisfying $\phi(f) = 1$. Notice 
that if $\Gamma_0$ is a normal subgroup of $\Gamma$ containing a nontrivial element having 
a fixed point, then by the first step of the proof there exists $h \in \Gamma_0$ with a 
repelling fixed point. By the second step we have $f \preceq h^n$ for $n\!\in\!\mathbb{N}$ 
large enough, and thus $\phi (h) \geq 1/n$. In particular, $\phi(\Gamma_0) \neq \{0\}$.

\vspace{0.21cm}

\noindent{\underbar{Fourth step.}} The action 
of $[\Gamma,\Gamma]$ on the line is free.\\

\vspace{0.1cm}

Indeed, $[\Gamma,\Gamma]$ being normal in $\Gamma$, if 
$\Gamma$ has a nontrivial element with a fixed point then 
$\phi \big( [\Gamma,\Gamma] \big) \neq 0$. Nevertheless, this is 
absurd, since $[\Gamma,\Gamma]$ is contained in the kernel of $\phi$.

Therefore, $[\Gamma,\Gamma]$ is semiconjugate to a group of 
translations, that is, there exists a group homomorphism 
$\phi_0\!\!: [\Gamma,\Gamma] \rightarrow (\mathbb{R},+)$ and a 
continuous non-decreasing surjective map $\varphi$ of the line 
such that \esp $\varphi(h(x)) = \varphi(x) + \phi_0(h)$ \esp 
for all $x \!\in\! \mathbb{R}$ and all $h \!\in\! [\Gamma,\Gamma]$.
We claim that $\phi_0 \big( [\Gamma,\Gamma] \big)$ is non 
discrete. If not, the conjugacies 
by elements in $\Gamma$ would preserve the generator of 
$\phi_0 \big( [\Gamma,\Gamma] \big) \! \sim \! \mathbb{Z}$, and so 
$[\Gamma,\Gamma]$ would be contained in the center of $\Gamma$. 
\index{center of a group} 
However, this is impossible, because $\Gamma$ 
contains elements having one fixed point.

\vspace{0.21cm}

\noindent{\underbar{Fifth step.}} End of the proof.\\

\vspace{0.1cm}

By the above, the image $\phi_0 \big( [\Gamma,\Gamma] \big)$ is dense in $\mathbb{R}$. 
It is then easy to see that the conjugacy $\varphi$ of the fourth step is unique up to 
composition by the right with elements in the affine group. Since $[\Gamma,\Gamma]$ is 
normal in $\Gamma$, for each $g \!\in\! \Gamma$ the homeomorphism $\varphi g \varphi^{-1}$ 
belongs to the affine group, and this finishes the proof. 
$\hfill\square$

\vspace{0.25cm}

Once again, we emphasize that neither H\"older's 
nor Solodov's theorem can be extended in a natural way 
to groups acting (even minimally and smoothly) on the circle so 
that every nontrivial element fixes at most two points \cite{Kov,Na-kova}.

\index{H\"older!theorem|)} \index{action!free|)}


\subsection{Translation numbers and quasi-invariant measures} 
\label{muylargo}
\index{measure!quasi-invariant|(}

\hspace{0.23cm} As in the case of actions on the circle, for a Radon 
\index{measure!Radon}
measure $\upsilon$ on the line we may consider the group $\Gamma_{\upsilon}$ of 
the homeomorphisms preserving it. For $g \in \Gamma_{\upsilon}$ we define 
its {\bf \textit{translation number}} with respect to $\upsilon$ by letting  
$$\tau_{v}(g) = \left \{ \begin{array} {l} 
v([x,g(x)[) \hspace{0.67cm} \mbox{ if } \esp g(x) > x,\\ 
0 \hspace{2.15cm} \mbox{ if } \esp g(x) = x,\\
- v([g(x),x[) \hspace{0.4cm} \mbox{ if } \esp g(x) < x. \end{array} \right.$$
It is easy to see that this number does not depend on the choice of 
$x\in\mathbb{R}$.\index{translation number|(}

The translation number satisfies many properties 
which are similar to those of the rotation number of circle 
homeomorphisms. For instance, for $g \!\in\! \Gamma_{\upsilon}$ one has 
\begin{equation}
\tau_{\upsilon}(g)=0 \hspace{0.05cm} \mbox{ if and only if }
g \mbox{ has a fixed point.}
\label{puntofijoplante}
\end{equation}
Indeed, if $\mathrm{Fix}(g) = \emptyset$ then the orbit of every point $x$ in 
the line is unbounded from both sides. Let us fix $x \in \mathbb{R}$, and let us 
assume that $g(x) > x$ (if this is not the case, then we may change $g$ by $g^{-1}$). 
If $\tau_{\upsilon}(g) = 0$, then letting $n$ go to infinity in the equality  
$\upsilon ([x,g^n(x)[) = \upsilon ([g^{-n}(x),x[) = 0$ we conclude 
that $\upsilon (]-\infty,+\infty[) = 0$, which is absurd. 
Conversely, 
if $\mathrm{Fix}(g)$ is nonempty, then by definition we have $\tau_{\upsilon}(g)=0.$

Let us remark that a stronger property holds for elements $g \in \Gamma_{\upsilon}$, 
namely
\begin{equation}
\mbox{if } \hspace{0.15cm} \mathrm{Fix}(g)
\neq \emptyset \hspace{0.15cm} \mbox{ then }
\hspace{0.1cm} supp (\upsilon) \subset \mathrm{Fix}(g).
\label{soporteplante}
\end{equation}
Indeed, if $supp (\upsilon)$ is not contained in $\mathrm{Fix}(g)$, 
then there is a positive $\upsilon$-measure set $A$ 
contained in a connected component of the complement of 
$\mathrm{Fix}(g)$ such that $A \cap g(A) = \emptyset$. At least one of the 
sets $\cup_{n \in \mathbb{N}} g^n(A)$ or $\cup_{n \in \mathbb{N}} g^{-n}(A)$ 
must be bounded, and therefore of finite $v$-measure. However, the $v$-measure 
of these sets equals $\sum_{n \in \mathbb{N}} \upsilon (A) = \infty.$

Notice that the function $\tau_{\upsilon} \! : \Gamma_{\upsilon}
\rightarrow \mathbb{R}$ is a group homomorphism from $\Gamma_{\upsilon}$ into 
$(\mathbb{R},+)$. This property will be very important for dealing with the 
problem of the uniqueness (up to a scalar factor) of the invariant Radon measure. 

\vspace{0.12cm}

\begin{lem} {\em If $\upsilon_1$ and $\upsilon_2$ are Radon measures which are invariant  
by a subgroup $\Gamma$ of \esp $\mathrm{Homeo}_{+}(\mathbb{R})$, then there exists 
$\kappa > 0$ such that the homomorphisms $\tau_{\upsilon_1}$ and \esp $\tau_{\upsilon_2}$ 
satisfy the relation $\tau_{\upsilon_1} = \kappa \tau_{\upsilon_2}$.}
\label{igualtau}
\end{lem}

\noindent{\bf Proof.} Due to (\ref{puntofijoplante}), the kernels of $\tau_{\upsilon_1}$ 
and $\tau_{\upsilon_2}$ coincide with $\Gamma_0 = \{ g \!: \mathrm{Fix}(g) \neq \emptyset\}$. 
We then dispose of two homomorphisms $\tau_1$ and $\tau_2$ from $\Gamma / \Gamma_0$ into 
$(\mathbb{R},+)$. Let us fix a point $x_0$ in $\mathrm{Fix}(\Gamma_0)$ (the existence of 
such a point is ensured by (\ref{soporteplante})). The group $\Gamma/\Gamma_0$ acts 
\index{action!free}
freely on the orbit $\Gamma(x_0)$, hence it admits an Archimedean 
\index{Archimedean property} ordering $\preceq$, namely the 
one given by $g_1\Gamma_0 \prec g_2 \Gamma_0 $ if $g_1(x_0) < g_2(x_0)$. Notice that both 
$\tau_{\upsilon_1}$ and $\tau_{\upsilon_2}$ preserve this order. If we fix an 
element $f \in \Gamma$ such that $\Gamma_0 \prec f \Gamma_0$, it is easy 
to see that for every $g \in \Gamma$ one has
$$\tau_i (g \Gamma_0) = \tau_i(f \Gamma_0) \cdot \lim\limits_{p \rightarrow \infty}
\left\{ \frac{q}{p}: \hspace{0.1cm} f^q \Gamma_0 \preceq g^p \Gamma_0 \prec f^{q+1}
\Gamma_0 \right\}.$$
Thus, $\tau_2(f) \cdot \tau_{v_1} = \tau_1(f)  \cdot \tau_{v_2} $, 
which concludes the proof. $\hfill\square$

\vspace{0.35cm}

In the case where $\Gamma$ preserves $\upsilon$ and $\tau_{\upsilon}(\Gamma)$ 
is trivial or isomorphic to $\mathbb{Z}$, one cannot expect having uniqueness 
(up to a scalar factor) of the invariant measure $\upsilon$. However, the case 
where $\tau_{\upsilon}(\Gamma)$ is dense in $\mathbb{R}$ is distinct.

\vspace{0.1cm}

\begin{prop} {\em If $\upsilon_1$ and $\upsilon_2$ are two Radon measures which 
are invariant by a subgroup $\Gamma$ of $\mathrm{Homeo}_{+}(\mathbb{R})$ so that  
that $\tau_{\upsilon_1}(\Gamma)$ and $\tau_{\upsilon_2}(\Gamma)$ are dense in 
$\mathbb{R}$, then there exists $\kappa > 0$ such that these measures satisfy 
the relation $\upsilon_1 = \kappa \esp \upsilon_2$.}
\label{igualmedida}
\end{prop}

\noindent{\bf Proof.} By the preceding lemma, after normalization  
we may assume that 
$\tau_{\upsilon_1} \!=\! \tau_{\upsilon_2}$. We will then show that 
$\upsilon_1 \!=\! \upsilon_2$. For this, first notice that none of these 
measures has atoms. Indeed, if $\upsilon_i \big( \{ x_0 \} \big)\!>\!0$ 
then every positive element in $\tau_{\upsilon_i}(\Gamma)$ 
would be greater than or equal to $\upsilon_i \big( \{ x_0 \} \big) > 0$,
which contradicts the fact that $\tau_{\upsilon_i}(\Gamma)$ is dense. A 
similar argument shows that the actions of $\Gamma$ on the supports of 
$\upsilon_1$ and $\upsilon_2$ are minimal. 

Now we show that the supports $supp(\upsilon_1)$ and $supp(\upsilon_2)$
are actually equal. Indeed, in case of non-equality there would be a point 
$x \!\in\! supp(\upsilon_i) \setminus supp(\upsilon_{i+1})$ 
(where $i \!\in\! \mathbb{Z} / 2 \mathbb{Z}$). 
By the density of the orbits on the supports, 
we could then choose $g \!\in\! \Gamma$ such that 
$g(x) \neq x$ and such that the $\upsilon_{i+1}$-measure of the interval of endpoints 
$x$ and $g(x)$ is zero. However, this would imply that $\tau_{\upsilon_{i+1}}(g) = 0$ 
and $\tau_{\upsilon_i}(g) \neq 0$, which is absurd.\\

To finish the proof of the equality between $\upsilon_1$ and $\upsilon_2$, we need to show 
that they give the same mass to intervals having endpoints in their common support. 
If $[x,y]$ is an interval of this type, we may choose $g_n \in \Gamma$ such that 
$g_n(x)$ converges to $y$. We then have 
\begin{eqnarray*}
\upsilon_1 \big( [x,y] \big) &=& \lim\limits_{n \rightarrow \infty}
\upsilon_1 \big( [x,g_n(x)] \big) \esp = \esp \lim\limits_{n \rightarrow
\infty} \tau_{\upsilon_1}(g_n)\\
&=& \lim\limits_{n \rightarrow \infty}
\tau_{\upsilon_2}(g_n) \esp 
= \esp \lim\limits_{n \rightarrow \infty} \upsilon_2 \big( [x,g_n(x)] \big) 
\esp = \esp \upsilon_2 \big( [x,y] \big), 
\end{eqnarray*} 
thus concluding the proof. $\hfill\square$

\vspace{0.25cm}

The preceding discussion shows how important is to know \textit{a priori} what 
conditions ensure the existence of an invariant Radon measure. The following 
result, due to Plante \cite{Pl3}, is an important issue in this direction.
\index{Plante!invariant measure theorem}

\vspace{0.1cm}

\begin{thm} {\em If \esp $\Gamma$ is a finitely generated 
virtually nilpotent subgroup of \esp $\mathrm{Homeo}_{+}(\mathbb{R})$, 
\index{group!nilpotent} then $\Gamma$ preserves a Radon measure on the line.}
\label{plante}
\end{thm}

\vspace{0.1cm}

We immediately state a corollary of this result which will be useful 
later on.

\vspace{0.1cm}

\begin{cor} {\em If \esp $\Gamma$ is a finitely generated virtually nilpotent subgroup 
of $\mathrm{Homeo}_{+}(\mathbb{R})$, then the action of the commutator subgroup 
$[\Gamma,\Gamma]$ has global fixed points.}
\label{conmutptofijo}
\end{cor}

\vspace{0.1cm}

Indeed, if $\upsilon$ is an invariant Radon measure, then the translation number of every 
element of $[\Gamma,\Gamma]$ with respect to $\upsilon$ is zero. The claim of the 
corollary then follows from (\ref{puntofijoplante}).

\vspace{0.1cm}

The original proof by Plante of Theorem \ref{plante} involves very interesting ideas 
related to growth of groups (see \S \ref{seccion-crecer}), the notion of pseudo-group 
(see \S \ref{psg}), and group amenability (see Appendix B), which apply in 
more general situations. Nevertheless, as an application of the methods 
from \S \ref{sec-minimal}, we will give a more direct proof based on 
the fact that nilpotent groups do not contain free semigroups.
\index{semigroup!free}

\begin{defn} Two elements $f,g$ in a group $\Gamma$ generate a {\bf {\em free semigroup}} 
if the elements of the form \esp $f^{n}g^{m_r}f^{n_r} \cdots g^{n_1}f^{n_1} g^{m}$, \esp 
where $n_j$ and $m_j$ are positive integers, $m \geq 0$, and $n \geq 0$, are two-by-two 
different for different choices of the exponents.
\label{def-semi-librre}
\end{defn}

\begin{small} \begin{ejer} Prove that virtually nilpotent groups do not contain free 
semigroups on two generators.

\noindent{\underbar{Hint.}} For nilpotent groups, use induction on the 
nilpotence degree.
\label{no-lib-nilp}
\end{ejer} \end{small}

To find free semigroups inside groups acting on the 
line, the following notion is quite appropriate. 

\begin{defn} Two orientation preserving homeomorphisms of the real line are 
{\bf \em{crossed}} \esp on an interval $]a,b[$ if one of them fixes $a$ and $b$ but no 
other point in $[a,b]$, while the other one sends either $a$ or $b$ into $]a,b[$. Here 
we allow the cases \esp $a = - \infty$ \esp or \esp $b = +\infty$.
\label{def-crossed-elements}\index{crossed elements}
\end{defn}

In Foliation Theory, 
\index{foliation} 
the notion of crossed elements corresponds to that of {\em resilient leaves} 
({\em feuilles ressort}, \esp in the French terminology): see for instance \cite{CC1}. 
However, the latter is more general, since it also applies to \textit{pseudo-groups} 
of homeomorphisms of one-dimensional manifolds ({\em c.f.,} Definition \ref{def-pseudo}). 
In addition, crossed elements are dynamically relevant since they somewhat correspond 
to one-dimensional versions of {\em Smale horseshoes} \cite{PalT}.
 
The next elementary criterion showing that certain semigroups 
are free is well-known.

\vspace{0.15cm}

\begin{lem} Every subgroup of $\mathrm{Homeo}_+(\mathbb{R})$ 
having crossed elements contains a free semigroup on two generators.
\label{semilibre}
\end{lem}

\noindent{\bf Proof.} Suppose that for two elements $f,g$ in 
$\mathrm{Homeo}_+(\mathbb{R})$ there exists an interval $[a,b]$ 
such that $\mathrm{Fix}(f) \cap [a,b] = \{a,b\}$ and $g(a) \! \in ]a,b[$ 
(the case where $g(b) \! \in ]a,b[$ is similar). Changing $f$ by its inverse 
if necessary, we may assume that $f(x) < x$ for every $x \! \in ]a,b[$. 
Let $c = g(a) \! \in ]a,b[$, and let $d'$ be a point in $]c,b[$. Since 
$gf^k(a) = c$ for every $k\!\in\!\mathbb{N}$, and since $gf^k(d')$ converges to 
$c$ as $k$ goes to infinity, for $k \in \! \mathbb{N}$ large enough the map $gf^k$ 
has a fixed point on $]a,d'[$. Let $n \!\in\! \mathbb{N}$ be such an integer, and 
let $d > c$ be the infimum of the set of fixed points of $gf^n$ in $]a,b[$. Let  
$m \!\in\! \mathbb{N}$ be large enough so that $f^m(d) < c$. Then the Positive 
Ping-Pong Lemma applied to the restrictions of $f^m$ and $gf^n$ to $[a,b]$ 
shows that the semigroup generated by these elements is free (see Exercise 
\ref{pingpong+}). $\hfill\square$

\index{Klein!Ping-Pong Lemma}

\vspace{0.3cm}

According to Exercise \ref{no-lib-nilp} and Lemma \ref{semilibre}, 
the following result corresponds to a generalization of Theorem 
\ref{plante} (compare \cite{bek,bekla,sol}).

\vspace{0.1cm}

\begin{prop} {\em Let $\Gamma$ be a finitely generated group of orientation 
preserving homeomorphisms of the real line. If \esp $\Gamma$ has no crossed 
elements, then $\Gamma$ preserves a Radon measure on $\mathbb{R}$.}
\label{radon-inv}
\end{prop}

\noindent{\bf Proof.} If $\Gamma$ has a global fixed point, then the claim is 
obvious: the Dirac measure with mass on any of these points is invariant by the 
action. Assume throughout that there is no global fixed point. According to Proposition 
\ref{para-la-recta}, $\Gamma$ preserves a non-empty minimal invariant closed 
set $\Lambda$. Moreover, by the comments after the proof of that proposition, 
there are three possibilities.

\vspace{0.15cm}

\noindent \underbar{Case (i).} $\Lambda' = \emptyset$.
 
\vspace{0.1cm}

In this case, $\Lambda$ coincides with the set of points of an increasing sequence 
$(y_n)_{n \in \mathbb{Z}}$ without accumulation points in $\mathbb{R}$. One then 
easily checks that the Radon measure \esp $\sum_{n \in \mathbb{Z}} \delta_{y_n}$ 
\esp is invariant by $\Gamma$.

\vspace{0.15cm}

\noindent \underbar{Case (ii).} $\partial \Lambda = \emptyset$.

\vspace{0.1cm}
 
In this case, the action of $\Gamma$ is minimal. We claim that this action is also 
free. Indeed, if not then there exist an interval in $\mathbb{R}$ of the form $[u,v[$ 
or $]u,v]$, and an element $g \!\in\! \Gamma$ fixing $]u,v[$ and with no fixed point 
inside. Since the action is minimal, there must be some $h \in \Gamma$ 
sending $u$ or $v$ inside $]u,v[$; 
however, this implies that $g$ and $h$ are crossed on $[u,v]$, 
thus contradicting our assumption. Now the action of $\Gamma$ being 
free, H\"older's theorem implies that $\Gamma$ 
\index{H\"older!theorem} \index{action!free} 
is topologically conjugate to a (in this case dense) group of translations. Pulling 
back the Lebesgue measure by this conjugacy, we obtain an invariant Radon measure 
for the action of $\Gamma$.

\vspace{0.15cm}

\noindent \underbar{Case (iii).} $\partial \Lambda = \Lambda' = \Lambda$.
 
\vspace{0.1cm}

Collapsing to a point the closure of each connected component of the complement 
of the ``local Cantor set'' $\Lambda$, we obtain a topological line with a 
$\Gamma$-action induced by semi-conjugacy. As in the second case, 
one easily checks that the induced action is free, hence it preserves 
a Radon measure. Pulling back this measure by the semi-conjugacy, one obtains a 
Radon measure on $\mathbb{R}$ which is invariant by the original action. $\hfill\square$

\vspace{0.2cm}

\begin{small}\begin{ejer} Let $\Gamma$ be a non-necessarily finitely generated group of 
homeomorphisms of the line without crossed elements.

\vspace{0.07cm}

\noindent (i) Prove directly (\textit{i.e.,} without using Proposition \ref{radon-inv}) 
that the set $\Gamma_0$ formed by the elements of $\Gamma$ having fixed points is 
a normal subgroup of $\Gamma$.

\vspace{0.07cm}

\noindent (ii) Prove that the group $\Gamma / \Gamma_0$ admits 
an Archimedean ordering.

\vspace{0.07cm}

\noindent (iii) Using (i) and (ii), give an alternative proof for Proposition \ref{radon-inv}.
\end{ejer} 

\begin{ejer} Let $\mathcal{G} = \{f_1,\ldots,f_k\}$ be a system of generators of a group 
$\Gamma$ of homeomorphisms of the line having no global fixed point and without crossed 
elements. Show that at least one of these generators does not have fixed points. 

\noindent{\underbar{Hint.}} Suppose for a contradiction that all the $f_i$'s have fixed points, 
and let $x_1 \!\in \mathbb{R}$ be any fixed point of $f_1$. If $f_2$ fixes $x_1$ then 
the point $x_2=x_1$ is fixed by both $f_1$ and $f_2$. If not, choose a fixed point 
$x_2 \!\in \mathbb{R}$ of $f_2$ such that $f_2$ does not fix any point between 
$x_1$ and $x_2$. Show that $x_2$ is still fixed by $f_1$. Continuing in this way, 
show that there is a point which is simultaneously fixed by all the $f_i$'s, 
thus contradicting the hypothesis.
\label{radon-inv-ejer}
\end{ejer}
\end{small}

The reader will easily check that the finite generation hypothesis is necessary for 
Proposition \ref{radon-inv} (see Exercise \ref{ejemplo-plante}). 
However, throughout the proof this 
hypothesis was only used for ensuring the existence of a non-empty minimal invariant 
closed set, which in its turn easily follows from the existence of an element without 
fixed points. Therefore, if such a condition is assumed {\em a priori}, then the theorem 
still holds (compare Exercise \ref{din-gen}). We record this fact in the 
case of Abelian groups as a proposition.

\vspace{0.1cm}

\begin{prop} {\em Every Abelian subgroup of $\mathrm{Homeo}_{+}(\mathbb{R})$ 
having elements without fixed points preserves a Radon measure on the line.}
\label{abelianoplante}
\end{prop}

\vspace{0.01cm}

\begin{small} \begin{ejer} In an alternative way, prove the proposition above 
by considering the action induced on the topological circle obtained as 
the quotient of the line by the action of an element without fixed points.
\end{ejer}

\begin{ejer} Give an example of a countable, Abelian, infinitely generated subgroup of 
$\mathrm{Homeo}_+(\mathbb{R})$ for which there is no invariant Radon measure 
on the line (see \cite{Pl2} in case of problems with this).
\label{ejemplo-plante}
\end{ejer}

\begin{obs} If a finitely generated group acts on the real line without global fixed 
points and preserving a Radon measure, then the corresponding translation number function 
provides us with a nontrivial homomorphism into $(\mathbb{R},+)$. However, there is 
a large variety of finitely generated orderable groups for which there is no such 
a homomorphism (compare \S \ref{super-witte-morris}). A concrete example is the 
preimage $\tilde{\mathrm{G}}$ in $\widetilde{\mathrm{Homeo}}_+(\mathbb{R})$ of 
Thompson's group G. Indeed, although $\tilde{\mathrm{G}}$ is not simple, it is 
a {\bf{\em perfect}} group, that is, 
\index{group!perfect}
it coincides with its first derived group \cite{CFP}. Another 
(historically important) example will be discussed in \S \ref{sec-th-th}. For 
these groups, all nontrivial actions on the line must have crossed elements.
\index{Thompson's groups!G}
\label{ojuela}
\end{obs}
\end{small}

\vspace{0.1cm}

Theorem \ref{plante} concerns nilpotent groups, and leads naturally 
to the case of solvable groups. \index{group!solvable}
To begin with, notice that if $t \neq 0$ and $\kappa \neq 1$, then 
the subgroup $\mathrm{Aff}_{+}(\mathbb{R})$ generated by $f(x) = x + t$ 
and $g(x) = \kappa x$ does not preserve any Radon measure on the line. However, since the 
elements in the affine group \index{affine group}
preserve the Lebesgue measure up to a factor, having in mind 
Proposition \ref{conjugaralafin} it is natural to ask for conditions ensuring that a 
solvable group of $\mathrm{Homeo}_{+}(\mathbb{R})$ leaves quasi-invariant a Radon measure 
on the line. Once again, we will only consider finitely generated groups, since there 
exist non finitely generated Abelian groups of homeomorphisms of the line which do not 
admit any nontrivial quasi-invariant Radon measure (compare Exercise \ref{ejemplo-plante}).

Solvable groups are constructed starting from Abelian groups by successive extensions.
We then begin with an elementary remark: if $\Gamma_0$ is a normal subgroup of a subgroup 
$\Gamma$ of $\mathrm{Homeo}_{+}(\mathbb{R})$ and preserves a Radon measure $\upsilon$, 
then for each $g \!\in\! \Gamma$ the measure $g(\upsilon)$ is also invariant by 
$\Gamma_0$. Indeed, for every $h \in \Gamma_0$ we have 
$$(ghg^{-1}) \big( g (\upsilon) \big) = g \big( h (\upsilon) \big) = g (\upsilon),$$
where $g(\upsilon)(A) = g_{*}(\upsilon)(A) = \upsilon \big( g^{-1}(A) \big)$ 
for every Borel set $A \subset \mathbb{R}$. 
By Lemma \ref{igualtau}, there exists a constant $\kappa(g)$ 
such that $\tau_{g(\upsilon)} = \kappa(g) \tau_{\upsilon}$. 


\begin{small}\begin{ejer}
Show that $\kappa (g)$ does not depend on the $\Gamma_0$-invariant Radon measure $v$.
\end{ejer}\end{small}


The next three lemmas deal with the problem of 
showing the existence of a quasi-invariant measure for a 
group starting from the quasi-invariance of some measure for a normal subgroup.

\vspace{0.1cm}

\begin{lem} {\em Let $\Gamma_0$ be a normal subgroup of a subgroup $\Gamma$ of 
\esp $\mathrm{Homeo}_{+}(\mathbb{R})$, and let $\upsilon$ be a Radon measure 
invariant by $\Gamma_0$. Suppose that $\tau_{\upsilon}(\Gamma_0) \neq \{0\}$. 
If $\kappa(\Gamma) = \{ 1 \}$ and $\Gamma/\Gamma_0$ is amenable, 
\index{group!amenable} then there 
exists a Radon measure which is invariant by $\Gamma$.}
\label{lema-uno}
\end{lem}

\noindent{\bf Proof.} If $\tau_{\upsilon}(\Gamma_0)$ is dense in $\mathbb{R}$, then 
Proposition \ref{igualmedida} shows that $\upsilon$ is invariant by $\Gamma$. Assume 
throughout that $\tau_{\upsilon}(\Gamma_0)$ is cyclic, and let us normalize 
$\upsilon$ so that $\tau_{\upsilon}(\Gamma_0)$ equals $\mathbb{Z}$. Let $h_0\in\Gamma_0$ 
be such that $\tau_{\upsilon}(h_0) = 1$, and let $\Gamma^*_{0}$ be the kernel of 
$\tau_{\upsilon}$. The set $\mathrm{Fix}(\Gamma_0)$ is non-empty (it contains the 
support of $\upsilon$) and is invariant by $\Gamma_0$; hence, it is unbounded in both 
directions. Moreover, $\Gamma^*_0$ coincides with 
$\{ h\in \Gamma_0\!: \mathrm{Fix}(h)\neq \emptyset\}$, 
which shows that $\Gamma_0^*$ is normal in $\Gamma$ (and not only in $\Gamma_0$).\\

We will now use a similar idea to that of the proof of Proposition \ref{abelianoplante}. 
We first claim that $\Gamma_0 / \Gamma_0^*$ is contained in the 
\index{center of a group} 
center of 
$\Gamma / \Gamma_0^*$. To prove this, it suffices to show that for every $g \!\in\! \Gamma$ 
and every $x \!\in\! \mathrm{Fix}(\Gamma_0^*)$ one has $g^{-1}h_0g(x) = h_0(x)$. Consider a 
measure $\upsilon_1$ giving mass $1$ to each point in the set $\{ h_0^n(x), \esp n\in\mathbb{Z}\}$.   
Notice that $\upsilon_1$ is invariant under $\Gamma_0$. From the hypothesis $\kappa(\Gamma) \!=\! 1$ 
we conclude that 
$$\upsilon_1 \big( [x,g^{-1}h_0g(x)[ \big) = \tau_{\upsilon_1}(g^{-1}h_0g) 
= \tau_{g (\upsilon_1)}(h_0) = \kappa(g) \esp 
\upsilon_1 \big( [x,h_0(x)[ \big)= \upsilon_1 \big( [x,h_0(x)[ \big),$$
which implies that $g^{-1}h_0g(x) \leq h_0(x)$. On the other hand, by considering 
the measure $\upsilon_2$ which gives mass $1$ to each point in the set $g^{-1}h_0^ng(x)$, and 
changing $h_0$ by $g^{-1}h_0g$, the same argument shows that $h_0(x) \leq g^{-1}h_0g(x)$.\\

Taking $x_0 \!\in\! \mathrm{Fix}(\Gamma_0^*)$, we notice that the group 
$\big( \Gamma/\Gamma_0^* \big)/\big( \Gamma_0/\Gamma_0^* \big) =
\Gamma/\Gamma_0$ acts on the quotient space \esp 
$\mathrm{Fix}(\Gamma_0^*)/ \langle h_0 \rangle (x_0)$, \esp 
which is compact. By the amenability hypothesis, $\Gamma/\Gamma_0$ preserves 
a probability measure on this space, which lifts to a $\Gamma$-invariant 
Radon measure on $\mathrm{Fix}(\Gamma_0^*)$. $\hfill\square$

\vspace{0.2cm}

\begin{lem} {\em Let $\Gamma_0$ be a normal subgroup of a subgroup $\Gamma$ of 
$\mathrm{Homeo}_{+}(\mathbb{R})$, and let $\upsilon$ be a Radon measure invariant 
by $\Gamma_0$. Suppose that $\tau_{\upsilon}(\Gamma_0)\neq \{ 0 \}$. If 
$\kappa(\Gamma) \neq \{1\}$, then $\upsilon$ is quasi-invariant by $\Gamma$.}
\label{lema-dos}
\end{lem}

\noindent{\bf Proof.} Let $h \in \Gamma_0$ and $g \in \Gamma$ be such that 
$\tau_{\upsilon}(h) \neq 0$ and $\kappa(g) \neq 1$. Since the equality 
$$\tau_{\upsilon}(g^{-n}h^mg^n) = \tau_{g^n (\upsilon)}(h^m) 
= m \kappa(g)^n \tau_{\upsilon}(h)$$
holds for every $m,n$ in $\mathbb{Z}$, the image $\tau_{\upsilon}(\Gamma)$ 
must be dense in $\mathbb{R}$. Proposition \ref{igualmedida} then shows 
that $\upsilon$ is quasi-invariant by $\Gamma$. $\hfill\square$

\vspace{0.2cm}

\begin{lem} {\em Let $\Gamma_0$ be a normal subgroup of a subgroup $\Gamma$ of 
$\mathrm{Homeo}_{+}(\mathbb{R})$, and let $\upsilon$ be a Radon measure 
quasi-invariant by $\Gamma_0$. Let $\Gamma_0^{*}$ be the subgroup of 
the elements $h \in \Gamma_0$ such that $h (\upsilon) = \upsilon$.
If $\tau_{\upsilon}(\Gamma_0^{*}) \neq \{ 0 \}$ and $\kappa(\Gamma_0)
\neq \{ 1 \}$, then $\upsilon$ is quasi-invariant by $\Gamma$.}
\label{lema-tres}
\end{lem}

\noindent{\bf Proof.} Since $\tau_{\upsilon}(\Gamma_0^{*}) \neq \{0\}$ and 
$\kappa(\Gamma_0) \neq \{ 1 \}$, it is easy to see that an element $h \in \Gamma_0$ 
belongs to $\Gamma_0^{*}$ if and only if either $\mathrm{Fix}(h) = \emptyset$ or 
$\mathrm{Fix}(h)$ is unbounded in both directions. Both conditions being stable 
under conjugacy, one concludes that $\Gamma_0^*$ is normal in $\Gamma$. It is also 
clear that $\kappa(\Gamma) \neq \{ 1 \}$. The claim then follows from the 
preceding lemma. $\hfill\square$

\vspace{0.3cm}

We are now ready to deal with the problem of the existence of a quasi-invariant Radon 
measure for a large family of solvable subgroups of $\mathrm{Homeo}_{+}(\mathbb{R})$.

\begin{thm} {\em Let $\Gamma$ be a solvable subgroup of $\mathrm{Homeo}_{+}(\mathbb{R})$. 
Suppose that $\Gamma$ admits a chain of subgroups $\{ id \} \!=\! \Gamma_n \sgn
\Gamma_{n-1} \sgn \ldots \sgn \Gamma_0 \!=\! \Gamma$ such that each $\Gamma_i$ 
is finitely generated and each quotient $\Gamma_{i-1} / \Gamma_i$ is Abelian. 
Then there exists a Radon measure which is quasi-invariant by $\Gamma$.}
\label{plantesoluble}
\end{thm}

\noindent{\bf Proof.} We will assume that $\Gamma$ does not preserve 
any Radon measure. We then have $\mathrm{Fix}(\Gamma) = \emptyset$, 
since otherwise the Dirac delta with mass on a global fixed 
point would be an invariant measure. Let $j > 0$ be the smallest 
index for which $\mathrm{Fix}(\Gamma_j) \neq \emptyset$. The Abelian and 
finitely generated group $\Gamma_{j-1}/\Gamma_j$ acts on the closed 
unbounded set $\mathrm{Fix}(\Gamma_j)$ by homeomorphisms which preserve 
the order. We leave to the reader the task of showing the existence of 
an invariant measure for this action by using (a slight modification 
of) Theorem \ref{plante}. This measure naturally induces a Radon 
measure $\upsilon$ on the line which is invariant by $\Gamma_{j-1}$ 
and whose support is contained in $\mathrm{Fix}(\Gamma_j)$. 
Notice that $\tau_{\upsilon}(\Gamma_{j-1}) \neq \{ 0 \}$,
since $\mathrm{Fix}(\Gamma_{j-1}) = \emptyset$.\\

Let $k > 0$ be the smallest index for which $\Gamma_k$ preserves a Radon measure. 
Slightly abusing of the notation, let us still denote this measure by $\upsilon$. 
As above, it is easy to see that $\tau_{\upsilon}(\Gamma_{k}) \neq \{ 0 \}$.
According to Lemma \ref{lema-uno}, we have $\kappa(\Gamma_{k-1})\neq \{ 1\}$, 
and by Lemma \ref{lema-dos} this implies that $\upsilon$ is quasi-invariant 
by $\Gamma_{k-1}$. Lemma \ref{lema-tres} then allows to prove by induction 
that $\upsilon$ is quasi-invariant by $\Gamma_{k-2}$, $\Gamma_{k-3}$, etc. Thus, 
at the end we conclude that $\upsilon$ is quasi-invariant by $\Gamma$. $\hfill\square$

\vspace{0.35cm}

The hypothesis of finite generation for each $\Gamma_i$ is satisfied by an important 
family of groups, namely that of {\bf{\em polycyclic groups}} (see Appendix A). 
In addition ,it may 
\index{group!polycyclic}
be weakened, as is shown below.

\begin{small} \begin{ejer} The dynamics of a subgroup $\Gamma$ of $\mathrm{Homeo}_{+}(\mathbb{R})$ 
is said to be {\bf \textit{boundedly generated}} if there exists a system of generators $\mathcal{G}$ 
of $\Gamma$ and a point $x_0 \in \Gamma$ such that $\{ h(x_0)\!: h \in \mathcal{G} \}$ is a bounded 
subset of the line. Prove that Theorems \ref{plante} and \ref{plantesoluble}, as well as 
Proposition \ref{radon-inv}, are true if the hypothesis of finite generation is replaced 
by a hypothesis of boundedly generated dynamics.
\label{din-gen}
\end{ejer} \end{small}
\index{translation number|)}
\index{measure!quasi-invariant|)}


\subsection{An application to amenable, orderable 
groups\index{group!orderable}\index{group!amenable}}
\label{super-witte-morris}

\hspace{0.45cm} In this section, we will use some previously developed ideas for giving 
a ``dynamical proof'' of the following algebraic result due to Witte-Morris  
\cite{Witam} (we recommend the lecture of Appendix B for the concept of 
amenability; see also Remark \ref{ojuela}, the example after Exercise 
\ref{th-analitico}, and Exercise \ref{nosecomo}).
\index{Witte-Morris!theorem on orderable amenable groups}

\vspace{0.04cm}

\begin{thm} {\em Every orderable, finitely generated, infinite, amenable 
group admits a nontrivial homomorphism into $(\mathbb{R},+)$.}
\label{witte-bien}
\end{thm}

\vspace{0.04cm}

This theorem settles an old problem in the theory of left-orderable groups.\footnote{A 
conjecture leading to the theorem above appears in the work by Linnell \cite{linelito}, 
\index{Linnell}
but quite remarkably it is already present in a different form in an old  seminal paper by 
\index{Thurston}
Thurston to be discussed in \S \ref{sec-th-th} (see \cite[Page 348]{Th2}).} 
We will follow essentially the same brilliant idea of Witte-Morris, but 
unlike \cite{Witam} we will avoid the use of the algebraic theory of 
$\mathcal{C}$-orderable \index{group!$\mathcal{C}$-orderable} groups. 

The first ingredient of the proof goes back to an idea 
independently due to Ghys 
\index{Ghys!space of orderings} and Sikora 
\index{Sikora}
(see for instance \cite{sikora}), which consists in the introduction 
of a {\bf \textit{space of orderings}} associated to an orderable 
group, and of a natural action of the group on it. More precisely, 
for a finitely generated orderable group $\Gamma$, let us denote 
by $\mathcal{O}(\Gamma)$ the set of all orderings on $\Gamma$. 
If we fix a finite system $\mathcal{G}$ of generators for $\Gamma$, then we 
may define the distance between \esp $\leq$ \esp and \esp $\preceq$ \esp 
in $\mathcal{O}(\Gamma)$ by putting \esp 
$dist(\leq,\preceq) = e^{-n}$, \esp where $n$ is the maximum non-negative 
integer such that the orderings $\leq$ and $\preceq$ coincide on the ball 
$B_{\mathcal{G}}(n)$ of radius $n$ in $\Gamma$ (see Appendix B for the notion 
of \textit{ball} inside a group). In other words, $n$ is the biggest non-negative 
integer such that for all $g,h$ in $B_{\mathcal{G}}(n)$ one has $g \leq h$ if and 
only if $g \preceq h$. If we also let \esp $dist(\preceq,\preceq) = 0$ \esp 
for every ordering $\preceq$, then it is easy to check that the thus defined 
function $dist$ is a distance on $\mathcal{O}(\Gamma)$ (which depends on 
$\mathcal{G}$). Actually, the resulting metric space is ultrametric and compact.  

\begin{small} \begin{obs} The structure of the space of orderings of an 
orderable group is interesting by itself. A quite elegant result by Tararin \cite{koko} 
completely describes all orderable groups admitting only finitely many orderings. 
If a group admits infinitely many orderings, then it necessarily admits uncountably many 
\cite{linnel,cnr,ordering2}. For higher-rank torsion-free Abelian groups \cite{sikora}, 
for non-Abelian torsion-free nilpotent 
\index{group!nilpotent} \index{group!torsion-free}groups \cite{ordering2}, 
and for non-Abelian free 
groups \cite{MC,ordering2}, the spaces of orderings are 
known to be homeomorphic to the Cantor set. However, there exist 
relevant examples of groups whose spaces of orderings are infinite 
but do contain isolated points \cite{adam-clay,dub,ordering2}. 
\index{Linnell} 
According to Linnell \cite{linnel}, one of the reasons for the interest 
in all of this concerns the structure of the semigroup formed by 
the positive elements, as stated in the following exercise. 
\end{obs} 

\begin{ejer} Prove that if $\preceq$ is a non-isolated point in the space of orderings of a 
finitely generated group, then its positive cone is not finitely generated as a semigroup.
\index{semigroup}
\end{ejer} \end{small}

The group $\Gamma$ acts (continuously) on $\mathcal{O}(\Gamma)$ by right multiplication: 
given an ordering $\preceq$ with positive cone $\Gamma_+$ and an element $f \!\in\! \Gamma$, 
the image of $\preceq$ under $f$ is the ordering $\preceq_{f}$ whose positive cone is the 
conjugate \esp $f \esp \Gamma_+ f^{-1}$ \esp of $\Gamma_+$ by $f$. In other words, one 
has $g \preceq_f h$ if and only if $fgf^{-1} \preceq fhf^{-1}$, which is equivalent 
to $gf^{-1} \preceq hf^{-1}$. 

We will say that an ordering $\preceq$ is {\bf \textit{right-recurrent}} if for every pair 
of elements $f,h$ in $\Gamma$ such that $f \!\succ\! id$, there exists $n \!\in\! \mathbb{N}$ 
satisfying $fh^n \succ h^n$. For instance, every bi-invariant ordering is right-recurrent. 
The following corresponds to the main step in the proof of Witte-Morris' theorem.

\vspace{0.1cm}

\begin{prop} {\em If \esp $\Gamma$ is a finitely generated, amenable, 
orderable group, then $\Gamma$ admits a right-recurrent ordering.}
\label{witte-fund}
\end{prop}

\vspace{0.1cm}

To prove this proposition we need the following weak form of the Poincar\'e 
Recurrence Theorem. 
\index{Poincar\'e!Recurrence Theorem}

\vspace{0.1cm}

\begin{thm} {\em If $T$ is a measurable map that preserves a probability 
measure $\mu$ on a space $\mathrm{M}$, then for every measurable subset 
$A$ of \esp $\mathrm{M}$ and $\mu$-almost every point $x \! \in \! A$ there 
exists $n \! \in \! \mathbb{N}$ such that $T^n(x)$ belongs to $A$.}
\label{rec-poincare}
\end{thm}

\noindent{\bf Proof.} The set of points in $A$ which do not come back to $A$ under 
iterates of $T$ is $B = A \setminus \cup_{n \in \mathbb{N}} T^{-n}(A)$. One easily 
checks that the sets $T^{-i}(B)$, with $i \!>\! 0$, are two-by-two disjoint. Since 
$T$ preserves $\mu$, these sets have the same measure, and since the total mass of 
$\mu$ equals $1$, the only possibility is that this measure equals zero. Therefore, 
$\mu(B)=0$, that is, $\mu$-almost every point in $A$ comes back to $A$ under some 
iterate of $T$. $\hfill\square$

\vspace{0.15cm}

\begin{small} \begin{ejer} 
Prove that, under the hypothesis of the preceding theorem, $\mu$-almost every 
point in $A$ comes back to $A$ under {\em infinitely many} iterates of $T$.
\label{poincare-general}
\end{ejer} \end{small}

\vspace{0.05cm}

\noindent{\bf Proof of Proposition \ref{witte-fund}.} By definition, if a finitely 
generated orderable group $\Gamma$ is amenable, then its action on $\mathcal{O}(\Gamma)$ 
preserves a probability measure $\mu$. We will show that $\mu$-almost every point in 
$\mathcal{O}(\Gamma)$ is right-recurrent. To do this, for each $g \!\in\! \Gamma$ let us consider  
the subset $A_g$ of $\mathcal{O}(\Gamma)$ formed by all of the orderings $\preceq$ on $\Gamma$ 
such that $g \!\succ\! id$. By the Poincar\'e Recurrence Theorem, for each $f \in \Gamma$ the set 
\esp $B_g(f) = A_g \setminus \cup_{n \in \mathbb{N}} f^{n}(A_g)$ \esp has null $\mu$-measure. 
Therefore, the measure of $B_g = \cup_{f \in \Gamma} B_g (f)$ is also zero, as well as the 
measure of $B = \cup_{g \in \Gamma} B_g$. Let us consider an arbitrary element $\preceq$ 
in the ($\mu$-full measure) set $A = \mathcal{O}(\Gamma) \setminus B$. Given $g \succ id$ and 
$f \in \Gamma$, from the inclusion $B_g(f) \subset B$ we deduce that $\preceq$ does not belong 
to $B_g(f)$, and thus there exists $n \!\in\! \mathbb{N}$ such that $\preceq$ belongs to $f^{n} (A_g)$. 
In other words, one has $g \succ_{f^{-n}} id$, that is, $gf^n \succ f^n$. Since $g \succ id$ and 
$f \in \Gamma$ were arbitrary, this shows the right-recurrence of $\preceq$. $\hfill\square$

\vspace{0.54cm}


\beginpicture

\setcoordinatesystem units <0.7cm,0.7cm>



\putrule from 0 0 to 0 8
\putrule from 0 0 to 8 0
\putrule from 8 0 to 8 8 
\putrule from 0 8 to 8 8

\put{$a$} at 0 -0.4
\put{$t(id)$} at 1.2 -0.4
\put{$c$} at 2 -0.4
\put{$c_m$} at 2.7 -0.4
\put{$c_m'$} at 3.3 -0.4
\put{$c_n$} at 6 -0.4
\put{$c_n'$} at 6.6 -0.4
\put{$d$} at 7.2 -0.4
\put{$b$} at 8 -0.4
\put{$c$} at -0.4 2 
\put{$f$} at 5.05 1.5 
\put{$g$} at 1.8 7 
\put{$h \!=\! h_n^k$} at 2.25 5.4 
\put{$h_n$} at 2.8 4.32   
\put{$h_m$} at 2.6 3.2 
\put{$h_m(c_n)$} at -0.87 3.65
\put{$h(a)$} at -0.55 4.32 

\setquadratic 

\plot 
0 2 1.6 2.55 2.7 2.7 /

\plot 
2.7 2.7 2.92 2.82 3 3 /

\plot 
3 3 3.1 3.2 3.3 3.3 /

\plot 
0 0 5.6 2.4 8 8 /

\plot 
0 2 1 6 2.3 8 /

\plot 
6 6 6.25 6.1 6.3 6.3 /

\plot 
6.3 6.3 6.36 6.5 6.6 6.6 /

\plot 
6 6 6.15 6.05 6.3 6.3 /

\plot 
6.3 6.3 6.42 6.52 6.6 6.6 /


\setlinear

\plot
0 2 
6 6 /

\plot 
0 4.2  
6 6 /

\plot 
6 6 
6.6 6.6 /

\plot 
3.3 3.3 
6 3.7 /

\plot 
6 3.7 
7.2 3.9 /

\plot 
6.6 6.6 
7.2 6.8 / 

\plot 
6.6 6.6 
7.2 6.7 /


\plot 3.3 3.3  3.5 3.3 /

\plot 3.7 3.3  3.9 3.3 /

\plot 4.1 3.3 4.3 3.3 /

\plot 4.5 3.3  4.7 3.3 /

\plot 4.9 3.3  5.1 3.3 /

\plot 5.3 3.3 5.5 3.3 /

\plot 5.7 3.3 5.9 3.3 /


\plot 
3.3 3.3  3.3 3.5 /

\plot 
3.3 3.7  
3.3 3.9 /

\plot
3.3 4.1 
3.3 4.3 /

\plot 
3.3 4.5
3.3 4.7 /

\plot 
3.3 4.9 
3.3 5.1 /

\plot
3.3 5.3 
3.3 5.5 /

\plot 
3.3 5.7 
3.3 5.9 /


\plot 
3.3 6  
3.5 6 /

\plot 
3.7 6  
3.9 6 /

\plot
4.1 6 
4.3 6 /

\plot 
4.5 6  
4.7 6 /

\plot 
4.9 6  
5.1 6 /

\plot
5.3 6 
5.5 6 /

\plot 
5.7 6 
5.9 6 /


\plot 
6 3.3 
6 3.5 /

\plot 
6 3.7   
6 3.9 /

\plot
6 4.1  
6 4.3 /

\plot 
6 4.5   
6 4.7 /

\plot 
6 4.9   
6 5.1 /

\plot
6 5.3  
6 5.5 /

\plot 
6 5.7 
6 5.9 /


\setdots

\plot 
0 0 
8 8 /

\putrule from 1.2 0 to 1.2 3.7 
\putrule from 2 0 to 2 2 
\putrule from 0 2 to 2 2 
\putrule from 2.7 0 to 2.7 2.7 
\putrule from 3.3 0 to 3.3 3.3 
\putrule from 6 0 to 6 3.3 
\putrule from 6.6 0 to 6.6 6.6 
\putrule from 0 6 to 3.3 6 
\putrule from 0 3.7 to 6 3.7 
\putrule from 7.2 0 to 7.2 7.2 


\put{Figure 13} at 4 -1.2 

\put{} at -5.4 0 

\endpicture


\vspace{0.54cm}

The right-recurrence for an ordering has dynamical 
consequences, as is shown below.

\vspace{0.15cm}

\begin{prop} {\em Let $\Gamma$ be a countable group admitting a right-recurrent ordering $\preceq$. If 
$(g_n)_{n \geq 0}$ is any numbering of $\Gamma$, then the dynamical realization of $\Gamma$ associated 
to $\preceq$ and this numbering is a subgroup of $\mathrm{Homeo}_+(\mathbb{R})$ without crossed elements.}
\label{conrad-yo}
\end{prop}

\noindent{\bf Proof.} The claim is obvious if $\Gamma$ is trivial. 
Therefore, in what follows, we will assume that $\Gamma$ is infinite. 
\index{crossed elements}
Let us suppose that 
there exist $f,g$ in $\Gamma$ and an interval $[a,b]$ such that (for the dynamical realizations one has) 
$\mathrm{Fix}(f) \cap [a,b] = \{a,b\}$ and $g(a) \! \in ]a,b[$ (the case where $g(b)$ belongs to $]a,b[$ 
is analogous). Changing $f$ by its inverse if necessary, we may assume that $f(x) \!<\! x$ for all 
$x \! \in ]a,b[$. As we have already observed after Theorem \ref{thorden}, there must exist some 
$g_i \!\in\! \Gamma$ such that $t(g_i)$ belongs to the interval $]a,b[$. Changing $g_i$ by $f^n g_i$ 
for a large $n > 0$ if necessary, we may assume that $t(g_i)$ is actually contained in $]a,c[$, where 
$c \!=\! g(a)$. Now by conjugating $f$ and $g$ by the element $g_i$ (and changing the points 
$a,b,c$ by their images under $g_i$), we may suppose that $t(id)$ belongs to $]a,c[$. 

Choose a point $d\! \in ]c,b[$. Since \esp $gf^n(a) \!=\! c$ \esp for all $n\!\in\!\mathbb{N}$, 
and since $gf^n(d)$ converges to \esp $c \!<\! d$ \esp when $n$ goes to infinity, if \esp 
$n\!\in\!\mathbb{N}$ \esp is big enough then the map $h_n=gf^n$ satisfies $h_n(a)\!>\!a$, 
\esp $h_n(d) \!<\! d$, \esp $\mathrm{Fix}(h_n) \cap ]a,d[ \subset [c_n,c_n'] \subset ]c,h_n(d)[$, 
\esp and $\{c_n,c_n'\} \subset \mathrm{Fix}(h_n)$, for some sequences of points $c_n$ and $c_n'$ 
converging to $c$ by the right (see Figure 13). Notice that the element $h_n$ is positive, since 
from $h_n \big( t(id) \big) > h_n (a) = c > t(id)$ one deduces that $t(h_n) > t(id)$, and by 
the construction of the dynamical realization this implies that $h_n \succ id$. Let us fix 
$m > n$ large enough so that the preceding properties hold for $h_m$ and $h_n$, and such that 
\esp $[c_m,c_m'] \!\subset ]c,c_n[$. \esp Let us fix $k \in \mathbb{N}$ large enough in order that 
\esp $h_n^k (a) > h_m(c_n)$, \esp and let us define $h \!=\! h_m^k$. For each $i \!\in\! \mathbb{N}$ 
one has $h^i \big( t(id) \big) \! \in \esp ]h_m(c_n),c_n[$, and thus $h_m h^i \big( t(id) \big) 
\!<\! h_m(c_n) \!<\! h(a) \!<\! h \big( t(id) \big)$. 
Therefore, $h_m h^i \prec h \prec h^i$ for all 
$i \!\in\! \mathbb{N}$. However, this is in contradiction 
with the hypothesis of right-recurrence for $\preceq$. $\hfill\square$ 

\vspace{0.5cm}

In the preceding proof, we used a property which is actually weaker than right-recurrence, 
namely for all positive elements $f,h$ in $\Gamma$ one has $fh^n \!\succ\! h$ 
for some $n \!\in\! \mathbb{N}$. This is called the {\bf \textit{Conrad property}}, 
and the groups which do admit an ordering satisfying it are said to be 
{\bf\em{$\mathcal{C}$-orderable}}. There exists a very 
rich literature on this property mostly from an algebraic viewpoint (see 
for instance \cite{conrad,koko} as well as Remark \ref{en-general}). The proof above 
together with Exercise \ref{equiv-corden} show that it has a natural dynamical 
counterpart: roughly, the Conrad property 
\index{Conrad!Conrad property|see{group, $\mathcal{C}$-orderable}} 
\index{group!$\mathcal{C}$-orderable}is the algebraic 
counterpart to the condition of non-existence of crossed elements 
for the corresponding action on the real line \cite{cnr,ordering2}. 


\vspace{0.2cm}

\begin{small} \begin{ejer}
Let $(x_n)$ be a dense sequence of points in the line, and let $\Gamma$ be a subgroup of 
$\mathrm{Homeo}_+(\mathbb{R})$. Prove that if $\Gamma$ has no crossed elements, then 
the order relation induced from $(x_n)$ as in the proof of Theorem \ref{thorden} 
satisfies the Conrad condition.
\label{equiv-corden} 
\end{ejer} 

\begin{ejer} Let $\preceq$ be an ordering on a group $\Gamma$.

\noindent (i) Show that if $\preceq$ has the Conrad property, then 
$fh^2 \succ h$ for all positive elements $f,h$ in $\Gamma$. 

\noindent{\underbar{Hint.}} Following \cite{ordering2}, suppose that the 
opposite inequality holds and show that for the positive elements $f$ and 
$g = fh$ in $\Gamma$ one has $fg^n \prec g$ for every $n \!\in\! \mathbb{N}$.

\noindent (ii) More generally, prove that for every positive integer $n$ 
one has $f h^{n} \succ h^{n-1}$.

\noindent (iii) Give examples of right-recurrent orderings $\preceq$ such 
that $fh^2 \!\prec\! h^2$ for some $f \succ id $ and $h \succ id$.
\end{ejer}

\begin{ejer} Give examples of $\mathcal{C}$-orderable groups which do not admit 
right-recurrent orderings (see \cite[Example 4.5]{Witam} in case of problems 
with this).
\end{ejer} 

\begin{ejer} Show that the Klein group 
$\Gamma = \langle f,g\!: fgf^{-1} \!=\! g^{-1} \rangle$ 
is not bi-orderable, though it admits right-recurrent orderings.
\index{Klein!group}
\end{ejer} \end{small}

\vspace{0.1cm}

\noindent{\bf Proof of Theorem \ref{witte-bien}.} If $\Gamma$ is amenable, 
finitely generated, and orderable, then Proposition \ref{witte-fund} 
provides us with a right-recurrent ordering on $\Gamma$. By Proposition 
\ref{conrad-yo}, the dynamical realization associated to this ordering and 
any numbering $(g_n)_{n \geq 0}$ of $\Gamma$ is a subgroup of 
$\mathrm{Homeo}_+(\mathbb{R})$ (isomorphic to $\Gamma$ and) 
without crossed elements. By Proposition \ref{radon-inv}, $\Gamma$ preserves 
\index{measure!Radon}
a Radon measure $\upsilon$ on the line. Let us now recall that, if $\Gamma$ is 
nontrivial, then its dynamical realizations have no global fixed point. 
Therefore, the translation number function with respect to $\upsilon$ is 
a nontrivial homomorphism into $(\mathbb{R},+)$ (see (\ref{puntofijoplante}) 
and (\ref{soporteplante})). $\hfill\square$ 
\index{translation number}

\vspace{0.12cm}

\begin{small} 
\begin{ejer} Give a faithful action of the Baumslag-Solitar group 
\index{group!Baumslag-Solitar} 
$\Gamma = \langle f,g: fgf^{-1} = g^2 \rangle$ on the line 
without crossed elements (see \cite{rivas} for more on this). 

\noindent{\underbar{Remark.}} Notice that $\Gamma$ embeds in the affine 
group by identifying $f$ with $x \mapsto 2x$, and $g$ with $x \mapsto x+1$.
\end{ejer}

\begin{ejer} 
Show that finitely generated subgroups of the group of piecewise affine homeomorphisms 
of the interval admit nontrivial homomorphisms into $(\mathbb{R},+)$ (see also Exercise 
\ref{no-F}). 

\noindent{\underbar{Remark.}} Recall that the group $\mathrm{PAff}_+([0,1])$ does 
not contain free subgroups on two generators ({\em c.f.,} Theorem \ref{arefad}). 
Since this is also the case of countable amenable groups, 
this turns natural the following question: Do finitely generated 
subgroups of $\mathrm{Homeo}_+ ([0,1])$ without free subgroups on two 
generators admit nontrivial homomorphisms into $(\mathbb{R},+)$? See 
\cite{linnell} for an interesting result pointing in the positive direction.
\end{ejer} \end{small}


\section{Invariant Measures and Free Groups}

\subsection{A weak version of the Tits alternative}
\label{sec-margulis}\index{Margulis!weak Tits alternative}
\index{group!free|(}

\hspace{0.45cm} A celebrated theorem by Tits establishes that every finitely generated 
subgroup of $\mathrm{GL}(n,\mathbb{R})$ either contains a free subgroup on two generators 
\index{Tits!alternative}
or is virtually solvable \index{group!solvable}
(see \cite{breuillard,br-gel2} for modern versions of this 
result). This dichotomy, known as the {\bf \textit{Tits alternative}}, does 
not hold for groups of circle homeomorphisms, as is shown by the following exercise.

\vspace{0.05cm}

\begin{small} \begin{ejer} Show that Thompson's  
\index{Thompson's groups!F} group $\mathrm{F}$ is not virtually 
solvable (however, recall that, by Theorem \ref{arefad}, 
F does not contain free subgroups on two generators).

\end{ejer} \end{small}

Nevertheless, a weak version for the alternative may be 
established. The following result was conjectured by Ghys 
\index{Ghys!weak Tits alternative} and proved by Margulis 
in \cite{Ma1}. We will develop the proof proposed some years later 
by Ghys himself in \cite{Gh1}. The idea of this proof will be pursued 
in a probabilistic setting in \S \ref{Tits-prob}.

\vspace{0.1cm}

\begin{thm} {\em If \esp $\Gamma$ is a subgroup of $\mathrm{Homeo}_+(\clo)$, 
then either there exists a probability measure on $\mathrm{S}^1$ which 
is invariant by $\Gamma$, or $\Gamma$ contains a free subgroup on two 
generators.}
\label{margulis}
\end{thm}

\vspace{-0.25cm}

Since virtually solvable groups are \index{group!amenable}
amenable ({\em c.f.}, Exercise \ref{clasemoy}), every action 
of such a group on the circle admits an invariant probability measure. Therefore, 
the theorem above may be considered as a weak form of the Tits alternative. However, the 
weaker alternative is not a dichotomy. Indeed, according to \S \ref{sec-witte}, the free 
group on two generators admits faithful actions on the interval, which induce actions on 
$\clo$ with a global fixed point, and hence with an invariant probability measure.

\vspace{0.05cm}

To prove Margulis' theorem we begin by noticing that, by Theorem \ref{invariante}, there 
are three distinct cases according to the type of minimal invariant closed set. In 
the case where $\Gamma$ has a finite orbit, there is obviously an invariant probability 
measure, namely the mean of the Dirac measures with mass on the points of the orbit. We 
claim that the case where there exists an exceptional minimal 
set\index{exceptional minimal set} $\Lambda$ reduces to 
the case where all the orbits are dense. Indeed, the associated action on the topological 
circle $\mathrm{S}^1_{\Lambda}$ is minimal. If there exists an invariant probability measure 
for this action, then it cannot have atoms, and its support must be total. Hence, it induces 
an invariant measure on $\mathrm{S}^1$ whose support is $\Lambda$. Moreover, if two elements 
in $\Gamma$ project into homeomorphisms of $\mathrm{S}^1_\Lambda$ generating a free group, 
then these elements generate a free subgroup of $\Gamma$.

\vspace{0.15cm}

By the discussion above, to prove Margulis' theorem we only need to consider 
the case where all the orbits are dense ({\em i.e.,} the minimal case).

\vspace{0.18cm}

\begin{defn} An action of a group $\Gamma$ by circle homeomorphisms is said to 
be {\bf \textit{equicontinuous}} if for every $\delta > 0$ there exists $\varepsilon > 0$ 
such that, if \esp $dist(x,y) < \varepsilon$, \esp then \esp $dist(g(x),g(y)) < \delta$ 
\esp for all $g \in \Gamma$. It is said to be {\bf \textit{expansive}} if for each 
$x \!\in\! \clo$ there exists an open interval $I$ containing $x$ and a sequence 
$(g_n)$ of elements in $\Gamma$ such that the length of the intervals $g_n(I)$ 
converges to zero.
\index{action!equicontinuous}\index{action!expansive}
\end{defn}

\vspace{0.01cm}

\begin{lem} {\em Every minimal action by circle 
homeomorphisms is either equicontinuous or expansive.}
\end{lem}

\noindent{\bf Proof.} If the action is not equicontinuous, then 
there exists $\delta> 0$ such that for all $n \in \mathbb{N}$
there exist $x_n,y_n$ in $\mathrm{S}^1$ and $g_n$ in the acting 
group $\Gamma$ such that $\big| ]x_n,y_n[ \big|$ goes to zero 
and $dist \big( g_n(x_n),g_n(y_n) \big) \geq \delta.$ Passing to 
a subsequence if necessary, we may assume that $g_n(x_n)$ (resp. 
$g_n(y_n)$) converges to some  $a \!\in\! \mathrm{S}^1$ (resp. 
$b \!\in\! \mathrm{S}^1$). Notice that 
$dist(a,b) \geq \delta$. By the compactness of $\mathrm{S}^1$, and 
since the action is minimal, there exist $h_1,...,h_k$ in $\Gamma$ such 
that $\mathrm{S}^1 \!=\! \cup_{i=1}^{k} h_i(]a,b[)$. Let $\varepsilon > 0$
be the Lebesgue number of this covering. If $dist(x,y) < \varepsilon$ 
then $x,y$ belong to $h_j(]a,b[)$ for some $j \!\in\! \{1,\ldots,k \}$. 
The points $g_n^{-1} h_j^{-1}(x)$ and 
$g_n^{-1} h_j^{-1}(y)$ then belong to $]x_n,y_n[$ 
for $n$ large enough, and hence the distance between them tends 
to zero. This shows that the $\varepsilon$-neighborhood of every 
point in the circle is ``contractable'' to a point by elements in 
$\Gamma$, and therefore the action is expansive. $\hfill\square$

\vspace{0.2cm}

\begin{small} \begin{ejer} 
Let $\Gamma$ be a group of circle homeomorphisms acting minimally. Show that the 
action of $\Gamma$ is expansive if and only if there exists $\varepsilon \!>\! 0$ 
such that, for all $x \neq y$ in $\clo$, there exists $g \in \Gamma$ satisfying 
\esp $dist (g(x),g(y)) \geq \varepsilon$. 
\end{ejer} \end{small} 

If the action of a group $\Gamma$ by circle homeomorphisms is 
equicontinuous then, by Ascoli-Arzela's theorem, the closure of 
$\Gamma$ in $\mathrm{Homeo}_{+}(\mathrm{S}^1)$ is compact. Therefore, 
by Proposition \ref{gatorot}, there exists a probability measure on 
$\mathrm{S}^1$ which is invariant by $\Gamma$ and allows conjugating  
$\Gamma$ to a group of rotations. In this way, to complete the proof 
\index{rotation!group} 
of Margulis' theorem it suffices to show the following proposition.

\vspace{0.2cm}

\begin{prop} {\em If the action of a group $\Gamma$ by circle homeomorphisms 
is minimal and expansive, then $\Gamma$ contains a free subgroup on two generators.}
\label{propmargulis}
\end{prop}

\vspace{0.15cm}

Margulis' argument consists in proving directly an equivalent claim to this proposition 
by using the Klein Ping-Pong Lemma. The proof by Ghys also uses this idea, but it 
is based on a preliminary study of the ``maximal domain of contraction'' which 
not only allows proving the proposition but provides useful additional  
information.

\vspace{0.05cm}

\begin{defn} If the action of a group $\Gamma$ on the circle is minimal and expansive, 
then we will say that it is {\bf \textit{strongly expansive}} if, for every open interval 
$I$ whose closure is not the whole circle, there exists a sequence of elements
$g_n \! \in \! \Gamma$ such that the length $|g_n(I)|$ converges to zero.
\index{action!strongly expansive}
\end{defn}

\vspace{0.05cm}

Not every minimal expansive action is strongly expansive: consider for instance 
the case of subgroups of $\mathrm{PSL}_k(2,\mathbb{R})$ acting minimally, with 
$k \geq 2$. In what follows, we will see that, in the minimal, expansive case, 
these ``finite coverings'' are the only obstructions to strong expansiveness.

\vspace{0.1cm}

\begin{lem} {\em If the action of a subgroup $\Gamma$ of $\mathrm{Homeo}_+(\clo)$ 
is minimal and expansive, then there exists a finite order homeomorphism 
$R \!: \clo \rightarrow \clo$ commuting with all of the elements in $\Gamma$ 
such that the action of \esp $\Gamma$ on the quotient circle $\clo \!/\! \sim$ 
obtained as the space of orbits of $R$ is strongly expansive.}
\label{defn-theta}
\end{lem}

\noindent{\bf Proof.} For each $x \!\in\! \clo$ let us define $R(x) \!\in\! \clo$ 
as the ``supremum'' of the points $y$ for which the interval $]x,y[$ is 
contractable, \textit{i.e.,} there exists a sequence $(g_n)$ of elements in 
$\Gamma$ such that the length $|]g_n(x),g_n(y)[|$ converges to zero. The fact 
that $R$ commutes with all of the elements in $\Gamma$ follows immediately from 
the definition. Moreover, $R$ is a monotonous function ``of degree 1''. 
\index{degree!of a map}
We claim that $R$ is a finite order homeomorphism.

To show that $R$ is ``strictly increasing'' we argue by contradiction, and we consider the 
(non-empty) set $\mathrm{Plan}(R)$ formed by the union of the interiors of the intervals 
on which $R$ is constant. Since $R$ centralizes $\Gamma$, this set is invariant under the 
action. Hence, by minimality, $\mathrm{Plan}(R)$ is the whole circle, and this implies that 
$R$ is a constant map, which is impossible. To show that $R$ is continuous, one may use a 
similar argument using the union $\mathrm{Salt}(R)$ of the interior of the intervals which 
are avoided by the image of $R$ (compare with the proof of Theorem \ref{po-ref}).

The rotation number of \index{rotation number} 
$R$ cannot be irrational. Indeed, if it were irrational and $R$ admitted 
an exceptional minimal set then, since $R$ centralizes $\Gamma$, this set 
would be invariant by $\Gamma$, which contradicts the minimality of the action. 
If $\rho(R)$ were irrational and the orbits of $R$ were dense, 
then the unique invariant probability measure by $R$ would be invariant by $\Gamma$, 
and hence $\Gamma$ would be topologically conjugate to a group of rotations (see 
Exercise \ref{uni-inv}), which contradicts the expansiveness of the action.

The homeomorphism $R$ has finite order. Indeed, since its rotation number is rational, 
it admits periodic points. The set $\mathrm{Per}(R)$ of these points is closed, 
non-empty, and invariant by $\Gamma$ (once again, the last property follows from the 
fact that $R$ centralizes $\Gamma$). From the minimality of the action of $\Gamma$, 
one concludes that $\mathrm{Per}(R)$ coincides with the whole circle. To finish  
the proof of the lemma just notice that, by the definition of $R$, the action 
of $\Gamma$ on $\clo\!/\!\sim$ is strongly expansive. $\hfill\square$

\vspace{0.35cm}

In the context of the preceding lemma, the order of the homeomorphism $R$ will be 
called the {\bf \textit{degree}} 
\index{degree!of a group of circle homeomorphisms} 
of $\Gamma$, and will be denoted by $d(\Gamma)$. 
Notice that $d(\Gamma)=1$ if and only if $R$ is the identity, that is, if the 
original action of $\Gamma$ is (minimal and) strongly expansive.

\vspace{0.15cm}

Let us come back to the proof of Proposition \ref{propmargulis}. As we already 
mentioned, the strategy consists in finding two elements whose action on $\clo$, 
up to a ``finite covering'', is similar to the action of a Schottky subgroup of  
$\mathrm{PSL}(2,\mathbb{R})$. Indeed, such a pair of 
\index{group!Schottky}
elements  generates a free group, as the famous Klein Ping-Pong Lemma shows.
\index{Klein!Ping-Pong Lemma}

\vspace{0.15cm}

\begin{lem} {\em Let $\mathrm{M}$ be a set and $\Gamma$ a group of bijections of $\mathrm{M}$. 
Assume that there exist non-empty disjoint subsets $A_0$ and  $A_1$ of $\mathrm{M}$, and elements 
$g_0$ and $g_1$ in $\Gamma$, such that $g_0^n(A_1) \subset A_0$ and $g_1^n(A_0) \subset A_1$ for 
all $n \in \mathbb{Z} \setminus \{0\}$. Then $g_0$ and $g_1$ generate a free group.}
\label{klein-lem}
\end{lem}

\vspace{-0.3cm}

\noindent{\bf Proof.} We need to show that every 
nontrivial word $W = g_0^{i_1} g_1^{j_1} \cdots g_0^{i_k} g_1^{j_k}$ 
which is {\bf \textit{reduced}} (\textit{i.e.,} such that the exponents 
are non-zero, possibly with the exception of the first and/or the last ones) 
represents a nontrivial element in $\Gamma$. To do this it suffices to notice that, 
if $i_1 = 0$ and $j_k \neq 0$ (resp. if $i_1 \neq 0$ and $j_k = 0$), then (the element 
represented by) $W$ sends $A_0$ into $A_1$ (resp. $A_1$ into $A_0$), and hence it cannot 
equal the identity. Moreover, if $i_1$ and $j_k$ are both non-zero, then by conjugating 
appropriately $g$ one obtains a word for which this argument may be applied. 
$\hfill\square$

\vspace{0.4cm}

If the action of a group $\Gamma$ on the circle is minimal and expansive, then 
we consider the associated action on the quotient circle $\clo\!/\!\!\sim$. This 
action is minimal and strongly expansive, and hence there exist two sequences of 
intervals $L_n, L_n'$ in $\clo\!/\!\sim$ converging to points $a,b$ respectively, 
in such a way that for each $n \!\in\! \mathbb{Z}$ there exists $f_n \!\in\! \Gamma$ 
so that $f_n(\clo \setminus L_n) \!=\! L_n'$. We can assume that $a \neq b$, 
since otherwise we may replace $f_n$ by $g f_n$ for some element 
$g \!\in\! \Gamma$ which does not fix $a'$ (such an element exists by minimality). 
We claim that there exists $f \!\in\! \Gamma$ such that 
$a' \!=\! f(a)$ and $b' \!=\! f(b)$ 
are different from both $a$ and $b$. Although this could be left as an exercise, 
it also follows from a very interesting (and surprisingly not well-known\footnote{It 
is interesting to remark that some of the arguments in the proof of the Tits alternative 
for linear groups where the properties of Zariski's topology are strongly used may 
be simplified using Newmann's lemma.}) lemma due to Newmann.
\index{Newmann's lemma}

\vspace{0.1cm}

\begin{lem} {\em No group can be written as a finite union 
of left classes with respect to infinite index subgroups.}
\end{lem}

\noindent{\bf Proof.} Suppose that $\Gamma = S_1 [g_1] \cup \ldots \cup S_k [g_k]$ 
is a decomposition of a group into left classes $[g_i]$, where each set $S_i$ is 
finite and the classes $[g_i]$ are taken with respect to different subgroups. We 
will show by induction on $k$ that a certain class $[g_i]$ has finite index 
in $\Gamma$.

If \esp $k\!=\!1$ \esp there is nothing to prove. Suppose that the 
claim holds for $k \leq n$, and let us consider a decomposition 
$\Gamma = S_1 [g_1] \cup \ldots \cup S_{n+1} [g_{n+1}]$ as above. 
If $S_1 [g_1] = \Gamma$ then $[g_1]$ has finite index in $\Gamma$. 
In the other case, there exists $g \in \Gamma$ such that 
$g [g_1] \cap S_1 [g_1] = \emptyset$. We then have 
$$g [g_1] \subset S_2 [g_2] \cup \ldots \cup S_{n+1} [g_{n+1}],$$
and thus
$$S_1 [g_1] \subset T_2 [g_2] \cup \ldots \cup T_{n+1} [g_{n+1}],$$
where $T_i = \{ h_1 g^{-1} h_i \!: h_1 \in S_1, h_i \in S_i\}$ 
is finite. Hence,
$$\Gamma = (T_2 \cup S_2) [g_2] \cup \ldots \cup (T_{n+1} \cup S_{n+1}) [g_{n+1}],$$
and by the induction hypothesis this implies the existence of an index 
$j \!\in\! \{ 2, \ldots , n+1 \} $ such that $[g_j]$ has finite 
index in $\Gamma$. $\hfill\square$

\vspace{0.3cm}

The lemma above allows to choose $f \in \Gamma$ and 
$n \in \mathbb{N}$ so that the intervals $L_n,L_n',f(L_n),$ and 
$f(L_n')$, are disjoint. Indeed, $\Gamma$ cannot be written as the union 
$\Gamma_{a,a} \cup \Gamma_{a,b} \cup \Gamma_{b,a} \cup \Gamma_{b,b}$, where 
$\Gamma_{c,c'}$ denotes the --perhaps empty-- left class of the elements in 
$\Gamma$ sending $c$ into $c'$. Therefore, there must exist $f \!\in\! \Gamma$ 
such that $\{f(a),f(b)\} \cap \{a,b\} = \emptyset$, and hence the claimed property  
holds for $n \!\in\! \mathbb{N}$ large enough. To complete the proof of Proposition 
\ref{propmargulis}, we fix such an $n$, and we let $g_0 = f_n$ and 
$g_1 = f g_0 f^{-1}$. If we denote by $J_0$, $I_0$, $J_1$, and $I_1$, 
the preimages under $R$ of $L_n$, $L_n'$, $f(L_n)$, and $f(L_n')$, respectively, 
then one easily checks that for all $k \in \mathbb{Z} \setminus \{0\}$ one has 
$g_0^k ( I_1 \cup J_1) \subset I_0 \cup J_0$ and 
$g_1^k (I_0 \cup J_0) \subset I_1 \cup J_1$. Thus, one may apply  
Lemma \ref{klein-lem} to conclude that $g_0$ and $g_1$ generate a free group.

\begin{small} 
\begin{ejer} Show that a group $\Gamma$ of circle homeomorphisms preserves a probability 
measure on $\clo$ if and only if for every pair of elements in $\Gamma$ there is a common 
invariant probability measure.
\end{ejer}

\begin{ejer} Prove the following ``positive version'' of the 
Ping-Pong Lemma: If $g_0$ and $g_1$ are two bijections of a set $\mathrm{M}$ 
for which there exist disjoint non-empty subsets $A_0$ and $A_1$ such that 
$g_0^n (A_1) \subset A_0$ and $g_1^n (A_0) \subset A_1$ for every {\em positive} 
integer $n$, then $g_0$ and $g_1$ generate a free semigroup ({\em c.f.,}  
Definition \ref{def-semi-librre}).
\index{semigroup!free}
\label{pingpong+}
\end{ejer} 

\begin{ejer} Prove that Thompson's group F satisfies no nontrivial {\bf {\em law}}, 
\index{law (in a group)} 
that is, for every nontrivial word $W$ in two letters, there exist $f,g$ in F such  
that $W(f,g)$ is different from the identity. (This result holds for most groups 
of piecewise homeomorphisms of the interval, as was first shown in \cite{BS}.)

\noindent{\underbar{Hint.}} Let $W = g_0^{i_1} g_1^{j_1} \cdots g_0^{i_k} g_1^{j_k}$ 
be a nontrivial reduced word, with $i_1 \!\neq\! 0$ and $j_k \!\neq\! 0$. Take two 
elements $f$ and $g$ having at least $2k$ fixed points $x_1,\ldots,x_{2k}$ and 
$y_1,\ldots,y_{2k}$, respectively, so that \esp 
$x_1 < y_1 < x_2 < y_2 < \ldots < x_{2k} < y_{2k}$. 
\esp Use a similar argument to that of the proof of the Ping-Pong Lemma to 
show that, if $p$ belongs to the ``middle open interval'' determined by these 
points, then $W (f^n,g^n)(p)$ is different from $p$ for $n$ big enough.
\label{no-ley}
\end{ejer} \end{small}
\index{group!free|)}


\subsection{A probabilistic viewpoint}
\label{Tits-prob}

\hspace{0.45cm} Let $\Gamma$ be a countable group of circle homeomorphisms and $p$ a 
probability measure on $\Gamma$ that is {\bf \textit{non-degenerate}} (in the sense that 
its support 
\index{support!of a measure} 
generates $\Gamma$ as a semigroup). 
\index{semigroup}
Let us consider the {\bf {\em diffusion 
operator}} \index{diffusion operator} defined on the space of continuous functions by 
\begin{equation}
D \psi (x) = \int_{\Gamma} \psi \big( g (x) \big) \thinspace dp (g).
\label{difundir}
\end{equation}
Let \thinspace $\nu \mapsto p * \nu$ \thinspace be the dual action of this 
operator on the space of probability measures on the circle (this dual action 
is also called {\bf\em{convolution}}). 
\index{convolution of measures} Such a measure 
is said to be {\bf \textit{stationary}} \index{measure!stationary} 
(with respect to $p$) if \esp $p * \mu = \mu$, \esp 
that is, if for every continuous function $\psi \! : \clo \rightarrow \mathbb{R}$ one has 
\begin{equation}\int_{\clo} \psi(x) \thinspace d \mu(x) =
\int_{\Gamma} \int_{\clo} \psi \big( g(x) \big) \thinspace
d \mu(x) \thinspace d p (g).
\label{caca-def}
\end{equation}
The existence 
of at least one stationary measure is ensured by the Kakutani Fixed Point Theorem 
\index{Kakutani!Fixed Point Theorem} 
\cite{yosida}; clearly, it may be also deduced from a simple argument using Birkhoff 
\index{Birkhoff!sums}
sums (\textit{i.e.,} using the technique of Bogolioubov and Krylov: see Appendix B).
\index{Bogolioubov-Krylov's theorem}

\vspace{0.1cm}

\begin{lem} {\em If the orbits by $\Gamma$ are dense, then $\mu$ has total support 
and no atoms. If $\Gamma$ admits a minimal invariant Cantor set, then this set 
coincides with the support of $\mu$, and $\mu$ has no atoms either.}
\label{soporte}
\end{lem}

\noindent{\bf Proof.} Let us first show that if $\mu$ has atoms then $\Gamma$ 
admits finite orbits (concerning this case, see Exercise \ref{genial}). 
Indeed, if $x$ is a point with maximal positive mass, then from the equality 
$$\mu(x) = \int_{\Gamma} \mu \big( g^{-1}(x) \big) \thinspace dp (g)$$
one concludes that $\mu \big( g^{-1}(x) \big) \!=\! \mu(x)$ for every $g$ in 
the support of $p$. Since $p$ is non-degenerate, by repeating this argument one 
concludes that the equality $\mu \big( g^{-1}(x) \big) \!=\! \mu(x)$ actually 
holds for every element $g \!\in\! \Gamma$. Since the total mass of $\mu$ 
is finite, the only possibility is that the orbit of $x$ is finite.

If the action of $\Gamma$ is minimal then the support of $\mu$ being a closed 
invariant set, it must be the whole circle. If $\Gamma$ admits an exceptional minimal 
set $\Lambda$ then, since this set is unique, it must be contained in the support of 
$\mu$. Therefore, to prove that $\Lambda$ and $supp (\mu)$ coincide, we need to verify 
that $\mu(I)\!=\!0$ for every connected component of $\clo \setminus \Lambda$. Now, 
if this were not the case, then choosing such a component $I$ with maximal measure one 
would conclude --by an argument similar to that of the case of finite orbits-- that 
the orbit of $I$ is finite. However, this is in contradiction with the fact that 
the orbits of the endpoints of $I$ are dense in $\Lambda$. $\hfill \square$

\vspace{0.35cm}

The existence of stationary measures allows establishing a nontrivial result of regularity 
after conjugacy for general group actions on the circle \cite{der,DKN}. We point out that, 
according to \cite{harrison}, such a result is no longer true in dimension greater than 1 
(see however Exercise \ref{conj-abs}).

\vspace{0.15cm}

\begin{prop} {\em Every countable subgroup of \esp $\mathrm{Homeo}_+(\clo)$ 
(resp. $\mathrm{Homeo}_+([0,1])$) is topologically conjugate 
to a group of bi-Lipschitz homeomorphisms.}
\index{Lipschitz!homeomorphism}
\label{abu}
\end{prop}

\noindent{\bf Proof.} Let us first consider the case of a countable subgroup $\Gamma$ of 
$\mathrm{Homeo}_+(\clo)$ whose orbits are dense. Endow this group with a non-degenerate 
{\bf\em{symmetric}} probability measure $p$, where {\em symmetric} means that $p(g) \!= p(g^{-1})$ 
for every $g \in \Gamma$. Let us consider an associated stationary measure $\mu$ on $\clo$. 
For each interval $I \subset \clo$ and each element $g \!\in\! supp (p)$ one has 
$$\mu(I) \esp = \sum_{h \in supp (p)} \mu \big( h^{-1}(I) \big) \thinspace p(h)
\esp \geq \esp \mu \big( g (I) \big) \thinspace p (g^{-1}),$$
and hence
\begin{equation}
\mu \big( g(I) \big) \leq \frac{1}{p (g)} \thinspace \mu(I).
\label{Lip}
\end{equation}
Let us now take a circle homeomorphism $\varphi$ sending $\mu$ 
into the Lebesgue measure. If $J$ is an arbitrary interval in 
$\clo$ then, by (\ref{Lip}), for every $g \in supp (p)$ one has 
$$|\varphi \!\circ\! g \!\circ\! \varphi^{-1} (J)| =
\mu \big( g\! \circ\! \varphi^{-1} (J) \big) \leq
\frac{1}{p(g)} \thinspace \mu \big( \varphi^{-1}(J) \big)
= \frac{1}{p(g)} \thinspace |J|.$$
Therefore, for every $g \in supp (p)$ the homeomorphism 
$\varphi \circ g \circ \varphi^{-1}$ is Lipschitz with constant 
\thinspace $1/p (g)$. \thinspace Since $p$ is non-degenerate, 
this implies the proposition in the minimal case.

If $\Gamma$ is an arbitrary countable subgroup of $\mathrm{Homeo}_+(\clo)$, then adding 
an irrational rotation and considering the generated group, the problem reduces to that 
of dense orbits. By the arguments above, the new group --and hence the original one-- 
is topologically conjugate to a group of bi-Lipschitz homeomorphisms. Finally, 
by identifying the endpoints of the interval $[0,1]$, each subgroup $\Gamma$ of 
$\mathrm{Homeo}_+([0,1])$ induces a group of circle homeomorphisms with a marked 
global fixed point. Therefore, if $\Gamma$ is countable then this new group is 
conjugate by an element $\varphi$ of $\mathrm{Homeo}_+(\clo)$ to a group of 
bi-Lipschitz circle homeomorphisms. To obtain a genuine conjugacy inside 
$\mathrm{Homeo}_+([0,1])$, it suffices to compose $\varphi$ with a rotation in 
such a way that the marked point of the circle is sent to itself. $\hfill\square$

\vspace{0.1cm}

\begin{small} \begin{ejer} Using the argument of the preceding proof, show that, for each 
$\varepsilon\!>\!0$, every homeomorphism of the circle or the interval is topologically 
conjugate to a Lipschitz homeomorphism whose derivative (well-defined at almost every 
point) is less than or equal to $1+\varepsilon$ \esp (compare 
Exercise \ref{mas-adelante}).

\noindent{\underbar{Remark.}} Using the well-known inequality $h_{top}(T) \leq d \log(Lip(T))$ 
for the topological entropy of Lipschitz maps $T$ on $d$-dimensional compact manifolds, the 
claim above gives a short and conceptual proof of the fact that the topological entropy 
of any homeomorphism of the circle or the interval is zero (see \cite{walters}).
\index{topological entropy}
\end{ejer}

\begin{ejer} Prove that every countable group of homeomorphisms of a compact manifold 
is topologically conjugate to a group of absolutely continuous homeomorphisms.

\noindent{\underbar{Hint.}} Use the classical Oxtoby-Ulam's theorem which establishes that every 
probability measure of total support and without atoms on a compact manifold is the image of 
the Lebesgue measure by a certain homeomorphism (see \cite{GNW} for a concise presentation 
of this result). Use also the fact that every compact manifold supports minimal countable 
group actions (see \cite{Fathi-H} for minimal, {\em finitely generated} group actions).
\index{Oxtoby-Ulam's theorem}
\label{conj-abs}
\end{ejer}\end{small}

\vspace{0.02cm}

The preceding definitions extend to any action on a measurable space M of a countable group $\Gamma$ 
provided with a probability measure $p$. For example, looking at the action (by left translations) 
of $\Gamma$ on itself, the convolution operator may be iterated: the $n^{\mathrm{th}}$-convolution 
of $p$ with itself will be denoted by $p^{*n}$. 

For the general case, we will denote by $\Omega$ the space of the sequences 
$(g_1, g_2, \ldots) \in \Gamma^{\mathbb{N}}$ endowed with the product measure 
$\mathbb{P} = p^{\mathbb{N}}$. If $\sigma$ is the (one side) 
{\bf {\em shift}} on $\Omega$, 
\index{shift}
\esp that is, $\sigma(g_1,g_2,\ldots) = (g_2,g_3,\ldots)$, then one easily 
checks that a probability measure $\mu$ on M is stationary with respect to $p$ if 
and only if the measure $\mathbb{P} \times \mu$ is invariant by the 
{\bf \textit{skew-product}} map \esp 
$T \!: \Omega \times \mathrm{M} \rightarrow \Omega \times \mathrm{M}$ \esp given by
$$T(\omega,x) = \big( \sigma(\omega),h_1(\omega)(x) \big) =
\big( \sigma(\omega), g_1 (x) \big), \quad \esp 
\omega = (g_1,g_2,\ldots).$$


\begin{small}\begin{ejer} A continuous function 
$\psi\!: \mathrm{M} \rightarrow \mathbb{R}$ is said 
\index{harmonic function}
to be {\bf \textit{harmonic}} if it is invariant by the diffusion, that is, $D(\psi) = \psi$. 
Prove the following version of the {\bf \textit{maximum principle}}: if $\psi$ is harmonic, 
then the set of points at which $\psi$ attains its maximum is invariant by $\Gamma$. 
Prove that the same holds for {\bf \textit{super-harmonic}} functions, that is, for 
functions $\psi$ satisfying $D \psi \geq \psi$.
\label{func-harmonica}
\end{ejer}\end{small} 

\vspace{-0.35cm}

We now concentrate on the case where M is a compact metric space. Following the 
seminal work of Furstenberg \cite{furst}, in order to study the evolution of the 
random compositions we consider the {\bf \textit{inverse process}} given by \thinspace 
$\bar{h}_n(\omega) = g_1 \cdots g_n$. \thinspace The main reason for doing this 
lies in the following observation: if $\psi$ is a continuous function defined on M 
and $\mu$ is a stationary measure, then the sequence of random variables 
$$\xi_n(\omega) = \int_{\mathrm{M}} \psi \thinspace d \big(g_1 \cdots g_n(\mu)\big)
= \int_{\mathrm{M}} \psi \esp d \big( \bar{h}_n(\omega)(\mu) \big)$$
is a martingale \index{martingale}
\cite{feller}. Indeed, for every $g_1,\ldots,g_n$ in $\Gamma$, 
an equality of type (\ref{caca-def}) applied to the function \esp 
$x \mapsto \psi \big( g_1 \cdots g_n (x) \big)$ \esp yields 
$$\int_{\Gamma} \int_{\mathrm{M}} \psi \esp d \big( g_1 \cdots g_n g (\mu) \big)
\esp dp(g) = \int_{\mathrm{M}} \psi \esp d \big( g_1 \cdots g_n (\mu) \big),$$
that is, \esp $\mathbb{E}(\xi_{n+1}|g_{n+1}) = \xi_n$.
By the Martingale Convergence Theorem, the sequence $\big( \xi_n(\omega) \big)$ 
converges for almost every $\omega \in \Omega$. Let $(\psi_k)$ be a dense 
sequence in the space of continuous functions on M. By the 
compactness of the space of probability measures on M, 
for a total probability subset $\Omega_0 \subset \Omega$, 
\index{weak-* topology}
the sequence 
$g_1 g_2 \cdots g_n (\mu) = \bar{h}_n(\omega)(\mu)$ converges (in the 
weak-* topology) to a probability $\omega(\mu)$. Moreover, the map 
$\omega \mapsto \omega(\mu)$ is well-defined at almost every $\omega$ 
and it is measurable (see \cite[page 199]{Ma2} for more details on this).

\vspace{0.1cm}

\begin{prop} {\em Let $\Gamma$ be a countable subgroup of $\mathrm{Homeo}_+(\clo)$ whose 
action is minimal. If the property of strong expansiveness is satisfied, then for almost 
every sequence $\omega \in \Omega_0$ the measure $\omega(\mu)$ is a Dirac measure.}
\index{action!strongly expansive}
\label{cont-top}
\end{prop}

\noindent{\bf Proof.} We will show that for every $\varepsilon\!\in ]0,1]$ there 
exists a total probability subset $\Omega_{\varepsilon}$ of $\Omega_0$ such that, 
for every $\omega \in \Omega_{\varepsilon}$, there exists an interval $I$ of length 
$|I| \leq \varepsilon$ such that $\omega(\mu)(I) \geq 1 - \varepsilon$. This allows 
us to conclude that, for all $\omega$ contained in the total probability subset 
\thinspace $\Omega^* = \cap_{n \in \mathbb{N}} \Omega_{1/n}$, 
\thinspace the measure $\omega(\mu)$ is a Dirac measure.

Let us then fix $\varepsilon \!>\! 0$. For each $n \in \mathbb{N}$ let us denote by 
$\Omega^{n,\varepsilon}$ the set of the sequences $\omega \!\in\! \Omega_0$ such that, 
for all $m \geq 0$ and every interval $I$ in $\clo$ of length $|I| \leq \varepsilon$, 
one has $\bar{h}_{n+m}(\omega)(\mu)(I) < 1-\varepsilon$. We need to show that the 
probability of $\Omega^{n,\varepsilon}$ is zero. We will begin by exhibiting 
a finite subset $\mathcal{G}_{\varepsilon}$ of $supp (p)$, as well as an integer 
$l \in \mathbb{N}$, such that for all $r \in \mathbb{N}$ and all $(g_1,\cdots,g_r) 
\in \Gamma^r$ there exist an interval $I$ of length $|I| \leq \varepsilon$, an integer 
$\ell \leq 2l$, and elements $f_1,\ldots,f_{\ell}$ in $\mathcal{G}_{\varepsilon}$, 
satisfying
\begin{equation}
g_1 \cdots g_r f_1 \cdots f_{\ell} (\mu) (I) \geq 1-\varepsilon.
\label{aprobar}
\end{equation}
To do this, let us fix two different points $a$ and $b$ in $\clo$, as well as 
an integer $q > 1 / \varepsilon$, and let us take $q$ different points $a_1,\ldots,a_q$ 
in the orbit of $a$ by $\Gamma$. For each $i \! \in \! \{1,\ldots,q\}$ let us choose 
$h_i \in \Gamma$ and an open interval $U_i$ containing $a_i$ in such a way that the 
$U_i$'s are two-by-two disjoint, $h_i(a) = a_i$, and $h_i(U) = U_i$ for some 
neighborhood $U$ of $a$ not containing $b$. Let us now consider a neighborhood  
$V$ of $b$ disjoint from $U$ and such that $\mu(\clo \setminus V) \geq 1 - \varepsilon$. 
By minimality and strong expansiveness, there exists $h \! \in \!\Gamma$ such that 
$h(\clo \setminus V) \subset U$. Now each element in $\{h_1,\ldots,h_q,h\}$ may be 
written as a product of elements in the support of $\mu$. This may be done in many 
different ways, but if we fix once and for all a particular choice for $h$ and each $h_i$, 
then the set $\mathcal{G}_{\varepsilon}$ of the elements in $supp(p)$ which are used is 
finite. Let $l$ the maximal number of factors appearing in these choices. To check 
(\ref{aprobar}) notice that, for $g = g_1 \cdots g_r$, the intervals $g(U_i)$ are 
two-by-two disjoint, and hence the length of at least one of them is bounded from 
above by $\varepsilon$. For such an interval $I = g(U_i)$ we have  
$$g_1 \cdots g_r h_i h (\mu) (I) = \mu \big( h^{-1} (U) \big)
\geq \mu(\clo \setminus V) \geq 1 - \varepsilon,$$
which concludes the proof of (\ref{aprobar}).

Now let $\rho = \min \{p (f) \!: f \in \mathcal{G}_{\varepsilon} \}$, and let us 
define $\Omega^{\varepsilon}_{n,m}$ as the set of $\omega \in \Omega_0$ 
such that, for every interval $I$ of length $|I| \leq \varepsilon$ and all 
$k \leq m$, one has $\bar{h}_{n+k}(\omega)(\mu)(I) < 1 - \varepsilon$. By 
(\ref{aprobar}) we have
$$\mathbb{P}(\Omega_{n,2 l t}^{\varepsilon})
\leq (1-\rho^{2 l})^t.$$ 
Passing to the limit as $t$ goes to infinity, this allows to conclude 
that $\mathbb{P}(\Omega^{n,\varepsilon}) = 0$, which 
finishes the proof.  $\hfill\square$

\vspace{0.3cm}

Using a well-known argument in the theory of random walks on groups, Proposition 
\ref{cont-top} allows us to prove a general uniqueness result for the stationary measure. 
Let us point out that, according to \cite{DK}, this result still holds in the 
much more general context of codimension-one foliations (the notion of stationary 
measure for this case is that of Garnett: see \cite{candel,garnett}).

\vspace{0.1cm}

\begin{thm} {\em Let $\Gamma$ be a countable group of circle homeomorphisms endowed with 
a non-degenerate probability measure $p$. If $\Gamma$ does not preserve any probability 
measure on $\clo$, then the stationary measure with respect to $p$ is unique.}
\label{unicidad}
\end{thm}

\noindent{\bf Proof.} First suppose that the action of $\Gamma$ is minimal 
and strongly expansive, and let $\mu$ be a probability measure on $\clo$ which 
is stationary with respect to $p$. For each $\omega \in \Omega$ such that \esp 
$\lim_{n\rightarrow \infty} \bar{h}_n(\omega) (\mu)$ \esp exists and is a Delta 
measure, let us denote by $\varsigma_{\mu} (\omega)$ the atom of the measure 
$\omega(\mu)$, \textit{i.e.,} the point in $\clo$ for which  
$\omega(\mu) = \delta_{\varsigma_{\mu}(\omega)}$.
The map $\varsigma_{\mu}\!: \Omega \rightarrow \clo$ is well-defined at  
almost every sequence and measurable. We claim that the measures $\mu$ and 
$\varsigma_{\mu}(\mathbb{P})$ coincide. Indeed, since $\mu$ is stationary, 
$$\mu = p^{*n} * \mu = \sum_{g \in \Gamma} p^{*n}(g) \thinspace g(\mu) =
\int_{\Omega} \bar{h}_n(\omega) (\mu) \thinspace d \mathbb{P}(\omega).$$
Therefore, passing to the limit as $n$ goes to infinity, we obtain 
$$\mu = \int_{\Omega} \lim_{n \rightarrow \infty} h_n(\omega) (\mu) \thinspace d
\mathbb{P}(\omega)
= \int_{\Omega} \delta_{\varsigma_{\mu}(\omega)} \thinspace d \mathbb{P}(\omega) 
= \varsigma_{\mu}(\mathbb{P}).$$

Now consider two stationary probabilities $\mu_1$ and $\mu_2$. The measure 
$\mu = (\mu_1 + \mu_2)/2$ is also a stationary probability, and the function 
$\varsigma_{\mu}$ satisfies, for $\mathbb{P}$-almost every $\omega \in \Omega$,
$$\frac{\delta_{\varsigma_{\mu_1}(\omega)} + \delta_{\varsigma_{\mu_2}(\omega)}}{2} 
= \delta_{\varsigma_{\mu}(\omega)}.$$
Clearly, this is impossible unless $\varsigma_{\mu_1}$ and $\varsigma_{\mu_2}$
coincide almost surely. Thus, by the claim of the first part of the proof, 
$$\mu_1 = \varsigma_{\mu_1}(\mathbb{P}) = \varsigma_{\mu_2}(\mathbb{P}) = \mu_2.$$ 

Suppose now that the action of $\Gamma$ is minimal and expansive, but not strongly 
\index{action!expansive}
expansive. Let $\mu$ be a stationary probability with respect to $p$. By Lemma 
\ref{defn-theta}, there exists a finite order homeomorphism $R\!: \clo \rightarrow \clo$  
which commutes with all the elements of $\Gamma$ and such that 
the action induced on the topological circle $\clo\!/\!\!\sim$ obtained as the space 
of orbits by $R$ is minimal and strongly expansive. For each $x \in \clo$ let us denote 
$\psi(x)\! = \! \mu \big( [x,R(x)[ \big)$. Since $\mu$ has no atom, $\psi$ is a continuous 
function, and since $R$ centralizes $\Gamma$, it is harmonic. Therefore, the set of 
the points at which $\psi$ attains its maximum value is invariant by $\Gamma$ (see  
Exercise \ref{func-harmonica}). Since the orbits by $\Gamma$ are dense, $\psi$ is 
constant; in other words, $\mu$ is invariant under $R$. On the other hand, $\mu$ 
projects into a stationary probability measure for the action of $\Gamma$ on 
$\clo \! / \!\! \sim$. By the first part of the proof, this projected measure 
is unique, and together with the $R$-invariance of $\mu$, this proves that 
$\mu$ is unique as well.

If $\Gamma$ admits an exceptional minimal set,
\index{exceptional minimal set} 
then this set coincides with the support of $\mu$. 
By collapsing to a point each connected component in its complementary set, one obtains a minimal 
action. If this action is expansive, then the arguments above give the uniqueness of the stationary 
measure. To complete the proof it suffices to notice that, in all the cases which have not 
been considered, $\Gamma$ preserves a probability measure of the circle. $\hfill \square$

\vspace{0.28cm}

Let us define the {\bf \textit{contraction coefficient}} $\mathrm{contr}(h)$ of a circle 
homeomorphism $h$ as the infimum of the numbers $\varepsilon > 0$ such that there exist 
closed intervals $I$ and $J$ in $\clo$ of length less than or equal to $\varepsilon$ 
such that $h(\overline{\clo \setminus I}) = J$. This notion allows us to give a ``topological 
version'' of Proposition \ref{cont-top} for the composition in the ``natural order''.

\vspace{0.12cm}

\begin{prop} {\em Under the hypothesis of Proposition {\em \ref{cont-top}}, for almost 
every sequence $\omega = (g_1,g_2,\ldots) \in \Omega$ the contraction coefficient of 
$h_n(\omega)\!=\! g_n \cdots g_1$ converges to zero as $n$ goes to infinity.}
\label{ultima}
\end{prop}

\noindent{\bf Proof.} Since $\mu$ has total support and no atoms, there exists a 
circle homeomorphism $\varphi$ sending $\mu$ into the Lebesgue measure. Hence, 
since the claim to be proved is invariant by topological conjugacy, we may 
assume that $\mu$ coincides with the Lebesgue measure.

Let $\bar{p}$ be the probability on $\Gamma$ defined by $\bar{p}(g)=p(g^{-1})$, and let 
$\bar{\Omega}$ be the probability space $\Gamma^{\mathbb{N}}$ endowed with the measure 
$\bar{p}^{\mathbb{N}}$. On this space let us consider the inverse process \esp 
$\bar{h}_n(\bar{\omega})=g_1 \cdots g_n$, \esp where 
$\bar{\omega}= (g_1,g_2,\ldots)$. \esp By Proposition 
\ref{cont-top}, for almost every $\bar{\omega} \in \bar\Omega$ and all 
$\varepsilon > 0$ there exists a positive integer $n(\varepsilon,\bar{\omega})$ 
such that, if $n \geq n(\varepsilon,\bar{\omega})$, then there exists a closed 
interval $I$ such that $\mu(I) \leq \varepsilon$ and $\bar{h}_n(\bar{\omega})(\mu)(I) 
\geq 1 -\varepsilon$. If we denote by $J$ the closure of $\clo \setminus
g_n^{-1} \cdots g_1^{-1}(I)$, \esp then one easily checks that 
$|I| = \mu(I) \leq \varepsilon$,
$$|J| = 1 - |g_n^{-1} \cdots g_1^{-1} (I)|
= 1 - \mu \big( \bar{h}_n (\bar{\omega})^{-1}(I) \big)
= 1 - \bar{h}_n(\bar{\omega})(\mu)(I) \leq \varepsilon,$$
and \esp $g_n^{-1} \cdots g_1^{-1} (\overline{\clo \setminus I}) \!=\! J$. 
Therefore,
\esp $\mathrm{contr} \big( g_n^{-1} \cdots g_1^{-1} \big) \!\leq\! \varepsilon$ \esp
for all 
$n \!\geq\! n(\varepsilon,\bar{\omega})$. The proof is then finished by noticing that 
the map \esp $(g_1,g_2,\ldots) \mapsto (g_1^{-1},g_2^{-1},\ldots)$ \esp identifies  
the spaces $(\bar{\Omega},\bar{p}^{\mathbb{N}})$ and $(\Omega,p^{\mathbb{N}})$. 
$\hfill\square$

\vspace{0.5cm}

The contraction coefficient is always realized, in the sense 
that for every circle homeomorphism $h$ there exist intervals $I,J$ such that 
\esp $\max \{ |I|,|J| \} \!=\! \mathrm{contr}(h)$ \esp and 
$h(\overline{\clo \setminus I}) \!=\! J$ (however, these intervals 
are not necessarily unique). As a consequence, and according to the preceding 
proof, for almost every $\omega \in \Omega$ we may choose two sequences of closed 
intervals $I_n(\omega)$ and $J_n(\omega)$ whose lengths go to zero and such that 
$h_n(\omega) \big( \overline{\clo \setminus I_n(\omega)} \big) \!=\! J_n(\omega)$ 
for all $n \in \mathbb{N}$. For any of these choices, the intervals $I_n(\omega)$ 
converge to the point $\varsigma_{\mu}(\omega)$.

\vspace{0.1cm}

\begin{small} \begin{ejer} Let $\Gamma$ be a group of circle homeomorphisms 
acting minimally.

\vspace{0.08cm}

\noindent (i) Show that if the $\Gamma$-action is strongly expansive, 
\index{centralizer}
then the centralizer of $\Gamma$ in $\mathrm{Homeo}_+(\clo)$ is trivial. 

\vspace{0.08cm}

\noindent (ii) More generally, show that if the $\Gamma$-action is expansive, 
then the only circle homeomorphisms commuting with all the elements of 
$\Gamma$ are the powers of the map $R$ from Lemma \ref{defn-theta}.
\end{ejer}

\begin{ejer} Using the preceding exercise, prove that if a product of groups 
$\Gamma = \Gamma_1 \times \Gamma_2$ acts on the circle, then at least 
one of the factors preserves a probability measure on $\clo$. 
\label{proba-para-productos}
\end{ejer}

\begin{ejer} The results in this section were obtained by the author 
in collaboration with Deroin and Kleptsyn in \cite{DKN} (partial results appear in 
\cite{victor}). Nevertheless, we must point out the existence of a prior work on this 
topic, namely the article \cite{antonov} by Antonov, where almost equivalent results 
\index{Antonov}\index{Deroin}\index{Kleptsyn}are 
stated (and proved) in a purely probabilistic language. 
The exercise below contains the essence of \cite{antonov}.

\vspace{0.1cm}

\noindent (i) Given a probability measure $p$ on a (non-necessarily finite) system of 
generators of a countable subgroup $\Gamma$ of $\mathrm{Homeo}_+(\clo)$ without 
finite orbits, consider the probability $\bar{p}$ on $\Gamma$ 
defined by $\bar{p}(g)~=~p(g^{-1})$ 
\esp (compare Proposition \ref{ultima}). If $\bar{\mu}$ is a measure on $\clo$ 
which is stationary with respect to $\bar{p}$, show that for every $x,y$ in $\clo$ 
the sequence of random variables
$$\xi_n^{x,y}(\omega) = \bar{\mu} ( [g_n \cdots g_1 (x), g_n \cdots g_1 (y)] )$$
is a martingale. \index{martingale}
In particular, this sequence converges almost surely. We want to show that, 
if there is no nontrivial circle homeomorphism centralizing $\Gamma$, then 
the corresponding limit equals $0$ or $1$.

\vspace{0.1cm}

\noindent (ii) Let $\nu$ be a stationary measure for the diagonal action of $\Gamma$ 
on the torus $\clo \!\times\! \clo$. Show that the distribution of $\xi_n^{x,y}$ with 
respect to $\mathbb{P} \times \nu$ coincides with that of $\xi_{n+1}^{x,y}$.

\vspace{0.1cm}

\noindent (iii) Using the relation (which holds for every square integrable martingale)
$$\mathbb{E} (\xi_{n+1}^2) = \mathbb{E} (\xi_n^2) + \mathbb{E} ( (\xi_{n+1} - \xi_n)^2),$$
conclude that for $\nu$-almost every point $(x,y)\in\clo\times\clo$ one has \esp 
$\xi_{n+1}^{x,y} = \xi_n^{x,y}.$

\vspace{0.1cm}

\noindent (iv) Defining $\psi ( (x,y) ) \!=\! \bar{\mu} ( [x,y] )$, conclude from 
(iii) that for $\nu$-almost every point $(x,y)$ in $\clo\times\clo$ and all 
$g \in \Gamma$ one has \esp $\psi ( (x,y) ) = \psi ( (g(x),g(y)) ).$

\vspace{0.1cm}

\noindent (v) Using (iv), prove that the support 
of the measure $\nu$ is contained in the diagonal.

\vspace{0.08cm}

\noindent{\underbar{Hint.}} Without loss of generality, one can 
restrict to the case where all of the orbits of $\Gamma$ on $\clo$ 
are dense. Suppose that $\alpha \! \in ]0,1[$ is such that the set
$$X_{\alpha,\varepsilon} =
\{(x,y) \in \clo \!\times\! \clo\!: \psi((x,y)) 
\in [\alpha-\varepsilon,\alpha+\varepsilon] \}$$
has positive $\nu$-measure for every $\varepsilon > 0$. Conclude from (iv)
that the normalized restriction of $\nu$ to $X_{\alpha,\varepsilon}$ is 
a stationary probability. Letting $\varepsilon$ go to $0$ and passing 
to a weak limit, this gives a stationary probability $\bar{\nu}$ on 
the torus concentrated on the set
$$\{ (x,y) \in \clo \!\times\! \clo\!: \esp \psi ( [x,y] )= \alpha\}.$$
The density of the orbits then allows defining a unique homeomorphism 
$R$ of $\clo$ that satisfies $\psi ( (x,R(x)) ) = \alpha$. Show that 
$R$ commutes with all the elements of $\Gamma$, thus contradicting 
the hypothesis at the end of item (i).

\vspace{0.1cm}

\noindent (vi) Fixing $x,y$ in $\clo$, consider the Dirac measure 
$\delta_{(x,y)}$ with mass on the point $(x,y) \in \clo\!\times\!\clo$. 
Show that every measure in the adherence of the set of measures 
$$\nu_n = \frac{1}{n} \sum_{j=0}^{n-1} p^{*j} * \delta_{(x,y)}$$
is a stationary probability concentrated on the diagonal. Finally, 
using the fact that $\xi_n^{x,y}$ converges almost surely, 
conclude that the limit of $\xi_n^{x,y}$ is equal to 0 or 1.
\end{ejer}

\begin{ejer} Prove that if $\Gamma$ is a subgroup of $\mathrm{Homeo}_+([0,1])$ 
without global fixed point in the interior and $p$ is a symmetric non-degenerate 
probability measure on $\Gamma$, then every probability measure on $[0,1]$ which 
is stationary (with respect to $p$) is supported at the endpoints of $[0,1]$.

\noindent{\underbar{Hint.}} Let $\mu$ be a stationary measure supported on $]0,1[$. 
Show first that $\mu$ has no atom. Then, by collapsing the connected components of 
the complement of the support and reparameterizing the interval, reduce the 
general case to that where $\mu$ coincides with the Lebesgue measure. One  
then has, for all $s\!\in ]0,1[$,
$$s = \mu ( [0,s] ) = \int_{\Gamma} \mu ( g^{-1}([0,s]) ) \thinspace dp (g)
= \int_{\Gamma} \frac{ \mu ( g([0,s]) ) +
\mu ( g^{-1}([0,s]) )}{2} \thinspace dp (g),$$
which gives
$$s = \int_{\Gamma} \frac{g(s) + g^{-1}(s)}{2} \thinspace dp (g).$$
Hence, by integrating between $0$ and an arbitrary point $t\!\in ]0,1[$,
\begin{equation}
t^2 = \int_{\Gamma} \int_{0}^{t} (g(s) + g^{-1}(s))
\thinspace ds \thinspace dp (g).
\end{equation}
Now by looking at Figure 14, conclude that for every homeomorphism 
$f$ of the interval and every $t\!\in \![0,1]$ one has
$$\int_0^t ( f(s) + f^{-1}(s) ) \thinspace ds \geq t^2,$$
where the equality holds if and only if $f(t) = t$. 
Conclude that $t$ is a global fixed point for the action, thus 
contradicting the hypothesis (see \cite{DKN} for an alternative 
and more conceptual proof using Garnett's ergodic theorem).

\noindent{\underbar{Remark.}} The hypothesis of symmetry for $p$ 
above is necessary, as is shown by \cite{ka-intervalo}.
\label{genial}
\end{ejer} \end{small}

\vspace{0.6cm}


\beginpicture

\setcoordinatesystem units <1cm,1cm>

\putrule from 0 0 to 0 5
\putrule from 5 0 to 5 5
\putrule from 0 5 to 5 5
\putrule from 0 0 to 5 0

\plot
0 0
0.7071 0.1
1 0.2
1.5811 0.5
2 0.8
2.236 1
2.5495 1.3
2.8482  1.6
3.1622 2
3.3911 2.3
3.6055 2.6
3.8729 3
4.1833 3.5
4.4721 4
4.7434 4.5
5 5 /

\plot
0 0
0.1 0.7071
0.2 1
0.5 1.5811
0.8 2
1 2.236
1.3 2.5495
1.6 2.8482
2 3.1622
2.3 3.3911
2.6 3.6055
3 3.8729
3.5 4.1833
4 4.4721
4.5 4.7434
5 5 /

\plot 0 0 0.1 0.1 /
\plot 0.2 0.2 0.3 0.3 /
\plot 0.4 0.4 0.5 0.5 /
\plot 0.6 0.6  0.7 0.7 /
\plot 0.8 0.8  0.9 0.9 /
\plot 1 1  1.1 1.1 /
\plot 1.2 1.2  1.3 1.3 /
\plot 1.4 1.4  1.5 1.5 /
\plot 1.6 1.6  1.7 1.7 /
\plot 1.8 1.8  1.9 1.9 /
\plot 2 2  2.1 2.1 /
\plot 2.2 2.2  2.3 2.3 /
\plot 2.4 2.4  2.5 2.5 /
\plot 2.6 2.6  2.7 2.7 /
\plot 2.8 2.8  2.9 2.9 /
\plot 3 3  3.1 3.1 /
\plot 3.2 3.2  3.3 3.3 /
\plot 3.4 3.4  3.5 3.5 /
\plot 3.6 3.6  3.7 3.7 /
\plot 3.8 3.8  3.9 3.9 /
\plot 4 4  4.1 4.1 /
\plot 4.2 4.2  4.3 4.3 /
\plot 4.4 4.4  4.5 4.5 /
\plot 4.6 4.6  4.7 4.7 /
\plot 4.8 4.8  4.9 4.9 /

\putrule from 3.5 0 to 3.5 0.2
\putrule from 3.5 0.4 to 3.5 0.6
\putrule from 3.5 0.8 to 3.5 1
\putrule from 3.5 1.2 to 3.5 1.4
\putrule from 3.5 1.6 to 3.5 1.8
\putrule from 3.5 2 to 3.5 2.2
\putrule from 3.5 2.4 to 3.5 2.6
\putrule from 3.5 2.8 to 3.5 3
\putrule from 3.5 3.2 to 3.5 3.4
\putrule from 3.5 3.6 to 3.5 3.8
\putrule from 3.5 4 to 3.5 4.2

\putrule from 0 3.5 to 0.2 3.5
\putrule from 0.4 3.5 to 0.6 3.5
\putrule from 0.8 3.5 to 1 3.5
\putrule from 1.2 3.5 to 1.4 3.5
\putrule from 1.6 3.5 to 1.8 3.5
\putrule from 2 3.5 to 2.2 3.5
\putrule from 2.4 3.5 to 2.6 3.5
\putrule from 2.8 3.5 to 3 3.5
\putrule from 3.2 3.5 to 3.4 3.5

\putrule from 0 4.18 to 0.2 4.18
\putrule from 0.4 4.18 to 0.6 4.18
\putrule from 0.8 4.18 to 1 4.18
\putrule from 1.2 4.18 to 1.4 4.18
\putrule from 1.6 4.18 to 1.8 4.18
\putrule from 2 4.18 to 2.2 4.18
\putrule from 2.4 4.18 to 2.6 4.18
\putrule from 2.8 4.18 to 3 4.18
\putrule from 3.2 4.18 to 3.4 4.18

\put{$f$} at 2.05 1.2
\put{$f^{-1}$} at 1.5 2.3
\put{$A$} at 2.5 0.5
\put{$A$} at 0.45 2.6
\put{$B$} at 2.9 2.3
\put{$\Delta$} at 3.15 3.725

\put{$0$} at 0 -0.3
\put{$t$} at 3.5 -0.3
\put{$1$} at 5 -0.3
\put{$t$} at -0.6 3.5
\put{$f^{-1}(t)$} at -0.6 4.18

\put{$A = \int_{0}^{t} f(s) \thinspace ds$} at 8 3.5
\put{$B = \int_0^t f^{-1}(s) \thinspace ds$ } at 8 2.8
\put{$A + B = t^2 + \Delta$} at 8 1.2


\setdots

\putrule from 0.7071 0.1 to 3.5 0.1
\putrule from 1 0.2 to 3.5 0.2
\putrule from 1.2 0.3 to 3.5 0.3
\putrule from 1.4 0.4 to 3.5 0.4
\putrule from 1.5811 0.5 to 3.5 0.5
\putrule from 1.8 0.6 to 3.5 0.6
\putrule from 1.9 0.7 to 3.5 0.7
\putrule from 2 0.8 to 3.5 0.8
\putrule from 2.15 0.9 to 3.5 0.9
\putrule from 2.236 1 to 3.5 1
\putrule from 2.35 1.1 to 3.5 1.1
\putrule from 2.45 1.2 to 3.5 1.2
\putrule from 2.5495 1.3 to 3.5 1.3
\putrule from 2.67 1.4 to 3.5 1.4
\putrule from 2.75 1.5 to 3.5 1.5
\putrule from 2.8482 1.6 to 3.5 1.6
\putrule from 2.9 1.7 to 3.5 1.7
\putrule from 3.04 1.8 to 3.5 1.8
\putrule from 3.1 1.9 to 3.5 1.9
\putrule from 3.1622 2 to 3.5 2
\putrule from 3.24 2.1 to 3.5 2.1
\putrule from 3.32 2.2 to 3.5 2.2
\putrule from 3.3911 2.3 to 3.5 2.3

\putrule from 0.1 0.7071 to 0.1 3.5
\putrule from 0.2 1 to 0.2 3.5
\putrule from 0.3 1.2 to 0.3 3.5
\putrule from 0.4 1.4 to 0.4 3.5
\putrule from 0.5 1.5811 to 0.5 3.5
\putrule from 0.6 1.8 to 0.6 3.5
\putrule from 0.7 1.9 to 0.7 3.5
\putrule from 0.8 2 to 0.8 3.5
\putrule from 0.9 2.15 to 0.9 3.5
\putrule from 1 2.236 to 1 3.5
\putrule from 1.1 2.35 to 1.1 3.5
\putrule from 1.2 2.45 to 1.2 3.5
\putrule from 1.3 2.5495 to 1.3 3.5
\putrule from 1.4 2.67 to 1.4 3.5
\putrule from 1.5 2.75 to 1.5 3.5
\putrule from 1.6 2.8482 to 1.6 3.5
\putrule from 1.7 2.9 to 1.7 3.5
\putrule from 1.8 3.04 to 1.8 3.5
\putrule from 1.9 3.1 to 1.9 3.5
\putrule from 2 3.1622 to 2 3.5
\putrule from 2.1 3.24 to 2.1 3.5
\putrule from 2.2 3.32 to 2.2 3.5
\putrule from 2.3 3.3911 to 2.3 3.5

\putrule from 2.3 3.5 to 2.3 3.3911
\putrule from  2.5 3.5 to 2.5 3.6055
\putrule from  2.7 3.5 to 2.7 3.72
\putrule from  2.9 3.5 to 2.9 3.8
\putrule from 3.1 3.5 to 3.1 3.9
\putrule from 3.3 3.5 to 3.3 4.1

\begin{tiny}
\put{$\cdot$} at 3.37 2.27
\put{$\cdot$} at 2.16 3.41
\end{tiny}

\put{Figure 14} at 4.25 -0.9

\begin{Large}
\put{$\Downarrow$} at 7.6 2
\end{Large}

\put{} at -2.3 0

\endpicture


\newpage
\thispagestyle{empty}
${}$
\newpage


\chapter{DYNAMICS OF GROUPS OF DIFFEOMORPHISMS}
\label{en-dos}


\section{Denjoy's theorem}
\label{denj-thm}\index{Denjoy!theorem|(}

\hspace{0.45cm} In \S \ref{teoria-poinc}, we have seen an example of a circle 
homeomorphism with irrational rotation number admitting an exceptional minimal set.
\index{exceptional minimal set} In this section, we will see that such 
a homeomorphism cannot be a diffeomorphism of certain regularity. This 
is the content of a classical and very important theorem due to Denjoy. 
To state it properly, we will denote by $\vl$ the group of $\ce^1$ 
circle diffeomorphisms whose derivatives have bounded variation. 
\index{variation of a function!total variation}
We will say that such a 
diffeomorphism is of class $\mathrm{C}^{1 + \mathrm{bv}}$, and we will denote 
by $V(f)$ the total variation 
of the logarithm of its derivative, that is, 
$$V(f) \esp \esp = \sup\limits_{a_0 < \ldots < a_n = a_0} \sum\limits_{i = 0}^{n-1}
\big| \log(f')(a_{i+1}) - \log(f')(a_i) \big| \esp = 
\esp \mathrm{var} \big( \log(f');\clo \big).$$
For an interval $I$ we will use the notation 
$V(f;I) = \mathrm{var} \big( \log(f')|_I \big)$.

\vspace{0.1cm}

\begin{thm} {\em If $f$ is a circle diffeomorphism of class $\mathrm{C}^{1 + \mathrm{bv}}$ with 
irrational rotation number, \index{rotation number}
then $f$ is topologically conjugate to the rotation of angle $\rho(f)$.}
\label{denj}
\end{thm}

\vspace{0.1cm}

Every circle diffeomorphism of class $\mathrm{C}^2$ or $\ce^{1+\mathrm{Lip}}$  
belongs to $\vl$. For instance, for every $f \in \dos$ one has 
\begin{equation}
V(f) = \int_{\mathrm{S}^1} \left| \frac{f''(s)}{f'(s)} \right| ds.
\label{vardos}
\end{equation}
The same formula holds for difeomorphisms of class $\mathrm{C}^{1+\mathrm{Lip}}$.   
(In this case the function $s \mapsto f''(s)/f'(s)$ is almost everywhere defined 
and essentially bounded.) However, let us point out that Denjoy's theorem does not 
longer hold in class $\mathrm{C}^{1+\tau}$ for any $\tau \!<\! 1$. This is a 
consequence of a construction due to Herman \cite{He} which we will 
\index{Herman}
reproduce in \S \ref{ejemplosceuno}.\index{Denjoy!counterexamples}

\vspace{0.1cm}

\begin{thm} {\em For every irrational angle $\theta$, every 
$\tau< 1$, and every neighborhood of the rotation $R_{\theta}$ 
in $\mathrm{Diff}_{+}^{1+\tau}(\mathrm{S}^1)$, there exists 
an element in this neighborhood with rotation number $\theta$ 
which is not topologically conjugate to $R_{\theta}$.}
\label{herman-contr}
\end{thm}

\vspace{0.1cm}

Before passing to the proof of Theorem \ref{denj}, we would like to give an 
``heuristic proof'' which applies 
in a particular but very illustrative case. Suppose that 
$f$ is a $\ce^1$ circle diffeomorphism with irrational rotation number and 
admitting a minimal invariant Cantor set $\Lambda$ so that the derivative 
of $f$ on $\Lambda$ is identically equal to $1$. Following \cite{har}, we will 
prove by contradiction that $f$ cannot be of class $\mathrm{C}^{1+\mathrm{Lip}}$   
(compare \cite{norton}; see also Exercise \ref{gugu}). For this, let us fix a connected 
component $I$ of the complementary set of $\Lambda$, and for each $n \in \mathbb{N}$ 
let us denote $I_n = g^n (I)$. Notice that $|I_{n+1}| = g'(p) |I_n|$ for some 
$p \in I_n$, and since the derivative of $g$ at the endpoints of $I_n$ equals 
$1$, we conclude that
$$\left| \frac{|I_{n+1}|}{|I_n|} - 1 \right| = |f'(p) - 1| \leq C |I_n|,$$
where $C$ is the Lipschitz constant of the derivative of $f$. We then have 
\index{Lipschitz!derivative}
\begin{equation}
|I_{n+1}| \geq |I_n| \esp (1 - C |I_n|).
\label{denjoy-heur}
\end{equation}
Without loss of generality, we may assume that the length of $I$ as well as that of all of 
its (positive) iterates by $f$ are less than or equal to $1 / 2C$. The contradiction 
we search then follows from the following elementary lemma by letting $d\!=\!1$ 
(the case $d > 1$ is naturally related to \S \ref{ejemplosceuno}4).

\vspace{0.1cm}

\begin{lem} {\em If \esp $d \!\in\! \mathbb{N}$ \esp and $(\ell_n)$ is a sequence 
of positive real numbers such that \esp $\ell_n\leq1 /C^{d} (1+1\!/\!d)^{d}$ \esp and 
\esp $\ell_{n+1} \geq \ell_n (1 - C \ell_n^{1/d})$ \esp for all $n \in \mathbb{N}$ and 
some positive constant $C$, then there exists $A>0$ such that \esp $\ell_n \geq A / n^d$ 
\esp for all $n \in \mathbb{N}$. In particular, if $d = 1$ then $S = \sum \ell_n$ 
diverges.}
\label{arg-norton}
\end{lem}

\noindent{\bf Proof.} The function $s \mapsto s (1 - Cs^{1\!/\!d})$ is increasing 
on the interval $\big[ 0, \big( \frac{1}{C(1+1\!/\!d)} \big)^{d} \big]$. Using 
this fact, the claim of the lemma easily follows by induction for the constant 
$A \!= \!\min \{\ell_1,d^d / 2^{d^2} C^d\}$. We leave the details to the 
reader. $\hfill \square$

\vspace{0.35cm}

We will give two distinct proofs of Denjoy's theorem (naturally, both of them will strongly 
use the combinatorial properties of circle homeomorphisms with irrational rotation number). 
The first one, due to Denjoy himself, consists in controlling the distortion produced by 
the diffeomorphism on the affine structure of the circle. As we will see later, this is 
also the main idea behind the proof of many other results in one-dimensional dynamics 
of group actions, as for instance those of Sacksteder and Duminy. In the second 
proof, we will concentrate on the distortion \index{control of distortion}
with respect to the projective structure of the circle. This 
technique has revealed useful in many other cases. For instance, 
it was used by Yoccoz 
\index{Yoccoz}
in \cite{Yo2} to extend Denjoy's theorem to $\mathrm{C}^{\infty}$ 
circle {\em homeomorphisms} whose critical points are not infinitely flat (in 
particular, for real-analytic circle homeomorphisms), and by Hu and Sullivan 
\index{Sullivan}
to obtain in \cite{HS} very fine results concerning the optimal regularity 
of Denjoy's theorem for diffeomorphisms. We will use a variation of this 
idea to give a general rigidity theorem in \S \ref{nose}.

\vspace{0.35cm}

\noindent{\bf First proof of Denjoy's theorem.} Assume for a contradiction that (the group 
generated by) a $\ce^{1+\mathrm{bv}}$ circle diffeomorphism $f$ admits an exceptional 
minimal set $\Lambda$. Let $I$ be a connected component of $\mathrm{S}^1 \setminus \Lambda$, 
let $x_0$ be an 
interior point of $I$, and let $\varphi$ be the semiconjugacy of $f$ to $R_{\theta}$,  
where $\theta = \rho(f)$. Without loss of generality, we may suppose that $\varphi(x_0)=0$. 
For $n \!\in\! \mathbb{N}$ let us consider the intervals $I_n$ and $J_n$ of the construction 
given at the end of \S \ref{teoria-poinc}, and let us define $I_n(f) = \varphi^{-1}(I_n)$ 
and $J_n(f) = \varphi^{-1}(J_n)$. Notice that $I$ is contained in each $I_n(f)$. 
From the combinatorial properties of $I_n$ and $J_n$ one concludes that:

\vspace{0.08cm}

\noindent{(i) the intervals in 
$\{ f^j(J_n(f))\!: j \!\in\! \{ 0,1,\ldots,q_{n+1}-1 \} \}$ 
cover the circle, and every point of $\clo$ is contained in at most two of them;}

\vspace{0.08cm}

\noindent{(ii) the interval $f^{q_{n+1}}(I_n(f))$ is contained in $J_n(f)$ for 
every $n \in \mathbb{N}$, and the same holds for $f^{-q_{n+1}}(I_{-n}(f))$.}

\vspace{0.08cm}

We claim that, for each $x,y$ in $I$ and each $n \in \mathbb{N}$, 
\begin{equation}
\big| (f^{q_{n+1}})'(x) \hspace{0.02cm} (f^{-q_{n+1}})'(y) \big| \geq e^{-2 V(f)}.
\label{denjoy}
\end{equation}
To show this, we let $\bar{y} = f^{-q_{n+1}}(y)$. Since $I$ is contained in 
$I_{-n}(f)$, by property (ii) we have $\bar{y} \in J_n(f)$. From 
(\ref{cociclolog}) and property (i) we obtain 
\begin{multline*}
\Big| \log\big(  (f^{q_{n+1}})'(x) \hspace{0.02cm} (f^{-q_{n+1}})'(y) \big) \Big| = 
\left| \log \bigg( \frac{(f^{q_{n+1}})'(x)}{(f^{q_{n+1}})'(\bar{y})} \bigg) \right| \leq \\
\leq \sum_{k=0}^{q_{n+1}-1} \big| \log(f')(f^k(x))-\log(f')(f^k(\bar{y})) \big| 
\leq \sum_{k=0}^{q_{n+1}-1} V \big( f;f^k(J_n(f)) \big) \leq 2 V(f),
\end{multline*}
from where one easily concludes inequality (\ref{denjoy}).

To close the proof notice that, for some $x,y$ in $I$,
$$\big| f^{q_{n+1}}(I) \big| \cdot \big| f^{-q_{n+1}}(I) \big| =  
(f^{q_{n+1}})'(x) \hspace{0.02cm} (f^{-q_{n+1}})'(y) \cdot |I|^2.$$
From inequality (\ref{denjoy}) one deduces that the right-hand expression  
is bounded from below by \esp $\exp\big( -2V(f) \big) \cdot |I|^2$. 
\esp However, this is absurd, since the intervals $f^k(I)$ 
are two-by-two disjoint for $k \in \mathbb{Z}$, and hence 
$\big| f^{q_{n+1}}(I) \big| \cdot \big| f^{-q_{n+1}}(I) \big|$ 
goes to zero as $n$ goes to infinity. $\hfill\square$

\vspace{0.35cm}

For the second proof we will use the following notation: given an interval 
$I=[a,b]$ and a diffeomorphism $f$ defined on $I$, let (compare (\ref{cochina}))
$$M(f;I) = \frac{|f(b)-f(a)|}{|b-a| \sqrt{f'(a)f'(b)}}.$$
The reader may easily check the relation (compare (\ref{cocicloschwarz}))
\begin{equation}
M(f \circ g;I) = M(g;I) \cdot M(f; g(I)).
\label{estunoyocoz}
\end{equation}
Moreover, from the existence of a point $c \in\! I$ such that 
\esp $|f(b) - f(a)| = f'(c) |b-a|$ \esp one concludes that 
\begin{equation}
\big| \log  M(f;I) \big| \leq
\frac{1}{2} \big[ \big| \log(f'(c)) - \log(f'(a)) \big|
\!+\! \big| \log(f'(c)) -\log(f'(b)) \big| \big] \leq \frac{V(f;I)}{2}.
\label{estdosyocoz}
\end{equation}

\vspace{0.25cm}

\noindent{\bf Second proof of Denjoy's theorem.} Suppose again that $f$ admits an 
exceptional minimal set $\Lambda$, and let us use the notation introduced at the 
beginning of the first proof. For each $n \in \mathbb{N}$ let us fix two points 
\hspace{0.01cm} $a \in f^{-(q_n + q_{n+1})}(I)$ \hspace{0.02cm} and  \hspace{0.02cm}
$b \in f^{-q_{n+1}}(I)$  \hspace{0.01cm} such that 
$$(f^{q_{n+1}})'(a) = \frac{\big| f^{-q_n}(I) \big|}{\big| f^{-(q_n+q_{n+1})}(I)\big|},  
\qquad \qquad (f^{q_{n+1}})'(b) = \frac{\big| I \big|}{\big|
f^{-q_{n+1}}(I) \big|}.$$
The open interval $J$ of endpoints $a$ and $b$ is contained in $J_n(f)$ (recall the 
combinatorial properties at the end of \S \ref{teoria-poinc}; in particular, see 
Figures 11 and 12). By (\ref{estunoyocoz}), (\ref{estdosyocoz}), and property (i) 
of the first proof, 
$$\big| \log M(f^{q_{n+1}};J) \big|
\leq \sum_{k=0}^{q_{n+1}-1} \big| \log M \big( f,f^k(J) \big) \big|
\leq \sum_{k=0}^{q_{n+1}-1} \frac{V \big( f;f^k(J) \big) }{2} \leq V(f).$$
We then obtain
$$\Big( \frac{\big| f^{q_{n+1}}(J) \big|}{ \big| J \big|} \Big)^2
\frac{1}{(f^{q_{n+1}})'(a) \hspace{0.04cm}
(f^{q_{n+1}})'(b)} \geq \exp \big( -2 V(f)\big),$$
that is,
$$\Big( \frac{\big| f^{q_{n+1}}(J) \big|}{ \big| J \big|} \Big)^2
\frac{\big| f^{-(q_n+q_{n+1})}(I) \big|}{\big| f^{-q_n}(I) \big|} \hspace{0.03cm}
\frac{\big| f^{-q_{n+1}}(I) \big|}{\big| I \big|} \geq \exp \big( -2 V(f) \big).$$
Since $\big| f^{q_{n+1}}(J) \big| \leq 1$ for all $n \! \in \! \mathbb{N}$ and  
$\big| f^{-(q_n+q_{n+1})}(I) \big|$ tends to zero as $n$ goes to infinity, this 
inequality implies that, for every $n$ large enough,
$$\big| f^{-q_{n+1}}(I) \big| \geq \big| f^{-q_n}(I) \big|.$$
Nevertheless, this contradicts the fact that $\big| f^{-q_n}(I) \big|$ 
converges to zero as $n$ goes to infinity. $\hfill\square$

\vspace{0.15cm}

\begin{small} \begin{ejer} 
The proof proposed below for Denjoy's theorem in class $\mathrm{C}^{1+\mathrm{Lip}}$ is 
quite interesting by the fact that it does not use any control of distortion type argument.

\vspace{0.08cm}

\noindent (i) Suppose that $f$ is a $\mathrm{C}^{1+\mathrm{Lip}}$ counterexample 
to Denjoy's theorem, and denote by $I$ one of the connected components 
of the complement of the exceptional minimal set. For each $j \in \mathbb{Z}$ 
choose a point $x_j \in f^j(I)$ such that \esp $| f^{j+1}(I) | = f'(x_j) \esp | f^j(I) |$. 
\esp Denoting by $C$ the Lipschitz constant of $f'$ and fixing $n \!\in\!\mathbb{N}$, 
from the inequality $| f'(x_j) - f'(x_{j-q_n}) | \leq C |x_j - x_{j-q_n}|$ conclude 
that, for some constant $\bar{C} > 0$ which is independent from $n$, and for  
every $j \in \{0,\ldots,q_{n}-1\}$,  
\begin{equation}
\frac{|f^{j+1}(I)| \cdot |f^{-q_n+j}(I)|}{|f^j(I)| \cdot |f^{-q_n+j+1}(I)|} 
\geq 1 - \bar{C} |x_j - x_{j-q_n}|.
\label{noseque}
\end{equation}
 
\vspace{0.08cm}

\noindent (ii) Using the combinatorial properties in \S \ref{teoria-poinc}, from the 
above inequality conclude the existence of a constant $D > 0$ such that, for every 
$n \in \mathbb{N}$,
\begin{equation}
\frac{|f^{q_n}(I)| \cdot |f^{-q_n}(I)|}{|I|^2} \geq D e^{-D}.
\label{nosequedos}
\end{equation}

\noindent{\underbar{Hint.}} The left-side expression coincides with 
$$\prod_{j=0}^{q_{n}-1} \frac{|f^{j+1}(I)| \cdot |f^{-q_n+j}(I)|}{|f^j(I)| 
\cdot |f^{-q_n+j+1}(I)|}.$$
Use inequality (\ref{noseque}) to estimate each factor of this product for $j$ large 
enough, and then use the elementary inequality $1 - x \geq e^{-2x}$ (which holds for 
$x > 0$ very small) as well as the fact that the intervals of endpoints $x_j$ and 
$x_{j-q_n}$, with $j \in \{0,\ldots,q_n-1\}$, are two-by-two disjoint.

\vspace{0.08cm}

\noindent (iii) Show that inequality (\ref{nosequedos}) cannot hold for all $n$.
\label{gugu}
\end{ejer} \end{small}

\vspace{0.1cm}

To construct smooth circle homeomorphisms admitting an exceptional minimal 
set, one may proceed in two distinct ways. The first consists in allowing 
critical points. In this direction, in \cite{Ha} one may find examples of 
$\mathrm{C}^{\infty}$ circle homeomorphisms with a single critical point and a 
minimal invariant Cantor set (by Yoccoz' theorem already mentioned \cite{Yo2}, 
such a homeomorphism cannot be real-analytic). The second issue consists in 
considering diffeomorphisms of differentiability class less than $\ce^2$. 
We will carefully study examples of this type in \S \ref{ejemplosceuno}.

Concerning the nature of the invariant Cantor set which may appear in differentiablity 
class less than $\ce^2$, let us first notice that any Cantor set contained in $\clo$ 
may be exhibited as the exceptional minimal set of a circle {\em homeomorphism}. 
However, if there is some regularity for the homeomorphism, then such a set 
must satisfy certain properties. For example, in \cite{McD} it is proved 
that the triadic Cantor set cannot appear as the exceptional minimal set 
of  any $\mathrm{C}^{1}$ circle diffeomorphism (see also \cite{norton-afin}). 
On the other hand, in \cite{NT} it is shown that every ``affine'' Cantor set 
may appear as the exceptional minimal set of a Lipschitz circle homeomorphism.
\index{Lipschitz!homeomorphism}

\vspace{0.1cm}

From the viewpoint of group actions, Denjoy's theorem may be reformulated 
by saying that there is no action of $(\mathbb{Z},+)$ by $\mathrm{C}^{1+\mathrm{bv}}$ 
circle diffeomorphisms which admits an exceptional minimal set. As we will see 
below, this property is shared by a large family of finitely generated groups. 

\vspace{0.15cm}

\begin{prop} {\em Let $\Gamma$ be a finitely generated group of 
$\ce^{1+\mathrm{Lip}}$ circle diffeomorphisms. If $\Gamma$ admits an 
exceptional minimal set, then it contains a free subgroup on two generators.}
\label{librecantor}
\end{prop}

\vspace{-0.35cm}

\noindent {\bf Proof.} We will show that $\Gamma$ acts expansively on the topological 
circle $\mathrm{S}^1_{\Lambda}$ associated to the exceptional minimal set $\Lambda$. By 
the proof of Theorem \ref{margulis}, 
\index{Margulis!weak Tits alternative}
this implies that $\Gamma$ contains a free group.
\index{group!free}

Suppose that the (minimal) action of $\Gamma$ on 
$\clo_{\Lambda}$ is not expansive. Then it is equicontinuous, 
\index{action!equicontinuous} and \index{action!expansive}
there 
\index{equicontinuous!action|see{action, equicontinuous}} 
is an invariant probability measure on 
$\clo_{\Lambda}$ with total support and no atoms
which allows conjugating $\Gamma$ to 
a group of rotations of $\mathrm{S}^1_{\Lambda}$. Let $g_1, \ldots, g_n$ be the 
generators of $\Gamma$. If each $g_i$ has finite order then every orbit in 
$\mathrm{S}^1_{\Lambda}$ is finite, which contradicts the minimality of 
$\Lambda$. On the other hand, if some generator $g_i$ has infinite order, 
then the rotation number of $g_i$ is irrational. Since $g_i$ has 
non-dense orbits (for instance, those contained in $\Lambda$), 
this contradicts Denjoy's theorem. $\hfill\square$

\vspace{0.3cm}

The proposition above applies to finitely generated, amenable groups. (Indeed, 
such a group cannot contain $\mathbb{F}_2$: see Appendix B). However, this 
may be shown in a more direct way ({\em i.e.,} without using Margulis' theorem).

\begin{small} \begin{ejer} Let $\Gamma$ be a finitely generated \index{group!amenable} 
amenable group of $\mathrm{C}^{1+\mathrm{bv}}$ circle diffeomorphisms. Suppose that 
$\Gamma$ admits an exceptional minimal set $\Lambda$, and consider an invariant 
probability measure $\mu$ for the action of $\Gamma$ on the topological circle 
$\clo_{\Lambda}$ associated to $\Lambda$. Prove that $\mu$ has total support 
and no atoms. Conclude that $\Gamma$ is semiconjugate to a group of rotations, 
and obtain a contradiction using Denjoy's theorem.
\index{rotation!group}
\label{a-referenciar}
\end{ejer} \end{small}

Proposition \ref{librecantor} still holds for non finitely generated groups  
provided that there exists a constant $C$ such that $V(g) \leq C$ for every 
generator $g$. Without such a hypothesis, it fails to hold for general non 
finitely generated subgroups of $\vl$, as is shown by the following example, 
essentially due to Hirsh \cite{Hir} (see also \cite{imanishi}).

\begin{small} \begin{ejem} Let $H\!: \mathbb{R} \!\rightarrow\! \mathbb{R}$ be 
a homeomorphism satisfying the properties (i) and (ii) of \S \ref{ghys-sergi} 
\index{Ghys!Ghys-Sergiescu's realization of Thompson's groups} 
and such that $\mathrm{Fix} (H) \!=\! \{ a, b \}$, with $a=0$ and $b\!>\!0$ small. 
This homeomorphism $H$ naturally induces a degree-$2$ circle map $\bar{H}$ (see 
Figure 10). Let us denote by $I$ the interval $]a,b[$. Let us now define 
$g_1 \!: \clo \rightarrow \clo$ by $g_1(x) = y$ if $\bar{H}(x) = \bar{H}(y)$ 
and $x \neq y$. In general, for each $n \in \mathbb{N}$ let us consider the 
degree-$2^n$ circle map $\bar{H}^{n}$, and let us define $g_n \!: \clo \rightarrow \clo$ 
by $g_n(x) = y$ if $\bar{H}^{n}(x) = \bar{H}^{n}(y)$ and $\bar{H}^{n}(x) \neq \bar{H}^{n}(y')$ 
for every $y' \! \in ]x,y[$.

The group $\Gamma \!=\! \Gamma_H$ generated by the $g_n$'s is Abelian and isomorphic to 
$\mathbb{Q}_2 (\clo)$. Its action on the circle is semiconjugate to that of a group of 
rotations, and its exceptional minimal set is 
$\Lambda \!=\! \clo \setminus \cup_{n \in \mathbb{N}} 
\cup_{i=0}^{2^n -1} g_n^i(C)$. Notice that if $H$ is a real-analytic diffeomorphism of 
the real line, then the elements in $\Gamma_H$ are real-analytic circle diffeomorphisms. 
The reader should notice the similarity between this construction and the last one of 
\S \ref{sec-minimal}. Indeed, the group $\Gamma_H$ is a subgroup of (the realization 
$\Phi_H (\mathrm{G})$) of Thompson's group G \index{Thompson's groups!G}
(though $\Phi_H(\mathrm{G})$ is not a group 
of real-analytic diffeomorphisms of $\clo$\esp!).
\label{ejemplohirsh}
\end{ejem} \end{small}

Denjoy's theorem may be also seen as a consequence of the next proposition, 
which is the basis for dealing with the relevant problem of the regularity 
of the {\bf \em{linearization}} \index{linearization} of smooth circle 
diffeomorphisms ({\em i.e.,} the conjugacy to the corresponding rotation).

\vspace{0.1cm}

\begin{prop} {\em Let $\theta \!\! \in ]0,1[$ be an irrational number, and let $p / q$ 
be one of the rational approximations of $\theta$ (i.e., $|\theta - p/q| \leq 1/q^2$). 
For a circle homeomorphism $f$ of rotation number $\theta$, 
\index{variation of a function!total variation}
denote by $\mu$ the unique probability measure on $\clo$ which is invariant by $f$. 
If $\psi: \clo \rightarrow \mathbb{R}$ is a (non-necessarily continuous) function 
with finite total variation, then for every $x \in \clo$ one has}
\begin{equation}
\left| \sum_{i=0}^{q-1} \psi \big( f^i(x)\big) - q \int_{\clo} \psi \esp d\mu \right|
\leq \mathrm{var} (\psi;\clo).
\label{denj-k}
\end{equation}
{\em In particular, if $f$ is differentiable and its derivative has bounded 
variation, then for every $x \in \clo$ and every $n \in \mathbb{N}$ one has}
\begin{equation}
\exp \big( -V(f) \big) \leq (f^{q_n})'(x) \leq \exp \big( V(f) \big).
\label{des-denjoy}
\end{equation}
\label{pden}
\end{prop}

\vspace{-0.6cm}

To prove this proposition, we will use the following combinatorial lemma.

\vspace{0.075cm} 

\begin{lem} {\em If \esp $\theta \! \in ]0,1[$ \esp is irrational and \esp $p,q$ \esp are 
two relatively prime integers such that $q \neq 0$ and $\big| \theta - p/q \big| \leq 1/q^2$, 
then for each $i \!\in\! \{0,\ldots,q-1\}$ the interval \esp $]i/q,(i+1)/q[$ \esp contains a 
unique point of the set 
$\big\{ j\theta \esp (\mathrm{mod } \esp 1)\!: j \in \{1,\ldots,q\} \big\}$.}
\end{lem}

\noindent{\bf Proof.} By the hypothesis, we have 
either \esp $0 < \theta - p/q < 1 / q^2$  \esp 
or \esp $- 1 / q^2 < \theta - p / q < 0$. \esp Both 
cases being analogous, let us consider only the first 
one. For every $j \! \in \! \{1,\ldots,q\}$ one has
$$0 < j \theta - \frac{jp}{q} < \frac{j}{q^2} \leq \frac{1}{q},$$
and hence \esp $j \theta \esp \mathrm{mod } \esp 1$ \esp belongs to  
$]jp/q,(jp+1)/q[$. The claim of the lemma follows from the fact that these 
intervals \esp $[jp/q,(jp+1)/q[$ \esp cover the circle without intersection. 
$\hfill\square$

\vspace{0.35cm}

\noindent{\bf Proof of Proposition \ref{pden}.} It is evident that the claim to 
be proved is equivalent to that the inequality 
$$\left| \sum_{i=1}^{q} \psi \big( f^i(x)\big) - q \int_{\clo} \psi \esp d\mu \right|
\leq \mathrm{var} (\psi;\clo)$$
holds for all $x\!\in\!\clo$ and every function of bounded variation $\psi$. 
To show this inequality, 
we consider the map from the circle into itself defined by
$$\varphi(y) = \int_x^y d\mu \quad \mathrm{mod } \esp 1.$$
From the relation $\varphi \circ f =  R_{\theta} \circ \varphi$ 
we deduce that, if we let $x_j \!=\! f^j (x)$, then we have 
$\varphi(x_j) \!=\! j\theta \esp \mathrm{mod} \esp 1$. 
By the preceding lemma, if for each $j \!\in\! \{1,\ldots,q\}$ we choose an interval 
$I_j \!= ]i_j/q,(i_j+1)/q[$ containing $\varphi(x_j)$, then the intervals which 
appear are two-by-two disjoint. Therefore, denoting by $J_j$ the interval 
$\varphi^{-1} (\bar{I}_j)$ (whose $\mu$-measure equals $1/q$), from the equality 
$$\left| \sum_{i=1}^{q} \psi \big( f^i(x)\big) - q \int_{\clo} \psi \esp d\mu \right|
= \left| \sum_{j=1}^q \Big( \psi\big(f^j(x)\big) - q \int_{J_j} \psi \esp d\mu
\Big)\right|$$
one concludes that 
\begin{small}
$$\left| \sum_{i=1}^{q} \psi \big( f^i(x)\big) \!-\! q \int_{\clo} \!\psi \esp d\mu \right|
\leq q \sum_{j=1}^q \!\int_{J_j} \!\!\big| \psi\big(f^j(x)\big) - \psi \big| \esp d\mu 
\leq \sum_{j=1}^q \sup_{y \in J_j} \big| \psi\big(f^j(x)\big) - \psi(y) \big| 
\!\leq\! \mathrm{var} (\psi,\clo).$$
\end{small}

The second part of the proposition follows from the first one applied to 
the function $\psi(x) = \log \big( f'(x) \big)$, thanks to the equality
\begin{equation}
\int_{\clo} \log (f') \esp d \mu = 0.
\label{intenula}
\end{equation}
However, we point out that this equality is not at all evident. 
For the proof, notice that from (\ref{denj-k}) one deduces 
that, for every $x \in \clo \!$,
$$\exp \big( -V(f) \big) \leq \frac{(f^{q_n})'(x)}{\exp\big( q_n \int_{\clo}
\log(f') d\mu \big)}
\leq \exp \big( V(f) \big).$$
Integrating with respect to Lebesgue measure, this allows to conclude that
$$\exp \big( -V(f) \big) \leq 
\frac{1}{\exp \big( q_n \int_{\clo} \log(f') d\mu \big)}
\leq \exp \big( V(f) \big),$$
which cannot hold for every $n \in \mathbb{N}$ 
unless (\ref{intenula}) holds (compare Exercise \ref{inte-cero}). $\hfill\square$

\vspace{0.3cm}

Relation (\ref{denj-k}) (resp. (\ref{des-denjoy})) is known as the {\bf {\em Denjoy-Koksma 
inequality}} (resp. {\bf \textit{Denjoy inequality}}). \index{Denjoy!inequality}
Let us point out that, similarly to Denjoy's theorem and equality (\ref{intenula}), 
these inequalities still hold for piecewise affine circle homeomorphisms, and 
the argument of the proofs are similar to those of the $\ce^{1+\mathrm{bv}}$ 
case (of course, in the piecewise affine case, one consider a lateral 
derivative instead of the usual one). \index{piecewise affine homeomorphism}

\begin{small}
\begin{obs} For every $\ce^{2+\tau}$ circle diffeomorphism $f$ of irrational 
rotation number, a stronger conclusion than that 
of Proposition \ref{pden} holds: The sequence 
$(f^{q_n})$ converges to the identity in the $\ce^1$ topology, where $p_n/q_n$ 
denotes the $n^{\mathrm{th}}$-rational approximation of $\rho(f)$ (see \cite{He,KO,Yo1}). 
\label{des-denj-conv}
\end{obs}

\begin{ejer} Prove that equality (\ref{intenula}) holds for every $\ce^1$ circle 
diffeomorphism by using the Birkhoff \index{Birkhoff!Ergodic Theorem}
Ergodic Theorem and the unique ergodicity of $f$ (see \cite{walters}).
\label{inte-cero}
\end{ejer} \end{small}

Notice that Denjoy's theorem implies that if $\Gamma$ is a subgroup of $\vl$ admitting 
an exceptional minimal set, then each of its elements has periodic points. Indeed, if 
$g \in \Gamma$ has no periodic point then its rotation number is irrational, which 
implies the density of the orbits by $g$ (and hence the density 
of the orbits by $\Gamma$). An  
interesting consequence of this fact is the rationality of the rotation number of 
the elements in Thompson's group G. Indeed, in \S \ref{sec-minimal} we constructed an 
action of G by $\ce^{\infty}$ circle diffeomorphisms which is semiconjugate to the 
standard action and admits an exceptional minimal set. The fact that $\rho(g)$ 
belongs to $\mathbb{Q}$ for each $g \in \mathrm{G}$ then follows from the above 
remark together with the invariance of the rotation number under semiconjugacy. 
Following \cite{liouse}, we invite the reader to develop an alternative and 
``more direct'' proof just using Denjoy's inequality (alternative proofs 
may be found in \cite{cal-thompson} and \cite{klep}).
\index{Thompson's groups!G}\index{Denjoy!theorem|)}

\begin{small} \begin{ejer} Given $g \in \mathrm{G}$ let us define $M_n \in \mathbb{N}$
and $N_n \in \mathbb{N} \cup \{0\}$ so that for each $n \in \mathbb{N}$ one has 
$g^n(0) = M_n / 2^{N_n}$, where $M_n$ is either odd or zero. Suppose that $0$ is 
not a periodic point for $g$, and thus $M_n \neq 0$ for every $n \in \mathbb{N}$.

\noindent (i) Prove that $N_n$ tends to infinity together with $n$.\\

\noindent{\underbar{Hint.}} Notice that for every $N \in \mathbb{N}$ the set of 
dyadic rational numbers with denominator less than or equal to $2^N$ is finite.
\index{dyadic!rational number}

\noindent (ii) Using the convergence of $N_n$ to infinity, conclude that 
$\lim_{n \rightarrow \infty} (g^n)'(0) \!=\! 0$ (where we consider 
the right derivative of the map).

\noindent{\underbar{Hint.}} Write $g(x) = 2^{\lambda(x)} + M(x) / 2^{N(x)}$ for 
some integer-valued uniformly bounded functions $\lambda$, $M$, and $N$, and 
show that for $n \!\in\! \mathbb{N}$ large enough one has 
$$M_{n+1} = M_n + 2^{N_n - \lambda(g^n(0)) - N(g^n(0))} M(g^n(0)), 
\qquad N_{n+1} = N_n - \lambda(g^n(0)).$$

\noindent (iii) Show that $\rho(g)$ is rational by applying Denjoy's inequality.
\end{ejer} \end{small}


\section{Sacksteder's theorem\index{Sacksteder!theorem|(}}
\label{psg}

\hspace{0.45cm} Let us begin by recalling the notion of 
\textit{pseudo-group} \index{pseudo-group} of homeomorphisms.

\begin{defn} A family $\Gamma = \{ g\!: dom(g) \rightarrow ran(g) \}$ of 
homeomorphisms between subsets of a topological space $\mathrm{M}$ is a 
{\bf \textit{pseudo-group}} if the properties below are satisfied:

\vspace{0.05cm}

\noindent{-- the domain $dom(g)$ and the image $ran(g)$ of each $g \in \Gamma$ 
are open sets;}

\vspace{0.05cm}

\noindent{-- if $g,h$ belong to $\Gamma$ and $ran(h) \!\subset\! dom(g)$, 
then $gh$ belongs to $\Gamma$;}

\vspace{0.05cm}

\noindent{-- if $g \in \Gamma$ then $g^{-1} \in \Gamma$;}

\vspace{0.05cm}

\noindent{-- the identity on $\mathrm{M}$ is an element of $\Gamma$;}

\vspace{0.05cm}

\noindent{-- if $g$ belongs to $\Gamma$ and $A \subset dom(g)$ is an open 
subset, then the restriction $g|_{A}$ of $g$ to $A$ is an element of $\Gamma$;}

\vspace{0.05cm}

\noindent{-- if $g : dom(g) \rightarrow ran(g)$ is a homeomorphism between 
open subsets of 
$\mathrm{M}$ such that for every $x \in dom(g)$ there exists a neighborhood 
$V_x$ such that $h|_{V_x}$ is an element of $\Gamma$, then $g$ belongs to 
$\Gamma$.}
\label{def-pseudo}
\end{defn}

The notions of orbit and invariant set for a pseudo-group are naturally defined. 
We will say that $\Gamma$ is {\bf \textit{finitely generated}} (resp. 
{\bf \textit{countably generated}}) 
if there exists a finite (resp. countable) subset 
$\mathcal{G}$ of $\Gamma$ such that every element in $\Gamma$ may be written 
as (the restriction to its domain of) a product of elements in $\mathcal{G}$. 
A measure $\mu$ on the Borel sets of $\mathrm{M}$ is invariant by $\Gamma$ 
if for every Borel set $A$ and every $g \!\in\! \Gamma$ one has 
$\mu \big( A \cap dom(g) \big) = \mu \big( g(A) \cap ran(g) \big)$. 
Finally, for $a \in \mathrm{M}$ we denote by $\Gamma(a)$ the orbit of $a$, 
and for each $p \in \Gamma(a)$ we define its {\bf \textit{order}} by letting 
$$\mathrm{ord}(p) = \min \{n \in \mathbb{N}:
\mbox{ there exist } g_{i_{j}} \in \mathcal{G} \mbox{ such that }
g_{i_n} \cdots  g_{i_1}(a) = p \}.$$

Pseudo-groups naturally appear in Foliation Theory: 
\index{foliation} after prescribing a set of open transversals, 
the holonomies along leaves form a pseudo-group under composition.
 But as we will see, even in the context of group actions, it is sometimes important to 
keep in mind the pseudo-group approach, since many results are proved by restricting 
the original action to certain open sets where the group structure is loosed. 


\subsection{The classical version in class $\mathbf{\ce^{1+\mathrm{Lip}}}$}
\label{chuartz}

\hspace{0.45cm} In this section, we will prove the most classical version of 
a theorem due to Sacksteder \cite{Sa} about the existence of elements with 
{\bf \em{hyperbolic fixed points}} in some pseudo-groups of diffeomorphisms 
of one-dimensional compact manifolds. 
\index{hyperbolic!fixed point}

\vspace{0.1cm}

\begin{thm} {\em Let $\Gamma$ be a pseudo-group of $\mathrm{C}^{1+\mathrm{Lip}}$ 
diffeomorphisms of a one-dimensional compact manifold generated by a finite family 
of elements $\mathcal{G}$. Suppose that for each generator $g \!\in\! \mathcal{G}$ 
there exists $\bar{g} \!\in\! \Gamma$ whose domain contains the closure of $dom(g)$ 
and that coincides with $g$ when restricted to $dom(g)$.\footnote{This is the 
so-called {\bf{\em compact generation property}} for the pseudo-group, whose 
relevance was cleverly noticed by Haefliger \cite{haefliger} (see Remark 
\ref{remark-haefliger} on this). It is always satisfied by the holonomy 
pseudo-group of a codimension-one foliation on a compact manifold.
\label{servira}} Suppose moreover that 
there exists an invariant closed set $\Lambda$ such that for some connected 
component $I \!\!= ]a,b[$ of its complement, either $a$ is an accumulation 
point of the orbit $\Gamma(a)$, or $b$ is an accumulation point of $\Gamma(b)$. 
Then there exist $p \in \Lambda$ and $h \in \Gamma$ such that 
$h(p) = p$ and $h'(p) < 1.$}
\label{sacksteder}
\end{thm}

\vspace{0.1cm}
\index{Haefliger}\index{pseudo-group!compactly generated}

A set $\Lambda$ satisfying the above hypotheses (with the only possible exception 
of the one concerning the regularity) is called a {\bf {\em local exceptional 
set}}.\index{local exceptional set} \index{Schwartz estimates}
To show the theorem, the main technical tool will be the next lemma, whose 
$\ce^{1+\mathrm{Lip}}$ version is essentially due to Schwartz \cite{Sc}. We 
give a general version including the 
$\ce^{1+\tau}$ case for future reference; 
a slightly refined $\ce^{1+\mathrm{Lip}}$ version will be discussed in 
\index{Duminy!estimates} \S \ref{duminy-fines2} (see also Exercise \ref{mucho}). 

\vspace{0.1cm}

\begin{lem} {\em Let $\Gamma$ be a pseudo-group of diffeomorphisms of class 
$\ce^{1+\tau}$ (resp. $\ce^{1+\mathrm{Lip}}$) of a one-dimensional, compact manifold. 
Suppose that there exist a finite subset $\mathcal{G}$ of $\Gamma$, a constant 
$S \!\in\! [1,\infty[$, and an open interval $I$, such that to each $g \in \mathcal{G}$ 
one may associate a compact interval $L_g$ contained in an open set where $g$ is defined  
in such a way that, for each $m \in \mathbb{N}$, there exists $g_{i_{m}} \!\in \mathcal{G}$ 
so that for all $n \!\in\! \mathbb{N}$ the element $h_n = g_{i_{n}} \cdots g_{i_{1}} 
\in \Gamma$ satisfies the following properties:

\vspace{0.1cm}

\noindent{-- \thinspace the interval $g_{i_{k-1}} \!\cdots g_{i_{1}}(I)$ is contained 
in $L_{g_{i_k}}$ \esp (where \esp $g_{i_{k-1}} \cdots g_{i_{1}} = Id$ for $k = 1$),}

\vspace{0.1cm}

\noindent{-- \thinspace one has the inequality} \esp \esp \esp 
$\sum_{k = 0}^{n-1} |g_{i_{k}} \cdots g_{i_{1}}(I)|^{\tau} \leq S \quad
\big( \mbox{resp. } \sum_{k = 0}^{n-1} |g_{i_{k}} \cdots g_{i_{1}}(I)| \leq S \big).$ 

\vspace{0.1cm}

\noindent Then there exists a positive constant $\ell = \ell(\tau,S,|I|;\mathcal{G})$ 
such that, if for some $n \!\in\! \mathbb{N}$ the interval $h_n(I)$ is contained in a 
$\ell$-neighborhood of $I$ and does not intersect the interior of $I$, then the 
map $h_n$ has a hyperbolic fixed point.}
\label{?}
\end{lem}

\noindent{\bf Proof.} To simplify, we will think of the case $\tau \!=\! 1$ as the 
Lipschitz case in order to treat it simultaneously 
with the H\"older case $\tau\!  \in ]0,1[$.  
\index{H\"older!derivative} \index{Lipschitz!derivative}
Let $\varepsilon > 0$ be a constant 
such that each $g \in \mathcal{G}$ is defined in a $2\varepsilon$-neighborhood of $L_g$. 
We fix $C \!>\! 0$ in such a way that, for all $g \in \mathcal{G}$ and all $x,y$ in the 
$\varepsilon$-neighborhood of $L_g$, one has
$$\big| \log(g'(x)) - \log(g'(y)) \big| \esp \leq \esp C \esp |x-y|^{\tau}.$$
We will show that the claim of the lemma is satisfied for
$$\ell \esp = \esp \min \left\{ \frac{|I|}{2 \exp(2^{\tau} C S )},
\frac{|I| \esp \varepsilon}{2 \exp(2^{\tau} C S )}, 
\frac{\varepsilon}{2} \right\}.$$

We denote by $J$ the $2\ell$-neighborhood of $I$, and we let  $I'$ (resp. $I''$) be the 
connected component of $J \setminus I$ to the right (resp. to the left) of $I$. We will 
show by induction that the following properties are simultaneously satisfied:

\vspace{0.1cm}

\noindent{$\mathrm{(i)}_{k}$ \hspace{0.25cm} $g_{i_{k}} \cdots g_{i_{1}}(I')$ 
is contained in the $\varepsilon$-neighborhood of $L_g$;}

\vspace{0.1cm}

\noindent{$\mathrm{(ii)}_{k}$ \hspace{0.25cm} $|g_{i_{k}} \cdots g_{i_{1}}(I')| \leq
|g_{i_{k}} \cdots g_{i_{1}}(I)|$;}

\vspace{0.1cm}

\noindent{$\mathrm{(iii)}_{k}$ \hspace{0.25cm} $\sup_{\{x,y\} \subset I \cup I'}
\frac{(g_{i_{k}} \cdots g_{i_{1}})'(x)}{(g_{i_{k}} \cdots g_{i_{1}})'(y)}
\leq \exp(2^{\tau} \thinspace C S )$.}

\vspace{0.1cm}

Condition (iii)$_0$ is obviously satisfied, while (i)$_0$ and (ii)$_0$ follow from the 
hypothesis $|I'| \!=\! 2\ell$ and the inequalities $\ell \!\leq\! \varepsilon / 2$ 
and  $\ell \!\leq\! |I| / 2$. Suppose that (i)$_j$, (ii)$_j$, and (iii)$_j$, are 
true for every $j \!\in\! \{0,\ldots,k-1\}$. In this case, for every $x,y$ in 
$I \cup I'$ we have
\begin{eqnarray*}
\left|\log\left(\frac{(g_{i_{k}} \!\cdots g_{i_{1}})'(x)}
{(g_{i_{k}} \!\cdots g_{i_{1}})'(y)}\right)\right|
\!\! &\leq& \!\! 
\sum_{j=0}^{k-1} \big| \log(g_{i_{j+1}}'(g_{i_{j}} \!\cdots g_{i_{1}}(x))) 
- \log(g_{i_{j+1}}'(g_{i_{j}} \!\cdots g_{i_{1}}(y))) \big|\\
\!\! &\leq& \!\!  
C \sum_{j=0}^{k-1} \big| g_{i_{j}} \cdots g_{i_{1}}(x) -
g_{i_{j,n}} \cdots g_{i_{1}}(y) \big|^{\tau}\\
\!\! &\leq& \!\! 
C \sum_{j=0}^{k-1} \big( |g_{i_{j}} \cdots g_{i_{1}}(I)| +
|g_{i_{j}} \cdots g_{i_{1}}(I')| \big)^{\tau}\\
\!\! &\leq& \!\! C \thinspace 2^{\tau} S.
\end{eqnarray*}
This inequality shows (iii)$_k$. Concerning (i)$_k$ and 
(ii)$_k$, notice that there exist $x \in I$ and $y \in I'$ such that
$$|g_{i_{k}} \cdots g_{i_{1}} (I)| = |I| \cdot (g_{i_{k}} \cdots g_{i_{1}})'(x),
\hspace{0.475cm} \hspace{0.475cm} |g_{i_{k}} \cdots g_{i_{1}}(I')|
= |I'| \cdot (g_{i_{k}} \cdots g_{i_{1}})'(y).$$
Hence, by (iii)$_k$,
$$\frac{|g_{i_{k}} \cdots g_{i_{1}}(I')|}{|g_{i_{k}} \cdots g_{i_{1}}(I)|}=
\frac{(g_{i_{k}} \cdots g_{i_{1}})'(x)}{(g_{i_{k}} \cdots g_{i_{1}})'(y)} 
\cdot \frac{|I'|}{|I|} \leq \exp(2^{\tau} C S) \frac{|I'|}{|I|},$$
which shows (i)$_k$ and (ii)$_k$ by the definition of $\ell$. Of 
course, similar properties also hold for $I''$.

Now suppose that $h_n(I)$ is contained in the $\ell$-neighborhood of the 
interval $I$ and does not intersect  the interior of $I$. Property (ii)$_n$ 
then gives $h_n(J) \subset J$. Moreover, if $h_n(I)$ is to the right 
(resp. to the left) of $I$, then $h_n(I \cup I') \subset I'$ (resp.
$h_n(I'' \cup I) \subset I''$). Both cases being analogous, we will only 
deal with the first one. Clearly, $h_n$ has a fixed point $x$ in $I'$, 
and we need to verify that this fixed point is hyperbolic (and 
contracting). To do this, chose $y \in I$ such that \esp 
$h_n'(y) = |h_n(I)| / |I| \leq \ell / |I|$. \esp 
By (iii)$_n$, 
$$h_n'(x) \leq h_n'(y) \exp(2^{\tau}CS) \leq 
\frac{\ell \exp(2^{\tau}CS)}{|I|} \leq \frac{1}{2},$$
and this finishes the proof. $\hfill\square$

\vspace{0.35cm}

An inequality of type $\sum_{k \geq 0} |g_{i_k} \cdots g_{i_1}(I)| \!\leq\! S$ obviously 
holds when the intervals $g_{i_k} \cdots g_{i_1}(I)$ are two-by-two disjoint. In 
this case, the issue of the preceding lemma lies on the possibility of controlling 
the distortion 
\index{control of distortion}
for the compositions in the $2\ell$-neighborhood of $I$, 
despite the fact that the images of this neighborhood are no longer disjoint.

\vspace{0.4cm}

\noindent{\bf Proof of Theorem \ref{sacksteder}.} Suppose that $b$ is an 
accumulation point of $\Gamma(b)$ (the case where $a$ is an accumulation point 
of $\Gamma(a)$ is analogous). Let $(h_n) \!=\! (g_{i_n} \cdots g_{i_1})$ be a 
sequence in $\Gamma$ such that $h_n(b)$ converges to $b$, such that $h_n(b) \neq b$ 
for every $n \!\in\! \mathbb{N}$, and such that the order of the point $h_n(b) \in \Gamma(b)$
is realized by $h_n$ (we leave to the reader the task of showing the existence of such a 
sequence). Clearly, the hypotheses of the preceding lemma are satisfied for the connected 
component $I$ of the complement of $\Lambda$ containing $b$  and every constant 
$S\! \geq \!1$ greater than the total length of the underlying one-dimensional manifold. 
We then conclude that, for $n$ large enough, the map $h_n \in \Gamma$ contracts 
(hyperbolically) into itself the interval $I'$ of length $\ell$ situated immediately to 
the right of $I$. Therefore, $h_n$ has an unique fixed point in $I$ (namely the point 
$p = \bigcap_{k \in \mathbb{N}} (h_n)^{k}(I')$), which is hyperbolic (see Figure 15). 
Remark finally that, since $\Lambda$ is a compact invariant set and $I'$ intersects 
$\Lambda$, the point $p$ must belong to it. $\hfill\square$

\vspace{0.54cm}

\beginpicture

\setcoordinatesystem units <0.98cm,0.98cm>

\putrule from -5.5 0 to 5.5 0

\putrule from -2.8 -0.015 to 2.8 -0.015

\putrule from -2.8 0.015 to 2.8 0.015

\putrule from 3 -0.015 to 3.5 -0.015

\putrule from 3 0.015 to 3.5 0.015

\circulararc -45 degrees from 0 0.5 center at 1.6 -3.3

\circulararc 132 degrees from 3.8 -0.5 center at 4.3 -0.3

\circulararc -55 degrees from 2.2 -0.5 center at -1.25 6

\plot 3.8 -0.5
3.9 -0.56 /

\plot 3.8 -0.5
3.8 -0.615 /

\plot 2.2 -0.5
2.12 -0.63 /

\plot 2.2 -0.5
2.03 -0.51 /

\plot 3.14 0.52
3 0.52 /

\plot 3.14 0.52
3.05 0.62 /

\put{$\Big|$} at -4.8 0
\put{$\Big|$} at 4.8 0
\put{$|$} at 3.722 0
\put{$|$} at 2.5 0
\put{$\Big($} at -2.8 0
\put{$\Big)$} at 2.8 0
\put{$($} at 3 0
\put{$)$} at 3.5 0
\put{$h_n$} at -3.8 -1.2
\put{$h_n$} at 4.4 -1.2
\put{$h_n$} at 1.3 1.1

\plot 3.615 1
3.615 0.3 /

\plot 3.615 0.3
3.55 0.45 /

\plot 3.615 0.3
3.68 0.45 /

\put{$I \!=\! [a,b]$} at -1.6 1
\put{Figure 15} at 0 -1.78
\put{} at -6.9 0

\small
\put{$h_n(I)$} at 3.25 -0.6
\put{$I$} at 0 -0.6
\put{$\bullet$} at 3.62 0
\put{$p$} at 3.615 1.3

\endpicture


\vspace{0.6cm}

From what precedes one immediately deduces the following.

\vspace{0.15cm}

\begin{cor} {\em If $\Gamma$ is a finitely generated subgroup of 
$\mathrm{Diff}^{1+\mathrm{Lip}}_+(\mathrm{S}^1)$ having an 
excep- tional minimal set $\Lambda$, then there exist \esp 
$p\!\in\!\Lambda$ \esp and \esp $h \!\in\! \Gamma$ \esp such 
that \esp $h(p)\!=\!p$ \esp and \esp $h'(p)\!<\!1$.}
\end{cor}

\vspace{0.1cm}

Notice that Denjoy's theorem for $\mathrm{C}^{1+\mathrm{Lip}}$ diffeomorphisms may be deduced 
from this corollary. Actually, Sacksteder's theorem is customarily presented as a generalization 
of Denjoy's \index{Denjoy!theorem} 
theorem for groups of diffeomorphisms and codimension-one foliations. Nevertheless, 
this approach is not appropriate from a dynamical point of view. Indeed, Denjoy's theorem 
strongly uses the very particular combinatorics of circle homeomorphisms with irrational 
rotation number, \index{rotation number}
which is quite distinct (even locally) from that of pseudo-groups acting 
without invariant probability measure.\footnote{The reader should notice that 
an analogous dichotomy was exploited when dealing with general orderable groups 
in \S \ref{super-witte-morris}.} In the next section we will develop this view to 
obtain finer results than those of this section. Here we content ourselves by showing 
how to extend Sacksteder's theorem to non-Abelian groups acting minimally on the 
circle. We begin by a clever remark due to Ghys 
\index{Ghys!remark on Saksteder's theorem}
and taken from \cite{conf}.

\vspace{0.2cm}

\begin{prop} {\em Let $\Gamma$ be a pseudo-group of $\mathrm{C}^{1+\mathrm{Lip}}$ 
diffeomorphisms of the circle or of a bounded interval, all of 
whose orbits are dense. Suppose 
that there exists an element $f \!\in\! \Gamma$ fixing a single point on a non-degenerate 
(non-necessarily open) interval containing this point. Then $\Gamma$ contains an element 
with a hyperbolic fixed point.}
\label{rem-ghys}
\end{prop}

\noindent{\bf Proof.} Let $a$ be the fixed point given by the hypothesis, and let 
$a' \!>\! a$ be such that the restriction of $f$ to an interval $[a,a']$ has no 
other fixed point than $a$ (if such a point $a'$ does not exist, then one may 
consider an interval of type $[a',a]$ with a similar property). Changing $f$ by 
$f^{-1}$ if necessary, we may assume that $f(y) < y$ for all $y \! \in ]a,a']$. 
Since the orbits by $\Gamma$ are dense, there must exist $h \in \Gamma$ such 
that $h(a) \! \in ]a,a'[$. Let $n \in \mathbb{N}$ be large enough so that 
$[a,f^n(a')] \subset dom(h)$ and $h([a,f^n(a')]) \subset ]a,a']$.\\

Let us consider the intervals $A_0 \!=\! [a,f^n(a')]$ and 
$B_0 \!=\! h([a,f^n(a')])$, and let us define by induction 
the sets $A_{j+1} \!=\! f^n(A_j \cup B_j)$ and $B_{j+1}\!=\! h(A_{j+1})$. 
For $\epsilon > 0$ very small, consider the pseudo-group generated by the restrictions 
of $f^n$ and $h$ to the intervals $]a-\varepsilon, h f^n(a')+\varepsilon[$ \esp and 
\esp $]a-\varepsilon,f^n(a')+\varepsilon[$, respectively. With respect to this 
pseudo-group, the Cantor set $\Lambda = \cap_{j \in \mathbb{N}} (A_j \cup B_j)$
satisfies the hypotheses of Sacksteder's theorem. The result then follows. 
$\hfill\square$

\vspace{0.4cm}

By H\"older's theorem, 
\index{H\"older!theorem} \index{rotation!group}
a group of circle homeomorphisms which is not semiconjugate to a group of rotations 
contains an element $f$ satisfying the hypothesis of the preceding proposition. 
Since every semiconjugacy from a group whose orbits are dense to a group of 
rotations is necessarily a conjugacy, as a corollary we obtain the following.

\vspace{0.2cm}

\begin{cor} {\em Let $\Gamma$ be a subgroup of 
$\mathrm{Diff}^{1+\mathrm{Lip}}_+(\mathrm{S}^1)$. Suppose that the orbits of 
$\Gamma$ are dense and that $\Gamma$ is not conjugate to a group of rotations 
(equivalently, suppose that $\Gamma$ acts minimally and without invariant 
probability measure). Then $\Gamma$ contains an element 
with a hyperbolic fixed point.}
\label{remarque-etiene}
\end{cor}

\vspace{0.2cm}

In the next section, we will see that this corollary still holds in class $\ce^{1}$.

\vspace{0.1cm}

\begin{small} \begin{ejer} The fact that some nontrivial holonomy is hyperbolic is 
relevant when studying stability and/or rigidity properties (see for instance \S \ref{sec-st}). 
However, knowing that some holonomy is nontrivial (but not necessarily hyperbolic) may be 
also useful. The next proposition corresponds to a clever remark (due to Hector) 
\index{Hector} 
in this direction. Before passing to its proof, the reader should check 
\index{Thompson's groups!G}
its validity for the modular group and for Thompson's group $\mathrm{G}$. 
\index{group!modular}

\begin{prop} {\em Let $\Gamma$ be a finitely generated group of $\ce^{1+\mathrm{Lip}}$ 
circle diffeomorphisms admitting an exceptional minimal set $\Lambda$. If $p$ is the 
endpoint of one of the connected components $I$ of \esp $\clo \setminus \Lambda$, 
then there exists $g \!\in\! \Gamma$ fixing $p$ and whose restriction to 
$I \cap V$ is nontrivial for every neighborhood $\esp \! V$ of $p$.}
\end{prop}

\noindent{\underbar{Hint.}} Suppose not and obtain a contradiction 
to the fact that for every neighborhood of $p$ there exist elements in 
$\Gamma$ having a hyperbolically \textit{repelling} fixed point inside.
\label{ejer-hector-int}
\end{ejer} \end{small}


\subsection{The $\mathbf{\ce^1}$ version for pseudo-groups}
\label{sec-pseudo}

\hspace{0.45cm} The main object of this section and the next one consists in formulating 
and proving several generalizations of Sacksteder's theorem in class $\ce^1$. Some closely 
related versions were obtained by Hurder in \cite{hurder1,hurder2,Hu,hurdito} 
\index{Hurder} 
via dynamical methods coming from the so-called foliated Pesin's theory. The optimal and 
definitive versions (obtained by the author in collaboration with Deroin and Kleptsyn in 
\cite{DKN}) make a strong use of probabilistic arguments.
\index{foliation}\index{Deroin}\index{Kleptsyn}

\vspace{0.1cm}

\begin{thm} {\em If $\Gamma$ is a pseudo-group of $\ce^1$ diffeomorphisms of a compact 
one-dimensional manifold without invariant probability measure, then $\Gamma$ has 
elements with hyperbolic fixed points.}
\label{sac-ceuno}
\end{thm}

\vspace{0.1cm}

We will show this theorem inspired by the proof of Proposition \ref{rem-ghys}. 
Let us then suppose that $\Gamma$ contains two elements $f$ and $h$ satisfying 
(compare Definition \ref{def-crossed-elements}):

\vspace{0.1cm}

\noindent{(i) \thinspace the domain of definition of $f$ contains an interval $[a,a'[$ 
so that \esp $f(a)\!=\!a$ \esp and \esp $f$ \esp topologically contracts towards $a$;}

\vspace{0.1cm}

\noindent{(ii) \thinspace $h$ is defined on a neighborhood of \esp $a$ \esp 
and \esp $h(a) \!\in ]a,a'[$.}

\vspace{0.1cm}

Put $c \!=\! h(a)$ and fix $d' \!\!\in ]c,a'[$. Replacing $f$ by $f^n$ for 
$n \!\in \!\mathbb{N}$ large enough if necessary, we may suppose that $f(d') < c$, 
that $f(d')$ belongs to the domain of definition of $h$, and that $hf (d')\!\in ]c,d'[$. 
The last condition implies in particular that $hf$ has fixed points in $]c,d'[$. 
Let $d$ be the first fixed point of $hf$ to the right of $c$, and let $b = f(d)$. The 
interval $I\!= ]b,c[$ corresponds to the ``first gap'' (\textit{i.e.}, the ``central'' 
connected component of the complement) of a Cantor set $\Lambda$ which is invariant 
by $f$ and $g = hf$ (see Figure 16).

\vspace{0.45cm}


\beginpicture

\setcoordinatesystem units <0.8cm,0.8cm>

\putrule from -7 0 to -2 0

\putrule from -7 -0.2 to -7 0.2

\putrule from -2 -0.2 to -2 0.2

\putrule from -5 -0.13 to -5 0.13

\putrule from -4 -0.13 to -4 0.13

\putrule from -5 -0.015 to -4 -0.015

\putrule from -5 0.015 to -4 0.015

\putrule from -6.2 -0.07 to -6.2 0.07

\putrule from -5.8 -0.07 to -5.8 0.07

\putrule from -6.2 -0.015 to -5.8 -0.015

\putrule from -6.2 0.015 to -5.8 0.015

\putrule from -3.2 -0.07 to -3.2 0.07

\putrule from -2.8 -0.07 to -2.8 0.07

\putrule from -3.2 -0.015 to -2.8 -0.015

\putrule from -3.2 0.015 to -2.8 0.015

\putrule from -5.5 1.1 to -3.5 1.1

\plot -5.5 1.1
-5.35 1.15 /

\plot -5.5 1.1
-5.35 1.05 /

\putrule from -6 0.7 to -3 0.7

\putrule from -6 0.3 to -6 0.7

\putrule from -3 0.3 to -3 0.7

\plot -3 0.3
-3.05 0.45 /

\plot -3 0.3
-2.95 0.45 /


\putrule from 2 -2 to 7 -2
\putrule from 2 3 to 7 3
\putrule from 2 -2 to 2 3
\putrule from 7 -2 to 7 3

\putrule from 2.02 0 to 2.02 1
\putrule from 1.98 0 to 1.98 1

\putrule from 1.85 0 to 2.15 0
\putrule from 1.85 1 to 2.15 1

\putrule from 2.02 -1.2 to 2.02 -0.8
\putrule from 1.98 -1.2 to 1.98 -0.8

\putrule from 2.02 1.8 to 2.02 2.2
\putrule from 1.98 1.8 to 1.98 2.2

\putrule from 1.9 -1.2 to 2.1 -1.2
\putrule from 1.9 -0.8 to 2.1 -0.8
\putrule from 1.9 1.8 to 2.1 1.8
\putrule from 1.9 2.2 to 2.1 2.2

\plot
2 1
2.125 1.05296
2.25 1.109
2.375 1.16056
2.5 1.216
2.625 1.26056
2.75 1.3009
2.875 1.35296
3 1.4 /

\plot
3 1.4
3.125 1.45296
3.25 1.509
3.375 1.56056
3.5 1.616
3.625 1.66056
3.75 1.7009
3.875 1.75296
4 1.8 /

\plot
4 1.8
4.125 1.85296
4.25 1.909
4.375 1.96056
4.5 2.016
4.625 2.06056
4.75 2.1009
4.875 2.15296
5 2.2 /

\plot
5 2.2
5.125 2.25296
5.25 2.309
5.375 2.36056
5.5 2.416
5.625 2.46056
5.75 2.5009
5.875 2.55296
6 2.6 /

\plot
6 2.6
6.125 2.65296
6.25 2.709
6.375 2.76056
6.5 2.816
6.625 2.86056
6.75 2.9009
6.875 2.95296
7 3 /


\plot
7 0
6.875 -0.05296
6.75 -0.109
6.625 -0.16056
6.5 -0.216
6.375 -0.26056
6.25 -0.3009
6.125 -0.35296
6 -0.4 /

\plot
6 -0.4
5.875 -0.45296
5.75 -0.509
5.625 -0.56056
5.5 -0.616
5.375 -0.66056
5.25 -0.7009
5.125 -0.75296
5 -0.8 /

\plot
5 -0.8
4.875 -0.85296
4.75 -0.909
4.625 -0.96056
4.5 -1.016
4.375 -1.06056
4.25 -1.1009
4.125 -1.15296
4 -1.2 /

\plot
4 -1.2
3.875 -1.25296
3.75 -1.309
3.625 -1.36056
3.5 -1.416
3.375 -1.46056
3.25 -1.5009
3.125 -1.55296
3 -1.6 /

\plot
3 -1.6
2.875 -1.65296
2.75 -1.709
2.625 -1.76056
2.5 -1.816
2.375 -1.86056
2.25 -1.9009
2.125 -1.95296
2 -2 /


\setdots

\putrule from 4 -2 to 4 3
\putrule from 5 -2 to 5 3
\putrule from 2 0 to 7 0
\putrule from 2 1 to 7 1

\putrule from 2 -1.2 to 4 -1.2
\putrule from 2 -0.8 to 5 -0.8
\putrule from 2 1.8 to 4 1.8
\putrule from 2 2.2 to 5 2.2

\plot 2 -2
7 3 /

\put{Figure 16} at 0 -2.6
\put{} at -8.6 0

\small

\put{$a$} at -7 -0.5
\put{$d$} at -2 -0.5
\put{$f$} at -4.5 1.3
\put{$h$} at -3.6 0.5
\put{$I\!= ]f(d),h(a)[$} at -4.5 -0.45


\put{$a$} at 1.7 -2.3
\put{$b$} at 4 -2.3
\put{$c$} at 5 -2.3
\put{$d$} at 7 -2.3

\put{$b$} at 1.5 0
\put{$c$} at 1.5 1
\put{$d$} at 1.7 3
\put{$I$} at 1.75 0.5
\put{$f$} at 4.7 -1.2
\put{$g$} at 4.7 1.8

\endpicture


\vspace{0.3cm}

\begin{prop} {\em With the above notation, the pseudo-group generated by 
$f$ and $g$ contains elements with hyperbolic fixed points in $\Lambda$.}
\label{resorte}
\end{prop}

\vspace{0.1cm}

The proof of this proposition is particularly simple in class $\ce^{1+\tau}$ 
(where $\tau \! > \! 0$). Moreover, it allows illustrating the usefulness of 
the probabilistic methods in the theory. The (slightly more technical) $\ce^1$ 
case will be discussed later. For an alternative ``purely deterministic" proof 
of a closely related result, we refer the reader to \cite{katok-per}.

\vspace{0.35cm}

\noindent{\bf Proof of Proposition \ref{resorte} in class $\ce^{1+\tau}\!$.} 
\esp On $\Omega \!=\! \{ f,g \}^{\mathbb{N}}$ let us consider the Ber- noulli 
measure $\mathbb{P}$ giving mass $1/2$ to each ``random variable" $f$, $g$. 
For $\omega\!=\!(g_1,g_2,\ldots)$ in $\Omega$ and $n \!\in\! \mathbb{N}$, let 
$h_n(\omega) \!=\! g_n \cdots g_1$, and let $h_0(\omega)\!=\!id$. Since 
the interval $I$ is wandering --in the sense that the intervals in 
the family \esp $\big\{ g_n \cdots g_1 (I)\!: n \geq 0, \thinspace 
(g_1,\ldots,g_n) \in \{f,g\}^n \big\}$ \esp are two-by-two disjoint--,
$$\sum_{n \geq 0} \sum_{(g_1,\ldots,g_n) \in \{ f, g \}^n} |g_n \cdots g_1 (I)|
\esp \leq \esp d\!-\!a \esp < \esp \infty.$$
Moreover, Fubini's theorem gives
\begin{small}
$$\sum_{n \geq 0} \sum_{(g_1,\ldots,g_n) \in \{ f, g \}^n} |g_n \cdots g_1 (I)|
= \sum_{n \geq 0} 2^n \Big( \int_{\Omega} |h_n(\omega)(I)|
\thinspace d\mathbb{P}(\omega) \Big) = 
\int_{\Omega} \Big( \sum_{n \geq 0} 2^n |h_n(\omega)(I)| \Big) d\mathbb{P}(\omega).$$
\end{small}This allows us to 
conclude that, for $\mathbb{P}$-almost every random sequence 
$\omega \!\in\! \Omega$, the series \esp $\sum  2^n |h_n (\omega)(I)|$ 
\esp converges. For each $B>0$ let us consider the set 
$$\Omega(B) \esp = \esp \big\{ \omega \in \Omega \!: \esp |h_n(\omega)(I)| 
\esp \leq \esp B/2^n \esp \mbox{ for all } \esp n \geq 0 \big\}.$$
The probability $\mathbb{P} [\Omega(B)]$ converges to $1$ as $B$ goes to 
infinity. In particular, we may fix $B$ so that $\mathbb{P}[\Omega(B)]>0$.
Notice that if $\omega$ belongs to $\Omega(B)$ then 
\begin{equation}
\sum_{n \geq 0} |h_n(\omega)(I)|^{\tau} \esp \leq \esp 
B^{\tau} \sum_{n \geq 0} \frac{1}{2^{n \tau}} 
\esp = \esp S \esp < \esp \infty.
\label{cota}
\end{equation}

\index{wandering interval}
Let us consider the interval $J' = [b-\ell,c+\ell]$ containing the wandering interval 
$I$, where \thinspace $\ell \!=\! \ell(\tau,S,|I|;\{f,g\})$ \thinspace is the constant 
which appears in Lemma \ref{?}. If $N \in \mathbb{N}$ is large enough, then both \esp 
$f^N g$ \esp and \esp $g^N f$ \esp send the whole interval $[0,1]$ into $J' \setminus I$. 
A direct application of the Borel-Cantelli Lemma then gives  
\index{Borel-Cantelli Lemma}
$$\mathbb{P}[h_n(\omega)(I) \subset J' \setminus I \mbox{ infinitely many times}] = 1.$$
If $\omega \!\in\! \Omega(B)$ and $n \!\in\! \mathbb{N}$ satisfy 
$h_n(\omega)(I) \subset J' \setminus I$, then Lemma \ref{?} shows that 
$h_n(\omega)$ has a hyperbolic fixed point $p$. Finally, since $\Lambda$ 
is invariant by the pseudo-group and the fixed point of $h_n(\omega)$ that we 
found contracts part of this set, $p$ must belong to $\Lambda$. $\hfill\square$

\vspace{0.45cm}

The proof of the general version (in class $\ce^1$) of Proposition \ref{resorte} 
needs some technical improvements for finding hyperbolic fixed points in the 
absence of control of distortion. The main idea lies in that, when one knows 
{\em a priori} that the dynamics is differentiably contracting somewhere, 
the continuity of the derivatives forces this contraction to persists 
(perhaps with a smaller rate) on a larger domain.

\vspace{0.15cm}

Keeping the notation of the preceding proof, fix once and for all a constant 
$\varepsilon \!\in ]0,1/2[$. We know that for $\mathbb{P}$-almost every 
$\omega \in \Omega$ there exists $B = B(\omega) \geq 1$ such that
\begin{equation}
|h_n(\omega)(I)| \leq \frac{B}{2^n}
\quad \mbox{ for all} \quad n \geq 0.
\label{int-exp}
\end{equation}

\vspace{0.02cm}

\begin{lem} {\em There exists a constant $\bar{C}$ depending only on $f$ 
and $g$ such that, if $\omega = (g_1,g_2,\ldots) \in \Omega$ satisfies
{\em (\ref{int-exp})}, 
then for all $x \in I$ and every integer $n \geq 0$ one has}
\begin{equation}
h_n(\omega)'(x) \leq \frac{B \bar{C}}{(2-\varepsilon)^n}.
\label{hyp-inf}
\end{equation}
\label{aaa}
\end{lem}

\vspace{-0.3cm}

\noindent{\bf Proof.} Let us fix $\varepsilon_0 \!>\! 0$ small enough so that 
for every pair of points $y,z$ in $[a,d]$ at distance less than or equal 
to $\varepsilon_0$ one has
\begin{equation}
\frac{f'(y)}{f'(z)} \leq \frac{2}{2-\varepsilon}, \qquad 
\qquad \frac{g'(y)}{g'(z)} \leq \frac{2}{2-\varepsilon}.
\label{continuity}
\end{equation}
It is easy to check the existence of $N \in\! \mathbb{N}$ such that, for every 
$\omega \!\in\! \Omega$ and every $i \!\geq\! 0$, the length of the interval 
$h_{N+i}(\omega)(I)$ is less than or equal to $\varepsilon_0$. We claim 
that (\ref{hyp-inf}) holds for $\bar{C} = \max\{ A,\bar{A} \}$, where 
$$A = \!\!\! 
\sup_{x\in I, n\leq N, \omega \in \Omega} \frac{h_n(\omega)'(x) \esp
(2-\varepsilon)^n}{B}, \quad \quad \bar{A} = \!\!\! 
\sup_{x,y \in I, \omega \in \Omega} \frac{h_N(\omega)'(x)}
{h_N(\omega)'(y) \esp |I|} \esp \Big( \frac{2-\varepsilon}{2} \Big)^N \!\!.$$
Indeed, if $n \!\leq\! N$ then (\ref{hyp-inf}) holds due to the condition \esp 
$\bar{C} \!\geq\! A$. \esp For $n \!>\! N$ let us fix $y \!=\! y(n) \!\in\! I$ 
such that \esp $|h_n(\omega)(I)| = h_n(\omega)'(y) \esp |I|$. For $x \!\in\! I$, 
the distance between $h_{N+i}(\omega)(x)$ and $h_{N+i}(\omega)(y)$ is smaller 
than or equal to $\varepsilon_0$ for all $i \geq 0$. Hence, by (\ref{continuity}),
$$\frac{h_n(\omega)'(x)}{h_n(\omega)'(y)} 
= \frac{h_N(\omega)'(x)}{h_N(\omega)'(y)}
\esp \frac{g_{N+1}'(h_N(\omega)(x))}{g_{N+1}'(h_N(\omega)(y))} \cdots
\frac{g_n'(h_{n-1}(\omega)(x))}{g_n'(h_{n-1}(\omega)(y))}
\leq \frac{h_N(\omega)'(x)}{h_N(\omega)'(y)} \esp
\Big( \frac{2}{2-\varepsilon} \Big)^{n-N}\!\!,$$
and therefore
\begin{small}
$$h_n(\omega)'(x)
\leq \frac{h_N(\omega)'(x)}{h_N(\omega)'(y)} \esp \frac{|h_n(\omega)(I)|}{|I|} 
\Big( \frac{2}{2-\varepsilon} \Big)^{n-N} \!\!
\leq \frac{h_N(\omega)'(x)}{h_N(\omega)'(y)} \esp \frac{B}{|I|\esp 2^n}
\Big( \frac{2}{2-\varepsilon} \Big)^{n-N} \!\!
\leq \frac{B\bar{C}}{(2-\varepsilon)^n},$$
\end{small}where the last inequality follows from the condition \esp 
$\bar{C} \geq \bar{A}$. \esp $\hfill\square$

\vspace{0.5cm}

To verify the persistence of the contraction beyond the interval $I$, we will use a 
kind of ``dual'' argument. Fix a small constant $\varepsilon_1 \!>\! 0$ such that, for 
every $y,z$ in $[a,d]$ at distance less than or equal to $\varepsilon_1$, one has
\begin{equation}
\frac{f'(y)}{f'(z)} \leq \frac{2-\varepsilon}{2-2\varepsilon}, \qquad 
\qquad \frac{g'(y)}{g'(z)} \leq \frac{2-\varepsilon}{2-2\varepsilon}.
\label{continuity-2}
\end{equation}

\vspace{0.2cm}

\begin{lem} {\em Let $C \geq 1$, $\omega = (g_1,g_2,\ldots) \in \Omega$, 
and $x \in [a,d]$, be such that}
\begin{equation}
h_n(\omega)'(x) \leq \frac{C}{(2-\varepsilon)^n} \quad \mbox{for all } n \geq 0.
\label{hypot}
\end{equation}
{\em If $y \in [a,d]$ is such that \esp $dist(x,y) \leq \varepsilon_1 / C$ 
\esp and \esp $n \geq 0$, then}
\begin{equation}
h_n(\omega)'(y) \leq \frac{C}{(2 - 2\varepsilon)^n}.
\label{tesis}
\end{equation}
\label{bbb}
\end{lem}

\vspace{-0.2cm}

\noindent{\bf Proof.} We will verify inequality (\ref{tesis}) by induction. For $n\!=\!0$ 
it holds due to the condition $C \! \geq \! 1$. Let us assume that it holds for every 
$j \!\in\! \{0,\ldots,n\}$, and let $y_j = h_{j}(\omega)(y)$ and $x_j = h_j(\omega)(y)$. 
Suppose that $y \leq x$ (the case $y \geq x$ is analogous). Each point $y_j$ belongs 
to the interval \esp $h_j(\omega)([x-\varepsilon_1/C,x])$. \esp By the 
induction hypothesis, 
$$\big| h_j(\omega)([x-\varepsilon_1/C,x]) \big| \leq \frac{C}{(2-2\varepsilon)^j}
\esp \big| [x-\varepsilon_1/C,x] \big| \leq C \esp \frac{\varepsilon_1}{C} =
\varepsilon_1.$$
By the definition of $\varepsilon_1$, for every $j \leq n$ one has \esp 
$g_{j+1}'(y_j) \leq g_{j+1}'(x_j) \esp \big( \frac{2-\varepsilon}{2-2\varepsilon}
\big)$. \esp Therefore, due to hypothesis (\ref{hypot}),
\begin{eqnarray*}
h_{n+1}(\omega)'(y) 
\esp \esp = \esp \esp 
g_1'(y_0) \cdots g_{n+1}'(y_n) 
&\leq& 
g_1'(x_0) \cdots g_{n+1}'(x_n) \esp \Big( \frac{2-\varepsilon}{2-2\varepsilon} \Big)^{n+1}\\
&\leq&
\frac{C}{(2-\varepsilon)^{n+1}} \esp
\Big( \frac{2-\varepsilon}{2-2\varepsilon} \Big)^{n+1} 
\esp \esp \esp \esp \leq 
\esp \esp \esp \esp \frac{C}{(2-2\varepsilon)^{n+1}},
\end{eqnarray*}
which concludes the inductive proof of (\ref{tesis}). $\hfill\square$

\vspace{0.3cm}

We can now complete the proof of Proposition \ref{resorte}. To do this, let us notice 
that, by Lemma \ref{aaa}, if $C$ is large enough then the probability of the set
$$\Omega(C,\varepsilon) = \Big\{ \omega \in \Omega : h_n(\omega)'(x) \leq
\frac{C}{(2-\varepsilon)^n}
\mbox{ for every } \esp n \geq 0 \esp \mbox{ and every } \esp x \in I \Big\}$$
is positive. Fix such a $C \!\geq\! 1$, let $\ell \!=\! \min \{\varepsilon_1/2 C, |I|/2 \}$, 
and let us denote by $J$ the $2\ell$-neighborhood of $I$. Lemma \ref{bbb} implies that, 
for every $\omega \in \Omega(C,\varepsilon)$, every $n \geq 0$, and every $y \in J$,
\begin{equation}
h_n(\omega)'(y) \leq \frac{C}{(2-2\varepsilon)^n}.
\label{mato}
\end{equation}
If $J'$ denotes the $\ell$-neighborhood of the interval $I$, then the probability of the set 
$$\big\{\omega: h_n(\omega)(I) \subset J' \setminus I \esp \mbox{ infinitely many times} \big\}$$  
equals 1. Since $\varepsilon < 1/2$, we deduce the existence 
of $\omega \in \Omega(C,\varepsilon)$ and $m \in \mathbb{N}$ such that 
$$h_m(\omega)(I) \subset J' \setminus I, \qquad \quad (2-2\varepsilon)^m > C.$$ 
By (\ref{mato}) and the inequality $\ell \leq |I|/2$, this implies that $h_m(\omega)$ 
sends $J$ into one of the two connected components of $J \setminus I$ in such a way 
that $h_m(\omega)$ has a fixed point in this component, which is necessarily hyperbolic 
and belongs to $\Lambda$. The proof of Proposition \ref{resorte} is then concluded.

\vspace{0.3cm}

From the discussion above we see that, in order to prove Theorem \ref{sac-ceuno}, it suffices 
to show that every pseudo-group of $\ce^1$ diffeomorphisms of a one-dimensional compact 
manifold without invariant probability measure must contain elements $f$ and $h$ which 
satisfy the properties (i) and (ii) at the beginning of this section (compare 
Proposition \ref{radon-inv}). Although this is actually true, the proof uses 
quite involved techniques related to the construction of the Haar measure, and  
we do not want to enter into the details here for avoiding overload the 
presentation: the interested reader may find the proof in \cite{DKN}. 

\begin{small}\begin{obs} 
An interesting issue of the proof above is the fact that it allows suppressing 
the hypothesis of compact generation for the pseudo-group (see Footnote 
\ref{servira} in \S \ref{chuartz}). 
\index{pseudo-group!compactly generated}
\label{remark-haefliger}
\end{obs}\end{small}

In the case of \textit{groups} acting on the circle without invariant probability 
measure, getting a pair of elements $f,h$ as before is relatively simple: for 
minimal actions this may be done by using H\"older's 
\index{H\"older!theorem} 
theorem, whereas for actions with an exceptional minimal 
set a similar argument applies to the action induced on the topological circle obtained 
after collapsing the connected components of its complement. However, 
in the next section we will see that, 
for groups of $\ce^1$ circle diffeomorphisms without invariant 
probability measure, a much stronger conclusion than the one of Sacksteder's theorem 
holds: these groups have elements having finitely many fixed points, all of them 
hyperbolic.

\vspace{0.12cm}

\begin{small} \begin{ejer} Anticipating the main result of the next section, 
show that every group $\Gamma$ of $\ce^{1+\mathrm{Lip}}$ circle diffeomorphisms 
without invariant probability measure contains elements with at least $2 d(\Gamma)$
hyperbolic fixed points, half of them contracting and half of them repelling.

\noindent{\underbar{Hint.}} Consider the minimal and strongly 
expansive case (the general case easily reduces to this one). 
\index{action!strongly expansive}\index{action!expansive}Prove 
that $\Gamma$ contains two elements $f,g$ such that for some 
interval $I \!=\![a,b]$ one has $g(I) \cup g^{-1}(I) \!\subset\! \clo \setminus I$  as well 
as $\lim_{n \rightarrow \infty} f^n(x) = b$ and $\lim_{n \rightarrow \infty} f^{-n}(x) = a$ 
for all $x \in \clo \setminus I$. Then, for $\varepsilon > 0$ very small, apply Schwartz' 
technique to the sequences of compositions $(f^{-n} g)$ and $(g^{-1}f^{n})$ over the 
intervals $[a,a+\varepsilon]$ and $g^{-1}([b-\varepsilon,b])$.
\index{Schwartz estimates}
\label{dos-hyp}
\end{ejer} \end{small}

\vspace{0.12cm}

Remark that the hypothesis of non-existence of invariant probability measure assumed 
throughout this section is necessary for the validity of Theorem \ref{sac-ceuno}. An easy 
(but uninteresting) example is given by (the group generated by) a diffeomorphism having 
fixed points, all of them parabolic. More interesting examples correspond to $\ce^1$ 
(and even $\ce^{1+\tau}$) circle diffeomorphisms with irrational rotation number and 
non-dense orbits (see \S \ref{ejemplosceuno} for the construction of such diffeomorphisms). 
These last examples show that the statement of Theorem \ref{sacksteder} is not valid 
for groups of $\ce^{1+\tau}$ circle diffeomorphisms. The $\ce^{1+\mathrm{bv}}$ case 
is very special. The proposition below is still true for 
pseudo-groups of $\ce^{1+\mathrm{bv}}$ diffeomorphisms having 
a minimal invariant Cantor set (compare \cite{ega}).

\vspace{0.25cm}

\begin{prop} {\em If \esp $\Gamma$ is a finitely generated group of $\ce^{1+\mathrm{bv}}$ 
circle diffeomorphisms preserving an exceptional minimal set, then $\Gamma$ has 
elements with hyperbolic fixed points.}
\label{cc-et}
\end{prop}

\noindent{\bf Proof.} It suffices to show that $\Gamma$ cannot preserve a probability measure, 
because this allows applying Theorem \ref{sac-ceuno}. However, if $\Gamma$ preserves 
a probability measure, then it is semiconjugate to a group of rotations. Since $\Gamma$ 
is finitely generated, at least one of its generators must have irrational rotation number. 
Nevertheless, the existence of such an element is in contradiction with Denjoy's theorem. $\hfill\square$

\vspace{0.5cm}

As we will see in the next section, under the hypothesis of Proposition \ref{cc-et} 
the group $\Gamma$ has elements having finitely many fixed points, all of them hyperbolic.

\vspace{0.12cm}

\begin{small}\begin{obs} We ignore whether Proposition \ref{cc-et} extends to groups of piecewise 
affine homeomorphisms. Actually, it seems to be unknown whether some version of Sacksteder's   
theorem holds in this context. This problem is certainly related to that of knowing what 
are the subgroups of $\mathrm{PAff}_+(\clo)$ which are conjugate to groups of $\ce^1$ circle 
diffeomorphisms (compare \S \ref{ghys-sergi}). A seemingly related problem consists in 
finding a purely combinatorial proof of Denjoy's theorem in the piecewise affine 
\index{piecewise affine homeomorphism} case.\index{Denjoy!theorem}
\label{no-se-que-pasa}
\end{obs}\end{small}


\subsection{A sharp $\mathbf{\ce^1}$ version via Lyapunov exponents}

\hspace{0.45cm} For groups of $\ce^1$ circle diffeomorphisms, it is possible to give a 
global and optimal version of Sacksteder's theorem. Recall that to every group $\Gamma$ 
of circle homeomorphisms without invariant probability measure, there is an associated 
\textit{degree} \esp $d(\Gamma) \!\in\! \mathbb{N}$ \esp (see \S \ref{sec-margulis}).
\index{degree!of a group of circle homeomorphisms}

\vspace{-0.2cm}

\begin{thm} {\em If $\Gamma$ is a finitely generated subgroup of $\mathrm{Diff}_+^1(\clo)$ 
without invariant probability measure, then $\Gamma$ contains elements having finitely many 
fixed points, all of them hyperbolic. More precisely, $\Gamma$ contains elements having 
$2d(\Gamma)$ fixed points, half of them hyperbolically contracting and half of them 
hyperbolically repelling.}
\label{morse-smale}
\end{thm}

\vspace{0.25cm}

The probabilistic methods will again be fundamental for the proof of this theorem. Let 
$\Gamma$ be a finitely generated subgroup of $\mathrm{Diff}_+^1(\clo)$ endowed with a 
non-degenerate finitely supported probability measure $p$. Let $\mu$ be a probability 
on $\clo$ which is stationary with respect to $p$, and let $T$ be the skew-product 
map from $\Omega \times \clo$ into itself given by
$$T(\omega,x)=\big( \sigma(\omega),h_1(\omega)(x) \big).$$ 
This map $T$ preserves $\mathbb{P}\!\times\!\mu$, 
and hence the Birkhoff Ergodic Theorem \cite{walters} 
\index{Birkhoff!Ergodic Theorem} applied to 
the function \thinspace $(\omega,x) \mapsto \log \big(h_1(\omega)'(x) \big)$ 
\thinspace shows that, for $\mu$-almost every point $x \!\in\! \clo$ and 
$\mathbb{P}$-almost 
every random sequence $\omega \in \Omega$, the following limit exists:
$$\lambda_{(\omega, x)}(\mu) =
\lim_{n \rightarrow \infty} \frac{\log \big( h_n(\omega)'(x) \big)}{n}.$$
This value will be called the {\bf {\em Lyapunov exponent}} of $(\omega,x)$. 
\index{Lyapunov exponent}

To continue with our discussion, the Kakutani Random Ergodic Theorem 
\index{Kakutani!Random Ergodic Theorem} below will be essential (the 
reader may find more details on this in \cite{furman,kifer}).

\vspace{0.1cm}

\begin{thm} {\em If the stationary measure $\mu$ cannot be written as a nontrivial 
convex combination of two stationary probabilities, then the map $T$ is ergodic with 
respect to $\mathbb{P} \!\times\! \mu$.}
\label{kak-th}
\end{thm}

\noindent{\bf Proof.} Let $\psi\!: \Omega \times \clo \rightarrow \mathbb{R}$ be 
a $T$-invariant integrable function. We need to show that $\psi$ is almost everywhere 
constant. For each $n \geq 0$ let $\psi_n$ be the conditional expectation of $\psi$ 
given the first $n$ entries of $\omega \!\in\! \Omega$. We have \esp 
$\psi_0(\omega,x) = \mathbb{E} (\psi_1) (\omega,x).$ \esp 
Since $\psi$ is invariant by $T$, for each $n \in \mathbb{N}$ we have  
$$\psi(\omega,x) = \psi \big( \sigma^n(\omega),h_n(\omega)(x) \big),$$
which allows easily to deduce that
\begin{equation}
\psi_n(\omega,x) = \psi_{n-1} \big( \sigma(\omega),h_1(\omega)(x) \big) 
= \ldots\ldots = \psi_0 \big( \sigma^n(\omega),h_n(\omega)(x) \big).
\label{a-iterar}
\end{equation}
We then conclude that 
$$\psi_0 (\omega,x) = \mathbb{E} (\psi_0) (\sigma(\omega),h_1(\omega)(x)).$$
Now, since $\psi_0(\omega,x)$ does not depend on $\omega$, we may 
write $\psi_0(\omega,x) = \bar{\psi}(x)$, and due to the above equality, 
$\bar{\psi}$ is invariant by the diffusion operator. We claim that this implies  
that the function $\bar{\psi}$ is $\mu$-almost everywhere constant. Before 
showing it, let us remark that, by (\ref{a-iterar}), this implies that all 
the functions $\psi_n$ are constant, and hence the limit function $\psi$ 
is also constant, thus finishing the proof of the ergodicity of $T$.

To prove that $\bar{\psi}$ is constant, let us suppose for a contradiction 
that for some $c \!\in\!\mathbb{R}$ the sets $S_{+} = \{\bar{\psi} \geq c\}$ 
and $S_{-} = \{\bar{\psi} < c\}$ have positive $\mu$-measure. The invariance 
of $\bar{\psi}$ and $\mu$ by the diffusion \index{diffusion operator}
allows easily to show that these sets are invariant by the action. This 
implies that \esp $\mu = \mu \mathcal{X}_{S_+} + \mu \mathcal{X}_{S_{-}}$ 
\esp is a nontrivial convex decomposition of $\mu$ into stationary 
probabilities, thus contradicting our hypothesis. $\hfill\square$

\vspace{0.15cm}

\begin{obs} A stationary measure which cannot be written as a nontrivial convex 
combination of two stationary probabilities is said to be {\bf \em{ergodic}}. 
Obviously, if the stationary measure is unique, then it is necessarily ergodic. 
In general, every stationary measure may be written as a ``wighted mean'' of  
ergodic ones (the so-called {\bf \em{ergodic decomposition}}; 
see for instance \cite{furman,kifer}).
\index{ergodic!decomposition}
\end{obs}

\vspace{0.15cm}

Whenever $T$ is ergodic, the Lyapunov exponent is 
$\mathbb{P} \times \mu$-almost everywhere equal to
$$ \lambda (\mu) = \int_{\Omega} \int _{\clo} \log \big( h_1(\omega)'(x) \big)
\thinspace d\mu(x) \thinspace d \mathbb{P} (\omega) = 
\int_{\Gamma} \int _{\clo} \log \big( g'(x) \big)
\thinspace d\mu(x) \thinspace dp(g).$$
If $p$ is moreover symmetric and the stationary measure 
$\mu$ is invariant by the action of $\Gamma$, then the Lyapunov exponent 
$\lambda_{(\omega,x)}(\mu)$ is equal to zero for almost every $(\omega,x)$. 
Indeed, if $\mu$ is invariant then the map from $\Omega \times \clo$ into itself 
defined by $(\omega,x) \mapsto \big(\omega,h_1(\omega)^{-1}(x)\big)$ preserves 
$\mathbb{P}\!\times\!\mu$. Due to the symmetry of $p$, this implies that 
$\lambda(\mu)$ coincides with 
$$\int_{\Omega} \int _{\clo} \!\log \big( h_1(\omega)'(x) \big)
\thinspace d\mu(x) \thinspace d \mathbb{P} (\omega)
\!=\! \int_{\Omega} \int _{\clo} \! 
\log \big( h_1(\omega)'(h_1(\omega)^{-1}(x)) \big)
\thinspace d\mu(x) \thinspace d \mathbb{P} (\omega)=$$
$$=\!-\!\int_{\Omega} \int _{\clo} \! \log \big( (h_1(\omega)^{-1})'(x) \big)
\thinspace d\mu(x) \thinspace d \mathbb{P} (\omega) \!=\! -\int_{\Gamma} \int_{\clo}
\! \log \big( (g^{-1})'(x) \big) \esp d\mu(x) \esp dp(g),$$
that is, $\lambda(\mu) = - \lambda(\mu)$. To summarize, if  
$\mu$ is invariant and ergodic, then its Lyapunov exponent equals zero. 
The general case may be deduced from this by an ergodic decomposition 
type argument.

\vspace{0.1cm}

The next proposition (which admits a more general version for codimension-one 
foliations: see \cite{DK}) may be seen as a kind of converse to the above remark.
\index{foliation}

\vspace{0.12cm}

\begin{prop} {\em If $p$ is non-degenerate, symmetric, and has finite 
support, and if $\Gamma$ does not preserve any probability 
measure on the circle, then the Lyapunov 
\index{measure!stationary} exponent of the 
unique\footnote{The uniqueness of the stationary measure 
comes from Theorem \ref{unicidad}.} 
stationary measure is negative.}
\label{neg}
\end{prop}


Before giving the proof of this proposition, we will explain how to obtain Theorem 
\ref{morse-smale} from it. Suppose for instance that the action of $\Gamma$ is minimal 
and strongly expansive 
\index{action!strongly expansive} 
(\textit{i.e.}, $d(\Gamma) = 1$: see \S \ref{sec-margulis}). 
In this case, we know that for every $\omega$ in a subset $\Omega^*$ of total 
probability of $\Omega$, the contraction coefficient $\mathrm{contr} (h_n(\omega))$ 
converges to 0 as $n$ goes to infinity. This means that for every $\omega \in \Omega^*$ 
there exist two closed intervals $I_n(\omega)$ and $J_n(\omega)$ whose lengths tend to 0 
and such that $h_n(\omega) \big( \overline{\clo \setminus I_n(\omega)} \big) = J_n(\omega)$. 
Moreover, $I_n(\omega)$ converges to a point $\varsigma_+(\omega)$, and the map 
$\varsigma_+\!\!: \Omega^* \rightarrow \clo$ is measurable. If $I_n(\omega)$ and 
$J_n(\omega)$ are disjoint, then the fixed points of $h_n(\omega)$ are in the interior 
of these two intervals. To show the uniqueness and the hyperbolicity of at least the 
fixed point in $J_n(\omega)$ (\textit{i.e.}, of the contracting fixed point), one uses 
the fact that the Lyapunov exponent is negative. However, two technical difficulties 
immediately appear: the intervals $I_n(\omega)$ and $J_n(\omega)$ are not necessarily 
disjoint, and the ``deviation'' from the local contraction rate depends on the 
initial point $(\omega,x)$. To overcome the first problem, it suffices to notice that the 
``times'' $n \!\in\! \mathbb{N}$ for which $I_n(\omega) \cap J_n(\omega) \neq \emptyset$ 
are rare (the density of this set of integers is generically equal to 0). The second 
difficulty is overcome by using the fact that the map $T$ is ergodic, which implies 
that almost every initial point $(\omega,x)$ will enter into a set where the 
deviation from the local contraction rate is well-controlled. 
Let us finally remark that, to show the uniqueness and the 
hyperbolicity of the repelling fixed point, one may apply an analogous argument 
to the compositions of the inverse maps in the opposite order (compare Exercise 
\ref{dos-hyp}). To do this, it suffices to consider the finite time distributions. 
Indeed, since the probability $p$ is supposed to be symmetric, the distributions 
of both Markov processes coincide for every finite time.
\index{Markov!process}

If the action is minimal and $d(\Gamma) > 1$, the preceding arguments may be 
applied ``after passing to a finite quotient''. Finally, if there is an exceptional 
minimal set, we may argue in an analogous way ``over this set''. Making these 
arguments formal is just a technical issue, and hence we will not develop the details 
here: the interested reader is referred to \cite{DKN}. However, we point out that from 
the arguments contained in \cite{DKN} one may deduce a more concise statement, namely 
if for $d \!\in\! \mathbb{N}$ we denote by $\mathrm{D}_d(\clo)$ the set of $\ce^1$ 
circle diffeomorphsims having exactly $2 d$ periodic points, all of them 
hyperbolic, then
$$\mathbb{P} \Big[ \lim_{n\rightarrow \infty} \frac{1}{N} \esp card
\big\{ n \! \in \! \{1,\ldots,N\} \!: \esp 
h_n(\omega) \in \mathrm{D}_{d(\Gamma)}(\clo) \big\}  = 1 \Big] = 1.$$

\vspace{0.02cm}

\begin{small}\begin{ejer}
Using some of the arguments above, prove that every group $\Gamma$ of circle homeomorphisms 
without invariant probability measure contains elements with at least $2d(\Gamma)$ fixed points. 
Give an example of such a group so that, for every element having fixed points, the number of 
these points is greater than $2d(\Gamma)$.\\

\noindent{\underbar{Remark.}} We ignore whether the group of piecewise affine circle homeomorphisms 
contains subgroups as the example asked for above (compare Remark \ref{no-se-que-pasa}).
\end{ejer}\end{small}

\vspace{0.1cm}

We now pass to the proof of Proposition \ref{neg}. Let  
$$\psi (x) = \int_{\Gamma} \log \big( g'(x) \big) \thinspace dp (g),$$
and let us suppose for a contradiction that \thinspace $\lambda(\mu) \geq 0$, 
\thinspace that is,
\begin{equation}
\int_{\clo} \psi(x) \thinspace d \mu(x) \geq 0.
\label{suppp}
\end{equation}
In this case, we will show that $\Gamma$ preserves a probability measure. 
For this, we will strongly use a lemma which is inspired from Sullivan's 
theory of foliated cycles \cite{Sullivan} (see also \cite{Gh}). Recall that the 
\index{Laplacian} \index{harmonic function}
{\bf \textit{Laplacian}} $\Delta \zeta$ of 
\index{Sullivan}
a continuous real-valued function $\zeta$ is defined 
by $\Delta \zeta = D \zeta - \zeta$, where $D$ denotes the diffusion operator (\ref{difundir}). 
A function $\psi$ is {\bf \textit{harmonic}} if its Laplacian is identically equal to 0 
({\em c.f.,} Exercise \ref{func-harmonica}).

\vspace{0.1cm}

\begin{lem} {\em Under hypothesis {\em (\ref{suppp})}, there exists 
a sequence of continuous functions $\zeta_n$ defined on the circle such 
that, for every integer $n\in \mathbb{N}$ and every point $x \in \clo$,}
\begin{equation}
\psi(x) + \Delta \zeta _n (x) \geq -\frac{1}{n}.
\label{casi-cociclo}
\end{equation}
\end{lem}

\noindent{\bf Proof.} Let us denote by $C(\clo)$ the space of continuous functions on the 
circle. Let $E$ be the subspace of functions arising as the Laplacian of some function in 
$C(\clo)$, and let $C_+$ be the convex cone in $C(\clo)$ formed by the non-negative functions. 
We need to show that if $\psi$ satisfies (\ref{suppp}), then its image by the projection 
$\pi \!\!: C(\clo) \rightarrow C(\clo) / \bar{E}$ is contained in $\pi(C_+)$. Suppose 
that this is not the case. Then Hahn-Banach's separation theorem provides us with 
a continuous linear functional 
$\bar{\mathrm{I}}: C (\clo) / \bar{E} \rightarrow \mathbb{R}$ such that 
$\bar{\mathrm{I}}(\pi(\psi)) < 0 \leq \bar{\mathrm{I}} (\pi(\eta))$ for every 
$\eta \in C_+$. Obviously, $\bar{\mathrm{I}}$ induces a continuous linear functional 
$\mathrm{I}\!: C (\clo) \rightarrow \mathbb{R}$ which is identically zero on $E$ 
and such that $\mathrm{I}(\psi) < 0 \leq \mathrm{I}(\eta)$ for every $\eta \in C_+$. 
We claim that there exists a constant $c \!\in\! \mathbb{R}$ such that \esp 
$\mathrm{I} = c \esp \mu$ \esp (we identify the probability measures with the linear 
functional induced on the space of continuous functions by integration). To show this 
let us begin by noticing that, since $\mathrm{I}$ is zero on $E$, for every 
$\zeta \in C (\clo)$ one has \esp 
$$\langle D \mathrm{I}, \zeta \rangle = \langle \mathrm{I}, D \zeta \rangle =
\langle \mathrm{I} ,\Delta \zeta + \zeta \rangle = \langle \mathrm{I} , \zeta
\rangle,$$
that is, $\mathrm{I}$ is invariant by the (dual operator to the) diffusion. 
Suppose that the Hahn decomposition of $\mathrm{I}$ may be expressed in the form \esp 
$\mathrm{I} = \alpha \nu_1 - \beta \nu_2,$ \esp where $\nu_1$ and $\nu_2$ are probability 
measures with disjoint support, $\alpha > 0$, and $\beta > 0$. In this case, the equality 
$D\mathrm{I} = \mathrm{I}$ and the uniqueness of the Hahn decomposition of $D\mathrm{I}$ 
show that $\nu_1$ and $\nu_2$ are also invariant by the diffusion. Consequently, 
\thinspace $\nu_1 = \nu_2 = \mu$, \thinspace which contradicts the fact that the 
supports of $\nu_1$ and $\nu_2$ are disjoint. The functional $\mathrm{I}$ may then 
be expressed in the form $\mathrm{I} = c \esp\! \nu$ for some probability measure $\nu$ 
on the circle, and the equality $\mathrm{I} = D \mathrm{I}$ implies that $\nu = \mu$.

Now notice that, since
$$0 > \mathrm{I}(\psi) = c \esp \mu(\psi) = c \esp \! \int_{\clo} \psi(x) \esp d\mu(x),$$
hypothesis (\ref{suppp}) implies that $\mu(\psi) > 0$ and $c < 0$. Nevertheless, since  
the constant function $1$ belongs to $C_+$, we must have \esp $c=\mathrm{I}(1) \geq 0$.
\esp This contradiction concludes the proof. $\hfill\square$

\vspace{0.3cm}

Coming back to the proof of Proposition \ref{neg} we begin by noticing 
that, adding a constant $c_n$ to each $\zeta_n$ if necessary, we 
may suppose that the mean of $\exp(\zeta_n)$ equals $1$ for every 
$n \!\in\! \mathbb{N}$. Let us consider the probability 
measures $\nu_n$ on the circle defined by 
$$\frac{d \nu_n (s)}{d Leb} =  \exp \big( \zeta_n (s) \big),$$
and let us fix a subsequence $\nu_{n_i}$ converging to a probability measure 
$\nu$ on $\clo$. We will show that this measure $\nu$ is invariant by $\Gamma$.

Let us first check that $\nu$ is stationary. 
To do this, notice that if we denote by 
$Jac_n(g)$ the Radon-Nikodym 
\index{Radon-Nikodym derivative}
derivative of $g \!\in\! \Gamma$ with respect to $\nu_n$, 
then relation (\ref{casi-cociclo}) implies that, for every $x \in \clo$, the 
value of 
$$\int_{\Gamma} \log \big( Jac_{n} (g) (x) \big) \thinspace dp (g)$$
equals
$$\int_{\Gamma} \log \big( g'(x) \big) \thinspace dp (g) +
\int_{\Gamma} \big[ \zeta_n(g(x)) - \zeta_n(x) \big] \thinspace dp (g) 
= \psi(x) + \Delta \zeta_n (x) \geq -\frac{1}{n}.$$
Now notice that, since the dual of the diffusion operator acts continuously on the 
space of probability measures on the circle, the sequence of measures $D \nu_{n_i}$ 
converges in the weak-* topology to the measure $D \nu$. 
\index{weak-* topology}
Moreover, the diffusion of 
$\nu_n$ is a measure which is absolutely continuous with respect to $\nu_n$, and 
whose density may be written in the form 
$$\frac{d \thinspace D \nu_n(x)}{d \thinspace \nu_n (x)} =
\int_{\Gamma} Jac_n(g^{-1}) (x) \thinspace dp (g) =
\int_{\Gamma} Jac_n(g)(x) \thinspace dp (g).$$
By the concavity of the logarithm, we have 
$$\frac{d \thinspace D \nu_{n_i}(x)}{d \thinspace \nu_{n_i} (x)} \geq
\exp \left( \int_{\Gamma} \log \big( Jac_{n_i} (g) (x)\big)
\thinspace dp (g) \right ) \geq \exp(-1/n_i),$$
that is,\thinspace $D \nu_{n_i} \geq \exp(-1/n_i) \thinspace \nu_{n_i}$. \thinspace
Passing to the limit we obtain $D \nu \geq \nu$, and since both $\nu$ and $D \nu$ 
have total mass 1, they must be equal. Hence, $\nu$ is stationary.

We may now complete the proof of the invariance of $\nu$. For this, fix an interval $J$ 
such that $\nu(J) > 0$, and consider the function $\psi_{n,J}\!: \Gamma\!\rightarrow ]0,1]$ 
defined by 
\thinspace $\psi_{n,J}(g) = \nu_n \big( g (J) \big)$. \thinspace
From 
$$\Delta \! \log (\psi_{n,J}) (id) =\! \int_{\Gamma} \!\log \!
\Big( \frac{\nu_{n} (g(J))}{\nu_n(J)} \Big) dp(g)
=\! \int_{\Gamma} \! \log \Big(\! \int_J \! Jac_{n}(g)(x) \thinspace 
\frac{d \nu_n(x)}{\nu_n(J)}\Big) dp(g)$$
one deduces that  
\begin{eqnarray*}
\Delta \log (\psi_{n,J}) (id) \!
&\geq& \int_{\Gamma} \Big( \int_J \log \big( Jac_{n}(g)(x) \big) \thinspace
\frac{d\nu_n(x)}{\nu_n(J)} \Big) \thinspace dp(g)\\
&=& \int_J \Big( \int_{\Gamma} \log \big( Jac_{n}(g)(x) \big) \thinspace dp(g) \Big)
\thinspace \frac{d\nu_n(x)}{\nu_n(J)} \esp \esp \geq \esp \esp -\esp\frac{1}{n}.
\end{eqnarray*}
Notice that this is true for every interval $J$ satisfying $\nu(J)>0$. Due to the 
relation \esp $\psi_{n,J}(gf)= \psi_{n,f(J)}(g),$ \esp this implies that the Laplacian 
of $\log(\psi_{n,J})$ is bounded from below by $-1/n$ at every element of $\Gamma$. 
Therefore, if $\psi_J$ is the limit of the sequence of the functions $\psi_{n_i,J}$, 
that is, $\psi_{J}(g) = \nu \big(g(J) \big)$, \esp then the function $\log (\psi_{J})$ 
is super-harmonic (\textit{i.e.}, its Laplacian is non-negative). On the other hand, 
since $p$ is symmetric, $\psi_{J}$ is harmonic. As a consequence, for every element 
$f \!\in\! \Gamma$ the inequalities below are forced to be equalities:
\[ \log (\psi_{J}) (f) \leq \int_{\Gamma} \log (\psi_{J}) (gf) \thinspace dp (g)
\leq \log \Big( \int_{\Gamma} \psi_{J}(gf)\thinspace d\mu(g)\Big)=\log (\psi_{J})
(f).\]
This implies that the function $\psi_J$ is constant. Now since the interval 
satisfying $\nu(J)\!>\!0$ was arbitrary, we deduce that the measure $\nu$ is 
invariant by all of the elements in $\Gamma$. Proposition \ref{neg} is then proved. 
\index{Sacksteder!theorem|)}


\section{Duminy's First Theorem: on the Existence of Exceptional Minimal 
Sets\index{exceptional minimal set}}

\subsection{The statement of the result}
\label{pres-dum}

\index{Duminy!first theorem|(}

\hspace{0.45cm} In general, the subgroups of non discrete groups which are generated by elements 
near the identity behave nicely. A classical result well illustrating this phenomenon 
is the \textit{Zazenh\"aus Lemma}: In every Lie group \index{Lie group}
there exists a neighborhood of the identity such that every discrete group generated 
by elements inside is nilpotent. \index{group!nilpotent} \index{Zazenh\"aus Lemma}

Another well-known result in the same direction 
is \textit{J\o rgensen's inequality} \cite{Be,Ka}: 
\index{J\o rgensen's inequality} 
If $f,g$ are elements in $\mathrm{PSL}(2,\R)$ 
generating a Fuchsian group of second kind, then 
$$|\mathrm{tr}^2(f)-4| + |\mathrm{tr}([f,g])-2| > 1,$$
where $\mathrm{tr}$ denotes the trace of the (equivalence class of the) corresponding 
matrix. In the same direction, Marden \cite{marden} showed the existence of  
\index{Marden's theorem}
a universal constant $\varepsilon_0 \!>\! 0$ such that for every non-elementary Fuchsian group 
$\Gamma$, every system of generators $\mathcal{G}$ of $\Gamma$, and every point $P$ in the 
Poincar\'e disk, \index{Poincar\'e!disk}
there exists $g \!\in\! \mathcal{G}$ satisfying \esp 
$dist (P,g(P)) \geq \varepsilon_0$. \esp This result was extended by Margulis to discrete 
\index{Margulis!inequality}
isometry groups of hyperbolic spaces of arbitrary dimension (\textit{Margulis' inequality} 
still holds for more general spaces of negative curvature: see \cite{BGS}).

For the case of groups of circle diffeomorphisms, there is a large variety of results in 
the same spirit, specially in the real-analytic case \cite{Bon,DFF,Gh3,Re1,Re2}. One of 
the main motivations is a beautiful theorem obtained by Duminy at the end of the 
seventies. Unfortunately, Duminy never published the proof of his result.

\vspace{0.2cm}

\begin{thm} {\em There exists a universal constant $V_0 > 0$ satisfying the following 
property: If $\Gamma$ is a subgroup of $\mathrm{Diff}_+^{1+\mathrm{bv}}(\clo)$ 
generated by a (non-necessarily finite) family $\mathcal{G}$ of diffeomorphisms such 
that at least one of them has finitely many periodic points and $V(g) < V_0$ for 
every $g \in \mathcal{G}$, then $\Gamma$ does not admit an exceptional minimal set.}
\label{thduminy}
\end{thm}

\vspace{0.2cm}

Notice that the condition $V(g) \!<\! V_0$ means that the elements of $\mathcal{G}$ are 
\index{rotation!group}
close to rotations: the equality $V(g) = 0$ is satisfied if and only if $g$ is a rotation. 
Concerning the hypothesis of existence of a generator with isolated periodic points, 
let us point out that it is ``generically'' satisfied \cite{MS}. It is very plausible 
that the theorem is still true without this assumption. This is known for instance in 
the real-analytic case: see Exercise \ref{dum-anal}.\\

Recall that if $g$ is a $\ce^{1+\mathrm{bv}}$ diffeomorphism and $I$ is an interval 
contained in the domain of definition of $g$, then $V(g;I)$ denotes the total variation 
of the logarithm of the derivative of the restriction of $g$ to $I$. Notice that 
\index{variation of a function!total variation}
\esp $V(g^{-1};I) = V(g;g^{-1}(I)).$ \esp 
This relation implies in particular that $V(g) = V(g^{-1})$ for a circle diffeomorphism $g$. Due 
to this, there is no loss of generality for the proof of Theorem \ref{thduminy} if we assume that 
$\mathcal{G}$ is symmetric. Moreover, for every circle diffeomorphism $g$ there must exist a 
point $p \in \clo$ such that $g'(p) = 1$, which allows to conclude that 
\begin{equation}
\inf_{x \in \clo} f'(x) \geq e^{-V(f)}, \qquad 
\qquad \sup_{x \in \clo} f'(x) \leq e^{V(f)}.
\label{minimo-maximo}
\end{equation}
From (\ref{vardos}) and  (\ref{minimo-maximo}) one deduces that if 
$g \!\in\! \dos$ satisfies $|g''(x)| < \delta$ for every $x \in \clo$, 
then $V(g) < 2 \pi \delta / (1 \!-\! 2 \pi \delta)$. This implies that, 
for some universal constant $\delta_0$, the conclusion of Theorem \ref{thduminy} applies  
to subgroups of $\dos$ generated by elements $g$ satisfying $|g''(x)| < \delta_0$ for every 
$x \!\in\! \clo$ (at least when one of the generators has finitely many periodic points). In 
the same way, the reader should have no problem in adapting the arguments below to groups 
of piecewise affine circle homeomorphisms generated by elements near rotations.

Before entering into the proof of Theorem \ref{thduminy}, we record two 
important properties. On the one hand, for every pair of intervals $I_1,I_2$ 
inside an interval $I$ contained in the domain of definition of $g$ one has  
\begin{equation}
\frac{|g(I_2)|}{|I_2|} \esp e^{-V(g;I)} \leq \frac{|g(I_1)|}{|I_1|}
\leq \frac{|g(I_2)|}{|I_2|} \esp e^{V(g;I)}.
\label{tresduminy}
\end{equation}
On the other hand, for every pair of $\ce^{1+\mathrm{bv}}$ diffeomorphisms 
$f,g$, the cocycle relation (\ref{cociclolog}) implies the inequality  
\begin{equation}
V(f \circ g; I ) \leq V(g;I) + V(f;g(I)).
\label{unoduminy}
\end{equation}


\subsection{An expanding first-return map}
\index{expanding!map|(}

\hspace{0.45cm} The proof of Theorem \ref{thduminy} is based on Duminy's work on 
semi-exceptional leaves of codimension-one foliations which will be studied later. 
The lemma below appears in one of his unpublished manuscripts (see however \cite{Na3}).

\vspace{0.15cm}

\begin{lem} {\em Let $a,b,b',c$ be points in the line such that 
$a\!<\!c\!<\!b$ and $a\!<\!b'\!<\!b.$ Let $f\!:[a,c] \rightarrow [a,b]$
and $g\!: [c,b] \rightarrow [a,b']$ be $\mathrm{C}^{1+\mathrm{bv}}$ 
diffeomorphisms such that $f(x) > x$ for every $x \neq a$, and $g(x) < x$ 
for every  $x$ (see Figure {\em 17}). Suppose that for some $m,n$ in 
$\mathbb{N}$ and $[u,v] \subset [c,b]$, the map $\bar{H} = g^{-n} \circ f^{-m}$ 
is defined on the whole interval $[u,v]$. Then one has the inequality}
\begin{equation}
\frac{\bar{H}(v)-\bar{H}(u)}{v-u} \leq
\frac{\bar{H}(v)-f^{-1}(\bar{H}(v))}{v-f^{-1}(v)}
\Big( 1- \frac{1}{\sup\limits_{x \in [a,c]} f'(x) } \Big) e^{V(f;[a,c])+V(g;[c,b])}.
\label{cin}
\end{equation}
\label{duminy}
\end{lem}

\vspace{-0.3cm}

\noindent{\bf Proof.} To simplify, let us denote $\bar{f} \!\!=\!\! f^{-1}$ and 
$\bar{g} \!=\! g^{-1}$ (see Figure 18). Recall that $V(\bar{f};[a,b]) \!=\! V(f;[a,c])$ 
and $V(\bar{g};[a,b']) \!=\! V(g;[c,b]).$ We then have to show the inequality
\begin{equation}
\frac{\bar{H}(v)-\bar{H}(u)}{v-u} \leq \frac{\bar{H}(v)-\bar{f}(\bar{H}(v))}
{v-\bar{f}(v)}
\left( 1- \inf\limits_{x \in [a,b]} \bar{f}'(x) \right)
e^{V(\bar{f};[a,b])+V(\bar{g};[a,b'])}.
\label{inversa}
\end{equation}
To do this, first remark that 
$$V \big( \bar{f}^m;[\bar{f}(v),v] \big) \leq
\sum\limits_{k=0}^{m-1}V \big( \bar{f};[\bar{f}^{k-1}(v),\bar{f}^{k}(v)] \big)
\leq V \big( \bar{f};[a,b] \big).$$
Since $\bar{f}(v) \leq c \leq u < v$, we have 
$$\frac{\bar{f}^m(v)-\bar{f}^m(u)}{v-u} 
\leq \frac{\bar{f}^m(v)-\bar{f}^{m+1}(v)}{v-\bar{f}(v)} \hspace{0.1cm}
e^{V(\bar{f}^m;[\bar{f}(v),v])} 
\leq \frac{\bar{f}^m(v)-\bar{f}^{m+1}(v)}{v-\bar{f}(v)}
\hspace{0.1cm} e^{V(\bar{f};[a,b])}.$$
Moreover,
$$\bar{f}^{m+1}(v)-a 
= \bar{f}(\bar{f}^m(v))-\bar{f}(a) \geq (\bar{f}^m(v)-a)
\inf\limits_{x \in [a,b]} \bar{f}'(x),$$
from where one obtains
$$\bar{f}^m(v) -\bar{f}^{m+1}(v) \leq (\bar{f}^m(v)-a)
\left( 1-\inf\limits_{x \in [a,b]} \bar{f}'(x) \right)\!.$$
We then deduce that
\begin{equation}
\frac{\bar{f}^m(v)-\bar{f}^m(u)}{v-u}
\leq \frac{\bar{f}^m(v)-a}{v-\bar{f}(v)} \left( 1-\inf\limits_{x \in [a,b]}
\bar{f}'(x) \right)
e^{V(\bar{f};[a,b])}.\\
\label{ses}
\end{equation}
Since \esp $a \!\leq\! \bar{f}^m(u) \!<\! \bar{f}^m(v) \!\leq\! c\! =\! \bar{g}(a)$ 
\esp and \esp $\bar{g}^n$ \esp is well-defined on the interval \esp 
$[\bar{f}^m(u),\bar{f}^m(v)]$, \esp analogous arguments show that 
\begin{equation}
\frac{\bar{g}^n(\bar{f}^m(v))-\bar{g}^n(\bar{f}^m(u))}{\bar{f}^m(v)-\bar{f}^m(u)}
\leq \frac{\bar{g}^n(\bar{f}^m(v))-\bar{g}^n(a)}{\bar{f}^m(v)-a}
\hspace{0.1cm} e^{V(\bar{g};[a,b'])}.
\label{sie}
\end{equation}
From ($\ref{ses}$) and ($\ref{sie}$) one deduces that \esp 
$(\bar{H}(v)-\bar{H}(u)) / (v-u)$ \esp is less than or equal to  
$$\frac{\bar{H}(v)-\bar{f}(\bar{H}(v))}{v-\bar{f}(v)} \esp  
\frac{\bar{H}(v) - \bar{g}(a)}{\bar{H}(v) - \bar{f}(\bar{H}(v))}
\left( 1 \!-\! \inf\limits_{x \in [a,b]} \bar{f}'(x) \right)
e^{V(\bar{f};[a,b])+V(\bar{g};[a,b'])}.$$
Inequality (\ref{inversa}) follows from this by using the 
fact that $\bar{f}(\bar{H}(v)) \leq c = \bar{g}(a)$. $\hfill\square$

\vspace{0.5cm}


\beginpicture

\setcoordinatesystem units <0.72cm,0.72cm>

\putrule from 0 0 to 5.1 0

\putrule from 0 0  to 0 5.1

\putrule from 5.1 0 to 5.1 5.1

\putrule from 1.723 0 to 1.723 5.1

\putrule from 0 5.1 to 5.1 5.1

\putrule from 0 5.1 to 5.1 5.1

\put{$a$} at 0 -0.3


\put{$b$} at 5.1 -0.3

\put{$c$} at 1.723 -0.3

\put{$c$} at -0.3 1.723

\put{$f$} at 0.7 3

\put{$g$} at 3.95 2.5

\put{Figure 17} at 2.5 -1

\plot
1.723 0
1.723 0.4 /

\plot
0  0
5.1  5.1 /

\setquadratic 

\plot 
1.723 0
2.9 2 
5.1 3.58 /

\plot
0 0
1.35 4 
1.723 5.1 /

\setlinear 

\put{$b'$} at  -0.3  3.58


\putrule from 10 0 to 15.1 0

\putrule from 10 0  to 10 5.1

\putrule from 15.1 0 to 15.1 5.1

\putrule from 10 1.723 to 15.1 1.723

\putrule from 10 5.1 to 15.1 5.1

\put{$a$} at 10 -0.3

\put{$b$} at 15.1 -0.3

\put{$c$} at 11.723 -0.3

\put{$c$} at 9.7 1.723

\put{$\bar{f}$} at 14.4 1

\put{$\bar{g}$} at 12.15 3.654

\put{Figure 18} at 12.5 -1

\plot
10 1.723
10.4 1.723 /

\plot
10  0
15.1  5.1 /

\setquadratic 

\plot 
10 1.723 
12 2.9 
13.58 5.1 /

\plot
10 0 
14 1.35 
15.1 1.723 /

\setlinear 

\put{$b'$} at 13.58 -0.3

\setdots
\putrule from 13.58 0 to 13.58 5.1
\putrule from 0 3.58 to 5.1 3.58
\putrule from 0 1.723 to 1.723 1.723
\putrule from 11.723 0 to 11.723 1.723

\put{} at -2 0


\endpicture


\vspace{0.5cm}

Given two maps $f,g$ as in Lemma \ref{duminy}, for each $x \!\in ]c,b]$ there exists a 
positive integer $n \!=\! n(x)$ such that $g^{n-1}(x)\!\in ]c,b]$ and $g^n(x)\!\in ]a,c]$. 
Analogously, for $y \!\in ]a,c]$ there exists $m = m(y)$ such that $f^{m-1}(y) \!\in ]a,c]$ 
and $f^{m}(y) \!\in ]c,b]$. We define the {\bf \textit{first-return map}}
$H \!: \esp ]c,b] \rightarrow ]c,b]$ by 
$$H (x) = f^{m(g^{n(x)}(x))} \circ g^{n(x)}(x).$$
Notice that the set of discontinuity points of this map is 
$\{ g^{-1}(f^{-j}(c)) \!\! : j \!\in\! \mathbb{N}\}$, and for every 
$\varepsilon > 0$ the intersection of this set with 
$[c+ \varepsilon,b]$ is finite. We let 
\begin{equation}
C(f) = \frac{\sup_{x \in [c,b]} (x-\bar{f}(x))}{\inf_{x \in [c,b]}(x-\bar{f}(x))}.
\label{graneme}
\end{equation}
\index{first-return map}

\begin{prop} {\em Under the hypothesis of Lemma {\em \ref{duminy}}, suppose moreover that}
\begin{equation}
\Big( 1- \frac{1}{\sup\limits_{x \in [a,c]} f'(x) } \Big)
\esp \! e^{ V(f;[a,c]) + V(g;[c,b]) }  < 1.
\label{verifie}
\end{equation}
{\em Then for every $\kappa \!>\! 1$ there exists $N \!\in\! \mathbb{N}$ such that 
$(H^{N})'(x) \!>\! \kappa$ for every $x$ at which $H^{N}$ is differentiable.}
\label{retorno}
\end{prop}

\noindent{\bf Proof.} Let $N \in \mathbb{N}$ be such that 
\begin{equation}
\Big( 1 - \frac{1}{\sup\limits_{x \in [a,c]} f'(x) } \Big)^{N}
e^{N(V(f;[a,c])+V(g;[c,b]))} < \frac{1}{\kappa \cdot C(f)}.
\label{contraria}
\end{equation}
We claim that each branch of the map $H^{N}$ is $\kappa$-expanding. To show this, let 
$x$ be a differentiability point of $H^{N}$. 
Let us fix an interval $[x-\varepsilon,x]$ contained in this component. Denote   
$[u_0,v_0] = [x-\varepsilon,x]$, and for $j \!\in\! \{ 1, \ldots, N-1 \}$ 
denote $[u_{j},v_{j}] \!=\! H^{j}([u_0,v_0])$. Inequality 
(\ref{cin}) applied to $\bar{H} \!=\! H^{-1}$ gives
$$\frac{v_j -u_j}{v_{j+1}-u_{j+1}} \leq \frac{v_j -
f^{-1}(v_j)}{v_{j+1}-f^{-1}(v_{j+1})}
\hspace{0.1cm} \Big( 1- \frac{1}{\sup\limits_{x \in [a,c]} f'(x) } \Big)
e^{ V(f;[a,c]) + V(g;[c,b]) }.$$
Taking the product from  \hspace{0.05cm} $j=0$ \hspace{0.05cm} to \hspace{0.05cm} 
$j=N-1$ \hspace{0.05cm} we obtain
$$\frac{v_0 -u_0}{v_N -u_N} \leq \frac{v_0 -f^{-1}(v_0)}{v_N -f^{-1}(v_N)}
\hspace{0.1cm}
\Big( 1- \frac{1}{\sup\limits_{x \in [a,c]} f'(x) } \Big)^N
e^{ N( V(f;[a,c]) + V(g;[c,b]) )}.$$
From this and (\ref{contraria}), one deduces that
$$\frac{(H^{N}(x)-H^{N}(x-\varepsilon))}{\varepsilon}  > \kappa.$$  
Since this inequality holds for every $\varepsilon > 0$ small enough, 
this allows to conclude that $(H^{N})'(x) \!>\! \kappa$. $\hfill\square$

\vspace{0.25cm}

\begin{prop} {\em If $f$ and $g$ satisfy the hypothesis of Proposition {\em \ref{retorno}}, 
then all the orbits of the pseudo-group \index{pseudo-group}
generated by these maps are dense in $]a,b[$.}
\label{retorno-dos}
\end{prop}

\noindent{\bf Proof.} Suppose that there exists an orbit which is non-dense in 
$]a,b[$. Then it is easy to see that the points $a$, $c$, and $b$, belong to the 
closure of this orbit. Among the connected components of the complementary set of 
this closure contained in $]c,b[$, let us choose one, say $]u,v[$, having maximal 
length. Then the proposition above shows that, for $N$ large enough, $H^N(]u,v[)$ 
has larger length, which is a contradiction. $\hfill\square$

\vspace{0.5cm}

At first glance, the above estimates may seem too technical. To understand them 
better, one may consider the particular case where the maps $f$ and $g$ are affine, 
say $f(x) = \lambda x$, with $\lambda > 1$, and $g(x) = \eta(x- 1 / \lambda)$.
In this case, for every $n \!\in\! \mathbb{N}$ the restriction of the first 
return map $H \!: [1 / \lambda,1] \rightarrow [1 / \lambda,1]$ to the interval 
$g^{-1} \big( [1 / \lambda^{n+1}, 1 / \lambda^{n}] \big)$ is given by $H(x) = f^n \circ g (x)$. 
Since for every $x \in g^{-1} \big( [1/\lambda^{n+1},1/\lambda^{n}] \big)$ one has 
$$\frac{1}{\lambda^{n}} \geq \eta \left( 1 - \frac{1}{\lambda} \right)
\geq \frac{1}{\lambda^{n+1}},$$
this gives
$$H'(x) = \eta \lambda^{n} \geq \frac{1}{\lambda - 1} = 
\frac{1}{\lambda(1 - 1 / \lambda)},$$
and since $C(f) = \lambda$, this inequality may be written in the form 
$$H'(x) \geq \frac{1}{C(f)
\Big( 1 - \frac{1}{\sup\limits_{x \in [0,1/\lambda]} f'(x)} \Big)}.$$
Now the values of $V(f)$ and $V(g)$ being equal to zero, the similarity between 
the last inequality and those appearing throughout this section becomes evident.

\begin{small} \begin{ejer}
State and show a proposition making formal the following claim: If the derivative of $f$ is 
near 1 and the total variation $V(f;[a,c])$ is small, then the constant $C(f)$ defined by 
(\ref{graneme}) is near 1.
\end{ejer} \end{small}

\begin{small} \begin{ejer} By modifying the example in \S \ref{sec-minimal} slightly, 
prove that Proposition \ref{retorno-dos} is optimal when imposing \textit{a priori} 
the condition $V(g)=0$; show that, in this case, the critical parameter corresponds 
to (an integer multiple of) $\log(2)$. Moreover, prove   
that the proposition is optimal under the condition $V(f) = V(g)$; show 
that, in this case, the critical parameter is (an integer multiple of) the 
logarithm of the golden number (see \cite{Na3} for more details on this).
\label{optimal}
\end{ejer} \end{small}

\index{expanding!map|)}


\subsection{Proof of the theorem}

\hspace{0.45cm} In order to 
illustrate the idea of the proof of Theorem \ref{thduminy}, suppose that a group 
acts on the circle preserving an exceptional minimal set $\Lambda$, and two generators 
$f,g$ are as in Figure 17 over an interval $[a,b]$ in $\clo$ whose interior intersects  
$\Lambda$. In this case, by an argument similar to that of Corollary \ref{retorno-dos}, 
we conclude that the opposite inequality to (\ref{verifie}) must be satisfied, that is,
\begin{equation}
\Big( 1- \frac{1}{\sup\limits_{x \in [a,c]} f'(x) } \Big)
e^{ V(f;[a,c]) + V(g;[c,b]) }  \geq 1.
\label{necesaria}
\end{equation}
By (\ref{minimo-maximo}), this inequality implies that 
$$( 1 - e^{-V(f)}) \hspace{0.1cm}  e^{V(f;[a,c]) + V(g;[c,b])} \geq 1.$$
If $V(f)$ and $V(g)$ are strictly smaller than some positive constant $V_0$, 
then this gives \esp $(1 - e^{-V_0}) \hspace{0.1cm} e^{2 V_0} > 1,$ 
\esp that is, 
$$e^{2V_0} - e^{V_0} - 1 > 0,$$ 
and hence \esp 
$V_0 > \log \big( (\sqrt{5}+1)/2 \big)$. \esp (The occurrence of the golden number here 
is not mysterious; actually, it seems to be related to the optimal constant for the 
theorem: see Exercise \ref{optimal}.)

\vspace{0.2cm}

Unfortunately, the proof for the general case involves a major technical problem:  
it is not always the case that there are generators for which one may directly apply the 
preceding argument. For instance, the modular group acts on $\clo$ admitting an exceptional 
minimal set, although its generators have finite order... This is the main reason why 
we assume the hypothesis that one of the generators has isolated periodic points.

\vspace{0.35cm}

\noindent{\bf Proof of Theorem \ref{thduminy}.} Suppose that $\Gamma$ admits an exceptional 
minimal set $\Lambda$. By Denjoy's theorem, 
\index{Denjoy!theorem} the set of periodic points of each element in 
$\Gamma$ is non-empty. By hypothesis, there exists a generator $g \!\in\! \mathcal{G}$
whose periodic points are isolated. Let us denote by $\mathrm{Per}(g)$ the set of 
these points, and let $\mathrm{P}(g) = \mathrm{Per}(g) \cap \Lambda$. Notice that 
$\mathrm{P}(g)$ is non-empty. Indeed, if the periodic points of $g$ have 
order $k$ and $p \in \Lambda$ is not a fixed point of $g^k$, then 
\esp $\lim_{j \rightarrow \infty} g^{jk}(p)$ \esp and \esp 
$\lim_{j \rightarrow \infty} g^{-jk}(p)$ \esp are fixed points of  
$g^k$ contained in $\Lambda$. 

We now claim that there exist 
$p \in \mathrm{P}(g)$ and $f \in \mathcal{G}$ such that 
$f(p) \in \mathrm{S}^1 \setminus \mathrm{P}(g).$ If not, 
the finite set $\mathrm{P}(g)$ would be invariant by 
$\Gamma$, thus contradicting the minimality of $\Lambda$.\\

Let $G \!=\! g^k \in \Gamma$, and let us denote by $u$ and $v$ the periodic points of $g$ 
immediately to the left and to the right of $f(p)$, 
respectively. The map $F= f \circ g^{k} \circ f^{-1}$ has a fixed point in $[u,v]$, namely $f(p)$. 
Let $a$ be the fixed point of this map to the left of $v$, and let $q$ be the fixed point to the 
right of $a$. Replacing $G$ by $G^{-1}$ and/or $F$ by $F^{-1}$ if necessary, we may suppose that 
\hspace{0.02cm} $G(x) < x$ \hspace{0.02cm} and \hspace{0.02cm} $F(x) > x$ \hspace{0.02cm} 
for every $x \!\in ]a,v[$ \esp (see Figure 19).\\

We now claim that, if $V_0$ is small enough, then the point $b=F(G^{-1}(a))$ 
belongs to the interval $]a,v[$. To show this we first notice that
$$V \big( F^{-1};[a,q] \big) = V \big( F;[a,q] \big) \leq \sum\limits_{j=0}^{k-1} 
V \big( f \circ g \circ f^{-1}; f \circ g^{j} \circ f^{-1} ([a,q]) \big)
\leq V(f \circ g \circ f^{-1}) < 3 V_0.$$
In the same way one obtains 
\esp $V \big( G^{-1};[u,v] \big) = V \big( G;[u,v] \big) < V_0.$ 
\esp Let $x_0 \! \in ]a,v[$ and $y_0 \! \in ]a,q[$ be such that 
$$(F^{-1})'(x_0) = \frac{F^{-1}(v)-a}{v-a}, 
\quad  \quad (F^{-1})'(y_0) = 1.$$
Clearly, we have 
\esp $\big| \mathrm{\mathrm{log}}(F^{-1})'(y_0) - \mathrm{log} (F^{-1})'(x_0) \big|
\leq V \big( F^{-1};[a,q] \big),$ \esp 
and hence 
\begin{equation}
F^{-1}(v) - a > e^{-3V_0}(v-a).
\label{onc}
\end{equation}
By an analogous argument one shows that
\begin{equation}
v - G^{-1}(a) > e^{-V_0}(v-a).
\label{doc}
\end{equation}
If $b$ were not contained in $]a,v[$ then $F^{-1}(v) \leq G^{-1}(a)$, 
and hence, by ($\ref{onc}$) and ($\ref{doc}$),
$$v-a \esp \geq \esp (F^{-1}(v)-a) + (v-G^{-1}(a)) 
\esp > \esp (e^{-V_0}+e^{-3V_0}) \hspace{0.1cm} (v-a).$$
Therefore, \esp $e^{3V_0} - e^{2V_0} > 1,$ \esp which 
is impossible if \esp $V_0 \leq \mathrm{log}(1.46557)$.\\

The elements $F$ and $G$ in $\Gamma$ are thus as in Figure 17 over the 
interval $[a,b]$. Now notice that for every $x \!\in\! [a,q]$ one has 
$$\big| \mathrm{log}(F')(x) - \mathrm{log}(F')(y_0) \big|
\leq V \big( F;[a,q] \big) < 3V_0,$$
and hence \esp $\sup_{x \in [a,q]} F'(x) < e^{3V_0}.$ \esp 
For $c = G^{-1}(a)$, inequality ($\ref{necesaria}$) 
applied to $F$ and $G$ yields 
$$3V_0 + V_0 \esp > \esp V \big( F;[a,c] \big) + V \big( G;[c,b] \big) \esp \geq \esp 
\mathrm{log} \Big( \frac{1}{1 \!-\! \frac{1}{\sup\limits_{x \in [a,c]} F'(x)}} \Big) 
\esp > \esp \mathrm{log} \!\left( \frac{1}{1-e^{-3V_0}} \right).$$
Therefore, \esp 
$e^{4V_0} - e^{V_0} > 1,$ \esp 
which is impossible if \esp $V_0 \leq \mathrm{log}(1.22074)$. $\hfill\square$

\vspace{0.35cm}

\beginpicture

\setcoordinatesystem units <0.75cm,0.75cm>

\putrule from 10 0 to 16 0
\putrule from 10 6 to 16 6
\putrule from 10 0 to 10 6
\putrule from 16 0 to 16 6
\putrule from 11 0 to 11 6
\putrule from 15 0 to 15 6
\putrule from 10 1 to 16 1
\putrule from 10 5 to 16 5

\putrule from 12.236 1 to 12.236 3.24
\putrule from 11 3.24 to 13.24 3.24
\putrule from 13.24 1 to 13.24 3.24

\plot
10 0
16 6 /

\put {$u$}  at 10     -0.3
\put {$a$}  at 11     -0.3
\put {$a$}  at 9.7     1
\put {$b$}  at 13.197 -0.3
\put {$v$}  at 15     -0.3
\put {$q$} at 16     -0.3
\put {$b$}  at 9.7 3.197

\put {$G$} at 12.95  1.4
\put {$F$} at 11.5  2.7

\put {Figure 19} at 13 -0.9


\plot
10 0
10.7071 0.1
11 0.2
11.5811 0.5
12 0.8
12.236 1
12.5495 1.3
12.8482  1.6
13.1622 2
13.3911 2.3
13.6055 2.6
13.8729 3
14.1833 3.5
14.4721 4
14.7434 4.5
15 5 /

\plot
11 1
11.065 1.3
11.16  1.6
11.24  1.8
11.333  2
11.44  2.2
11.56 2.4
11.6933 2.6
12  3
12.36 3.4
12.56  3.6
12.7733  3.8
13 4
13.4933  4.4
13.8983 4.7
14.3333 5
14.96 5.4
15.2933 5.6
16 6 /

\put{} at 3.8 0

\setdots
\putrule from 13.24 0 to 13.24 1
\putrule from 11 2.236 to 12.236 2.236
\putrule from 10 3.24 to 11 3.24

\endpicture


\vspace{0.23cm}

\begin{small} \begin{ejer} Prove that for every group of circle homeomorphisms which is not semiconjugate 
to a group of rotations, and for every system of generators $\mathcal{G}$, there exists an element in the 
ball $B_{\mathcal{G}}(4)$ of radius $4$ which is not semiconjugate to a rotation. Conclude 
that Duminy's theorem holds for groups of real-analytic diffeomorphisms generated by elements $h$ 
satisfying $V(h) < V_0$ for some small $V_0\!>\!0$, even in the case where the hypothesis of existence 
of a generator with isolated periodic points is not satisfied (see \cite{Na3} for more details on this).
\label{dum-anal}
\end{ejer} 

\begin{ejer} Prove that for every $g \!\in\! \mathrm{PSL}(2,\mathbb{R})$ one has 
\begin{equation}
V(g) = 4 \esp dist (O, g(O)),
\label{dum-yam}
\end{equation}
where $O$ denotes the point $(0,0)$ in the Poincar\'e 
disk, and \esp $dist$ \esp stands for the hyperbolic distance.  
Conclude that, for subgroups of $\mathrm{PSL}(2,\mathbb{R})$, the claim of Duminy's theorem 
\index{Margulis!inequality}
follows from either Marden's theorem or Margulis' inequality discussed at the beginning of 
\S \ref{pres-dum}.

\noindent\underbar{Remark.} By successive applications of J\o rgensen's inequality, 
\index{J\o rgensen's inequality}\index{Marden's theorem}and 
using some results on the classification of Fuchsian groups, Yamada computed the value of the 
best constant for Marden's theorem. Thanks to (\ref{dum-yam}), a careful reading of \cite{yamada} 
then allows showing that the best constant $V_0$ of Duminy's theorem for subgroups of 
$\mathrm{PSL}(2,\mathbb{R})$ is $4 \log(5/3)$, and this value is critical for the actions 
of $\mathrm{PSL}(2,\mathbb{Z})$ (with respect to its canonical system of generators).
\label{adolfo} \index{Duminy!first theorem|)}
\end{ejer} \end{small}


\section{Duminy's Second Theorem: on the Space of Semi-Exceptional Orbits}
\label{super-du}

\index{Duminy!second theorem|(}

\subsection{The statement of the result}

\hspace{0.45cm} The aim of this section is to present another great result 
by Duminy on the structure ``to the infinity'' of some special orbits of 
\index{pseudo-group} 
pseudo-groups of diffeomorphisms of one-dimensional 
manifolds. Roughly, in class $\ce^{1+\mathrm{Lip}}$, 
semi-exceptional orbits 
are forced to have infinitely many ends. A proof of this result 
appears in \cite{CC2}. In what follows, we will give the 
original (and unpublished) remarkable proof by Duminy.

We begin by recalling some basic notions. Given a compact subset $K$ of $\mathrm{M}$, 
we denote the set of connected components of $\mathrm{M} \setminus K$ by 
$\mathcal{E}_K$, and we endow this set with the discrete topology. 
The {\bf {\em space of ends}} 
\index{space of ends} 
of M is the inverse limit of the spaces $\mathcal{E}_K$, 
where $K$ runs over all compact subsets of $\mathrm{M}$. For 
locally compact separable metric spaces, there is an equivalent, but 
more concrete, definition. Fix a sequence $(K_n)$ of compact subsets of M such that 
$K_n \!\subset\! K_{n+1}$ for every $n \!\in\! \mathbb{N}$ and $\mathrm{M}=\cup_{n \in
\mathbb{N}}K_n$, and consider a sequence of connected components $C_n$ of 
$\mathrm{M} \setminus K_n$ such that $C_{n+1} \! \subset \! C_n$ for every 
$n \in \mathbb{N}$. If $(K_n')$ and $(C_n')$ satisfy the same properties, we say 
that $(C_n)$ is equivalent to $(C_n')$ if for every $n \in \mathbb{N}$ there exists 
$m \in \mathbb{N}$ such that $C_n \subset C_m'$ and $C_n' \subset C_m$. Then the space 
of ends of M identifies with the set of equivalence classes of this relation, 
and it is (compact and) metrizable: after fixing the 
sequence $(K_n)$, the distance between (the ends determined by) $(C_n^1)$ and 
$(C_n^2)$ is $e^{-k}$, where $k$ is the minimal integer for which $C_k^1 \neq C_k^2$.

Let $\Gamma$ be a finitely generated 
pseudo-group of homeomorphisms of some space, and let $p$ be a 
point in this space. Fix a (symmetric) system of generators $\mathcal{G}$ for $\Gamma$. The 
{\bf {\em Cayley graph}}\footnote{Perhaps the appropriate terminology should be
{\em Schreier graph}.} 
\index{Cayley graph}\index{Schreier graph}(with respect 
to $\mathcal{G}$) of the orbit $\Gamma(p)$ is the graph whose vertexes are the points 
in $\Gamma(p)$, so that two vertexes $p'$ and $p''$ are connected by an edge (of 
length 1) if and only if there exists $g \!\in\! \mathcal{G}$ such that $g(p')=p''$. 
We define the {\bf {\em space of ends of the orbit of p}} as the space of ends 
$\mathcal{E}(\Gamma(p))$ 
of the graph $Cay_{\mathcal{G}}(\Gamma(p))$. It is not difficult to check that, 
although the Cayley graph depends on $\mathcal{G}$, the associated space of ends is independent 
of the choice of the system of generators. To simplify, in what follows we will omit the 
reference to $\mathcal{G}$ when this causes no confusion, and for each $n \!\in\! \mathbb{N}$ 
we will denote by $B(p,n)$ the ball of center $p$ and radius $n$ in $Cay (\Gamma(p))$ 
(compare \S \ref{muylargo}).

An orbit of a pseudo-group of homeomorphisms of a one-dimensional manifold is said to 
be {\bf \textit{semi-exceptional}}\index{semi-exceptional orbit} 
if it corresponds to the orbit of a point $p$ in a 
local exceptional\index{local exceptional set} 
set $\Lambda$ so that $p$ is an endpoint of one of the connected 
components of the complement $\Lambda^c$ of $\Lambda$ 
(see \S \ref{chuartz}). \textit{Duminy's 
second theorem} may then be stated as follows.

\vspace{0.05cm}

\begin{thm} {\em Every semi-exceptional orbit of a finitely generated pseudo-group of 
$\ce^{1+\mathrm{Lip}}$ diffeomorphisms of a one-dimensional, compact manifold has 
infinitely many ends.}
\label{superduminy}
\end{thm}


It is unknown whether this result holds for \textit{all} the orbits of local exceptional 
sets $\Lambda$ as above. (This is known in the real-analytic case, according to a result 
due to Hector.) 
\index{Hector} 
To deal with semi-exceptional orbits, the basic idea consists in 
associating, to each point $p$ in such an orbit, the connected component of $\Lambda^c$ 
having $p$ as an endpoint. To be more concrete, we let $\mathcal{G}^* =
\cup_{k \in \mathbb{N}} \mathcal{G}^k$, and to each $\hat{g} = (h_1,\ldots,h_k) \in
\mathcal{G}^*$ we associate the element $g = h_k \cdots h_1$. For each 
connected component $I$ of $\Lambda^c$ on which $g$ is defined, we denote 
$$S(\hat{g} ; I) \esp = \esp |I| + \sum_{i=1}^{k-1} |h_i \cdots h_1(I)|.$$
For any two points $p',p''$ belonging to some semi-exceptional orbit $\Gamma(p)$, we 
denote the set of the elements $\hat{g} \in \mathcal{G}^*$ such that $g(p')=p''$ 
by $\mathcal{G}^*(p',p'')$. Then we consider the connected component 
$I$ of $\Lambda^c$ having $p'$ as endpoint, and we let
$$S(p',p'') \esp = \esp \inf 
\{ S(\hat{g};I) \!: \esp \hat{g} \in \mathcal{G}^*(p',p''), 
\esp I \subset dom(g) \}.$$ 
In the case where there is no $\hat{g} \in \mathcal{G}^*(p',p'')$ 
for which $I \subset dom(g)$, we let \esp $S(p',p'') \!=\! \infty$.

\vspace{0.1cm}

\begin{lem} {\em If $p'$ and $p''$ are points of a semi-exceptional orbit $\Gamma(p)$ 
for which there exists an element $f \!\in\! \Gamma$ such that $f(p') \! \neq \! p'$ 
and $f(p'') \! \neq \! p''$, then the sequences of points $p_n' \!=\! f^n(p')$ and 
$p_n'' \!=\! f^n(p'')$ determine ends of $\Gamma(p)$, and these ends coincide 
if and only if \esp $\lim_{n \rightarrow \infty} S(p_n',p_n'') = 0$.}
\label{dist-lem}
\end{lem}

\noindent{\bf Proof.} Since $f(p') \!\neq\! p'$, the points $f^n(p')$ are two-by-two 
different. Therefore, the sequence $(f^n(p'))$ escapes from the compact sets of 
$Cay(\Gamma(p))$. Since the distance in $Cay (\Gamma(p))$ between $p_n'$ 
and $p_{n+1}'$ is independent of $n$, this easily shows that $(p_n')$
determines an end of the orbit of $p$; the same applies to $(p_{n}'')$.

Now fix a constant $\delta \!>\! 0$ such that, if $h \! \in \! \mathcal{G}$ and 
$x \!\in\! dom(h) \cap \Lambda$, then the $2\delta$-neighborhood of $x$ is contained 
in the domain of $h$. Fix also $m_0 \in \mathbb{N}$ so that $B(p,m_0)$ contains all the 
points in $\Gamma(p)$ which are endpoints of some connected component $\Lambda^c$ of 
length greater than or equal to $\delta$. For each $n \!\in\! \mathbb{N}$ we denote by 
$I_n'$ (resp. $I_n''$) the connected component of $\Lambda^c$ whose closure contains 
$p_n'$ (resp. $p_n''$). From the hypotheses $f(p') \!\neq\! p'$ and $f(p'') \!\neq\! p''$,  
one deduces that the intervals $I_n'$ (resp. $I_n''$) are two-by-two disjoint, and 
hence the sum of their lengths is finite (in particular, $|I_n'|$ and $|I_n''|$ tend 
to 0 as $n$ goes to infinity). Since the space of ends does not depend on the choice 
of the system of generators, without loss of generality we may assume that $f$ 
belongs to $\mathcal{G}$.

Suppose that, for every $n \!\in\! \mathbb{N}$ large enough, every path $\hat{g}$ 
from $p_n'$ to $p_n''$ (\textit{i.e.,} such that $g(p_n') = p_n''$) satisfies 
$I_n' \not\subset dom(g)$. If this is the case, these paths must intersect 
$B(p,m_0)$, from where one deduces that the ends determined by $(p_n')$ 
and $(p_n'')$ are different, and 
$$\lim_{n \rightarrow \infty} S(p_n',p_n'') = \infty.$$

In what follows, suppose that for infinitely many $n \!\in\! \mathbb{N}$ there 
exists $\hat{g} \in \mathcal{G}^*(p_n',p_n'')$ such that $I_n' \subset dom(g)$. 
One then easily checks that the same is true for {\em every} $n$ large enough. 
A concatenation type argument then shows that, for all $m \geq n$ very large,
$$S(p_m',p_m'') \esp < \esp 
S(p_n',p_n'') + \sum_{j=n}^m \big( |I_j'| + |I_j''| \big).$$
This inequality obviously implies that $\big( S(p_n',p_n'') \big)$ 
is a Cauchy sequence, and hence converges (to a finite limit).

Fix an arbitrary $\varepsilon > 0$, and denote the sum of the lengths of the connected 
components of $\Lambda^c$ by $S$. 
Clearly, we may choose a finite family $\mathcal{I}_{\varepsilon}$ 
of these components so that the sum of their lengths is greater than $S\!-\!\varepsilon$. 
Let $m_{\varepsilon} \in \mathbb{N}$ be such that $B(p,m_{\varepsilon})$ contains all the 
points in $\Gamma(p)$ situated in the closure of some interval in $\mathcal{I}_{\varepsilon}$. 
On the one hand, if 
the ends determined by $(p_n')$ and $(p_n'')$ coincide, then for every $n$ large enough 
there exists a path $\hat{g} \in \mathcal{G}^*$ from $p_n'$ to $p_n''$ avoiding 
$B(p,m_{\varepsilon})$. This implies that \esp 
$S(p_n',p_n'') \leq S(\hat{g};I_n') \leq \varepsilon$, \esp 
and since $\varepsilon > 0$ was arbitrary, this shows that \esp $\lim_{n \rightarrow 
\infty} S(p_n',p_n'')=0$. On the other hand, 
if the ends we consider are different, then there exists an 
integer $m_0'$ such that, for every $n \!\in\! \mathbb{N}$ large enough, each path from $p_n'$ 
to $p_n''$ passes through $B(p,m_0')$. Denoting by $\varepsilon_0$ the minimum among the 
lengths of the connected components of $\Lambda^c$ having a vertex in $B(p,m_0')$ as an 
endpoint, we conclude that \esp $S(p_n',p_n'') > \varepsilon_0$ \esp for all $n$ large 
enough. Therefore, \esp $\lim_{n \rightarrow \infty} S(p_n',p_n'') > 0$. $\hfill\square$

\vspace{0.3cm}

To show Theorem \ref{superduminy} we will argue by contradiction. Assuming that the ends 
determined by two sequences as those in the preceding lemma coincide, we will use the 
equivalent condition \esp $\lim_{n \rightarrow \infty} S(p_n',p_n'') \!=\! 0$ \esp 
to find paths $\hat{g}_n$ from $p_n'$ to $p_n''$ with arbitrarily 
small distortion. \index{control of distortion}
We will then conclude the proof by applying Corollary \ref{retorno-dos}  
to $f^{-1}$ and $g_n^{-1}$. However, we immediately point out that, 
throughout the proof, we will need to overcome several nontrivial 
technical difficulties.

\vspace{0.05cm}

\begin{small} \begin{obs} The orbits associated to a free action of $\mathbb{Z}$ 
\index{action!free}
(resp. $\mathbb{Z}^d$, for $d \!\geq\! 2$) 
clearly have two ends (resp. one end). Since there exist free 
non minimal actions of $\mathbb{Z}^d$ by $\ce^{1+\tau}$ circle diffeomorphisms for 
some positive $\tau$ (see \S \ref{ejemplosceuno}), this shows that the 
$\ce^{1+\mathrm{Lip}}$ regularity hypothesis is necessary for  Duminy's second 
theorem. However, we ignore whether the theorem still holds in lower 
differentiability assuming \textit{a priori} that there is no 
invariant probability measure. 
\end{obs} \end{small}


\subsection{A criterion for distinguishing two different ends}
\label{duminy-fines2}

\hspace{0.45cm} To fix notation, in what follows we will consider a semi-exceptional orbit 
$\Gamma(p)$ associated to a finitely generated pseudo-group $\Gamma$ of $\ce^{1+\mathrm{Lip}}$ 
diffeomorphisms. We denote by $\Lambda$ the local exceptional set, we choose a finite and symmetric 
system of generators $\mathcal{G}$ for $\Gamma$, and we denote by $C$ a simultaneous Lipschitz 
\index{Lipschitz!derivative}
constant for the logarithm of the derivative of each element $h \in \mathcal{G}$ on a 
$\delta$-neighborhood of $\Lambda \cap dom(h)$ (where $\delta > 0$).

By Sacksteder's theorem, 
\index{Sacksteder!theorem} 
there exists $f \!\in\! \Gamma$ with a fixed point $a \!\in\! \Lambda$ 
such that $f'(a) < 1$. Since the orbit of $p$ intersects every open interval containing 
points of $\Lambda$, it approaches $a$ at least from one side. In what follows, we will 
assume that $\Gamma(p)$ approaches $a$ by the right, since the other case is analogous. 
Let $b > a$ be very near to $a$ so that
$$0 \esp < \esp \inf_{x \in [a,b]} f'(x) \esp = \esp m(f) \leq M(f) \esp 
= \esp \sup_{x \in [a,b]} f'(x) \esp < \esp 1.$$
Given $\varepsilon > 0$ we may choose \esp $b_{\varepsilon} \!\!\in ]a,b[$ \esp such that 
\esp $C (b_{\varepsilon} - a) / (1 - M(f)) \leq \varepsilon / 2$; \esp thus  
\begin{equation}
V \big( f^n;[a,b_{\varepsilon}] \big) \leq C \sum_{k=0}^{n-1} \big| f^k
([a,b_{\varepsilon}]) \big| \leq C \esp \! (b_{\varepsilon} - a) 
\sum_{k=0}^{n-1} M(f)^{k} \leq \frac{\varepsilon}{2}.
\label{dist-todo-ene}
\end{equation}


\begin{prop} {\em Let $\varepsilon \!>\! 0$ be such that 
$\varepsilon \!\leq\! \frac{1}{2} \log \big( \frac{1}{1-m(f)} \big)$. Suppose that 
there exist connected components $]x,y[$ and $]u,v[$ of the complementary set of 
$\Lambda$ such that:

\vspace{0.1cm}

{\noindent {\em (i)} \esp $a < x < y < u < v < f(b_{\varepsilon})$,}

\vspace{0.1cm}

{\noindent {\em (ii)} \esp $\{x,u \}$ \esp (resp. $\{y,v\}$) 
\esp is contained in $\Gamma(p)$,}

\vspace{0.1cm}

{\noindent {\em (iii)} \esp one has the inequalities}}
$$\frac{f^{-1}(u)-u}{f^{-1}(u)-x} \exp(2\varepsilon) \leq \frac{v-u}{y-x}
\leq \frac{u-a}{x-a} \exp(-2\varepsilon).$$

\vspace{0.1cm}

\noindent{Then the sequences $(f^n(x))$ and $(f^n(u))$ (resp. $(f^n(y))$ 
and $(f^n(v))$) determine two different ends of $\Gamma (p)$.}
\label{distinguir}
\end{prop}

\vspace{0.25cm}

To show this proposition we will need to slightly refine for the 
$\ce^{1+\mathrm{Lip}}$ case the Schwartz' estimates 
\index{Schwartz estimates} from \S \ref{chuartz}.

\vspace{0.1cm}

\begin{lem} {\em Let $I$ be a connected component of the complementary set of $\Lambda$, 
and let $J$ be an interval containing $I$. Suppose that for some $\lambda > 1$ and 
$n \!\in\! \mathbb{N}$ there exists $\hat{g} = (h_1, \ldots, h_n) \in \mathcal{G}^n$ 
such that}
$$S(\hat{g};I) \leq \inf \Big\{ \frac{\log(\lambda)}{\lambda C},
\frac{\delta}{\lambda} \Big\},
\qquad  \qquad \frac{|J|}{|I|} \leq \lambda \exp \big( -\lambda C
S(\hat{g};I) \big).$$
{\em \noindent{Then one has:}}

\vspace{0.1cm}

\noindent{(i) \esp $J \subset dom(g)$,}

\vspace{0.1cm}

\noindent{(ii) \esp $V(g;J) \leq \lambda C S(\hat{g};I)$,}

\vspace{0.1cm}

\noindent{(iii) \esp $|g(J)| \leq \lambda |g(I)|$.}
\label{dum-sc}
\end{lem}

\noindent{\bf Proof.} By induction. For $n\!=\!1$ we have 
$\hat{g} \!=\! g \!\in\! \mathcal{G}$ and $S(\hat{g};I) = |I|$. \esp Property (i) 
follows from the fact that $J$ contains $I$ and \esp 
$|J| \!\leq\! \lambda |I| \!=\! \lambda S(\hat{g};I) \leq \delta$. \esp 
For (ii) notice that \esp $V(g;J) \!\leq\! C |J| \!\leq\! \lambda C S(\hat{g};I)$.
\esp Concerning (iii) notice that, by hypothesis and (\ref{tresduminy}),
$$\frac{|g(J)|}{|g(I)|} \leq \frac{|J|}{|I|} e^{V(g;J)} \leq \lambda.$$

\index{Duminy!estimates}

Suppose now that the claim holds up to $k \!\in\! \mathbb{N}$, and consider 
$\hat{g} \!=\! (h_1,\ldots,h_{k+1})$ in $\mathcal{G}^{k+1}$ satisfying the hypothesis. 
If (i), (ii), and (iii), do hold for $(h_1,\ldots,h_k) \in \mathcal{G}^k$, then:

\vspace{0.1cm}

\noindent{-- property (i) holds for $\hat{g}$, since $h_k \cdots h_1$ is 
defined on $J$, the interval $h_k \cdots h_1 (J)$ intersects $\Lambda$, and}
$$|h_k \cdots h_1(J)| \leq \lambda |h_k \cdots h_1(I)| 
< \lambda S(\hat{g};I) \leq \delta;$$

\vspace{0.1cm}

\noindent{-- property (ii) holds for $\hat{g}$, since by the inductive hypothesis we have}
\begin{eqnarray*}
V(g;J) 
&\leq& V \big( h_{k+1};h_k \cdots h_1(J) \big) + V \big(h_k \cdots h_1;J \big)\\
&\leq& C \big| h_k \cdots h_1(J) \big| + \lambda C S \big( (h_k,\ldots,h_1);I \big)
\esp \esp \leq \esp \esp \lambda C S (\hat{g};I);
\end{eqnarray*}

\noindent{-- property (iii) holds for $\hat{g}$, since 
from (ii) and the hypothesis of the lemma it follows that}
$$\mbox{ } \qquad \qquad \qquad \qquad \hspace{1.25cm}
\frac{|g(J)|}{|g(I)|} \leq \frac{|J|}{|I|} \exp \big( V(g;J) \big) \leq \lambda. 
\hspace{4.25cm} \square$$

\vspace{0.2cm}

Notice that the parameters $S(\hat{g};I)$ and $\lambda$ occur simultaneously in 
the hypothesis of the preceding lemma. In particular, if $S(\hat{g};I)$ is small, 
then the claim of the lemma asserts that the distortion is controlled in a 
neighborhood of $I$ whose length is large with respect to that of $|I|$. 
This clever remark by Duminy will be fundamental in the sequel.

\vspace{0.38cm}

\noindent{\bf Proof of Proposition \ref{distinguir}.} We will give the proof only 
for the case where $\{x,u \}$ is contained in $\Gamma(p)$, since the other case 
is analogous. Suppose for a contradiction that the ends determined by the sequences
$\big( f^n(x) \big)$ and $\big( f^n(u) \big)$ coincide, and let us denote 
$x_n \!=\! f^n(x)$, $y_n \!=\! f^n(y)$, $u_n \!=\! f^n(u)$, and $v_n \!=\! f^n(v)$. 
By Lemma \ref{dist-lem}, $S_n=S(x_n,u_n)$ converges to $0$ as $n$ goes to infinity. 
Hence, by definition, for each $n \!\in\! \mathbb{N}$ there exists $\hat{g}_n \!\in\!
\mathcal{G}^*$ such that \esp $g_n (x_n) \!=\! u_n$ \esp and \esp $S(\hat{g}_n;[x_n,y_n]) 
\leq 2 S_n$. \esp Letting \esp $\lambda_n = 1 / \sqrt{S_n}$ \esp we see that, for $n$ 
large enough, one has \esp $\lambda_n > 1$ \esp and
$$2 S_n \leq \min \Big\{ \frac{\log(\lambda_n)}{\lambda_n C},
\frac{\delta}{\lambda_n} \Big\}.$$
For these integers $n$ let
$$\alpha_n = \lambda_n \exp(-2 \lambda_n C S_n), \quad
x_n' = x_n - \alpha_n (y_n - x_n), \quad y_n' = y_n + \alpha_n (y_n-x_n).$$
Remark that $\alpha_n$ goes to infinity together with $n$. By Lemma 
\ref{dum-sc}, $g_n$ is defined on the whole interval $[x_n',y_n']$, and
\begin{equation}
V \big( g_n;[x_n',y_n'] \big) \leq 2 \lambda_n C S_n.
\label{dist-a-cero}
\end{equation}

\vspace{0.25cm}

\noindent{\underbar{Claim (i).} There exists an integer $N_1$ such that, 
if $n \geq N_1$, then $[a,f^{-1}(u_n)]$ is contained in $[x_n',y_n']$.}

\vspace{0.1cm}

Indeed, by (\ref{dist-todo-ene}),
$$\frac{x_n-a}{y_n-x_n} = \frac{f^n(x)-f^n(a)}{f^n(y)-f^n(x)} \leq
\frac{x-a}{y-x} \exp \big( V(f^n;[a,b_{\varepsilon}]) \big) \leq \frac{x-a}{y-x}
\exp(\varepsilon / 2).$$
Hence, if we choose $n$ so that \esp 
$\alpha_n \geq e^{\varepsilon/2}(x-a)/(y-x)$, \esp
then \esp $(x_n-a)/(y_n-x_n) \leq \alpha_n$, \esp that is, 
\esp \esp $x_n' \leq a$. Analogously, if $n$ is such that \esp
$\alpha_n \geq e^{\varepsilon/2}(f^{-1}(u)-y)/(y-x),$ \esp 
then the inequality 
$$\frac{f^{-1}(u_n) - y_n}{y_n - x_n} \leq \frac{f^{-1}(u) - y}{y - x}
\exp(\varepsilon/2)$$
implies that \esp $y_n' \geq f^{-1}(u_n)$.

\vspace{0.35cm}

\noindent{\underbar{Claim (ii).} There exists an integer $N_2 \geq N_1$ such 
that, if $n \geq N_2$, then $g_n(t) > t$ for all $t \in [a,f^{-1}(u_n)]$.}

\vspace{0.1cm}

To show this claim notice that, from (\ref{dist-todo-ene}) 
and hypothesis (iii) of the proposition it follows that  
$$\frac{f^{-1}(u_n) - u_n}{f^{-1}(u_n) - x_n} \exp(\varepsilon) \leq 
\frac{v_n-u_n}{y_n-x_n} \leq \frac{u_n - a}{x_n - a} \exp(-\varepsilon).$$
From (\ref{dist-a-cero}) one concludes that \esp $V(g_n;[x_n',y_n'])$ \esp 
converges to 0 as $n$ goes to infinity. 
Take $N_2 \geq N_1$ so that $V \big( g_n;[x_n',y_n'] \big) \leq \varepsilon / 2$ 
for every $n \geq N_2$. For such an $n$ and every $t \in [x_n',y_n']$ we have 
$$g_n'(t) \leq \frac{v_n - u_n}{y_n - x_n} \exp(\varepsilon / 2) < \frac{u_n -
a}{x_n - a},$$
$$g_n'(t) \geq \frac{v_n - u_n}{y_n - x_n} \exp(-\varepsilon / 2)
> \frac{f^{-1}(u_n) - u_n}{f^{-1}(u_n) - x_n}.$$
By the Mean Value Theorem, the first inequality shows that \esp $g_n(t)\!>\! t$ 
\esp for every $t \!\in\! [a,x_n]$, whereas the second one gives $g_n(t) \!> t$ 
for all $t \in [x_n,f^{-1}(u_n)]$ (see Figure 20).

\vspace{0.35cm}

Now fix $n \!\geq\! N_2$. From $g_n(a) \!<\! g_n(x_n) \!=\! u_n$ it follows 
that $f^{-1}g_n(a) \!<\! f^{-1}(u_n) \!<\! b_{\varepsilon}$. Let us consider 
the restrictions of $f$ and $g_n$ to the interval $[a,f^{-1}g_n(a)]$. We have  
\begin{eqnarray*}
V \big( f;[a,f^{-1} g_n(a)] \big) + V \big( g_n;[a,g_n^{-1}f^{-1}g_n(a)] \big) 
&\leq& 
V \big( f;[a,b_{\varepsilon}] \big) + V \big( g_n;[x_n',y_n'] \big)\\ 
&\leq& \varepsilon \esp \\
&<& \esp \log \Big( \frac{1}{1 - m(f)} \Big).
\end{eqnarray*}
However, this is in contradiction with Corollary \ref{retorno-dos} applied to the 
restrictions of $f^{-1}$ and $g_n^{-1}$ to the interval $[a,f^{-1}g_n(a)]$. $\hfill\square$

\vspace{1cm}


\beginpicture

\setcoordinatesystem units <1cm,1cm>

\circulararc 28 degrees from 0.5 0.4
center at -1.8 15

\circulararc -29 degrees from -0.3 -0.2
center at 16 -10

\putrule from -0.3 0.5 to 6.5 0.5
\putrule from 1 0.5 to 1 5.3

\plot 1 0.5 6.5 6 /

\putrule from 1.4 0.49 to 1.95 0.49
\putrule from 1.4 0.51 to 1.95 0.51

\putrule from 2.7 0.49 to 3.32 0.49
\putrule from 2.7 0.51 to 3.32 0.51


\setdots
\putrule from -0.3 0 to -0.3 0.5
\putrule from 1.4 2.2 to 5.55 2.2
\putrule from 1.4 0.5 to 1.4 2.2
\putrule from 5.55 0.5 to 5.55 5.05
\putrule from 2.7 0.5 to 2.7 2.2
\putrule from 6 0.5 to 6 6.2
\putrule from 1.95 0.5 to 1.95 2.82
\putrule from 1.95 2.82 to 3.32 2.82
\putrule from 3.32 0.5 to 3.32 2.82

\plot 1 0.5 1.4 2.2 /
\plot 1.4 2.2 5.55 5.05 /

\put{$a$} at 0.8 0.7
\put{$x_n'$} at -0.23 0.8
\put{$f^{-1}(u_n)$} at 5.45 0.2
\put{$x_n$} at 1.4 0.2
\put{$y_n$} at 1.95 0.2
\put{$u_n$} at 2.7 0.2
\put{$v_n$} at 3.32 0.2
\put{$y_n'$} at 6.35 0.8
\put{${\Large f}$} at 6.8 2.5
\put{${\Large g_n}$} at 2.9 4.4

\put{} at -3.15 -0.73
\put{Figure 20} at 3.6 -0.2

\endpicture


\vspace{0.25cm}

\subsection{End of the proof}

\hspace{0.45cm} To complete the proof of Theorem \ref{superduminy}, it suffices to find two 
intervals $]x,y[$ and $]u,v[$ satisfying the hypothesis of Proposition \ref{distinguir}. 
However, we immediately point out that this is a nontrivial issue, and it is very 
illustrative to follow the final steps (\textit{i.e.,} the proofs of claims (iii) 
and (iv)) in the proof below for the pseudo-group illustrated by Figure 16. 
We continue considering a positive constant 
$\varepsilon \leq \frac{1}{2} \log \big( \frac{1}{1 - m(f)} \big)$.

\vspace{0.25cm}

\noindent{\underbar{Claim (i).}} There exist $c_{\varepsilon} \!\!\in
]a,f(b_{\varepsilon})[$ and $g_{\varepsilon} \!\in \Gamma$ such that 
$[a,c_{\varepsilon}] \subset dom(g_{\varepsilon})$, \esp 
$g_{\varepsilon} (c_{\varepsilon}) = c_{\varepsilon}$, \esp 
$g_{\varepsilon} (t) > t$ \esp for all \esp $t \in [a,c_{\varepsilon}[$, 
\esp and \esp $V(g_{\varepsilon};[a,c_{\varepsilon}]) \leq \varepsilon$.

\vspace{0.15cm}

Indeed, since the orbit of $a \!\in\! \Lambda$ intersects the 
interval $]a,f(b_{\varepsilon})[$, we may choose $g \in \Gamma$ 
and $c \!\in ]a,f(b_{\varepsilon})[$ so that $[g(a),g(c)] \subset 
]a,f(b_{\varepsilon})[$ and $V(g;[a,c]) \leq \varepsilon / 2$. 
For $k$ large enough we have $f^k g (c) < c$, and since $f^k g(a) > a$, 
the map $g_{\varepsilon} = f^k g$ has fixed points in $]a,c[$. If we denote 
by $c_{\varepsilon}$ the first of these points to the right of $a$, then from 
(\ref{dist-todo-ene}) one deduces that \esp $V(g_{\varepsilon};[a,c_{\varepsilon})] 
\leq V \big( g;[a,c_{\varepsilon}] \big) + V \big( f^k;[a,b_{\varepsilon}] \big) 
\leq \varepsilon$.

\vspace{0.4cm}

\noindent{\underbar{Claim (ii).}} One has the inequality \esp
$f'(a) + g_{\varepsilon}'(c_{\varepsilon}) < e^{\varepsilon}$.

\vspace{0.15cm}

To show this claim, first notice that
$$V \big( f;[a,c_{\varepsilon}] \big) + V \big( g_{\varepsilon};[a,c_{\varepsilon}]
\big) \leq \frac{3 \varepsilon}{2} < \log \Big( \frac{1}{1 - m(f)} \Big).$$
Since the orbits of $\Lambda$ are not 
dense in $\Lambda \esp \cap \esp ]a,c_{\varepsilon}[$,
by aplying Corollary \ref{retorno-dos} to $f^{-1}$ and $g_{\varepsilon}^{-1}$
we conclude that $f(c_{\varepsilon}) \!<\! g_{\varepsilon} (a)$, which obviously 
implies that \esp $m(f) + m(g_{\varepsilon}) < 1$ \esp (where $m(g_{\varepsilon}) 
= \inf_{x \in [a,c_{\varepsilon}]} g_{\varepsilon}'(x)$). We thus  conclude that 
the value of \esp $f'(a) + g_{\varepsilon}'(c_{\varepsilon})$ \esp 
is bounded from above by  
$$m(f) \exp \big(V(f;[a,c_{\varepsilon}]) \big)
+ m(g_{\varepsilon}) \exp \big( V(g_{\varepsilon};[a,c_{\varepsilon}]) \big) \esp 
\leq \esp e^{\varepsilon} \big( m(f) + m(g_{\varepsilon}) \big) \esp 
< \esp e^{\varepsilon}.$$

In what follows, we will denote \esp $\alpha \!=\! f'(a)$ \esp and 
\esp $\beta \!=\! g_{\varepsilon}' (c_{\varepsilon})$ \esp (so that 
$\alpha + \beta < e^{\varepsilon}$), and we will impose to $\varepsilon 
> 0$ the extra condition (which holds for $\varepsilon>0$ small enough)
$$1 - \alpha e^{\varepsilon} \esp > \esp (e^{\varepsilon} - \alpha)^2.$$

\vspace{0.25cm}

\noindent{\underbar{Claim (iii).}} For each \esp $n \!\in\! \mathbb{N}$ 
\esp one may choose \esp $\kappa(\varepsilon,n) > 0$ \esp so that \esp 
$\lim_{\varepsilon \rightarrow 0} \esp \kappa(\varepsilon,n) > 0$ \esp and, 
for every $i \!\in\! \{1,\ldots,n\}$ and every $t \!\in\! [a,c_{\varepsilon}]$,
\begin{equation}
f^{n-i} g_{\varepsilon} f^i(t) - f^{n-i+1} g_{\varepsilon} f^{i-1}(t)
\esp \geq \esp (c_{\varepsilon} - a) \esp \kappa(\varepsilon,n).
\label{cuadro}
\end{equation}

\vspace{0.15cm}

Indeed, from \esp $\alpha e^{-\varepsilon} \leq f'(t) \leq \alpha e^{\varepsilon}$
\esp and \esp $f(a)=a$ \esp we conclude that, for every $t \in [a,c_{\varepsilon}]$,
$$a + \alpha e^{-\varepsilon}(t-a) \leq f(t) \leq a + \alpha e^{\varepsilon}(t-a).$$
Analogously, from \esp $\beta e^{-\varepsilon} \!\leq\!g_{\varepsilon}'(t) \!\leq\! \beta
e^{\varepsilon}$ \esp and \esp $g_{\varepsilon}(c_{\varepsilon}) \!=\! c_{\varepsilon}$ 
\esp we deduce that 
$$c_{\varepsilon} - \beta e^{\varepsilon} (c_{\varepsilon} - t) \leq
g_{\varepsilon}(t) \leq
c_{\varepsilon} - \beta e^{-\varepsilon} (c_{\varepsilon} - t).$$
Hence,
\begin{eqnarray*}
g_{\varepsilon}f(t) - fg_{\varepsilon}(t)
\!\! &\geq& \!\!
g_{\varepsilon} ( a+\alpha e^{-\varepsilon}(t-a) ) -
              ( a+\alpha e^{\varepsilon} (g_{\varepsilon}(t)-a) )\\
&\geq& \!\! c_{\varepsilon} \!-\! \beta e^{\varepsilon} 
( c_{\varepsilon} - a - \alpha e^{-\varepsilon} (t-a))
- ( a + \alpha e^{\varepsilon} 
(c_{\varepsilon} - \beta e^{-\varepsilon} (c_{\varepsilon}-t)-a) )\\
&=& \!\! (c_{\varepsilon} - a) [1 - \beta (e^{\varepsilon} - \alpha) - \alpha
e^{\varepsilon}]\\
&>& \!\! (c_{\varepsilon} - a) 
[1 - \alpha e^{\varepsilon} - (e^{\varepsilon} - \alpha)^2].
\end{eqnarray*}
Letting $\kappa(\varepsilon,n) =
(\alpha e^{-\varepsilon})^n \esp [1 - \alpha e^{\varepsilon} - (e^{\varepsilon} -
\alpha)^2] > 0$
we then deduce that, for every $i \!\in\! \{1,\ldots,n\}$ 
and every $t \in [a,c_{\varepsilon}]$,
\begin{eqnarray*}
f^{n-i} g_{\varepsilon} f^i (t) - f^{n-i+1} g_{\varepsilon} f^{i-1} (t)
\!\!&\geq&\!\! (\alpha e^{-\varepsilon})^{n-i}
                      \big[ g_{\varepsilon}f (f^{i-1}(t)) - f g_{\varepsilon}
(f^{i-1}(t)) \big]\\
&>& \!\! (\alpha e^{-\varepsilon})^{n-i} (c_{\varepsilon} - a) 
[1 - \alpha e^{\varepsilon} - (e^{\varepsilon} - \alpha)^2]\\
&>& \!\!(c_{\varepsilon} - a) \esp \kappa(\varepsilon,n).
\end{eqnarray*}
Finally, notice that
$$\lim_{\varepsilon \rightarrow 0} \kappa(\varepsilon,n) 
= \alpha^n [1 - \alpha - (1-\alpha)^2] > 0.$$

The preceding construction may be carried out for every $\varepsilon > 0$ very small. 
Let us then fix $n \!\in\! \mathbb{N}$, and let us impose to $\varepsilon>0$ the 
supplementary conditions
\begin{equation}
1 + \kappa(\varepsilon,n) \geq e^{2 (n+2) \varepsilon}, \qquad 
\qquad  1 - \kappa(\varepsilon,n) \leq e^{-2 (n+2) \varepsilon}.
\label{ast}
\end{equation}
Let $I$ be a connected component of $\Lambda^c$ contained in 
$]a,c_{\varepsilon}[$, and for each $i \!\in\! \{1,\ldots,n\}$ let 
$I_i = f^{n-i} g_{\varepsilon} f^i (I).$

\vspace{0.4cm}

\noindent{\underbar{Claim (iv).}} For all $i \!<\! j$, the intervals $]x,y[=\! I_i$ 
and $]u,v[ =\! I_j$ satisfy the conditions in Proposition \ref{distinguir}.

\vspace{0.15cm}

Conditions (i) and (ii) easily follow from the construction. To check 
condition (iii) let us first notice that, for every $k\!\in\!\{1,\ldots,n\}$,
$$(\alpha e^{-\varepsilon})^n \beta e^{-\varepsilon} \leq \frac{|I_k|}{|I|}
\leq (\alpha e^{\varepsilon})^n \beta e^{\varepsilon},$$
and therefore
$$e^{-2(n+1)\varepsilon} \leq \frac{|I_i|}{|I_j|} \leq e^{2(n+1)\varepsilon}.$$
Hence, by (\ref{cuadro}) and (\ref{ast}),
\begin{eqnarray*}
\frac{u-a}{x-a} \cdot \frac{y-x}{v-u} 
&=& \Big( \frac{u-x}{x-a} + 1 \Big) \esp \frac{y-x}{v-u} \\
&\geq& \Big( \frac{(c_{\varepsilon}-a) \esp \kappa(\varepsilon,n)}{x-a} + 1 \Big)
\esp \frac{y-x}{v-u} \\
&\geq& \big( 1 + \kappa(\varepsilon,n) \big) 2^{-2(n+1) \varepsilon} \\
&\geq& e^{2\varepsilon},
\end{eqnarray*}
that is,
$$\frac{u-a}{x-a} \exp(-2\varepsilon) \geq \frac{v-u}{y-x}.$$
Analogously,
\begin{eqnarray*}
\frac{f^{-1}(u)-u}{f^{-1}(u)-x} \!\cdot\! \frac{y-x}{v-u} 
&=& \Big( 1 \!-\! \frac{u-x}{f^{-1}(u)-x} \Big)
\frac{y-x}{v-u}\\
&\leq& \Big( 1 \!-\! \frac{(c_{\varepsilon} - a) \esp
\kappa(\varepsilon,n)}{f^{-1}(u)-x} \Big) \frac{y-x}{v-u} 
\leq \big( 1 - \kappa(\varepsilon,n) \big) e^{2(n+1)\varepsilon} 
\leq e^{-2\varepsilon},
\end{eqnarray*}
and hence
$$\frac{f^{-1}(u) - u}{f^{-1}(u) - x} \exp(2\varepsilon) \leq \frac{v-u}{y-x}.$$
This concludes the verification of condition (iii) of Proposition 
\ref{distinguir}, thus finishing the proof of Duminy's second theorem.

\vspace{0.11cm}

\begin{small} \begin{ejer} Show that, for the pseudo-group generated by the maps $f$ and $g$ 
in Figure 16, the Cayley graph associated to the semi-exceptional orbit is a tree. Using this 
fact, show directly that the corresponding space of ends is infinite. For the orbits in $\Lambda$ 
which are not semi-exceptional show that, despite the fact that they may not have the tree 
structure, they still have infinitely many ends. 
\index{Duminy!second theorem|)} 
\index{tree}
\label{finir}
\end{ejer} \end{small}


\section{Two Open Problems}

\hspace{0.45cm} There exist two major open questions concerning the dynamics of groups 
of $\ce^{1+\mathrm{Lip}}$ circle diffeomorphisms, namely the ergodicity of minimal 
actions, and the zero Lebesgue measure for exceptional minimal sets. In what 
follows, we will make an overview of some partial results, and we will 
explain some of the difficulties lying behind them.


\subsection{Minimal actions}
\label{sec-minimal-inv}

\hspace{0.45cm} Let us begin by making more precise the question of ergodicity 
for minimal actions. For this, recall that an action is {\bf{\em ergodic}} if the 
measurable invariant sets have null or total measure. In what follows, the ergodicity 
will be always considered with respect to the (normalized) Lebesgue measure.
\index{ergodic!action|see{action, ergodic}}
\index{action!ergodic|(}

\vspace{0.25cm}

\noindent{\underbar{\bf Problem 1.}} Let $\Gamma$ be a finitely generated 
subgroup of $\lip$ all of whose orbits are dense. Is the action of $\Gamma$ 
on $\clo$ ergodic~? 

\vspace{0.25cm}

An analogous question for germs at the origin of analytic diffeomorphisms of the complex 
plane has an affirmative answer \cite{gomez-mont,ilyashenko}. This is one of the reasons 
why one should expect that the answer to the problem above is positive. (See 
Exercise \ref{ejer-victor} for another evidence.) This is the 
case for instance when $\Gamma$ has an element with irrational rotation number, 
according to the result below obtained independently by Herman for the $\ce^2$ case,  
and by Katok for the general $\ce^{1+\mathrm{bv}}$ case (see \cite{He} and 
\cite{katok-erg}, respectively).  
\index{Herman} \index{Katok}Here we reproduce 
Katok's proof, which strongly uses the combinatorial structure 
of the dynamics of irrational rotations: see \S \ref{teoria-poinc}. 

\vspace{0.18cm}

\begin{thm} {\em If a diffeomorphism $g \in \vl$ has irrational rotation number, 
then the action of (the group generated by) $g$ on the circle is ergodic.}
\label{katokherman}
\end{thm}

\noindent{\bf Proof.} Recall that a {\bf\em{density point}} 
of a measurable subset $A \subset \clo$ 
\index{density point}
is a point $p$ such that  
$$\frac{ Leb \big( A \cap [p-\varepsilon,p+\varepsilon] \big) }{
Leb \big( [p-\varepsilon,p+\varepsilon] \big)}$$ 
converges to $1$ as $\varepsilon$ goes to zero. A classical theorem by Lebesgue 
asserts that, in any measurable set, almost every point is a density point.

Let $A \!\subset\! \clo$ be a $\Gamma$-invariant measurable subset. Assuming that the Lebesgue 
measure of $A$ is positive, we will show that this measure is actually total. For this, let 
us consider a density point $x_0$ of $A$. Let $\varphi$ be the conjugacy between $g$ 
and the rotation of angle $\rho(g)$ such that $\varphi(x_0) = 0$. As in the proof 
of Denjoy's theorem, let $I_n(g) \!=\! \varphi^{-1}(I_n)$ and 
$J_n(g) \!=\! \varphi^{-1}(J_n)$. Notice that $x_0$ belongs to 
$J_n(g)$, and $|J_n(g)|$ converges 
to zero as $n$ goes to infinity. Moreover, the intervals $g^k \big( J_n(g) \big)$,
$k \in \{0,\ldots,q_{n+1}-1 \}$, cover the circle, and each point in $\clo$ is contained 
in at most two of these intervals. Using inequality (\ref{tresduminy}) it is not 
difficult to conclude that, for every $k \in \{ 0, \ldots, q_{n+1}-1 \}$,
\begin{eqnarray*}
\frac{Leb\big( g^k(J_n(g) \setminus A)\big) }{Leb\big( g^k(J_n(g)) \big)}
&\leq& \exp \big( V(g^k;J_n(g)) \big) \cdot
\frac{Leb\big( J_n(g) \setminus A \big)}{Leb\big( J_n(g)\big)}\\
&\leq& \exp \big( 2 V(g) \big) \cdot
\frac{Leb\big( J_n(g) \setminus A \big)}{Leb\big( J_n(g)\big)}.
\end{eqnarray*}
Therefore, by the invariance of $A$,
\begin{eqnarray*}
Leb (\clo \setminus A) &\leq&
\sum_{k=0}^{q_{n+1}-1} Leb \big( g^k(J_n(g) \setminus A) \big)\\
&\leq& \exp \big( 2 V(g) \big) \hspace{0.02cm}
\frac{Leb\big( J_n(g) \setminus A \big)}{Leb\big( J_n(g)\big)}
\sum_{k=0}^{q_{n+1}-1} Leb \big( g^k(J_n(g)) \big)\\
&\leq& 2 \hspace{0.03cm} \exp \big( 2 V(g) \big) \hspace{0.02cm}
\frac{Leb\big( J_n(g) \setminus A \big)}{Leb\big( J_n(g)\big)}.
\end{eqnarray*}
Since $x_0$ is a density point of $A$, the value of \esp 
$Leb\big( J_n(g) \setminus A \big) / Leb\big( J_n(g)\big)$ 
\esp converges to $0$ as $n$ goes to infinity, from where one 
easily concludes that \esp $Leb(\clo \setminus A) = 0$. $\hfill\square$

\vspace{0.15cm}

\begin{small}\begin{obs} Although every $\ce^{1+\mathrm{bv}}$ circle diffeomorphism with irrational 
rotation number is topologically conjugate to the corresponding rotation, in ``many'' cases the 
conjugating map is {\em singular}, even in the real-analytic case (see the discussion at the 
begining of \S \ref{poquito}). Thus, the preceding theorem does not follow from Denjoy's theorem.
\label{singular-linearization}
\end{obs}\end{small}

\vspace{0.02cm}

According to Proposition \ref{gatomedible}, every group of circle homeomorphisms 
acting minimally and preserving a probability measure is topologically conjugate 
to a group of rotations. Moreover, if such a group is finitely generated, then  
at least one of its generators must have irrational rotation number. From Theorem 
\ref{katokherman} one then deduces the following.

\vspace{0.1cm}

\begin{cor} {\em If a finitely generated subgroup of 
$\mathrm{Diff}^{1+\mathrm{bv}}_+(\clo)$ acts minimally and 
preserves a probability measure, then its action on $\clo$ is ergodic.}
\label{general-erg}
\end{cor}

\vspace{0.15cm}

The differentiability hypothesis is necessary for this result. 
(The hypothesis of finite generation is also necessary: see 
\cite{DKN-nuevo} for an example illustrating this.) Indeed, in \cite{OR} 
the reader may find examples of $\ce^1$ circle diffeomorphisms with irrational 
rotation number acting minimally but not ergodically. It is very plausible that, 
by refining the construction method of \cite{OR}, one may provide analogous 
examples in class $\ce^{1+\tau}$ for every $\tau \!\in ]0,1[$.

\vspace{0.15cm}

Due to Theorem \ref{katokherman}, for groups of $\ce^{1+\mathrm{Lip}}$ circle 
diffeomorphisms acting minimally, the ergodicity question arises when the rotation 
number of each of its elements is rational. Under such a hypothesis, the general  
answer to this question is unknown. Nevertheless, there exist very important 
cases where the ergodicity is ensured. 

\vspace{0.15cm}

\begin{defn} Given a subgroup $\Gamma$ of $\mathrm{Diff}_{+}^1(\clo)$, a point 
$p \in \clo$ is {\bf \em{expandable}} if there exists $g \in \Gamma$ such that 
$g'(p) > 1$. The action is {\bf{\em differentiably expanding}} if for every 
point $p \in \clo$ is expandable.
\index{expandable point}
\index{action!differentiably expanding}
\index{expanding!action|see{action, differentiably expanding}}
\label{expergo}
\end{defn}

\vspace{0.15cm}

As we will see, if the action of a finitely generated subgroup of 
$\mathrm{Diff}_{+}^{1+\tau}(\clo)$ is differentiably expanding and minimal,  
with $\tau > 0$, then it is necessarily ergodic. Notice, however, that the hypothesis 
of this claim is not invariant under smooth conjugacy, although the conclusion is. 
To state an {\em a priori} more general result which takes into account all of 
this, let us fix some notation. Given a finitely generated subgroup $\Gamma$ of 
$\mathrm{Diff}_+^1(\clo)$, let $\mathcal{G}$ be a finite symmetric system 
of generators. For each $n \geq 1$ let (compare Appendix B)
$$B_{\mathcal{G}}(n) \esp = \esp \{ h \in \Gamma: \esp h = h_m \cdots h_1
\quad \mbox{for some } h_i \in \mathcal{G} \mbox{ and } m \leq n\},$$
and for each $x \in \clo$ let
$$\lambda(x) \esp = \esp \limsup_{n \rightarrow \infty}
\left( \max_{h \in B_{\mathcal{G}}(n)} \frac{\log \big( h'(x) \big)}{n} \right).$$
Notice that this number is always finite, since it is bounded from above by 
$$\sup_{h \in \mathcal{G}, y \in \clo} \log \big( h'(y) \big).$$
For each $\lambda \!>\! 0$ \esp let \esp 
$E_{\lambda}(\Gamma) = \{ x \!\in\! \clo\!: \lambda(x) \!\geq\! \lambda \}.$ 
\esp The {\bf{\em exponential set}} 
\index{exponential set of an action} 
$E(\Gamma)$ of the action is defined as the union of the sets $E_{\lambda}(\Gamma)$ 
with $\lambda \!>\! 0$; its complement $S(\Gamma)$ is called the 
{\bf{\em sub-exponential set}}.
\index{sub-exponential set of an action} 
Notice that every $E_{\lambda}(\Gamma)$, as well as $E(\Gamma)$ and $S(\Gamma)$, are
Borel sets. Moreover, the function $x \mapsto \lambda(x)$ is invariant by the 
$\Gamma$-action. Therefore, the sets $E_{\lambda}(\Gamma)$, $E(\Gamma)$, and 
$S(\Gamma)$, are invariant by $\Gamma$. Finally, if (the action of) $\Gamma$ 
is differentiably expanding, then $\lambda(x) > 0$ for every $x \in \clo$.

\vspace{0.22cm}

\begin{thm} {\em Let $\Gamma$ be a finitely generated subgroup of 
$\mathrm{Diff}_+^{1+\tau}(\clo)$ whose action is minimal, where $\tau > 0$. 
If the exponential set of $\Gamma$ has positive Lebesgue measure, 
then it has full measure, and the action is ergodic.}
\label{principal}
\end{thm}

\vspace{0.01cm}

\begin{small} \begin{obs} Both hypotheses above (the facts that $\tau \!>\! 0$ and 
that the group is finitely ge- nerated) are necessary for the validity of the 
theorem. For a detailed discussion of this, see \cite{DKN-nuevo}.
\end{obs} 

\begin{ejem} Theorem \ref{principal} allows showing that if a subgroup of $\lip$ 
satisfies the hypothesis of Duminy's first theorem and has no finite orbit, then its 
action is ergodic. Indeed, let $\Gamma$ be a subgroup of $\lip$ acting minimally and  
generated by a family $\mathcal{G}$ of diffeomorphisms satisfying $V(f) \!<\! V_0$ for 
every $f \!\in\! \mathcal{G}$. If we assume that the set $\mathrm{Per}(g)$ of periodic 
points is finite for at least one element $g \!\in\! \mathcal{G}$, then there must exist 
$p \!\in\! \mathrm{Per}(g)$ and $f \!\in\! \mathcal{G}$ so that $f(p)$ belongs to $\clo 
\setminus \mathrm{Per}(g).$  Using some of the arguments in the proof of Theorem 
\ref{thduminy}, 
\index{Duminy!first theorem} 
starting with $f$ and $g$ it is possible to create elements $F, G$ in $\Gamma$ 
which are as in Figure 17 over an interval $[a,b]$ of $\clo$ in such a way that 
$V(F;[a,c])$ and $V(G;[c,b])$ are small, where $c = G^{-1}(a)$. Fix $\varepsilon > 0$
very small, and let $\{ I_1, \ldots, I_n \}$ be a family of intervals covering 
the circle so that there exist $h_1, \ldots, h_n$ in $\Gamma$ satisfying  
$h_i(x) \!\in\! [c+\varepsilon,b]$ for every $x \in I_i$. Let $C$ be 
the constant defined by $C^{-1} \!=\! \inf \{ h_i'(x) \!: x \in I_1 
\cup \ldots \cup I_n \}$, and let $N$ be a sufficiently large integer so that 
each branch of the return map $H^N$ induced by $F$ and $G$ is $C$-expanding 
(see Proposition \ref{retorno}). For $g_i = H^N h_i \in \Gamma$ one has 
$g_i'(x) > 1$ for every $x \in I_i$ (where we consider the right derivative 
in case of discontinuity). Thus, the action of $\Gamma$ is differentiably 
expanding, and hence its ergodicity follows from Theorem \ref{principal}.
\end{ejem}\end{small}

Unfortunately (or perhaps fortunately), there exist minimal actions for which 
the exponential set has zero Lebesgue measure. For instance, this is the case  
for the (standard action of the) modular group, as well as for the smooth, 
minimal actions of Thompson's group $\mathrm{G}$: see \cite{DKN-nuevo}.

\begin{small}
\begin{ejer} Give a precise statement and show a result in the following 
spirit: If $\Gamma$ admits a continuous family of representations $\Phi_t$ 
in $\mathrm{Diff}_+^{1+\mathrm{bv}}(\clo)$ so that all the orbits by 
$\Phi_0 (\Gamma)$ are dense and $\Phi_t(\Gamma)$ admits an exceptional 
minimal set for each $t \!> \!0$, then the action of $\Phi_0 (\Gamma)$ 
is not differentiably expanding. 
\label{claro-no}
\end{ejer}

\begin{obs} Recall that if $\Gamma$ is a non-Abelian countable group of circle 
homeomorphisms acting minimally, then its action is ergodic with respect to 
every stationary measure (this follows directly from Theorems \ref{unicidad} 
and \ref{kak-th}). It is then natural to ask whether in the case of groups of 
{\em diffeomorphisms}, there always exists a probability distribution on the 
group so that the corresponding stationary measure is absolutely continuous with 
respect to the Lebesgue measure (compare Remark \ref{singular-linearization}). 
However, this relevant problem seems to be very hard. A partial result from 
\cite{DKN-nuevo} points in the negative direction: for the cases of the 
modular group and Thompson's group G, the stationary measure associated to 
any finitely supported probability distribution on the group is singular. 
Actually, for these cases the exponential set of the action has zero Lebesgue 
measure, but its mass with respect to the stationary measure is total.
\end{obs}\end{small}

\index{measure!stationary}

Despite the preceding discussion, for the actions of the modular group and Thompson's 
group $\mathrm{G}$ already mentioned --as well as for most ``interesting'' actions in 
the literature-- the set of points which are non-expandable is finite and is made 
up of isolated fixed points of certain elements. Under such a hypothesis we can 
give the following general result from \cite{DKN-nuevo}, which covers Theorem 
\ref{principal} at least in the $\ce^{1+\mathrm{Lip}}$ case (for a complete 
proof of Theorem \ref{principal}, we refer the reader to \cite{yo-brasil}). 

\vspace{0.15cm}

\begin{thm} {\em Let $\Gamma$ be a finitely generated subgroup of $\ce^{1+\mathrm{Lip}}$ 
circle diffeomorphisms whose action is minimal. Assume that for every non-expandable point 
$x\in \clo$ there exist $g_+,g_-$ in $\Gamma$ such that $g_+(x)=g_-(x)=x$, and such that 
$g_+$ (resp. $g_-$) has no fixed point in some interval $]x,x+\varepsilon[$ (respectively, 
$]x-\epsilon,x[$). Then the action of \esp $\Gamma$ is ergodic.}
\label{principalazo}
\end{thm}

\vspace{0.15cm}

To show this result we will use some slight modifications of Schwartz 
\index{Schwartz estimates} 
estimates from \S \ref{chuartz} that we leave as exercises to the reader. 

\vspace{0.02cm}

\begin{small} \begin{ejer} Given two intervals $I,J$ and a $\ce^1$ map $f\!: I \rightarrow J$ 
which is a diffeomorphism onto its image, define the {\bf {\em distortion coefficient}} of 
$f$ on $I$ as 
$$ \varkappa (f;I) \esp = \esp 
\log \Big( \frac{\max_{x \in I} f'(x)}{\min_{y \in I} f'(y)} \Big), $$
and its {\bf \em{distortion norm}} as
$$ \eta (f;I) \esp = 
\sup_{\{x,y\} \subset I} \frac{\log \big( \frac{f'(x)}{f'(y)}\big)}{|f(x)-f(y)|} 
\esp = \esp \max_{z \in J} \left| \big( \log ( (f^{-1})')\big)'(z) \right|.$$
\index{control of distortion}

\vspace{-0.03cm}

\noindent (i) Show that the distortion coefficient is subadditive under composition. 

\vspace{0.08cm}

\noindent (ii) Show that $\varkappa(f,I) \leq C_f |I|,$ where the constant $C_f$ equals 
the maximum of the absolute value of the derivative of the function $\log(f')$. Conclude  
that if $\mathcal{G}$ is a subset of $\mathrm{Diff}^{1+\mathrm{Lip}}_+(S^1)$ such that the 
set $\{|(\log(f'(x)))'|\!: f \in \mathcal{G}, \esp x \in \clo \}$ is bounded, then there 
exists a constant $C_{\mathcal{G}}$ (depending only on $\mathcal{G}$) such that, for 
every interval $I$ in the circle and every $f_1,\dots,f_n$ in $\mathcal{G}$, 
$$ \varkappa(f_n \circ \cdots \circ f_1;I) \esp \leq \esp  
C_{\mathcal{G}} \sum_{i=0}^{n-1} |f_i \circ \cdots \circ f_1 (I)|.$$

\vspace{0.08cm}

\noindent (iii) Under the assumptions in (ii), let $I$ be an interval of $\clo$\!, 
and let $x_0$ be a point in $I$. Denoting $F_i = f_i \circ \cdots \circ f_1$, 
\esp $I_i = F_i(I)$, \esp and \esp $x_i = F_i(x_0)$, \esp show that
\begin{equation}\label{eq:dist}
\exp \Big( -C_{\mathcal{G}}\sum_{j=0}^{i-1} |I_j| \Big) \cdot \frac{|I_i|}{|I|} 
\esp \leq \esp 
F_i'(x_0) 
\esp \leq \esp 
\exp \Big( C_{\mathcal{G}}\sum_{j=0}^{i-1} |I_j| \Big) \cdot \frac{|I_i|}{|I|},
\end{equation}
\begin{equation}\label{eq:lsum}
\sum_{i=0}^n |I_i| (
\esp \leq  \esp  
|I| \exp \Big( C_{\mathcal{G}} \sum_{i=0}^{n-1} |I_i| \Big) \sum_{i=0}^{n} F_i'(x_0).
\end{equation}

\vspace{0.08cm}

\noindent (iv) Still under the conditions in (ii), show that if for  
$x_0\!\in\!S^1$ we denote $S = \sum_{i=0}^{n-1} F_i'(x_0)$, then 
for every $\delta \leq \log(2) / 2 C_{\mathcal{G}} S$ one has 
\begin{equation}\label{baba}
\varkappa \big( F_n;]x_0 - \delta/2,x_0 + \delta/2[ \big) 
\esp \leq \esp 2 C_{\mathcal{G}} S \delta.
\end{equation}
\label{mucho}
\end{ejer} 

\vspace{-0.75cm}

\begin{ejer} The action of a group of circle homeomorphisms $\Gamma$ is said to be 
{\bf {\em conservative}} if for every measurable subset $A \subset \clo$ of positive 
Lebesgue measure there exists $g \neq id$ in $\Gamma$ such that $Leb ( A \cap g(A) ) > 0$. 
Show that every ergodic action is conservative. Conversely, following the steps below, 
show that if the action of a finitely generated group of $\ce^{1+\mathrm{Lip}}$ circle 
diffeomorphisms is minimal, then it is conservative.

\vspace{0.08cm}

\noindent (i) Assuming that the action is not conservative, show that the 
set $\sum (\Gamma)$ formed by the points $y \in \clo$ for which the series  
$\sum_{g \in \Gamma} g'(y)$ converges is non-empty. 

\vspace{0.08cm} 

\noindent (ii) Using (\ref{baba}), from (i) deduce that $\sum (\Gamma)$ is 
open, and using the minimality of the action conclude that this set actually 
coincides with the whole circle.
  
\vspace{0.08cm}

\noindent (iii) Obtain a contradiction by choosing an infinite order element 
$g \!\in\! \Gamma$ and using the equality 
$$\sum_{i=0}^{n-1} \int_{\clo} (g^i)'(y) \esp dy \esp = \esp 2 \pi \esp \! n,$$
\label{ejer-victor}
\end{ejer}\end{small}

\vspace{-0.31cm}

From now on, we will assume that $\Gamma$ is a subgroup of $\lip$ satisfying the hypothesis 
of Theorem \ref{principalazo}. In this case, the set of non-expandable points is closed, 
since it coincides with \esp $\bigcap_{g\in \Gamma} \{x\!: g'(x)\le 1\}$. \esp Together 
with the lemma below, this implies that this set is actually finite. 

\vspace{0.15cm}

\begin{lem}\label{p:isolated}
{\em The set of non-expandable points is made up of isolated points.}
\end{lem}

\noindent{\bf Proof.} For each non-expandable point $y \!\in\! \clo$,  we will find an interval 
of the form $]y,y+\delta[$ which does not cointain any non-expandable point. The reader 
will notice that a similar argument provides us with an interval of the form 
$]y - \delta',y[$ sharing this property.

By hypothesis, there exist $g_+ \!\in G$ and $\varepsilon > 0$ such
that $g_+(y)=y$ and such that~$g_+$ has no other fixed point in
$]y,y+\varepsilon[$. Changing $g_+$ by its inverse if necessary, we
may assume~$y$ to be a right topologically repelling fixed point of
$g_+$. Let us consider the point $\bar{y} = y + \varepsilon/2 \in 
]y,y+\varepsilon[$, and for each integer $k \geq 0$ let $\bar{y}_k =
g_+^{-k}(\bar{y})$ and $J_k = ]\bar{y}_{k+1},\bar{y}_{k}[$. Taking
$a=\bar{y}_1, b=\bar{y}$, and applying (\ref{eq:dist}), we conclude 
the existence of $k_0 \in \mathbb{N}$ such that $(g_+^k)'(x) \geq 2$ 
for all $k\geq k_0$ and all $x \!\in\! J_k$. Hence, each $J_k$ with $k \geq k_0$ 
contains no non- expandable point, which clearly implies that the same holds 
for the interval $]y,\bar{y}_{k_0}[$. $\hfill\square$

\vspace{0.5cm}

According to the proof above, for each non-expandable point $y \in \clo$ one can find  
an interval $I_y^+ \!= ]y,y+\delta^{+}[$, a positive integer $k_0^+$, and an element 
$g_+ \in \Gamma$ having~$y$ as a right topologically repelling fixed point and with 
no other fixed point than $y$ in the closure of $I_{y}^+$ such that, if for 
$x \in I_{y}^+$ we take the smallest integer $n \geq 0$ satisfying  
$g_+^n (x) \notin I_y^+$, \esp then $(g_+^{n+k_0^+})'(x) \geq 2$. 
\esp Obviously, one can also find an interval $I_{y}^-\!=]y-\delta^{-},y[$), 
a positive integer $k_0^-$, and an element $g_- \in \Gamma$, satisfying 
analogous properties. We then let
$$U_{y} = I_{y}^+ \cup I_{y}^- \cup \{y\}.$$
By definition (and continuity), for every expandable point $y$ there exist
$g = g_y \in \Gamma$ and a neighborhood $V_y$ of $y$ such that \esp 
$\inf_{x \in V_y} g'(x) > 1$. The sets 
$\{U_y \!: y \mbox{ non-expandable} \}$ 
and $\{V_y \!: y \mbox{ expandable} \}$ form an open cover 
of the circle, from which we can extract a finite sub-cover
$$\{U_y \!: y \mbox{ non-expandable} \} \esp \bigcup \esp \{V_{y_1},\dots, V_{y_k} \}.$$
Let
$$\lambda \esp = \esp 
\min \Big\{2, \inf_{x \in V_{y_1}} g_{y_1}'(x), \dots , 
\inf_{x \in V_{y_k}} g_{y_k}'(x) \Big\}.$$
Since $\lambda$ is the minimum among finitely 
many numbers greater than~$1$, we have $\lambda>1$.

\vspace{0.3cm}

\begin{lem} {\em For every point $x \in S^1$ either the set of derivatives 
$\{g'(x) \!: g \in \Gamma \}$ is unbounded, or $x$ belongs to the orbit of 
some non-expandable point.}
\end{lem}

\noindent{\bf Proof.} If $x \!\in\! S^1$ is expandable, then 
it belongs to some of the sets $I_y^{\pm}$ or $V_{y_j}$, and 
there exists a map $g \!\in\! \Gamma$ such that $g'(x) \ge \lambda$. 
Similarly, the image point $g(x)$ is either non-expandable, or there exists 
$h \!\in\! \Gamma$ such that $h' \big( g(x) \big) \geq \lambda$. By repeating 
this procedure we see that if we do not fall into a non-expandable point by 
some composition, then we can always continue expanding by a factor at least 
equal to $\lambda$ by some element of $\Gamma$. Therefore, for each point 
not belonging to the orbit of any non-expandable point, the set of derivatives 
$\{g'(x) \!: g \in \Gamma \}$ is unbounded. Since for each $x$ in the orbit of 
a non-expandable point this set is obviously bounded, this proves the lemma. 
$\hfill\square$

\vspace{0.5cm}

For the proof of Theorem \ref{principalazo} we will use 
the so-called ``expansion argument'' (essentially due to Sullivan), 
\index{Sullivan} 
which is one of the most important techniques for showing the ergodicity of dynamical 
systems having some hyperbolic behavior. Let $A \subset S^1$ be an invariant
measurable set of positive Lebesgue measure, and let~$a$ be a density point of~$A$ not 
belonging to the orbit of any non-expandable point (notice that, since there are only 
finitely many non-expandable points and $\Gamma$ is countable, such a point $a$ 
exists). The idea now consists in ``blowing-up'' a very small neighborhood of 
$a$ by using a well chosen sequence of compositions so that the distortion is 
controlled and the length of the final interval is ``macroscopic'' ({\em i.e.,} 
larger than some prescribed positive number). Since $a$ is a density point of $A$, this 
final interval will mostly consist of points in $A$, and by the minimality of the action, 
this will imply that the measure of the set $A$ is very close to 1. Finally, by performing 
this procedure starting with smaller and smaller neighborhoods of $a$, this will yield 
the total measure for the invariant set $A$.
 
In our context, the expansion procedure will work by applying the ``exit-maps'' 
$g_{\pm}^{n+k_0^{\pm}}$ to points in $I_y^+ \cup I_y^-$, and the maps $g_{y}$ to points 
in the neighborhoods $V_{y}$. To simplify the control of distortion estimates, in what 
follows we consider a symmetric generating system $\mathcal{G}$ of $\Gamma$ containing 
the elements of the form $g_+$ and $g_{-}$. 

\vspace{0.1cm}

\begin{lem}
{\em There exists a constant $C_1>0$ such that, for every expandable 
point $x \!\in\! \clo$, one can find $f_1,\dots,f_n$ in $\mathcal{G}$ 
such that \esp $(f_n\circ\dots\circ f_1)'(x)\ge \lambda$ \esp and}
\begin{equation}
\label{eq:sum-ld}
\frac{\sum_{j=1}^n (f_j\circ\dots\circ f_1)'(x)}{(f_n\circ\dots\circ f_1)'(x)} \le C_1.
\end{equation}
\end{lem}

\noindent{\bf Proof.} 
A compactness type argument reduces the general case to studying (arbitrarily small) 
neighborhoods of non-expandable points. To find a small interval of the form 
$]y,y+\delta[$ formed by points satisfying the desired property, 
take the corresponding element 
$g_+$ having $y$ as a right-repelling fixed point and no fixed point in $]y,\bar{y}]$, 
and for each $k \geq 0$ let $y_k = g_+^{-k}(\bar{y})$ and $J_k = g^{-k} (J_0)$, where 
$J_0 = [\bar{y}_1,\bar{y}[$. We know that for some $k_0^+$ we have $(g_+^n)'(x) \geq \lambda$ 
for all $n \geq k_0^+$ and all $x \!\in\! J_n$. For each $x$ in the interval 
$I_y^+ \!= ]y,\bar{y}_{k_0^+}[$ take $n \in \mathbb{N}$ so that $x \!\in\! J_n$. 
Then the estimates in Exercise \ref{mucho} show that
$$\sum_{j = 1}^n (g_+^j)'(x) \esp \leq \esp \exp(C_{\mathcal{G}}) \cdot 
\frac{\sum_{j=1}^{n} |g_+^j (J_n)|}{|J_n|} \esp \leq \esp 
\frac{\exp(C_{\mathcal{G}})}{|J_n|},$$ 
$$(g_+^n)'(x) \esp \geq \esp \exp(-C_{\mathcal{G}}) \cdot 
\frac{|g^n(J_n)|}{|J_n|} \esp = \esp \exp(-C_{\mathcal{G}}) \cdot \frac{|J_0|}{|J_n|}.$$
Therefore, letting $f_1,\ldots,f_n$ be all equal to $g_+$,
$$\frac{\sum_{j=1}^n (f_j\circ\dots\circ f_1)'(x)}{(f_n\circ\dots\circ f_1)'(x)}
\esp = \esp \frac{\sum_{j = 1}^n (g_+^j)'(x)}{(g_+^n)'(x)}
\leq \frac{\exp \esp \!(2 C_{\mathcal{G}})}{|J_0|}.$$
Moreover, for $n \geq k_0^+$ we have \esp 
$(f_n\circ\dots\circ f_1)'(x) = (g_+^n)'(x) \geq 2 \geq \lambda$. \esp 
Analogous arguments may be applied to the interval \esp $[y-\delta,y[$ 
\esp associated to $y$. Since there are only finitely many 
non-expandable points, this concludes the proof. $\hfill\square$

\vspace{0.3cm}

\begin{lem}\label{l:derivatives-quotient} 
{\em There exists a constant $C_2$ such that, 
for every point $x$ which does not belong to the
orbit of any non-expandable point and every $M > 1$, 
one can find $f_1,\ldots,f_n$ in $\mathcal{G}$ such 
that $(f_n \circ \dots \circ f_1)'(x) \geq M$ and}
\begin{equation}\label{eq:sums-ld-2}
\frac{\sum_{j=1}^n (f_j\circ\dots\circ f_1)'(x)}{(f_n\circ\dots\circ f_1)'(x)} 
\esp \leq \esp C_2.
\end{equation}
\end{lem}

\noindent{\bf Proof.} Starting with $x_0 = x$ we let
$$
\widetilde{F}_k \esp = \esp f_{k,n_k} \circ\dots\circ f_{k,1}, 
\qquad x_k = \widetilde{F}_k (x_{k-1}),$$
where the elements $f_{k,j} \in \mathcal{G}$ satisfy
$$\frac{\sum_{j=1}^{n_k} (f_{k,j} \circ\dots\circ 
f_{k,1})'(x_{k-1})}{(f_{k,n_k}\circ\dots\circ f_{k,1})'(x_{k-1})}
\le C_1, \qquad (f_{k,n_k} \circ\dots\circ f_{k,1})'(x_{k-1}) \ge \lambda.$$
If we perform this procedure \esp $K \geq \log(M) / \log(\lambda)$ \esp 
times, \esp then for the compositions \esp 
$F_k = \widetilde{F}_k \circ \dots \circ \widetilde{F}_1$ \esp we obtain
$$F_K'(x) \esp = \esp \prod_{k=1}^K (f_{k,n_k}\circ\dots\circ f_{k,1})'(x_{k-1}) 
\ge \lambda^K \geq M.$$

To estimate the left-side expression of~\eqref{eq:sums-ld-2}  
notice that, for $y = f_n\circ\dots\circ f_1 (x)$, 
\begin{equation}\label{eq:sum-inverses}
\frac{\sum_{j=1}^n (f_j\circ\dots\circ f_1)'(x)}{(f_n\circ\dots\circ
f_1)'(x)} \esp = \esp \sum_{j=1}^n (f_{j+1}^{-1}\circ \dots f_n^{-1})'(y).
\end{equation}
If we denote $y \!=\! F_K (x)$, then using~\eqref{eq:sum-inverses} we see 
that the left-side expression of~\eqref{eq:sums-ld-2} is equal to
\begin{multline*}
\sum_{k=1}^K \sum_{j=1}^{n_k} \left((f_{k,j+1}^{-1} \circ \dots
\circ f_{k,n_k}^{-1})\circ (\widetilde{F}_{k+1}^{-1}\circ\dots\circ
\widetilde{F}_K^{-1})\right)' (y) \esp = 
\\
= \esp \sum_{k=1}^K (\widetilde{F}_{k+1}^{-1}\circ\dots\circ
\widetilde{F}_K^{-1})' (y) \cdot \sum_{j=1}^{n_k} (f_{k,j+1}^{-1}
\circ \dots \circ f_{k,n_k}^{-1})'(x_{k}) \esp \leq
\\
\leq \esp \sum_{k=1}^{K} \frac{1}{\lambda^{K-k}} \cdot C_1 
\esp \le \esp \frac{C_1}{1-\lambda^{-1}}. \quad
\end{multline*}

\vspace{-0.9cm}

$\hfill\square$

\vspace{0.35cm}

\begin{lem}\label{p:final}
{\em For certain $\varepsilon > 0$ the following property holds: for 
every point $x$ not belonging to the orbit of any non-expandable point, 
there exists a sequence  $V_k$ of neighborhoods of $x$ converging to $x$ 
such that, to each $k \in \mathbb{N}$, one may associate a sequence of 
elements $f_1,\ldots,f_{n_k}$ in $\mathcal{G}$ satisfying 
$|f_{n_k} \circ \dots \circ f_1 (V_k)| = \varepsilon$ 
and $\varkappa (f_{n_k} \circ \dots \circ f_1 ; V_k) \leq \log (2)$.}
\end{lem}

\noindent{\bf Proof.} We will check the conclusion for \hspace{0.01cm}
$\varepsilon = \log(2) / 2 C_{\mathcal{G}} C_2$. \hspace{0.01cm} To do 
this, fix $M > 1$, and consider the sequence of compositions $f_n \circ
\dots \circ f_1$ associated to $x$ and $M$ provided by the preceding 
lemma. Denoting $\bar{F}_n = f_n \circ \dots \circ f_1$ and
$y=\bar{F}_n(x)$, for the neighborhood 
$V = \bar{F}_n^{-1}(]y-\varepsilon/2,y+\varepsilon/2[))$ of $x$ we have 
$$\varkappa ( \bar{F}_n ; V ) 
\esp = \esp \varkappa \big( \bar{F}^{-1}_n ; ]y-\varepsilon/2,y+\varepsilon/2[ \big)
\esp = \esp \varkappa \big( f_1^{-1} \circ \dots \circ f_n^{-1};
]y-\varepsilon/2,y+\varepsilon/2[ \big).$$
The estimates in Exercise \ref{mucho} 
then show that the distortion coefficient of the 
composition $f_1^{-1} \circ \dots \circ f_n^{-1}$ is bounded 
from above by $\,\log (2)\,$ in a neighborhood of~$y$ of radius
\esp $r = \log(2) /4 C_{\mathcal{G}} S,$ \esp where
$$S = \sum_{j=1}^{n} (f_{j+1}^{-1} \circ \dots \circ f_n^{-1})'(\bar{F}_n (x_0)).$$
Now according to~\eqref{eq:sums-ld-2} and~\eqref{eq:sum-inverses}, we have 
$$
S = \frac{\sum_{j=1}^{n} (f_j\circ\dots\circ f_1)'(x)}{(f_n \circ\dots\circ
f_1)'(x)}  \le C_2.
$$
Thus, $\varepsilon / 2 \leq r$, which yields the desired estimate for the
distortion. Finally, since 
$$| V | = \big| \bar{F}_n^{-1} (]y-\varepsilon/2,y+\varepsilon/2[) \big| 
\leq \frac{\varepsilon \hspace{0.01cm} \exp \big(
\varkappa \big( \bar{F}_n^{-1}; ]y-\varepsilon/2,y+\varepsilon/2[ \big) \big)}{
(f_n \circ \dots \circ f_1)'(x)} \leq \frac{2 \varepsilon}{M},$$ and
since the last expression tends to zero as $M$ goes to infinity,
this concludes the proof of the lemma.
$\hfill\square$

\vspace{0.35cm}

To complete the proof of Theorem \ref{principalazo}, recall that for our  
\index{density point}
invariant subset $A \subset S^1$ of positive Lebesgue measure we found 
a density point $a$ not belonging to the orbit of any non-expandable 
point. Fix $\delta > 0$. By the preceding lemma, for some $\varepsilon > 0$ 
there exists a sequence $V_k$ of neighborhoods of $a$ converging to $a$ such 
that, for each $k \in \mathbb{N}$, there exist elements $f_1,\ldots,f_{n_k}$ 
in $\mathcal{G}$ satisfying $|f_{n_k} \circ \dots \circ f_1 (V_k)| = \varepsilon$ 
and $\varkappa (f_{n_k} \circ \dots \circ f_1 ; V_k) \leq \log (2)$. Passing to a 
subsequence if necessary, we may assume that the sequence of intervals  
$f_{n_k} \circ \cdots \circ f_{1} (V_k)$ converges to some interval 
$V$ of length $\varepsilon$. Due to the minimality of the 
action, the proof will be finished when showing that \esp 
$Leb \big( (\clo \setminus A) \cap V \big)=0.$ \esp 
To do this first observe that, by the invariance of $A$, 
\begin{eqnarray*}
\frac{Leb \big( (\clo \setminus A) \cap f_{n_k} \circ \cdots \circ f_{1} (V_k) \big)}
{\varepsilon}  
&=& \frac{Leb \big( (\clo \setminus A) \cap f_{n_k} \circ \cdots \circ f_{1} (V_k) \big)}
{Leb(f_{n_k} \circ \cdots \circ f_{1} (V_k))}\\ 
&\leq& \exp \big( \varkappa (f_{n_k} \circ \dots \circ f_1 ; V_k) \big) 
\cdot \frac{Leb \big( (\clo \setminus A) \cap V_k \big)}{Leb(V_k)} \\
&\leq& \frac{2 \esp Leb \big( (\clo \setminus A) \cap V_k \big)}{Leb(V_k)}. 
\end{eqnarray*} 
Now notice that the first expression in this inequality converges to 
\esp $Leb \big( (\clo \setminus A) \cap V \big) / \varepsilon$, \esp while 
the last expression converges to 0, because $a$ is a density point of $A$. 

\index{action!ergodic|)}


\subsection{Actions with an excepcional minimal set}
\label{sec-excep}

\hspace{0.45cm} Another major open question in the theory 
concerns the Lebesgue measure of exceptional minimal sets.

\vspace{0.35cm}

\noindent{\underbar{\bf Problem 2.}} If $\Gamma$ is a finitely generated subgroup 
of $\mathrm{Diff}_+^{1+\mathrm{lip}}(\clo)$ having an exceptional minimal set 
$\Lambda$, is the Lebesgue measure of $\Lambda$ necessarily equal to zero~?\\

\vspace{0.35cm}

A closely related problem concerns the finiteness for the number of orbits of 
connected components of $\clo \setminus \Lambda$ (if this is always true, this 
could be considered as a kind of generalization of the classical Ahlfors 
Finiteness Theorem \cite{ah-lect}). 
\index{Ahlfors Finiteness Theorem}
In the same spirit, a conjecture by Dippolito suggests that the action of $\Gamma$ on 
$\Lambda$ should be topologically conjugate to the action of a group of piecewise 
affine homeomorphisms \cite{Di}.
\index{Dippolito}

The Lebesgue measure of exceptional minimal sets is zero in the case of Fuchsian groups, 
but the techniques involved in the proof cannot be adapted to the case of general 
diffeomorphisms (see for instance \cite{nicholls}). In what follows, we will concentrate 
on a particular class of dynamics, namely that of \textit{Markovian} minimal sets. 
\index{Markov!minimal set}
Roughly speaking, over these sets the dynamics is conjugate to a {\em subshift}, 
and this information 
strongly simplifies the study of the combinatorics. 
Following \cite{CC3,matsu}, we will see that, for 
this case, the answer to Problem 2 is positive (the same 
holds for the finiteness of the connected components of the 
complement of exceptional minimal sets modulo the action \cite{CC4}).

\vspace{0.05cm}

Let $P \!=\! \big( p_{i j}\big)$ be a $k\!\times\!k$ {\bf \em{incidence matrix}} 
(\textit{i.e.,} a matrix with entries 0 and 1). Let us consider the space 
$\Omega = \{1,\ldots,k\}^{\mathbb{N}}$, and the subspace $\Omega^*$ formed by 
the {\bf\textit{admissible sequences}}, that is, by the 
$\omega \!=\! (i_1,i_2,\ldots) \!\in\! \Omega$ such that $p_{i_j i_{j+1}} = 1$
for every $j \in \mathbb{N}$. We endow this subspace with the topology induced 
from the product topology on $\Omega$, and we consider the dynamics of the 
shift $\sigma\!: \Omega^* \rightarrow \Omega^*$ on it.
\index{shift}

\vspace{0.1cm}

\begin{defn} Let $\mathcal{G} \!=\! \{g_1,\ldots,g_k\}$ be a family of 
homeomorphisms defined on open bounded intervals. If $\{I_1,\ldots,I_k\}$ 
is a family of closed intervals such that 
$I_j \!\subset\! ran(g_j)$ for every $j \! \in \!\{1,\ldots,k\}$, then we say 
that $\mathcal{S} = (\{I_1,\ldots,I_k\},\mathcal{G},P)$ is a 
{\bf\em{Markov system}} for the {\bf\textit{Markov pseudo-group}} 
$\Gamma_{\mathcal{S}}$ generated by the $g_i$'s if the properties below are satisfied:

\vspace{0.1cm}

\noindent (i) \esp $ran(g_i) \cap ran(g_j) = \emptyset$ \esp for every $i \neq j$,

\vspace{0.1cm}

\noindent (ii) \esp if $p_{ij} = 1$ (resp. $p_{ij} = 0$), then $I_j \subset dom (g_i)$ 
and $g_i (I_j) \subset I_i$ (resp. $I_j \cap dom(g_i) = \emptyset$).
\end{defn}

\vspace{0.1cm}

For each $\omega \! = \! (i_1,i_2,\ldots) \in \Omega^*$ and 
each $n \in \mathbb{N}$, the domain of definition of the map 
$\bar{h}_n(\omega) = g_{i_1} \cdots g_{i_n}$ contains the 
interval $g_{i_n}^{-1}(I_{i_n})$. Let 
$I_n(\omega) = \bar{h}_n(\omega) \big( g_{i_n}^{-1}(I_{i_n}) \big)$ and 
$$\Lambda_{\mathcal{S}} = \bigcap_{n \in \mathbb{N}} \bigcup_{\omega \in \Omega^*}
I_n(\omega).$$
Notice that if \esp 
$T\!: \cup_{i=1}^k ran(g_i) \!\rightarrow\! \cup_{i=1}^k dom(g_i)$ \esp 
denotes the map whose restriction to each set $ran(g_i)$ coincides with $g_i^{-1}$, 
then the restriction of $T$ to $\Lambda_{\mathcal{S}}$ is naturally {\bf{\em semiconjugate}} 
to the shift \esp $\sigma\!: \Omega^* \! \rightarrow \! \Omega^*$, \esp in the sense 
that each $x \!\in\! \Lambda_{\mathcal{S}}$ uniquely determines an admissible sequence  
$\varphi(x) = \omega = (i_1,i_2,\ldots)$ such that $x$ belongs to $I_n (\omega)$ for 
each $n$, and the equality \esp $\varphi (T(x)) = \sigma (\varphi(x))$ \esp holds.

\vspace{0.15cm}

\begin{defn} A set $\Lambda$ which is invariant by the action of a pseudo-group 
is said to be {\bf\textit{Markovian}} if there exists an open interval $L$ intersecting 
$\Lambda$ so that $L \cap \Lambda$ equals $\Lambda_{\mathcal{S}}$ for a Markov system 
$\mathcal{S}$ defined on $L$. 
\end{defn}

\vspace{0.15cm}

\index{pseudo-group}

In this chapter, we have already studied an example of a Markov pseudo-group, 
namely the one illustrated by Figure 16. Indeed, the properties in the definition 
are satisfied for the elements $g_1$, $g_2$ corresponding in Figure 16 to $f$, 
$g$, respectively, where $I_1 \!=\! [a,b]$ and $I_2 \! = \! [c,d]$ (we consider 
an open interval slightly larger than $[a,d]$ as the domain of definition 
of $g_1$ and $g_2$). In this particular case, we have $p_{i,j} \!=\! 1$ 
for all $i,j$ in $\{1,2\}$, and thus $\Omega^* = \Omega$.

\vspace{0.12cm}

\begin{thm} {\em If $\Lambda$ is a Markovian local exceptional set for a Markov 
pseudo-group of $\ce^{1+\mathrm{Lip}}$ diffeomorphisms of a one-dimensional 
compact manifold, then the Lebesgue measure of $\Lambda$ is zero.}
\label{th-cc}
\end{thm}

\vspace{0.12cm}

In the case where the maps $g_i$ may be chosen 
(uniformly) differentiably contracting, this result 
still holds in class $\ce^{1+\tau}$ (compare Theorem \ref{th-hur}), but not in class 
$\ce^1$ (see Example \ref{ejem-bowen}). Actually, in the former case a stronger result 
holds \cite{PalT}: The {\em Hausdorff dimension} of $\Lambda$ is less than $1$. 
Nevertheless, this is no longer true in the non-contracting Markovian case, 
even in the real-analytic context \cite{rams}.

We will give the proof of Theorem \ref{th-cc} only for the case of the Markov system 
discussed above and illustrated by Figure 16: the reader should have no problem with  
adapting the arguments below to the general case by taking care of some technical 
details of combinatorial nature. First notice that, for every $x \!\in\! I \!=\! [b,c]$ 
and every element $g = g_{i_n} \cdots g_{i_1} \in \Gamma$ (where each 
$g_{i_j}$ belongs to $\mathcal{G} = \{g_1,g_2\}$), one has 
$$\big| \!\log(g'(b)) - \log(g'(x))  \big| \leq C  
\sum_{j=0}^{n-1}  \big| g_{i_j} \cdots g_{i_1} (b) - g_{i_j} \cdots g_{i_1} (x) \big|
\leq C  \sum_{j=0}^{n-1}  \big| g(I) \big|  \leq  C (d-a).$$
Choosing \esp $x \!\in\! I$ \esp so that \esp $g'(x) = |g(I)| / |I|$, \esp from 
this inequality we conclude that \esp $g'(b) \leq |g(I)| \esp e^{C (d-a)} / |I|$, 
\esp and therefore  
$$\sum_{g \in \Gamma} g'(b) \esp \leq \esp 
\frac{e^{C (d-a)}}{|I|} \sum_{g \in \Gamma} \big| g(I) \big| 
\esp \leq \esp 
\frac{e^{C (d-a)} (d-a)}{|I|} 
\esp = \esp S.$$ 
Letting $\ell = \log(2) / 4CS$, from the estimates in Exercise \ref{mucho} 
it follows that for every $g,h$ in $\Gamma$ and every $x \!\in\! [a,d]$ 
in the $\ell$-neighborhood of $h(b)$ one has \esp 
$g'(x) / 2 \leq g'(b) \leq 2 g'(x)$, \esp and hence 
\begin{equation}
\sum_{g \in \Gamma} \esp \big| g \big( [h(b)-\ell,h(b)+\ell] \big) \big| \esp 
\esp \leq \esp \esp 2 \ell \esp \sum_{g \in \Gamma} \esp \sup_{|x-h(b)| \leq \ell} g'(x) 
\esp \esp \leq \esp \esp 4 \ell S.
\label{rec-min-eq}
\end{equation}
Fix $r \!\in\! \mathbb{N}$ so that for every 
$(i_1,\ldots,i_r) \!\in\! \{ 1,2 \}^r$ the length of the interval 
\esp $g_{i_1} \cdots g_{i_r} ([a,b])$ is smaller than or equal to $\ell$.  
Notice that these intervals cover $\Lambda$, and hence for every $n \!\in\! \mathbb{N}$ 
the same holds for the intervals of the form $g_{i_1} \cdots g_{i_n} h([a,d])$, where 
$(i_1,\ldots,i_n)$ ranges over $\{1,2\}^n$, and $h$ ranges over the elements of 
{\em length} $r$ (that is, the elements of type $g_{i_1} \cdots g_{i_r}$, 
with $(i_1,\ldots,i_r) \in \{1,2\}^r$). Therefore, 
$$Leb(\Lambda) \esp \esp \esp \leq 
\sum_{(i_1,\ldots,i_n) \in \{1,2\}^n} \esp \esp  
\sum_{(h_{1},\ldots,h_{r}) \in \{g_1,g_2\}^r} \!\!\!
\big| g_{i_1} \cdots g_{i_n} h_{1} \cdots h_{r}([a,d]) \big|,$$
and since for each $(h_1,\ldots,h_r) \in \{g_1,g_2\}^r$ the point 
$h_1 \cdots h_r(b)$ belongs to $h_1 \cdots h_r ([a,d])$, from 
this we conclude that 
$$Leb(\Lambda) \esp\esp\esp 
\leq \sum_{(h_1,\ldots,h_r) \in \{g_1,g_2\}^r} \esp \esp \esp 
\sum_{\lh (g) = n} \!\!\!
\big| g \big( [h_1 \cdots h_r (b) - \ell,h_1 \cdots h_r (b) + \ell] \big) \big|.$$
Letting $n$ go to infinity in this inequality, from the convergence of the series in 
(\ref{rec-min-eq}) one concludes that \esp $Leb(\Lambda) = 0$, \esp thus finishing the proof. 

\vspace{0.35cm}

We refer the reader to \cite{CC2} and \cite{inaba} for the realization of 
Markov pseudo-groups as holonomy pseudogroups of codimension-one foliations.
\index{foliation} Let us point out, however, that 
\index{holonomy} it is not difficult to construct exceptional minimal 
sets which are not induced by Markovian systems: see \cite{CC4}. 

\begin{small} \begin{ejer}
Completing Exercise \ref{finir}, show that if an exceptional minimal set 
$\Lambda$ is Markovian, then \textit{all} the orbits in $\Lambda$ have 
infinitely many ends (see \cite{CC3} in case of problems with this). 
\end{ejer} \end{small}

To close this section, let us briefly discuss another approach to Problem 2, which is 
closely related to our approach to the ergodicity question. Let us begin by pointing 
out that Theorem \ref{principal} has a natural analogue (due to Hurder 
\index{Hurder} 
\cite{Hur}) in the present context.

\vspace{0.1cm}

\begin{thm} {\em Let $\Gamma$ be a subgroup of \esp \! $\mathrm{Diff}_+^{1+\tau}(\clo)$ 
admitting an exceptional minimal set $\Lambda$. If \esp $\tau \!>\! 0$, then 
the set \esp $\Lambda \cap E(\Gamma)$ \esp has null Lebesgue measure.}
\label{th-hur}
\end{thm}

\vspace{0.14cm}

The necessity of the hypothesis $\tau \! > \! 0$ for this result is 
illustrated by the following classical example due to Bowen \cite{bowen}.

\begin{small} \begin{ejem} As the reader can easily check, the maps $f,g$ 
in Figure 16 may be chosen so that the following supplementary conditions 
are verified (to put ourselves in the Markovian context, we denote again 
$g_1 \!=\! f$, $g_2 \!=\! g$, and $I\! =\! [b,c]$, with $a\!=\!0$ and 
$d\!=\!1$).

\vspace{0.03cm}

\noindent (i) There exists a sequence of positive real numbers 
$\ell_n$ (where $n \geq 0$) such that $|I| = \ell_0$, 
$$\sum_{n \geq 0} \ell_n < 1, \qquad \lim_{n \rightarrow \infty}
\frac{\ell_{n+1}}{\ell_n} = 1,$$
and such that if for $n \!\in\! \mathbb{N}$ and  
$(g_{i_1},\ldots,g_{i_n}) \!\in\! \{g_1,g_2\}^n$ we denote 
$I_{i_1,\ldots,i_n} \!=\! g_{i_1} \cdots g_{i_n} (I)$, then 
\esp $|I_{i_1,\ldots,i_n}| = \ell_n / 2^n$.

\vspace{0.03cm}

\noindent (ii) Each $g_i$ is smooth over $I$ and over each $I_{i_1,\ldots,i_n}$, 
its derivative at the endpoints of these intervals equals $1/2$, and 
$$\lim_{n \rightarrow \infty} \esp 
\max_{(i_1,\ldots,i_n) \in \{1,2\}^n} 
\esp \sup_{x \in I_{i_1,\ldots,i_n}} \esp 
\Big| g_i'(x) - \frac{1}{2} \Big| \esp = \esp 0.$$

\vspace{0.03cm}

\noindent (iii) Each $g_i$ is differentiable on a neighborhood of $a=0$ and $d=1$, with 
derivative identically equal to $1/2$ to the left (resp. to the right) of $a$ (resp. $d$).

\vspace{0.03cm}

With these conditions, it is not difficult to show that $g_1$ and $g_2$ are of class 
$\ce^1$ over the whole interval $[a,d] = [0,1]$, and their derivatives equal $1/2$ 
over the Markov minimal set $\Lambda$. Now for the Lebesgue measure of $\Lambda$ 
we have
$$Leb(\Lambda) \esp = \esp 1 - |I| - \sum_{n \geq 1} \sum_{(i_1,\ldots,i_n)}
|I_{i_1,\ldots,i_n}|\esp = \esp 1 - \sum_{n \geq 0} \ell_n \esp > \esp 0.$$
Notice that, using this example, it is not difficult to create a 
subgroup of $\mathrm{Diff}_+^1(\clo)$ with two generators and having an 
exceptional minimal set of positive Lebesgue measure entirely contained 
in the exponential set of the action.
\label{ejem-bowen}
\end{ejem}

\begin{ejer} After reading \S \ref{ejemplosceuno}, show that for every $\tau > 0$ 
there exist finitely generated Abelian subgroups of $\mathrm{Diff}^{1+\tau}_+(\clo)$ 
admitting an exceptional minimal set of positive measure. Prove directly ({\em i.e.,} 
without using Theorem \ref{th-hur}) that, for these examples, the exceptional minimal 
set is contained (up to a null measure set) in the sub-exponential set of the action.
\index{sub-exponential set of an action}
\end{ejer}

\begin{ejem} Following an idea which seems to go back to Ma\~n\'e, 
given $\tau < 1$ let us consider a $\ce^{1+\tau}$ circle diffeomorphism with 
irrational rotation number and admitting a minimal invariant Cantor set $\Lambda_0$ 
of positive Lebesgue measure (see the preceding exercise). Fix a connected component 
$I$ of the complement of $\Lambda_0$, and let $\bar{H}: \clo \rightarrow \clo$ 
be a degree-2 
\index{degree!of a map}
map which coincides with $f$ outside $I$ (\textit{i.e.,} the map $\bar{H}$ 
makes an ``extra turn'' around the circle over $I$). Using $\bar{H}$, 
it is easy to construct a pseudo-group with two generators admitting an 
exceptional minimal set $\Lambda$ containing $\Lambda_0$. Since the 
Lebesgue measure of $\Lambda_0$ is positive, this is 
also the case of $\Lambda$.
\label{mane}
\end{ejem}

\begin{ejer} Prove that if $\Gamma$ is a finitely generated group of real-analytic 
circle diffeomorphisms admitting an exceptional minimal set $\Lambda$ then, excepting 
at most countably many points, $\clo \setminus \Lambda$ is contained in the sub-exponential 
set $S(\Gamma)$.

\vspace{0.07cm}

\noindent{\underbar{Remark.}} Quite possibly, the same holds 
for subgroups of $\mathrm{Diff}_+^{1+\mathrm{bv}}(\clo)$.\\

\vspace{0.07cm}

\noindent{\underbar{Hint.}} Use the analyticity to show that, discarding at most 
countably many points, if $x$ belongs to $\clo \setminus \Lambda$ then there exists 
an open interval $I_x$ containing $x$ in its interior such that the images of $I_x$ 
by the elements in the group are two-by-two disjoint. Then use the relations 
$$1 \esp \geq \esp \sum_{g \in \Gamma} |g(I_x)| \esp = \esp \sum_{g \in \Gamma} \int_{I_x} 
g'(y) dy \esp = \esp \int_{I_x} \Big( \sum_{g \in \Gamma} g'(y) \Big) dy.$$
\end{ejer} \end{small}


\section{On the Smoothness of the Conjugacy between Groups of Diffeomorphisms}
\label{poquito}

\hspace{0.45cm} The problem of the smoothness for the conjugacy between groups of circle 
diffeomorphisms is technical and difficult. The relevant case of free actions is 
already extremely hard. Nevertheless, this case is very well understood thanks to  
the works of Siegel, Arnold, Herman, Moser, Yoccoz, Khanin, Katznelson and Orstein, 
\index{Yoccoz}\index{Herman}among others. For $r > 2$, the main technical tool for 
obtaining the differentiability of the conjugacy\footnote{Notice that from Exercise 
\ref{uni-inv} one easily concludes that two different conjugacies as above differ 
by an Euclidean rotation. Due to this, there is no ambiguity when speaking about 
the differentiability of the conjugacy.} of a $\mathrm{C}^r$ circle 
diffeomorphism $f$ having irrational rotation number \index{rotation number} 
$\rho(f)$ to the rotation of angle $\rho(f)$ 
corresponds to the Diophantine nature of $\rho(f)$. Roughly, the conjugacy is forced 
to be smooth if the approximations of $\rho(f)$ by rational numbers are ``bad'', in 
the sense that they are slow with respect to the denominator of the approximating rational. 
This corresponds to one of the main issues of the so-called Small Denominators Theory. 
\index{small denominators}

In what follows, we will study the case of non-free actions by circle diffeomorphisms, which 
is essentially different. There is no systematic theory for this case and, perhaps, the only 
definitive results correspond to Theorem \ref{con-jugando} to be studied in detail in the 
next section, and Proposition \ref{conj-caca} in the next chapter. However, we point out  
that these results are much simpler than those of the theory of small denominators. 
Slightly stronger results hold in some particular cases. For instance, Sullivan 
\index{Sullivan}
showed that every 
topological and absolutely continuous conjugacy between Fuchsian groups of first kind is 
real-analytic, and he further obtained analogous results on the smoothness over the exceptional 
minimal set for conjugacies between Fuchsian groups of second kind \cite{Su1,Su2}. On the 
other hand, using classical results on the existence and uniqueness of absolutely continuous 
invariant measures for expanding maps of the interval, the author showed in \cite{Na3} that if 
two groups of $\mathrm{C}^r$ circle diffeomorphisms, with $r \geq 2$, satisfy the hypothesis 
of Duminy's first theorem (\textit{i.e.,} if they are generated by elements near rotations), 
and if their orbits are dense, then every topological and absolutely continuous conjugacy 
between them is a $\mathrm{C}^r$ diffeomorphism (see \cite{Re1,Re2} for a stronger 
result in the real-analytic context).


\subsection{Sternberg's linearization theorem and $\mathbf{\ce^1}$ conjugacies
\index{linearization}}
\index{Sternberg!linearization theorem|(}
\label{sec-st}

\hspace{0.23cm} Let $f$ and $g$ be $\ce^r$ diffeomorphisms of a neighborhood of 
the origin in the line into their images. If $f$ and $g$ fix the origin, we say 
that they are {\bf {\em equivalent}} if there exists $\varepsilon > 0$ such that 
$f|_{]-\varepsilon,\varepsilon[} = g|_{]-\varepsilon,\varepsilon[}$. Modulo 
this equivalence relation, the class of $g$ will be denoted by $[g]$. The set of 
classes forms a group with respect to the composition of representatives, that is, 
$[f][g] = [f \circ g]$. This group, called the group of {\bf \em{germs}} of 
$\mathrm{C}^r$ diffeomorphisms of the line which fix the origin and preserve 
orientation, will be denoted by $\mathcal{G}^r_{+}(\mathbb{R},0)$.
\index{germ of a diffeomorphism}

Notice that the derivative $g^{(i)}(0)$ of order $i \!\leq\! r$ at the origin is well-defined 
for all $[g] \in \ger$. We will say that $[g]$ is {\bf {\em hyperbolic}} if $g'(0) \neq 1$. 
\index{hyperbolic!germ}

\vspace{0.1cm}

\begin{lem} {\em If $[g] \!\in\! \ger$ is hyperbolic, with $1 \leq r < \infty$, then 
there exists $[h]$ in $\ger$ such that $h'(0) = 1$ and $(hgh^{-1})^{(i)}(0) = 0$ 
for every $2 \leq i \leq r$.}
\label{yococito}
\end{lem}

\noindent{\bf Proof.} Let $g(x) = ax + a_2x^2 + \ldots + a_r x^r + \mathrm{o} (x^r)$ 
be the Taylor series expansion of $g$ about the origin. Let us formally write 
$$\hat{h}(x) = x + b_2x^2 + \ldots + b_rx^r + \ldots,$$
and let us try to find the coefficients $b_i$ so that 
\begin{equation}
\hat{h} \circ g = M_{a}  \circ \hat{h},	
\label{aconjugar}
\end{equation}
where $M_a$ denotes the multiplication by $a \!=\! g'(0)$. If we identify the coefficients of 
$x^2$ in both sides of this equality we obtain $a_2 \!+\! b_2a^2 \!=\! ab_2$, and therefore 
$b_2 \!=\! a_2/(a-a^2)$. In general, assuming that $b_2, \ldots, b_{i-1}$ are already known 
(where $i \leq r$), from (\ref{aconjugar}) one obtains 
$$Q_i(b_1, \ldots, b_{i-1}, a_1, \ldots, a_r) + b_ia^i = ab_i$$
for some polynomial $Q_i$ in $(r+i-1)$ variables. Thus, 
$$b_i = \frac{Q_i(b_1, \ldots, b_{i-1}, a_1, \ldots, a_r)}{a - a^i}.$$
Now let $h$ be a $\mathrm{C}^r$ diffeomorphism defined on a neighborhood of the origin and 
such that $h(0)=0$, $h'(0) = 1$, and $h^{(i)}(0) = i! \esp b_i$ for every $2 \leq i \leq r$. 
By reversing the above computations, it is easy to see that $[h] \in \ger$ satisfies 
the desired properties. $\hfill\square$

\vspace{0.35cm}

Sternberg's theorem appears in the form below in \cite{Yo1} (see \cite{St} for the 
original, weaker form). Let us point out that an analogous result holds for germs of 
real-analytic diffeomorphisms \cite{CG}.

\vspace{0.1cm}

\begin{thm} {\em Let $g \!\in\! \ger$ be a hyperbolic germ, where \esp $2 \!\leq\! r \!\leq\! \infty$. 
\esp If we denote \esp $a \!=\! g'(0)$, \esp then there exists $[h] \!\in\! \ger$ such that \esp 
$h'(0) \!=\! 1$ \esp and \esp $h(g(x)) \!=\! a h(x)$ \esp for every $x$ near the origin. Moreover, 
if $[h_1] \!\in\! \mathcal{G}^1_{+}(\mathbb{R},0)$ satisfies the last two properties, then $[h_1]$ 
belongs to $\ger$, and $[h_1] = [h].$}
\end{thm}

\noindent{\bf Proof.} Let us first consider the case $r < \infty$. By the lemma above, to obtain 
a conjugacy we may suppose that $g^{(i)}(0) = 0$ for each $i \in \{ 2,\ldots,r \}$. Changing $g$ by 
$g^{-1}$ if necessary, we may suppose moreover that $g'(0) = a < 1$. Let $0 < \delta < 1$ be such 
that the domain of definition of $g$ contains the interval $[-\delta,\delta]$. Let us define 
$$C(\delta) = \sup \{ |g^{(r)}(t)| \!: t \in [-\delta,\delta]\}.$$
A direct application of the Mean Value Theorem shows that, for all  
$t \!\in\! [-\delta,\delta]$, one has $|g'(t)| \leq a + C(\delta)$ 
and $|g^{(i)}(t)| \leq C(\delta)$ for every $i \in \{ 2, \ldots, r \}$.

Let \esp $E_\delta$ \esp be the space of real-valued functions \esp $\psi$ 
\esp of class \esp $\ce^r$ \esp on \esp $[-\delta,\delta]$ \esp satisfying 
\esp $\psi(0) \!=\! \psi'(0)\! = \!\ldots\! =\! \psi^{(r)}(0)\! =\! 0$. \esp 
Endowed with the norm
$$\| \psi \| = \sup \{ |\psi^{(r)}(t)| \!: t \in [-\delta,\delta] \},$$
$E_{\delta}$ becomes a Banach space. For \esp $\psi \in E_\delta$, \esp 
$t \in [-\delta,\delta]$, \esp and \esp $i \in \{ 0, \ldots, r \}$, one has 
$$|\psi^{(i)}(t)| \leq \frac{t^{r-i}}{(r-i)!} \| \psi \|.$$
Let us consider the linear operator $S_\delta \!: E_\delta \rightarrow E_\delta$ defined by 
$S_\delta(\psi) = ( \psi \circ g ) / a$. We claim that if $\delta>0$ is small enough, then 
$S_\delta$ is a contraction (that is, $\|S_\delta\| < 1$). Indeed, it is easy to check that 
$$(\psi \circ g)^{(r)} = \sum_{i=1}^{r} \big( \psi^{(i)} \circ g \big)
\cdot Q_i \big( g', \ldots, g^{(r-i+1)} \big),$$
where $Q_i$ is a polynomial with positive coefficients in $(r + 1 - i)$ 
variables and $Q_r(x) = x^r$. For $t \in [-\delta,\delta]$ this yields 
$$ \big| (\psi \circ g)^{(r)}(t) \big| \leq K(\delta) \hspace{0.1cm} \|\psi\|,$$ 
where 
$$K(\delta) = \big( a + C(\delta) \big)^r +
\sum_{i=1}^{r-1}\frac{\delta^{r-i}}{(r-i)!}
Q_i \big( a+C(\delta), \ldots, a+C(\delta) \big).$$
Notice that $K(\delta)$ tends to $a^r$ as $\delta$ goes to zero. 
Moreover, \esp $\|S_{\delta}\| \leq K(\delta) / a$. \esp Since  
\esp $r \!\geq\! 2$ \esp and $a < 1$, \esp this shows that the value 
of $\|S_\delta\|$ is smaller than $1$ for $\delta$ small enough.\\

Now fix $\delta > 0$ such that $\|S_\delta\| < 1$. Notice that the restriction of the map 
$x \mapsto g(x) - ax$ to the interval $[-\delta,\delta]$ defines an element $\psi_1$ in 
$E_\delta$. Since $S_\delta$ is a contraction, the equation (in the variable $\psi$) 
$$S_\delta(\psi) + a^{-1} \psi_1 = \psi$$
has a unique solution $\psi_0 \in E_\delta$. For $h = Id + \psi_0$ we have 
\begin{eqnarray*}
h (g (x)) \esp = \esp g(x) + \psi_0 \circ g(x)  
\!&=&\!  \psi_1(x) + ax + \psi_0 \circ g(x)\\
\!&=&\!  \psi_1(x) + ax + aS_\delta(\psi_0)(x) 
\esp = \esp a\psi_0(x) + ax \esp = \esp ah(x).
\end{eqnarray*}
Since $\psi_0 \!\in\! E_\delta$, one has $h(0) \!=\! 0$ and $h'(0) \!=\! 1$. Therefore, 
if $r \!<\! \infty$, then $[h]$ is a germ in $\ger$ satisfying the required properties.

If $h_1 \in \mathcal{G}_{+}^{1}(\mathbb{R},0)$ verifies the same properties, then
$$ahh_1^{-1}(t) \esp = \esp h \circ g \circ h_1^{-1} (t) \esp = \esp h h_1^{-1} (at)$$
for every $t$ near the origin. From this one deduces that 
$$hh_1^{-1} (t) \esp = \esp \lim\limits_{n \rightarrow \infty} \frac{hh_1^{-1}(a^nt)}{a^n}
\esp = \esp t \lim\limits_{n \rightarrow \infty} \frac{hh_1^{-1}(a^n t)}{a^n t} \esp  
= \esp t \esp (hh_1^{-1})'(0) \esp = \esp t,$$
which shows the uniqueness.

Finally, the case $r = \infty$ easily follows from the uniqueness 
already established for each finite $r \geq 2$. $\hfill\square$

\vspace{0.35cm}

A direct consequence of the preceding theorem is the following result.

\vspace{0.1cm}

\begin{cor} {\em Let $g_1$ and $g_2$ be elements in $\ger$, where $2 \!\leq\! r \!\leq\! \infty$. 
Suppose that $[\varphi] \!\in\! \mathcal{G}_{+}^1(\mathbb{R},0)$ conjugates $g_1$ to $g_2$, that 
is, on a neighborhood of the origin one has $\varphi \circ g_1 = g_2 \circ \varphi$. If $g_1$ is 
hyperbolic, then $[\varphi]$ is an element of $\ger$.}
\label{corolarito}
\end{cor}

\noindent{\bf Proof.} Let us begin by noticing that $[g_2]$ is also hyperbolic, 
because $g_2'(0) = g_1'(0)$. If $h \in \ger$ conjugates $g_2$ to its linear part, then 
$h \circ \varphi$ conjugates $g_1$ to its linear part. If $b = (h \varphi)'(0)$, then 
$M_{1/b} h \varphi$ still conjugates $g_1$ to its linear part; moreover, it satisfies 
$(M_{1/b} h \varphi)'(0) \!=\! 1$. By the uniqueness of such a conjugacy, we have 
$[M_{1/b} h \varphi] \in \ger$, and therefore $[\varphi] \in \ger$. $\hfill\square$

\vspace{0.45cm}

Sternberg's linearization theorem still holds in class $\ce^{1+\tau}$ for every 
$\tau > 0$ (the proof consists in a straightforward extension of the preceding 
arguments; see also \cite{chap}). However, the theorem is no longer true in 
class $\ce^1$, as is illustrated by the next example due to Sternberg himself.
\index{Sternberg!example}

\begin{small}\begin{ejem}
Let us consider the map $g$ defined on an open interval about the origin by $g(0)=0$ 
and \esp $g(x) = a x (1 - 1/\log(x))$ \esp for $x \neq 0$, where $0 < a < 1$. It is 
easy to see that $[g]$ is the germ of a $\ce^1$ diffeomorphism satisfying $g'(0) = a$. 
Despite this fact, $[g]$ is not conjugate to the germ of the linear map $M_a$ by 
any germ of a Lipschitz homeomorphism 
\index{Lipschitz!conjugacy}
(in particular, $[g]$ is not $\ce^1$ conjugate to $[M_a]$). 
To show this, fix a constant $\bar{a} \in ]0,1[$ so that $g'(x) \geq \bar{a}$ for every 
$x$ near the origin, and choose a local homeomorphism $h$ fixing $0$ and such that for 
these points $x$ one has $g(x) = h M_a h^{-1}(x)$. For every $k \in \mathbb{N}$ one 
has $g^k (x) = h \big( a^k h^{-1}(x) \big)$, and hence
\begin{equation}
\frac{g^k(x)}{a^k} \esp = \esp \frac{h(a^k h^{-1}(x))}{a^k}.
\label{act}
\end{equation}
If $h$ were Lipschitz with constant $\bar{C}$ then, on the one hand, 
the right-hand member in the above equality would be bounded from 
above by $\bar{C} \esp \! h^{-1}(x)$ for every $k \!\in\! \mathbb{N}$. 
On the other hand, from the definition of $g$ one easily concludes that 
$$\frac{g^k(x)}{a^k} \esp = \esp x \prod_{i=0}^{k-1} \left( 1 -
\frac{1}{\log \big( g^i(x) \big)} \right).$$
Since \esp $g^i(x) \geq \bar{a}^i x$, \esp there exists an integer 
$i_0 \geq 1$ such that, for some positive constant $C$ and every $i \geq i_0$, 
$$\log \big( g^i(x) \big) \esp \geq \esp 
i \log(\bar{a}) + \log(x) \esp \geq \esp -\frac{i}{C}.$$
From this it follows that
$$\prod_{i=i_0}^{k-1} \left( 1 - \frac{1}{\log \big( g^i(x) \big)} \right) 
\esp \geq \esp \prod_{i=i_0}^{k-1} \left( 1 + \frac{C}{i} \right).$$
Since the last product diverges as $k$ goes to infinity, we conclude that the 
left-hand member in (\ref{act}) is unbounded, thus giving a contradiction
\label{ejst}
\end{ejem}

\begin{ejer} Show that the germ at the origin in Example \ref{ejem-bowen} is 
(hyperbolic and) non linearizable by the germ of any $\ce^1$ diffeomorphism.
\end{ejer}\end{small}

\vspace{0.06cm}

As an application of Sternber's theorem and Corollary 
\ref{remarque-etiene}, we now reproduce (with a much simpler 
proof than the original one) a result obtained by Ghys and 
Tsuboi in \cite{GT}.
\index{Ghys}\index{Tsuboi}

\vspace{0.1cm}

\begin{thm} {\em Let $\Phi_1$ and $\Phi_2$ be two representations of a group $\Gamma$
in $\mathrm{Diff}_{+}^{r}(\clo)$, where $2 \leq r \leq \infty$. Suppose that their 
orbits are dense and that $\Phi_1(\Gamma)$ is not topologically conjugate to a group 
of rotations. If a $\mathrm{C}^1$ circle diffeomorphism $\varphi$ conjugates 
$\Phi_1$ and $\Phi_2$, then $\varphi$ is a $\mathrm{C}^r$ diffeomorphism.}
\label{con-jugando}
\end{thm}

\noindent{\bf Proof.} By Corollary \ref{remarque-etiene}, there exist $g \!\in\! \Gamma$ 
and $x \!\in\! \clo$ such that $\Phi_1(g)(x) \!=\! x$ and $\Phi_1(g)'(x) \!\neq\! 1$. 
In local coordinates we obtain a hyperbolic germ of a diffeomorphism fixing the 
origin. The point $\varphi(x) \!\in\! \clo$ is fixed by 
$\Phi_2(g) \!=\! \varphi \circ \Phi_1(g) \circ \varphi^{-1} \in \Gamma_2$.
Therefore, $\varphi$ induces a conjugacy between hyperbolic germs. Corollary 
\ref{corolarito} then implies that $\varphi$ is of class $\mathrm{C}^r$ in 
a neighborhood of $x$. Since $\varphi \circ \Phi_1(h) = 
\Phi_2(h) \circ \varphi$ holds for every  $h \in \Gamma$, the set 
of points around which $\varphi$ is of class $\mathrm{C}^r$ 
is invariant by $\Phi_1(\Gamma)$, and since the orbits by 
$\Phi_1(\Gamma)$ are dense, $\varphi$ must be of class $\mathrm{C}^r$ on the 
whole circle. To show that $\varphi^{-1}$ is of class $\mathrm{C}^r$ on $\clo$, 
it suffices to interchange the roles of $\Phi_1$ and $\Phi_2$. $\hfill\square$

\vspace{0.1cm}

\begin{small} \begin{obs} The theorem above still holds for bi-Lipschitz conjugacies, provided 
that we assume \textit{a priori} \esp that the actions are ergodic (see Theorem \ref{japon1}). 
Moreover, the theorem is also true for conjugacies between groups of $\ce^{1+\tau}$ 
circle diffeomorphisms which are non-Abelian and act minimally.
\label{st-obs}
\end{obs}

\begin{ejer} Prove that Theorem \ref{con-jugando} still holds (in class $\ce^2$ or $\ce^{1+\tau}$) for 
groups admitting an exceptional minimal set\index{exceptional minimal set} 
and which are non semiconjugate to groups of rotations.\\
\index{rotation!group}

\noindent{\underbar{Hint.}} Apply the same argument as before having in 
mind that every orbit must accumulate on the minimal invariant Cantor set.

\noindent{\underbar{Remark.}} This result does not generalize 
to bi-Lipschitz conjugacies: see Theorem \ref{japon2}.
\end{ejer}
\end{small}

\index{Sternberg!linearization theorem}


\subsection{The case of bi-Lipschitz conjugacies}
\label{conja-sec}

\hspace{0.45cm} The main object of this section is to extend Theorem \ref{con-jugando} 
to bi-Lipschitz 
\index{Lipschitz!conjugacy} 
conjugacies under the hypothesis of ergodicity for the action (see 
\S \ref{sec-minimal-inv} for a discussion about this hypothesis).

\vspace{0.1cm}

\begin{thm} {\em Let $\Phi_1$ and $\Phi_2$ be two representations of a finitely generated 
group $\Gamma$ into $\mathrm{Diff}_{+}^{r}(\clo)$, where $2 \!\leq\! r \leq \infty$. 
Suppose that $\Phi_1(\Gamma)$ acts minimally and ergodically. If a bi-Lipschitz circle 
homeomorphism $\varphi$ conjugates $\Phi_1$ and $\Phi_2$, then $\varphi$ is a $\mathrm{C}^1$ 
diffeomorphism. Moreover, if $\Phi_1(\Gamma)$ is not topologically conjugate 
to a group of rotations, then $\varphi$ is a $\ce^r$ diffeomorphism.}
\label{japon1}
\end{thm}

\vspace{0.1cm}

The proof of this result uses an equivariant version of a classical lemma of 
cohomological flavor due to Gottschalk and Hedlund. 
\index{Gottschalk-Hedlund Lemma}
Let M be a compact metric space and $\Gamma$ a finitely 
generated group acting on it by homeomorphisms. A {\bf {\em cocycle}} 
\index{cocycle} 
associated to this action (compare \S \ref{def-T}) is a map 
$c\!: \Gamma \times \mathrm{M} \rightarrow \mathbb{R}$ such that for each 
fixed $f \!\in\! \Gamma$ the map \esp $x \mapsto c(f,x)$ \esp is continuous, 
and such that for every $f,g$ in $\Gamma$ and every $x \in \mathrm{M}$ one has
\begin{equation}
c(fg,x) = c(g,x) + c(f,g(x)).
\label{coc-GH}
\end{equation}
With this notation, the Equivariant Gottschalk-Hedlund 
Lemma may be stated as follows.

\vspace{0.15cm}

\begin{lem} {\em Suppose that $\Gamma$ is finitely generated and its action 
on $\mathrm{M}$ is minimal. Then the following are equivalent:}

\vspace{0.12cm}

\noindent{(i) \em there exist $x_0 \in \mathrm{M}$ and a constant 
$C_0 > 0$ such that $|c(f,x_0)| \leq C_0$ for every $f \in \Gamma$,}

\vspace{0.12cm}

\noindent{(ii) \em there exists a continuous function $\phi: \mathrm{M} 
\rightarrow \mathbb{R}$ such that $c(f,x) = \phi (f(x)) - \phi(x)$ for all 
$f \in \Gamma$ and all $x \in \mathrm{M}$.}
\end{lem}

\noindent{\bf Proof.} If the second condition is satisfied, then (i) follows from
$$|c(f,x_0)| \leq |\phi(f(x_0))| + |\phi(x_0)|\leq 2 \| \phi \|_{C^0}.$$

Conversely, let us suppose that the first condition holds. For each $f \in \Gamma$ 
consider the homeomorphism $\hat{f}$ of the space $\mathrm{M} \times \mathbb{R}$ 
defined by $\hat{f}(x,t) = (f(x),t+c(f,x))$. It is easy to see that the cocycle 
relation (\ref{coc-GH}) implies that this defines a group action of $\Gamma$ on 
$\mathrm{M} \times \mathbb{R}$, in the sense that $\hat{f} \hat{g} = \widehat{fg}$ 
for all $f,g$ in $\Gamma$. Moreover, condition (i) implies that the orbit of 
the point $(x_0,0)$ under this action is bounded; in particular, its closure 
is a (non-empty) compact, invariant set. One may then apply the Zorn Lemma 
\index{Zorn Lemma}
to deduce the existence of a non-empty, minimal invariant, compact 
subset $\Lambda$ of $\mathrm{M} \times \mathbb{R}$. We claim that  
$\Lambda$ is the graph of a continuous real-valued function on M.

First of all, since the action of $\Gamma$ on M is minimal, the projection 
of $\Lambda$ into M is the whole space. Moreover, if $(\bar{x},t_1)$ and 
$(\bar{x},t_2)$ belong to $\Lambda$ for some $\bar{x} \in \mathrm{M}$ and 
some $t_1 \neq t_2$, then this implies that $\Lambda \cap \Lambda_t \neq 
\emptyset$, where $t = t_2 - t_1$ and $\Lambda_t = \{(x,s+t) \!: \esp (x,s) 
\in \Lambda \}$. Notice that the action of $\Gamma$ on $\mathrm{M} \times \mathbb{R}$ 
commutes with the map $(x,s) \mapsto (x,s+t)$; in particular, $\Lambda_t$ is also 
invariant. Since $\Lambda$ is minimal, this implies that $\Lambda = \Lambda_t$. 
One then concludes that \esp $\Lambda = \Lambda_t = \Lambda_{2t} = \ldots$, 
which is impossible because $\Lambda$ is compact and $t \!\neq\! 0$.

We have then proved that, for every $x \in \mathrm{M}$, the set $\Lambda$ contains 
exactly one point of the form $(x,t)$. Putting $\phi(x) = t$, one obtains a 
function from M into $\mathbb{R}$. This function must be continuous, since 
its graph (which coincides with $\Lambda$) is compact. 

Finally, since the graph of $\phi$ is invariant by the action, 
for all $f \in \Gamma$ and all $x \in \mathrm{M}$ the point $\hat{f}(x,\phi(x)) 
= (f(x),\phi(x) + c(f,x))$ must be of the form $(f(x),\phi(f(x)))$, 
which implies that \esp $c(f,x) = \phi(f(x)) - \phi(x)$. $\hfill\square$

\vspace{0.5cm}

The following corresponds to a ``measurable version'' of the preceding lemma.

\vspace{0.15cm}

\begin{lem} {\em Let $\mathrm{M}$ be a compact metric space and $\Gamma$ a finitely generated 
group acting on it by homeomorphisms. Suppose that the action of $\Gamma$ on $\mathrm{M}$ is 
minimal and ergodic with respect to some probability measure $\mu$, and let $c$ be a cocycle 
associated to this action. If $\phi$ is a function in 
$\mathcal{L}^{\infty}_{\mathbb{R}}(\mathrm{M},\mu)$ 
such that for all $f \in \Gamma$ and $\mu$-almost every $x \in \mathrm{M}$ one has} 
\begin{equation}
c(f,x) = \phi (f(x)) - \phi(x),
\label{cacha}
\end{equation}
{\em then there exists a continuous function $\tilde{\phi} \!: 
\mathrm{M} \rightarrow \mathbb{R}$ which coincides $\mu$-a.e. with 
$\phi$ and such that, for all $f \in \Gamma$ and all $x \in \mathrm{M}$,}
\begin{equation}
c(f,x) = \tilde{\phi} (f(x)) - \tilde{\phi} (x).
\label{cacha2}
\end{equation}
\label{cle}
\end{lem}

\vspace{-0.5cm}

\noindent{\bf Proof.} Let $\mathrm{M}_0$ be the set of points for which (\ref{cacha}) 
does not hold for some $f \in \Gamma$. Since $\Gamma$ is finitely generated, 
$\mu(\mathrm{M}_0) = 0$. Let $\mathrm{M}_1'$ the complementary set of 
the essential support 
\index{support!of a function}
of $\phi$, and let $\mathrm{M}_1 = \cup_{f \in \Gamma} f(\mathrm{M}_1')$. Take a point 
$x_0$ in the full measure set $\mathrm{M} \setminus (\mathrm{M}_0 \cup \mathrm{M}_1)$. 
Equality (\ref{cacha}) then gives \esp 
$|c(f,x_0)| \leq 2 \| \phi \|_{\mathcal{L}^\infty}$ \esp for all $f \in \Gamma$. 
By the preceding lemma, there exists a continuous function 
$\tilde{\phi}\!: \mathrm{M} \rightarrow \mathbb{R}$ such that (\ref{cacha2}) 
holds for {\em every} \esp $x$. This implies that $\mu$-a.e. one has  
$$\tilde{\phi} \circ f - \tilde{\phi} = \phi \circ f - \phi,$$
hence 
$$\tilde{\phi} - \phi = (\tilde{\phi} - \phi) \circ f.$$
Since the action of $\Gamma$ on $\mathrm{M}$ is assumed to be $\mu$-ergodic, 
the difference $\tilde{\phi} - \phi$ must be $\mu$-a.e. equal to a constant 
$C$. To conclude the proof, just change $\tilde{\phi}$ by $\tilde{\phi} - C$. 
$\hfill\square$

\vspace{0.28cm}

We may now pass to the proof of Theorem \ref{japon1}. First notice that, if $\varphi$ is a 
bi-Lipschitz circle homeomorphism conjugating the actions $\Phi_1$ and $\Phi_2$ of our group 
$\Gamma$, then $\varphi$ and $\varphi^{-1}$ are almost everywhere differentiable, and their 
derivatives belong to $\mathcal{L}^{\infty}_{\mathbb{R}}(\clo,Leb)$. Therefore, the function 
$x \mapsto \log(\varphi'(x))$ is also in $\mathcal{L}^{\infty}_{\mathbb{R}}(\clo,Leb)$. The 
relation \esp $\Phi_1 (f) = \varphi^{-1} \circ \Phi_2(f) \circ \varphi$ \esp gives almost 
everywhere the equality 
$$\log(\Phi_1(f)'(x)) = \log(\varphi'(x)) - \log(\varphi'(\Phi_1(f)(x))) 
+ \log(\Phi_2(f)'(\varphi(x))).$$
Putting $\phi = - \log(\varphi')$ and $c(f,x) = \log(\Phi_1(f)'(x)) - 
\log(\Phi_2(f)'(\varphi(x)))$ this yields, for all $f \in \Gamma$ 
and almost every $x \in \mathrm{S}^1$, 
$$c(f,x) = \phi(\Phi_1(f)(x)) - \phi(x).$$
Using the relation \esp $\Phi_2 (g) = \varphi \circ \Phi_1 (g) \circ \varphi^{-1}$, 
\esp one easily checks the cocycle identity
$$c(fg,x) = c(g,x) + c(f,\Phi_1(g)(x)).$$
Since the action $\Phi_1$ is supposed to be ergodic, Lemma \ref{cle} 
ensures the existence of a continuous function $\tilde{\phi}$ which 
coincides almost everywhere with $\phi$ and such that (\ref{cacha2}) 
holds for every $x$ and all $f$. By integration, one concludes that 
the derivative of $\varphi$ is well-defined everywhere and coincides 
with $\exp(-\tilde{\phi})$. In particular, $\varphi$ is of 
class $\ce^1$. Changing the roles of $\Phi_1$ and 
$\Phi_2$, one concludes that $\varphi$ is a $\ce^1$ diffeomorphism. 
Finally, to show that $\varphi$ is a $\ce^r$ diffeomorphism in the case where 
the actions are non conjugate to actions by rotations, it suffices to apply 
Theorem \ref{con-jugando}. $\hfill\square$

\vspace{0.1cm}

\begin{small}
\begin{ejer} Let $\Gamma$ be a group of $\ce^1$ circle diffeomorphisms and $\varphi$ 
a bi-Lipschitz circle homeomorphism. Suppose that for every $f \in \Gamma$ and almost 
every $x \in \clo$ one has
$$\log ( f'(x) ) - \log ( f'(\varphi(x)) ) =
\log ( \varphi'(x) ) - \log ( \varphi'(f(x)) ).$$
Show that $\varphi$ {\bf \em{centralizes}} $\Gamma$ 
(that is, it commutes with every element of $\Gamma$).
\index{centralizer}
\end{ejer}

\begin{ejer} Given a circle diffeomorphism $f$, let $F$ be the diffeomorphism 
of $\mathbb{R}^2 \setminus \{ O \}$ defined in polar coordinates by 
$$F(r,\theta) = \left( \frac{r}{\sqrt{f'(\theta)}}, f(\theta) \right).$$

\vspace{0.08cm}

\noindent (i) Show that, if $f$ is of class $\ce^2$ and its rotation 
number is irrational, then $f$ is $C^1$ conjugate to $R_{\rho(f)}$ 
if and only if there exist positive $r_1,r_2$ such that 
$F^n \big( B( O, r_1) \big) \subset B \big( O, r_2 \big)$ for all $n \geq 0$.

\vspace{0.08cm}

\noindent (ii) Show that $F$ preserves the Lebesgue measure, and that for every $\theta$ 
the parallelogram generated by the vectors $F(1,\theta)$ and $\frac{d}{d\theta} F(1,\theta)$ 
has area 1. Conclude that there exists a real-valued function $\psi$ defined on the circle 
such that, for every $\theta \in \clo$, 
$$F(1,\theta) + \psi \esp \frac{d^2}{d \theta^2} F(1,\theta) = 0.$$

\noindent (iii) Endow the circle with the canonical projective structure ({\em c.f.,} 
Exercise \ref{estructura-proyectiva}), and identify each angle $\theta \in [0,2\pi]$ 
with the corresponding point in the projective space (so that $\theta$ and $\theta + \pi$ 
identify with the same point for $\theta \in [0,\pi]$). 
\index{projective!space} 
Show that, in the corresponding projective 
coordinates, the value of the function $\psi$ above is given by 
$$\psi = \frac{S(f)}{2} - 1.$$

\noindent{\underbar{Remark.}} This corresponds to the starting point of a 
proof using Sturm-Liouville's theory of a beautiful theorem due to Ghys: For every  
diffeomorphism of the circle (endowed with the canonical projective structure) 
there exist at least four points at which the Schwarzian derivative vanishes 
(see \cite{taba}). 
\index{Ghys!four vertex theorem for the Schwarzian derivative}  
\index{derivative!Schwarzian}
\end{ejer}
\end{small}

\vspace{0.1cm}

Remark that a conjugacy of an action to itself corresponds to a homeomorphism which 
centralizes the action. Moreover, if $\Phi_1$ and $\Phi_2$ are actions by $\ce^r$ 
circle diffeomorphisms which are conjugate by some $\ce^r$ diffeomorphism $\varphi_0$, 
and if $\varphi$ is another bi-Lipschitz homeomorphism which conjugates them, then the 
(bi-Lipschitz) homeomorphism $\varphi_0^{-1} \varphi$ centralizes the action $\Phi_1$. 
This is the reason why it is so important to deal with the study of centralizers before 
\index{centralizer} 
studying the general problem of conjugacies (see however Exercise \ref{no-anal-jap}). In this 
direction, one may show that Theorem \ref{japon1} is far from being true for non minimal actions.

\vspace{0.2cm}

\begin{thm} {\em Let $\Gamma$ be a finitely generated group of $\ce^{1+\mathrm{Lip}}$ 
circle diffeomorphisms whose action is not minimal. If the restrictions of the 
stabilizers of intervals are either trivial or infinite cyclic, then 
there exist bi-Lipschitz homeomorphisms which are not $\ce^1$ and which 
commute with all of the elements in $\Gamma$. Moreover, these homeomorphisms 
may be chosen non differentiable on every non-empty open interval in $\clo$.}
\label{japon2}
\end{thm}

\vspace{0.2cm}

This theorem is based on a very simple construction related to the techniques of the 
next chapter. This is the reason why we postpone its proof to \S \ref{lem-kopp}. For 
the moment, let us remark that the hypothesis on the stabilizers is not too strong. For 
instance, it is always satisfied in the real-analytic case (this result is due to Hector; 
\index{Hector}
a complete proof may be found in the Appendix of \cite{Na-kova}). Obviously, it is also 
satisfied by many other interesting non real-analytic actions, as for instance those 
of Thompson's group G. Nevertheless, without such a hypothesis, bi-Lipschitz 
conjugacies are forced to be smooth in many cases.

\vspace{0.05cm}

\begin{small} \begin{ejer}
Give examples of finitely generated groups of $\ce^{\infty}$ circle diffeomorphisms which are 
conjugate by some bi-Lipschitz homeomorphism though there is no smooth conjugacy between them.\\

\noindent{\underbar{Remark.}} These examples may be constructed either having finite 
orbits or an exceptional minimal set. However, we ignore whether there exists groups 
of {\em real-analytic diffeomorphisms} having the desired property.
\label{no-anal-jap}
\end{ejer} \end{small}


\newpage
\thispagestyle{empty}
${}$
\newpage


\chapter{STRUCTURE AND RIGIDITY VIA DYNAMICAL METHODS}
\label{algebra}

\section{Abelian Groups of Diffeomorphisms}
\label{seccion-kopell}

\subsection{Kopell's lemma}
\label{lem-kopp}\index{Kopell!lemma|(}

\hspace{0.45cm} For a group of homeomorphisms of an interval, the circle, or the real 
line, we will say that $]a,b[$ is an {\bf {\em irreducible component}} 
\index{irreducible component}
for the action if  
it is invariant and does not contain strictly any invariant interval. Let us denote 
by $\mathrm{Diff}_{+}^{1+\mathrm{bv}}([0,1[)$ the group of $\mathrm{C}^{1+\mathrm{bv}}$ 
diffeomorphisms of $[0,1[$, {\em i.e.}, the group of $\ce^1$ diffeomorphisms $f$ of 
$[0,1[$ such that the total variation of the logarithm of the derivative is 
finite on each compact interval $[a,b] \subset [0,1[$. Recall that 
\index{variation of a function!total variation}
this variation is denoted by $V(f;[a,b])$, that is,
$$V \big( f;[a,b] \big) \esp \esp = 
\sup\limits_{a = a_0 \! < \! a_1 \! < \! \cdots < \! a_n = b}
\sum_{i=1}^{n} \big| \log(f')(a_i) - \log(f')(a_{i-1}) \big|.$$
The group $\mathrm{Diff}_{+}^{1+\mathrm{bv}}(]0,1])$ is defined in a similar 
way. Notice that every element $f$ in $\mathrm{Diff}_+^2([0,1[)$ belongs to 
$\mathrm{Diff}_{+}^{1+\mathrm{bv}}([0,1[)$. Indeed, for $0 \leq a < b < 1$ one has
$$V \big( f;[a,b] \big) \esp = \esp \int_a^b \big| (\log(f'))'(s) \big| ds 
\esp = \esp \int_a^b \left| \frac{f''(s)}{f'(s)} \right| ds.$$
The important result below is stated as a theorem, although it is widely known (and we will 
refer to it) as Kopell's lemma, since it corresponds to the first lemma in the thesis of 
Kopell \cite{Kop}. Notwithstanding, we provide a much simpler proof, 

\vspace{0.2cm}

\begin{thm} {\em Let $f$ and $g$ be commuting diffeomorphisms of the interval $[0,1[$ 
or $]0,1]$. \esp Suppose that $f$ is of class $\mathrm{C}^{1+\mathrm{bv}}$ and $g$ 
of class $\mathrm{C}^1$. If $f$ has no fixed point in $]0,1[$ and $g$ has at least 
one fixed point in $]0,1[$, then $g$ is the identity.}
\end{thm}

\noindent{\bf Proof.} We will give the proof only for the case of the interval $[0,1[$, 
since the case of $]0,1]$ is analogous. Replacing $f$ by $f^{-1}$ if necessary, we 
may suppose that $f(x) < x$ for every $x \! \in ]0,1[$. Let $b \! \in ]0,1[$ be a 
fixed point of $g$. For each $n \in \mathbb{Z}$ denote $b_n = f^n(b)$, and denote 
$a = b_1 = f(b)$. Since $g$ commutes with $f$, it must fix all the intervals 
$[b_{n+1},b_n]$, and hence 
$$g'(0) \esp = \esp 
\lim_{n \rightarrow \infty} \frac{g(b_n) - g(0)}{b_n - 0} \esp = \esp 1.$$ 
Let $\delta = V(f;[0,b])$. If $u$ and $v$ belong to $[a,b]$, then
$$\left| \log \Big( \frac{(f^n)'(v)}{(f^n)'(u)} \Big) \right| \leq
\sum\limits_{i=1}^{n} \Big| \log (f')(f^{i-1}(v)) -
\log (f')(f^{i-1}(u)) \Big| \leq \delta.$$
Let $u \!=\! x \!\in\! [a,b]$ and $v \!=\! f^{-n} g f^n (x) \!=\! g(x) \!\in\! [a,b]$. 
Using the equality 
$$g'(x) = \frac{(f^{n})'(x)}{(f^{n})'(f^{-n} g f^n(x))}
\hspace{0.1cm} g'(f^n(x)) = \frac{(f^{n})'(x)}{(f^{n})'(g(x))}
\hspace{0.1cm} g'(f^n(x)),$$
and passing to the limit as $n$ goes to infinity, it follows that 
$\sup_{x \in [a,b]} g'(x) \leq e^{\delta}$. Now remark that this remains true 
if we replace $g$ by $g^j$ for any $j \in \mathbb{N}$ (this is due to the fact that 
the constant $\delta$ depends only on $f$). Therefore, 
$$\sup_{x \in [a,b]} (g^j)'(x) \leq e^{\delta}.$$
Since $g$ fixes $a$ and $b$, this is not possible unless the restriction of $g$ 
to $[a,b]$ is the identity. Finally, since $f$ and $g$ commute, this implies  
that $g$ is the identity on the whole interval $[0,1[$. $\hfill\square$

\vspace{0.3cm}

The preceding theorem allows us to conclude that, for every 
$f \!\in \!\mathrm{Diff}_{+}^{1+\mathrm{bv}}([0,\!1[)$ without fixed points in $]0,1[$, 
its centralizer in $\mathrm{Diff}_{+}^1([0,1[)$ acts freely on $]0,1[$. H\"older's theorem 
then implies that this 
\index{H\"older!theorem}\index{action!free}centralizer 
\index{centralizer}
is semiconjugate to a group of translations. In fact, 
if this group of translations is dense, then the semiconjugacy 
is actually a conjugacy, as is stated in the next proposition.

\vspace{0.1cm}

\begin{prop} {\em Let $\Gamma$ be a subgroup of $\mathrm{Diff}_+^1([0,1[)$ which 
is semiconjugate to a dense subgroup of the group of translations. If $\Gamma$ 
contains an element of class $\mathrm{C}^{1+\mathrm{bv}}$ without fixed 
points in $]0,1[$, then the semiconjugacy is a topological conjugacy.}
\label{no-semiconjugado}
\end{prop}

\noindent{\bf Proof.} Suppose that $\Gamma$ is a subgroup of $\mathrm{Diff}^1_+([0,1[)$ which 
is semiconjugate to a dense subgroup of the group of translations without being conjugate to 
it. Let $f \!\in\! \Gamma$ be the element given by the hypothesis. Without loss of generality, 
we may suppose that $f$ is sent to the translation $T_{-1}\!: x \!\mapsto\! x-1$ by 
the homomorphism induced by the semiconjugacy. In particular, one has $f(x) < x$ for 
every $x \!\in ]0,1[$. Choose an interval $[a,b]$ non reduced to a point which is sent 
to a single point by the semiconjugacy, and which is maximal for this property. By the 
choice of $[a,b]$, there exists an increasing sequence $(n_i)$ of positive integers 
such that for every $i \!\in\! \mathbb{N}$ there exists $\bar{f}_i \in \Gamma$ verifying, 
for every $n \in \mathbb{N}$,
$$\bar{f}_i^{n_i}(f^n(a)) \geq f^{n+1}(a),
\qquad \bar{f}_i^{n_i+1}(f^n(a)) < f^{n+1}(a),$$
$$\bar{f}_i^{n_i}(f^n(b)) \geq f^{n+1}(b),
\qquad \bar{f}_i^{n_i+1}(f^n(b)) < f^{n+1}(b).$$
Let $a_n \!=\! f^n(a)$ and $b_n \!=\! f^n(b)$. Passing to 
the limit as $n$ goes to infinity in the inequalities
$$\frac{\bar{f}_i^{n_i+1}(a_n)}{a_n} < \frac{f(a_n)}{a_n}
\leq \frac{\bar{f}_i^{n_i}(a_n)}{a_n}$$
we obtain
\begin{equation}
(\bar{f}_i'(0))^{n_i +1} \esp \leq \esp 
f'(0) \esp \leq \esp (\bar{f}_i'(0))^{n_i}.
\label{contderi}
\end{equation}
For $n \geq 0$ the intervals $f^n(]\bar{f}_i(a),b[)$ are two-by-two disjoint. If we 
denote $\delta = V(f;[0,b])$, then for every $u,v$ in $]\bar{f}_i(a),b[$ one has
\begin{equation}
\left| \log \Big( \frac{(f^n)'(v)}{(f^n)'(u)} \Big) \right| \esp \leq \esp 
\sum\limits_{i=1}^{n} \Big| \log (f')(f^{i-1}(v)) - \log (f')(f^{i-1}(u)) \Big| 
\esp \leq \esp \delta.
\label{clasiquita}
\end{equation}
Passing to the limit as $n$ goes to infinity in the inequality
$$|\bar{f}_i([a,b])| \esp = \esp |f^{-n}\bar{f}_if^n ([a,b])|
\esp \geq \esp \frac{\inf_{u \in [\bar{f}_i(a),b]}(f^n)'(u)}
{\sup_{v \in [\bar{f}_i(a),b]}(f^n)'(v)} \cdot
\inf_{x \in [f^n(a),f^n(b)]} \bar{f}_i'(x) \cdot |[a,b]|,$$
and using the estimates (\ref{contderi}) and (\ref{clasiquita}),  
we obtain, for some positive constant $C$ and every $i \in \mathbb{N}$,
$$|\bar{f}_i([a,b])| \geq e^{-\delta} (f'(0))^{1/n_i} \cdot |[a,b]| \geq C.$$
However, this inequality is absurd, since $|\bar{f}_i([a,b])|$ 
obviously converges to zero as $i$ goes to infinity. $\hfill\square$

\vspace{0.1cm}

\begin{small}\begin{ejer} Using Denjoy's theorem, prove that Proposition \ref{no-semiconjugado} 
holds for every finitely generated subgroup of $\mathrm{Diff}_+^{1+\mathrm{bv}}(]0,1[)$ whose 
center contains an element without fixed points (compare Example \ref{ejemplohirsh}).
\end{ejer}\end{small}

From the preceding proposition one immediately deduces the following.

\vspace{0.05cm}

\begin{cor} {\em If $f$ is an element in $\mathrm{Diff}_{+}^{1+\mathrm{bv}}([0,1[)$ without fixed 
points in $]0,1[$, then its centralizer in $\mathrm{Diff}_{+}^1([0,1[)$ is topologically 
conjugate to a group of translations of the line.}
\label{flujocelog}
\end{cor}

\begin{small} \begin{ejer}
Prove the ``strong version'' of Kopell's lemma (in class  $\ce^{1+\mathrm{Lip}}$) 
given by the next proposition (see either \cite{level} or \cite{firmo2} in 
case of problems with this).

\begin{prop} {\em Let $f\!:\! [0,1[ \rightarrow [0,1[$ be a 
$\ce^{1+\mathrm{Lip}}$ diffeomorphism such that $f(x) \!<\! x$ 
for all $x \!\! \in ]0,1[$. Fix a point $a \! \in ]0,1[$, and for each 
$n \!\in\! \mathbb{N}$ let $g_n\!\!: [f(a),a] \!\rightarrow\! [f(a),a]$ 
be a diffeomorphism tangent to the identity at 
the endpoints. If $g\!: ]0,1] \!\rightarrow ]0,1]$ is such that its 
restriction to $[f^{n+1}(a),f^n(a)]$ coincides with $f^n g_n f^{-n}$ 
for every $n \!\in\! \mathbb{N}$, then $g$ extends to a $\mathrm{C}^1$ 
diffeomorphism of $[0,1[$ if and only if $(g_n)$ converges to the 
identity in the $\mathrm{C}^1$ topology.}
\end{prop}
\label{kopell-fuerte-lip}
\end{ejer}

\begin{ejer} Prove the following ``real-analytic version'' of Kopell's lemma 
(due to Hector): If 
\index{Hector}
$f$ and $g$ are real-analytic diffeomorphisms that may be written in the form 
$f(x) = x + a_ix^i + \cdots$ and $g(x) = x + b_jx^j + \cdots$ for $|x| \leq \varepsilon$, 
with $j > i$ and $f(x) < x$ for every small positive $x$, then the sequence of maps $f^{-n}gf^n$ 
converges uniformly to the identity on an interval $[0,\varepsilon']$ (see \cite{Na-kova} for 
an application of this claim).
\end{ejer} \end{small}

To close this section, we give a proof of Theorem \ref{japon2}. 
Let us begin by considering a diffeomorphism $f$ of class $\ce^{1+\mathrm{bv}}$ 
of an interval $I \!=\! [a,b]$ such that $f(x) < x$ for every 
$x \! \in ]a,b[$. Fix an arbitrary point $c \!\in ]a,b[$, 
and consider any bi-Lipschitz 
\index{Lipschitz!conjugacy}
homeomorphism $h$ from $[f(c),c]$ into itself. Extending $h$ to $]a,b[$ 
so that it commutes with $f$, and then putting $h(a)\!=\!a$ and $h(b)\!=\!b$, we obtain 
a well-defined homeomorphism of $[a,b]$ (which we still denote by $h$). We claim that 
if $C$ is a bi-Lipschitz constant for $h$ in $[f(c),c]$, then $C e^V$ is a bi-Lipschitz 
constant for $h$ on $[a,b]$, where $V\!=\!V(f;[a,b])$. Indeed, let us consider 
for instance a point $x \!\in\! f^n([f(c),c])$ for some $n \geq 0$, and such that $h$ is 
differentiable on $f^{-n}(x) \in [f(c),c]$, with derivative less than or equal to $C$ 
(notice that this is the case for almost every point $x \!\in\! [f^{n+1}(c),f^n(c)]$). 
From the relation $h \!=\! f^n h f^{-n}$ we obtain 
\begin{equation}
h'(x) = h'(f^{-n}(x)) \cdot \frac{(f^n)'(hf^{-n}(x))}{(f^n)'(f^{-n})(x)} \leq
C \cdot \frac{(f^n)'(hf^{-n}(x))}{(f^n)'(f^{-n})(x)}.
\label{unilita}
\end{equation}
Letting $y = f^{-n}(x) \in [f(c),c]$ and $z = h(y) \in [f(c),c]$, 
and arguing as in the proof of Kopell's lemma, this gives 
$$\left| \log \Big( \frac{(f^n)'(z)}{(f^n)'(y)} \Big) \right| = 
\left| \log \Big( \frac{\prod_{i=0}^{n-1} f'(f^i(z))}{\prod_{i=0}^{n-1} f'(f^i(y))}
\Big) \right| \leq \sum_{i=0}^{n-1} \Big| \log(f'(f^i(z))) - \log(f'(f^i(y))) \Big| 
\esp \leq \esp V.$$
Introducing this inequality in (\ref{unilita}) we deduce that $h'(x) \!\leq\! C e^{V}$. 
Since $x$ was a generic point, this shows that the Lipschitz constant of $h$ is bounded 
from above by $Ce^V$. Obviously, a similar argument shows that this bound also holds  
for the Lipschitz constant of $h^{-1}$.

\vspace{0.15cm}

For the proof of Theorem \ref{japon2} we will use a similar construction. To simplify, 
we will only prove the first claim of the theorem, leaving to the reader the task of 
proving the second claim concerning the existence of bi-Lipschitz centralizers 
which are non differentiable on any open set.

Let us begin by recalling that if $\Gamma$ is a group of $\ce^{1+\mathrm{Lip}}$ circle 
diffeomorphisms preserving an exceptional minimal set, then the stabilizer of every 
connected component $]a,b[$ of the complement 
of this set is nontrivial ({\em c.f.,} Exercise 
\ref{ejer-hector-int}). By the hypothesis of the theorem, the restriction to $]a,b[$ 
of this stabilizer is either trivial or infinite cyclic. In the first case, we define 
$h$ as being any bi-Lipschitz and non differentiable homeomorphism of $[a,b]$. In the 
second case, let $f \!\in\! \Gamma$ be so that its restriction to $]a,b[$ generates 
the restriction of the corresponding stabilizer. Let $[\bar{a},\bar{b}] \subset [a,b]$ 
be such that $f(x) \neq x$ for every $x \in ]\bar{a},\bar{b}[$, and such that 
$f(\bar{a}) = \bar{a}$ and $f(\bar{b}) = \bar{b}$. Changing $f$ by $f^{-1}$ 
if necessary, we may assume that $f^n(x)$ converges to $\bar{a}$ as $n$ goes to 
infinity for every $x \in [\bar{a},\bar{b}[$. Arguing as in the case of the 
interval above, fix a point $\bar{c} \in ]\bar{a},\bar{b}[$, and consider any 
bi-Lipschitz and non differentiable homeomorphism $h$ of $[f(\bar{c}),\bar{c}]$. 
This homeomorphism extends in a unique way to $[a,b]$ so that it commutes with the 
restriction of $f$ to $[\bar{a},\bar{b}]$ and coincides with the identity on 
$I \setminus [\bar{a},\bar{b}]$.

By the hypothesis on the stabilizers, it is not difficult to check that there exists a 
unique extension of $h$ to a circle homeomorphism (which we still denote by $h$) which 
commutes with (every element of) $\Gamma$ and which coincides with the identity on the 
complement of $\cup_{g \in \Gamma} \hspace{0.1cm} g(]a,b[)$. We claim that this 
extension is bi-Lipschitz. More precisely, fixing a finite system of generators 
$\mathcal{G} \!=\! \{g_1,\ldots,g_k\}$ for $\Gamma$, letting $V$ be the maximum 
among the total variation of the logarithm of the derivatives of the $g_i$'s, and 
denoting by $C$ the bi-Lipschitz constant of $h$ on $[a,b]$, 
we claim that the globally defined homeomorphism $h$ has a 
bi-Lipschitz constant less than or equal to $Ce^{kV}$. The proof of this 
claim is similar to that of the case of the interval. Let us fix for instance 
$x \!\in\! \cup_{g \in \Gamma} \hspace{0.1cm} (g(I) \setminus I)$, and let us 
try to estimate the value of $h'(x)$. For this, let us consider the smallest 
$n \!\in\! \mathbb{N}$ for which there exists an element 
$g \!=\! g_{i_n} \ldots g_{i_1} \!\in\! \Gamma$, with $g_{i_j}$ in $\mathcal{G}$, 
such that $g(x) \in I$. The minimality of $n$ implies that the intervals
$I, g_{i_n}^{-1}(I),g_{i_{n-1}}^{-1}g_{i_n}^{-1}(I), \ldots,g_{i_1}^{-1} 
\cdots g_{i_n}^{-1}(I)$ have disjoint interior. Using the relation  \esp 
$h = g^{-1} h g$ \esp we obtain, for a generic $x \!\in\! g^{-1}(I)$,
\begin{equation}
h'(x) = h'(g(x)) \cdot \frac{g'(x)}{g'(h(x))} \leq C \cdot
\frac{g'(x)}{g'(y)},
\label{jime}
\end{equation}
where $y = h(x) \in g^{-1}(I)$. Now using the fact that the total variation of the 
logarithm of the derivative of each $g_i$ is bounded from above by $V$, we obtain
\begin{small}
$$\left| \log \Big( \frac{g'(x)}{g'(y)} \Big) \! \right|
\!\leq\! \sum_{j = 0}^{n-1} \!\big| \! \log(g_{i_{j+1}}'(g_{i_j} \cdots g_{i_1})(x))
- \log(g_{i_{j+1}}'(g_{i_j} \cdots g_{i_1})(y)) \big|  
\!\leq\! \sum_{i=1}^k \! V \big( \log(g_i'); \clo \big) \!\leq\!  kV.$$
\end{small}From (\ref{jime}) we conclude  
that $h'(x) \leq C e^{kV}$, as we wanted to check.

Finally, let us consider the case where $\Gamma$ admits finite orbits. If $\Gamma$ is finite, 
then we consider any bi-Lipschitz and non differentiable diffeomorphism which commutes with 
its generator. If $\Gamma$ is infinite, H\"older's theorem 
\index{H\"older!theorem}
implies that the action cannot be free. Let $I$ be 
a connected component of the complementary set of the union of the finite orbits, and let 
$f \!\in\! \Gamma$ be so that its restriction to $I$ generates the restriction of the 
corresponding stabilizer. If we proceed as before with $I$ and $f$, we may construct 
a bi-Lipschitz and non differentiable homeomorphism centralizing $\Gamma$. We leave 
the details to the reader.
\index{Kopell!lemma|)}


\subsection{Classifying Abelian group actions in class $\ce^2$}

\hspace{0.45cm} We may now give a precise description of the Abelian groups of 
$\mathrm{C}^{1+\mathrm{bv}}$ diffeomorphisms of one-dimensional manifolds. The case 
of the interval is quite simple. Indeed, by Corollary \ref{flujocelog}, the restriction 
of such a group to each irreducible component 
\index{irreducible component}
is conjugate to a group of translations. 
The case of the circle is slightly more complicated. We begin with a lemma which is 
interesting by itself (compare \S \ref{sec-margulis} and Exercise \ref{a-referenciar}).

\vspace{0.15cm}

\begin{lem} {\em If $\Gamma$ is an amenable subgroup of $\mathrm{Homeo}_{+}(\clo)$, 
then either $\Gamma$ is semiconjugate to a group of rotations or it contains a 
finite index subgroup having fixed points.}
\label{moyen-medida}
\end{lem}

\vspace{-0.35cm}

\noindent{\bf Proof.} Since $\Gamma$ is amenable, \index{group!amenable} it must preserve 
a probability measure $\mu$ on $\clo$. If the orbits by $\Gamma$ are dense, then $\mu$ 
has total support, and by reparameterizing the circle one easily checks that $\Gamma$ is 
topologically conjugate to a group of rotations. 
\index{rotation!group} 
If there exists a minimal invariant 
Cantor set, then the support of $\mu$ coincides with this set, and this allows 
semiconjugating $\Gamma$ to a group of rotations. Finally, if there is a finite 
orbit, then the elements in $\Gamma$ preserve the cyclic order of the points in 
this orbit. The stabilizer of these points is then a finite index subgroup of 
$\Gamma$ having fixed points. $\hfill\square$

\vspace{0.45cm}

From the preceding lemma (and its proof) it follows that, if $\Gamma$ is an Abelian subgroup 
of $\mathrm{Diff}_+^{1+\mathrm{bv}}(\clo)$, then it is either semiconjugate to a group 
of rotations, or a finite central extension of an at most countable  product of 
Abelian groups acting on disjoint intervals. By Corollary \ref{flujocelog}, the factors 
of this product are conjugate to groups of translations on each irreducible component. 
Finally, recall that for finitely generated subgroups of $\mathrm{Diff}_+^{1+\mathrm{bv}}(\clo)$, 
every semiconjugacy to a group of rotations is necessarily a conjugacy; however, this 
is no longer true for non finitely generated groups ({\em c.f.,} Example 
\ref{ejemplohirsh}).

\begin{small}
\begin{ejer} Prove that every virtually Abelian group of $\ce^{1+\mathrm{bv}}$ 
diffeomorphisms of the interval or the circle is actually Abelian.
\label{virt-abel}
\end{ejer}

\begin{ejer}
Prove that every Abelian subgroup of $\mathrm{Diff}^{1+\mathrm{bv}}_{+}(\mathbb{R})$
preserves a Radon measure on the line (compare Proposition \ref{abelianoplante}).
\index{measure!Radon}
\label{plantecdos}
\end{ejer} \end{small}


\subsection{Szekeres' theorem}
\label{szek}
\index{Szekeres' theorem|(}

\hspace{0.45cm} In class $\ce^2$, the homomorphism given by Corollary \ref{flujocelog} is 
necessarily surjective. This follows from the result below, due to Szekeres \cite{sz-th} 
(see also \cite{Yo1}). To state it properly, 
given a non-empty non-degenerate interval $[a,b[$ and $2 \leq r \leq \infty$, let us 
denote by $\mathrm{Diff}^{r,\Delta}_{+}([a,b[)$ the subset of $\mathrm{Diff}^r_{+}([a,b[)$ 
formed by the elements $f$ such that $f(x) \neq x$ for all $x \! \in ]a,b[$. For simplicity, 
we let $\infty\!-\!1 = \infty$.
\index{centralizer}

\vspace{0.1cm}

\begin{thm} {\em For every $f \in \mathrm{Diff}^{r,\Delta}_{+}([0,1[)$ there exists 
a unique vector field $X_f$ on $[0,1[$ without singularities in $]0,1[$ satisfying:}

\vspace{0.07cm}

\noindent{(i) {\em $X_f$ is of class $\mathrm{C}^{r-1}$ on $]0,1[$ and 
of class $\mathrm{C}^1$ on $[0,1[$;}}

\vspace{0.07cm}

\noindent{(ii) {\em if $f^{\mathbb{R}} \!=\! \{f^t \! : t \in \mathbb{R} \}$ 
is the flow associated to this vector field, then $f^1 \!=\! f$;}}

\vspace{0.07cm}

\noindent{(iii) {\em the centralizer of $f$ in $\mathrm{Diff}_{+}^{1}([0,1[)$ 
coincides with $f^{\mathbb{R}}$.}}
\label{sz}
\end{thm}

\vspace{0.1cm}

This theorem is particularly interesting if the germ of $f$ at the origin is non hyperbolic, 
since otherwise we may use Sternberg's linearization theorem. Actually, in the hyperbolic 
case, the claim before Example \ref{ejst} allows extending Theorem 
\ref{sz} to class $\mathrm{C}^{1+\tau}$ for every $\tau \!>\! 0$.

\begin{small} \begin{ejer} Given $\lambda < 0$, consider a $\ce^1$ vector 
field $X \!=\! \varrho \esp \! \frac{\partial}{\partial x}$ on $[0,1]$ 
such that for every $x$ small enough one has \esp 
$\varrho(x) = \lambda x \big( 1 - 1 / \log(x) \big).$ \esp 
Prove that $X$ is not $\ce^1$ linearizable, \index{linearization} 
that is, there is no $\ce^1$ diffeomorphism conjugating $X$ 
to its linear part $\lambda x \frac{\partial}{\partial x}$ 
(compare Exercise \ref{ejst}).

\noindent{\underbar{Hint.}} Show that if $f$ 
is the time 1 of the flow associated to $X$, then
$$f(x) = e^{\lambda} \esp x \esp \left( \frac{1 - \log(x)}{1 - \log(f(x))} \right)$$
for all $x$ near the origin. Using this, show more generally that $f$ 
is not conjugate to its linear part by any bi-Lipschitz homeomorphism.
\index{Lipschitz!conjugacy}
\end{ejer} \end{small}

We will give the proof of Theorem \ref{sz} only for the case $r\!=\!2$ (the extension to 
bigger $r$ is straightforward). First notice that 
we may assume that $f(x) > x$ for all $x \!\in ]0,1[$: the vector field associated to 
a diffeomorphism which is topologically contracting at the origin may be 
obtained by changing the sign of the vector field associated to its inverse. 

\vspace{0.1cm}

\begin{lem} {\em Let $f\!:[0,1[ \rightarrow \! [0,1[$ be a $\mathrm{C}^2$ diffeomorphism such 
that $f(x) > x$ for every $x \! \in ]0,1[$. If $X (x) = \varrho(x) \frac{\partial}{\partial x}$ 
is a $\mathrm{C}^1$ vector field on $[0,1[$ and $a \! \in ]0,1[$, then $X$ is associated to 
$f$ if and only if the following conditions are satisfied:}

\vspace{0.07cm}

\noindent{(i) \em the function $\varrho$ is strictly positive on $]0,1[$;}

\vspace{0.07cm}

\noindent{(ii) \em for every $x \in [0,1[$ one has $\varrho (f(x)) = f'(x) \varrho(x)$;}

\vspace{0.07cm}

\noindent{(iii) \em one has $\int_{a}^{f(a)} \frac{ds}{\varrho(s)} = 1$.}
\label{caract-camp}
\end{lem}

\noindent{\bf Proof.} The first two conditions are clearly necessary. 
Concerning the third one, notice that if $X$ is associated to $f$, 
then for every $t \geq 0$ we have 
$$\frac{d}{dt} \int_{a}^{f^t(a)} \frac{ds}{\varrho (s)}
= \frac{d f^t(a)}{dt} \cdot \frac{1}{\varrho(f^t(a))} = 1.$$ 
Thus, 
$$\int_{a}^{f^t(a)} \frac{ds}{\varrho (s)} = t,$$
which implies (iii) by letting 
$t = 1$. Conversely, if condition (ii) is satisfied, then 
one easily checks that for every $x \!\in ]0,1[$ the derivative of the function 
$x \mapsto \int_{x}^{f(x)} \frac{ds}{\varrho(s)}$ is identically equal to zero. 
Hence, for every $x \!\in ]0,1[$,
$$\int_x^{f(x)} \frac{ds}{\varrho(s)} = \int_a^{f(a)} \frac{ds}{\varrho(s)} = 1.$$
Therefore, if we denote by $\hat{f}$ the diffeomorphism obtained by integrating $X$ up 
to time $1$, then $\hat{f}(0) \!=\! f(0) \!=\! 0$, and for all $x \!\in ]0,1[$ one has 
$$1 = \int_{x}^{\hat{f}(x)} \frac{ds}{\varrho(s)} = \int_{x}^{f(x)}
\frac{ds}{\varrho(s)}.$$
From this one easily concludes that $\hat{f} \!=\! f$, 
that is, $X$ is associated to $f$. $\hfill\square$

\vspace{0.2cm}

To construct the vector field $X$, we begin by considering the ``discrete difference''
$\Delta = \Delta_f$ defined by $\Delta(x) = f(x) - x$. Although $\Delta(f(x))$
and $\Delta(x) f'(x)$ do not coincide, the ``error'' has a nice expression. 
Indeed, if we define 
\begin{equation}
\Theta(x) = \Theta_f(x) = \log(f'(x)) - \log \left[ \int_{0}^{1} f'(x+t\Delta (x))dt
\right],
\label{def-theta}
\end{equation}
then from the equality 
$$f^2(x) = f(x) + \Delta(x) \int_{0}^{1} f'(x+t\Delta (x)) dt$$
it follows that
\begin{equation}
\Delta(f(x)) \exp(\Theta(x)) = \Delta(x) f'(x).
\label{delta-eq}
\end{equation}
For later estimates we will further use the second order Taylor series expansion  
$$f^2(x) = f(x) + \Delta(x) f'(x) + \Delta(x)^2 \int_0^1 (1-t)f''(x+t\Delta(x))dt,$$
which allows checking the equality
\begin{equation}
\Theta(x) = \log \left( 1 - \frac{\Delta(x)^2}{\Delta(f(x))}
\int_0^1 (1-t) f'' \big( x+t\Delta(x) \big) dt \right).
\label{theta-eq}
\end{equation}

The function $\varrho$ corresponding to $X$ will be obtained by multiplying 
$\Delta$ by the sum of the successive errors under the iteration, modulo some 
normalization.  To be more concrete, let us formally define a function 
$\Sigma$ by letting $\Sigma (0) = 0$ and, for $x \! \in ]0,1[$,
$$\Sigma (x) = \sum_{n > 0} \Theta(f^{-n}(x)).$$
This function satisfies the formal equality
\begin{equation}
\Sigma(f(x)) = \Sigma(x) + \Theta(x).
\label{sigma-eq}
\end{equation}
Therefore, if we define the field 
$Y(x) \!=\! \Delta(x) \exp(\Sigma (x)) \frac{\partial}{\partial x}$,
from (\ref{delta-eq}) and (\ref{sigma-eq}) one concludes that
$$Y(f(x)) = \Delta(f(x))\exp\! \big( \Sigma (f(x))\big) \frac{\partial}{\partial x} 
= \frac{\Delta(x) f'(x)}{\exp (\Theta(x))} \exp \big( \Theta(x) + \Sigma (x) \big)
\frac{\partial}{\partial x} = f'(x) Y(x).$$
The vector field $X$ will be then of the form $X = c Y$ for 
some normalizing constant $c = c(f)$ for which 
$$\int_a^{f(a)} \frac{ds}{X(s)} = 1.$$
Remark that the condition imposed by this equality corresponds to
\begin{equation}
c = \int_a^{f(a)} \frac{dx}{\Delta(x) \exp(\Sigma (x))}.
\label{condon}
\end{equation}
Suppose that $f'(0) > 1$, and denote $\lambda = f'(0)$. For 
$a \sim 0$ one has $\Sigma (x) \sim 0$ for $x \leq f(a)$, and 
$$\Delta(x) = x \left( \frac{f(x)}{x} - 1 \right) \sim x (\lambda - 1).$$
Hence, according to (\ref{condon}), we must have 
$$c = \int_a^{f(a)} \frac{dx}{\Delta(x) \exp(\Sigma (x))} \sim \int_a^{f(a)}
\frac{dx}{x(\lambda - 1)} = \frac{\log(f(a)/a)}{\lambda - 1}
\sim \frac{\log(\lambda)}{\lambda - 1},$$
and therefore the appropriate choice is \esp
$c(f) = \log(\lambda)/(\lambda - 1)$. \esp
When $\lambda$ tends to $1$, this expression converges to $1$, and thus 
in the case $f'(0)=1$ the normalizing constant to be considered should be 
\esp \!$c(f) = 1$. \!\esp We will now verify that the definition we have 
just sketched is pertinent and provides a vector field associated to $f$.

\vspace{0.1cm}

\begin{prop} {\em Let $f$ be a diffeomorphism satisfying the hypothesis of Lemma} 
\ref{caract-camp}. {\em If we define $X \!=\! \varrho \frac{\partial}{\partial x}$ 
by $\varrho(x) \!=\! c(f) \Delta(x) \exp(\Sigma(x))$, then $X$ is a $\ce^1$ vector 
field associated to $f$.}
\end{prop}

\noindent{\bf Proof.} We need to show three claims: the vector field $X$ is 
well-defined (in the sense that the series defining the function $\Sigma$ converges), 
it is of class $\ce^1$ on $[0,1[$, and $f$ is the time 1 of the flow associated to it.

\vspace{0.2cm}

\noindent{\underbar{First step.}} The convergence of $\Sigma(x)$.

\vspace{0.1cm}

Since $f$ is of class $\mathrm{C}^1,$ for each $x \!\in ]0,1[$
there exists $y = y(x) \in [x,f(x)]$ such that 
$$\int_{0}^{1} f' \big( x + t\Delta(x) \big) dt = f'(y).$$
From this equality and (\ref{def-theta}) one deduces that 
$$\big| \Theta(x)| = |\log(f'(x)) - \log(f'(y)) \big| \leq C |y - x| \leq C \Delta(x),$$
where $C$ is the Lipschitz constant of the function $\log(f')$ on $[x,f(x)]$ 
\index{Lipschitz!derivative}
(which equals the supremum of the function $|f''| / |f' |$ on this interval). One 
then concludes that the series corresponding to $\Sigma (x)$ converges for every 
$x \!\in ]0,1[$. Moreover, $|\Sigma (x)| \leq C x$, which implies that the 
function $\Sigma$ extends continuously to $[0,1[$ by letting $\Sigma(0) = 0$.

\vspace{0.2cm}

\noindent{\underbar{Second step.}} The differentiability of $\varrho$.

\vspace{0.1cm}

To show that $\varrho$ is differentiable at the origin, it suffices to notice that 
$$\lim_{t \rightarrow 0} \frac{\varrho(t)}{t} = c \hspace{0.05cm}
\lim_{t \rightarrow 0} \exp \big( \Sigma(t) \big) \hspace{0.05cm} \lim_{t
\rightarrow 0} \frac{\Delta(t)}{t} = c \hspace{0.05cm} \Delta'(0) = \log(\lambda).$$
To check the differentiability at the interior we will use an indirect argument. 
Denoting by $L_X$ the Lie derivative along the vector field $X$, from the relation 
\esp $\varrho \big( f(x) \big) = \varrho(x) f'(x)$ \esp it follows that \esp 
$L_X(\Theta \circ f^{-1}) = (L_X \Theta) \circ f^{-1}$. \esp Now from 
(\ref{theta-eq}) it is easy to conclude that, for some constant $C>0$,
$$\varrho(x) \esp \Theta'(x) \leq C \esp \Delta(x).$$
Therefore, the series \esp $\sum_{n>0} (L_X \Theta) \circ f^{-n}$ \esp converges,
and its value equals $L_X (\Sigma) (x)$. This shows that $\Sigma$, and hence $X$, 
is of class $\ce^1$ on $]0,1[$. Moreover, for every $x \in ]0,1[$,
$$\varrho'(x) = \varrho(x) \esp \frac{\Delta'(x)}{\Delta(x)} +
\sum_{n>0} L_X (\varrho) \circ f^{-n} (x).$$
Letting $x$ go to the origin, one readily concludes from this relation 
that $\varrho'(x)$ converges to $c\esp\Delta'(0) = \log(\lambda)$, 
which shows that $X$ is of class $\ce^1$ on the whole interval $[0,1[$.

\vspace{0.2cm}

\noindent{\underbar{Third step.}} The time 1 of the flow.

\vspace{0.1cm}

We need to check that condition (iii) holds for every $x \!\! \in ]0,1[$. 
Now, as we have already seen in the proof of Lemma \ref{caract-camp}, the 
equality $\varrho(f(x)) = \varrho(x) f'(x)$ implies that the function 
$x \mapsto \int_{x}^{f(x)} \frac{ds}{\varrho(s)}$ is constant.

Suppose first that $f$ is tangent to the identity at the origin, that is, $\lambda=1$. 
In this case, one has $\Delta'(0)=0$. Hence, for every $\varepsilon > 0$ there 
exists $\delta > 0$ such that, if $a < \delta$ and $t \in [a,f(a)]$, then 
$$|\Delta(t) - \Delta(a)| < \varepsilon (t-a), \qquad 
\qquad 1-\varepsilon < \frac{1}{\exp(\Sigma(t))} < 1+\varepsilon.$$
The first of these inequalities implies that \esp 
$|\Delta(t) - (f(a)-a)| < \varepsilon (b-a),$ \esp 
whereas the second one gives
$$\left| \frac{\Delta(t)}{\varrho(t)} - 1 \right| < \varepsilon.$$
One then obtains, for $a < \delta$,
$$\int_a^{f(a)} \frac{du}{\varrho(u)} > (1-\varepsilon) \int_a^{f(a)}
\frac{du}{\Delta(u)} > \frac{1-\varepsilon}{1+\varepsilon},$$
$$\int_a^{f(a)} \frac{du}{\varrho(u)} < (1+\varepsilon) \int_a^{f(a)}
\frac{du}{\Delta(u)}
< \frac{1+\varepsilon}{1-\varepsilon}.$$
Nevertheless, since the function $x \mapsto \int_x^{f(x)} \frac{ds}{\varrho(s)}$
is constant, letting $a$ go to the origin as $\varepsilon$ goes to 
zero one deduces that this constant equals $1$.

Suppose now that $\lambda > 1$. In this case, 
$\Delta'(0) = c-1$, and hence for $t$ small one has
$$|\Delta(t) - t(c-1)| < \varepsilon t, \quad \quad
\left| \frac{\lambda \Delta(t)}{\varrho(f)(t)} - 1 \right| < \varepsilon.$$
From this it follows that
\begin{eqnarray*}
\int_{a}^{f(a)} \frac{ds}{\varrho(s)}
&<& (1+\varepsilon) \!\int_{a}^{f(a)} \! \frac{ds}{\lambda \Delta(s)}
\esp\esp = \esp\esp \frac{1+\varepsilon}{\lambda} \!\int_a^{f(a)} \! \frac{ds}{\Delta(s)}
\esp\esp < \esp\esp \frac{1+\varepsilon}{\lambda} \!\int_a^{f(a)} \!\!\! \frac{ds}{(c-1-\varepsilon) s}\\
&=& \frac{1+\varepsilon}{\lambda(c-1-\varepsilon)} \log \left( \frac{f(a)}{a} \right)
\esp\esp = \esp\esp \frac{1+\varepsilon}{\lambda(c-1-\varepsilon)} \log \left( 1+\frac{\Delta(a)}{a}
\right)\\
&<& \frac{1+\varepsilon}{\lambda} \cdot \frac{\log(c+\varepsilon)}{c-1-\varepsilon}
\esp\esp = \esp\esp (1+\varepsilon) \frac{\log(c+\varepsilon)}{\log(c)} \cdot
\frac{c-1}{c-1-\varepsilon}.
\end{eqnarray*}
Passing to the limit this yields 
$$\int_{a}^{f(a)} \frac{ds}{\varrho(s)} \leq 1.$$
A similar argument shows the opposite inequality, 
thus concluding the proof. $\hfill\square$

\vspace{0.25cm}

\begin{small} \begin{ejer}
Show that if the original diffeomorphism is of class $\mathrm{C}^{1+\mathrm{Lip}}$, then 
the preceding construction provides us with a Lipschitz vector field associated to it.
\end{ejer} \end{small}

\vspace{0.01cm}

We leave to the reader the task of checking that, if the diffeomorphism $f$ is of class 
$\mathrm{C}^r$ for some $r \geq 2$, then the associated vector field $X$ is of class 
$\mathrm{C}^{r-1}$ on $]0,1[$. However, $X$ may fail to be twice differentiable at the 
origin. Moreover, if $f$ is a diffeomorphism of the interval $[0,1]$, then the vector field 
$X$ defined on $[0,1[$ extends continuously to $[0,1]$ by letting $X(1) = 0$, but this 
extension is non differentiable in ``most of the cases'' (see \cite[Chapters IV and V]{Yo1} 
for more details on this; see also \cite{elena}). 
Anyway, the uniqueness and the smoothness of $X$ allows establishing 
an interesting result of regularity for the conjugacy between diffeomorphisms.

\vspace{0.08cm}

\begin{prop} {\em Let $r \geq 2$, and let $f_1,f_2$ be $\ce^r$ diffeomorphisms of an interval 
which is closed at least from one side. If $f_1,f_2$ have no fixed point in the interior 
of this interval, then the restriction to the interior of every (if any) $\ce^1$ 
diffeomorphism conjugating them is a $\ce^r$ diffeomorphism.}
\label{conj-caca}
\end{prop}
\index{Szekeres' theorem|)}


\subsection{Denjoy counterexamples}
\label{ejemplosceuno}

\hspace{0.45cm} The regularity $\ce^{1+\mathrm{Lip}}$ (or $\ce^{1+\mathrm{bv}}$) 
\index{Denjoy!counterexamples|(} 
is necessary for many of the dynamical results of the preceding chapter, as well as 
for some of the algebraic results of this one. Before passing to the construction 
of ``counterexamples'' to some of them, we need to recall the notion of 
{\bf \em{modulus of continuity}}. 
\index{modulus of continuity}

\begin{defn} Given a homeomorphism 
$\omega\!: [0,1] \!\rightarrow \! [0,\omega(1)]$, we say that a function 
$\psi\!: [0,1] \!\rightarrow\! \mathbb{R}$ is $\omega$-continuous if there 
exists $C \in \mathbb{R}$ such that, for every $x \neq y$ in $[0,1]$,
$$\left| \frac{\psi(x) - \psi(y)}{\omega(|x-y|)} \right| \esp \leq \esp C.$$
\end{defn}

Let us denote the supremum of the left-side expression by $\|\psi\|_{\omega}$, and 
let us call it the $\ce^{\omega}$-norm of $\psi$. The main interest on the notion 
of $\omega$-continuity comes from the obvious fact that every sequence of functions 
$\psi_n$ defined on $[0,1]$ and such that 
$$\sup_{n \in \mathbb{N}} \| \psi_n \|_{\omega} < \infty$$
is equicontinuous. 
\index{equicontinuous!sequence of functions}
In a certain way, having a uniformly bounded modulus of continuity for 
a sequence of functions is a kind of quantitative (and sometimes 
easy to handle) criterion for establishing equicontinuity. 

\begin{small} \begin{ejem}
For $\omega(s) = s^{\tau}$, where $0\!<\!\tau\!<\! 1$, the notions 
of $\omega$-continuity and $\tau$-H\"older continuity coincide.  
\index{H\"older!continuity}
\label{wholder}
\end{ejem}

\begin{ejem}
For $\omega(s) = s$, the notion of $\omega$-continuity corresponds to that of 
Lipschitz continuity.
\index{Lipschitz!continuity}
\label{wlipschitz}
\end{ejem}

\begin{ejem}
Given $\varepsilon \!\geq\! 0$, let $\omega \!=\! \omega_{\varepsilon}$ be such 
that $\omega_{\varepsilon}(s) \!=\! s [\log(1/s)]^{1+\varepsilon}$ for $s$ small. 
If a map is $\omega_{\varepsilon}$-continuous, then it is $\tau$-continuous for every 
$0 < \tau < 1$. Indeed, one easily checks that 
$$s\log \left( \frac{1}{s}\right)^{1+\varepsilon}
\leq \esp C_{\varepsilon,\tau} s^{\tau}, \quad \mbox{ where } \quad 
C_{\varepsilon,\tau}
= \frac{1}{e^{1+\varepsilon}} \left( \frac{1+\varepsilon}{1-\tau}
\right)^{1+\varepsilon}.$$
Notice that the map $s \mapsto s \log(1/s)^{1+\varepsilon}$ is not Lipschitz. Therefore, 
$\omega_{\varepsilon}$-continuity for a function does not imply that the function 
is Lipschitz.
\label{we}
\end{ejem}

\begin{ejem}
A modulus of continuity $\omega$ satisfying $\omega(s) = 1/\log(1/s)$ for every 
$s$ small is weaker than any modulus $s \mapsto s^{\tau}$, where $\tau > 0$.
\label{wlog}
\end{ejem} \end{small}

For our construction, one of the main problems will consist in controlling the modulus of 
continuity for the derivatives of 
maps obtained by fitting together infinitely many diffeomorphisms defined 
on small intervals. To do this, the following elementary lemma will be quite useful.

\begin{lem} {\em Let $\{I_n\!: n \in \mathbb{N}\}$ be a family of closed intervals in 
$[0,1]$ (resp. in $\clo$) having disjoint interiors and such that the complement of 
their union has zero Lebesgue measure. Suppose that $\varphi$ is a homeomorphism of 
$[0,1]$ such that its restrictions to each interval $I_n$ are $C^{1+\omega}$ 
diffeomorphisms which are $C^1$-tangent to the identity at the endpoints of 
$I_n$ and whose derivatives have $\omega$-norms bounded from above by a constant 
$C$. Then $\varphi$ is a $C^{1+\omega}$ diffeomorphism of the whole interval 
$[0,1]$ (resp. of $\clo$), and the $\omega$-norm of its derivative is less 
than or equal to $2C$.}
\label{pegar}
\end{lem}

\noindent{\bf Proof.} We will just consider the case of the interval, since that of 
the circle is analogous. Let $x\!<\!y$ be points in $\cup_{n \in \mathbb{N}} I_n$. 
If they belong to the same interval $I_n$ then, by hypothesis,
$$\left| \frac{\varphi'(y) - \varphi'(x)}{\omega(y-x)} \right| \leq C.$$
Now take $x \in I_i = [x_i,y_i]$ and $y \in I_j = [x_j,y_j]$, with 
$y_i \leq x_j$. In this case we have
\begin{eqnarray*}
\left| \frac{\varphi'(y)\!-\!\varphi'(x)}{\omega(y-x)} \right| 
&=& \left|\frac{(\varphi'(y)\!-\!1)+(1\!-\!\varphi'(x))}{\omega(y-x)}\right| 
\esp \esp \leq \esp \esp 
\left|\frac{\varphi'(y)\!-\!\varphi'(x_j)}{\omega(y-x)} \right|
+ \left| \frac{\varphi'(y_i)\!-\!\varphi'(x)}{\omega(y-x)}\right| \\
&\leq& C \left[ \frac{\omega(y-x_j)}{\omega(y-x)} +
\frac{\omega(y_i-x)}{\omega(y-x)} \right] \esp \esp \leq \esp \esp 2C.
\end{eqnarray*}
The map $x\!\mapsto \!\varphi'(x)$ is therefore uniformly continuous on the dense 
subset $\cup_{n\in\mathbb{N}} I_n$, and hence extends to a continuous function 
defined on $[0,1]$ whose derivative has $\ce^{\omega}$-norm bounded from above by 
$2C$. Since $I \setminus \cup_{n\in\mathbb{N}} I_n$ has zero measure, this function 
must coincide (at every point) with the derivative of $\varphi$. $\hfill\square$

\vspace{0.45cm}

Because of the lemma above, it would be nice to dispose 
of a ``good'' family of diffeomorphisms between intervals.

\vspace{0.05cm}

\begin{defn} A family  $\{ \varphi_{a,b}\!: [0,a]\!\rightarrow
[0,b]; \esp a > 0, \esp b > 0 \}$ of homeomorphisms is said to be 
{\bf{\em equivariant}} if \esp $\varphi_{b,c} \circ \varphi_{a,b} 
= \varphi_{a,c}$ \esp for all positive $a, b, c$.
\label{def-equiv}\index{equivariant family of diffeomorphisms}
\end{defn}

Given an equivariant family and two intervals  
$I\!\!=\!\![x_1,x_2]$ and $J\!\!=\!\![y_1,y_2]$, let us denote by 
$\varphi(I,J)\!:I \!\rightarrow\! J$ the homeomorphism defined by 
$$\varphi(I,J)(x) = \varphi_{x_2-x_1,y_2-y_1}(x-x_1) + y_1.$$
Notice that $\varphi(I,I)$ must coincide with the identity.

Perhaps the simplest equivariant family of diffeomorphisms is the one formed by the linear 
maps $\varphi_{a,b}(x) = bx/a$. Nevertheless, it is clear that this family is not appropriate 
for fitting maps together in a smooth way. To overcome this difficulty, let us introduce a 
general procedure for constructing equivariant families. Given a family of homeomorphisms 
$\{ \varphi_{a} \!: \esp ]0,a[ \rightarrow \mathbb{R}; \esp a>0 \}$, let us define 
$\varphi_{a,b} = \varphi_{b}^{-1} \circ \varphi_a \!: \hspace{0.1cm} ]0,a[ \rightarrow ]0,b[$. 
We have 
$$\varphi_{b,c} \circ \varphi_{a,b} = (\varphi_c^{-1} \circ \varphi_b) \circ
(\varphi_b^{-1} \circ \varphi_a) = \varphi_c^{-1} \circ \varphi_a = \varphi_{a,c}.$$
Thus, extending $\varphi_{a,b}$ continuously to the whole 
interval $[0,a]$ by setting $\varphi_{a,b}(0) = 0$ and 
$\varphi_{a,b}(a)=b$, we obtain an equivariant family. 
For some of our purposes, a good equivariant family is obtained via this procedure 
using the maps $\varphi_a \!: \esp ]0,a[ \rightarrow \mathbb{R}$ defined by 
\begin{equation}
\varphi_a (x) = - \frac{1}{a} \esp \mathrm{ctg} \left( \frac{\pi x}{a} \right).
\label{yoccocito}
\end{equation}
The associated equivariant family was introduced by Yoccoz. The elements of 
this family satisfy remarkable differentiability properties that we now discuss.
\index{Yoccoz}

Letting $u = \varphi_a (x)$, we have 
$$\varphi_{a,b}'(x) = (\varphi_b^{-1})' (\varphi_a(x)) \cdot \varphi_a'(x) 
= \frac{(\varphi_b^{-1})'(u)}{(\varphi_a^{-1})'(u)}
= \frac{u^2 + 1/a^2}{u^2 + 1/b^2}.$$
Notice that if $x \!\rightarrow\! 0$ (resp. $x \!\rightarrow\! a$), 
then $u \!\rightarrow\! -\infty$ (resp. $u \!\rightarrow\! +\infty$) 
and $\varphi_{a,b}'(x) \!\rightarrow\! 1$. Therefore, the map $\varphi_{a,b}$ 
extends to a $\ce^1$ diffeomorphism from $[0,a]$ into $[0,b]$ which is tangent 
to the identity at the endpoints. Moreover, for $a \!\geq\! b$ (resp. $a \!\leq\! b$), 
the function $u \!\mapsto\! \frac{u^2 + 1/a^2}{u^2 + 1/b^2}$ attains its minimum 
(resp. maximum) value at $u \!=\! 0$; this value being equal to $b^2/a^2$, we have 
$$\sup_{x \in [0,a]} |\varphi_{a,b}'(x) - 1| = \left| \frac{b^2}{a^2} - 1 \right|.$$

For the second derivative of $\varphi_{a,b}$ we have 
\begin{eqnarray*}
\big| \varphi_{a,b}''(x) \big| &=& \frac{d \varphi_{a,b}'(x)}{du} 
\cdot \big| \frac{du}{dx} \big| \esp \esp = \esp \esp 
\frac{\big|2u(u^2+1/b^2) - 2u(u^2+1/a^2)\big|}{(u^2+1/b^2)^2} \esp \pi (u^2 + 1/a^2) \\
&=& \pi \frac{u^2 + 1/a^2}{(u^2 + 1/b^2)^2}
\left| \Big[ 2u \Big( \frac{1}{b^2} - \frac{1}{a^2} \Big) \Big] \right| 
\esp \esp = \esp \esp \pi \frac{u^2 + 1/a^2}{u^2 + 1/b^2} 
\cdot \frac{|2u(1/b^2 - 1/a^2)|}{u^2 + 1/b^2}.
\end{eqnarray*}
From this it follows that $\varphi_{a,b}$ is a $\ce^2$ diffeomorphism satisfying  
$\varphi_{a,b}''(0) = \varphi_{a,b}''(a)=0$. The inequality $\frac{2|u|}{u^2
+t^2}\leq \frac{1}{t}$ applied to $t = 1/b$ yields
$$\left| \varphi_{a,b}''(x) \right| \leq
\pi \frac{u^2 + 1/a^2}{u^2 + 1/b^2} \left| \frac{1}{b^2} - \frac{1}{a^2} \right|b.$$
For $a \leq b$, this implies that 
$$\left| \varphi_{a,b}''(x) \right| \esp \esp \leq \esp \esp 
\pi \frac{b^2}{a^2}\left( \frac{b^2 - a^2}{a^2b^2} \right)b \esp \esp 
= \esp \esp \frac{\pi b}{a^2} \left( \frac{b^2}{a^2}-1 \right)\!.$$
Hence, if \esp $a \leq b \leq 2a$ \esp then 
$$\left| \varphi_{a,b}''(x) \right| \leq 6\pi \left| \frac{b}{a} - 1 \right|
\frac{1}{a}.$$
Analogously, if \esp $2b \geq a \geq b$ \esp then 
$$\left| \varphi_{a,b}''(x) \right| \leq
\frac{\pi}{b} \left(1 - \frac{b^2}{a^2} \right) \leq 2\pi \left| \frac{b}{a}-1
\right| \frac{1}{b}
\leq 4\pi \left| \frac{b}{a} - 1 \right| \frac{1}{a}.$$
Therefore, in both cases we have 
\begin{equation}
\left| \varphi_{a,b}''(x) \right| \leq 6 \pi \left| \frac{b}{a} - 1 \right|
\frac{1}{a}.
\label{tres}
\end{equation}
The last inequality together with the exercise below show that the family of maps 
$\varphi_{a,b}$ is nearly optimal.

\begin{small} \begin{ejer}
Let $\varphi\!: [0,a] \rightarrow [0,b]$ be a $\ce^2$ diffeomorphism. Suppose 
that $\varphi'(0)\!=\!\varphi'(a)\!=\!1$. Show that there exists a point 
$s \! \in ]0,a[$ such that
$$\big| \varphi''(s) \big| \geq \frac{2}{a} \left| \frac{b}{a}-1\right|.$$
\end{ejer} \end{small}

With the only exception of Exercise \ref{caqui}, in what follows the notation 
$\varphi_{a,b}$ and $\varphi_{I,J}$ 
will be only used for denoting the maps in Yoccoz' family. Without 
loss of generality, we will suppose that the function $s\!\mapsto\! \omega(s)/s$ 
is decreasing. (Notice that the moduli from Examples \ref{wholder}, \ref{we}, 
and \ref{wlog}, may be modified far from the origin so that they satisfy 
this property.) Under this assumption we have the following.

\vspace{0.1cm}

\begin{lem} {\em If $a>0$ and $b>0$ are such that \esp $a/b \leq 2$, 
\esp $b/a \leq 2$, \esp and
$$\left| \frac{b}{a} - 1 \right| \frac{1}{\omega(a)} \leq C,$$
then the $\mathrm{C}^{\omega}$-norm of $\varphi_{a,b}'$ 
is less than or equal to $6\pi C$.}
\label{simple}
\end{lem}

\noindent{\bf Proof.} According to 
(\ref{tres}), for every $x \in [0,a]$ one has 
$$|\varphi_{a,b}''(x)| \leq \frac{6 \pi C \omega(a)}{a}.$$
For every $y < z$ in $[0,a]$ there exists $x \in [y,z]$ satisfying  
$\varphi_{a,b}'(z) - \varphi'(y) = \varphi''_{a,b}(x) (z-y)$. Since 
the function $s \mapsto \omega(s)/s$ is decreasing 
and \esp $z\!-\!y \leq a$, \esp this implies that 
$$\left| \frac{\varphi_{a,b}'(z) - \varphi'_{a,b}(y)}{\omega(z-y)} \right| \esp 
= \esp |\varphi''_{a,b}(x)| \left| \frac{z-y}{\omega(z-y)} \right| \esp 
\leq \esp |\varphi''_{a,b}(x)| \left| \frac{a}{\omega(a)} \right| \esp 
\leq \esp 6 \pi C,$$
which shows the lemma. $\hfill\square$

\vspace{0.3cm}

On the basis of the preceding discussion, we can now give a conceptual construction 
of the so-called \textit{Denjoy counterexamples} (\textit{i.e.,} we can give a 
proof of Theorem \ref{herman-contr}). The method of proof will be used later 
on for smoothing many other group actions on the interval and the circle.
\index{Herman}

\vspace{0.42cm}

\noindent{\bf Proof of Theorem \ref{herman-contr}.} Slightly more generally, for all 
$\varepsilon > 0$, every irrational angle $\theta$, and every $k \in \mathbb{N}$, 
we will exhibit a $\ce^{1+\omega_{\varepsilon}}$ circle diffeomorphism with rotation 
number $\theta$ whose derivative has $\ce^{\omega_{\varepsilon}}$-norm bounded from 
above by $2/[\log(k)]^{\varepsilon/2}$. To do this let us fix $x_0 \!\in \!\clo$, 
and for each $n \!\in\! \mathbb{Z}$ let us replace each point 
$x_n = R_{\theta}^n(x_0)$ of its orbit by an interval $I_n$ of length 
$$\ell_n = \frac{1}{(|n|+k)[\log(|n|+k)]^{1+\varepsilon/2}}.$$
The original rotation induces a homeomorphism $h\!=\!\bar{R}_{\theta}$ 
of a circle of length $\bar{\ell}_k = \sum_{n \in \mathbb{Z}} \ell_n$ be letting 
$\bar{R}_{\theta}(x) = \varphi(I_n,I_{n+1}) (x)$ for $x \!\in\! I_n$ and by extending 
continuously to $\clo$. We now check that $\bar{R}_{\theta}$ is a diffeomorphism 
satisfying the desired properties. For this, we need to estimate the value of 
expressions of type
\begin{equation}
\left| \frac{\ell_{n+1}}{\ell_n} - 1 \right| \frac{1}{\omega_{\varepsilon}(\ell_n)}.
\label{acagar}
\end{equation}
We will make the explicit computations for $n \geq 0$, leaving the 
case $n \leq 0$ to the reader. Expression (\ref{acagar}) equals 
$$\left|
\frac{(n+k)[\log(n+k)]^{1+\varepsilon/2}}{(n+k+1)[\log(n+k+1)]^{1+\varepsilon/2}}-1
\right| \frac{(n+k)[\log(n+k)]^{1+\varepsilon/2}}
{[\log(n+k) + (1+\varepsilon/2)\log \log (n+k)]^{1+\varepsilon}},$$
which is bounded from above by 
$$\left|
\frac{(n+k)[\log(n+k)]^{1+\varepsilon/2}}{(n+k+1)[\log(n+k+1)]^{1+\varepsilon/2}}-1
\right| \frac{n+k}{[\log(n+k)]^{\varepsilon/2}}.$$
By applying the Mean Value Theorem to the function 
$s \mapsto s[\log(s)]^{1+\varepsilon/2}$, we obtain 
the following upper bound for the latter expression: 
$$\frac{[\log(n+k+1)]^{1+\varepsilon/2} +
(1+\varepsilon/2)[\log(n+k+1)]^{\varepsilon/2}}
{[\log(n+k+1)]^{1+\varepsilon/2} \esp [\log(n+k)]^{\varepsilon/2}} 
\esp \leq \esp \frac{2}{[\log(n+k)]^{\varepsilon/2}}.$$
The diffeomorphism $h$ that we constructed acts on a circle of length  
$\bar{\ell}_k \sim 2/\varepsilon [\log(k)]^{\varepsilon/2}$. Therefore, to 
obtain a diffeomorphism of the unit circle, we must renormalize the circle, 
say by an affine map $\varphi$. Notice that this procedure does not increase the  
$\ce^{\omega_{\varepsilon}}$-norm for the derivative. Indeed, from the equality 
\esp $\bar{R}_{\theta}' = (\varphi\circ h\circ\varphi^{-1})' =
\bar{\ell}_k^{-1} \cdot (h'\circ\varphi^{-1}) \cdot  L_k$ \esp 
one deduces that
\begin{eqnarray*}
\frac{|\bar{R}_{\theta}'(x) - \bar{R}_{\theta}'(y)|}
{\omega_{\varepsilon}(|x-y|)} 
&=& \frac{|(\varphi\circ h\circ\varphi^{-1})'(x) - (\varphi\circ
h\circ\varphi^{-1})'(y)|}
{\omega_{\varepsilon}(|\varphi^{-1}(x) - \varphi^{-1}(y)|)} \cdot
\frac{\omega_{\varepsilon}(|\varphi^{-1}(x)-\varphi^{-1}(y)|}{\omega_{\varepsilon}(|x-y|)}\\
&\leq& \frac{|h'(\varphi^{-1}(x)) - h'(\varphi^{-1}(y))|}  
{\omega_{\varepsilon}(|\varphi^{-1}(x) - \varphi^{-1}(y)|)}
\esp \esp \leq \esp \esp \frac{2}{[\log(k)]^{\varepsilon/2}},
\end{eqnarray*}
where the first inequality comes from the fact that 
$\omega_{\varepsilon}$ is increasing and $\bar{\ell}_k < 1$. $\hfill\square$

\vspace{0.43cm}

The preceding construction might suggest that there exist Denjoy counterexamples 
for any modulus of continuity for the derivative 
weaker than the Lipschitz one. However, there exist subtle 
obstructions related to the Diophantine nature of the rotation number which are not 
completely understood: see \cite[Chapter X]{He} for a partial result on this.

\vspace{0.12cm}

\begin{small} \begin{ejer} Show that there is no Denjoy counterexample of class 
$\mathrm{C}^{1+\omega_0}$ whose derivative is identically equal to $1$ on the minimal 
invariant Cantor set, where $\omega_0 (s) = s \log(1/s)$ (compare Lemma \ref{arg-norton} 
and the argument before it).

\noindent{\underbar{Remark.}} It seems to be unknown whether the claim above 
is still true without the hypothesis on the derivative over the invariant Cantor 
set. Actually, the general problem of finding the weakest modulus of continuity  
ensuring the validity of Denjoy's theorem seems interesting. 
\end{ejer} \end{small}

\vspace{0.02cm} 

We would now like to use the preceding technique for constructing smooth actions on 
the circle 
of higher-rank Abelian groups with an exceptional minimal set. However, Proposition 
\ref{arg-norton} suggests that we should find obstructions in differentiability 
classes less than $\ce^{1+\mathrm{Lip}}$. Indeed, this will be the main object 
of the next section where we will show, for instance, that for every $d \geq 2$ and all 
$\varepsilon > 0$, every free action of $\mathbb{Z}^d$ by $\ce^{1+1\!/\!d+\varepsilon}$ 
circle diffeomorphisms is minimal.\footnote{Actually, a stronger result holds: Every free 
$\mathbb{Z}^d$-action by circle diffeomorphisms is minimal provided that the generators 
$f_i$, $i \in \{1,\ldots,d\}$, are respectively of class $\ce^{1+\tau_i}$, and \esp 
$\tau_1 + \cdots + \tau_d > 1$ \esp (see \cite{nuestro}).} 
According to \cite{TsP} 
\index{Tsuboi!conjecture}
(see also \cite{tsuboi-fourier}), the hypothesis $\varepsilon\!>\!0$ should be superfluous 
for this result. More generally, the theorem should be true for actions by diffeomorphisms 
with bounded $d$-variation ({\em c.f.,} Exercise \ref{cuadr}). Here we will content ourselves 
with showing that, in regularity less than $\ce^{1+1/d}$, this result is no longer true.

\vspace{0.1cm}

\begin{prop} {\em For every $\varepsilon > 0$ and every positive integer $d$, 
there exist free $\mathbb{Z}^d$-actions by $\ce^{1+1\!/\!d-\varepsilon}$ 
circle diffeomorphisms admitting an exceptional minimal set.}
\end{prop}

\noindent{\bf Proof.} Notice that the case $d\! = \!1$ corresponds to Theorem 
\ref{herman-contr}. However, the construction involves simultaneously all 
the cases, and produces examples of groups of circle diffeomorphisms 
of class $\ce^{1+\omega}$,  where \esp 
$\omega(s) = s^{1\!/\!d} \esp [\log(1/s)]^{1\!/\!d + \varepsilon}$.

Let us begin by fixing a rank-$d$ torsion-free subgroup of $\mathrm{SO}(2,\mathbb{R})$, 
and let $\theta_1, \ldots,\theta_d$ be the angles of the generators. Let $m \geq d-1$ 
be an integer number and $p$ a point in $\clo$. For each $(i_1,\ldots,i_d) \in
\mathbb{Z}^d$ let us replace the point $R_{\theta_1}^{i_1} \cdots
R_{\theta_d}^{i_d}(p)$ by an interval $I_{i_1,\ldots,i_d}$ of length  
$$\ell_{(i_1,\ldots,i_d)} \esp = \esp 
\frac{1}{\big(|i_1|+\cdots+|i_d|+m\big)^d \thinspace
\big[\log(|i_1|+\cdots+|i_d|+m)\big]^{1+\varepsilon}}.$$
This procedure induces a new circle (of length \esp 
$\bar{\ell}_{m} \!\leq\! 2^d / \varepsilon [\log(m)]^{\varepsilon}\! (d-1)!$), 
\esp upon which the rotations $R_{\theta_j}$ induce homeomorphisms $f_j$ 
satisfying, for every $x \in I_{i_1,\ldots,i_j,\ldots,i_d}$,
$$f_j(x) = \varphi(I_{i_1,\ldots,i_j,\ldots,i_d},I_{i_1,\ldots,1+i_j,\ldots,i_d})(x).$$
By the equivariance properties of the $\varphi(I,J)$'s, these homeomorphisms $f_j$ 
commute between them. The verification that each $f_i$ is of class $\ce^{1+\omega}$ 
is analogous to the proof of Theorem \ref{herman-contr}, and we leave it to the 
reader. Finally, notice that by renormalizing the circle and letting $m$ go 
to infinity, each $f_i$ converges (in the $\ce^{1+\omega}$ topology) to the 
corresponding rotation $R_{\theta_i}$. $\hfill\square$

\vspace{0.4cm}

An analogous procedure leads to counterexamples to Kopell's lemma.  
\index{Kopell!lemma} 
Slightly more generally, for each integer $d \geq 2$ and each $\varepsilon > 0$, there 
exist $\ce^{1+1\!/\!(d-1)-\varepsilon}$ diffeomorphisms $f_1,\ldots,f_{d}$ of $[0,1]$ 
and disjoint open intervals $I_{n_1,\ldots,n_{d}}$ disposed in $]0,1[$ according to 
the lexicographic ordering so that, for every $(n_1,\ldots,n_d) \! \in \! \mathbb{Z}^d$ 
\index{lexicographic ordering}
and every $j\! \in \! \{1,\ldots,d\}$,
$$f_j (I_{n_1,\ldots,n_j,\ldots,n_{d}}) = I_{n_1,\ldots,n_j-1,\ldots,n_{d}}.$$
To construct these diffeomorphisms, a natural procedure goes as follows. Given 
an integer $m \geq d-1$, for each $(i_1,\ldots,i_{d}) \in \mathbb{Z}^{d}$ let 
$$\ell_{i_1,\ldots,i_{d}} \esp = \esp \frac{1}{\big(|i_1|+\cdots+|i_{d}|+m\big)^d 
\thinspace \big[\log(|i_1|+\cdots+|i_{d}|+m)\big]^{1+\varepsilon}}.$$
Let us inductively define 
\thinspace $\ell_{i_1,\ldots,i_{j-1}} = \sum_{i_j \in \mathbb{Z}}
\ell_{i_1,\ldots,i_{j-1},i_j}$. \thinspace Let \thinspace
$[x_{i_1,\ldots,i_j,\ldots,i_d},y_{i_1,\ldots,i_j,\ldots,i_d}]$
\thinspace be the interval whose endpoints are  
$$x_{i_1,\ldots,i_j,\ldots,i_{d}} = \sum_{i_1' < i_1} \ell_{i_1'} +
\sum_{i_2' < i_2} \ell_{i_1,i_2'} + \ldots +
\sum_{i_{d}' < i_{d}} \ell_{i_1,\ldots,i_{d},i_{d}'},$$
$$y_{i_1,\ldots,i_j,\ldots,i_{d}} =
x_{i_1,\ldots,i_j,\ldots,i_{d}} + \ell_{i_1,\ldots,i_{d}},$$
and let $f_j$ be the diffeomorphism of $[0,1]$ that restricted to each 
$[x_{i_1,\ldots,i_j,\ldots,i_{d}},y_{i_1,\ldots,i_j,\ldots,i_{d}}]$ coincides with 
$$\varphi \big( [x_{i_1,\ldots,i_j,\ldots,i_{d}},y_{i_1,\ldots,i_j,\ldots,i_{d}}],
[x_{i_1,\ldots,i_j-1,\ldots,i_{d}},y_{i_1,\ldots,i_j-1,\ldots,i_{d}}] \big).$$
Using some of the estimates in this section, one readily checks that the obtained 
$f_j$'s are of class $\ce^{1+\omega}$, where 
$\omega(s) = s^{1\!/\!d} \esp [\log(1/s)]^{1\!/\!d+\varepsilon}$. In particular, 
these maps are $\ce^{1+1\!/\!d - \varepsilon}$ diffeomorphisms. However, 
notice that, by this 
method, we have not achieved the optimal regularity $\ce^{1+1\!/\!(d-1)-\varepsilon}$ 
that we claimed. Actually, a direct consequence of Proposition \ref{arg-norton} is 
that by this procedure it is impossible to reach the class $\ce^{1+1\!/\!d}$. To 
reach the optimal regularity $\ce^{1+1\!/\!(d-1)-\varepsilon}$, it is necessary to 
avoid many of the tangencies to the identity of the maps. To do this we will use 
the original technique by Pixton \cite{Pix}, following the brilliant presentation 
of \index{Tsuboi!Pixton-Tsuboi actions} Tsuboi \cite{TsP}.

Let $X = \varrho \esp \frac{\partial}{\partial x}$ be a $\ce^{\infty}$ vector field 
on $[0,1]$ satisfying \esp $|\varrho'(x)| \leq 1$ \esp for all $x \!\in\! [0,1]$, and 
such that \esp $\varrho(x) = x$ \esp for every \esp $x \!\in\! [0,1/3]$, \esp and \esp 
$\varrho(x) = 0$ \esp for every \esp $x \!\in\! [1/2,1].$ Let $\varphi^t(x)$ be the 
flow associated to $X$, that is, the solution of the differential equation 
$$\frac{d \varphi^t (x)}{dt}(x) = \varrho \big( \varphi^t(x)\big), \qquad
\varphi^0(x) = x.$$
For every positive $a,b$ and each $t \geq 0$, the diffeomorphism 
$x \mapsto b \esp \varphi^t (x/a)$ sends the interval $[0,a]$ onto $[0,b]$. Moreover, 
its derivative is identically equal to $b/a$ in $[a/2,a]$, and equals $be^t/a$ 
at the origin. Given $a'\!<\!0\!<\!a$ and $b'\!<\!0\!<\!b$, let us define a 
diffeomorphism $\varphi_{b',b}^{a',a}\!: [0,a]\rightarrow [0,b]$ by letting
$$\varphi^{a',a}_{b',b}(x) \esp = \esp 
b \esp \varphi^{\log(b'a/a'b)} \Big( \frac{x}{a} \Big).$$
Notice that for every $c'<0<c$ one has the equivariance property
\begin{equation}
\varphi^{b',b}_{c',c} \circ \varphi^{a',a}_{b',b} \esp = \esp \varphi^{a',a}_{c',c},
\label{abc}
\end{equation}
which is analogous to that in Definition \ref{def-equiv}. Nevertheless, the fact that 
in the family $\{\varphi^{a',a}_{b',b}\}$ four (and not only two) parameters are involved 
allows obtaining a better control for the modulus of continuity of the derivatives (with 
the mild cost of having to suppress the tangencies to the identity at the endpoints).

\vspace{0.15cm}

\begin{lem} {\em Letting \esp $C = \sup_{y \in [0,1]} \varrho''(y)$,  
\esp for every $x \in [0,a]$ one has the inequality}
\begin{equation}
\left| \frac{\partial}{\partial x} \log \Big( \frac{\partial
\varphi^{a',a}_{b',b}}{\partial x}
\Big)(x) \right| \leq \frac{C}{a} \esp \Big| \frac{b' a}{a' b} - 1 \Big|.
\label{est-pix}
\end{equation}
\end{lem}

\noindent{\bf Proof.} From \esp $d \varphi^t(x) / dt = \varrho \big(
\varphi^t (x) \big)$ \esp one concludes that 
$$\frac{d}{dt} \frac{\partial \varphi^t}{\partial x}(x) = \frac{\partial}{\partial x}
\varrho \big( \varphi^t(x) \big) = \frac{\partial \varrho}{\partial x} \big(
\varphi^t(x) \big)
\cdot \frac{\partial \varphi^t}{\partial x}(x),$$
and hence 
\begin{equation}
\frac{d}{dt} \log \Big( \frac{\partial \varphi^t}{\partial x} \Big) (x) =
\frac{\partial \varrho}{\partial x} \big( \varphi^t(x) \big).
\label{ddt}
\end{equation}
On the other hand, from $|\varrho'(x)| \leq 1$ it follows that 
$$\left| \log \Big( \frac{\partial \varphi^t}{\partial x} \Big)(x) \right|
= \left| \int_0^t \frac{d}{ds} \log \Big( \frac{\partial \varphi^s}{\partial x}
\Big) (x) \esp ds\right| = \left| \int_0^t \frac{\partial \varrho}{\partial x} 
\big( \varphi^s(x)\big) \esp ds \right| \leq \int_0^t ds = t,$$
and therefore
\begin{equation}
\frac{\partial \varphi^t}{\partial x}(x) \esp \leq \esp  e^t.
\label{ddx}
\end{equation}
Since
$$\log \Big( \frac{\partial \varphi^{a',a}_{b',b}}{\partial x} \Big)(x) \esp = 
\esp \log \Big(\frac{b}{a} \Big) + \log \Big( \frac{\partial}{\partial x} 
\varphi^{\log(b'a/a'b)} \Big) \big(\frac{x}{a}\big),$$
using (\ref{ddt}) and (\ref{ddx}) we obtain 
\begin{eqnarray*}
\left| \frac{\partial}{\partial x} \log \Big( \frac{\partial \varphi^{a',a}_{b',b}
}{\partial x} \Big)(x) \right|
\!&=&\! \left| \frac{\partial}{\partial x} \log \Big( \frac{\partial}{\partial x}
\varphi^{\log(b'a/a'b)} \Big)
\big( \frac{x}{a} \big) \right|\\
&=&\! \left| \int_0^{\log(b'a/a'b)} \frac{d}{dt} \frac{\partial}{\partial x} \log
\Big( \frac{\partial \varphi^t}{\partial x} \Big) \big( \frac{x}{a} \big) dt \right|\\
&=&\! \frac{1}{a} \left| \int_0^{\log(b'a/a'b)} \frac{\partial}{\partial x}
\frac{d}{dt} \log
\Big( \frac{\partial \varphi^t}{\partial x} \Big) \big( \frac{x}{a} \big) \cdot
\frac{d \varphi^t}{dt} \big( \frac{x}{a} \big) \esp dt \right|\\
&=&\! \frac{1}{a} \left| \int_0^{\log(b'a/a'b)} \frac{\partial}{\partial x}
\frac{\partial \varrho}{\partial x} \Big( \varphi^t \big( \frac{x}{a} \big) \Big)
\cdot \frac{d \varphi^t}{dt} \big( \frac{x}{a} \big)  \esp dt \right|,
\end{eqnarray*}
and hence 
$$\log \Big( \frac{\partial \varphi^{a',a}_{b',b}}{\partial x} \Big)(x)
\esp \leq \esp \frac{C}{a} \left| \int_{0}^{\log(b'a/a'b)} e^t dt \right| 
\esp = \esp \frac{C}{a} \esp \left| \frac{b'a}{a'b} - 1 \right|.$$

\vspace{-0.85cm}

$\hfill\square$

\vspace{0.9cm}

Notice that 
\begin{small}\begin{equation}
\frac{\partial \varphi^{a',a}_{b',b}}{\partial x} (0) = \frac{b}{a} \cdot
\frac{\partial}{\partial x}
\varphi^{\log(b'a/a'b)}(0) = \frac{b}{a} \cdot \exp \big( \log(b'a/a'b) \big) =
\frac{b'}{a'}, \qquad
\frac{\partial \varphi^{a',a}_{b',b}}{\partial x} (a) = \frac{b}{a}.
\label{abd}
\end{equation}\end{small}Therefore, the 
diffeomorphisms $\varphi^{a',a}_{b',b}$ are not necessarily tangent to the identity 
at the endpoints. This allows fitting together them with certain freedom. A careful 
choice of the (infinitely many) parameters $a'\!<\!0$ and $b'\!<\!0$ will then lead to the 
optimal regularity. In what follows, we will just sketch the (quite technical) construction 
of commuting diffeomorphisms of the interval which are $\ce^{2-\varepsilon}$ counterexamples 
to Kopell's lemma. For the general case of $\mathbb{Z}^d$-actions on the interval we refer 
to \cite{TsP} (see also \cite{edo}).

Let us fix $k \in \mathbb{N}$ and 
$\varepsilon > 0$, and for each pair of integers $m,n$ let 
$$\ell_{m,n} = \frac{1}{\big( |m|^2 + |n|^2 + k^2 \big)^{1+\varepsilon}}.$$ 
Let $\{I_{m,n} \!=\! [x_{m,n},y_{m,n}] \}$ be a family of intervals placed inside a  
closed interval $I$ respecting the lexicographic ordering and so that their union 
has full measure in $I$, with $|I_{m,n}| = \ell_{m,n}$. Define two homeomorphisms 
$f,g$ of $I$ by letting, for each $x \in I_{m,n}$,
$$f(x) = \varphi^{\ell_{m,n-1},\ell_{m,n}}_{\ell_{m-1,n-1},\ell_{m-1,n}}
( x- x_{m,n} ) + x_{m-1,n},$$
$$g(x) = \varphi^{\ell_{m,n-1},\ell_{m,n}}_{\ell_{m,n-2},\ell_{m,n-1}} (x - x_{m,n})
+ x_{m,n-1}.$$
Properties (\ref{abc}) and (\ref{abd}) imply that $f$ and $g$ are commuting $\ce^1$ 
diffeomorphisms of $I$. Using (\ref{est-pix}), the reader should be able to check  
that $f$ and $g$ are actually $\ce^{2-\varepsilon'}$ diffeomorphisms, where 
$\varepsilon'>0$ depends only on $\varepsilon > 0$ and goes to zero together 
with $\varepsilon$. However, we must warn that this is not an easy 
task, mainly because of the absence of an analogue to Lemma 
\ref{pegar} for this case: see \cite{TsP} for the details.
\index{Denjoy!counterexamples|)}

\begin{small}\begin{ejer} 
Slightly abusing of the notation, define \esp $\varphi_{a,b} \!: [0,a] \rightarrow [0,b]$ 
\esp by letting \esp $\varphi_{a,b}(x) = \varphi^{-1,a/2}_{-1,b/2} (x)$ \esp and  
\esp $\varphi_{a,b}(x) \!=\! b \!-\! \varphi^{-1,a/2}_{-1,b/2}(x)$ 
\esp for $x$ in $[a/2,a]$ and $[0,a/2]$, respectively. Check that 
$\{\varphi_{a,b}\}$ is an equivariant family of diffeomorphisms tangent to 
the identity at the endpoints which satisfy (compare inequality (\ref{tres}))
$$\left| \frac{\partial}{\partial x} \log \Big( \frac{\partial
\varphi_{a,b}}{\partial x} \Big) (x)
\right| \leq \frac{2 C}{a} \esp \left| \frac{a}{b} - 1 \right|.$$
\label{caqui}
\end{ejer}\end{small}

\subsection{On intermediate regularities}
\label{sec-inter}

\vspace{0.2cm}

\hspace{0.45cm} The next result (obtained by the author in collaboration 
with Deroin and Kleptsyn in \cite{DKN}) may be considered an extension 
(up to some parameter $\varepsilon > 0$) of Denjoy's theorem for free 
\index{action!free} 
actions of higher-rank Abelian groups on the circle. We ignore whether,  
under the same hypothesis, the corresponding actions are ergodic 
(compare Theorem \ref{katokherman}).
\index{Denjoy!generalized theorem|(}\index{Deroin}\index{Kleptsyn}

\vspace{0.1cm}

\begin{thm} {\em If $d$ is an integer greater than or equal to $2$ 
and $\varepsilon\! >\! 0$, then every free action of $\mathbb{Z}^d$ 
by $\ce^{1+1\!/\!d+\varepsilon}$ circle diffeomorphisms is minimal.}
\label{denjoy-gen}
\end{thm}

\vspace{0.15cm}

Before giving the proof of this result, we would like to explain the main idea, which 
is inspired by the famous {\em Erd\"os principle} and shows once again how fruitful the 
probabilistic methods are in the theory. Suppose for a contradiction that $I$ is a wandering 
interval for the dynamics of a rank-2 Abelian group of $\ce^{1+\tau}$ circle 
diffeomorphisms generated by elements $g_1$ and $g_2$. Letting 
$\mathbb{N}_0 \!=\! \mathbb{N} \cup \{0\}$, the family  
$$\big\{ g_1^m g_2^n(I) \!: 
(m,n) \in \mathbb{N}_0 \times \mathbb{N}_0, \esp m+n \leq k \big\}$$
consists of $(k+1)(k+2)/2$ intervals. Since these intervals are two-by-two disjoint, 
we should expect that, ``typically'', their length 
has order $1/k^2$. Hence, for a ``generic'' random 
sequence $I,h_1(I),h_2(I) \ldots$, where either $h_{n+1} = g_1 h_n$ \esp or 
\esp $h_{n+1} = g_2  h_n$, for $\tau > 1/2$ we should have
$$\sum_{k \geq 1} \big| h_k(I) \big|^{\tau} \esp \leq \esp 
C \sum_{k \geq 1} \frac{1}{k^{2\tau}} \esp < \esp \infty.$$
Now notice that the left-side 
expression of this inequality corresponds to the series which allows 
controlling the distortion for the successive compositions. More precisely, if this 
sum is finite, then Lemma \ref{?} should allow us to find elements with hyperbolic 
fixed points, thus contradicting the hypothesis of freeness for the action.

To formalize the idea above, we need to make the random nature of the compositions 
more precise. For this, let us consider the Markov process on 
$\mathbb{N}_0 \times \mathbb{N}_0$ with transition probabilities 
\begin{small}\begin{equation}
p \big( (m,n) \rightarrow (m+1,n) \big) = \frac{m+1}{m+n+2},
\qquad p \big( (m,n) \rightarrow (m,n+1) \big) = \frac{n+1}{m+n+2}.
\label{paso}
\end{equation}
\end{small}This Markov process induces a probability measure $\mathbb{P}$ on 
the corresponding space of paths $\Omega$. One easily shows that, starting 
from the origin, the probability of arriving at the point $(m,n)$ in $k$ 
steps equals $1/(k+1)$ (resp. 0) if \thinspace $m+n=k$ \thinspace 
(resp. \thinspace $m+n \neq k$).
\index{Markov!process}

\vspace{0.45cm}


\beginpicture

\setcoordinatesystem units <1cm,1cm>

\putrule from 1.5 0 to 1.5 6
\putrule from 3 0 to 3 6
\putrule from 4.5 0 to 4.5 6

\putrule from 0 1.5 to 6 1.5
\putrule from 0 3 to 6 3
\putrule from 0 4.5 to 6 4.5

\putrule from 0 0 to 0 6
\putrule from 0 0 to 6 0

\putrule from 3 -0.6 to 5 -0.6
\putrule from -0.9 2 to -0.9 4

\plot 5 -0.6
4.85 -0.65 /

\plot 5 -0.6
4.85 -0.55 /

\plot -0.9 4
-0.95 3.85 /

\plot -0.9 4
-0.85 3.85 /

\put{Figure 21} at 3.9 -1.3


\plot 1 0
0.85 0.05 /

\plot 1 0
0.85 -0.05 /

\plot 1 1.5
0.85 1.55 /

\plot 1 1.5
0.85 1.45 /

\plot 1 3
0.85 2.95 /

\plot 1 3
0.85 3.05 /

\plot 1 4.5
0.85 4.55 /

\plot 1 4.5
0.85 4.45 /



\plot 2.5 0
2.35 0.05 /

\plot 2.5 0
2.35 -0.05 /

\plot 2.5 1.5
2.35 1.55 /

\plot 2.5 1.5
2.35 1.45 /

\plot 2.5 3
2.35 2.95 /

\plot 2.5 3
2.35 3.05 /

\plot 2.5 4.5
2.35 4.55 /

\plot 2.5 4.5
2.35 4.45 /



\plot 4 0
3.85 0.05 /

\plot 4 0
3.85 -0.05 /

\plot 4 1.5
3.85 1.55 /

\plot 4 1.5
3.85 1.45 /

\plot 4 3
3.85 2.95 /

\plot 4 3
3.85 3.05 /

\plot 4 4.5
3.85 4.55 /

\plot 4 4.5
3.85 4.45 /


\plot 5.5 0
5.35 0.05 /

\plot 5.5 0
5.35 -0.05 /

\plot 5.5 1.5
5.35 1.55 /

\plot 5.5 1.5
5.35 1.45 /

\plot 5.5 3
5.35 2.95 /

\plot 5.5 3
5.35 3.05 /

\plot 5.5 4.5
5.35 4.55 /

\plot 5.5 4.5
5.35 4.45 /



\plot 0 1
-0.05 0.85 /

\plot 0 1
0.05 0.85 /

\plot 1.5 1
1.45 0.85 /

\plot 1.5 1
1.55 0.85 /

\plot 3 1
2.95 0.85 /

\plot 3 1
3.05 0.85 /

\plot 4.5 1
4.55 0.85 /

\plot 4.5 1
4.45 0.85 /



\plot 0 2.5
-0.05 2.35 /

\plot 0 2.5
0.05 2.35 /

\plot 1.5 2.5
1.45 2.35 /

\plot 1.5 2.5
1.55 2.35 /

\plot 3 2.5
2.95 2.35 /

\plot 3 2.5
3.05 2.35 /

\plot 4.5 2.5
4.45 2.35 /

\plot 4.5 2.5
4.55 2.35 /



\plot 0 4
-0.05 3.85 /

\plot 0 4
0.05 3.85 /

\plot 1.5 4
1.45 3.85 /

\plot 1.5 4
1.55 3.85 /

\plot 3 4
2.95 3.85 /

\plot 3 4
3.05 3.85 /

\plot 4.5 4
4.45 3.85 /

\plot 4.5 4
4.55 3.85 /



\plot 0 5.5
-0.05 5.35 /

\plot 0 5.5
0.05 5.35 /

\plot 1.5 5.5
1.45 5.35 /

\plot 1.5 5.5
1.55 5.35 /

\plot 3 5.5
2.95 5.35 /

\plot 3 5.5
3.05 5.35 /

\plot 4.5 5.5
4.45 5.35 /

\plot 4.5 5.5
4.55 5.35 /


\setdots

\plot
0 0
6 6 /

\small

\put{$(0,0)$} at -0.4 -0.4

\put{$p \big( (m,n) \rightarrow (m,n+1) \big) \geq \frac{1}{2}$} at 3 6.5

\put{$n \geq m \hspace{0.3cm} \Longrightarrow $} at 3 6.9

\put{$p \big( (m,n) \rightarrow (m+1,n) \big) \geq \frac{1}{2}$} at 7 3.58

\put{$m \geq n \hspace{0.3cm} \Longrightarrow $} at 7 3.98

\put{$1/2$} at 0.6 -0.3
\put{$1/3$} at 0.6 1.2
\put{$1/4$} at 0.6 2.7
\put{$2/3$} at 2.1 -0.3
\put{$2/4$} at 2.1 1.2
\put{$3/4$} at 3.6 -0.3

\put{$1/2$} at -0.4 0.55
\put{$1/3$} at 1.1 0.55
\put{$1/4$} at 2.6 0.55
\put{$2/3$} at -0.4 2.05
\put{$2/4$} at 1.1 2.05
\put{$3/4$} at -0.4 3.55

\put{$g_1$} at 4 -0.8
\put{$g_2$} at -1.1 3

\put{$\bullet$} at 0 0

\put{} at -2.7 0

\endpicture


\vspace{0.65cm}

To prove Theorem \ref{denjoy-gen} in the case $d\!=\!2$, 
let $g_1$ and $g_2$ be two commuting $\ce^{1+\tau}$ circle 
diffeomorphisms generating a rank-2 Abelian group. 
\index{semigroup} 
The semigroup $\Gamma^+$ generated by $g_1$ and $g_2$ identifies with  
$\mathbb{N}_0 \times \mathbb{N}_0$; therefore, the Markov process already described 
induces a ``random walk'' on $\Gamma^+$. In what follows, we will identify $\Omega$ with 
the corresponding space of paths over $\Gamma^+$. For every $\omega \in \Omega$ and all 
$n \in \mathbb{N}$, we denote by $h_n(\omega)\in \Gamma^+$ the product of the first $n$ 
entries of $\omega$. In other words, for $\omega = (g_{i_1},g_{i_2},\ldots) \in \Omega$ 
we let $h_n(\omega) = g_{i_n} \cdots g_{i_1}$ (where $g_{i_j}$ equals either $g_1$ or 
$g_2$), and we let $h_0(\omega) = id$.

If the action of $\Gamma \!=\! \langle g_1, g_2 \rangle \!\sim \!\mathbb{Z}^2$ is free 
then, by H\"older's theorem and \S \ref{medinv}, the restriction of the rotation number 
function to $\Gamma$ is a group homomorphism. This implies that the rotation numbers 
\index{rotation number} 
$\rho(g_1)$ and $\rho(g_2)$ are {\bf {\em independent}} over 
the rationals, in the sense that for all $(r_0,r_1,r_2) \in \mathbb{Q}^3$ distinct 
from $(0,0,0)$ one has \esp $r_1 \rho(g_1) + r_2 \rho(g_2) \neq r_0$. \esp Indeed, 
otherwise one could find nontrivial elements with rational rotation number. 
These elements would then have periodic points, and since they are of 
infinite order, this would contradict the freeness of the action.

Now suppose for a contradiction that the action of $\Gamma$ is non minimal. In this case, 
there exists a minimal invariant Cantor set for the action. Moreover, every connected 
component $I$ of the complement of this set is wandering for the dynamics. 
\index{wandering interval}

\vspace{0.1cm}

\begin{lem} {\em If $\tau > 1/2$ then the series \thinspace $\ell_{\tau} (\omega)  
= \sum_{n \geq 0} |h_n(\omega)(I)|^{\tau}$ \thinspace converges 
for $\mathbb{P}$-almost every $\omega \in \Omega$.}
\label{2d}
\end{lem}

\noindent{\bf Proof.} We will show that, if $\tau > 1/2$, 
then the expectation (with respect to $\mathbb{P}$) of the 
function $\ell_{\tau}$ is finite, which obviously implies the claim of the lemma. 
For this, first notice that since the arrival probabilities in $k$ steps are equally 
distributed over the points at simplicial distance $k$ from the origin,
$$\mathbb{E}(\ell_{\tau})
= \mathbb{E} 
\Big( \sum_{k \geq 0} \big| h_k(\omega)(I) \big|^{\tau} \Big)
= \sum_{k \geq 0} \mathbb{E} \big( \big| h_k(\omega)(I) \big|^{\tau} \big)
= \sum_{k \geq 0} \sum_{m+n=k} \frac{\big| g_1^m g_2^n(I) \big|^{\tau}}{k+1}.$$
By H\"older's inequality,  
\index{H\"older!inequality}
$$\sum_{m+n=k} \frac{\big| g_1^m g_2^n(I) \big|^{\tau}}{k+1} \leq
\left( \sum_{m+n=k} \big| g_1^m g_2^n(I) \big| \right)^{\tau}
\left( (k+1) \cdot \frac{1}{(k+1)^{1/(1-\tau)}}\right)^{1-\tau},$$
and hence
$$\mathbb{E}(\ell_{\tau}) \leq \sum_{k \geq 0}
\frac{\left(\sum_{m+n=k} \big| g_1^m g_2^n(I) \big|\right)^{\tau}}{(k+1)^{\tau}}.$$
By applying H\"older's inequality once again we obtain
$$\mathbb{E}(\ell_{\tau}) \leq
\left[ \sum_{(m,n) \in \mathbb{N}_0 \times \mathbb{N}_0}
\big| g_1^m g_2^n(I) \big| \right]^{\tau} \left[ \sum_{k \geq 1} 
\left( \frac{1}{k^{\tau}} \right)^{\frac{1}{1-\tau}} \right]^{1-\tau}.$$
Since $\tau > 1/2$, the series
$$\sum_{k \geq 1} \left( \frac{1}{k^{\tau}} \right)^{\frac{1}{1-\tau}} =
\sum_{k \geq 1} \frac{1}{k^{\tau/(1-\tau)}}$$
converges, and since the intervals of type $g_1^mg_2^n(I)$ are two-by-two 
disjoint, this shows the finiteness of $\mathbb{E}(\ell_{\tau})$. $\hfill\square$

\vspace{0.2cm}

\begin{small} \begin{ejer} Given $\alpha > 0$, we say that 
a function $\psi$ has bounded {\bf {\em $\alpha$-variation}} if 
\index{variation of a function!$\alpha$-variation}
$$\sup_{a_0<a_1<\ldots<a_n} \sum_{i=1}^n \big| \psi(a_i) - \psi(a_{i-1})
\big|^{\alpha} < \infty.$$
Suppose that $g_1$ and $g_2$ are commuting circle diffeomorphisms 
which generate a rank-$2$ Abelian group acting freely and non minimally. Let 
$x,y$ be points in a connected component of the complement of the minimal invariant 
Cantor set. Prove that if $g_1'$ and $g_2'$ have bounded $\alpha$-variation for 
some $\alpha < 2$, then the function $V_{x,y}: \Omega \rightarrow \mathbb{R}$ 
defined by
$$V_{x,y}(\omega) \esp = \esp \sum_{k \geq 0} \big| g_{i_{k+1}}'(h_k(\omega)(x)) -
g_{i_{k+1}}'(h_k(\omega)(y)) \big|$$
has finite expectation.
\label{cuadr}
\end{ejer} \end{small}

From Lemma \ref{2d} one deduces that, if $S > 0$ is big enough, then the probability of the set 
$\Omega(S) \!=\! \big\{ \omega \in \Omega \!: \ell_{\tau}(\omega) \leq S \big\}$ is positive (actually, 
$\mathbb{P}[\Omega(S)]$ converges to $1$ as $S$ goes to infinity). Fix such an $S \! > \! 0$, and 
let \thinspace $\ell \!=\! \ell(\tau,S,|I|;\{g_1,g_2\})$ \thinspace be the constant of Lemma \ref{?}. 
Let us finally consider the open interval \esp $L'$ \esp of length $|L'| = \ell$ which is next to 
$I$ by the right. We claim that
\begin{equation}
\mathbb{P} \big[ \omega \in \Omega\!: h_n(\omega)(I) \not\subset L'
\mbox{ for all } n \in \mathbb{N} \big] \esp = \esp 0.
\label{total}
\end{equation}
To show this, recall that the action of $\Gamma$ is 
semiconjugate, but non conjugate, to an action by rotations. 
\index{rotation!group} 
Therefore, if we ``collapse'' each connected 
component of the complement of the minimal invariant Cantor set $\Lambda$, then we obtain a 
topological circle $\clo_{\Lambda}$ on which $g_1$ and $g_2$ induce minimal homeomorphisms. 
On the other hand, the interval $L'$ becomes an interval $U$ of non-empty interior in
$\clo_{\Lambda}$. Since the rotation numbers of $g_1$ and $g_2$ are irrational, there must 
exist $N \in \mathbb{N}$ so that, after collapsing, $g_1^{-1}(U), \ldots, g_1^{-N}(U)$ 
cover $\clo_{\Lambda}$, and the same holds for $g_2^{-1}(U),\ldots,g_2^{-N}(U)$. On the 
original circle $\clo$ this implies that, for every connected component $I_0$ of 
$\clo \setminus \Lambda$, there exist $n_1$ and $n_2$ in $\{1,\ldots,N\}$ so that 
\thinspace $g_1^{n_1}(I_0) \subset L'$ \thinspace and \thinspace $g_2^{n_2}(I_0) \subset L'$.

From the definition of $N$ it follows immediately that, for every 
integer $k \geq 0$, the conditional probability of the event 
$$\big[ g_1^i h_k(\omega)(I) \not\subset L' \thinspace \mbox { and }
\thinspace g_2^i h_k(\omega)(I) \not\subset L' 
\mbox{ for all } i \in \{1,\ldots,N \} \thinspace\big]$$ 
given that $\thinspace h_j(\omega)(I) \not\subset L' 
\mbox{ for all } j \!\in\! \{1,\ldots,k \}$ is zero. Remark now the following 
elementary property which follows directly from (\ref{paso}): the probabilities 
of ``jumping'' to the right (resp. upwards) of the Markov process are greater 
than or equal to $1/2$ under (resp. over) the diagonal (see Figure 21). 
Together with what precedes, this implies that
\begin{equation}
\mathbb{P} \big[ h_{k+i} (\omega)(I) \not\subset L', 
i \!\in \!\{1,\ldots,N \} \thinspace
\big| \thinspace h_j(\omega)(I) \not\subset L', 
j \!\in \!\{1,\ldots,k \} \big] \esp \leq \esp 1 - \frac{1}{2^N}.
\label{total2}
\end{equation}
As a consequence, for every $r \in \mathbb{N}$,
$$\mathbb{P} \big[ h_n(\omega)(I) \not\subset L', 
n \!\in\! \mathbb{N}  \big]
\leq \mathbb{P} \big[ h_n(\omega)(I) \not\subset L', 
i \!\in\! \{ 1, \ldots, rN \} \big] \esp 
\leq \esp \Big( 1 \!-\! \frac{1}{2^N} \Big)^r,$$
from where (\ref{total}) follows by letting $r$ 
go to the infinity.

\vspace{0.1cm}

To conclude the proof of Theorem \ref{denjoy-gen} (in the case $d=2$) remark that, if 
$\omega \!\in\! \Omega(S)$ and $n \!\in\! \mathbb{N}$ satisfy $h_n(\omega)(I) \!\subset\! L'$, 
then Lemma \ref{?} allows finding a hyperbolic fixed point for $h_n(\omega) \in \Gamma^+$, 
which contradicts the freeness of the action.

\vspace{0.21cm}

The proof of Theorem \ref{denjoy-gen} for $d > 2$ is analogous to that given for the case 
$d = 2$: assuming the existence of a wandering interval, one may consider the Markov process 
on $\mathbb{N}_0^d$ with transition probabilities 
$$p \big( (n_1,\ldots,n_i,\ldots,n_d) \longrightarrow
(n_1,\ldots,1+n_i,\ldots,n_d) \big) \esp = \esp \frac{1+n_i}{n_1 + \cdots + n_d + d}.$$
Once again, the arrival probabilities in $k$ steps are equally distributed over the points 
at simplicial distance $k$ from the origin. This allows controlling the 
distortion \index{control of distortion} of almost every random sequence 
(that is, an analogous of Lemma \ref{2d} holds for $\tau \!>\! 1/d$). Moreover, 
one easily checks that each $(n_1,\ldots,n_d) \in \mathbb{N}_0^d$ is the starting point of 
at least one half-line for which the transition probabilities between two adjacent vertexes 
is greater than or equal to $1/d$ (it suffices to follow the direction of the 
$i^{\mathrm{th}}$-coordinate for which $n_i$ attains its maximum value). This allows 
obtaining an inequality which is analogous to (\ref{total2}) (whose right-hand member 
will be equal to \thinspace $1\!-\!1/d^N$ \thinspace for some big integer $N$). Such an  
inequality implies property (\ref{total}), which thanks to the control of distortion 
allows using Lemma \ref{?} for finding elements with hyperbolic fixed points, thus 
contradicting the freeness of the action.

\vspace{0.15cm}

\begin{small} \begin{ejer} Show the following generalization of Theorem \ref{denjoy-gen}: 
if $\tau \!>\! 1/d$ for some $d \in \mathbb{N}$ and $\Gamma$ is a subgroup of 
$\mathrm{Diff}_+^{1+\tau}(\clo)$ which is semiconjugate to a group of rotations 
and contains $d$ elements with rotation numbers independent over the rationals, then 
the semiconjugacy is a conjugacy. In particular, if $f$ is a Denjoy counterexample 
of class $\ce^{1+\tau}$, then its centralizer in $\mathrm{Diff}_+^{1+\tau}(\clo)$ cannot 
contain elements $f_1,\ldots,f_{d-1}$ such that $\rho(f_1),\ldots,\rho(f_{d-1})$ 
and $\rho(f)$ are independent over the rationals.
\index{centralizer}

\noindent{\underbar{Hint.}} A careful reading of the proof of Theorem \ref{denjoy-gen} shows 
that the commutativity between the generators is needed only on the minimal invariant set.
\label{no-conm-1}
\end{ejer}
\end{small}

\vspace{0.05cm}

Similarly to the case of Denjoy's theorem, Kopell's lemma may be extended to the 
case of Abelian groups of interval diffeomorphisms with intermediate regularity.
\index{Kopell!generalized lemma}
The next result was also established in \cite{DKN}.

\vspace{0.1cm}

\begin{thm} {\em Let $d\!\geq\!2$ be an integer and $\varepsilon \!>\! 0$. 
Let $f_1,\ldots,f_{d+1}$ be $\ce^1$ commuting diffeomorphisms of $[0,1[$ 
for which there exist disjoint open intervals $I_{n_1,\ldots,n_{d}}$ 
placed on $[0,1]$ according to the lexicographic ordering and such 
that, for all $(n_1,\ldots,n_d) \! \in \! \mathbb{Z}^d$ and all 
$i\! \in\! \{1,\ldots,d\}$,
$$f_i(I_{n_1,\ldots,n_i,\ldots,n_{d}}) = I_{n_1,\ldots,n_i-1,\ldots,n_{d}},  
\qquad f_{d+1}(I_{n_1,\ldots,n_d}) = I_{n_1,\ldots,n_d}.$$ 
If $f_1,\ldots,f_d$ are of class $\ce^{1+1\!/\!d+\varepsilon}\!,$ then the action 
of $f_{d+1}$ on the union of the intervals $I_{n_1,\ldots,n_d}$ is trivial.}
\label{kop-gen}
\end{thm}

\noindent{\bf Proof.} We will only deal with the case $d = 2$, leaving to the reader the 
task of adapting the arguments to the general case. Let \esp $\tau = 1/2 + \varepsilon$, 
\esp let us identify with $\mathbb{N}_0 \times \mathbb{N}_0$ the semigroup $\Gamma^+$ 
generated by the elements $f_1$, $f_2$ in 
$\mathrm{Diff}_+^{1+\tau} \left( [0,1] \right)$, 
and let us consider again the Markov process given by (\ref{paso}). If 
we denote $I=I_{0,0} \!= ]a,b[$ and for each $\omega \in \Omega$ we let
$$\ell_{\tau}(\omega) = \sum_{i \geq 0} \big| h_i(\omega)(I) \big|^{\tau}\!,$$
then the argument of the proof of Lemma \ref{2d} allows showing that, since 
$\varepsilon > 0$, the function $\ell_{\tau}\!: \Omega  \rightarrow \mathbb{R}$ 
is almost surely finite.

Let us consider a $\tau$-H\"older constant $C$ for \esp $\log(f_1')$ \esp 
\index{H\"older!derivative}
and \esp $\log(f_2')$ \esp 
over $[0,b]$. For each $\omega = (f_{j_1},f_{j_2},\ldots) \in \Omega$, 
each $n,k$ in $\mathbb{N}$, and each $x \!\in\! I$, the equality \esp 
$f_3^n = h_k(\omega)^{-1} f_3^n h_k(\omega)$ \esp 
implies that $\log \big( (f_3^n)'(x)\big)$ coincides with 
$$\log \big( (f_3^n)' (h_k(\omega)(x)) \big) +  
\sum_{i=1}^k \big[ \log \big( f_{j_i}'(h_{i-1}(\omega)(x)) \big) 
- \log \big( f_{j_i}'(f_3^n \circ h_{i-1}(\omega)(x)) \big) \big],$$
and hence
\vspace{-0.2cm}
\begin{eqnarray*}
\big| \log ( (f_3^n)'(x) ) \big| 
&\leq& \big| \log \big( (f_3^n)' (h_k(\omega)(x)) \big) \big| \esp + \esp 
C \sum_{i=1}^k \big| h_{i-1}(\omega)(x) - f_3^n \circ h_{i-1}(\omega)(x)
\big|^{\tau} \\
&\leq& \big| \log \big( (f_3^n)' (h_k(\omega)(x)) \big) \big| 
\esp + \esp C \sum_{i=1}^k |h_{i-1}(\omega)(I)|^{\tau}.
\end{eqnarray*}
Choosing $\omega \in \Omega$ so that $\ell_{\tau}(\omega) = S$
is finite, the above inequality implies that 
$$\big| \log \big( (f_3^n)'(x)\big) \big| \esp \leq \esp 
\big| \log \big( (f_3^n)' (h_k(\omega)(x)) \big) \big| + CS.$$
Notice that the sequence $\big( h_k(\omega)(x) \big)$ converges to a (necessarily parabolic) 
fixed point of $f_3$ (actually, for almost every $\omega \in \Omega$ this sequence converges 
to the origin). Letting $k$ go to infinity we conclude that 
\thinspace $\big| \log \big( (f_3^n)'(x)\big) \big| \leq CS.$ \thinspace Therefore, 
$(f_3^n)'(x) \leq \exp(CS)$ for all $x \in I$ and all $n \in \mathbb{N}$, which implies 
that the restriction of $f_3$ to the interval $I$ is the identity. By the commutativity 
hypothesis, the same is true over all the intervals $I_{n_1,n_2}$, which concludes 
the proof. $\hfill\square$

\vspace{0.2cm}

\begin{small} \begin{ejer} Prove that Theorem \ref{kop-gen} still holds if 
$f_1,\ldots,f_d,f_{d+1}$ are $\ce^1$ commuting diffeomorphisms of the interval 
and the $(d\!-\!\varepsilon)$-variation of the derivatives of $f_1,\ldots,f_d$ 
are finite for some $\varepsilon > 0$ ({\em c.f.,} Exercise \ref{cuadr}).
\end{ejer}

\begin{ejer} One may give a version of Theorem \ref{kop-gen} for which the commutativity 
hypothesis between the $f_i$'s is weaker (compare Exercise \ref{no-conm-1}). More 
precisely, let $d\!\geq\!2$ be an integer, and let $f_1,\ldots,f_{d+1}$ be $\ce^1$ 
diffeomorphisms of the interval $[0,1[$ (which do not necessarily commute). Suppose that 
there exist disjoint open intervals $I_{n_1,\ldots,n_{d+1}}$ disposed on $]0,1[$ respecting the 
lexicographic ordering and so that, for every $(n_1,\ldots,n_{d+1}) \! \in \! \mathbb{Z}^{d+1}$,
$$f_i(I_{n_1,\ldots,n_i,\ldots,n_{d+1}}) = I_{n_1,\ldots,n_i-1,\ldots,n_{d+1}}
\quad \mbox{for all} \quad i\! \in\! \{1,\ldots,d+1\}.$$
Prove that $f_1,\ldots,f_d$ cannot be all simultaneously of 
class $\ce^{1+1\!/\!d+\varepsilon}$ for any $\varepsilon \!>\! 0$.

\noindent{\underbar{Hint.}} Use similar arguments to those of the proof of Theorem 
\ref{kop-gen}. After establishing a (generic) control for the distortion, apply 
the technique of the proof of Proposition \ref{kopelcito} to conclude.
\label{no-conm-2}
\end{ejer}
\end{small}

The combinatorial hypothesis of Theorem \ref{kop-gen} might seem strange. Nevertheless, 
according to \cite{ts-annals}, the actions of $\mathbb{Z}^{d+1}$ on $[0,1]$ 
satisfying this hypothesis are very interesting from a cohomological viewpoint. 
Notice, however, that it is not difficult to construct actions of 
$\mathbb{Z}^{d+1}$ by $\ce^{2-\varepsilon}$ diffeomorphisms of the 
interval without global fixed points in the interior. To do this, let us 
consider for instance $d$ diffeomorphisms $f_2,\ldots,f_{d+1}$ of an interval 
$[a,b]\! \subset ]0,1[$ whose open supports are disjoint, 
\index{support!of a map} 
and let $f_1$ be a diffeomorphism 
of $[0,1]$ without fixed point in the interior and sending $]a,b[$ into a disjoint interval. 
Extending $f_2,\ldots,f_{d+1}$ to the whole interval $[0,1]$ so that they commute with 
$f_1$, one obtains a faithful action of $\mathbb{Z}^{d+1}$ by homeomorphisms of the 
interval and without global fixed point in $]0,1[$. Clearly, the methods from 
\S \ref{ejemplosceuno} allow smoothing this action up to class 
$\ce^{2-\varepsilon}$ for every $\varepsilon > 0$.

It is also interesting to remark the existence of actions by diffeomorphisms of the interval 
which at the interior are free but non conjugate to actions by translations. According to 
Proposition \ref{no-semiconjugado}, such a phenomenon cannot appear in class $\ce^{1+\mathrm{bv}}$. 
However, using the methods from \S \ref{ejemplosceuno}, it is possible to construct examples of 
$\mathbb{Z}^2$-actions by $\ce^{3/2 - \varepsilon}$ diffeomorphisms of $[0,1]$ which are 
free on $]0,1[$ but do admit wandering intervals. Once again, the regularity \thinspace 
$\ce^{3/2- \varepsilon}$ \thinspace is optimal for the existence of these examples. 
The next result from \cite{DKN} should be compared with Proposition 
\ref{no-semiconjugado} and Corollary \ref{flujocelog}.
\index{wandering interval}

\vspace{0,1cm}

\begin{thm} {\em Let $\Gamma$ be a subgroup of $\mathrm{Diff}_+^{1+\tau}([0,1[)$ 
isomorphic to $\mathbb{Z}^d$, with $\tau > 1/d$ and $d \geq 2$. If the restriction 
to $]0,1[$ of the action of \esp $\Gamma$ is free, then it is minimal and 
topologically conjugate to the action of a group of translations.}
\label{ult-gen}
\end{thm}

\noindent{\bf Proof.} Once again, we will give the complete proof just for the case 
$d=2$. Let $f_1$ and $f_2$ be the generators of a group $\Gamma \!\sim\! \mathbb{Z}^2$ of 
$\ce^{1+\tau}$ diffeomorphisms of $[0,1[$ acting freely on $]0,1[$. Changing some of these 
elements by their inverses if necessary, we may assume that they topologically contract 
towards the origin. Suppose that the action of $\Gamma$ on $]0,1[$ is not conjugate 
to an action by translations. In this case, a simple control of distortion argument of 
hyperbolic type shows that all of the elements in $\Gamma$ are tangent to the identity 
at the origin. Indeed, suppose for a contradiction that there exist wandering intervals and 
an element $f \in \Gamma$ such that $f'(0) < 1$. Fix $\lambda < 1$ and $c \! \in ]0,1[$ 
so that $f'(x) \leq \lambda$ for all $x \!\in \![0,c[$, and fix a maximal open wandering 
interval $I\!\!=]a,b[$ contained in $]0,c[$. If we denote by $L$ the interval $[f(b),b]$, 
then
$$\sum_{n \geq 0} \big| f^n(L) \big|^{\tau} 
\leq |L|^{\tau} \sum_{n \geq 0} \lambda^{n \tau}
= \frac{|L|^{\tau}}{1-\lambda^{\tau}} = \bar{S}.$$
Consequently, \thinspace $(f^n)'(x) / (f^n)'(y) \leq \exp(C \bar{S})$
\thinspace for all $x,y$ in $L$, where $C>0$ is a $\tau$-H\"older constant for 
$\log(f')$ over $[0,c]$. This estimate allows applying the arguments of the 
proof of Proposition \ref{no-semiconjugado}, thus obtaining a contradiction.\footnote{Remark 
that this argument just uses the hypothesis $\tau > 0$. This is related to the fact that 
Sternberg's linearization theorem still holds in class $\ce^{1+\tau}$ 
for any positive $ \tau$ (see Remark \ref{st-obs}). 
Indeed, since the centralizer of a nontrivial linear germ 
coincides with the 1-parameter group of 
linear germs, this prevents the existence of wandering intervals for the dynamics 
we are dealing with.}
\index{germ of a diffeomorphism}

Let us now identify the semigroup $\Gamma^+$ generated by $f_1$ and $f_2$ with 
$\mathbb{N}_0 \times \mathbb{N}_0$, and let us consider once again our Markov 
process on it. If $\tau > 1/2$, then the proof of Lemma \ref{2d} shows the 
finiteness of the expectation of the function 
$$\omega \mapsto \ell_{\tau}(\omega) 
= \sum_{k \geq 0} \big| h_k(\omega)(I) \big|^{\tau}.$$
Choose a large enough $S>0$ so that the set \thinspace 
$\Omega(S) = \{\omega \in \Omega\!: \ell_{\tau}(\omega) \leq S \}$ \thinspace has 
positive probability, and let \thinspace $\bar{\ell} = |I| / \exp(2^{\tau}CS)$. 
\thinspace The first part of the proof of Lemma \ref{?} shows that if $I''$ denotes 
the interval which is next to $I$ by the left and has length $\bar{\ell}$, then for 
every $x,y$ in $J = \bar{I}'' \cup \bar{I}$, every $\omega \in \Omega(S)$, and 
every $n \in \mathbb{N}$,
\begin{equation}
\frac{h_n(\omega)'(x)}{h_n(\omega)'(y)} \leq \exp(2^{\tau} C S).
\label{controlado}
\end{equation}
Since the interval $I$ is not strictly contained in any other open wandering 
interval, there must exist $h \in \Gamma$ such that $h(I) \subset I''$ (and 
hence \thinspace $|h(I)| < |I| \thinspace \exp(- 2^{\tau} C S)$). \thinspace
Fix an arbitrary point $y \!\in\! I$. Since $h'(0)=1$ and $h_n(\omega)(y)$ 
converges to the origin for all $\omega \in \Omega(S)$, from the equality 
$$h'(y) = \frac{h_n(\omega)'(y)}{h_n(\omega)'(h(y))} \thinspace \esp 
h' \big( h_n(\omega)(y) \big)$$
and (\ref{controlado}) one deduces that \thinspace
$h'(y) \!\geq\! \exp(-2^{\tau} C S).$ \thinspace
By integrating this inequality we obtain 
$|h(I)| \geq |I| \thinspace \exp(-2^{\tau} C S)$,
which is in contradiction with the choice of $h$. $\hfill\square$

\vspace{0.1cm}

\begin{small} \begin{obs} Let $f_1,\ldots,f_d$ be the generators of a group 
$\Gamma \!\sim\! \mathbb{Z}^d$ acting freely by homeomorphisms of $]0,1[$ so 
that the action is not conjugate to an action by translations. By identifying 
the points in the orbits of $f_1$, the maps $f_2,\ldots,f_d$ become the 
generators of a rank--$(d\!-\!1)$ Abelian group of circle homeomorphisms 
which is semiconjugate, but non conjugate, to a group of rotations: Theorem 
\ref{denjoy-gen} then implies that the $f_i$'s cannot be all of class 
$\ce^{1 \! / \! (d-1) + \varepsilon}$. Notice that this argument uses only 
the regularity of the $f_i$'s at the interior; in this context, the obstruction 
appears in class $\ce^{1 \! / \! (d-1)}$. However, according to Theorem \ref{ult-gen}, 
for interval diffeomorphisms of $[0,1[$ the obstruction already appears in class 
$\ce^{1 \! / \! d}$. Obviously, the difference lies in the differentiability 
of the maps at the origin. This actually plays an important role along the proof.
\end{obs}

\begin{ejer} Extend Theorem \ref{ult-gen} to subgroups of $\mathrm{Diff}_+^{1+\tau}([0,1[)$
that are semiconjugate (on $]0,1[$) to groups of translations but do not act freely 
on the interior (compare Exercises \ref{no-conm-1} and \ref{no-conm-2}).
\end{ejer}
\end{small}
\index{Denjoy!generalized theorem|)}


\section{Nilpotent Groups of Diffeomorphisms}

\subsection{Plante-Thurston's theorems}

\hspace{0.45cm} Given a nilpotent group \index{group!nilpotent} $\Gamma$, let us denote 
$$\{ id \} =\Gamma_{k}^{\mathrm{nil}} \sgn \Gamma_{k-1}^{\mathrm{nil}} \sgn \ldots \sgn
\Gamma_{1}^{\mathrm{nil}} \sgn \Gamma_0^{\mathrm{nil}} =\Gamma$$
the {\bf {\em central series}} of $\Gamma$, that is, 
\index{central series} 
$\Gamma_{i+1}^{\mathrm{nil}} = [\Gamma,\Gamma_i^{\mathrm{nil}}]$ 
and $\Gamma_{k-1}^{\mathrm{nil}} \neq \{ id \}$. Remark 
that the subgroup $\Gamma_{k-1}^{\mathrm{nil}}$ 
is contained in the center of $\Gamma$. 
\index{center of a group}
The next result corresponds to a (weak version of a) theorem 
\index{Plante-Thurston's theorem|(}
due to Plante and Thurston \cite{Pl3,PT}.

\vspace{0.07cm}

\begin{thm} {\em Every nilpotent subgroup of 
$\mathrm{Diff}_{+}^{1+\mathrm{bv}} ([0,1[)$ is Abelian.}
\label{chancho}
\end{thm}

\noindent{\bf Proof.} Without loss of generality, we may suppose that $]0,1[$ is an irreducible 
component for the action of $\Gamma$. We will show that the restriction of $\Gamma$ to $]0,1[$ 
is free, which allows concluding the proof using H\"older's theorem. 
\index{H\"older!theorem}\index{action!free}Suppose for a contradiction that there 
exists a nontrivial $f \in \Gamma$ whose set of fixed points in $]0,1[$ is non-empty, 
and let $x_0 \! \in ]0,1[$ be a point in the boundary of this set. Fix a nontrivial 
element $g$ in the center of $\Gamma$. We claim that $g(x_0) = x_0$. For otherwise, 
replacing $g$ by its inverse if necessary, we may assume that $g(x_0) < x_0$. Let 
$$a' = \lim\limits_{n \rightarrow +\infty} g^n(x_0), \qquad
b' = \lim\limits_{n \rightarrow -\infty} g^n(x_0).$$
Notice that $[a',b'[ \subset [0,1[$. Moreover, the restriction of $g$ to $]a',b'[$ has no 
fixed point. Since $f$ and $g$ commute, each $g^n(x_0)$ is a fixed point of $f$, and hence 
$f$ fixes $a'$ and $b'$. We then obtain a contradiction by applying Kopell's lemma to the 
restrictions of $f$ and $g$ to the interval $[a',b'[$. Therefore, $g(x_0) = x_0$. 
\index{Kopell!lemma}

Since $g$ was nontrivial and $g(x_0)=x_0$, the boundary in $]0,1[$ of the set of fixed 
points of $g$ is non-empty. Fix a point $x_1 \!\in \esp\!]0,1[$ in this boundary, and let $h$ 
be an arbitrary element in $\Gamma$. Since $g(x_1) = x_1$ and $gh=hg$, the same argument 
as before shows that $h(x_1) \!=\! x_1$. Therefore, the intervals $[a',x_1[$ and $[x_1,b'[$ 
are invariant by $\Gamma$, and this contradicts the fact that $]0,1[$ was an 
irreducible component. $\hfill\square$
\index{irreducible component}

\vspace{0.15cm}

\begin{small} \begin{ejer} Prove that Theorem \ref{chancho} still holds for 
nilpotent groups of germs of $\ce^{1+\mathrm{bv}}$ diffeomorphisms of 
the line fixing the origin.
\label{PT-ger}\label{germ of a diffeomorphism}
\end{ejer} \end{small}

Let us now consider the case of the circle. Since nilpotent groups are amenable, if $\Gamma$ 
is a nilpotent subgroup of $\mathrm{Diff}_+^{1+\mathrm{bv}}(\clo)$ then it must preserve a 
probability measure $\mu$ on $\clo$. The rotation number function $g \mapsto \rho(g)$ 
is then a group homomorphism from $\Gamma$ into $\mathbb{T}^1$ (see \S \ref{medinv}). 
If the rotation number 
\index{rotation number}
of an element $g \in \Gamma$ is irrational, then $\mu$ is conjugate 
to the Lebesgue measure. The group $\Gamma$ is therefore conjugate 
to a group of rotations; 
\index{rotation!group} 
in particular, $\Gamma$ is Abelian.

Suppose now that the rotation number of every element in $\Gamma$ is rational. In this case, 
the support of $\mu$ is contained in the intersection of the set of periodic points of these 
elements. If $\Gamma$ is not Abelian then we can take $f,g$ in $\Gamma_{k-2}^{\mathrm{nil}}$ 
so that $h = [f,g] \in \Gamma_{k-1}^{\mathrm{nil}}$ is nontrivial. Notice that $h$ is 
contained in the center of $\Gamma$; moreover, $\rho(h) = 0$, and hence $h$ has fixed 
points. From the equality $f^{-1}g^{-1}fg = h$ one obtains $g^{-1} f g = f h$.  
Hence, $g^{-1} f^{n} g = f^n h^n$, and through successive conjugacies by $g$ one  
concludes that $g^{-m}f^{n}g^{m} = f^{n}h^{mn}$ for all 
$m \in \mathbb{N}$. It follows that
\begin{equation}
h^{mn} \esp = \esp f^{-n}g^{-m}f^{n}g^{m}.
\label{conmutadortri}
\end{equation}
If $x_0$ belongs to the support of $\mu$ and $m,n$ are positive integers so that $x_0$ 
is a fixed point for both $f^n$ and $g^m$, then (\ref{conmutadortri}) shows that 
$h(x_0) = x_0$. Let us consider 
the restriction to $[x_0,x_0 + 1[$ of the group generated by $f^n$, $g^m$, and $h$. 
Since this group is nilpotent, this restriction must be Abelian, and hence from 
(\ref{conmutadortri}) one concludes that $h^{mn}$ is the identity over 
$[x_0,x_0+1[$. Therefore, $h$ itself is the identity, which is a 
contradiction. We have then proved the following result.

\vspace{0.1cm}

\begin{thm} {\em Every nilpotent subgroup of 
$\mathrm{Diff}_+^{1+\mathrm{bv}}(\clo)$ is Abelian.}
\label{casi-ultimo}
\end{thm}

Finally, in the case of the line we have the following result (recall that a group $\Gamma$ 
is {\bf {\em metabelian}} if its first derived subgroup $\Gamma' = [\Gamma,\Gamma]$ is Abelian).
\index{group!metabelian}

\vspace{0.13cm}

\begin{thm} {\em Every nilpotent subgroup of 
$\mathrm{Diff}_{+}^{1+\mathrm{bv}}(\mathbb{R})$ is metabelian.}
\label{ultimo}
\end{thm}

\vspace{0.13cm}

This result is a direct consequence of Corollary 
\ref{conmutptofijo} and the next proposition.

\vspace{0.13cm}

\begin{prop} {\em Let $\Gamma$ be a nilpotent subgroup of $\mathrm{Diff}_{+}^{1+\mathrm{bv}}(\mathbb{R})$. 
If every element in $\Gamma$ fixes at least one point in the line, then $\Gamma$ is Abelian.}
\end{prop}

\noindent{\bf Proof.} Fix a nontrivial element $g$ in the center of $\Gamma$, and let $b$ be a point 
in the boundary of the set of fixed points of $g$. We claim that $b$ is fixed by every element in 
$\Gamma$. Indeed, suppose for a contradiction that  $f \!\in\! \Gamma$ is such that $f(b) \neq b$. 
Replacing $f$ by $f^{-1}$ if necessary, we may suppose that $f(b) \!<\! b$. 
By hypothesis, at least one of the 
sequences $\big( f^n(b) \big)$ or $\big( f^{-n}(b) \big)$ converges to a fixed point 
$a \!\in\! \mathbb{R}$ of $f$. Since both cases are analogous, let us only deal with the 
first one. Notice that $f^n(b)$ is a fixed point of $g$ for every $n \!\in\! \mathbb{N}$, 
and hence $g(a) = a$. On the other hand, $f$ has no fixed point in $]a,b[$. Therefore, 
letting $a' = \lim_{n \to \infty} f^{-n}(b)$ we obtain a contradiction by applying 
Kopell's lemma to the restrictions of $f$ and $g$ to the interval $[a,a'[$.\\

The above discussion implies that every element $f \!\in\! \Gamma$ fixes 
the intervals $]\!-\!\infty,b]$ and $[b,\infty[$. The conclusion of the 
proposition then follows from Theorem \ref{chancho}. $\hfill\square$

\vspace{0.4cm}

It is important to point out that Theorem \ref{ultimo} does not mean that the nilpotence 
degree of every nilpotent subgroup of $\mathrm{Diff}^{1+\mathrm{bv}}_+(\mathbb{R})$ is less 
than or equal to 2. Actually, $\mathrm{Diff}^{\infty}_+(\mathbb{R})$ contains nilpotent 
subgroups of arbitrary degree of nilpotence. The following example, taken from \cite{FF}, 
uses ideas which are very similar to those of \cite{Pl2}.
\index{degree!of nilpotence}

\begin{small} \begin{ejem}
Let us consider a $\ce^{\infty}$ diffeomorphism $g\!: [0,1] \rightarrow [0,1]$ 
that is infinitely tangent to the identity at the endpoints (that is, 
$g'(0) \!=\! g'(1) \!=\! 1$, and all the higher derivatives at 
$0$ and $1$ are zero). Let $f$ be the translation $f(x)=x-1$ on the line. 
For each integer $n \!\geq\! 0$ let $k(n,m)=\esp$${m+n-1}\choose{n}$, and let  
$h_n$ be the diffeomorphism of the line defined by
$$h_n(x) = g^{k(n,m)}(x-m) + m \qquad \mbox{for } \esp x \in [m,m+1[.$$
We leave to the reader the task of showing that the maps $h_n$ commute between 
them, and that $[f,h_n] = f^{-1}h_n^{-1}fh_{n} = h_{n-1}$ for every $n \geq 1$, 
while $[f,h_0] = Id$. One then easily concludes that the group $\Gamma_n$ 
generated by $f$ and $h_0, \ldots, h_n$ is nilpotent with degree of 
nilpotence $n\!+\!1$.
\label{nilp-locos}
\end{ejem} 

\begin{ejer}
Using Theorem \ref{ultimo}, show directly the following weak version of Theorem 
\ref{wittefuerte}: \index{Witte-Morris!theorem of non orderability of lattices}
for $n \geq 4$, every action of a finite index subgroup of 
$\mathrm{SL}(n,\mathbb{Z})$ by $\ce^{1+\mathrm{bv}}$ diffeomorphisms of the 
line is trivial.
\end{ejer} \end{small}
\index{Plante-Thurston's theorem|)}
\index{Plante!Plante-Thurston's theorem|see{Plante-Thurston's theorem}}


\subsection{On the growth of groups of diffeomorphisms}
\label{seccion-crecer} \index{growth!of groups|(}

\hspace{0.45cm} Another viewpoint of Plante-Thurston's 
\index{Plante-Thurston's theorem} theorem is related to the notion of 
{\em growth} of groups. Recall that for a group provided with a finite system of generators, 
the {\bf {\em growth function}} is the one that assigns, to each positive integer $n$, the 
number of elements that may be written as a product of no more than $n$ generators and 
their inverses. One says that the group has {\bf {\em polynomial growth}} or 
{\bf {\em exponential growth}}, if 
its growth function has the corresponding asymptotic behavior. These notions 
do not depend on the choice of the system of generators (see Appendix B). 

\begin{small} \begin{ejer} Prove that if a group contains a free 
semigroup on two generators, then it has exponential growth.
\index{semigroup!free|(}

\noindent{\underbar{Remark.}} There exist groups of exponential growth 
without free semigroups on two generators: see for instance \cite{olsha}.
\label{chanta-chanta}
\end{ejer}

\begin{ejem}
The construction of groups with {\bf \em intermediate growth} 
(that is, neither polynomial nor exponential growth) is a nontrivial and 
fascinating topic. The first examples were given by Grigorchuk in \cite{grigor-deg}. 
One of them, namely the so-called {\em first Grigorchuk group} $\hat{\mathrm{H}}$, 
may be seen either as a group acting on the binary rooted 
\index{tree}\index{Grigorchuk's examples}tree $\mathcal{T}_2$, or 
as a group acting isometrically on the Cantor set $\{0,1\}^{\mathbb{N}}$. (These 
points of view are essentially the same, since the boundary at infinity of 
$\mathcal{T}_2$ naturally identifies with $\{0,1\}^{\mathbb{N}}$.) Here we 
record the explicit definition. For more details, we strongly 
recommend the lecture of the final chapter of \cite{harpe}.

Using the convention $(x_1,(x_2,x_3,\ldots)) = (x_1,x_2,x_3,\ldots)$ for $x_i \! \in \! \{0,1\}$, 
the generators of $\hat{\mathrm{H}}$ are the elements $\hat{a},\hat{b},\hat{c}$, and $\hat{d}$, 
whose actions on sequences in $\{0,1\}^{\mathbb{N}}$ are recursively defined by \esp 
$\hat{a}(x_1,x_2,x_3,\ldots) = (1-x_1,x_2,x_3,\ldots)$ \esp and 
$$\hat{b}(x_1,x_2,x_3,\ldots) = \left \{ \begin{array} {l}
(x_1, \hat{a}(x_2,x_3, \ldots )), \hspace{0.1cm} x_1 = 0,\\
(x_1, \hat{c}(x_2,x_3,\ldots)), \hspace{0.11cm} x_1 = 1, \end{array} \right.$$
$$\hat{c}(x_1,x_2,x_3,\ldots) = \left \{ \begin{array} {l}
(x_1, \hat{a}(x_2,x_3, \ldots )), \hspace{0.1cm} x_1 = 0,\\
(x_1, \hat{d}(x_2,x_3,\ldots)), \hspace{0.1cm} x_1 = 1, \end{array} \right.$$
$$\hat{d}(x_1,x_2,x_3,\ldots) = \left \{ \begin{array} {l}
(x_1, x_2,x_3, \ldots ), \hspace{0.52cm} x_1 = 0,\\
(x_1, \hat{b}(x_2,x_3,\ldots)), \hspace{0.13cm} x_1 = 1. \end{array} \right.$$
In this way, the action on $\mathcal{T}_2$ of the element $\hat{a} \!\in\! \hat{\mathrm{H}}$ 
consists in permuting the first two edges (together with the trees which are rooted at 
the second level vertexes). The elements $\hat{b}$, $\hat{c}$, and $\hat{d}$, fix these 
edges, and their actions on the higher levels of $\mathcal{T}_2$ are illustrated below.
\label{grigor-ejem}
\end{ejem} \end{small}

\vspace{0.02cm}


\beginpicture

\setcoordinatesystem units <0.873cm,0.873cm>


\plot
0 0
0.7 -0.8 /
\plot
0 0
-0.7 -0.8 /
\plot
0.7 -0.8
1 -1.5 /
\plot
0.7 -0.8
0.4 -1.5 /
\plot
-0.7 -0.8
-0.4 -1.5 /
\plot
-0.7 -0.8
-1 -1.5 /
\plot
1 -1.5
1.15 -1.9 /
\plot
1 -1.5
0.85 -1.9 /
\plot
0.4 -1.5
0.55 -1.9 /
\plot
0.4 -1.5
0.25 -1.9 /
\plot
-0.4 -1.5
-0.25 -1.9 /
\plot
-0.4 -1.5
-0.55 -1.9 /
\plot
-1 -1.5
-1.15 -1.9 /
\plot
-1 -1.5
-0.85 -1.9 /
\put{........................} at 0 -2.3


\plot
3.5 0
4.2 -0.8 /
\plot
3.5 0
2.8 -0.8 /
\plot
4.2 -0.8
4.5 -1.5 /
\plot
4.2 -0.8
3.9 -1.5 /
\plot
2.8 -0.8
3.1 -1.5 /
\plot
2.8 -0.8
2.5 -1.5 /
\plot
4.5 -1.5
4.65 -1.9 /
\plot
4.5 -1.5
4.35 -1.9 /
\plot
3.9 -1.5
4.05 -1.9 /
\plot
3.9 -1.5
3.75 -1.9 /
\plot
3.1 -1.5
3.25 -1.9 /
\plot
3.1 -1.5
2.95 -1.9 /
\plot
2.5 -1.5
2.35 -1.9 /
\plot
2.5 -1.5
2.65 -1.9 /
\put{........................} at 3.5 -2.3


\plot
0 0
0.7 -0.8 /
\plot
0 0
-0.7 -0.8 /
\plot
0.7 -0.8
1 -1.5 /
\plot
0.7 -0.8
0.4 -1.5 /
\plot
-0.7 -0.8
-0.4 -1.5 /
\plot
-0.7 -0.8
-1 -1.5 /
\plot
1 -1.5
1.15 -1.9 /
\plot
1 -1.5
0.85 -1.9 /
\plot
0.4 -1.5
0.55 -1.9 /
\plot
0.4 -1.5
0.25 -1.9 /
\plot
-0.4 -1.5
-0.25 -1.9 /
\plot
-0.4 -1.5
-0.55 -1.9 /
\plot
-1 -1.5
-1.15 -1.9 /
\plot
-1 -1.5
-0.85 -1.9 /
\put{........................} at 0 -2.3


\plot
7 0
7.7 -0.8 /
\plot
7 0
6.3 -0.8 /
\plot
7.7 -0.8
8 -1.5 /
\plot
7.7 -0.8
7.4 -1.5 /
\plot
6.3 -0.8
6.6 -1.5 /
\plot
6.3 -0.8
6 -1.5 /
\plot
8 -1.5
8.15 -1.9 /
\plot
8 -1.5
7.85 -1.9 /
\plot
7.4 -1.5
7.55 -1.9 /
\plot
7.4 -1.5
7.25 -1.9 /
\plot
6.6 -1.5
6.45 -1.9 /
\plot
6.6 -1.5
6.75 -1.9 /
\plot
6 -1.5
5.85 -1.9 /
\plot
6 -1.5
6.15 -1.9 /
\put{........................} at 7 -2.3

\plot
10.5 0
11.2 -0.8 /
\plot
10.5 0
9.8 -0.8 /
\plot
11.2 -0.8
11.5 -1.5 /
\plot
11.2 -0.8
10.9 -1.5 /
\plot
9.8 -0.8
10.1 -1.5 /
\plot
9.8 -0.8
9.5 -1.5 /
\plot
11.5 -1.5
11.65 -1.9 /
\plot
11.5 -1.5
11.35 -1.9 /
\plot
10.9 -1.5
11.05 -1.9 /
\plot
10.9 -1.5
10.75 -1.9 /
\plot
10.1 -1.5
10.25 -1.9 /
\plot
10.1 -1.5
9.95 -1.9 /
\plot
9.5 -1.5
9.35 -1.9 /
\plot
9.5 -1.5
9.65 -1.9 /
\put{........................} at 10.5 -2.3


\put{$\hat{a}$} at 0 0.6
\put{$\curvearrowright$} at 0 0.2
\put{$\longleftrightarrow$} at 0 -0.7

\put{$\hat{b}$} at 3.5 0.6
\put{$\curvearrowright$} at 3.5 0.2
\put{$\hat{a}$} at 2.8 -0.2
\put{$\hat{c}$} at 4.3 -0.2
\put{$\curvearrowright$} at 2.8 -0.5
\put{$\curvearrowright$} at 4.3 -0.5

\put{$\hat{c}$} at 7 0.6
\put{$\curvearrowright$} at 7 0.2
\put{$\hat{a}$} at 6.3 -0.2
\put{$\hat{d}$} at 7.8 -0.17
\put{$\curvearrowright$} at 6.3 -0.5
\put{$\curvearrowright$} at 7.8 -0.5

\put{$\hat{d}$} at 10.5 0.6
\put{$\curvearrowright$} at 10.5 0.2
\put{$id$} at 9.8 -0.2
\put{$\hat{b}$} at 11.3 -0.17
\put{$\curvearrowright$} at 9.8 -0.5
\put{$\curvearrowright$} at 11.3 -0.5

\put{Figure 22} at 5.2 -2.8

\small
\put{$\bullet$} at 0 0
\put{$\bullet$} at 3.5 0
\put{$\bullet$} at 7 0
\put{$\bullet$} at 10.5 0
\put{$\bullet$} at 0.7 -0.8
\put{$\bullet$} at 4.2 -0.8
\put{$\bullet$} at 7.7 -0.8
\put{$\bullet$} at 11.2 -0.8
\put{$\bullet$} at -0.7 -0.8
\put{$\bullet$} at 2.8 -0.8
\put{$\bullet$} at 6.3 -0.8
\put{$\bullet$} at 9.8 -0.8

\put{$\bullet$} at 1 -1.5
\put{$\bullet$} at 4.5 -1.5
\put{$\bullet$} at 1 -1.5
\put{$\bullet$} at 11.5 -1.5
\put{$\bullet$} at 10.9 -1.5
\put{$\bullet$} at -0.4 -1.5
\put{$\bullet$} at 3.1 -1.5
\put{$\bullet$} at 6.6 -1.5
\put{$\bullet$} at 9.5 -1.5

\put{$\bullet$} at 8 -1.5
\put{$\bullet$} at 0.4 -1.5
\put{$\bullet$} at 3.9 -1.5
\put{$\bullet$} at 7.4 -1.5
\put{$\bullet$} at 10.1 -1.5
\put{$\bullet$} at -1 -1.5
\put{$\bullet$} at 2.5 -1.5
\put{$\bullet$} at 6 -1.5
\put{$\bullet$} at 9.5 -1.5

\put{} at -2.5 0

\endpicture


\vspace{0.3cm}

A celebrated theorem of Gromov \index{Gromov!theorem on polynomial growth groups}
establishes that a group has polynomial growth if and only if 
it is almost nilpotent, {\em i.e.}, if it contains a finite index nilpotent 
\index{group!nilpotent} subgroup. Due to this 
and Exercise \ref{virt-abel}, Plante-Thurston's theorem 
\index{Plante-Thurston's theorem} 
may be reformulated by saying that every 
subgroup of polynomial growth of $\mathrm{Diff}_+^{1+\mathrm{bv}}([0,1[)$ is Abelian. Actually, 
this was the original statement of the theorem, which is prior to Gromov's theorem (see 
Exercise \ref{original-PL}). Following \cite{ceuno}, in what follows we will generalize 
this statement to groups of sub-exponential growth, and even more generally to groups 
without free semigroups on two generators.

\vspace{0.1cm}

\begin{thm} {\em Every finitely generated subgroup of $\mathrm{Diff}_+^{1+\mathrm{bv}}([0,1[)$ 
without free semigroups on two generators is Abelian.}
\label{mio}
\end{thm}

\vspace{0.01cm}

\begin{small} \begin{ejem} The {\bf {\em wreath product}} 
\index{wreath product} \esp 
$\Gamma = \mathbb{Z} \wr \mathbb{Z} = \mathbb{Z} \ltimes \oplus_{\mathbb{Z}} \mathbb{Z}$
\esp naturally acts on the interval: it suffices to identify the generator of the 
$\mathbb{Z}$-factor 
in $\Gamma$ with a homeomorphism $f$ of $[0,1]$ satisfying $f(x) \!<\! x$ for all 
$x \! \in \esp \!]0,1[$, and the element $(\ldots,0,1,0,\ldots)$ of the second factor in 
$\Gamma$ to a homeomorphism $g$ satisfying $g(x) \neq x$ for all $x\!\in ]f(a),a[$ and 
$g(x)=x$ for all $x \in [0,1] \setminus [f(a),a]$, where $a$ is some point in $]0,1[$. 
This action may be smoothed up to class $\ce^{\infty}$ (see Example \ref{one} 
for more details on this construction). Notice that $\Gamma$ is a metabelian, 
\index{group!metabelian}
non virtually 
\index{Milnor-Wolf's theorem}
nilpotent group. According to a classical theorem by Rosenblatt \cite{ros} (which generalizes 
prior results of Milnor \cite{Mil} and Wolf \cite{Wo}; see also \cite{br-ros}), every 
solvable \index{group!solvable}
group which is non virtually nilpotent contains free semigroups in two generators. 
Therefore, $\Gamma$ contains such a semigroup. Actually, looking at the action on 
the line, one may easily verify that the semigroup generated by $f$ and $g$ is free.
\label{wreath-exp}
\end{ejem} \end{small}

To show Theorem \ref{mio} we will need a version of Kopell's lemma which does 
not use the commutativity hypothesis too strongly. This version will be 
very useful in \S \ref{solvito} for the study of solvable groups of 
diffeomorphisms of the interval. The proof below was taken from \cite{CC1}. 
\index{Kopell!lemma}

\vspace{0.1cm}

\begin{lem} {\em Let $h_1$ and $h_2$ be two diffeomorphisms from the interval $[0,1[$ into 
their images that fix the origin. Suppose that $h_1$ is of class $\mathrm{C}^{1+\mathrm{bv}}$ 
and $h_2$ of class $\mathrm{C}^1$. Suppose moreover that $h_1(x) < x$ for all $x \! \in ]0,1[$, 
that $h_2(x_0) = x_0$ for some $x_0 \! \in ]0,1[$, and that for each $n \!\in\! \mathbb{N}$ 
the point $x_n = h_1^n(x_0)$ is fixed by $h_2$. Suppose finally that $h_2(y) \geq z > y$ 
(resp. $h_2(y) \leq z < y$) for some $y,z$ in $]x_1,x_0[$. Then there exists 
$N \!\in\! \mathbb{N}$ such that $h_2(h_1^n(y)) < h_1^n(z)$ \esp (resp. 
$h_2(h_1^n(y)) > h_1^n(z)$) \esp for all $n \geq N$.}
\label{kopelcito}
\end{lem}

\noindent{\bf Proof.} Notice that the sequence $(x_n)$ tends to zero as $n \!\!\in\!\! \mathbb{N}$ 
goes to infinity. Let $\delta \!=\! V(h_1;[0,x_0])$. For all $u,v$ in $[x_1,x_0]$ and all 
$n \!\in\! \mathbb{N}$ one has (compare inequality (\ref{clasiquita}))
\begin{equation}
\left| \log \Big( \frac{(h_1^n)'(u)}{(h_1^n)'(v)} \Big) \right| \leq \delta.
\label{gatito}
\end{equation}
Since $h_2$ fixes each point $x_n$, one necessarily has $h_2'(0) = 1$. Suppose that the 
claim of the lemma is not satisfied by some $y < z$ in $]x_1,x_0[$ (the other case may be 
reduced to this one by changing $h_2$ by $h_2^{-1}$). Let $\kappa$ be a constant such that 
$$1 < \kappa < 1 + \frac{z-y}{e^{\delta}(y-x_1)}.$$
Fix $N \in \mathbb{N}$ large enough so that 
\begin{equation}
\quad h_2'(w) \leq \kappa \hspace{0.1cm} \esp \mbox{ and } \hspace{0.1cm} \esp
(h_2^{-1})'(w) \leq \kappa \esp \mbox{ for all } \hspace{0.1cm} w \in
{[x_{N+1},x_N]}.
\label{controlcapel}
\end{equation}
For some $n \geq N$ we have $h_2(h_1^n(y)) \geq h_1^n(z)$. Let 
$y_n = h_1^n(y)$ and $z_n = h_1^n(z)$. By the Mean Value Theorem, 
there must exist points $u, v$ in $[x_1,x_0]$ such that
$$\frac{h_2(y_n) - y_n}{y_n - x_{n+1}} \esp \geq \esp 
\frac{z_n - y_n}{y_n - x_{n+1}} \esp = \esp 
\frac{(h_1^n)'(u) \esp (z-y)}{(h_1^n)'(v) \esp (y-x_1)}.$$
From inequality (\ref{gatito}) it follows that
$$\frac{h_2(y_n) - y_n}{y_n - x_{n+1}} \esp \geq \esp 
\frac{z-y}{e^{\delta}(y-x_1)} \esp > \esp \kappa -1.$$
One then deduces that, for some $w \in [x_{n+1},y_n]$,
$$h_2'(w) \esp = \esp \frac{h_2(y_n)-h_2(x_{n+1})}{y_n - x_{n+1}} 
\esp = \esp \frac{h_2(y_n) - x_{n+1}}{y_n - x_{n+1}} \esp > \esp \kappa,$$
which contradicts (\ref{controlcapel}). $\hfill\square$

\vspace{0.35cm}

In \S \ref{muylargo}, we introduced an elementary criterion (namely Lemma \ref{semilibre}) 
for detecting free semigroups inside groups of interval homeomorphisms. We now give 
another criterion for the smooth case which will be fundamental for the proof of 
Theorem \ref{mio}. Example \ref{wreath-exp} well illustrates the lemma below.

\vspace{0.1cm}

\begin{lem} {\em Let $f$ and $g$ be elements in $\mathrm{Homeo}_+([0,1[)$ such that $f(x) < x$ 
for all $x \! \in ]0,1[$. Suppose that there exists an interval $[a,b]$ contained in $]0,1[$ 
such that $g(x) < x$ for every $x \! \in ]a,b[ $, $g(a)=a$, $g(b)=b$, and $f(b) \leq a$. Let 
$[y,z] \subset ]a,b[$ be an interval so that either $g(y) \geq z$ or $g(z) \leq y$. Suppose 
that $g$ fixes all the intervals $f^n([a,b])$, and that there exists $N_0 \in \mathbb{N}$ such 
that $g$ has fixed points inside $f^n(]a,b[)$ for every $n \geq N_0$. If the group generated 
by $f$ and $g$ has no crossed elements, then the semigroup generated by these elements 
is free.}
\label{log-semilibre}
\end{lem}

\noindent{\bf Proof.} Changing $[a,b]$ by its image under some positive iterate of $f$, 
we may assume that $g$ has no fixed point inside $]a,b[$, but it has fixed points in 
$f^n(]a,b[)$ for all $n \geq 1$. We need to show that any two different words $W_1$ 
and $W_2$ in positive powers of $f$ and $g$ represent different homeomorphisms. 
Up to conjugacy, we may suppose that these words are of the form \esp 
$W_1 = f^{n}g^{m_r}f^{n_r} \cdots g^{m_1}f^{n_1}$ \esp and \esp 
$W_2 = g^{q}f^{p_s}g^{q_s} \cdots f^{p_1}g^{q_1}$, \esp where 
$m_j,n_j,p_j,q_j$ are positive integers, $n \geq 0$, and $q \geq 0$ (with $n > 0$
if $r=0$, and $p>0$ if $s=0$). Since the group generated by $f$ and $g$ 
\index{crossed elements}
has no crossed elements, according to Proposition \ref{radon-inv} this 
group must preserve a Radon measure $v$ on $]0,1[$. Notice that 
\index{measure!Radon}
$$\tau_{v}(W_1) = (n_1+\ldots+n_r+n) \esp \tau_{v}(f), 
\qquad \tau_{v}(W_2) = (p_1+\ldots+p_s) \esp \tau_{v}(f).$$
Since $\tau_v (f) \neq 0$, if the values of $(n_1+\ldots+n_r+n)$ and $(p_1+\ldots+p_s)$
are different then $\tau_v (W_1) \neq \tau_v (W_2)$, 
and hence $W_1 \neq W_2$. Now assuming that these values are both equal to 
some $N \in \mathbb{N}$, we will show that the elements $f^{-N} W_1$ and $f^{-N} W_2$ 
are different. To do this, notice that
\begin{eqnarray*}
f^{-N} W_1
\!\!\!\!&=&\!\!\!\! 
f^{-N} f^n g^{m_r} f^{n_r} \cdots g^{m_1} f^{n_1} \\
&=&\!\!\!\! f^{-N} f^n g^{m_r} f^{n_r} \!\cdots g^{m_2} f^{n_1+n_2}
\big( f^{-n_1} g^{m_1} f^{n_1} \big) \\
&=&\!\!\!\! f^{-N} f^n g^{m_r} f^{n_r}\! \cdots g^{m_3} f^{n_1+n_2+n_3}
\big( f^{-(n_1+n_2)} g^{m_2} f^{n_1+n_2} \big)
\big( f^{-n_1} g^{m_1} f^{n_1} \big) \\
&\vdots& \\
&=& \!\!\!\! \big( f^{-(N-n)} g^{m_r} f^{N-n}\big)\! \cdots
\big( f^{-(n_1+n_2)} g^{m_2} f^{n_1+n_2} \big)
\big( f^{-n_1} g^{m_1} f^{n_1} \big)
\end{eqnarray*}
and
\begin{eqnarray*}
f^{-N} W_2\!\!\!
&=& \!\!\! f^{-N} g^{q} f^{p_s} g^{q_s} \cdots f^{p_1} g^{q_1} \\
&=& \!\!\! f^{-N} g^{q} f^{p_s} g^{q_s} \cdots f^{p_3} g^{q_3} f^{p_1+p_2}
\big( f^{-p_1} g^{q_2} f^{p_1} \big) g^{q_1} \\
&\vdots& \\
&=& \!\!\! \big(f^{-N} g^q f^N \big) \cdots
\big( f^{-p_1} g^{q_2} f^{p_1} \big) g^{q_1}.
\end{eqnarray*}
Since the group generated by $f$ and $g$ does not contain crossed 
elements, and since all of the maps \esp
$f^{-(N-n)} g^{m_r} f^{N-n}, \ldots, 
f^{-(n_1+n_2)} g^{m_2} f^{n_1+n_2},$ \esp 
$f^{-n_1} g^{m_1} f^{n_1}\!,$ \esp and \esp 
$f^{-N} g^q f^N$,$\ldots$, $f^{-p_1} g^{q_2} f^{p_1}$, \esp 
have fixed points inside $]a,b[$, these maps must fix $]a,b[$. On the other hand, $g^{q_1}$ fixes 
the interval $]a,b[$, but it does not have fixed points in its interior. Therefore, if $\bar{v}$ 
is a Radon measure on $]a,b[$ which is invariant by the group generated by (the restrictions 
\index{measure!Radon}
to $]a,b[$ of) all of these maps (including $g^{q_1}$), then
\esp $\tau_{\bar{v}}(f^{-N}W_1) = 0$ \esp and \esp
$\tau_{\bar{v}}(f^{-N}W_2) = \tau_{\bar{v}}(g^{q_1}) \neq 0$. 
\esp This shows that $f^{-N} W_1 \neq f^{-N} W_2$. $\hfill\square$

\vspace{0.45cm}

We may now proceed to the proof of Theorem \ref{mio}. Let $\Gamma$ be a finitely generated 
subgroup of $\mathrm{Diff}^{1+\mathrm{bv}}_+([0,1[)$ without free semigroups on two generators. 
To show that $\Gamma$ is Abelian, without loss of generality we may suppose that $\Gamma$ has 
no global fixed point in $]0,1[$. By Proposition \ref{radon-inv}, $\Gamma$ preserves a Radon 
measure $v$ on $]0,1[$; moreover, Exercise \ref{radon-inv-ejer} provides us with an element 
$f \!\in\! \Gamma$ such that $f(x) \!<\! x$ for all $x \!\in\! \mathbb{R}$. Let $\Lambda$ 
be the set of points in $]0,1[$ which are globally fixed by the action of the derived group 
$\Gamma' \!\!=\!\! [\Gamma,\Gamma]$. The set $\Lambda$ is non-empty, since it contains the 
support of $v$. If $\Lambda$ coincides with $]0,1[$, then $\Gamma$ is Abelian. Suppose 
now that $\Lambda$ is strictly contained in $]0,1[$ and that the restriction of $\Gamma'$ 
to each connected component of \esp $]0,1[ \esp\!\setminus \Lambda$ \esp is free. 
\index{H\"older!theorem} \index{action!free}By H\"older's theorem, 
the restriction of $\Gamma'$ to each of these connected components is Abelian, and hence 
$\Gamma$ is metabelian. By Rosenblatt's theorem \cite{br-ros,ros}, $\Gamma$ is virtually 
nilpotent, and by Plante-Thurston's theorem, $\Gamma$ is virtually Abelian. Finally, 
from Exercise \ref{virt-abel} one concludes that, in this case, $\Gamma$ is Abelian.

It remains the case where the action of $\Gamma'$ on some connected component $I$ of the 
complement of $\Lambda$ is not free. The proof of Theorem \ref{mio} will be then finished 
by showing that, in this situation, $\Gamma$ contains free semigroups on two generators. 
For this, let us consider an element $h \in \Gamma'$ and an interval $]y,z[$ strictly contained 
in $I$ so that $h$ fixes $]y,z[$, but no point in $]y,z[$ is fixed by $h$. We claim that there 
must exist an element $g \in \Gamma'$ sending $]y,z[$ into a disjoint interval (contained in 
$I$). Indeed, if this is not the case then, since $\Gamma$ has no crossed elements, every 
element in $\Gamma'$ must fix $]y,z[$, and hence the points $y$ and $z$ are contained in 
$\Lambda$, thus contradicting the fact that $]y,z[$ is strictly contained in the 
connected component $I$ of $]0,1[ \setminus \Lambda$.

The element $g \in \Gamma'$ fixes the connected component $f^n(I)$ of 
$]0,1[ \setminus \Lambda$ for every $n \geq 0$. Moreover, since $\Gamma$
has no crossed elements, for each $n \geq 0$ the intervals $f^n(]y,z[)$ 
and $gf^n(]y,z[)$ either coincide or are disjoint. Lemma \ref{kopelcito} 
(applied to $h_1 = f$ and $h_2 = g$, with $x_0$ being the right endpoint 
of $I$) implies the existence of an integer $N_0 \in \mathbb{N}$ such 
that $g$ fixes the interval $f^n(]y,z[)$ for every $n \geq N_0$ (see 
Figure 23). By considering the dynamics of $f$ and $g$ on the open convex 
closure of $\cup_{n \in \mathbb{Z}} g^n(]y,z[)$, it follows from Lemma 
\ref{log-semilibre} that the semigroup generated by these elements 
is free. This concludes the proof of Theorem \ref{mio}.

\vspace{0.15cm}

\begin{small} \begin{ejer} Show that subgroups 
of $\mathrm{Diff}_+^{1+\mathrm{bv}}(\clo)$ (resp.  
$\mathrm{Diff}_+^{1+\mathrm{bv}}(\mathbb{R})$) without 
free semi- groups on two generators are Abelian (resp. metabelian)   
(compare Theorems \ref{casi-ultimo} and \ref{ultimo}).
\end{ejer} \end{small}

\vspace{0.25cm}

\beginpicture

\setcoordinatesystem units <0.75cm,0.75cm>

\putrule from -0.9 0 to 14.1 0

\begin{small}
\put{$w$} at 7 -0.6
\put{$z$} at 11.5 -0.6
\put{$f^{N_0}(w)$} at 1.2 -0.6
\put{$f^{N_0}(z)$} at 4.9 -0.6

\putrule from 7 -1.4 to 11.5 -1.4
\putrule from 1.2 -1.4 to 4.8 -1.4
\put{$($} at 1.2 -1.4
\put{$($} at 7 -1.4
\put{$)$} at 4.8 -1.4
\put{$)$} at 11.5 -1.4
\put{$I$} at 9.25 -1.1
\put{$f^{N_0}(I)$} at 3 -1.1

\put{$($} at 1.2 0
\put{$($} at 7 0
\put{$)$} at 4.8 0
\put{$)$} at 11.5 0
\end{small}

\begin{Large}
\put{$0$} at -1.4 -0.9
\put{$1$} at 13.6 -0.9
\put{$($} at -1.4 0
\put{$)$} at 13.6 0
\end{Large}

\begin{tiny}
\put{$g(u)$} at 9.5 -0.35
\put{$g(v)$} at 10.5 -0.35
\put{$u$} at 7.5 -0.35
\put{$v$} at 8.5 -0.35
\put{$f^{N_0}(u)$} at 1.8 -0.35
\put{$f^{N_0}(v)$} at 3.2 -0.35
\put{$($} at 7.5 0
\put{$($} at 9.5 0
\put{$($} at 1.7 0
\put{$)$} at 10.5 0
\put{$)$} at 8.5 0
\put{$)$} at 3.3 0
\end{tiny}

\circulararc 38 degrees from 9 1.2
center at 6 -7.5

\circulararc 105 degrees from 10 0.25
center at 9 -0.5

\circulararc 180 degrees from 2.9 0.25
center at 2.4 0.25

\put{$g$} at 9 0.4
\put{$g$} at 2.4 0.4
\put{$f^{N_0}$} at 6.1 1.3
\put{$\bullet$} at 8 0.25
\put{$\bullet$} at 9 1.2
\put{$\bullet$} at 1.9 0.25

\plot
2.9 0.25
2.76 0.41 /

\plot
2.9 0.25
2.97 0.45 /

\plot
3 1.2
3.32 1.2 /

\plot
3 1.2
3.28 1.4 /

\plot
10 0.25
9.8 0.35 /

\plot
10 0.25
9.95  0.45 /

\put{Figure 23} at 5.8 -2

\put{} at -2.4 0

\endpicture


\vspace{0.15cm}

\begin{small} \begin{ejer} Throughout this exercise, $\Gamma$ will always 
be a finitely generated group of polynomial growth of degree $k$.
\index{degree!of polynomial growth for a group}

\noindent (i) Prove that every finitely generated subgroup of $\Gamma$ 
has polynomial growth of degree less than or equal to $k$.

\noindent (ii) Prove that if \esp \esp
$\Gamma_n \subset \Gamma_{n-1} \subset \ldots \subset \Gamma_0 = \Gamma$ \esp \esp 
is a series of subgroups such that for every $i \in \{0,\ldots,n-1\}$ there exists 
a nontrivial group homomorphism $\phi_i\!: \Gamma_i \rightarrow (\mathbb{R},+)$ 
satisfying $\phi_i (\Gamma_{i+1}) = \{0\}$, then $n \le k$.

\noindent{\underbar{Hint.}} For each $i \in \{0,\ldots,n-1\}$ choose an element 
$g_i \in \Gamma_{i}$ such that $\phi_i (g_i) \neq 0$, and consider the words 
of the form \esp $W = g_1^{m_1} g_2^{m_2} \cdots g_n^{m_n}$.

\noindent (iii) Suppose that $\Gamma$ is a subgroup of $\mathrm{Homeo}_+([0,1])$. 
Using (ii) and the translation number homomorphism, show that $\Gamma$ is solvable 
with degree of solvability less than or equal to $k$. Applying Milnor-Wolf's 
\index{degree!of solvability|(}
theorem \cite{Mil,Wo}, 
\index{Milnor-Wolf's theorem} 
conclude that $\Gamma$ is virtually nilpotent.

\noindent (iv) Suppose now that $\Gamma$ is a subgroup of 
$\mathrm{Diff}_+^{1+\mathrm{bv}}([0,1[)$. Using Theorem \ref{chancho} 
and Exercise \ref{virt-abel}, conclude that $\Gamma$ is Abelian.

\noindent (v) Prove that the claim in (iv) still holds if $\Gamma$ is contained 
in the group of germs $\mathcal{G}_{+}^{1+\mathrm{bv}}(\mathbb{R},0)$.

\noindent{\underbar{Hint.}} Instead of using the translation number 
homomorphism, use Thurston's stability theorem from \S \ref{sec-th-th}.
\label{original-PL} \index{translation number}
\end{ejer}

\begin{ejer} The claim in item (ii) below may be proved either by a direct 
argument or by using the Positive Ping-Pong Lemma (see \cite{br-ros}).
\index{Klein!Ping-Pong Lemma}

\noindent (i) Under the hypothesis of Lemma \ref{log-semilibre}, show that the group 
generated by the family of elements $\{f^ngf^{-n} \! : n \in \mathbb{Z} \}$ 
is not finitely generated.\\

\noindent (ii) Prove that if $f,g$ are elements in a group so that 
the subgroup generated by $\{f^ngf^{-n} \!\! : n \in \mathbb{Z} \}$ is not 
finitely generated, then the semigroup generated by $f$ and $gfg^{-1}$ is free.

\noindent (iii) Use (ii) to give an alternative proof to Theorem \ref{mio}.
\end{ejer}

\begin{obs} We ignore whether Theorem \ref{mio} extends (at least in the case 
of sub-exponential growth) to groups of {\em germs} of $\ce^{1+\mathrm{bv}}$ 
diffeomorphisms of the line fixing the origin (compare Exercise \ref{PT-ger}).
\index{germ of a diffeomorphism}
\end{obs} \end{small}


\subsection{Nilpotence, growth, and intermediate regularity}

\hspace{0.45cm} As in the cases of Denjoy's theorem and Kopell's lemma, the results from the 
preceding two sections are no longer true in regularity less than $\ce^{1+\mathrm{Lip}}$. 
Indeed, using the first of the techniques in \S \ref{ejemplosceuno}, Farb and 
Franks provided in \cite{FF} a 
$\ce^1$ realization of each group $\mathrm{N}_n$ from Exercise \ref{columnas} 
(notice that $\mathrm{N}_n$ is nilpotent for every $n$ and non-Abelian for $n \!\geq\! 3$). 
On the other hand, according to a result by Malcev (see \cite{Ra}), every finitely 
generated torsion-free \index{group!torsion-free} nilpotent \index{group!nilpotent} 
group embeds into some $\mathrm{N}_n$. All of this shows the existence of a large 
variety of $\ce^1$ counterexamples to Plante-Thurston's theorem. 
\index{Plante-Thurston's theorem}

In fact, using the second method of construction from \S \ref{ejemplosceuno}, 
it is possible to show that the canonical action of $N_n$ on the interval 
may be smoothed until reaching any differentiability class less than 
$\ce^{1+1/(n-2)}$. Notice that, by Theorem \ref{kop-gen}, this 
action cannot exceed this regularity. 
It is then tempting to conjecture some general result which relates 
the nilpotence degree and the optimal regularity for finitely generated 
torsion-free nilpotent group actions on the interval. However, one may show 
that the groups from Example \ref{nilp-locos} can be realized as groups of 
$\ce^{1+\tau}$ diffeomorphisms of $[0,1]$ for any $\tau < 1$ (see \cite{edo}). 
Actually, it seems that the right parameter giving the obstruction for sharp 
smoothing corresponds to the degree of solvability rather than the degree 
of nilpotence.
 
The following general question remains, therefore, open: Given a finitely generated 
torsion-free nilpotent group $\Gamma$, what is the best regularity for faithful 
actions of $\Gamma$ by diffeomorphisms of the (closed) interval~?

\begin{small}
\begin{obs} For each pair of positive integers $m<n$, the group $\mathrm{N}_m$ 
naturally embeds into $\mathrm{N}_n$. If we denote by N the union of all of these groups, 
then one easily checks that N is orderable, countable, and non finitely generated. Moreover, 
N contains all torsion-free finitely generated nilpotent groups. It would be interesting 
to check whether the methods from \S \ref{ejemplosceuno} allow showing that N is 
isomorphic to a group of $\ce^1$ diffeomorphisms of the interval.
\end{obs}
\end{small}

Summarizing the remarks above, every finitely generated torsion-free nilpotent group 
may be embedded into $\mathrm{Diff}_+^{1+\tau}([0,1])$ for some $\tau < 1$ small enough. 
It is then natural to ask whether some other groups of sub-exponential growth may appear 
as groups of $\ce^{1+\tau}$ diffeomorphisms of the interval. According to Exercise 
\ref{chanta-chanta}, the following result (due to the author \cite{ceuno}) 
provides a negative answer to this.

\vspace{0.082cm}

\begin{thm} {\em For all $\tau > 0$, every finitely generated group 
of $\ce^{1+\tau}$ diffeomorphisms of the interval $[0,1[$ which does 
not contain free semigroups on two generators is virtually nilpotent.}
\label{growth-thm}
\end{thm}

\vspace{-0.2cm}

\begin{small}
\begin{obs}
In recent years, another notion has come to play 
an important role in the theory of growth of groups, 
namely the notion of {\bf {\em uniformly exponential growth}}. 
\index{growth!of groups}
One says that a finitely generated 
group $\Gamma$ has uniformly exponential growth if there exists a constant $\lambda>1$ such 
that for every (symmetric) finite system of generators $\mathcal{G}$ of $\Gamma$, the 
number $L_{\mathcal{G}}(n)$ of elements in $\Gamma$ that may be written as a product 
of at most $n$ elements in $\mathcal{G}$ satisfies
$$\lim_{n \rightarrow \infty} \frac{\log \big( L_{\mathcal{G}}(n) \big)}{n} \geq \lambda.$$
Although in many situations groups of exponential growth are forced to have uniform exponential 
growth (see for instance \cite{br-ros,Br-Gel,eskin-oh,osin2}), there exist groups 
for which these properties are not equivalent \cite{wilson}. We ignore whether 
this may happen for groups of $\ce^{1+\tau}$ diffeomorphisms of the interval.
\end{obs}
\end{small}

The condition $\tau > 0$ in Theorem \ref{growth-thm} is necessary. To show this, below we 
sketch the construction of a subgroup of $\mathrm{Diff}_+^1([0,1])$ having intermediate 
growth. This subgroup is closely related to the group $\hat{\mathrm{H}}$ from Example 
\ref{grigor-ejem}. Indeed, although $\hat{\mathrm{H}}$ is a torsion group (the reader may easily 
check that the order of every element is a power of $2$; see also \cite{harpe}), starting with 
$\hat{\mathrm{H}}$ 
one may create a {\em torsion-free} group of intermediate growth. This group was introduced 
in \cite{grigor2}. Geometrically, the idea consists in replacing $\mathcal{T}_2$ by a rooted 
tree whose vertexes have (countable) infinite valence. More precisely, we consider the group 
$\bar{\mathrm{H}}$ acting on the space $\Omega = \mathbb{Z}^{\mathbb{N}}$ and generated 
by the elements $\bar{a},\bar{b},\bar{c}$, and $\bar{d}$, recursively defined by 
\esp $\bar{a}(x_1,x_2,x_3, \ldots) = (1+x_1,x_2,x_3, \ldots)$ \esp and
$$\bar{b}(x_1,x_2,x_3,\ldots) = \left \{ \begin{array} {l}
(x_1, \bar{a}(x_2,x_3, \ldots )), \hspace{0.1cm} x_1 \mbox{ even},\\
(x_1, \bar{c}(x_2,x_3,\ldots)), \hspace{0.12cm} x_1 \mbox{ odd}, \end{array} \right.$$
$$\bar{c}(x_1,x_2,x_3,\ldots) = \left \{ \begin{array} {l}
(x_1, \bar{a}(x_2,x_3, \ldots )), \hspace{0.1cm} x_1 \mbox{ even},\\
(x_1, \bar{d}(x_2,x_3,\ldots)), \hspace{0.1cm} x_1 \mbox{ odd}, \end{array} \right.$$
$$\bar{d}(x_1,x_2,x_3,\ldots) = \left \{ \begin{array} {l}
(x_1, x_2,x_3, \ldots ), \hspace{0.55cm} x_1 \mbox{ even},\\
(x_1, \bar{b}(x_2,x_3,\ldots)), \hspace{0.1cm} x_1 \mbox{ odd}. \end{array} \right.$$

The group $\bar{\mathrm{H}}$ preserves the lexicographic ordering on $\Omega$, 
from where 
\index{lexicographic ordering}
one concludes that it is orderable, 
\index{group!orderable} 
and hence according to \S \ref{sec-witte}, 
it may be realized as a group of homeomorphisms of the interval (compare \cite{GrM}). 
To explain this better, we give below a concrete realization of 
$\bar{\mathrm{H}}$ as a group of bi-Lipschitz homeomorphisms of $[0,1]$ 
\index{Lipschitz!homeomorphism}
(compare Proposition \ref{abu}).

\vspace{0.1cm}

\begin{small} \begin{ejem} Given $C > 1$, let $(\ell_i)_{i \in \mathbb{Z}}$ 
be a sequence of positive real numbers such that 
$$\sum_{i \in \mathbb{Z}} \ell_i \!=\! 1, \qquad 
\max \left\{ \frac{\ell_{i+1}}{\ell_i}, \frac{\ell_i}{\ell_{i+1}} \right\}
\esp \leq \esp C \quad \mbox {for all} \quad i \in \mathbb{Z}.$$
Let $I_i$ be the interval \esp $]\sum_{j<i}\ell_j,\sum_{j\leq i}\ell_j[$, \esp and let 
$f\!: [0,1] \rightarrow [0,1]$ be the homeomorphism sending affinely each interval $I_i$ onto 
$I_{i+1}$. Let us denote by $g$ the affine homeomorphism sending $[0,1]$ onto $\bar{I}_0$, 
and by $\lambda = 1/\ell_0$ the (constant) value of its derivative. Let $a,b,c$, and $d$, 
be the maps recursively defined on a dense subset of $[0,1]$ by $a(x) = f(x)$ 
and, for $x \in I_i$,
$$b(x) = \left \{ \begin{array} {l}
f^i g a g^{-1} f^{-i}(x), \hspace{0.1cm} i \mbox{ even},\\
f^i g c g^{-1} f^{-i}(x), \hspace{0.1cm} i \mbox{ odd}, \end{array} \right.$$
$$c(x) = \left \{ \begin{array} {l}
f^i g a g^{-1} f^{-i}(x), \hspace{0.1cm} i \mbox{ even},\\
f^i g d g^{-1} f^{-i}(x), \hspace{0.1cm} i \mbox{ odd}, \end{array} \right.$$
$$d(x) = \left \{ \begin{array} {l}
x, \hspace{1.985cm} i \mbox { even},\\
f^i g b g^{-1} f^{-i}(x), \hspace{0.1cm} i \mbox{ odd}. \end{array} \right.$$
We claim that $a, b, c$, and $d$, are bi-Lipschitz homeomorphisms with Lipschitz constant 
bounded from above by $C$ (notice that $C$ may be chosen as near to $1$ as we want). 
Indeed, this is evident for $a$, while for $b$, $c$, and $d$, this may be easily 
checked by induction. For instance, if $x \in I_i$ for some odd integer $i$, then
$$b'(x) = \frac{(f^i)'(gag^{-1}f^{-i}(x))}{(f^i)'(f^{-i}(x))} \cdot
\frac{g'(ag^{-1}f^{-i}(x))}{g'(g^{-1}f^{-i}(x))} \cdot a'(g^{-1}f^{-i}(x)),$$
and since \esp $g'|_{[0,1]} \!\equiv\! \lambda$ \esp and \esp 
$(f^i)'|_{I_0} \!\equiv\! \ell_i/\ell_0$, \esp we have \esp 
$b'(x) \!=\! a'(g^{-1}f^{-i}(x)) \!\leq\! C$. \esp It is 
geometrically clear that the group generated by $a,b,c$, 
and $d$, is isomorphic to $\bar{\mathrm{H}}$.
\label{primerlisa}
\end{ejem}

\vspace{0.01cm}

\begin{ejer} Give an example of a finitely generated group $\Gamma$ of homeomorphisms of 
the interval and/or the circle for which there exists a finite system of generators 
and a constant $\delta > 0$ such that, for every topological conjugacy of $\Gamma$ 
to a group of bi-Lipschitz homeomorphisms, at least one of these generators or its 
inverse has Lipschitz constant bigger than or equal to $1 + \delta$.

\noindent{\bf Remark.} It seems to be interesting to determine whether $\Gamma$ above 
may be a group of $ \ce^1$ diffeomorphisms without crossed elements ({\em c.f.}, 
Definition \ref{def-crossed-elements}).
\label{mas-adelante}\index{crossed elements}
\end{ejer} \end{small}

The preceding idea is inappropriate to obtain an embedding of $\bar{\mathrm{H}}$ 
into $\mathrm{Diff}_+^1([0,1])$, since the derivatives of the maps that are produced 
have discontinuities at each ``level'' of the action of H. Next we provide another 
realization along which it will be essential to ``renormalize'' the geometry 
at each step. We begin by fixing an equivariant family of homeomorphisms \esp 
$\{ \varphi_{u,v} \!: \esp u>0,v>0 \}.$ For each $n \in \mathbb{N}$ 
and each $(x_1,\ldots,x_n)$ in $\mathbb{Z}^n$, 
let us consider a non-degenerate closed interval 
$I_{x_1,\ldots,x_n}\!=\![u_{x_1,\ldots,x_n},w_{x_1,\ldots,x_n}]$ 
and an interval $J_{x_1,\ldots,x_n}\!=\![v_{x_1,\ldots,x_n},w_{x_1,\ldots,x_n}]$, 
both contained in some interval $[0,T]$. Suppose that the 
conditions below hold (see Figure 24):

\vspace{0.08cm}

\noindent{(i) $\sum_{x_1 \in \mathbb{Z}} |I_{x_1}| = T$,}

\vspace{0.08cm}

\noindent{(ii) $u_{x_1,\ldots,x_n} < v_{x_1,\ldots,x_n} \leq w_{x_1,\ldots,x_n}$, 
\esp so that $J_{x_1,\ldots,x_n} \subset I_{x_1,\ldots,x_n}$,}

\vspace{0.08cm}

\noindent{(iii) $w_{x_1,\ldots,x_{n-1},x_n} = u_{x_1,\ldots,x_{n-1},1+x_n}$,}

\vspace{0.08cm}

\noindent{(iv) $\lim_{x_n \rightarrow -\infty} u_{x_1,\ldots,x_{n-1},x_n} =
u_{x_1,\ldots,x_{n-1}}$,}

\vspace{0.08cm}

\noindent{(v) $\lim_{x_n \rightarrow \infty} u_{x_1,\ldots,x_{n-1},x_n} =
v_{x_1,\ldots,x_{n-1}}$,}

\vspace{0.08cm}

\noindent{(vi) $\lim_{n \rightarrow \infty} \sup_{(x_1,\ldots,x_n) \in \mathbb{Z}^n}
|I_{x_1,\ldots,x_n}| = 0$.}

\vspace{0.32cm}

\beginpicture
\setcoordinatesystem units <0.8cm,0.8cm>

\putrule from 0 0 to 12 0
\putrule from 0 1.5 to 5 1.5
\putrule from 7.2 1.5 to 12 1.5
\putrule from 5.3 0.7 to 7.85 0.7
\putrule from 9.6 0.7 to 12 0.7
\putrule from 2.17 -0.5 to 3.2 -0.5

\put{$J_{x_1,\ldots,x_n}$} at 8.8 0.7

\begin{small}
\put{$I_{x_1,\ldots,x_n,x_{n\!+\!1}}$} at 2.6 -1
\end{small}

\begin{tiny}
\put{$u_{x_1\!,\!\ldots\!,\!x_n\!,\!x_{n\!+\!1}}$} at 1.45 0.5
\put{$w_{x_1\!,\!\ldots\!,\!x_n\!,\!x_{n\!+\!1}}$} at 3.35 0.5
\put{$|$} at 1.8 0
\put{$|$} at 2.8 0
\put{$|$} at 1.8 -0.5
\put{$|$} at 2.8 -0.5
\end{tiny}

\begin{Large}
\put{$u_{x_1,\ldots,x_n}$} at -0.8 -0.9
\put{$w_{x_1,\ldots,x_n}$} at 11.6 -0.9
\put{$|$} at -0.75 0
\put{$|$} at 11.25 0
\put{$v_{x_1,\ldots,x_n}$} at 4.6 -0.9
\put{$|$} at 4.65 0
\put{$I_{x_1,\ldots,x_n}$} at 5.45 1.5
\end{Large}

\put{Figure 24} at 5.12 -1.9

\put{} at -2.2 2.2

\endpicture


\vspace{0.5cm}

\noindent Notice that
\begin{equation}
|J_{x_1,\ldots,x_n}| \esp + \sum_{x_{n+1} \in \mathbb{Z}}
|I_{x_1,\ldots,x_n,x_{n+1}}| \esp = \esp |I_{x_1,\ldots,x_n}|.
\label{defjota}
\end{equation}
Let us denote by $\bar{\mathrm{H}}_n$ the stabilizer in $\bar{\mathrm{H}}$ of the 
$n^{\mathrm{th}}$-level of the tree $\mathcal{T}_{\infty}$. For each $n \in \mathbb{N}$ 
we will define four homeomorphisms $a_n, b_n, c_n$, and $d_n$, so that the group generated by 
them is isomorphic to $\bar{\mathrm{H}} / \bar{\mathrm{H}}_n$. For this, let us consider 
the group homomorphisms $\phi_0$ and $\phi_1$ from the subgroup 
$\langle \bar{b}, \bar{c}, \bar{d} \rangle \subset \bar{\mathrm{H}}$ 
into $\bar{\mathrm{H}}$ defined by 
$$\phi_0(\bar{b}) = \bar{a}, \esp \phi_0(\bar{c}) = \bar{a}, \esp \phi_0(\bar{d}) = id, 
\quad \mbox{  } \quad \phi_1(\bar{b}) = \bar{c}, \esp \phi_1(\bar{c}) = \bar{d},
\esp \phi_1(\bar{d}) = \bar{b}.$$

\noindent{\underbar{Definition of $a_n$.}}

\vspace{0.08cm}

\noindent{-- If $p \in J_{x_1,\ldots,x_i}$ for some $i \!<\! n$, we let  
$a_n(p) \!=\! \varphi(J_{x_1,x_2,\ldots,x_i},J_{1+x_1,x_2,\ldots,x_i})(p)$.}

\vspace{0.08cm}

\noindent{-- For $p \in I_{x_1,\ldots,x_n}$, we let 
$a_n(p) = \varphi(I_{x_1,x_2,\ldots,x_n},I_{1+x_1,x_2,\ldots,x_n})(p)$.}

\vspace{0.2cm}

\noindent{\underbar{Definition of $b_n$.}}

\vspace{0.08cm}

Suppose that \esp $p \! \in ]0,1[$ \esp belongs to $I_{x_1,\ldots,x_n}$, and denote by \esp 
$(\bar{x}_1,\ldots,\bar{x}_n) \in \{0,1\}^n$ \esp the sequence obtained after reducing 
each entry modulo 2.

\vspace{0.08cm}

\noindent{-- If $\phi_{\bar{x}_1}(\bar{b})$, $\phi_{\bar{x}_2} \phi_{\bar{x}_1}(\bar{b})$, 
$\ldots$, $\phi_{\bar{x}_n} \ldots \phi_{\bar{x}_2} \phi_{\bar{x}_1}(\bar{b})$ are 
well-defined, we let $b_n(p) = p$.}

\vspace{0.08cm}

\noindent{-- In the other case, we denote by \esp $i=i(p) \leq n$ \esp the smallest 
integer for which \esp 
$\phi_{\bar{x}_i} \ldots \phi_{\bar{x}_2} \phi_{\bar{x}_1}(\bar{b})$ \esp 
is not defined, and we let}

\newpage

$$b_n(p) = \left \{ \begin{array} {l}
p, \esp \esp \esp \esp \esp \esp \hspace{4.8cm} p \!\in\! J_{x_1,\ldots,x_j}, \esp j < i,\\ 
\varphi(J_{x_1,\ldots,x_i,\ldots,x_j},J_{x_1,\ldots,1+x_i,\ldots,x_j})(p), 
\esp \esp \esp p \!\in\! J_{x_1,\ldots,x_i,\ldots,x_j}, \esp 
i \!\leq \!j\! <\! n,\\
\varphi(I_{x_1,\ldots,x_i,\ldots,x_n},I_{x_1,\ldots,1+x_i,\ldots,x_n})(p), 
\esp\esp\esp p \!\in\! I_{x_1,\ldots,x_n}. 
\end{array} \right.$$

The definitions of $c_n$ and $d_n$ are similar to that of $b_n$. Clearly, the maps 
$a_n,b_n,c_n$, and $d_n$, extend to homeomorphisms of $[0,T]$. The fact that they generate 
a group isomorphic to $\bar{\mathrm{H}}/\bar{\mathrm{H}}_n$ follows from the equivariance 
properties of the maps $\varphi(I,J)$.   
Condition (vi) implies moreover that the sequences of maps $a_n$, $b_n$, $c_n$, and $d_n$, 
converge to homeomorphisms $a,b,c$, and $d$, respectively, which generate a group 
isomorphic to $\bar{\mathrm{H}}$. We leave the details to the reader.

\vspace{0.1cm}

\begin{small} \begin{ejem}
Given a sequence $(\ell_i)_{i \in \mathbb{Z}}$ of positive real numbers such that 
$\sum \ell_i = 1$, let us define $|I_{x_1,\ldots,x_n}|$ and $|J_{x_1,\ldots,x_n}|$ by 
$$|J_{x_1,\ldots,x_n}| = 0, \qquad  \quad
|I_{x_1,\ldots,x_n}| = \ell_{x_1} \ldots \ell_{x_n}.$$
If we carry out the preceding construction (for $T \!=\! 1$) using the equivariant family 
of affine maps $\varphi([0,u],[0,v])(x) \!=\! vx/u$, then we reobtain the inclusion of 
$\bar{\mathrm{H}}$ in the group of bi-Lipschitz homeomorphisms of the interval from 
Example \ref{primerlisa} (under the same hypothesis $\ell_{i+1}/\ell_i \!\leq\! C$ 
and $\ell_i/\ell_{i+1} \!\leq\! C$ for all $i \in \mathbb{Z}$).
\end{ejem} \end{small}

\vspace{0.1cm}

The details for the rest of the construction are straightforward and we leave them to the reader 
(see \cite{ceuno} in case of problems). Let $\omega$ be a modulus of continuity satisfying $\omega(s) 
= 1/\log(1/s)$ for $s \leq 1/2e$, and such that the map $s \mapsto \omega(s)/s$ is decreasing. 
Fix a constant $C>0$, and for each $k \in \mathbb{N}$ let $T_k = \sum_{i \in \mathbb{Z}}
\frac{1}{(|i|+k)^2}$. Consider an increasing sequence $(k_n)$ of positive integers, and for 
$n \!\in\! \mathbb{N}$ and $(x_1,\ldots,x_n) \!\in\! \mathbb{Z}^n$ let
$$|I_{x_1,\ldots,x_n}| \esp = \esp \frac{1}{(|x_1|+ \ldots |x_n|+k_n)^{2n}}.$$
Using the equivariant family induced by (\ref{yoccocito}), our general method provides us with  
subgroups of $\mathrm{Diff}_+^{1}([0,T_{k_1}])$ generated by elements $a_n,b_n,c_n$, and $d_n$, 
and isomorphic $\bar{\mathrm{H}}/\bar{\mathrm{H}}_n$. The goal then consists in controlling 
the $\ce^{\omega}$-norm for the derivative of these diffeomorphisms. Now, if the sequence 
$(k_n)$ satisfies some ``rapid growth'' conditions  
\index{growth!of a sequence} 
(which depend only on $C$), it is not 
difficult to verify that these norms are bounded from above by $C$ for all $n \in \mathbb{N}$ 
(and the same holds for the derivatives of the inverses of these diffeomorphisms). The 
corresponding sequences are therefore equicontinuous, 
\index{equicontinuous!sequence of functions}
and it is easy to see that they converge 
to $\omega$-continuous maps having $\ce^{\omega}$-norm bounded from above by $C$. These maps 
correspond to the derivatives of C$^{1+\omega}$ diffeomorphisms $a,b,c,d$ (and their 
inverses), which generate a group isomorphic to $\bar{\mathrm{H}}$. Since this group 
acts on the interval $[0,T_{k_1}]$, to obtain an action on $[0,1]$ we may conjugate by the 
affine map $g \!: [0,1] \rightarrow [0,T_{k_1}]$. Since for $k_1$ large enough one has $T_{k_1} \leq 1$,  
this conjugacy procedure does not increase the $\ce^{\omega}$-norm for the derivatives. Finally, notice  
that since for every diffeomorphism $f$ of $[0,1]$ there exists a point at which the derivative 
equals $1$, if the $\ce^{\omega}$-norm of $f'$ is bounded from above by $C$ then
$$\sup_{x \in [0,1]} |f'(x) - 1| \leq C \esp \omega(1).$$
As a consequence, for $C$ small the preceding realization of $\bar{\mathrm{H}}$ 
has its generators near the identity with respect to the $\ce^{1+\omega}$ topology.

\begin{small} \begin{ejer}
Inside Grigorchuk-Maki's group $\bar{H}$ there are ``a lot'' of commuting elements. 
More precisely, for each $d \in \mathbb{N}$ one may choose $d\!+\!1$ elements in 
$\bar{H}$ satisfying a combinatorial property similar to that of Theorem \ref{kop-gen} 
(take for instance $f_1 = \bar{a}^{-2}$, $f_2 = \bar{b}^{-2}$, $f_3 = \bar{a}^{-1}
\bar{b}^{-2} \bar{a}$, etc). In this way, to verify that the natural action of $\bar{H}$ 
is not semiconjugate to an action by $\ce^{1+\tau}$ diffeomorphisms, it suffices to apply 
Theorem \ref{kop-gen} for $d > 1 / \tau$. The items below allow showing the same claim 
via more elementary methods (see \cite{ceuno} for more details).\\

\vspace{0.03cm}

\noindent (i) Prove that if $h$ is a $\ce^{1+\tau}$ diffemorphism of a closed interval 
$[u,v]$, and $C$ denotes the H\"older constant of $h'$, then for every $x \in [u,v]$ 
\index{H\"older!derivative}
one has \esp $|h(x) - x| \leq C \esp |v-u|^{1+\tau}$.

\vspace{0.03cm}

\noindent (ii) Using (i) prove directly ({\em i.e.}, without using any probabilistic argument) 
that Theorem \ref{kop-gen} holds for $d \geq d(\alpha)$, where $d(\alpha)$ is the minimal 
integer greater than or equal to $2$ for which \esp $\alpha (1+\tau)^{d-2} \geq 1$.

\vspace{0.03cm}

\noindent (iii) Using (ii) conclude that the canonical action of $\bar{H}$ is 
not semiconjugate to an action by $\ce^{1+\tau}$ diffeomorphisms for any $\tau > 0$.
\end{ejer} \end{small}

To conclude this section, we notice that the element $\bar{a}^2$ belongs to the 
center of 
\index{center of a group}
$\bar{\mathrm{H}}$. Having in mind the realization of $\bar{\mathrm{H}}$ 
as a group of $\ce^1$ diffeomorphisms of the interval, the study of centralizers in 
$\mathrm{Diff}_+^1([0,1])$ of diffeomorphisms (or homeomorphisms) of the interval 
becomes natural. The next ``weak Kopell lemma'', due to Bonatti, Crovisier, 
and Wilkinson, is an interesting issue in this direction.
\index{centralizer}
\index{Bonatti-Crovisier-Wilkinson's lemma}

\vspace{0.2cm}

\begin{prop} {\em If $h$ is a homeomorphism of $[0,1[$ without fixed 
points in $]0,1[$, then the group of $\mathrm{C}^1$ diffeomorphisms 
of $[0,1[$ that commute with $h$ has no crossed elements.}
\label{facil}
\end{prop}

\noindent{\bf Proof.} Suppose for a contradiction that $f$ and $g$ are $\ce^1$ 
diffeomorphisms of $[0,1[$ which commute with $h$, and that they are crossed 
\index{crossed elements}
on some interval $[x,y] \subset [0,1]$. As in the proof of Lemma \ref{semilibre}, 
we may restrict ourselves to the case where $f(x)=x$, $f(y)\! \in ]x,y[$, $g(x)\! \in ]x,y[$, 
and $g(y)=y$. Moreover, changing $f$ and/or $g$ by some of their iterates, we may suppose 
that $g(x) > f(y)$. All of these properties are preserved under conjugacy, and hence 
$f$ and $g$ satisfy $f(h^n(x)) = h^n(x)$, $f(h^n(y)) \! \in \esp ]h^n(x),h^n(y)[$, 
$g(h^n(x)) \! \in \esp ]h^n(x),h^n(y)[$, $g(h^n(y)) = h^n(y)$, and $g(h^n(x)) > f(h^n(y))$, 
for all $n \in \mathbb{Z}$. If $h(z) < z$ (resp. if $h(z) > z$) for all $z \! \in ]0,1[$, 
then the sequences $(h^n(x))$ and $(h^n(y))$ (resp. $(h^{-n}(x))$ and $(h^{-n}(y))$)
converge to the origin. Since $f$ and $g$ are of class $\ce^1$, this implies that 
$f'(0)=g'(0)=1$. On the other hand, since $g(h^n(x)) > f(h^n(y))$, there must 
exist a sequence of points $z_n \! \in \esp ]x_n,y_n[$ such that, for each 
$n \!\in\! \mathbb{Z}$, either $f'(z_n) \! < \! 1/2$ or $g'(z_n) \! < \! 1/2$. 
However, this contradicts the continuity of the derivatives 
of $f$ and $g$ at the origin. $\hfill\square$ 
\index{growth!of groups|)} \index{semigroup!free|)}


\section{Polycyclic Groups of Diffeomorphisms}
\label{policiclicos-suaves}
\index{group!polycyclic|(}

\hspace{0.45cm} Unlike nilpotent groups, inside the group of $\ce^2$ diffeomorphisms 
of the interval there is a large variety of non-Abelian {\em polycyclic} groups. 
To see this, first notice that the affine group contains many polycyclic groups 
(see Exercise \ref{a-clasif}). On the other hand, the affine group is conjugate 
to a subgroup of $\mathrm{Diff}_+^{\infty}([0,1])$. Indeed, inspired by Exercise 
\ref{truco-de-tsuboi}, 
\index{Tsuboi!M\"uller-Tsuboi conjugacy} 
let us fix a constant $0 < \varepsilon < 1/2$, and 
let us consider two $\ce^{\infty}$ diffeomorphisms 
\index{affine group}
$\varphi_1\!\!: \esp ]0,1[ \rightarrow \mathbb{R}$ 
and $\varphi_2\!\!: \esp ]0,1[ \rightarrow ]0,1[$ such that 
$$\varphi_1(x) = - \frac{1}{x} \quad \mbox{and} \quad
\varphi_2(x) = \exp \left( - \frac{1}{x} \right)
\quad \mbox{ for } x \! \in ]0,\varepsilon],$$
$$\quad \varphi_1(x) = \frac{1}{1-x} \quad \mbox{and} 
\quad \varphi_2(x) =  1 - \exp \left( \frac{1}{x-1} \right)
\quad \mbox{ for } x \in [1-\varepsilon,1[.$$
Consider now the vector fields on the line $Y_1 = \frac{\partial}{\partial x}$ and 
$Y_2 = x\frac{\partial}{\partial x}$, which generate the Lie algebra of the affine group. 
One readily checks that the vector fields $X_j = \varphi^{*}(Y_j)$, where $j \in \{1,2\}$ 
and $\varphi = \varphi_1 \circ \varphi_2^2$, extend to $\mathrm{C}^{\infty}$ vector fields 
on $[0,1]$ which are (zero and) infinitely flat at the endpoints. Due to this, for each 
$g \!\in\! \mathrm{Aff}_{+}(\mathbb{R})$ the map $\varphi^{-1} \circ g \circ \varphi$ 
is a $\mathrm{C}^{\infty}$ diffeomorphism of $[0,1]$ infinitely tangent to the 
identity at the endpoints.

\begin{small}\begin{ejer}
Show that every polycyclic subgroup of the affine group is isomorphic to 
$\mathbb{Z}^k \ltimes \mathbb{Z}^n$ for some non-negative integers $n,k$.
\label{a-clasif}
\end{ejer}\end{small}

The construction above was first given by 
Plante \cite{Pl1,Pl2} who showed that, under 
\index{Plante!examples of solvable groups}
a mild assumption, polycyclic groups of $\ce^2$ diffeomorphisms of the interval are forced 
to be conjugate to subgroups of the affine group: see Theorem \ref{thm-moriyama}, item (i). 
\index{Moriyama}
(Plante's results were used by Matsumoto to classify certain codimension-one 
\index{foliation}
foliations: see \cite{matsumotito}.) 
\index{Matsumoto}
It was Moriyama \cite{moriyama} who provided another kind of examples 
and gave a complete classification. To introduce his examples, given any integer 
$n \!\geq\! 2$ and a matrix $A$ in $\mathrm{SL}(n,\mathbb{Z})$ with an eigenvalue 
$\lambda \! \in ]0,1[$, let $(t_1,\ldots,t_n) \in \mathbb{R}^n$ be an eigenvector 
associated to it. Given an interval $[a,b] \!\subset ]0,1[$, let $X$ be a non-zero 
vector field on $[a,b]$ which is infinitely tangent to zero at the endpoints. 
Consider a diffeomorphism $f$ of $[0,1[$ topologically 
contracting towards the origin such that $f(b) = a$ and such that, if we 
extend $X$ to $[0,1[$ by imposing the condition $f_{*}(X) = \frac{1}{\lambda} X$, 
then the resulting $X$ is of class $\ce^{\infty}$ on $[0,1[$ (to ensure the 
existence of such a vector fied, one can either give an 
explicit construction or use Exercise \ref{reg-no-facil}). Let 
$\{\varphi^t \!: t \in \mathbb{R} \}$ be the flow associated to $X$. 
For $i \in \{1,\ldots,n\}$ let $f_i$ be the diffeomorphism $\varphi^{t_i}$. 
We claim that the group $\Gamma$ generated by $f,f_1,\ldots,f_n$ is polycyclic. 
Indeed, from 
$$f^{-1} f_i f \esp = \esp \varphi^{\lambda t_i} \esp = \esp 
\varphi^{\sum_j \! a_{i,j} t_j} \esp = \esp \prod_{j=1}^n (\varphi^{t_j})^{a_{i,j}} 
\esp = \esp \prod_{j = 1}^{n} f_j^{a_{i,j}}$$
one easily concludes that the Abelian subgroup generated by $f_1,\ldots,f_n$ 
is normal. Moreover, $\Gamma'$ coincides with the subgroup consisting of 
the elements of the form $f_1^{m_1} f_2^{m_2} \cdots f_n^{m_n}$, where 
$(m_1,\ldots.m_n)$ belongs to \esp $(A - Id) (\mathbb{Z}^n)$. \esp 
Therefore, $\Gamma$ is polycyclic.

\vspace{0.1cm}

\begin{small}\begin{ejer} Given $k \!\geq\! 1$, let $f \in \mathrm{Diff}^{k+1}_+([0,1[)$ be 
a diffeomorphism topologically contracting towards the origin and such that $f'(0) = 1$. 
Let $X$ be a continuous vector field on $[0,1[$ which is of class $\ce^k$ on $]0,1[$. Show 
that if $f_*(X) = \frac{1}{\lambda} X$ for some $\lambda \! \in ]0,1[$, then $X$ is of 
class $\ce^k$ on $[0,1[$ (see \cite{moriyama} in case of problems). 
\label{reg-no-facil}
\end{ejer}\end{small}

\vspace{0.08cm}

To state Moriyama's theorem properly, recall that to every polycyclic group $\Gamma$ there is an 
associated {\bf{\em nilradical}} $N(\Gamma)$ which corresponds to the maximal nilpotent normal subgroup 
of $\Gamma$ (see Appendix A). For example, in the preceding example, $N(\Gamma)$ coincides with the 
Abelian group generated by $f_1,\ldots,f_n$. Moreover, though in general $N(\Gamma)$ does not 
necessarily contain the derived subgroup of $\Gamma$, there exists a finite index 
subgroup $\Gamma_0$ of $\Gamma$ such that $[\Gamma_0,\Gamma_0] \subset N(\Gamma)$. 
\index{nilradical}

\vspace{0.1cm}

\begin{thm} {\em Let $\Gamma$ be a polycyclic subgroup of \esp 
$\mathrm{Diff}_+^{1+\mathrm{bv}}([0,1[)$ having no global fixed 
point in $]0,1[$. The following two possibilities may occur:}

\vspace{0.05cm}

\noindent (i) {\em if $N(\Gamma)$ has no global fixed point in $]0,1[$, 
then $\Gamma$ is conjugate to a subgroup of the affine group;}

\vspace{0.05cm}

\noindent (ii) {\em if $N(\Gamma)$ has a global fixed point in $]0,1[$, then $\Gamma$ 
is a semidirect product $\mathbb{Z}^n \rtimes \mathbb{Z}$.}
\label{thm-moriyama}
\end{thm}

\vspace{0.1cm}

Actually, in case (ii) the dynamics can be fully described, and 
it is very similar to that of the second example previously exhibited.

Here we just show that polycyclic subgroups of $\mathrm{Diff}_+^{1+\mathrm{bv}}([0,1[)$ are 
metabelian. (Theorem \ref{thm-moriyama} will then follow from the classification of solvable 
subgroups of $\mathrm{Diff}_+^{1+\mathrm{bv}}([0,1[)$: see \S \ref{seccion-metabeliano}.) 
Actually, this is quite easy. Indeed, let $\Gamma$ be a non-Abelian polycyclic subgroup of 
$\mathrm{Diff}_+^{1+\mathrm{bv}}([0,1[)$. By Plante-Thurston's theorem, the nilradical $N(\Gamma)$ 
is Abelian (and its restriction to each of its irreducible components is conjugate to a 
group of translations). Hence, if $\Gamma_0$ is a finite index subgroup of $\Gamma$ such 
that $[\Gamma_0,\Gamma_0] \subset N(\Gamma)$, then $\Gamma_0$ is metabelian. The fact 
that $\Gamma$ itself is metabelian then follows from Exercise \ref{indice-finito}.
\index{Plante-Thurston's theorem}

We conclude this section with a clever result due to Matsuda \cite{matsuda}. 
\index{Matsuda} 
It turns out that, in most 
of the cases, the algebraic properties of subgroups of $\mathrm{Diff}_+^2 ([0,1[)$ and 
$\mathrm{Diff}_+^2 ([0,1])$ are similar. However, in what concerns polycyclic subgroups, 
the situation is rather special. We ignore whether the result below extends to 
the class $\ce^{1+\mathrm{bv}}$.  

\vspace{0.1cm}

\begin{thm} {\em Every polycyclic subgroup of $\mathrm{Diff}_+^{1+\mathrm{Lip}}([0,1])$ 
without global fixed points in $]0,1[$ is topologically conjugate to a subgroup of the 
affine group.}
\label{thm-matsuda}
\end{thm}

\vspace{0.1cm}

The proof of this theorem strongly uses Thurston's stability theorem. 
The discussion is thus postponed to \S \ref{sec-th-th}.

\index{group!polycyclic|)}


\section{Solvable Groups of Diffeomorphisms}
\label{solvito}

\subsection{Some examples and statements of results}

\hspace{0.45cm} The discussion in the preceding section turns natural the question of 
knowing whether every solvable subgroup of $\mathrm{Diff}_{+}^{1+\mathrm{bv}}([0,1[)$ 
is metabelian. However, the examples below (inspired by \cite[Chapitre V]{godbillon}) 
show that there exist solvable groups of diffeomorphisms of the interval with arbitrary 
solvability degree having dynamics very different from that of the affine group. For 
this, it suffices to take successive extensions by $\mathbb{Z}$ in an appropriate way.

\begin{small} \begin{ejem} Let $f$ be a $\mathrm{C}^{\infty}$ diffeomorphism of $[0,1]$ 
without global fixed point in $]0,1[$ and topologically contracting towards $0$. 
Suppose that $f$ is the time 1 of the flow associated to a $\mathrm{C}^{\infty}$ vector 
field $[0,1]$ which is infinitely flat at the endpoints. Fix $a \! \in ]0,1[$, and 
consider a vector field $X$ on $[f(a),a]$ with zeros only at $f(a)$ and $a$, and 
which is infinitely flat at these points. Extend $X$ by letting $X(x)=0$ for 
$x \in [0,1] \setminus [f(a),a]$, thus obtaining a $\mathrm{C}^{\infty}$ 
vector field on $[0,1]$ which is infinitely flat at $0$ and $1$.

Let $g$ by the $\mathrm{C}^{\infty}$ diffeomorphism obtained by integrating the vector field $X$ 
(up to time $1$), and let $\Gamma$ the group generated by $f$ and $g$. Clearly, $\Gamma$ has 
no fixed point in $]0,1[$. We claim that $\Gamma$ is solvable with order of solvability 
equal to $2$. Indeed, let $\Gamma^{*}$ be the Abelian subgroup of $\Gamma$ formed by the elements 
fixing the points $f^n(a)$ (with $n \in \mathbb{Z}$), so that their restrictions to each interval 
$[f^{n+1}(a),f^n(a)]$ are contained in the group generated by the restriction of the element 
$f^ngf^{-n}$ to this interval. Clearly, $\Gamma^{*}$ is an Abelian normal subgroup of $\Gamma$ 
containing the derived group $\Gamma'$; moreover, the quotient $\Gamma / \Gamma^{*}$ identifies 
with $(\mathbb{Z},+)$, thus showing the claim.

In what follows, we will successively apply the preceding idea to obtain, for each 
$k \geq 2$, a solvable group $\bar{\Gamma}_k$ of $\ce^{\infty}$ diffeomorphisms 
of $[2-k,k-1]$ with order of solvability equal to $k$ and generated by elements 
$f_{1,k}, \ldots, f_{k,k}$ each of which is the time 1 of a flow associated to a 
$\ce^{\infty}$ vector field that is infinitely flat at $2-k$ and $k-1$. To do this, 
we argue inductively. For $k=2$ we let $\bar{\Gamma}_2 = \Gamma$, where $f_{1,2} = g$  
and $f_{2,2} = f$. Suppose now that we have already constructed the group $\Gamma_{k}$, 
and let us consider the vector fields $X_{i,k}\!: [2-k,k-1] \rightarrow \mathbb{R}$ 
corresponding to the $f_{i,k}$'s. Let us consider a vector field $X_{k+1,k+1}$ on 
$[1-k,k]$ which is of the form $\varrho \frac{\partial}{\partial x}$ for some $\ce^{\infty}$ 
function $\varrho\!: [1-k,k] \rightarrow \mathbb{R}$ that is negative at the interior and 
infinitely flat at the endpoints. After multiplying this vector field by a scalar factor  
if necessary, we may suppose that the time 1 of the associated flow is a diffeomorphism 
$f_{k+1,k+1}$ of $[1-k,k]$ satisfying $f_{k+1,k+1}(k-1) = 2 - k$.

For $i \in \{1,\ldots,k\}$ we extend $X_{i,k}$ to a vector field $X_{i,k+1}$ by letting 
$X_{i,k+1}(x) = 0$ for every point $x \in [1-k,k] \setminus [2-k,k-1]$. The vector fields 
thus obtained are of class $\mathrm{C}^{\infty}$ on $[1-k,k]$ and infinitely flat at $1-k$ 
and $k$. Let $f_{i,k+1}$, $i \in \{1, \ldots,k\}$, be the time $1$ of the associated 
flows. We claim that the group $\bar{\Gamma}_{k+1}$ generated by the $f_{i,k+1}$'s, 
$i \in \{1, \ldots,k+1\}$, is solvable with solvability degree $k+1$. Indeed, the 
stabilizer in $\Gamma$ of $[2-k,k-1]$ is a normal subgroup $\bar{\Gamma}^{*}$ which 
identifies with a direct sum of groups isomorphic to $\bar{\Gamma}_k$, and by the induction 
hypothesis the latter group is solvable with solvability degree $k$. On the other hand, 
$\bar{\Gamma}_{k+1}/\bar{\Gamma}^{*}$ identifies with $(\mathbb{Z},+)$, and the derived 
group of $\Gamma$ is contained in $\bar{\Gamma}^{*}$. Using this fact, one easily concludes 
the claim. Notice that $\bar{\Gamma}_{k+1}$ has no fixed point in the interior of $[1-k,k]$.
\label{one}
\end{ejem}

\begin{ejem} We next improve the first step of the preceding example using the 
construction of the beginning of \S \ref{policiclicos-suaves}. We will thus obtain 
a more interesting family of solvable subgroups of $\mathrm{Diff}_{+}^{\infty}([0,1])$ 
with solvability degree 3. For this, let us consider two $\ce^{\infty}$ vector fields 
$X_1$ and $X_2$ defined on the interval $[1/3,2/3]$ which are infinitely flat at the 
endpoints and whose flows induce a conjugate of the affine group. Let us denote by 
$g$ and $h$ the $\mathrm{C}^{\infty}$ diffeomorphisms obtained by integrating up 
to time 1 the vector fields $X_1$ and $X_2$, respectively.

Fix a $\ce^{\infty}$ diffeomorphism $f\!: [0,1] \rightarrow [0,1]$ without fixed point at the 
interior and satisfying $f(2/3) = 1/3$. For $a=2/3$ let us consider a sequence (to be fixed) 
of positive real numbers $(t_n)_{n \in \mathbb{Z}}$, and let us extend by induction the 
definition of the $X_j$ by letting
$$X_j(x)= t_n X_j(f^{-1}(x)) \big/ (f^{-1})'(x), \qquad x \in [f^{n+1}(a),f^n(a)], \quad
n \geq 1,$$
$$X_j(x)= t_n X_j(f(x)) \big/ f'(x), \qquad x \in [f^{n+1}(a),f^n(a)], \quad n \leq -1.$$
Suppose that $\Pi_{i=1}^{n} t_i \!\rightarrow\! 0$ as $n \!\rightarrow\! +\infty$,  
that $\Pi_{i=n}^{0} t_i \! \rightarrow \! 0$ as $n \!\rightarrow\! -\infty$, and 
that these convergences are very fast (to fix the ideas, let us suppose that the speed  
of convergence is super-exponential; 
compare Exercise \ref{kopell-fuerte-lip}). In this case, it is not difficult to check that 
the $X_j$'s extend to $\mathrm{C}^{\infty}$ vector fields on $[0,1]$ which are zero and 
infinitely flat at the endpoints. Let $\Gamma$ be the group generated by $f$, $g$, and 
$h$. Obviously, the restriction of $\Gamma$ to the interior of $[0,1]$ is not semiconjugate 
to a subgroup of the affine group, and $\Gamma$ has no fixed point in $]0,1[$. We claim 
that $\Gamma$ is solvable of degree of solvability $3$. Indeed, let $\Gamma^{*}$ 
be the metabelian subgroup of $\Gamma$ formed by the elements fixing 
the points $f^n(a)$, so that their restrictions to the interior of each interval 
$[f^{n+1}(a),f^n(a)]$ are contained in the conjugate of the affine group generated 
by $X_1$ and $X_2$. If we respectively denote by $g^t_{[f^{n+1}(a),f^n(a)]}$ and 
$h^t_{[f^{n+1}(a),f^n(a)]}$ the flows associated to the restrictions of $X_1$ and $X_2$ to 
$[f^{n+1}(a),f^n(a)]$ (with $t \in \mathbb{R}$), then for every $n \in \mathbb{N}$ one has
$$f^{-1} \circ g_{[f^{n+1}(a),f^n(a)]} \circ f \!=\! g^{t_n}_{[f^n(a),f^{n-1}(a)]}, 
\quad f^{-1} \circ h_{[f^{n+1}(a),f^n(a)]} \circ f \!=\! h^{t_n}_{[f^n(a),f^{n-1}(a)]}.$$
The subgroup $\Gamma^{*}$ is therefore normal in $\Gamma$. Moreover, the quotient 
$\Gamma / \Gamma^{*}$ identifies with $(\mathbb{Z},+)$, and $\Gamma^*$ contains the derived 
group of $\Gamma$. Starting from this one easily concludes our claim.

As in the preceding example, for each integer $k \geq 2$ one may take successive 
extensions of $\Gamma$ by $(\mathbb{Z},+)$ for obtaining solvable subgroups of 
$\mathrm{Diff}_{+}^{\infty}([1-k,k])$ with solvability degree $k+2$.
\label{two}
\end{ejem} 


\begin{ejer} Prove that, in the second of the above examples, 
the convergences \esp $\Pi_{i=1}^{n} t_i \rightarrow 0$ \esp and \esp 
$\Pi_{i=n}^{0} t_i \rightarrow 0$ \esp are necessary. More precisely, 
prove the following proposition.

\vspace{0.02cm}

\begin{prop} {\em Let $g \!:\![0,1[ \rightarrow [0,g(1)[$ be a $\ce^2$ diffeomorphism, 
and let $a<b$ be two fixed points of $g$ in $]0,1[$ so that $g$ has no fixed point in 
$]a,b[$. Let $f\!:\![0,b[ \rightarrow [0,f(b)[$ be a $\mathrm{C}^2$ diffeomorphism with 
no other fixed point on $[0,b]$ than $0$, and such that $f(b) \leq a$. Suppose that 
there exists a sequence $(t_n)_{n \in \mathbb{N}}$ of positive real numbers so that 
for all $n \!\in\! \mathbb{N}$ one has 
$f^{-1} \circ g_{[f^n(a),f^n(b)[} \circ f = g_{[f^{n-1}(a),f^{n-1}(b)[}^{t_n}$.
Then the value of $\Pi_{i=1}^n t_i$ converges to zero as $n$ goes to infinity.}
\label{casikopell}
\end{prop}
\end{ejer} \end{small}

\vspace{0.12cm}

The method of construction of Examples \ref{one} and \ref{two} is essentially 
the only possible one for creating solvable groups of interval diffeomorphisms. 
Let us begin by giving a precise version of this fact in the metabelian case.

\vspace{0.08cm}

\begin{thm} {\em If $\Gamma$ is a metabelian subgroup of 
$\mathrm{Diff}_{+}^{1+\mathrm{bv}}([0,1[)$ without fixed points in $]0,1[$, then $\Gamma$ 
is either conjugate to a subgroup of the affine group, or a semidirect product between 
$(\mathbb{Z},+)$ and a subgroup of a product (at most countable) of groups which are 
conjugate to groups of translations}.
\label{clas-meta}
\end{thm}

\vspace{0.08cm}

We will give the proof of this theorem in the next section. A complete classification 
of solvable subgroups of $\mathrm{Diff}_{+}^{1+\mathrm{bv}}([0,1[)$ may be obtained  
by using the same ideas throug a straightforward inductive argument \cite{Nasol}. To 
state the theorem in the general case, let us denote by $r(1)$ the family of groups that are 
conjugate to groups of translations, and by $r(2)$ the family of groups that are either conjugate 
to non-Abelian subgroups of the affine group, or a semidirect product between $(\mathbb{Z},+)$ 
and a subgroup of an at most countable product of nontrivial groups of translations. For $k > 2$ 
we define by induction the family $r(k)$ formed by the groups that are a semidirect product 
between $(\mathbb{Z},+)$ and a subgroup of an at most countable product of groups in 
$\mathcal{R}(k-1) = r(1) \! \cup \! \ldots \! \cup \! r(k-1)$, so that at least one 
of the factors does not belong to $\mathcal{R}(k-2)$.

\vspace{0.08cm}

\begin{thm} {\em Let $\Gamma$ be a solvable subgroup of \esp 
$\mathrm{Diff}_{+}^{1+\mathrm{bv}}([0,1[)$ without fixed point 
in $]0,1[$. If the solvability degree of $\Gamma$ equals 
$k \geq 2$, then $\Gamma$ belongs to the family $r(k)$.}
\label{clas-gen}
\end{thm}

\vspace{0.08cm}

The preceding classification allows obtaining interesting rigidity results. For instance, 
the normalizer of a solvable group of diffeomorphisms of the interval is very similar 
to the original group, as is established in the next result from \cite{Na-subexp}.

\vspace{0.03cm}

\begin{thm} {\em Let $\Gamma$ be a solvable subgroup of \esp 
$\mathrm{Diff}_{+}^{1+\mathrm{bv}}([0,1[)$ of solvability 
degree $k \geq 1$ and without fixed point in $]0,1[$. If 
$\mathcal{N}(\Gamma)$ denotes its normalizer in 
$\mathrm{Diff}_+^1([0,1[)$, then one of the 
following possibilities occurs:}

\vspace{0.07cm}

\noindent{(i) {\em if $k > 1$ then $\mathcal{N}(\Gamma)$ is solvable with solvability degree $k$;}}

\vspace{0.07cm}

\noindent{(ii) {\em if $k = 1$ and $\Gamma$ is not infinite cyclic, then $\mathcal{N}(\Gamma)$ 
is topologically conjugate to a (perhaps non-Abelian) subgroup of the affine group;}}

\vspace{0.07cm}

\noindent{(iii) {\em if $k = 1$ and $\Gamma$ is infinite cyclic, then $\mathcal{N}(\Gamma)$ 
is topologically conjugate to a subgroup of the group of translations.}}
\end{thm}

\vspace{0.03cm}

The classification of solvable groups of circle diffeomorphisms reduces, thanks 
to Lemma \ref{moyen-medida}, to that of the case of the interval. Using 
Theorem \ref{clas-gen}, the reader should easily verify the following.

\vspace{0.1cm}

\begin{thm} {\em Let $\Gamma$ be a solvable subgroup of 
$\mathrm{Diff}_+^{1+\mathrm{bv}}(\clo)$. If the degree 
of solvability of $\Gamma$ equals $k+1$, then one 
of the following possibilities occurs:}

\vspace{0.07cm}

\noindent{(i) {\em $\Gamma$ is topologically semiconjugate 
\index{rotation!group} to a group of rotations;}}

\vspace{0.07cm}

\noindent{(ii) {\em there exists a non-empty finite subset $F$ of $\clo$ which 
is invariant by $\Gamma$ so that the derived group $\Gamma'$ fixes each point 
of $F$, and the restriction of $\Gamma'$ to each one of the connected components 
of $\clo \setminus F$ belongs to the family $\mathcal{R}(k)$.}}
\end{thm}

\vspace{0.1cm}

The case of the line involves some extra difficulties. However, using some of 
the results from \S \ref{muylargo}, one may obtain the description below.

\vspace{0.1cm}

\begin{thm} {\em Let $\Gamma$ be a solvable subgroup of 
$\mathrm{Diff}_{+}^{1+\mathrm{bv}}(\mathbb{R})$. If the 
solvability degree of $\Gamma$ equals $k \geq 1$, then 
one of the following possibilities occurs:}

\vspace{0.07cm}

\noindent{(i) {\em $\Gamma$ is topologically semiconjugate 
to a subgroup of the affine group;}}

\vspace{0.07cm}

\noindent{(ii) {\em $\Gamma$ is a subgroup of a product (at most 
countable) of groups in the family $\mathcal{R}(k)$, so that at 
least one of the factors does not belong to $\mathcal{R}(k-1)$;}}

\vspace{0.07cm}

\noindent{(iii) {\em $\Gamma$ belongs either to $r(k)$ or $r(k+1)$.}}
\label{solrecta}
\end{thm}

\vspace{0.1cm}

Notice that if $\Gamma$ is solvable and semiconjugate to a group of affine transformations 
without being conjugate to it, then the second derived group $\Gamma''$ acts fixing a 
countable family of disjoint open intervals whose union is dense. Theorems \ref{clas-meta} 
and \ref{clas-gen} then allow describing the dynamics of $\Gamma''$. On the other hand, the 
fact that $\Gamma$ may belong to $r(k+1)$ when its solvability degree is $k$ is quite natural, 
since unlike the case of the interval, in the case of the line it is possible to produce 
{\em central} extensions of nontrivial groups, even in class $\ce^{\infty}$.

\vspace{0.2cm}


\subsection{The metabelian case}
\label{seccion-metabeliano}
\index{group!metabelian}

\vspace{0.2cm}

\hspace{0.45cm} The elementary lemma below will play a fundamental role in what follows.

\vspace{0.1cm}

\begin{lem} {\em The normalizer in $\mathrm{Homeo}_{+}(\mathbb{R})$ of every 
dense subgroup of the group of translations is contained in the affine group.}
\label{normalizando}
\end{lem}

\noindent{\bf Proof.} Up to a scalar factor, the only Radon measure on the line that 
is invariant by a dense group of translations is the Lebesgue measure. The normalizer 
\index{measure!quasi-invariant} \index{measure!Radon}of such a group 
leaves quasi-invariant this measure, and hence the conclusion of the 
lemma follows from Proposition \ref{conjugaralafin}. $\hfill\square$

\vspace{0.5cm}

The next lemma is an improved version of the preceding 
one in class $\mathrm{C}^{1+\mathrm{bv}}$.

\vspace{0.1cm}

\begin{lem} {\em If \esp $\Gamma$ is an Abelian subgroup of 
$\mathrm{Diff}_+^{1+\mathrm{bv}}([0,1[)$ without fixed points 
in $]0,1[$, then its normalizer in $\mathrm{Diff}_{+}^1([0,1[)$ 
is conjugate to a subgroup of the affine group.}
\label{normalizar-en-cedos}
\end{lem}

\noindent{\bf Proof.} Let $g \! \in \! \Gamma$ be a nontrivial element. 
By Corollary \ref{flujocelog}, the hypothesis of non-existence of fixed 
points in $]0,1[$ is equivalent to that $\Gamma$ is contained in a 1-parameter 
group $g^{\mathbb{R}} = \{g^t \!: t \in \mathbb{R} \}$ which is conjugate to the 
group of translations and such that $g^1 = g$; moreover, the centralizer in 
$\mathrm{Diff}_+^1([0,1[)$ of every nontrivial element of $\Gamma$ is contained 
in $g^{\mathbb{R}}$. If the 
image of $\Gamma$ by the conjugacy is a dense subgroup of $(\mathbb{R},+)$, then 
Lemma \ref{normalizando} implies that 
the image by this conjugacy of the normalizer $\mathcal{N}$ 
of $\Gamma$ in $\mathrm{Homeo}_+([0,1])$ is contained in the affine group. 
On the other hand, we claim that 
if $\{ t \in \mathbb{R} \!: g^t \in \Gamma \}$ is infinite cyclic, then 
$\mathcal{N}$ equals $g^{\mathbb{R}}$. Indeed, if $k$ 
is a positive integer such that $g^{1/k}$ is the generator of $\{ t \in \mathbb{R} \!: 
g^t \in \Gamma \}$, then for every $h \!\in\! \mathcal{N}$ there exists positive integers $n,m$ 
such that $hg^{1/k}h^{-1} \!=\! (g^{1/k})^n$ and $h^{-1}g^{1/k} h \!=\! (g^{1/k})^m$. One then has 
$$(g^{1/k})^{mn} = ((g^{1/k})^{m})^n =
(h^{-1}g^{1/k}h)^n = h^{-1}(g^{1/k})^{n}h = g^{1/k},$$
from where one obtains $m \!=\! n \!=\! 1$. This implies that the elements in $\mathcal{N}$ commute 
with $g^{1/k}$, and hence $\mathcal{N}$ is contained in $g^{\mathbb{R}}$. Since $g^{\mathbb{R}}$ 
centralizes (and hence normalizes) $\Gamma$, this shows that 
$\mathcal{N} = g^{\mathbb{R}}$. $\hfill\square$

\vspace{0.4cm}

We may now pass to the proof of Theorem \ref{clas-meta}. By Corollary \ref{flujocelog}, 
if $\Gamma$ is an Abelian subgroup of $\mathrm{Diff}_+^{1+\mathrm{bv}}([0,1[)$ without 
fixed points in $]0,1[$, then it is conjugate to a group of translations. Now let 
$\Gamma$ be a metabelian and non-commutative subgroup of 
$\mathrm{Diff}_+^{1+\mathrm{bv}}([0,1[)$ without fixed points in $]0,1[$. If there exists 
$g \in \Gamma'$ such that the orbits by $g$ accumulate at $0$ and $1$, then the Abelian 
group $\Gamma'$ is contained in the topological flow associated to $g$. Therefore, 
$\Gamma'$ acts without fixed point in $]0,1[$, and the preceding lemma implies 
that $\Gamma$ is conjugate to a (non commutative) subgroup of the affine group.

In what follows, suppose that every $g \in \Gamma'$ has fixed points in $]0,1[$. Kopell's  
lemma then easily implies that $\Gamma'$ must have global fixed points in $]0,1[$. 
\index{irreducible component} 
Let $]a,b[$ be an irreducible component of $\Gamma'$. Notice that, for every $h \in \Gamma$, 
the interval $h(]a,b[)$ is also an irreducible component of $\Gamma'$. In particular, if 
$h(]a,b[) \neq ]a,b[$, then $h(]a,b[) \cap ]a,b[ = \emptyset$.

If $a \!=\! 0$ or $b \!=\! 1$, then every $f \!\in\! \Gamma$ fixes $]a,b[$, which contradicts the 
hypothesis of non-existence of fixed points in $]0,1[$. Therefore, $[a,b]$ is contained in $]0,1[$. 
If $f \!\in\! \Gamma$ fixes $]a,b[$, then Lemma \ref{normalizar-en-cedos} shows that the 
restriction of $f$ to $]a,b[$ is affine in the coordinates induced by $\Gamma'$. The 
case of the elements that do not fix $]a,b[$ is more interesting.

\vspace{0.25cm}

\noindent{\underbar{Claim (i).}} \esp If $f$ is an element of $\Gamma$ that does not fix $]a,b[$, 
and if $u$ and $v$ are the fixed points of $f$ to the left and to the right of $]a,b[$, respectively,
then the interval $]u,v[$ is an irreducible component of $\Gamma$ (that is, $u = 0$ and $v = 1$).

\vspace{0.15cm}

Suppose not, and let $\bar{f} \in \Gamma$ an element which does not fix $]u,v[$. Replacing 
$f$ by $f^{-1}$ if necessary, we may suppose that $f(x) > x$ for all $x \! \in ]u,v[$. One then 
has $f(a) \geq b$. For $n \!\in\! \mathbb{N}$ the element $f^{-1}\bar{f}^{-n}f\bar{f}^n$ 
belongs to $\Gamma'$ and, therefore, fixes the points $a$ and $u$. Therefore, for 
every $n \in \mathbb{N}$,
\begin{equation}
f \bar{f}^n(u) = \bar{f}^n f(u) = \bar{f}^n(u), \qquad
f\bar{f}^n(a)=\bar{f}^nf(a) \geq \bar{f}^n(b).
\label{paraellema}
\end{equation}
One has $\bar{f}f\bar{f}^{-1} = f\bar{g}$ for some $\bar{g} \in \Gamma'$. On the other hand, 
$f\bar{g}$ has no fixed point in $]u,v[$ and fixes $u$ and $v$. Therefore, $f$ does not 
have fixed points in $]\bar{f}(u),\bar{f}(v)[$, and fixes $\bar{f}(u)$ and $\bar{f}(v)$. 
This shows that $]\bar{f}(u),\bar{f}(v)[ \cap ]u,v[ = \emptyset$. Replace $\bar{f}$ 
by $\bar{f}^{-1}$ if necessary so that $\bar{f}(u) < u$. Notice that the sequence 
$(\bar{f}^n(u))$ tends to a fixed point of $\bar{f}$. One then easily checks that 
relations (\ref{paraellema}) contradict Lemma \ref{kopelcito} (applied to the 
elements $h_1=\bar{f}$ and $h_2=f$ with respect to the points $x_0 \!=\! \bar{f}^{-1}(u)$, 
$y \!=\! a$, and $z \!=\! b \!\leq\! f(a)$). This concludes the proof of the claim.

\vspace{0.1cm}

Denote by $\Gamma^*$ the normal subgroup of $\Gamma$ formed by the elements fixing the 
irreducible components of $\Gamma'$. Since Claim (i) holds for every irreducible component 
$]a,b[$ of $\Gamma'$, an element of $\Gamma$ belongs to $\Gamma^*$ if and only if it fixes 
at least one of the irreducible components of $\Gamma'$. The restriction of $\Gamma^*$ to 
each irreducible component of $\Gamma'$ is affine in the induced coordinates. Remark that 
$\Gamma^*$ may admit irreducible components contained in the complement of the union of 
the irreducible components of $\Gamma'$. However, the restriction of $\Gamma^*$ to such 
a component is Abelian, and hence conjugate to a subgroup of the group of translations. We 
then conclude that $\Gamma^*$ is a subgroup of a product (at most countable) of groups 
conjugate to groups of affine transformations. Moreover, the quotient group 
$\Gamma/\Gamma^*$ acts freely 
\index{action!free} 
on the set $\mathrm{Fix}(\Gamma')$. Fix an irreducible component $]a,b[$ 
of $\Gamma'$, and define an order relation $\preceq$ on $\Gamma/\Gamma^*$ 
by $f_1\Gamma^* \prec f_2\Gamma^*$ when $f_1(]a,b[)$ is to the left of 
$f_2(]a,b[)$. This relation is total, bi-invariant and Archimedean. 
\index{Archimedean property} 
The argument of the proof of H\"older's theorem 
\index{H\"older!theorem} 
then shows that the group $H = \Gamma / \Gamma^*$ is naturally isomorphic to a 
subgroup of $(\mathbb{R},+)$. Notice that $H$ is nontrivial, since $\Gamma$ does 
not fix $]a,b[$. Proposition \ref{no-semiconjugado} then implies the claim below.

\vspace{0.25cm}

\noindent{\underbar{Claim (ii).}} \esp $H$ is an infinite cyclic group.

\vspace{0.2cm}

The proof of Theorem \ref{clas-meta} is then finished by the claim below.

\vspace{0.25cm}

\noindent{\underbar{Claim (iii).}} \esp $\Gamma^*$ is a subgroup of a product 
of groups that are conjugate to groups of translations.

\vspace{0.2cm}

To show this, fix an element $g \!\in\! \Gamma'$ whose restriction to an irreducible 
component $]a,b[$ of $\Gamma'$ has no fixed point. Assume for a contradiction that 
there exists $h \!\in\! \Gamma$ fixing $]a,b[$ so that the restrictions of $g$ and 
$h$ to this interval generate a non-Abelian group. Without loss of generality, 
we may suppose that the fixed point $a$ of $h$ is topologically repelling 
by the right. Let $f \!\in\! \Gamma$ be an element such that $f \Gamma^*$ generates 
$\Gamma / \Gamma^*$ and  $f(b) \!\leq \! a$. For every $n \!\in\! \mathbb{N}$ the 
element $f^{-n}hf^n$ fixes the interval $]a,b[$, and hence its restriction therein equals 
the restriction of $hg^{t_n}$ for some $t_n \!\in \!\mathbb{R}$, where $g^{\mathbb{R}}$ 
stands for the flow associated to the restriction of $g$ to $[a,b[$. For 
$\delta \!=\! V(f;[0,b]) \!>\! 0$, inequality (\ref{clasiquita})
allows showing that, for every $x \! \in ]a,b[$,  
\begin{small}
\begin{equation}
(hg^{t_n})'(x) = (f^{-n} h f^n)'(x) \leq 
\frac{(f^n)'(x)}{(f^n)'(f^{-n}hf^n(x))} \! \sup_{y \! \in ]0,f^n(b)[} h'(y)
\leq e^{\delta} \!\!\! \sup_{y \! \in ]0,f^n(b)[} h'(y),
\label{estimar en ce cero}
\end{equation}\end{small}and
\begin{equation}
(hg^{t_n})'(x) \esp \geq \esp 
e^{- \delta} \inf_{y \! \in ]0,f^n(b)[} h'(y).
\label{estt}
\end{equation}
From (\ref{estimar en ce cero}) one concludes that 
$$\sup_{x \! \in ]a,b[} (g^{t_n})'(x) \esp \leq \esp e^{\delta} \cdot
\frac{\sup_{y \! \in ]0,b[} h'(y)}{\inf_{y \! \in ]a,b[} h'(y)}.$$
This clearly implies that $(|t_n|)$ is a bounded sequence. Take a subsequence $(t_{n_k})$ 
converging to a limit $T \!\in\! \mathbb{R}$. Since $h$ fixes each interval $[f^n(a),f^n(b)]$, 
one has $h'(0)\!=\!1$. By integrating (\ref{estimar en ce cero}) and (\ref{estt}) we obtain, 
for every $k$ large enough and every $x \! \in ]a,b[$,
$$(x-a)/2e^{\delta} \esp \leq \esp 
hg^{t_{n_k}}(x) - a  \esp \leq \esp 2e^{\delta} (x-a),$$
and passing to the limit as $k$ goes to infinity we conclude that
$$(x-a)/2e^{\delta} \esp \leq \esp 
hg^{T}(x) - a \esp \leq \esp 2e^{\delta} (x-a).$$
Now notice that the argument above still holds when replacing $h$ by $h^j$ for 
any $j \!\in\! \mathbb{N}$ (indeed, the constant $\delta$ depends only on $f$). 
Therefore, for every $x \! \in ]a,b[$ and all $j \in \mathbb{N}$,
$$(x-a)/2e^{\delta} \esp \leq 
\esp (hg^{T})^j(x) - a \esp \leq \esp 2e^{\delta} (x-a),$$
which is impossible for $x$ near $a$ since the latter point is fixed by 
$hg^T$ and topologically repelling by the right. This concludes the proof. 

\begin{small}\begin{ejer} Show that if $\Gamma$ is a subgroup of 
$\mathrm{Diff}_+^{1+\mathrm{bv}} ([0,1[)$ containing a finite index 
subgroup $\Gamma_0$ which is metabelian (resp. solvable with solvability 
degree $k$), then $\Gamma$ itself is metabelian (resp. solvable with 
solvability degree $k$) (compare Exercise \ref{virt-abel}). 
\label{indice-finito}
\end{ejer}

\begin{ejer} Show that every polycyclic subgroup of 
$\mathrm{Diff}_+^{1+\mathrm{bv}}([0,1[)$ is strongly 
polycyclic ({\em c.f.}, Appendix A).
\end{ejer}\end{small}

It is worth noticing that if a solvable subgroup of $\mathrm{Diff}_+^{1+\mathrm{bv}}([0,1[)$ 
has solvability degree greater than $2$, then it necessarily contains nontrivial elements having 
infinitely many fixed points in every neighborhood of the origin; in particular, these elements 
cannot be real-analytic. As a consequence, every solvable group of real-analytic 
diffeomorphisms of the interval is topologically conjugate to a subgroup of the 
affine group (provided that there is no global fixed point at the interior). 
The reader may find this and many other related results in \cite{Nak} (see also \cite{wilk,Nasol}). 
Let us finally point out that the results already described for solvable subgroups of 
$\mathrm{Diff}_+^{1+\mathrm{bv}}([0,1[)$ still hold --and may be refined-- for subgroups 
of the group $\mathrm{PAff}_+([0,1[)$ of piecewise affine homeomorphisms 
\index{piecewise affine homeomorphism}
of the interval.

\begin{small} \begin{ejer} Prove that if $\Gamma$ is a nontrivial subgroup of 
$\mathrm{PAff}_+([0,1[)$ acting freely on $]0,1[$, then $\Gamma$ is infinite 
cyclic. Using this, prove that every finitely generated metabelian subgroup 
of $\mathrm{PAff}_+([0,1])$ is isomorphic to a semidirect product between 
$(\mathbb{Z},+)$ and a {\em direct sum} (at most countable) of groups isomorphic 
to $(\mathbb{Z},+)$ acting on two-by-two disjoint intervals. State and show 
a general result of classification for finitely generated solvable subgroups of 
$\mathrm{PAff}_+([0,1])$ (see \cite{Na-subexp} for more details; see also 
\cite{bleak}).
\end{ejer} \end{small}


\subsection{The case of the real line}

\hspace{0.45cm} For a solvable group of solvability degree $k$, we denote by 
$\{ id \} = \Gamma_k^{\sol} \sgn \ldots \sgn \Gamma_{0}^{\sol} = \Gamma$ its derived 
series, that is, $\Gamma_{i}^{\sol} = [\Gamma_{i-1}^{\sol},\Gamma_{i-1}^{\sol}]$  
for every $i \!\in\! \{ 1, \ldots, k\}$, with $\Gamma_{k-1}^{\sol} \neq \{id\}$. 
The result below should be compared with Theorem \ref{plantesoluble}.

\vspace{0.1cm}

\begin{prop} {\em Let $\Gamma$ be a solvable subgroup of $\mathrm{Homeo}_{+}(\mathbb{R})$ of 
solvability degree $k$. If there exists an index $i \leq k$ for which $\Gamma_i^{\mathrm{sol}}$ 
preserves a Radon measure $v_i$ such that the associated translation number homomorphism 
$\tau_{v_i}$ satisfies $\tau_{v_i}(\Gamma_i^{\sol}) \neq \{0\}$, then there exists a 
Radon measure on the line which is quasi-invariant by $\Gamma$.}
\label{ahorasi}
\end{prop}

\noindent{\bf Proof.} Let $j < k$ be the smallest index for which $\Gamma_j^{\sol}$ 
preserves a Radon measure, and let $v_j$ be this measure. Using the hypothesis and 
(\ref{soporteplante}), it is not difficult to check that $\tau_{v_j}(\Gamma_j^{\sol}) \neq \{0\}$. 
By Lemma \ref{lema-uno} one has $\kappa(\Gamma_{j-1}^{\sol}) \neq \{1\}$, and by Lemma \ref{lema-dos} 
this implies that $v_j$ is quasi-invariant by $\Gamma_{j-1}^{\sol}$. Lemma \ref{lema-tres} 
then shows that $v_j$ is quasi-invariant by $\Gamma_{j-2}^{\sol}$, $\Gamma_{j-3}^{\sol}$, 
$\ldots$ , and hence 
\index{measure!quasi-invariant} by $\Gamma$. $\hfill\square$

\vspace{0.25cm}

We now give the proof of Theorem \ref{solrecta} (assuming Theorem \ref{clas-gen}). Fix a solvable 
subgroup $\Gamma$ of $\mathrm{Diff}^{1+\mathrm{bv}}_{+}(\mathbb{R})$ with solvability degree 
$k \!\geq\! 2$. If $\Gamma$ has global fixed points, then Theorems \ref{clas-meta} and 
\ref{clas-gen} imply that $\Gamma$ is a subgroup of a product (at most countable) of 
groups in the family $\mathcal{R}(k)$ so that at least one of the factors does not 
belong to $\mathcal{R}(k\!-\!1)$. Assume throughout that $\Gamma$ has no 
global fixed point. We begin with a proposition which should be compared with 
the end of \S \ref{muylargo}.

\vspace{0.1cm}

\begin{prop} {\em Every solvable subgroup of 
$\mathrm{Diff}_+^{1+\mathrm{bv}}(\mathbb{R})$ 
leaves quasi-invariant a Radon measure on the line.}
\label{quasiinvariante}
\end{prop}

\noindent{\bf Proof.} Let us consider the smallest index $j$ 
for which $\Gamma_j^{\sol}$ has a global fixed point. Notice 
that, since $\Gamma$ has no fixed point, $j$ must be positive. 
There are two distinct cases.

\vspace{0.25cm}

\noindent{\underbar{First case.}} The index $j$ equals $k$.

\vspace{0.15cm}

We claim that $\Gamma_{k-1}^{\sol}$ preserves a Radon measure so that the translation 
number homomorphism is non identically zero. Indeed, since $\Gamma_{k-1}^{\sol}$ is 
Abelian, using Kopell's lemma one \index{Kopell!lemma} \index{translation number}easily 
shows the existence of elements in $\Gamma_{k-1}^{\sol}$ 
without fixed points, and hence Proposition \ref{abelianoplante} implies the existence 
of a Radon measure invariant by $\Gamma_{k-1}$ such that the translation number 
homomorphism is nontrivial. The existence of a Radon measure $v$ on the line which 
is quasi-invariant by $\Gamma$ then follows from Proposition \ref{ahorasi}.

\vspace{0.25cm}

\noindent{\underbar{Second case.}} The index $j$ is less than $k$.

\vspace{0.15cm}

Fix an irreducible component $]p_j,q_j[$ of $\Gamma_j^{\sol}$.

\vspace{0.1cm}

\noindent{\underbar{Claim (i).}} If $\bar{f}_1$ and $\bar{f}_2$ are elements of 
$\Gamma_{j-1}^{\sol}$ that have fixed points but do not fix $]p_j,q_j[$, then the 
fixed points in $\mathbb{R} \cup \{-\infty,+\infty\}$ of these elements which are 
next to $]p_j,q_j[$ by the left coincide. The same holds for the fixed points 
which are next by the right.

\vspace{0.1cm}

To show this, let $p$ and $q$ the fixed points of $\bar{f}_1$ to the left and to the 
right of $[p_j,q_j]$, respectively. Suppose that $\bar{f}_2$ does not fix $[p,q]$. For 
each $n \in \mathbb{Z}$ there exists $\bar{g} \in \Gamma_{j}^{\sol}$ so that 
$\bar{f}_2^n \bar{f}_1 \bar{f}_2^{-n} = \bar{f}_1 \bar{g}$. Since $\bar{f}_1 \bar{g}$ 
has no fixed point in $]p,q[$ and fixes $p$ and $q$, the element $\bar{f}_1$ fixes 
the interval $\bar{f}_2^n(]p,q[)$ and has no fixed point inside it. One then deduces 
that the intervals $\bar{f}_2^n(]p,q[)$ are two-by-two disjoint for $n \!\in\! \mathbb{Z}$. 
Changing $\bar{f}_2$ by its inverse if necessary, we may suppose that $\bar{f}_2^n(p_j)$ 
tends to a fixed point (in $\mathbb{R}$) of $\bar{f}_2$ as $n$ goes to the 
infinity. We then obtain a contradiction by applying the arguments of the 
proof of Claim (i) from \S \ref{seccion-metabeliano}.

Let us define the interval $[p_{j-1}^*,q_{j-1}^*]$ as being equal to $[p_j,q_j]$ if 
every element in $\Gamma_{j-1}^{\sol}$ which does not fix $[p_j,q_j]$ has no fixed 
point. Otherwise, let us choose an element $f_j \in \Gamma_{j-1}^{\sol}$ having 
fixed points but which does not fix $[p_j,q_j]$, and define $p_{j-1}^*$ and 
$q_{j-1}^*$ as the fixed points of $f_j$ to the left and to the right of 
$[p_j,q_j]$, respectively. Finally, denote by $\Gamma_{j-1}^*$ the stabilizer of 
$[p_{j-1}^*,q_{j-1}^*]$ in $\Gamma_{j-1}$. The group $\Gamma_{j-1}^*$ is normal 
in $\Gamma_{j-1}$, because it is formed by the elements in $\Gamma_{j-1}$ 
having fixed points in the line.

\vspace{0.15cm}

\noindent{\underbar{Claim (ii).}} the group $\Gamma_{j-1}^{\sol}$ 
preserves a Radon measure and has no fixed point.

\vspace{0.1cm}

Indeed, since Claim (i) holds for every irreducible component of $\Gamma_j$, the action 
of $\Gamma_{j-1}^{\sol} / \Gamma_{j-1}^*$ on $\mathrm{Fix}(\Gamma_{j-1}^{*})$ is free. 
The argument of proof of H\"older's theorem then shows that $\Gamma_{j -1}^{\sol}$ fixes 
a Radon measure whose support is contained in $\mathrm{Fix}(\Gamma_{j-1}^{*})$. 

Recall that, by the definition of $j$, the subgroup $\Gamma_{j-1}^{\sol}$ has no fixed 
point (which together with (i) implies the existence of elements therein without fixed points 
points). The translation number function with respect to the $\Gamma_{j-1}^{\sol}$-invariant  
Radon measure is therefore non identically zero. Proposition \ref{ahorasi} then allows 
concluding the proof of the proposition in the second case. $\hfill\square$

\vspace{0.35cm}

We can now conclude the proof of Theorem \ref{solrecta}. For this, fix a Radon measure 
$v$ which is quasi-invariant by $\Gamma$. Suppose first that $v$ has no atom. Consider 
the equivalence relation $\sim$ which identifies two points if they belong to the closure 
of the same connected component of the complement of the support of $v$. The quotient space 
$\mathbb{R} /\!\!\sim$ is a topological line upon which the group $\Gamma$ acts naturally 
by homeomorphisms. Since $v$ has no atom, it induces a $\Gamma$-quasi-invariant non-atomic 
Radon measure $\overline{v}$ on $\mathbb{R}/\!\!\sim$ whose support is total. 
Proposition \ref{conjugaralafin} then implies that the latter 
action is topologically conjugate to an action by affine transformations. Therefore, 
the original action of $\Gamma$ on the line is semiconjugate to an action by 
affine transformations.

It remains the case where $v$ has atoms. First notice that $\Gamma'$ preserves $v$, and 
hence the second derived group $\Gamma''$ fixes each atom of $v$. This argument shows 
in particular that, in this case, the index $j$ considered in the proof of Proposition 
\ref{quasiinvariante} equals either $1$ or $2$. We will see below that it is actually 
equal to $1$.

Denote by $\Gamma_{v}$ the subgroup of $\Gamma$ formed by the elements fixing $v$, and denote 
by $\Gamma_{v}'$ its derived group. The elements in $\Gamma_{v}'$ fix the atoms in $v$. Denote 
by $\Gamma^*_{v}$ the normal subgroup of $\Gamma$ formed by the elements fixing the irreducible 
components of $\Gamma'_{v}$. The arguments of proof of Claims (i) and (ii) in Proposition 
\ref{quasiinvariante} show that $\Gamma_{v} / \Gamma_{v}^*$ is isomorphic to a nontrivial subgroup 
$H$ of $(\mathbb{R},+)$. Notice that $H$ cannot be dense, for otherwise there would 
be atoms of $v$ of the same mass accummulating on some points in the line, 
thus contradicting the 
fact that $v$ is a Radon measure. Therefore, $H$ is infinite cyclic, and since 
$\Gamma$ acts by automorphisms of $H$, the latter group must 
preserve $v$. Therefore, $\Gamma = \Gamma_{v}$ (this shows 
that $j \!=\! 1$). We then conclude that $\Gamma$ is an extension 
(actually, a semidirect product) by $(\mathbb{Z},+)$ of a solvable 
subgroup of a product of groups of diffeomorphisms of closed intervals. 
This concludes the proof of Theorem \ref{solrecta}.


\section{On the Smooth Actions of Amenable Groups\index{group!amenable}}

\hspace{0.45cm} On the basis of the previous sections, it would be natural 
to continue the classification of groups of interval and circle 
diffeomorphisms by trying to describe 
the amenable ones. However, this problem seems to be 
extremely difficult. Actually, the only relevant result in this direction 
\index{Witte-Morris!theorem on orderable amenable groups}
is Theorem \ref{witte-bien}, 
which due to Thurston's stability theorem from \S \ref{sec-th-th} 
does not provide any information in the 
case of groups of diffeomorphisms. As a matter of example, let us recall that, 
according to \S \ref{ghys-sergi}, \index{Thompson's groups!F}
Thompson's group F embeds into $\mathrm{Diff}_+^{\infty}([0,1])$. However, 
the problem of knowing whether F is amenable 
has been open for more than 30 years$\esp\!!$ 
In what follows, we will give some partial results on the classification of a 
particular but relevant family of amenable subgroups of $\mathrm{Diff}_+^2([0,1])$.

Recall that amenability is stable under {\bf {\em elementary operations}}, that is, it is 
preserved by passing to subgroups or quotients, and by taking extensions or unions (see Exercise 
\ref{clasemoy}). Following \cite{osin}, we denote by $\mathrm{SA}$ the smallest family of 
amenable groups which is closed with respect to these operations and which contains the 
groups of sub-exponential growth (see Exercise \ref{sub-prom}). 

A result in \cite{Nasol} asserts that  
\index{growth!of groups}\index{group!metabelian}every group of real-analytic 
diffeomorphisms of $[0,1[$ contained in $\mathrm{SA}$ is metabelian. However, 
$\mathrm{Diff}_+^{\infty}([0,1])$ contains interesting subgroups from SA. The 
construction below should be compared with those in \cite{Pl2}.

\begin{small} \begin{ejem} Given points $a\!<\!c\!<\!d\!<\!b$ in $]0,1[$, let $f$ be 
a $\mathrm{C}^{\infty}$ diffeomorphism of $[a,b]$ which is infinitely tangent to the 
identity at the endpoints and such that $f(c) = d$. Extend $f$ to $[0,1]$ by letting 
$f(x)=x$ for $x \!\notin ]a,b[$. Let $g$ be a $\mathrm{C}^{\infty}$ diffeomorphism of 
$[0,1]$ which is infinitely tangent to the identity at the endpoints, with a single  
fixed point at the interior, and such that $g(c)=a$ and $g(d)=b$. Denote by $\Gamma$ 
the subgroup of $\mathrm{Diff}^{\infty}_+([0,1])$ generated by $f$ and $g$. For each 
$n \!\in\! \mathbb{N}$, the subgroup $\Gamma_n$ of $\Gamma$ generated by 
$\{f_i \!=\! g^{-i}fg^i\!: |i| \!\leq\! n \}$ is solvable with solvability degree 
$2n\!+\!1$. Moreover, the subgroup $\Gamma^* \!=\! \cup_{n\in\mathbb{N}}\Gamma_n$ is 
normal in $\Gamma$, and the quotient $\Gamma/\Gamma^*$ is isomorphic to $(\mathbb{Z},+)$ 
(a generator is $g\Gamma^*$). The group $\Gamma$ is therefore finitely generated and 
non solvable, \index{group!solvable}
and belongs to the family SA (see \cite{ana,gromov-ggd} for an interesting 
property concerning its F\o lner sequences).
\index{F\o lner sequence}
\label{moy-no}
\end{ejem} \end{small}

Despite the above example, the subgroups of $\mathrm{Diff}_+^{1+\mathrm{bv}}([0,1[)$ 
which do belong to the family SA may be partially classified. For this, let us define 
by transfinite induction the subfamilies $\mathrm{SA}_{\alpha}$ of SA by:

\vspace{0.1cm}

\noindent{(i) $\mathrm{SA}_1$ is the family of countable groups all 
of whose finitely generated subgroups have sub-exponential growth.}

\vspace{0.1cm}

\noindent{(ii) If $\alpha$ is not a limit ordinal, then $\mathrm{SA}_{\alpha}$
is the family of groups $\Gamma$ obtained either as a quotient, as a subgroup, or 
as an extension of groups in $\mathrm{SA}_{\alpha - 1}$. The latter means that 
$\Gamma$ contains a normal subgroup $G$ in $\mathrm{SA}_{\alpha - 1}$ so that 
the quotient $\Gamma / G$ also belongs to $\mathrm{SA}_{\alpha -1}$.}

\vspace{0.1cm}

\noindent{(iii) If $\alpha$ is a limit ordinal, then $\mathrm{SA}_{\alpha}$ 
is the family of groups obtained as union of groups in the union 
\esp $\bigcup \mathrm{SA}_{\beta}$, \esp with $\beta < \alpha$.}

\vspace{0.1cm}

A group $\Gamma$ belongs to $\mathrm{SA}$ if and only if it belongs to $\mathrm{SA}_{\alpha}$ 
for some ordinal $\alpha$ \cite{osin}. For example, the group from example \ref{moy-no} belongs 
to $\mathrm{SG}_{\alpha+1}$, where $\alpha$ is the first infinite ordinal. The next result 
(whose proof is based on those of \S \ref{solvito}) appears in \cite{Na-subexp} (see also 
\cite{brin-subexp} for a more accurate version in the piecewise affine case).

\vspace{0.12cm}

\begin{thm} {\em Every subgroup of $\mathrm{Diff}_+^{1+\mathrm{bv}}([0,1[)$ 
which belongs to $\mathrm{SA}_{\alpha}$ for some finite $\alpha$ is 
solvable with solvability degree less than or equal to $\alpha$.}
\index{degree!of solvability|)}
\label{two2}
\end{thm}

\vspace{0.12cm}

Theorem \ref{two2} only applies to particular families of amenable groups. Indeed, 
from \cite{bart} it follows that the family SA is smaller than that of amenable 
groups. This produces a breaking point in our classification of subgroups of 
$\mathrm{Diff}_+^{1+\mathrm{bv}}(\clo)$. This is the reason why in the next 
chapter we will follow an ``opposite direction'' in our study: we will show 
that, due to some internal (mostly cohomological) properties, certain 
groups cannot act on the interval or the circle in a nontrivial way. 

\newpage
\thispagestyle{empty}
${}$
\newpage


\chapter{RIGIDITY VIA COHOMOLOGICAL METHODS}
\label{cap-cohom}


\section{Thurston's Stability Theorem\index{Thurston!stability theorem}}
\label{sec-th-th}

\hspace{0.45cm} The result below corresponds to a weak version of a theorem due to 
Thurston. For the general version related to the famous Reeb's stability theorem for 
foliations, we refer to the remarkable paper \cite{Th2}. We also strongly 
recommend to look at Exercise \ref{thc}.
\index{foliation}

\vspace{0.15cm}

\begin{thm} \textit{Let $\Gamma$ be a finitely generated group. If $\Gamma$ does not 
admit any nontrivial homomorphism into $(\mathbb{R},+)$, then every representation 
$\Phi\!: \Gamma \rightarrow \mathrm{Diff}_{+}^{1}([0,1[)$ is trivial.\footnote{Equivalently, 
for every nontrivial finitely generated subgroup $\Gamma$ of $\mathrm{Diff}_{+}^{1}([0.1])$, 
the first {\em cohomology space} $H^1_{\mathbb{R}}(\Gamma)$ is nontrivial.}}
\label{ththurston}
\end{thm}

\index{cohomology space}

\vspace{0.15cm}

For the proof of this theorem, given a subset $B$ of an arbitrary group $\Gamma$ and 
$\varepsilon \geq 0$, we will say that a function $\phi\!: B \rightarrow \mathbb{R}$ 
is a {\bf \textit{$(B,\varepsilon)$-homomorphism}} (into $\mathbb{R}$) 
if for every $g,h$ in $B$ such that $gh \in B$ one has 
$$\big| \phi(g) + \phi(h) - \phi(gh) \big| \leq \varepsilon$$ 
Fix a finite and symmetric system of generators $\mathcal{G}$ of $\Gamma$. We will 
say that $\phi$ above is {\bf \em{normalized}} if \esp 
$\max_{g \in \mathcal{G}} |\phi(g)| = 1.$ \esp For simplicity, we let 
$$\nabla\phi(g,h) = \phi(g) + \phi(h) - \phi(gh).$$ 
Notice that a function  
$\phi\!: \Gamma \rightarrow \mathbb{R}$ is a $(\Gamma,0)$-homomorphism if and 
only if $\nabla\phi$ is identically zero, which is equivalent to that $\phi$ is a  
homomorphism from $\Gamma$ 
into $(\mathbb{R},+)$. Recall finally that $B_{\mathcal{G}}(k)$ denotes 
the set of elements in $\Gamma$ which may be written as a product of at most $k$ 
elements in $\mathcal{G}$ (see Appendix B).

\vspace{0.08cm}

\begin{lem} \textit{If for each $k \in \mathbb{N}$ there exists a normalized 
$\big( B_{\mathcal{G}}(k),1/k \big)$-homomorphism, then there exists a 
nontrivial homomorphism from $\Gamma$ into $(\mathbb{R},+)$.}
\label{mejor-aun}
\end{lem}

\noindent{\bf Proof.} Given $k \!\in\! \mathbb{N}$, let $\phi_k$ be a normalized 
$\big( B_{\mathcal{G}}(k),1/k \big)$-homomorphism. It is easy to check that 
$|\phi_k(g)| \leq k(1+\varepsilon)$ for every $g \!\in\! B_{\mathcal{G}}(k)$. 
In particular, there exists a subsequence 
$(\phi_k^1)$ of $\phi_k$ such that $\phi^1_k|_{B_{\mathcal{G}}(1)}$ converges 
(punctually) to a normalized function from $B_{\mathcal{G}}(1)$ into $\mathbb{R}$. 
Arguing by induction, for each  
$i \in \mathbb{N}$ we may find a subsequence $(\phi_{k}^{i+1})_{k \in \mathbb{N}}$ 
of $(\phi_k^i)_{k \in \mathbb{N}}$ so that $\phi^{i+1}_{k} |_{B_{\mathcal{G}}(i+1)}$ 
converges to a normalized function from $B_{\mathcal{G}}(i+1)$ into $\mathbb{R}$. 
The sequence $(\phi_k^k)$ converges (punctually) to a function
$\phi\!: \Gamma \!\rightarrow\! \mathbb{R}$. By construction, $\phi$ is a 
homomorphism into $(\mathbb{R},+)$, which is necessarily nontrivial since it is 
normalized. $\hfill\square$

\vspace{0.35cm}

Now let $\Phi\!: \Gamma \rightarrow \mathrm{Diff}_+^1([0,1])$ be a nontrivial 
representation, and let $x \!\in\! [0,1[$ be a point in the boundary of the set 
of points which are fixed by $\mathcal{G}$ (and hence by $\Gamma$). For each 
$y \!\in\! [0,1[$ which is not fixed by $\mathcal{G}$ we consider the function
$$\phi_y(g) = \frac{1}{C(y)} \big[ \Phi(g)(y)-y \big],$$
where $C(y)=\max_{g \in \mathcal{G}} \big| \Phi(g)(y)-y \big|.$
The lemma below shows that, if \esp 
$\Phi(g)'(x) \!=\! 1$ \esp for every $g \!\in\! \Gamma$, 
then for $y$ near to $x$ (and not fixed by $\mathcal{G}$) the function $\phi_y$ 
behaves infinitesimally like a homomorphism from $\Gamma$ into $(\mathbb{R},+)$.

\vspace{0.1cm}

\begin{lem} \textit{Under the above conditions, suppose moreover 
that $\Phi(g)'(x) \!=\! 1$ for every $g \in \Gamma$. Then for each 
$n \in \mathbb{N}$ and each $\varepsilon > 0$ there exists $\delta > 0$ 
such that, if $|x-y| < \delta$, then $\phi_y|_{B_{\mathcal{G}}(n)}$ is a 
normalized $\big( B_{\mathcal{G}}(n),\varepsilon \big)$-homomorphism.}
\label{lemthurston}
\end{lem}

\noindent{\bf Proof.} 
For $k \in \mathbb{N}$ and $\varepsilon' > 0$ we inductively define 
$$\lambda_0(0,\varepsilon') = 0, \qquad 
\lambda_0(k+1,\varepsilon') = 1+ \lambda_0(k,\varepsilon')(1+\varepsilon').$$
Let $\epsilon' > 0$ be small enough so that \esp 
$\varepsilon' \lambda_0(n,\varepsilon') \leq \varepsilon,$ \esp 
and let \esp $\delta' > 0$ be such that \esp 
$| \Phi(g)'(y) - 1 | \leq \varepsilon'$ \esp for all 
$g \in B_{\mathcal{G}}(n)$ and every point $y \in [0,1]$ satisfying 
$|x-y| \leq \delta'$. Finally, let $\delta \! \in ]0,\delta'[$ be such that, 
if $|x-y| \leq \delta$, then \esp $|\Phi(g)(y)-x| \leq \delta'$ \esp for 
every $g \in B_{\mathcal{G}}(n)$. We claim that this parameter $\delta$ 
verifies the claim of the lemma.

Let us first show that, for every $k \leq n$ 
and every \esp $g \!\in\! B_{\mathcal{G}}(k)$,
\begin{equation}
|\phi_y(g)| \leq \lambda_0(k,\varepsilon').
\label{thcua}
\end{equation}
Indeed, this inequality is evident for $k=0$ and $k=1$. Let us suppose that it holds 
for $k=i$. If $g$ is an element in $B_{\mathcal{G}}(i+1)$, then $g=h_1h_2$ for 
some $h_1 \!\in\! \mathcal{G}$ and $h_2 \!\in\! B_{\mathcal{G}}(i)$. Thus, 
\begin{equation}
|\phi_y(g)| \esp = \esp \frac{1}{C(y)} \big| \Phi(g)(y)-y \big|
      \esp \leq \esp \frac{1}{C(y)} \big| \Phi(h_1)\Phi(h_2)(y) - \Phi(h_2)(y) \big| 
                        + |\phi_y(h_2) |.
\label{thcinc}
\end{equation}
Notice that from the equality  
$$\Phi(h_1)\Phi(h_2)(y)-\Phi(h_2)(y) \esp = \esp 
  \int_{y}^{\Phi(h_2)(y)} \!\! \big[ \Phi(h_1)'(s)-1 \big] ds \esp 
         + \esp \big[ \Phi(h_1)(y)-y \big]$$
one deduces that
\begin{small}
$$\frac{1}{C(y)} \big| \Phi(h_1)\Phi(h_2)(y) - \Phi(h_2)(y) \big| \esp \leq
\esp \max\limits_{|s-x| \leq \delta'} \big| \Phi(h_1)'(s)-1 \big| \cdot
\frac{1}{C(y)} \big| \Phi(h_2)(y)-y \big| + 1.$$ 
\end{small}Using (\ref{thcinc}) and the induction hypothesis we conclude that
$$|\phi_y(g)| \leq \varepsilon' \lambda_0(i, \varepsilon') + 1 +
\lambda_0(i,\varepsilon') = \lambda_0(i+1,\varepsilon'),$$
which finishes the proof of (\ref{thcua}).\\

Let us now estimate the value of $\nabla\phi_y$. For $h_1, h_2$ 
in $B_{\mathcal{G}}(n)$ such that $h_1 h_2$ belongs to 
$B_{\mathcal{G}}(n)$ we have
$$\nabla\phi_y(h_1,h_2)
= \frac{1}{C(y)} \big[ \Phi(h_1)(y)-y+\Phi(h_2)(y)-y-\Phi(h_1h_2)(y)+y \big] ,$$
that is
\begin{eqnarray*}
\nabla\phi_y(h_1,h_2)
&=& -\frac{1}{C(y)} \big[ \Phi(h_1)\Phi(h_2)(y)-\Phi(h_2)(y) - (\Phi(h_1)(y) -y)
\big] \\
&=& -\frac{1}{C(y)}\int_{y}^{\Phi(h_2)(y)} \big[ \Phi(h_1)'(s)-1 \big] ds.
\end{eqnarray*}
Thus, 
$$|\nabla\phi_y(h_1,h_2)| \esp 
\leq \esp \Big| \frac{1}{C(y)}[\Phi(h_2)(y)-y] \Big|
\cdot \!\sup\limits_{|s-x| \leq \delta'}|\Phi(h_1)'(s) - 1|
\esp \leq \esp \lambda_0(n,\varepsilon') \esp \! \varepsilon' 
\esp \leq \esp \varepsilon.$$
Therefore, $\phi_y|_{B_{\mathcal{G}}(n)}$ is a normalized 
$\big( B_{\mathcal{G}}(n),\varepsilon \big)$-homomorphim. $\hfill\square$

\vspace{0.4cm}

\noindent{\bf Proof of Theorem \ref{ththurston}.} Suppose that 
$\Phi\!: \Gamma \rightarrow \mathrm{Diff}_+^1([0,1[)$ is a nontrivial 
representation, and let $x \in [0,1[$ be a point in the boundary of 
the set of fixed points of $\Gamma$. The function 
\esp $g \mapsto \log \big( \Phi(g)'(x) \big)$ \esp 
is a group homomorphism from $\Gamma$ into $(\mathbb{R},+)$. 
By hypothesis, this homomorphism must be trivial (compare 
Exercise \ref{truco-de-tsuboi}). We are hence under the 
hypothesis of Lemma \ref{lemthurston}, which together with Lemma 
\ref{mejor-aun} implies that $\Gamma$ admits a nontrivial homomorphism 
into $(\mathbb{R},+)$. $\hfill\square$

\vspace{0.4cm}

Notice that the argument above only involves the behavior of the maps near 
a global fixed point. It is then easy to see that Thurston's theorem still 
holds (with the very same proof) for subgroups of the group of germs of 
diffeomorphism $\mathcal{G}_{+}^{1}(\mathbb{R},0)$.
\index{germ of a diffeomorphism}

\vspace{0.02cm}

\begin{small}\begin{ejer} Prove that Thurston's stability theorem does not hold 
for non finitely generated countable groups of diffeomorphisms of the interval.

\noindent{\underbar{Hint.}} Consider the first derived group F' of 
Thompson's group F, \index{Thompson's groups!F} and use the fact 
that F' is simple \cite{CFP}.\index{group!simple}
\end{ejer}

\begin{ejer} Prove that Thurston's stability theorem holds for (perhaps non 
finitely generated) groups of germs of {\em real-analytic} diffeomorphisms.

\noindent{\underbar{Hint.}} Following an argument due to Haefliger (and prior to 
Thurston's theorem), analyze the coefficients corresponding to the Taylor series 
expansions of the germs about the origin.
\label{th-analitico}\index{Haefliger}
\end{ejer} \end{small}

We now give an example (also due to Thurston, but rediscovered some 
years later by Bergman \cite{bergman} in the context of orderable groups) 
\index{Thurston!Example} 
showing that Thurston's stability theorem does not 
extend to actions by homeomorphisms (see also Remarks \ref{ojuela} 
and \ref{en-general}). Let $\tilde{G}$ be 
the group of presentation \esp $\tilde{G} = \langle f,g,h : \esp f^2=g^3=h^7 = fgh \rangle.$ 
\esp We leave to the reader the task of showing that every homomorphism from $\tilde{G}$ 
into $(\mathbb{R},+)$ is trivial. (Actually, 
$\tilde{G} = [\tilde{G},\tilde{G}]$, that is, $\tilde{G}$ is a {\em perfect} group.) 
\index{group!perfect} This is quite natural, since $\tilde{G}$ is the 
fundamental group of a homology 3-sphere, that is, a closed 3-dimensional 
manifold with trivial homology but not homeomorphic to the 3-sphere \cite{Th1,Th2}.

To construct a nontrivial $\tilde{G}$-action on $[0,1]$, 
let us consider the tessellation of the Poincar\'e disk 
\index{Poincar\'e!disk} by hyperbolic triangles of angles 
$\pi/2$, $\pi/3$, and $\pi/7$. The preimage in $\widetilde{\mathrm{PSL}}(2,\mathbb{R})$ 
of the subgroup of $\mathrm{PSL}(2,\mathbb{R})$ preserving this tessellation is a 
group isomorphic to $\tilde{G}$. Since $\widetilde{\mathrm{PSL}}(2,\mathbb{R})$ acts 
on $\widetilde{\mathrm{S}^1} = \mathbb{R}$, viewing the interval $[0,1]$ as a 
two-point compactification of $\mathbb{R}$ we obtain a faithful 
$\tilde{G}$-action on $[0,1]$. 
Notice that, by Thurston's theorem, the latter action cannot be smooth: the obstruction 
to the differentiability localizes around each endpoint of the interval.
\index{homology sphere} 

\begin{small}\begin{ejer} Using the fact that $fgh$ belongs to the center of $\tilde{G}$, 
\index{center of a group}
as well as Propositions \ref{radon-inv} and \ref{facil}, prove that the 
$\tilde{G}$-action above is not conjugate to an action by $\ce^1$ diffeomorphisms 
of the interval $[0,1[$ without using Thurston's stability theorem.
\end{ejer}

\begin{ejer} Following \cite{cal-forcing}, 
consider the group $\hat{G}$ with presentation 
$$\langle f_1\!,\!g_1\!,\!h_1\!,\!f_2,\!g_2,\!h_2\!: 
f_1^2\!=\!g_1^3\!=\!h_1^7 \!=\! f_1g_1h_1, f_2^{-1}f_1f_2\!=\!\!f_1^2, 
g_2^{-1}g_1g_2\!=\!g_1^2, h_2^{-1}h_1h_2\!=\!h_1^2 \rangle.$$
Show that $\hat{G}$ contains a copy of $\tilde{G}$ and  
acts faithfuly by circle homeomorphisms, but it 
does not embed into $\mathrm{Diff}_+^1(\clo)$. 
\index{action!faithful}
\end{ejer}

\begin{obs} The family of finitely generated groups of homeomorphisms of the interval containing 
finitely generated subgroups that do not admit nontrivial homomorphisms into $(\mathbb{R},+)$ 
is quite large. Indeed, an important result in the theory of orderable 
\index{group!orderable} groups asserts that 
a group (of arbitrary cardinality) is \index{group!locally indicable}
{\bf {\em locally indicable}} ({\em i.e.,} all 
of its finitely generated subgroups admit nontrivial homomorphisms into $(\mathbb{R},+)$) if 
and only if it is $\mathcal{C}$-orderable 
\index{group!$\mathcal{C}$-orderable} 
(see \S \ref{super-witte-morris} for the 
notion of $\mathcal{C}$-order, and see \cite{ordering2} for an elementary 
proof of this result). Therefore, although all finitely generated groups of interval 
homeomorphisms are topologically conjugate to groups of Lipschitz homeomorphisms of $[0,1]$ 
\index{Lipschitz!homeomorphism}
({\em c.f.}, Proposition \ref{abu}), many of them do not embed into $\mathrm{Diff}_+^1([0,1[)$ 
because they fail to be locally indicable. However, this is not the only algebraic 
obstruction: there exist finitely generated, locally indicable groups without 
(faithful) actions by $\ce^1$ diffeomorphisms of the interval~! (see 
Remark \ref{thu-na-rem}).
\label{en-general}
\end{obs}\end{small}

A nice consequence of what precedes is that every countable group of $\ce^1$ diffeomorphisms of 
the interval admits a faithful action by homeomorphisms of the interval without crossed elements. 
Indeed, by Theorem \ref{ththurston}, such a group $\Gamma$ is locally indicable, and hence by 
Remark \ref{en-general}, it is $\mathcal{C}$-ordenable. Proposition \ref{conrad-yo} (and the 
comments after its proof) then implies that the dynamical realization of any Conradian ordering 
on $\Gamma$ is a subgroup of $\mathrm{Homeo}_+(\mathbb{R})$ without crossed elements (compare 
Exercise \ref{no-F}).

\index{crossed elements}

\vspace{0.02cm}

\begin{small}\begin{obs} 
As was cleverly shown by Calegari in \cite{cal-thurston}, the technique 
of proof of Thurston's stability theorem (but not its conclusion~!) may be 
used for showing that certain actions of locally indicable groups (as for 
example the free group $\mathbb{F}_2$) on the closed interval are non $\ce^1$ 
smoothable. This is the case for instance if the set of fixed points 
\index{Calegari}
of the generators accumulate at the endpoints and their commutator has no fixed 
point inside. Let us point out, however, that these actions have (plenty of) crossed 
elements ({\em c.f.,} Definition \ref{def-crossed-elements}). 
\end{obs}

\begin{obs} Theorem \ref{ththurston} may be also used for studying group actions on higher 
dimensional manifolds. For instance, using the idea of its proof, in \cite{cal-th} it 
is shown that the group of the $\ce^1$ diffeomorphisms of the closed disk fixing all 
points in the boundary is orderable.
\end{obs}

\begin{ejer} By combining Thurston's stability theorem with the result in Exercise 
\ref{proba-para-productos}, prove the following result from \cite{parwa}: If 
$\Gamma = \Gamma_1 \times \Gamma_2$ is a finitely generated group such that 
no finite index normal subgroup of $\Gamma_1$ and $\Gamma_2$ has a nontrivial  
homomorphism into $(\mathbb{R},+)$, then for every action of $\Gamma$ by $\ce^1$ 
circle diffeomorphisms the image of either $\Gamma_1$ or $\Gamma_2$ is finite.
\label{ejer-parwa}
\end{ejer}

\begin{ejer} A topological group is said to be {\bf {\em compactly generated}} 
if there exists a compact subset such that every element may be written as a product 
of elements therein. Extend Thurston's stability theorem to continuous representations 
in $\mathrm{Diff}_+^1([0,1[)$ of compactly generated groups.
\index{group!compactly generated}
\label{compactamente-generado}
\end{ejer}

\begin{ejer} The goal of this exercise consists in giving an alternative proof of  
Thurston's stability theorem using a technique from \cite{th-1} and \cite{th-2}. Let us 
consider a finitely generated subgroup $\Gamma$ of $\mathrm{Diff}^1_+([0,1[)$, and to 
simplify let us suppose that $\Gamma$ has no fixed point inside $]0,1[$ and that all 
of its elements are tangent to the identity at the origin (as we have already seen, the 
general case easily reduces to this one). Let $\mathcal{G}\!=\!\{ h_1, \ldots, h_k \}$ 
be a finite family of generators for $\Gamma$. For each $f \in \Gamma$ let us define 
the {\bf{\em displacement function}} $\Delta_f$ by $\Delta_f(x)=f(x)-x$. Notice that 
$(\Delta_f)'(0)=0$ for every $f \in \Gamma$.\\

\vspace{0.1cm}

\noindent{(i) Prove that for every $x \geq 0$ and every $f,g$ in $\Gamma$, 
there exist $y,z$ in $[0,x]$ such that}
$$\Delta_{fg}(x) = \Delta_f(x) + \Delta_g(x) + (\Delta_f)'(y) \esp \Delta_g(x), 
\qquad \Delta_{f^{-1}(x)} = - \Delta_f (x) - (\Delta_f)'(z) \esp \Delta_{f^{-1}}(x).$$

\vspace{0.08cm}

\noindent{(ii) Fixing an strictly decreasing sequence of points 
$x_n$ converging to the origin, for each $n \in \mathbb{N}$ let us choose 
$i_n \!\in\!\{1,\ldots,k\}$ so that \esp $|\Delta_{h_{i_n}}(x_n)| \geq |\Delta_{h_j}(x_n)|$ 
\esp for every $j \!\in\!\{1,\ldots,k\}$. Passing to a subsequence if necessary, we may 
assume that $i_n$ is constant (say, equal to $1$ after reordering the indexes), and that 
each of the $k$ sequences $( \Delta_{h_i} (x_n) /\Delta_{h_1}(x_n))$ converges to a 
limit $\phi_i$ (less than or equal to $1$) as $n$ goes to infinity. Use the equalities 
in (i) for showing that the map \esp $h_i \longmapsto \phi_i$ \esp extends to a 
normalized homomorphism from $\Gamma$ into $(\mathbb{R},+).$}
\label{thc}
\end{ejer}

\begin{ejer} Following \cite{muller} and \cite{Ts}, consider a diffeomorphism $\varphi$ 
from $]0,1[$ into itself such that $\varphi(s) = \exp(-1/s)$ for $s > 0$ small enough. 
Prove that if $k \geq 0$ and $f: [0,1[ \rightarrow [0,1[$ is a $\mathrm{C}^k$ diffeomorphism, 
then (the extension to $[0,1[$ of) $\varphi^{-1} \circ f \circ \varphi$ (resp. 
$\varphi^{-2} \circ f \circ \varphi^2$) is a $\mathrm{C}^k$ diffeomorphism with 
derivative $1$ (resp. tangent to the identity up to order $k$) at the origin.
\index{Tsuboi!M\"uller-Tsuboi conjugacy}
\label{truco-de-tsuboi}
\end{ejer} \end{small}


We close this section by sketching the proof of Theorem \ref{thm-matsuda}. Let 
$\Gamma$ be a polycyclic subgroup of $\mathrm{Diff}_+^{1+\mathrm{Lip}}([0,1])$. 
By \S \ref{policiclicos-suaves}, we know that $\Gamma$ is metabelian. According 
to Theorem \ref{clas-meta}, we need to show that $\Gamma$ cannot be a semidirect 
product between $\mathbb{Z}$ and a subgroup of a product of groups of translations 
as described in the proof of that theorem. Assume for a contradiction that 
this is the case. Then 
it is easy to see that the nilradical $N(\Gamma)$ of $\Gamma$ is torsion-free Abelian (say, 
isomorphic to $\mathbb{Z}^k$ for some $k \!\in\! \mathbb{N}$) and coincides with 
\index{nilradical}
the subgroup formed by the elements having fixed points in $]0,1[$. Moreover, if 
$\{f_1,\ldots,f_k\}$ is a system of generators for $N(\Gamma)$, then the group $\Gamma$ is 
generated by these elements and a certain $f \in \Gamma$ without fixed points in $]0,1[$. 

Let $]a,b[$ be an irreducible component of $N(\Gamma)$; we must necessarily have 
$[a,b] \subset ]0,1[$. If we identify $N(\Gamma)$ with $\mathbb{Z}^k$ via the map 
\esp $(m_1,\ldots,m_k) \mapsto g = f_1^{m_1} \cdots f_k^{m_k} \in N$, \esp 
we readily see that $f$ naturally induces an element $f_{*}$ in $\mathrm{Hom} 
(\mathrm{Hom}(\mathbb{Z}^k,\mathbb{R}),\mathrm{Hom}(\mathbb{Z}^k,\mathbb{R}))$, 
\esp namely \esp $f_{*} (\phi) (g) = \phi(fgf^{-1})$. \esp Viewing each element 
in $\mathrm{Hom}(\mathbb{Z}^k,\mathbb{R})$ as the restriction of a homomorphism 
defined on $\mathbb{R}^k$, the homomorphism $f_*$ yields a linear map 
$M \!\in\! \mathrm{GL} (k,\mathbb{R})$. 

Let $\phi_0 \!: \mathbb{Z}^k \rightarrow \mathbb{R}$ be the homomorphism obtained via 
Thurston's stability theorem applied to the restriction of $N(\Gamma) \sim \mathbb{Z}^k$ 
to $[a,b[$. From the construction of $\phi_0$, it is easy to see that if $(g_n)$ is a 
sequence of elements in $N(\Gamma)$ whose restrictions to $[a,b]$ converge to the 
identity in the $\ce^1$ topology, then $\phi_0 (g_n)$ converges to zero as 
$n$ goes to infinity. Now recall the 
following important fact from Exercise \ref{kopell-fuerte-lip}: 
For each $g \!\in\! N(\Gamma)$ 
both sequences $(f^{n}gf^{-n})$ and $(f^{-n}gf^{n})$ converge to the identity in 
the $\ce^1$ topology when restricted to $[a,b]$. Therefore, both $f_*^n (\phi_0)$ 
and $f^{-n}_* (\phi_0)$ punctually converge to the zero homomorphism. This means 
that for the linear map $M$ there exists a non-zero vector $v$ so that both 
$M^n (v)$ and $M^{-n} (v)$ converge to zero, which is clearly impossible. 
This contradiction concludes the proof. 
\index{group!polycyclic}


\section{Rigidity for Groups with Kazhdan's Property (T)}
\label{nose}

\subsection{Kazhdan's property (T)}
\label{def-T}\index{Kazhdan's property (T)|(}

\hspace{0.45cm} Let $\Gamma$ be a countable group and $\Psi: \Gamma \rightarrow U(\mathcal H)$ 
a unitary representation of $\Gamma$ in some (real) Hilbert space $\mathcal H$. We  say tthat 
$c\!: \Gamma \rightarrow \mathcal H$ is a {\bf {\em cocycle}} with respect to $\Psi$ 
\index{cocycle}
if for every $g_1,g_2$ in $\Gamma$ one has 
$$c(g_1g_2) = c(g_1) + \Psi(g_1)c(g_2).$$  
We say that a cocycle $c$ is a {\bf {\em coboundary}} 
\index{coboundary}
if there exists $K \in \mathcal H$ so that, for every $g \!\in\! \Gamma$, 
$$c(g) =  K - \Psi(g)K.$$
We denote the space of cocycles by $Z^1(\Gamma,\Psi)$ and the subspace of coboundaries 
by $B^1(\Gamma,\Psi)$. The quotient $H^1(\Gamma,\Psi) = Z^1(\Gamma,\Psi) / B^1(G,\Psi)$ 
is the {\bf \textit{first cohomology space}} of $\Gamma$ (with values in $\Psi$). 
\index{cohomology space}

\vspace{0.05cm}

\begin{defn} A countable group $\Gamma$ satisfies the {\bf \textit{Kazhdan's Property {\em (T)}}} 
(or, to simplify, is a {\em Kazhdan group}) if for every unitary representation $\Psi$ of 
$\Gamma$, the first cohomology $H^1(\Gamma,\Psi)$ is trivial.
\label{kazhdandefinicion}
\end{defn}

\vspace{0.05cm}

To give a geometrical insight of this definition, let us recall that 
every {\bf\em{isometry}}\footnote{By an isometry of a Hilbert space 
$\mathcal{H}$ we mean a {\em surjective} map $A$  from $\mathcal{H}$ 
into itself satisfying $\| A(K_1) - A(K_2) \| = \|K_1 -K_2\|$ \esp 
for all $K_1, K_2$ in $\mathcal{H}$.} 
of a Hilbert space is the composition of a unitary transformation and a translation. 
Indeed, the group of  isometries of $\mathcal{H}$ is the semidirect product between 
$U(\mathcal{H})$ and $\mathcal{H}$.

Let $\Psi\!: \Gamma \!\rightarrow\! U(\mathcal H)$ be a unitary representation and 
$c\!: \Gamma \!\rightarrow\! \mathcal{H}$ a map. If for each $g \!\in\! \Gamma$ we 
define the isometry $A(g) \!=\! \Psi(g) + c(g),$ then one easily checks that the 
equality $A(g_1)A(g_2) \!=\! A(g_1g_2)$ holds for every $g_1,g_2$ in $\Gamma$ 
if and only if $c$ is a cocycle associated to $\Psi$. In this case, the 
correspondence $g \mapsto A(g)$ defines an action by isometries.

Suppose now that $K \!\in\! \mathcal{H}$ is a fixed point for the isometric action 
associated to a cocycle $c\!: \Gamma \!\rightarrow\! \mathcal{H}$. Then we have, 
for every $g \!\in\! \Gamma$, 
$$\Psi(g)K + c(g) = A(g)K = K.$$
Therefore, \esp $c(g) = K - \Psi(g)K,$ \esp that is, \esp 
$c$ \esp is a coboundary. Conversely, it is easy to check that if \esp $c$ 
\esp is a coboundary, then there exists an invariant vector for the associated 
isometric action. We then have the following geometric interpretation of Definition 
\ref{kazhdandefinicion}: A group has Kazhdan's property (T) if and only if every 
action by isometries of a Hilbert space has an invariant vector.

\begin{small} \begin{ejem} Every finite group has Kazhdan's property (T). To see this, it 
suffices to remark that the mean along any orbit of an action by isometries on a Hilbert 
space is an invariant vector for the action.
\end{ejem}

\begin{ejer}
A vector $K_{\mathcal{C}} \in \mathcal{H}$ is the {\bf{\em geometric center}} of a 
subset $\mathcal{C}$ of $\mathcal{H}$ if $K_{\mathcal{C}}$ minimizes the function 
$K \mapsto \sup_{\bar{K} \in \mathcal{C}} \| K - \bar{K} \|$. Show that if $\mathcal{C}$ 
is bounded, then it has a unique center. Conclude that if $\mathcal{C}$ is invariant 
under an action by isometries, then its center remains fixed by the action. 
Use this to prove the following {\em Center Lemma} (due to Tits): 
\index{Tits!Center Lemma} 
An action by isometries has a fixed point if and only if its orbits are bounded. 
(In particular, a group has property (T) if and only if the orbits associated to 
its isometric actions on Hilbert spaces are bounded.)

\noindent{\underbar{Remark.}} The same argument applies to spaces on which 
the distance function satisfies a {\em convexity property}, as for instance 
simplicial trees or spaces with non-positive curvature \cite{Hae}.
\index{tree}
\label{centro-tits}
\end{ejer}
\end{small}

Every group with property $\mathrm{(T)}$ is finitely generated, according to a result 
due to Kazhdan himself (actually, this was one of his motivations for introducing 
property (T)). We give below an alternative argument due to Serre based on 
the result from Example \ref{arbolito}.

Let $\Gamma$ be a countable group, and let  
$\Gamma_1 \subset \Gamma_2 \subset \ldots \subset \Gamma_n \subset \ldots$
be a sequence of finitely generated subgroups whose union is $\Gamma$. Let us 
consider the oriented simplicial tree $\mathcal{T}$ whose vertexes are the left 
classes of $\Gamma$ with respect to the $\Gamma_n$'s, and such that between two 
vertexes $[g]$ in $\Gamma / \Gamma_n$ and $[h]$ in $\Gamma / \Gamma_{n+1}$ there 
is an edge (oriented from $[g]$ to $[h]$) if $g \!\in\! [h]$. The group $\Gamma$ 
naturally acts on $\mathcal{T}$ preserving the orientation of the edges. By 
Example \ref{arbolito}, if $\Gamma$ has property $(\mathrm{T})$ then there 
exists a vertex $[g]$ which is fixed by this action. If $n \!\in\! \mathbb{N}$
is such that $[g]$ represents the class of $g$ with respect to $\Gamma_n$, 
this implies that $\Gamma = \Gamma_n$, and hence $\Gamma$ is finitely generated.

\begin{small} \begin{ejem} If $\Gamma$ has property $\mathrm{(T)}$ and is amenable,  
then it is finite (see Appendix B for the notion of amenability as well as some notation). 
In particular, the only Abelian groups having Kazhdan's property are the finite ones. 
Actually, we will see below that every finitely generated 
amenable group satisfies the so-called 
{\bf{\em Haagerup's property}}, 
\index{Haagerup's property} 
that is, it acts by isometries of a Hilbert space in a {\bf{\em geometrically proper}} 
way (in the sense that $\|c(g)\|_{\mathcal{H}}$ goes to infinity together with 
$\ell ength(g)$). Clearly, such a group cannot have property (T) unless it is 
finite (see \cite{CJV} for more on this relevant property). 

Let $\mathcal{G}$ be a finite system of generators of $\Gamma$, and let 
$(A_n)$ be a F\o lner sequence associated to $\mathcal{G}$ so that
\begin{equation}
\frac{card(\partial A_n)}{card(A_n)} \leq \frac{1}{2^n}.
\label{folnerexp}
\end{equation}
\index{F\o lner sequence}Let us consider the Hilbert space 
$\mathcal{H} = \oplus_{n \geq 1} \hspace{0.1cm} n \ell^2_{\mathbb{R}}(\Gamma)$
defined by 
$$\mathcal{H} = \Big\{ K = (K_1, \ldots, K_n, \ldots), \hspace{0.1cm}
K_n \in \ell^2_{\mathbb{R}}(\Gamma), \hspace{0.1cm}
\sum_{n \geq 1} n^2 \|K_n\|_{\ell^2_{\mathbb{R}}(\Gamma)}^2 < \infty \Big\}.$$
Given $g,h$ in $\Gamma$, for $K_n \!\in\! \ell^2_{\mathbb{R}}(\Gamma)$ 
let \esp $\Psi_n(g) K_n (h) = K_n(g^{-1} h)$, \esp and let 
$\Psi\!: \Gamma \!\rightarrow\! U(\mathcal{H})$ 
be the {\em regular representation} given by \esp 
$\Psi(g)K \!=\! (\Psi_1(g)K_1, \!\Psi_2(g)K_2, \ldots)$. \esp 
Finally, for each $g \in \Gamma$ let 
$c(g) = (c(g)_1, c(g)_2, \ldots) \in \mathcal{H}$ be defined by 
$$c(g)_n = \frac{1}{\sqrt{card(A_n)}} \big( \mathcal{X}_{A_n} -
\mathcal{X}_{g(A_n)} \big),$$
where $\mathcal{X}$ denotes the characteristic function. By 
(\ref{folnerexp}), for every $g \in \mathcal{G}$ we have
$$\big\| c(g) \big\|_{\mathcal{H}}^2 \esp = \esp 
\sum_{n \geq 1} n^2 \Big\| \frac{\mathcal{X}_{A_n} -
\mathcal{X}_{g(A_n)}}{\sqrt{card(A_n)}} \Big\|_{\ell^2_{\mathbb{R}}(\Gamma)}^2 
\esp \leq \esp \sum_{n \geq 1}n^2 \left( \frac{card(\partial A_n)}{card (A_n)}
\right) \esp \leq \esp \sum_{n \geq 1} \frac{n^2}{2^n} < \infty.$$
The cocycle relation \esp $c(g_1g_2) = c(g_1) + \Psi(g_1)c(g_2)$ \esp can 
be easily checked. We claim that the function 
$c \!: \Gamma \rightarrow \mathcal{H}$ is 
geometrically proper. To show this, fix an integer $k$ bigger than 
the diameter of $A_n^{-1} \!=\! \{h^{-1}\!: h \in A_n\}$. 
If an element $g \in \Gamma$ belongs to 
$B_{\mathcal{G}}(k) \setminus B_{\mathcal{G}}(k-1)$, then 
$g(A_n) \cap A_n = \emptyset$, and therefore 
$$\Big\|\frac{\mathcal{X}_{A_n} -
\mathcal{X}_{g(A_n)}}{\sqrt{card(A_n)}} \Big\|_{\ell^2_{\mathbb{R}}(\Gamma)}^2 \geq 2.$$
Thus, $\|c(g)\|_{\mathcal{H}}^2 \geq 2n^2$, and this shows our claim.
\label{proprecito}
\end{ejem}

\begin{ejem} Group homomorphic images and finite extensions of Kazhdan groups also 
have property $\mathrm{(T)}$. The proof is easy and we leave it to the reader.
\label{easyto}\end{ejem} \end{small}

As a consequence of the example above, every Kazhdan 
group satisfies the hypothesis of Thurston's stability theorem. 
\index{Thurston!stability theorem}
Indeed, the image of a Kazhdan group by a homomorphism into $(\mathbb{R},+)$ is Abelian 
and has property $\mathrm{(T)}$. Thus, it must be finite, and hence trivial.

\begin{small}\begin{ejer}
Show directly that no nontrivial subgroup $\Gamma$ of $(\mathbb{R},+)$ has 
property $(\mathrm{T})$. For this, consider the representation of $\Gamma$ 
by translations on $\mathcal{L}^2_{\mathbb{R}}(\mathbb{R},Leb)$ and the 
associated cocycle \esp 
$c(g) = \mathcal{X}_{[0,\infty[} - \mathcal{X}_{[g(0),\infty[},$ 
\esp where $\mathcal{X}$ denotes the characteristic function and 
$g\in\Gamma$. Show that this cocycle is not a coboundary.
\label{formal-1}
\end{ejer}\end{small}

\begin{small} \begin{ejer} Prove directly that every finitely generated 
subgroup of $\mathrm{Diff}_+^{1+\tau}([0,1])$ with property (T) is 
trivial, where $\tau >0$.

\noindent{\underbar{Hint.}} Passing to a quotient if necessary, one may assume that the 
underlying group $\Gamma$ has no fixed point in $]0,1[$. Consider the Radon measure 
\index{measure!Radon} $\mu$ on $]0,1[$ defined by $d\mu = dx/x$, and let $\Psi$ 
be the {\bf{\em regular representation}} of $\Gamma$ on 
$\mathcal{H} = \mathcal{L}^2_{\mathbb{R}}([0,1],\mu)$ defined by 
$$\Psi(g) K(x) = K(g^{-1}(x)) \Big[ \frac{dg^{-1}}{d\mu}(x) \Big]^{1/2}.$$
For each $g \in \Gamma$ let
\begin{equation}
c(g) \esp = \esp 1 - \Big[ \frac{dg^{-1}}{d\mu}(x) \Big]^{1/2}.
\label{cociclo-de-contacto}
\end{equation}
Check the cocycle relation $c(gh) = c(g) + \Psi(g)c(h)$. 
Using the fact 
that $\tau$ is positive, prove that $c(g)$ belongs to $\mathcal{H}$. Show that 
if $c$ is cohomologically trivial, then there exists $K \in \mathcal{H}$ such that 
the measure $v$ on $]0,1[$ whose density function (with respect to $\mu$) is 
$x \mapsto [1 - K(x)]^2$ is invariant by $\Gamma$. Conclude the proof  
using some of the results from \S \ref{muylargo}.
\label{thurtito}
\end{ejer}
\end{small}

Nontrivial examples of Kazhdan groups are lattices in 
\index{lattice}
(connected) simple Lie groups of (real) rank greater than $1$. 
\index{Lie group!simple} 
(Recall that the (real) {\bf\em{rank}} of a Lie group is the dimension of the 
maximal Abelian subalgebra upon which the adjoint 
representation is diagonalizable (over 
$\mathbb{R}$), and that a discrete subgroup of a locally compact topological 
group is a {\bf{\em lattice}} if the quotient space has finite Haar measure.) 
For instance, $\mathrm{SL}(3,\mathbb{Z})$ has property $\mathrm{(T)}$, since 
the rank of the simple Lie group 
$\mathrm{SL}(3,\mathbb{R})$ is 2, and $\mathrm{SL}(3,\mathbb{Z})$ 
is a lattice inside. For a detailed discussion of this, see \cite{HV}.

Different kinds of groups with property $(\mathrm{T})$ have 
been constructed by many people. In particular, in his seminal work on random groups  
\index{Gromov!random groups} 
\cite{Gromov-random} (see also \cite{Zuk}), Gromov shows 
that ``generic'' finitely presented groups have Kazhdan's 
property. Therefore, a theorem which is true for Kazhdan 
groups is somehow valid for ``almost every group".

To close this section, we give two simple results on obstructions to 
property (T) which are relevant for us. They concern Thompson's group 
G and the Neretin's groups to be defined below. 
\index{Neretin's groups}
The reader who is eager for the connexions with groups of circle 
diffeomorphisms may skip this discussion and pass directly to 
the next section.

\vspace{0.1cm}

\begin{prop} {\em Thompson's group $\mathrm{G}$ does not have Kazhdan's property} (T).
\label{farley}
\end{prop}

\noindent{\bf Proof.} We will give a (modified version of a) nice argument due to 
Farley \cite{Fa}, which actually shows that $\mathrm{G}$ has Haagerup's property 
({\em c.f.,} Example \ref{proprecito}). 
\index{Haagerup's property}\index{Farley} 

Denote by $\mathrm{G}_0$ the subgroup of $\mathrm{G}$ formed by the elements whose 
restrictions to the subinterval $[0,1/2]$ of $\clo$ is the identity. Let us consider 
the Hilbert space $\mathcal{H} = \ell^2_{\mathbb{R}}(\mathrm{G} / \mathrm{G}_0)$. 
The group $\mathrm{G}$ naturally acts by isometries of $\mathcal{H}$ by letting 
\esp $\Psi(g) K ([h]) = K([g^{-1}h])$ \esp for all $g,h$ in G and all 
$K \!\in\! \mathcal{H}$.

Given a dyadic interval $I \subset \clo$, choose $g_I \in \mathrm{G}$ sending 
\index{dyadic!interval} 
$[0,1/2]$ into $I$ affinely. Notice that the class $[g_I]$ of 
$g_I$ modulo $\mathrm{G}_0$ does not depend on this choice.   
For each $g \!\in\! \mathrm{G}$ let $c(g) \!\in\! \mathcal{H}$ 
be the function defined by 
$$c(g) = \sum \left( \delta_{[g_I]} - \delta_{[g g_I]} \right),$$
where the sum extends over the set of all dyadic intervals $I$, 
and $\delta_{[g_I]}$ is the characteristic function of the set 
$\{ [g_I] \}$. Each function $c(g)$ has finite support, since for every 
$g \in \Gamma$ one has \esp $[g_{g(I)}] = [g \esp \! g_I]$ \esp 
for $|I|$ small enough (the restriction of $g$ to small 
dyadic intervals is affine). Therefore, $c(g)$ belongs to $\mathcal{H}$.
Moreover, it is not difficult to see that $c$ satisfies the cocycle 
identity with respect to the unitary representation $\Psi$.
\index{Thompson's groups!G}

To compute $\|c(g)\|$, notice that from the definition one easily deduces 
that $\|c(g)\|^2$ equals two times the number of dyadic intervals $I$ such 
that either the restriction of $g$ to $I$ is not affine, or the image of 
$I$ by $g$ is not a dyadic interval. It easily follows from this fact that 
$\|c(g)\|$ tends to infinity as $\ell ength(g)$ goes to infinity. 
$\hfill\square$

\vspace{0.4cm}

\index{Neretin's groups}
In order to introduce Neretin's groups, let us denote the 
homogeneous simplicial tree of valence $p+1$ 
by $\mathcal{T}^p$. Let $\sigma$ be a marked vertex 
of $\tepe$, which will be called the {\bf {\em origin}}. A (possibly empty) 
subtree $\mathcal{A}$ of $\mathcal{T}^p$ 
\index{tree} 
is {\bf{\em complete}} if it is connected, compact, 
and each time two edges in $\mathcal{A}$ have a common 
vertex, all of the edges containing this vertex are included in 
$\mathcal{A}$. Notice that the complement of a complete subtree 
is either all $\tepe$ or a finite union of rooted trees.

Given a pair of complete subtrees $\mathcal{A}, \mathcal{B}$ of $\mathcal{T}^p$, 
we denote by $\ner(\mathcal{A},\mathcal{B})$ the set of bijections from 
$\overline{\tepe \setminus \mathcal{A}}$ \hspace{0.02cm} onto \hspace{0.01cm} 
$\overline{\tepe \setminus \mathcal{B}}$ sending each connected component of 
$\overline{\tepe \setminus \mathcal{A}}$ isometrically onto a connected component 
of $\overline{\tepe \setminus \mathcal{B}}$. If $g$ belongs to 
$\ner(\mathcal{A},\mathcal{B})$, then $g$ induces a homeomorphism 
of $\partial \tepe$, which we will still denote by $g$. Notice that 
$\partial \tepe$ may be endowed with a natural metric: the distance between 
$x,y$ in $\partial \tepe$ is given by \esp $\db(x,y) =  p^{-n}$, \esp where 
$n$ is the distance \esp $dist$ \esp (on $\tepe$) between 
$\sigma$ and the geodesic joining $x$ and $y$.

\vspace{0.1cm}

\begin{defn} The group of homeomorphisms of $\partial \tepe$ induced by elements 
in some $\mathcal{N}^p(\mathcal{A},\mathcal{B})$ is called {\bf{\em Neretin's group}} 
(or {\bf{\em spheromorphisms group}}), and denoted $\ner$.
\end{defn}

\vspace{0.05cm}

Roughly, $\ner$ is the group of ``germs at infinity'' of isometries of 
$\partial \tepe$. Notice that if $g_1 \!\in\! \mathcal{N}^p(\mathcal{A},\mathcal{B})$ 
and $g_2 \!\in \!\mathcal{N}^p(\mathcal{A'},\mathcal{B'})$ induce the same element of 
$\ner$, then they coincide over $\tepe \setminus (\mathcal{A} \cup \mathcal{A}')$. 
A representative $\tilde{g} \!\in\! \mathcal{N}^p(\mathcal{A},\mathcal{B})$ of 
$g \in \ner$ will be said to be {\bf{\em maximal}} if its domain of definition 
contains the domain of any other representative of $g$. Each element in $\ner$ 
has a unique maximal representative. For $g \!\in\! \ner$, let us denote by 
$\mathcal{A}_g$  (resp. $\mathcal{B}_g$) the closure of the complement of the 
domain of definition of $\tilde{g}$ (resp. of $\tilde{g}^{-1}$). Notice that 
$\mathcal{A}_g \!=\! \mathcal{B}_{g^{-1}}$. The group $\mathrm{Isom}(\tepe)$ 
of the (extensions to the boundary of the) isometries of $\tepe$ is a subgroup of 
$\ner$. An element $g \!\in\! \ner$ comes from an element in 
$\mathrm{Isom}(\tepe)$ if and only 
if $\mathcal{A}_g \!=\! \mathcal{B}_g \!=\! \emptyset$. Notice finally that 
Thompson's group G may be also seen as a subgroup of $\mathcal{N}^p$.

Neretin's groups appear naturally in the  $p$-adic context. 
Indeed, if $p$ is a prime number, then the group of diffeomorphisms 
of the projective line over $\mathbb{Q}_p$ (that is, of the ``$p$-adic circle'') 
embeds into $\ner$. In a certain sense, $\ner$ is a combinatorial analogue of the 
group of circle diffeomorphisms. For further information on this, we refer to 
the nice survey \cite{yuner}. Here we will content ourselves with proving 
the following result from \cite{Na2}.

\vspace{0.15cm}

\begin{prop} {\em Let $\Gamma$ be a subgroup of $\ner$. If $\Gamma$ has 
property $\mathrm{(T)}$, then there exists a finite index subgroup $\Gamma_0$ 
of $\Gamma$ so that the boundary of $\tepe$ decomposes into finitely many balls 
which are fixed by $\Gamma_0$, and $\Gamma_0$ acts isometrically on each of them.}
\label{propner}
\end{prop}

\noindent{\bf Proof.} For each vertex $a \neq \sigma$ let $A_a$ be the subtree 
rooted at $a$ and ``pointing to the infinity'' ({\em i.e.,} in the opposite 
direction to that of the origin). Let us choose one of the $p+1$ subtrees 
of $\tepe$ rooted at the origin, and slightly abusing of the notation,  
let us denote it by $A_{\sigma}$. Let $\ner_{\sigma}$ be the subgroup of $\ner$ 
formed by the elements which fix the boundary at infinity $\partial A_{\sigma}$ 
of $A_{\sigma}$ and act isometrically (with respect to the metric $\db$) 
on it. The group $\ner$ has a natural unitary action $\Psi$ on the Hilbert 
space $\mathcal{H} = \ell^2_{\mathbb{R}}(\ner / \ner_{\sigma})$, namely 
\esp $\Psi(g) K ([h]) = K([g^{-1} h]).$ 

To each vertex $a$ of $\tepe$ we associate the left-class  
$\phi_{a} \in \ner / \ner_{\sigma}$ defined by $\phi_a = [h]$, where 
$h \!\in\! \ner$ is an element whose maximal representative sends $A_{\sigma}$ 
into $A_{a}$ isometrically (with respect to the metric $dist$). Given $g \in \ner$ 
we let
$$c(g) \esp = \esp \sum g(\delta_{_{\phi_a}}) - 
\sum \delta_{_{\phi_b}},$$
where $\delta_{_{\phi_a}}$ is the characteristic function of the set 
$\{ {\phi_a} \}$. Notice that, in the expression above, most of the terms 
cancel, and only finitely many remain; as a consequence, the function 
$c(g)$ belongs to $\mathcal{H}$. Moreover, one easily checks that $c$ 
is a cocycle with respect to $\Psi$.

For $g \in \Gamma \setminus \mathrm{Isom}(\tepe)$, let $d = d(g)$ be the distance 
between $\sigma$ and $\mathcal{A}_g$. Let us consider a geodesic $\gamma$ joining 
$\sigma$ to the vertex in $\mathcal{A}_g$ for which the distance to $\sigma$ 
is maximal, and let 
$a_1,a_2,\ldots,a_{d-1}$ be the vertexes in the interior of $\gamma$. One easily 
checks that for every vertex $b$ in $\tepe$ and all $i \!\in\! \{ 1, \ldots, d\!-\!1 \}$ 
one has $g(\delta_{_{\phi_{a_i}}}) \neq \delta_{_{b}}$, from where one deduces that 
$||c(g)||^2 \geq d-1$. 
 
Let now $g \in \Gamma \cap \mathrm{Isom}(\tepe)$. Let 
$d'  \!=\! d'(g)$ be the distance between $\sigma$ and $\tilde{g}(\sigma)$, and let 
$\gamma$ be the geodesic joining these points. Denoting by $a'_1, \ldots, a'_{d'-1}$ 
the vertexes in the interior of $\tilde{g}^{-1}(\gamma)$, one easily checks that 
for all $i \in \{ 1, \ldots, d'-1 \}$ and every vertex $b$ of $\tepe$ one has 
$g(\delta_{_{\phi_{a'_i}}}) \neq \delta_{_{b}}$ (this is because 
$A_{\tilde{g}(a'_i)}$ does not coincide with the image of $A_{a'_i}$ 
by $\tilde{g}$, since they point toward different directions~!). 
One thus concludes that $||c(g)||^2 \geq d'-1$.

If $\Gamma$ has property $(\mathrm{T})$, then the function $g \mapsto ||c(g)||$ 
must be bounded. Therefore, there exists an integer $N > 0$ such that 
$dist(\sigma,\mathcal{A}_g) \leq N$ for all 
$g \in \Gamma \setminus \mathrm{Isom}(\tepe)$ 
$\big($and hence 
$dist(\sigma,\mathcal{B}_g) \leq N$ for all
$g \in \Gamma \setminus \mathrm{Isom}(\tepe) \big)$, and such that for all 
$g \in \Gamma \cap \mathrm{Isom}(\tepe)$ one has $dist(\sigma,\tilde{g}(\sigma))
\leq N$. The proposition easily follows from these facts. $\hfill\square$

\index{Kazhdan's property (T)|)}


\subsection{The statement of the result}

\hspace{0.45cm} There is a lot of evidence suggesting that any (reasonable) one-dimensional 
structure on a space is an obstruction to actions of Kazhdan group on it. For instance, 
we will see below that for every action of a group with property (T) by isometries of 
a simplicial tree, there exists a global fixed point. The proof that we present 
is essentially due to Haglund, Paulin, and Valette. An easy modification allows 
showing that the same result holds for actions by isometries on real trees. 
The results are originally due to Alperin and Watatani \cite{Al,Wa}.

\begin{small} \begin{ejem} Let $\mathcal{T}$ be an oriented simplicial 
tree (we do not assume that the valence of the edges is finite). Let us denote 
by $\overrightarrow{edg}(\mathcal{T})$ the set of open (oriented) edges of 
$\mathcal{T}$. For each $\vec{\Upsilon} \in \overrightarrow{edg}(\mathcal{T})$, let 
$ver(\vec{\Upsilon})$ be the set of the vertexes in $\mathcal{T}$ which are connected to 
$\vec{\Upsilon}$ by a geodesic whose initial segment is $\vec{\Upsilon}$ (with the corresponding 
orientation). For each vertex $v \in \mathcal{T}$ let us denote by $\overrightarrow{edg}(v)$ 
the set of the oriented edges $\vec{\Upsilon}$ for which $v \in ver(\vec{\Upsilon})$.

Let us consider the Hilbert space 
$\mathcal{H} = \ell^2_{\mathbb{R}}(\overrightarrow{edg}(\mathcal{T}))$. 
Let $\Gamma$ be a subgroup of the group of (orientation preserving) 
isometries of $\mathcal{T}$. Fix a vertex $\sigma$ 
in $\mathcal{T}$, and define a unitary representation 
$\Psi$ of $\Gamma$ on $\mathcal{H}$ by letting \esp 
$\Psi(g)K(\vec{\Upsilon}) \! = \! 
K(g^{-1}(\vec{\Upsilon}))$. \esp For each \esp\! $g \!\in\! \Gamma$ \esp\! let 
\!\esp $c(g) \!: \overrightarrow{edg}(\mathcal{T}) \rightarrow \mathbb{R}$ 
\esp\! be defined by 
$$c(g) = \mathcal{X}_{\overrightarrow{edg}(\sigma)} -
\mathcal{X}_{\overrightarrow{edg}(g(\sigma))},$$
where $\mathcal{X}$ stands for the characteristic function. 
It is easy to see that $c(g)$ has finite support, 
\index{support!of a function} 
and hence belongs to $\mathcal{H}$. Moreover, the 
correspondence $g \mapsto c(g)$ is a cocycle with respect to $\Psi$.

If $\Gamma$ has Kazhdan's property (T), then the cocycle $c$ is a coboundary. In other words, 
there exists a function $K \in \mathcal{H}$ such that for every $g \in \Gamma$ one has \esp 
$c(g) = K - \Psi(g)K$. \esp In particular, $\|c(g)\| \leq 2 \hspace{0.05cm} \|K\|$. On 
the other hand, it is easy to see that \esp $\|c(g)\|^2 = 2 \hspace{0.1cm} dist(\sigma,g(\sigma))$. 
\esp One then deduces that the orbit of $\sigma$ by the $\Gamma$-action stays inside a bounded 
subset of $\mathcal{T}$. We leave to the reader the task of showing that this implies the 
existence of a vertex which is fixed by the action (see \cite{HV,Se} in case of 
problems; compare also Exercise \ref{centro-tits}).  

If we have a group action on a general simplicial tree, then by taking the first barycentric 
subdivision we may reduce to the oriented case. If the group has property (T), then one 
concludes that it fixes either a vertex or the middle point of an edge of the original tree. 
\index{tree}
\label{arbolito}
\end{ejem} \end{small}

\vspace{0.1cm}

In \S \ref{sec-witte}, we have seen that no finite index subgroup of $\mathrm{SL}(3,\mathbb{Z})$ 
acts by homeomorphisms of the line in a nontrivial way. According to \cite{Wi} (see also 
\cite{nota-al-CRAS}), 
\index{Witte-Morris!theorem of non orderability of lattices}this 
is also true for many other lattices in Lie groups having 
property (T), but it is still unknown whether there exists 
a nontrivial Kazhdan group action on the line. (Equivalently, it 
is unknown whether there exist nontrivial, orderable Kazhdan groups.) 
\index{group!orderable}

More accurate results are known in the case of the circle. In particular, 
a theorem of Ghys \cite{Gh2} (a closely related version was independently 
and simultaneously obtained by Burger and Monod \cite{BM}) asserts that if 
$\Phi\!: \Gamma \!\rightarrow\! \mathrm{\mathrm{Diff}}_{+}^1(\mathrm{S}^1)$ 
is a representation of a lattice $\Gamma$ in a simple Lie group of rank greater 
than $1$, then the image $\Phi(\Gamma)$ is finite. For a nice discussion of 
this result (as well as a complete proof for the case of lattices in 
$\mathrm{SL}(3,\mathbb{R})$) we strongly recommend \cite{Gh1}.

\vspace{0.03cm}

\begin{small} \begin{ejer}
The main step in Ghys' proof consists in showing that, 
\index{Ghys!theorem for lattice actions}
for every action by circle {\em homeomorphisms} 
of a lattice $\Gamma$ in a higher-rank simple Lie group, 
there is an invariant probability measure. Show that this 
suffices for proving the theorem (compare Exercise \ref{ejer-parwa}).\\

\vspace{0.08cm}

\noindent{\underbar{Hint.}} Recall that if there is an invariant probability measure, 
then the rotation number function is a group homomorphism into $\mathbb{T}^1$ (see 
\S \ref{medinv}). Using the fact that $\Gamma$ has Kazhdan's property, show that the 
image by this homomorphism is finite. Conclude that the orbit of every point in the 
support of the invariant measure is finite. Finally, apply Thurston's stability theorem.
\label{etiene}
\end{ejer}

\begin{obs} It is not difficult to extend Theorem \ref{wittefuerte} and show that finite 
index subgroups in $\mathrm{SL}(3,\mathbb{Z})$ do not admit nontrivial actions by circle 
homeomorphisms. Indeed, this result (also due to Witte-Morris) may be proved by combining 
the claim at the beginning of the preceding exercise and Theorem \ref{wittefuerte}, though 
it was originally proved using Margulis' normal subgroup theorem \cite{Ma2}. 
\index{Margulis!normal subgroup theorem} 
Let us point out that the proof is quite easy in class $\ce^2$. Indeed, the involved lattices 
contain nilpotent subgroups which are not virtually Abelian. Hence, by Plante-Thurston's 
theorem, their actions by $\ce^2$ circle diffeomorphisms have infinite kernel. By Margulis' 
normal subgroup theorem, these actions must have finite image. 
\end{obs}
\end{small}
\index{Witte-Morris!theorem of non orderability of lattices} 
\index{Plante-Thurston's theorem}

\vspace{0.05cm}

The theorem below, obtained by the author in \cite{Na1} (inspired by \cite{Rez}), 
may be thought of as an analogue for Kazhdan groups of the results discussed 
above under a supplementary regularity condition.

\vspace{0.2cm}

\begin{thm} \textit{Let $\Phi\!: \Gamma \rightarrow \mathrm{\mathrm{Diff}}_{+}^{1+\tau}(\mathrm{S}^1)$ 
be a representation of a countable group, with $\tau > 1/2$. If $\Gamma$ has Kazhdan's property 
{\em (T)}, then $\Phi(\Gamma)$ is finite.} 
\label{principalT}
\end{thm}

\vspace{0.15cm}

Kazhdan's property may be also considered for non discrete groups, mainly for locally 
compact ones. (In this case, the representations and cocycles involved in the definition 
should be continuous functions.) The reader can easily check that the technique of proof of 
Theorem \ref{principalT} still applies in this general context. Another Kazhdan group may 
then arise as a group of circle diffeomorphisms, namely $\mathrm{SO}(2,\mathbb{R})$ 
(compare \S \ref{Lie}).

Notice that Theorem \ref{principalT} shows again that Thompson's group G does not have 
property (T) (see Proposition \ref{farley}). Indeed, according to \S \ref{ghys-sergi}, 
G may be realized as a group of $\ce^{\infty}$ circle diffeomorphisms (actually, it 
suffices to use the $\ce^{1+\mathrm{Lip}}$ realization from \S \ref{sec-th}). 
However, our technique of proof does not lead to Haagerup's 
property for G shown in \S \ref{def-T}.
\index{Haagerup's property}

\vspace{0.02cm}

\begin{small}\begin{ejer} Show that every Kazhdan subgroup of $\mathrm{PSL}(2,\mathbb{R})$ 
is finite by following the steps below (see \cite[\S 2.6]{HV} in case of problems). 

\noindent (i) Consider the action of the M\"obious group on the space M of 
non-oriented geodesics in the Poincar\'e disk $\mathbb{D}$. By identifying M 
with the set of classes $\mathrm{PSL} (2,\mathbb{R}) / G_{\gamma}$ (where 
$G_{\gamma}$ is the stabiliser of a geodesic $\gamma$), use the Haar 
measure on $\mathrm{PSL} (2,\mathbb{R})$ to endow M with a Radon measure $v$.
\index{measure!Radon}

\vspace{0.07cm}

\noindent (ii) For each $P,Q$ in $\mathbb{D}$, let M$_{P,Q}$ be the 
set of geodesics which ``separate'' $P$ and $Q$ ({\em i.e.}, which 
divide $\mathbb{D}$ into two hemispheres so that $P$ and $Q$ do not 
belong to the same one). Show that there exists a constant $C$ such that 
$v$-measure of M$_{P,Q}$ equals \esp $C dist(P,Q)$ \esp for all \esp 
$P,Q$, where $dist$ stand for the hyperbolic distance on $\mathbb{D}$. 

\vspace{0.07cm}

\noindent (iii) Consider the Hilbert space $\mathcal{H} = \mathcal{L}^2_{\mathbb{R}} (M,v)$  
and the unitary action $\Psi$ of $\mathrm{PSL}(2,\mathbb{R})$ on $\mathcal{H}$ given 
by \esp $\Psi(f) K (\gamma) = K (f^{-1} \gamma)$. \esp Show that the function 
$c$ defined on $\mathrm{PSL}(2,\mathbb{R})$ by 
$c(f) \esp = \esp \mathcal{X}_{\mathrm{M}_{O,f(O)}}$ 
belongs to $\mathcal{H}$, where $O = (0,0) \in \mathbb{D}$ and $\mathcal{X}$ stands for 
the characteristic function. Show moreover that $c$ is a cocycle with respect to $\Psi$.
\index{M\"obius group}

\vspace{0.07cm}

\noindent (iv) Using the construction in (iii), show that if $\Gamma$ is a Kazhdan 
subgroup of $\mathrm{PSL}(2,\mathbb{R})$, then all of its elements are elliptic. 
Conclude that $\Gamma$ is finite by using the result from Exercise \ref{asi-no-mas}. 
\label{se-acabo}
\end{ejer}\end{small}


\subsection{Proof of the theorem}
\label{pruebabacan}

\hspace{0.45cm} To simplify, we will denote by $\bar{g}$ the diffeomorphism 
$\Phi(g^{-1})$. Recall that the Liouville measure $Lv$ on 
\index{measure!Liouville}
$\clo \times \clo$ has density function (see \S \ref{liouv}) 
$$(r,s) \mapsto \frac{1}{ 4 \esp \mathrm{sin}^2 \big( \frac{r-s}{2} \big)}.$$ 

Let $\mathcal{H}=\mathcal{L}^{2,\Delta}_{\mathbb{R}}(\clo \times \clo,Lv)$ be 
the subspace of $\mathcal{L}^2_{\mathbb{R}}(\mathrm{S}^1 \times \mathrm{S}^1,Lv)$ 
formed by the functions $K$ satisfying \esp $K(x,y) = K(y,x)$ \esp for almost 
every $(x,y) \in \clo \times \clo$. Let $\Psi$ be the unitary 
representation of $\Gamma$ on $\mathcal{H}$ given by 
$$\Psi(g)K (r,s) = K(\bar{g}(r),\bar{g}(s))
\cdot [Jac(\bar{g})(r,s)]^{\frac{1}{2}},$$
where $Jac(\bar{g})(r,s)$ is the Jacobian (with respect to $Lv$) 
of the map $(r,s) \mapsto (\bar{g}(r),\bar{g}(s))$: 
$$ Jac(\bar{g})(r,s) =
\frac{\mathrm{sin}^2 \left( \frac{r-s}{2} \right) }
{\mathrm{sin}^2 \left( \frac{\bar{g}(r)-\bar{g}(s)}{2} \right) }
\cdot \big[ \bar{g}'(r) \bar{g}'(s) \big].$$
For each $g \in \Gamma$ let us consider the function (compare (\ref{cochina}))
\begin{equation}
c(g)(r,s) = 1 - [Jac(\bar{g})(r,s)]^{\frac{1}{2}}.
\label{coc}
\end{equation}
One may formally check the relation 
$c(g_1g_2) \!=\! c(g_1) \!+\! \Psi(g_1)c(g_2)$ 
(compare (\ref{cocicloschwarz})). Indeed, this {\bf {\em Liouville cocycle}} $c$ 
corresponds to the ``formal coboundary'' of the constant function $1$, which does 
not belong to $\mathcal{L}^{2}_{\mathbb R}(\clo \times \clo,Lv)$.\footnote{The 
reader will readily notice that the cocycles from Examples \ref{proprecito}, 
\ref{formal-1}, \ref{thurtito}, and \ref{arbolito}, that of Exercise 
\ref{se-acabo}, and those of the proofs of Propositions \ref{farley} 
and \ref{propner}, also arise as formal coboundaries.} The main issue here is 
that, if $\Phi$ takes values in $\mathrm{\mathrm{Diff}}_{+}^{1+\tau}(\mathrm{S}^1)$ 
for some $\tau \!>\! 1/2$, then 
$c(g)$ is a 
true cocycle with values in $\mathcal{H}.$ This is the content of the proposition 
below, essentially due to Segal \cite{PS} (see also \cite{Rez}).
\index{cocycle!Liouville}

\vspace{0.2cm}

\begin{prop} \textit{If $\tau \!>\! 1/2$, then $c(g)$ belongs to 
$\mathcal{L}^{2,\Delta}_{\mathbb R}(\mathrm{S}^1 \times \mathrm{S}^1,Lv)$ 
for all \esp $g \!\in\! \Gamma$.}
\end{prop}

\noindent{\bf Proof.} Notice that for a certain continuous function 
$K_1: [0,2\pi] \times [0,2\pi] \rightarrow \mathbb{R}$ one has 
$$\frac{1}{\left| \mathrm{sin} \left( \frac{r-s}{2} \right) \right|}
= 2 \left[ \frac{1}{|r-s|} + K_1(r,s) \right].$$
Therefore, to prove that $c(g)$ is in 
$\mathcal{L}^{2,\Delta}_{\mathbb R}(\mathrm{S}^1 \times \mathrm{S}^1, Lv)$,
we need to verify that the function 
$$(r,s) \esp \mapsto \esp 
\frac{[\bar{g}'(r)\bar{g}'(s)]^{\frac{1}{2}}}{|\bar{g}(r)-\bar{g}(s)|}
- \frac{1}{r-s}$$
belongs to $\mathcal{L}^{2}_{\mathbb{R}}(\mathrm{S}^1
\times \mathrm{S}^1,Leb).$ Now for all $r,s$ in $\mathrm{S}^1$ such 
that $|r-s|< \pi$, there exists $t \in \clo$ in the shortest segment 
joining them so that $|\bar{g}(r) - \bar{g}(s)| = \bar{g}'(t)|r-s|$. 
We then have 
\begin{eqnarray*}
\left|\! \frac{[\bar{g}'(r)\bar{g}'(s)]^{\frac{1}{2}}}{|\bar{g}(r)\!-\!\bar{g}(s)|}
- \frac{1}{|r\!-\!s|} \!\right|
&=& 
\frac{1}{|r-s|\bar{g}'(t)} \cdot
\left| [\bar{g}'(r)\bar{g}'(s)]^{\frac{1}{2}} - \bar{g}'(t)\right| \\
&=& 
\frac{1}{|r-s|\bar{g}'(t)} \cdot \frac{|\bar{g}'(r)\bar{g}'(s) -
\bar{g}'(t)^2|}{[\bar{g}'(r)\bar{g}'(s)]^{\frac{1}{2}}
+ \bar{g}'(t)}.\\
&\leq& 
\frac{1}{2 \inf (\bar{g}')^2 |r\!-\!s|} 
\big[ |\bar{g}'(r) \!-\! \bar{g}'(t)|\bar{g}'(s) +
\bar{g}'(t)|\bar{g}'(s) \!-\! \bar{g}'(t)| \big].
\end{eqnarray*}
Since $\bar{g}'$ is $\tau$-H\"older continuous 
({\em c.f.,} Example \ref{wholder}), this gives  
\index{H\"older!derivative}
$$\left| \frac{[\bar{g}'(r)\bar{g}'(s)]^{\frac{1}{2}}}{\bar{g}(r)-\bar{g}(s)} -
\frac{1}{r-s}\right| \leq \frac{|\bar{g}'|_{\tau} \sup (\bar{g}')}{2|r-s|
\inf (\bar{g}')^2} \big[ |r-t|^{\tau} + |s-t|^{\tau} \big] 
\leq C \esp |r-s|^{\tau - 1},$$
where the constant $C$ does not depend on $(r,s)$. The proof is finished by noticing 
that, since $\tau \!>\! 1/2$, the function $(r,s) \!\mapsto\! |r-s|^{\tau - 1}$ belongs 
to $\mathcal{L}^2_{\mathbb{R}}(\mathrm{S}^1 \times \mathrm{S}^1,Leb)$. $\hfill\square$

\vspace{0.5cm}

If the group $\Gamma$ has property (T) and $\tau > 1/2$, then the cocycle 
(\ref{coc}) is a coboundary. In other words, there exists a function 
$K \in \mathcal{L}^{2,\Delta}_{\mathbb{R}}(\clo \times \clo,Lv)$ such 
that, for every $g \in \Gamma$ and almost every 
$(r,s) \in \clo \times \clo$, one has 
$$1 - [Jac(\bar{g})(r,s)]^{\frac{1}{2}} = K(r,s) - K(\bar{g}(r),\bar{g}(s))
\cdot [Jac(\bar{g})(r,s)]^{\frac{1}{2}},$$
that is,
$$ [ 1 - K(\bar{g}(r),\bar{g}(s))]^2 \cdot Jac(\bar{g})(r,s) = [1-K(r,s)]^2.$$
We thus conclude the following.

\vspace{0.2cm}

\begin{prop} \textit{Let $\Phi\!: \Gamma \rightarrow \mathrm{Diff}_{+}^{1+\tau}(\mathrm{S}^1)$
be a representation, with $\tau>1/2.$ If $\Gamma$ has property $\mathrm{(T)}$, then there 
exists $K \in \mathcal{L}^{2,\Delta}_{\mathbb{R}}(\clo \times \clo,Lv)$ such that $\Gamma$ 
preserves the geodesic current $L_K$ given by \index{geodesic current}}
$$\frac{d \hspace{0.1cm} L_K}{d \hspace{0.1cm} Lv}
= \big[ 1 - K(r,s) \big]^2.$$
\label{medida invariante}
\end{prop}

Since $L_K$ is absolutely continuous with respect to the Lebesgue measure on 
${\V}^1 \times {\V}^1 \setminus \Delta$, the property below is evident:
\begin{equation}
L_K([a,a] \times [b,c]) = 0  \hspace{0.15cm} \mbox{ for 
all } \hspace{0.15cm} a < b \leq c < a. \label{cero}
\end{equation}
On the other hand, the fact that $K$ is an square integrable 
function implies that 
\begin{equation}
L_K \left( [a,b[ \times ]b,c] \right) =
\infty \hspace{0.15cm} \mbox{ for all } \hspace{0.15cm} a < b < c < a.\\
\label{infinito}
\end{equation}
Indeed, noticing that 
$$L_K([a,x]\times[y,c])^{\frac{1}{2}} \esp = \esp 
   \left( \int_a^x\!\!\int_y^c [ 1 \!-\! K(r,s)]^2 \hspace{0.05cm}
   d \hspace{0.05cm} Lv  \right)^{\frac{1}{2}}
   \esp \geq \esp \left( \int_a^x\!\!\int_y^c \!
   \frac{dr \hspace{0.1cm}ds}{4 \hspace{0.05cm} \mathrm{sin}^2(\frac{r-s}{2})}
   \right)^{\frac{1}{2}} \!\!- \|K\|_2,$$
equality (\ref{infinito}) easily follows from the relation 
$$\int_a^x\!\!\int_y^c \frac{dr \hspace{0.1cm} ds}
{4 \hspace{0.05cm}\mathrm{sin}^2(\frac{r-s}{2}) } \esp 
= \esp \log \big( [e^{ia},e^{ix},e^{iy},e^{ic}] \big)$$
and from the fact that the value of the cross-ratio \esp 
$[e^{ia},e^{ix},e^{iy},e^{ic}]$ \esp goes to infinity as 
$x$ and $y$ tend to $b$ \hspace{0.01cm} (with $a<x<b<y<c$).

We will say that a geodesic current satisfying properties (\ref{cero}) 
and (\ref{infinito}) is {\bf{\em stable}}. Proposition 
\ref{medida invariante} implies the following.

\vspace{0.2cm}

\begin{prop} \textit{Let $\Phi: \Gamma \rightarrow
\mathrm{Diff}^{1+\tau}_{+}(\mathrm{S}^1)$ be a representation,
with $\tau > 1/2$. If $\Gamma$ has property $\mathrm{(T)}$, then there 
exists a stable geodesic current that is invariant by $\Gamma$.} \label{rigidez}
\end{prop}

\vspace{0.2cm}

The measure $L_K$ is not necessarily fully supported, that is, there may be nontrivial 
intervals $[a,b]$ and $[c,d]$ for which \esp $L_K([a,b] \times [c,d]) = 0.$ \esp This 
may lead thinking that the group of circle homeomorphisms preserving $L_K$ is not 
necessarily well behaved. Nevertheless, we will see that the stability properties 
of $L_K$ lead to rigidity properties for this group.

\vspace{0.1cm}

\begin{lem} \textit{If a circle homeomorphism preserves a stable geodesic 
current and fixes three different points, then it is the identity.}
\label{puntosfijos}
\end{lem}

\noindent{\bf Proof.} Suppose that a homeomorphism \esp $f \!\neq\! Id$ \esp fixes at 
least three points and preserves a stable geodesic current $L$, and let $I \!= ]a,b[$ 
be a connected component of the complement of the set of fixed points of $f$. Notice 
that $a$ and $b$ are fixed points of $f$. Let $c \! \in ]b,a[$ be another fixed point 
of $f$. Since $f$ does not fix any point in $]a,b[$, for each $x \! \in ]a,b[$ the 
sequence $(f^i(x))$ converges to either $a$ or $b$. Both cases being analogous, 
let us only consider the second one. Then $f^{-i}(x)$ converges to $a$ 
as $i$ goes to infinity. This yields 
$$L \big( [a,x] \times [b,c] \big) =
L \big( [a,f(x)] \times [b,c] \big),$$
and therefore $L \big( [x,f(x)] \times [b,c] \big) \!=\! 0$.
Since $x \! \in ]a,b[$ was arbitrary, 
$$L \big( ]a,b[ \times [b,c] \big) =
\sum_{i \in \mathbb{Z}}
L \big( [f^i(x),f^{i+1}(x)] \times [b,c] \big) =0.$$
However, this is in contradiction with (\ref{infinito}). $\hfill\square$

\vspace{0.35cm}

The next proposition is a direct consequence of the preceding one. The reason 
for using the notation $\Phi(g)$ instead of $\bar{g}$ here will be clear shortly.

\vspace{0.15cm}

\begin{prop} \textit{Let $\Phi: \Gamma \rightarrow
\mathrm{Diff}^{1+\tau}_{+}(\mathrm{S}^1)$ be a representation,
with $\tau > 1/2$. If $\Gamma$ has property $\mathrm{(T)}$ and 
$g \in \Gamma$ is such that $\Phi(g)$ fixes three points, 
then $\Phi(g)$ is the identity.}
\label{trespuntosfijos}
\end{prop}

\vspace{0.15cm}

The pretty argument of the proof of the next proposition 
was kindly communicated to the author by Witte-Morris. 
\index{Witte-Morris!3-fold covering trick}

\vspace{0.2cm}

\begin{prop} \textit{Let $\Phi: \Gamma \rightarrow
\mathrm{Diff}^{1+\tau}_{+}(\mathrm{S}^1)$ be a representation,
with $\tau > 1/2$. If $\Gamma$ has property 
$\mathrm{(T)}$ and $g \in \Gamma$ is such that $\Phi(g)$ has a 
fixed point, then $\Phi(g)$ is the identity.}
\label{unpuntofijo}
\end{prop}
\vspace{-0.4cm}
\noindent{\bf Proof.} Let us consider the 3-fold circle covering $\hat{\mathrm{S}}^1$. 
Upon this covering acts (by $\mathrm{C}^{1+\tau}$ diffeomorphisms) a degree-3 
central extension of $\Gamma$. More precisely, there exists a group 
\index{degree!of a group extension}
$\hat{\Gamma}$ containing a (central) normal subgroup isomorphic to 
$\mathbb{Z} / 3 \mathbb{Z}$, such that the quotient is isomorphic to $\Gamma$. 
Since $\Gamma$ has property (T) and $\mathbb{Z} / 3\mathbb{Z}$ is finite, 
$\hat{\Gamma}$ also has property ($\mathrm{T}$) 
({\em c.f.,} Example \ref{easyto}). If $g \in \Gamma$ is such that $\Phi(g)$ 
fixes a point in the original circle, then one of its preimages in $\hat{\Gamma}$
fixes three points in $\hat{\mathrm{S}}^1$ by the induced action. Since $\hat{\mathrm{S}}^1$ 
identifies with the circle, using the preceding proposition one easily deduces that 
$\Phi(g)$ is the identity. $\hfill\square$

\vspace{0.5cm}

Is is now easy to complete the proof of Theorem \ref{principalT}. Indeed, by the proposition above, 
the action of the group $\Phi(\Gamma)$ on $\clo$ is free. By H\"older's theorem, this group is Abelian. 
\index{H\"older!theorem} \index{action!free}But since $\Phi(\Gamma)$ still satisfies property (T), 
it is forced to be finite.

\vspace{0.12cm}

\begin{small} \begin{obs} By Exercise \ref{conjcorestable}, every group of circle homeomorphisms 
preserving a stable geodesic current is conjugate to a subgroup of the M\"obius group. 
\index{M\"obius group} Therefore, 
if $\tau \!>\! 1/2$ and $\Phi\!: \Gamma \rightarrow \mathrm{Diff}_{+}^{1+\tau}(\clo)$ is a 
representation whose associate cocycle (\ref{coc}) is a coboundary, then the image $\Phi(\Gamma)$ 
is topologically conjugate to a subgroup of $\mathrm{PSL}(2,\mathbb{R})$. This allows 
giving an alternative end of proof for Theorem \ref{principalT} by 
using the result from Exercise \ref{se-acabo}.
\label{pillo}
\end{obs}

\begin{ejer} Prove that for every function 
$K \!\in\! \mathcal{L}^{2,\Delta}_{\mathbb{R}}(\clo \times \clo,Lv)$,    
the group $\Gamma_{L_K}$ of the circle homeomorphisms preserving 
\index{quasisymmetric homeomorphism} 
$L_K$ is uniformly quasisymmetric (see Exercise \ref{quasisim}).
\end{ejer}

\begin{obs} A result from \cite{lp} allows extending Theorem \ref{principalT} 
to actions by $\ce^{3/2}$ diffeomorphisms. The theorem is perhaps true even  
for actions by $\ce^1$ diffeomorphisms (compare Exercise \ref{thurtito} and 
the comments before it; compare also Remark \ref{thu-na-rem}). 
However, an extension to general actions by {\em homeomorphisms} is 
unclear: see for instance Example \ref{interesnia}. A particular case which 
is interesting by itself is that of piecewise affine circle homeomorphisms.  
\index{piecewise affine homeomorphism} 
\label{interesnita}
\end{obs}

\begin{obs} By combining Theorem \ref{principalT} with the two-dimensional version of Thurston's  
stability theorem \cite{Th2} one may show that, if $\tau \!>\! 1/2$, then every countable group 
\index{Thurston!stability theorem}
of $\ce^{1+\tau}$ diffeomorphisms of the closed disk or the closed annulus with property (T) is 
finite. It is unknown whether this is true for the open disk and/or annulus, as well as for 
compact surfaces of non-positive Euler characteristic (notice that the sphere $\mathrm{S}^2$ 
supports a faithful action of the Kazhdan group $\mathrm{SL}(3,\mathbb{Z})$). 
\label{pierre}\end{obs}

\begin{ejer} Following \cite{He}, let us denote $\mathrm{D}^{\infty} =
\mathrm{Diff}^{\infty}_+(\clo) \setminus int(\rho^{-1}(\mathbb{Q}/\mathbb{Z}))$, 
where $\rho$ is the rotation number function. 
\index{rotation number}
The space $\mathrm{D}^{\infty}$ is closed in $\mathrm{Diff}_+^{\infty}(\clo)$, and 
hence it is a Baire space. Prove that the infinite cyclic group generated by a 
generic $g \in \mathrm{D}^{\infty}$ is not uniformly quasisymmetric, and thus 
the set $\{\|c(g^n)\|\!: n \in \mathbb{Z}\}$ is unbounded. (Recall that a property 
is {\bf {\em generic}} when it is satisfied on a $G_{\delta}$-dense set; recall 
moreover that, in a Baire space, countable intersections of $G_{\delta}$-dense 
sets is still a $G_{\delta}$-dense set).

\vspace{0.02cm}

\noindent{\underbar{Hint.}} The set of $g \in \mathrm{D}^{\infty}$ satisfying 
\esp $\sup_{[a,b,c,d] = 2} \hspace{0.17cm} \sup_{n \in \mathbb{N}}
\hspace{0.25cm} [g^n(a),g^n(b),g^n(c),g^n(d)] = \infty \label{unica}$ 
\esp is a $G_{\delta}$-set. Show that this set contains all infinite 
order circle homeomorphism of rational rotation number 
having at least three periodic points. Then use the fact 
that these homeomorphisms are dense in $\mathrm{D}^{\infty}$.

\vspace{0.02cm}

\noindent{\underbar{Remark.}} According to Remark \ref{des-denj-conv}, 
for every $g \!\in\! \mathrm{D}^{\infty}$ of irrational rotation 
number, the sequence $(g^{q_n})$ converges to the identity 
in the $\ce^1$ topology, where $p_n/q_n$ is the $n^{\mathrm{th}}$-rational 
approximation of $\rho(g)$. Moreover, if $k \geq 2$, then for a generic subset 
of elements $g$ in $\mathrm{D}^{\infty}$ there exists an increasing sequence 
of integers $n_i$ such that $g^{n_i}$ converges to the identity in the $C^k$ 
topology. Indeed, letting $dist_k$ be a metric inducing the $\ce^k$ topology on 
$\mathrm{Diff}_+^k(\clo)$, the set of $g \!\in\! \mathrm{D}^{\infty}$ satisfying 
$$\liminf_{n > 0} \hspace{0.15cm} dist_{k} (g^n,Id) = 0$$
is a $G_{\delta}$-set; moreover, this set contains all $g$ for which 
$\rho(g)$ verifies a Diophantine condition \cite{He,KO,Yo1}. One thus 
concludes that, for a generic $g \!\in\! \mathrm{D}^{\infty}$, the sequence 
$( \|c(g^n)\| )$ is unbounded, but it has the zero vector as an accumulation point. 
Furthermore, the orbits of the corresponding isometry $A(g) = \Psi(g) + c(g)$ 
are unbounded and {\bf{\em recurrent}}, in the sense that all of their points 
are accumulation points (see \cite{c-t-v} for more on unbounded, recurrent 
actions by isometries on Hilbert spaces).
\label{rec-frappend}
\end{ejer}

\begin{ejer} In class $\ce^{1+\mathrm{Lip}}$, Kopell's lemma 
\index{Kopell!lemma}
may be shown by using the Liouville cocycle. More precisely, let $f,g$ be 
commuting $\ce^{1+\mathrm{Lip}}$ diffeomorphisms of $[0,1]$. Consider the 
``Liouville measure'' $\overline{Lv}$ over $[0,1[ \times [0,1[$ whose density 
function is $(r,s) \mapsto 1/(r-s)^2$.

\noindent (i) Show that the functions $c(f)$, $c(g)$ defined by 
$$c(f)(r,s) = 1 - \big[ Jac(f^{-1})(r,s) \big]^{\frac{1}{2}}, \qquad 
c(g)(r,s) = 1 - \big[ Jac(g^{-1})(r,s) \big]^{\frac{1}{2}},$$
are square integrable over all compact subsets of $[0,1[ \times [0,1[$. 

\noindent (ii) Show that, for every $a \! \in ]0,1[$, there exists a constant 
$C = C(a,f)$ such that for every interval $[b,c]$ contained in $[0,a]$ one has 
$$\int_b^c \int_b^c \big\| c(f) \big\|^2 \esp \leq \esp C \esp \! |c-b|^2.$$

\noindent (iii) Suppose that $g$ fixes a point $a \!\in ]0,1[$ and that $f(x)\!<\!x$ 
for every $x \!\in ]0,1[$. Using the preceding inequality and the relation  \esp 
$c(g^k) = c(f^{-n} g^k f^n) = c(f^{-n}) + \Psi(f^{-n})c(g^k) + \Psi(f^{-n}g) c(f^n),$ 
\esp conclude that the value of 
$$\int_{f^2(a)}^a \int_{f^2(a)}^{a} \big\| c(g^k) \big\|^2 $$
is uniformly bounded (independently of $k \in \mathbb{N}$).

\noindent (iv) By an argument similar to that of the proof of 
Proposition \ref{trespuntosfijos}, show that the restriction of $g$ to 
$[f^2(a),a]$ \esp (and hence to the whole interval $[0,1[$) \esp is the identity.
\label{kopell-liouv}
\end{ejer}

\begin{ejer} Let $\mu$ be the (finite) measure on the boundary of $\tepe$ giving mass 
$p^{-n}$ to each ball of radius $p^{-n}$, where $n \geq 1$ (see \S \ref{def-T}). Given 
a homeomorphism $g$ of $\partial \tepe$, let us denote by $g'$ its Radon-Nikodym 
derivative with respect to  
\index{Radon-Nikodym derivative} 
$\mu$ (whenever it is defined). A natural ``Liouville measure'' $Lv$ on 
\esp $\partial \tepe \! \times \! \partial \tepe$ \esp is given by 
$$d Lv(x,y) = \frac{d\mu(x) \times d\mu(y)}{\db(x,y)^2}.$$
Consider the Hilbert space $\mathcal{H}$ formed by the functions 
$K$ in $\mathcal{L}^2_{\mathbb{R}}(\partial \tepe\!\times\!\partial \tepe,Lv)$
satisfying \esp $K(x,y) = K(y,x)$ \esp  
\index{Neretin's groups}
for almost every $(x,y) \in \partial \tepe\!\times\!\partial \tepe$. Consider 
the unitary representation $\Psi$ of $\ner$ on $\mathcal{H}$ given by 
$$\Psi(g) K (x,y) \esp = \esp 
K ( g^{-1} (x), g^{-1} (y) ) \cdot [ Jac (g^{-1})(x,y) ]^{1/2},$$
where $Jac(g^{-1})(x,y)$ denotes the Jacobian (with respect to the 
Liouville measure) of the map $(x,y) \mapsto (g^{-1}(x),g^{-1}(y))$.

\noindent(i) Prove that the extension to the boundary of every isometry $f$ of $\tepe$ 
satisfies the equality \esp $\db(x,y) f'(x)f'(y) \!=\! \db ( f(x),f(y) )$
\esp for all points $x, y$ in $\partial \tepe$.\\

\noindent(ii) Conclude that for each $g \in \ner$ the function 
$c(g)\!: \partial \tepe\!\times\!\partial \tepe \rightarrow \mathbb{R}$ 
given by 
$$c(g) (x,y) = 1 - \big[ Jac(g^{-1})(x,y) \big]^{1/2}$$
belongs to $\mathcal{H}$.\\

\noindent(iii) Show that $c$ satisfies the cocycle relation with respect to $\Psi$.\\

\noindent(iv) Prove that if $\Gamma$ is a subgroup of $\ner$ satisfying property 
$(\mathrm{T})$, then there exists $K \in \mathcal{H}$ such that $\Gamma$ 
preserves the ``geodesic current'' $L_K$ given by 
$$dL_K = \big[ 1 - K(x,y) \big]^2 d Lv.$$

\noindent(v) Conclude that there exists a compact subset 
$\mathcal{C} = \mathcal{C}(K)$ of $\tepe$ such that, for 
every element $g \in \Gamma \setminus \mathrm{Isom}(\tepe)$, 
either $\mathcal{A}_g \cap \mathcal{C} \neq \emptyset$ 
or $\mathcal{B}_g \cap \mathcal{C} \neq \emptyset$.\\

\noindent(vi) Use (v) to give an alternative proof for Proposition \ref{propner}.
\end{ejer} 
\end{small}


\subsection{Relative property (T)\index{Kazhdan's property (T)!relative} 
and Haagerup's property}

\hspace{0.45cm} Motivated by the preceding section, it is natural to ask whether finitely 
generated subgroups of $\mathrm{Diff}_+^{1+\tau}(\clo)$ necessarily satisfy Haagerup's 
property when $\tau > 1/2$ ({\em c.f.,} Example \ref{proprecito}). This is for instance 
the case of discrete subgroups of $\mathrm{PSL}(2,\mathbb{R})$ (see 
Exercise \ref{se-acabo}) as well 
as of Thompson's group G (see the proof of Proposition \ref{farley}). One of the main 
difficulties for this question relies on the fact that very few examples of groups 
satisfying neither Haagerup's nor Kazhdan's property are known. Indeed, most of the 
examples of groups without Haagerup's \index{Haagerup's property} property actually 
satisfy a weak version of property (T), namely the relative property (T) defined below.

\vspace{0.1cm}

\begin{defn}
If $\Gamma$ is a locally compact group and $\Gamma_0$ is a subgroup of $\Gamma$, then 
the pair $(\Gamma,\Gamma_0)$ has the {\bf {\em relative property}} (T) if for every 
(continuous) representation by isometries of $\Gamma$ on a Hilbert space, there 
exists a vector which is invariant by $\Gamma_0$.
\end{defn}

\vspace{0.1cm}

A relevant example of a pair satisfying the relative property (T) 
is $(\mathrm{SL}(2,\mathbb{Z}) \ltimes \mathbb{Z}^2, \esp \mathbb{Z}^2)$. The reader 
will find more examples, as well as a discussion of this notion, in \cite{HV,CJV}. 
Let us point out, however, that in most of the examples in the literature, 
if neither $\Gamma$ nor $\Gamma_0$ have property (T), 
then $\Gamma_0$ (contains a cocompact subgroup which) is 
normal in $\Gamma$ (see however \cite{cornulier}). Under such a hypothesis, the 
result below may be considered as a small generalization of Theorem \ref{principalT}. 
Its interest relies on Example \ref{interesnia}.

\vspace{0.25cm}

\begin{thm} {\em Let $\Gamma$ be a subgroup of $\mathrm{Diff}_+^{1+\tau}(\clo)$, with 
$\tau\!>\!1/2$. Suppose that $\Gamma$ has a normal subgroup $\Gamma_0$ such that the 
pair $(\Gamma,\Gamma_0)$ has the relative property} (T). {\em Then either $\Gamma$ 
is topologically conjugate to a group of rotations or $\Gamma_0$ is finite.}
\index{rotation!group} 
\label{no-se-extiende}
\end{thm}

\noindent{\bf Proof.} Let us use again the technique of the proof of 
Theorem \ref{principalT}. 
The Liouville cocycle induces an isometric representation of $\Gamma$ on
$\mathcal{L}^{2,\Delta}_{\mathbb{R}}(\clo \times \clo,Lv)$. If $(\Gamma,\Gamma_0)$ 
has the relative property (T), then this representation admits a $\Gamma_0$-invariant 
vector, and the arguments of the preceding section (see Remark \ref{pillo}) show that 
$\Gamma_0$ is topologically conjugate to a subgroup of the 
\index{M\"obius group} M\"obius group.

Relative property (T) is stable under finite, central extensions. Using this 
fact, and by means of the 3-fold covering trick in the proof of Proposition 
\ref{unpuntofijo}, one easily shows that $\Gamma_0$ is actually topologically 
\index{Witte-Morris!3-fold covering trick}
conjugate to a subgroup of the group of rotations (see Remark \ref{pillo}).
If $\Gamma_0$ is infinite, then it is conjugate to a {\em dense} subgroup 
of the group of rotations. Now it is easy to check that the normalizer in 
$\mathrm{Homeo}_+(\clo)$ of every dense subgroup of $\mathrm{SO}(2,\mathbb{R})$ 
coincides with the whole group of rotations 
(compare Exercise \ref{uni-inv}). Therefore, if $\Gamma_0$ is 
infinite, then $\Gamma$ is topologically conjugate to a subgroup of 
$\mathrm{SO}(2,\mathbb{R})$, which concludes the proof. $\hfill\square$

\vspace{0.05cm}

\begin{small}\begin{ejem}
Theorem \ref{no-se-extiende} does not extend to actions by homeomorphisms:  
the group $\mathrm{SL}(2,\mathbb{Z}) \ltimes \mathbb{Z}^2$ acts faithfully 
by circle homeomorphisms. Indeed, starting from the canonical action of 
$\mathrm{PSL}(2,\mathbb{R})$, one can easily realize $\mathrm{SL}(2,\mathbb{R})$ 
as a group of real-analytic circle diffeomorphisms. Let $p \!\in\! \clo$ be a point 
whose stabilizer under the corresponding 
$\mathrm{SL}(2,\mathbb{Z})$-action is trivial. Replace each 
point $f(p)$ of the orbit of $p$ by an interval $I_f$ (where 
$f \in \mathrm{SL}(2,\mathbb{Z})$) in such a way that the total sum of these 
\index{action!faithful} 
intervals is finite. Doing this, we obtain a topological circle $\clo_p$ provided 
with a faithful $\mathrm{SL}(2,\mathbb{Z})$-action (we use affine transformations 
for extending the maps in $\mathrm{SL}(2,\mathbb{Z})$ to the intervals $I_f$). 

Let $I \!=\! I_{id}$ be the interval corresponding to the point $p$, and let 
$\{\varphi^t \! : t \in \mathbb{R}\}$ be a nontrivial topological flow on $I$. 
Choose any real numbers $u,v$ which are linearly independent over the rationals, 
and let $g = \varphi^u$ and $h = \varphi^v$. Extend $g,h$ to $\clo_p$ by letting 
$$g(x) = f^{-1} \big( g^a h^c (f(x)) \big), \qquad 
h(x) = f^{-1} \big( g^b h^d (f(x)) \big),$$ 
for 
\begin{equation*}
f = \left(
\begin{array}
{cc}
a & b  \\
c & d  \\
\end{array}
\right) \esp \in \esp \mathrm{SL}(2,\mathbb{R})
\end{equation*}
and $x \!\in\! I_{f^{-1}}$.   
(For $x$ in the complement of the union of the $I_f$'s, we simply 
put $g(x)=h(x)=x$.) The reader will easily check that the group generated 
by $g,h$, and the copy of $\mathrm{SL}(2,\mathbb{Z})$ acting on $\clo_p$, 
is isomorphic to $\mathrm{SL}(2,\mathbb{Z}) \ltimes \mathbb{Z}^2$.
\label{interesnia}
\end{ejem}

\begin{ejem}
If $\mathbb{F}_2$ is a free subgroup of finite index in $\mathrm{SL}(2,\mathbb{Z})$, 
then the pair $( \mathbb{F}_2 \ltimes \mathbb{Z}^2, \mathbb{Z}^2)$ still has the relative 
property (T). Moreover, $\mathbb{F}_2 \ltimes \mathbb{Z}^2$ 
is an orderable group (it is actually locally indicable: 
\index{group!locally indicable} 
see Remark \ref{en-general}), thus it acts faithfully by homeomorphisms of the interval. 
(Notice that no explicit action arises as the restriction of the action constructed in 
the preceding example; however, a faithful action may be constructed by following a 
similar procedure.) Since $\mathbb{F}_2 \ltimes \mathbb{Z}^2$ is not bi-orderable, let 
us formulate the following question: Does there exist a finitely generated, bi-orderable 
\index{group!bi-orderable} 
group without Haagerup's property~?
\index{group!orderable}
\end{ejem}

\begin{obs} The actions constructed in the preceding examples are not $\ce^1$ smoothable. 
Actually, for every non solvable subgroup $H$ of $\mathrm{SL}(2,\mathbb{Z})$, the group 
$H \ltimes \mathbb{Z}^2$ does not embed in neither $\mathrm{Diff}^1_+(\clo)$ nor 
$\mathrm{Diff}^1_+(\mathbb{R})$. The proof of this fact does not rely on 
cohomological properties, and it is mostly based on dynamical methods: 
see \cite{na-th}.
\label{thu-na-rem}
\end{obs}
\end{small}


\section{Super-rigidity for Higher-Rank Lattice Actions}

\subsection{Statement of the result}
\index{super-rigidity theorem!for actions on the circle|(}

\hspace{0.45cm} Kazhdan's property fails to hold for some higher-rank 
{\em semisimple} Lie groups and their lattices. Indeed, some of these groups may 
act by circle diffeomorphisms. However, these actions are quite particular, as was 
first shown by Ghys in \cite{Gh3}. To state Ghys' theorem properly, recall 
\index{Ghys!theorem for lattice actions} 
that a lattice $\Gamma$ in a Lie group $G$ is said to be {\bf {\em irreducible}} if 
there are no normal subgroups $G_1$ and $G_2$ generating $G$ so that $G_1 \cap G_2$ 
is contained in the center of $G$ (which is supposed to be finite) and such that 
$(\Gamma \cap G_1) \cdot (\Gamma \cap G_2)$ has finite index in $\Gamma$.
\index{Lie group!semisimple}\index{lattice}

\vspace{0.1cm}

\begin{thm} {\em Let $G$ be a (connected) semisimple 
Lie group of rank greater than or equal to $2$, and let 
$\Gamma$ be an irreducible lattice in $G$. If $\Phi$ is a homomorphism from 
$\Gamma$ into the group of $\mathrm{C}^1$ circle diffeomorphisms, then either 
$\Phi$ has finite image, or the associated action is semiconjugate to a finite 
covering of an action obtained as the composition of the following homomorphisms:}

\vspace{0.08cm}

\noindent{-- {\em the inclusion of \esp $\Gamma$ in $G$,}}

\vspace{0.08cm}

\noindent{-- {\em a surjective homomorphism from $G$ into 
$\mathrm{PSL}(2,\mathbb{R})$},}

\vspace{0.08cm}

\noindent{-- {\em the natural action of $\mathrm{PSL}(2,\mathbb{R})$ 
on the circle.}}
\label{gros-et}
\end{thm}

\vspace{0.12cm}

For the proof of this result, Ghys starts by examinating the case 
of some ``standard'' higher-rank semisimple Lie groups 
($\mathrm{SL}(n,\mathbb{R})$, $\mathrm{Sp}(2r,\mathbb{R})$,
$\mathrm{SO}(2,q)$, $\mathrm{SU}(2,q)$, and 
$\mathrm{PSL}(2,\mathbb{R}) \times \mathrm{PSL}(2,\mathbb{R})$), 
and then he uses the 
classification of general semisimple Lie groups \cite{knap}. 
Notice that the former four cases correspond to higher-rank {\em simple} 
Lie groups (the involved lattices satisfy Kazhdan's property...). 
For the latter case (which is dynamically more interesting) Ghys proves 
that, up to a semiconjugacy and a finite covering, every homomorphism 
$\Phi \!: \Gamma \rightarrow \mathrm{Diff}_+^1(\clo)$ factors through the 
projection of $\Gamma$ into one of the factors, and then by the projective 
action of this factor on the circle. To understand this case better, 
the reader should have in mind as a fundamental example the embedding 
of $\mathrm{PSL}(2,\mathbb{Z}(\sqrt{2}))$ as a lattice in 
$\mathrm{PSL}(2,\mathbb{R}) \times \mathrm{PSL}(2,\mathbb{R})$ 
through the map
$$  \left(
\begin{array}
{cc}
\!\!a_1 \!+\! b_1 \sqrt{2} \!\!&\!\! a_2 \!+\! b_2 \sqrt{2} \!\! \\
\!\!a_3 \!+\! b_3 \sqrt{2} \!\!&\!\! a_4 \!+\! b_4 \sqrt{2} \!\! \\
\end{array}
\right) \!\longmapsto\!
\left(\! \left(
\begin{array}
{cc}
\!\!a_1 \!+\! b_1 \sqrt{2} \!\!&\!\! a_2 \!+\! b_2 \sqrt{2}  \\
\!\!a_3 \!+\! b_3 \sqrt{2} \!\!&\!\! a_4 \!+\! b_4 \sqrt{2}  \\
\end{array}
\right)\!,
\! \left(
\begin{array}
{cc}
\!\!a_1 \!-\! b_1 \sqrt{2} \!\!&\!\! a_2 \!-\! b_2 \sqrt{2}  \\
\!\!a_3 \!-\! b_3 \sqrt{2} \!\!&\!\! a_4 \!-\! b_4 \sqrt{2}  \\
\end{array}
\right) \!\right).$$

To extend Theorem \ref{gros-et}, we need to generalize the notion of ``rank'' 
to an arbitrary lattice. Many tempting definitions have been already proposed 
(see for instance \cite{eber}). Here we will deal with the perhaps simplest 
one, which has been successfully exploited (among others) by Shalom \cite{Sha}. 
The higher-rank hypothesis corresponds in this general framework to a 
commutativity hypothesis inside the ambient group. In more concrete 
terms, the ``general setting'' that we consider --which is actually 
that of \cite{Sha}-- is the following one:
\index{lattice}

\vspace{0.1cm}

\noindent -- $G = G_1 \times \cdots \times G_k$ 
is a locally compact, {\em compactly generated}  
\index{group!compactly generated} topological group (see Exercise 
\ref{compactamente-generado}), with $k \geq 2$, and $\Gamma$ 
is a finitely generated cocompact lattice inside;

\vspace{0.1cm}

\noindent -- the projection $pr_i(\Gamma)$ of $\Gamma$ into each factor $G_i$ is dense;

\vspace{0.1cm}

\noindent -- in the case where each $G_i$ is an algebraic linear group over a local field 
\cite{Ma2}, we also allow the possibility that $\Gamma$ be a non-cocompact lattice  in $G$.

\vspace{0.1cm}

We point out that, in the last case, $\Gamma$ is necessarily finitely generated. 
This follows from important theorems due to Kazhdan and Margulis \cite{Ma2}. 
\index{Margulis!normal subgroup theorem} 
On the other hand, the second condition is similar to 
the irreducibility hypothesis in Theorem \ref{gros-et}.

In the Introduction of \cite{Sha}, the reader may find many other motivations, as 
well as some relevant references, concerning the general setting that we consider. 
We should mention that examples of {\em nonlinear} lattices satisfying the first 
two properties have been constructed in \cite{BeGl,BMoz,remi}. For these lattices, 
as well as for the {\em linear} ones, ({\em i.e.}, those which embed into 
$\mathrm{GL}(n,\mathbb{K})$ for some field $\mathbb{K}$), the next 
super-rigidity theorem for actions on the circle (obtained by the 
author in \cite{Nasup}; see also \cite{Na-cociclo}) holds.

\vspace{0.1cm}

\index{Thurston!Plante-Thurston's theorem|see{Plante-Thurston's theorem}}

\begin{thm} {\em In the preceding context, let 
$\Phi\!: \Gamma \rightarrow \mathrm{Diff}_+^{1+\tau}(\clo)$ 
be a homomorphism such that $\Phi(\Gamma)$ does not preserve any probability 
measure on the circle. If $\tau > 1/2$, then either $\Phi(\Gamma)$ is conjugate 
to a subgroup of $\mathrm{PSL}(2,\mathbb{R})$, or $\Phi$ is semiconjugate 
to a finite covering of an action obtained as a composition of:}

\vspace{0.07cm}

\noindent{-- {\em the injection of $\Gamma$ in $G$,}}

\vspace{0.07cm}

\noindent{-- {\em the projection of $G$ into one of the factors $G_i$,}}

\vspace{0.07cm}

\noindent{-- {\em an action $\Phi$ of $G_i$ by circle homeomorphisms.}}
\label{primero1}
\end{thm}

\vspace{0.1cm}

The hypothesis that $\Phi(\Gamma)$ does not fix any probability measure 
on $\clo$ may be suppressed if the first cohomology (with real values) 
of every finite index normal subgroup of $\Gamma$ is trivial. Let us 
point out that, according to \cite{Sha}, this condition is fulfilled if 
$\mathrm{H}^1_{\mathbb{R}}(G)$ is trivial (this is the case for instance if 
the $G_i$'s are semisimple algebraic linear groups over local fields \cite{Ma2}).

\vspace{0.1cm}

\begin{cor} {\em Let $\Gamma$ be a finitely generated 
lattice satisfying the hypothesis of the general setting, and let 
$\Phi \!: \Gamma \rightarrow \mathrm{Diff}_+^{1+\tau}(\clo)$ be a homomorphism, 
with $\tau > 1/2$. If $\mathrm{H}^1_{\mathbb{R}}(\Gamma_0) = \{0\}$ for every 
finite index normal subgroup of $\Gamma$, then the conclusion of 
Theorem} \ref{primero1} {\em still holds}.
\label{corolarito2}
\end{cor}

\vspace{0.1cm}

Thanks to the results above, to understand the actions of irreducible lattices 
in higher-rank semisimple Lie groups by circle diffeomorphisms we may use the 
classification from \S \ref{Lie}. We can then obtain refined 
versions under any one of the hypotheses below:

\vspace{0.1cm}

\noindent{(i) the kernel of $\Phi$ is finite and the orbits by $\Phi(\Gamma)$ 
are dense;

\vspace{0.1cm}

\noindent{(ii) the kernel of $\Phi$ is finite and $\Phi$ takes 
values in the group of real-analytic circle diffeomorphisms;

\vspace{0.1cm}

\noindent{(iii) every normal subgroup of $\Gamma$ either is finite or has finite 
index (that is, $\Gamma$ satisfies Margulis' normal subgroup theorem 
\index{Margulis!normal subgroup theorem} 
\cite{Ma2}).}

\vspace{0.1cm}

\begin{thm} {\em Suppose that the hypothesis of Corollary} \ref{corolarito} 
{\em are satisfied, and that at least one of the preceding 
hypotheses} (i), (ii), {\em or} (iii), {\em is satisfied. If the image $\Phi(\Gamma)$ 
is infinite then, up to a topological semiconjugacy and a finite covering, $\Phi(\Gamma)$ 
is a non metabelian subgroup of $\mathrm{PSL}(2,\mathbb{R})$.}
\label{ultimito}
\end{thm}

\vspace{0.1cm}

\index{group!metabelian}
Hypothesis (i), (ii), or (iii), allows avoiding the degenerate case where $G$ 
has infinitely many connected components and its action on the circle factors, 
modulo topological semiconjugacy, through the quotient by the 
connected component of the identity. Let us 
point out that hypothesis (iii) is satisfied when each $G_i$ is a linear 
algebraic simple group, as well as for the lattices inside products of groups 
of isometries of trees constructed in \cite{BMoz} (see \cite{Bshalom} for a 
general version of this fact).

One can show that the hypothesis of connectedness for $G$ in Theorem \ref{ultimito} 
can be weakened into that no $G_i$ is discrete. This is not difficult but rather 
technical: see \cite{burger-nuevo} for a good discussion on this point.

The proofs we give for the results of this section strongly use 
a super-rigidity theorem for isometric actions on Hilbert spaces 
(due to Shalom) which will be discussed in the next section. This theorem 
has been generalized in \cite{lp} to isometric actions on $\mathcal{L}^p$ spaces. This 
allows obtaining $\ce^{1+\tau}$ versions of the preceding results for every $\tau \!>\! 0$. 
Actually, by means of rather different techniques, these results have been 
recently extended to actions by homeomorphisms by Bader, Furman, and Shaker 
in \cite{ba-gen} (using ``boundary theory''), and by Burger in \cite{burger-nuevo} 
(using bounded cohomology). However, to avoid overloading the presentation, we 
will not discuss these extensions here.


\subsection{Cohomological super-rigidity}

\hspace{0.45cm} Let us again consider a unitary action $\Psi$ of a group 
$\Gamma$ on a Hilbert space $\mathcal{H}$, where $\Gamma$ now may be non 
discrete, but it is supposed to be locally compact and 
compactly generated.
\index{Shalom}
\index{super-rigidity theorem!for reduced cohomology}

\begin{defn} One says that $\Psi$ {\bf \em{almost has invariant vectors}} 
if there exists a sequence of unitary vectors 
$K_n \!\in \mathcal{H}$ such that, for every compact subset $\mathcal{C}$ of 
$\Gamma$, the value of \esp $\sup_{g \in \mathcal{C}} \| K_n - \Psi(g)K_n \|$ 
\esp converges to zero as $n$ goes to infinity.
\end{defn}

\begin{defn} A cocycle $c\!: \Gamma \rightarrow \mathcal{H}$ associated to $\Psi$ is an 
{\bf {\em almost coboundary}} (or it is {\bf{\em almost cohomologically trivial}}) if 
there exists a sequence of coboundaries $c_n$ such that, for every compact subset 
$\mathcal{C}$ of $\Gamma$, the value of \esp $\sup_{g \in \mathcal{C}} \| c_n(g) - c(g) \|$ 
\esp converges to zero as $n$ goes to infinity.
\end{defn}

\begin{small} \begin{ejer} Show that the cocycle from Example \ref{thurtito} 
is almost cohomologically trivial for every finitely generated subgroup of  
$\mathrm{Diff}_+^{1+\tau}([0,1])$, where $\tau > 0$. 

\noindent{\underbar{Hint.}} Consider the coboundary $c_n$ associated to the 
function $K_n(x) = \mathcal{X}_{[1/n,1]}(x)$.
\end{ejer} \end{small}

\index{coboundary}

The elementary lemma below, due to Delorme \cite{HV}, appears to be fundamental 
for studying almost coboundaries. 
\index{Delorme's lemma}

\vspace{0.14cm}

\begin{lem} {\em If $\Psi$ does not almost have fixed vectors, then every 
cocycle which is an almost coboundary is actually cohomologically trivial.}
\end{lem}

\noindent{\bf Proof.} Let $\mathcal{G}$ be a compact generating set for $\Gamma$. 
By hypothesis, there exists $\varepsilon > 0$ such that, for every 
$K \in \mathcal{H}$,
\begin{equation}
\sup_{h \in \mathcal{G}} \| K - \Psi(h)K \| \geq \varepsilon \| K \|.
\label{guich}
\end{equation}
Since $c$ is almost cohomologically trivial, there must 
exist a sequence $(K_n)$ in $\mathcal{H}$ such that 
\esp $c(g) \!=\! \lim_{n \rightarrow +\infty}(K_n - \Psi(g)K_n)$ 
\esp for every $g \!\in\! \Gamma$.
\esp Inequality (\ref{guich}) then gives \esp 
$M=\sup_{h\in\mathcal{G}} \| c(h) \| \geq \varepsilon \lim\sup_n \| K_n \|$,
\esp and thus \esp \esp $\lim\sup_n \|K_n\| \leq M/\varepsilon$. \esp 
Therefore, 
$$\big\| c(g) \big\| \leq \limsup_{n \rightarrow \infty} 
\big( \|K_n\| \!+\! \|\Psi(g) K_n\| \big) \leq 
\frac{2M}{\varepsilon}.$$ 
Since this holds for every $g \!\in\! \Gamma$, 
the cocycle $c$ is uniformly bounded. By the 
\index{Tits!Center Lemma}
Tits Center Lemma ({\em c.f.}, Exercise \ref{centro-tits}), 
it is cohomologically trivial. $\hfill\square$

\vspace{0.35cm}

We now give a version of Shalom's super-rigidity theorem \cite{Sha}, which plays a central 
role in the proof of Theorem \ref{primero1}. Let us point out, however, that we will not use 
this result in its full generality. Indeed, in our applications we will always reduce to 
the case where the corresponding unitary representations do not almost have invariant 
vectors, and for this case the super-rigidity theorem becomes much more elementary.

\vspace{0.15cm}

\begin{thm} {\em Let $G \!=\! G_1 \times \cdots \times G_k$ be a 
locally compact, compactly generated 
topological group, and let $\Gamma$ be a lattice inside satisfying the hypothesis of 
our general setting. Let $\Psi \! : \Gamma \!\rightarrow\! U(\mathcal{H})$ be a unitary 
representation which does not almost have invariant vectors. If $c$ is a cocycle associated 
to $\Psi$ which is not an almost coboundary, then $c$ is cohomologous to a cocycle  
$c_1\!+\!\ldots\!+\!c_k$ such that each $c_i$ is a cocycle taking values in a 
$\Psi(\Gamma)$-invariant subspace $\mathcal{H}_i$ upon which the isometric action 
$\Psi + c_i$ continuously extends to an isometric action of $G$ 
which factors through $G_i$.}
\label{sup-sh}
\end{thm}

\vspace{0.1cm}

This remarkable result was obtained by Shalom inspired by the proof of Margulis' normal 
subgroup theorem. It is based on the general principle that commuting isometries of 
Hilbert spaces are somehow ``degenerate" (Exercises \ref{mon} and \ref{fati} well 
illustrate this fact). Instead of providing a proof (which may be found in \cite{Sha}), 
we have preferred to include two examples where one may appreciate (some of) its 
consequences ({\em c.f.,} Exercises \ref{sh-ejer} and \ref{liouv-red}). For the 
first one we give a useful elementary lemma for extending group 
homomorphisms from a lattice to the ambient group. 
To state it properly, recall that a topological group $H$ is
{\bf{\em sequentially complete}} if every sequence $(h_n)$
in $H$ such that \hspace{0.02cm} 
$\lim_{m,n\rightarrow +\infty}h_m^{-1}h_n=id_{H}$
\hspace{0.02cm} converges to a limit in $H$.
\index{group!sequentially complete}

\vspace{0.12cm}

\begin{lem} {\em Let $G$ and $\Gamma$ be groups satisfying the hypothesis 
of our general setting. Let $\Phi\!: \Gamma \rightarrow H$ be a group homomorphism, 
where $H$ is a sequentially complete, Hausdorff, topological group. Suppose that 
there exists $i \!\in\! \{1,\ldots,k\}$ such that for every sequence $(g_n)$ in 
$\Gamma$ satisfying $\lim_{n \rightarrow +\infty} pr_i(g_n) = id_{G_i}$, one has 
$\lim_{n \rightarrow +\infty} \Phi(g_n) = id_{H}$. Then $\Phi$ extends to a 
continuous homomorphism from $G$ into $H$ which factors through $G_i$.}
\label{extender}
\end{lem}

\noindent{\bf Proof.} For $g \!\in\! G_i$ take an arbitrary sequence $(g_n)$ 
in $\Gamma$ so that $pr_i(g_n)$ converges to $g$. By hypothesis, letting 
$h_n \!=\! \Phi(g_n)$ we have \esp 
$\lim_{m,n\rightarrow +\infty} h_m^{-1}h_n \!=\! id_{H}$. \esp Let us define 
$\hat{\Phi}(g) \!=\! \lim_{n \rightarrow +\infty} h_n$. This definition 
is pertinent, because of the hypothesis on the topology of $H$. Moreover, 
it does not depend on the chosen sequence, and it is easy to see that the 
map $\hat{\Phi}$ thus defined is a continuous homomorphism from $G$ which 
factors through $G_i$. We leave the details to the reader. $\hfill\square$

\vspace{0.15cm}
\index{action!non-elementary (on a tree)}

\begin{small} \begin{ejer} Assuming Theorem \ref{sup-sh}, and following the 
steps below, show the following super-rigidity theorem for actions on trees: 
If $\Gamma$ is a lattice satisfying the hypothesis of the general setting 
and $\Phi$ is a non-elementary action of $\Gamma$ by isometries of a simplicial 
tree $\mathcal{T}$, then there exists a $\Gamma$-invariant subtree over which 
the action extends to $G$ and factors through one of the $G_i$'s (recall that the 
action is {\bf{\em non-elementary}} if there is no vertex, edge, or {\em point 
at the infinity}, which is fixed). 

\noindent \underbar{Remark.} The result above (contained in 
\index{super-rigidity theorem!for actions on trees} \cite{Sha}) 
still holds for \esp {\em non-cocompact} \esp irreducible lattices. 
The proof for this general case appears in \cite{monod-shalom}, 
\index{tree}
and uses a super-rigidity theorem for bounded cohomology.
\index{super-rigidity theorem!for bounded cohomology}

\vspace{0.08cm}

\noindent (i) By means of a barycentric subdivision, reduce the general 
case to that where no element fixes an edge and interchanges its vertexes.

\vspace{0.08cm}

\noindent(ii) Suppose that there exist a vertex $v_0$ in $\mathcal{T}\!$ and an index 
$i \! \in \! \{1,\ldots,k\}$ verifying the following condition: for every sequence 
$(g_n)$ in $\Gamma$ such that \esp $\lim_{n \rightarrow +\infty} pr_i(g_n) = id_{G_i}$, 
\esp the vertex $v_0$ is a fixed point of $\Phi(g_n)$ for every $n$ large enough. 
Prove that the claim of the theorem is true.

\noindent{\underbar{Hint.}} The set of vertexes verifying the above condition is 
contained in a $\Gamma$-invariant tree over which Lemma \ref{extender} applies.

\vspace{0.08cm}

\noindent (iii) Using Delorme's lemma, prove that the hypothesis that the action is 
non-elementary implies that the associated regular representation $\Psi$ of $\Gamma$ 
on \esp $\mathcal{H} \!=\! \ell^2_{\mathbb{R}}(\overrightarrow{edg}(\mathcal{T}))$ 
does not almost have invariant vectors (see Example \ref{arbolito}).

\vspace{0.08cm}

\noindent (iv) Let $\mathcal{H}_i$ be the subspace given by Theorem \ref{sup-sh}. 
Define $edg^* (\mathcal{T})$ as the set of edges in $\mathcal{T}$ 
modulo the equivalence relation which identifies $\vec{\Upsilon}_1$ with 
$\vec{\Upsilon}_2$ if \esp $\psi(\vec{\Upsilon}_1) = \psi(\vec{\Upsilon}_2)$ 
\esp for every function $\psi$ in 
$\mathcal{H}_i$. Consider the subset $edg^*_0(\mathcal{T})$ of the classes in 
$edg^*(\mathcal{T})$ over which at least one of the functions in $\mathcal{H}_i$ 
is non-zero. Show that every class in $edg^*_0(\mathcal{T})$ is finite.

\vspace{0.08cm}

\noindent (v) Let $[\vec{\Upsilon]} \in edg^*_0(\mathcal{T})$ be one of the finite classes 
above, and let $\Gamma^*$ be its stabiliser in $\Gamma$. The subgroup $\Gamma^*$ fixes 
the geometric center of the set of edges in $[\vec{\Upsilon]}$, and thus (i) implies 
that $\Gamma^*$ fixes a vertex $v_0$. Show that $v_0$ satisfies the hypothesis in (ii).

\noindent{\underbar{Hint.}} Given a sequence $(g_n)$ in $\Gamma$ such that 
$\lim_{n \rightarrow +\infty} pr_i(g_n) = id_{G_i}$, show that $g_n$ belongs to 
$\Gamma^*$ for all $n$ large enough. For this, argue by contradiction and consider 
separately the cases where the classes of the $g_n$'s with respect to $\Gamma^*$ are 
equal or different, keeping in mind the fact that the functions in $\mathcal{H}_i$ 
are square integrable and non-zero over the edges in $[\vec{\Upsilon]}$.
\label{sh-ejer}
\end{ejer}

\begin{ejer} Recall that there exist nontrivial commuting $\ce^{1+\tau}$ 
diffeomorphisms $f,g$ of the interval $[0,1]$, where $1/2 \!<\! \tau \!<\! 1$, 
so that the set of fixed points of $f$ (resp. $g$) in $]0,1[$ is empty (resp. 
non-empty and discrete): see \S \ref{ejemplosceuno}. Show 
that the restriction of the Liouville cocycle 
$c\!: \Gamma \!\rightarrow \! 
\mathcal{L}^{2,\Delta}_{\mathbb{R}}([0,1] \times [0,1],\overline{Lv}) = \mathcal{H}$ 
to the group $\Gamma \!\sim\! \mathbb{Z}^2$ 
generated by $f$ and $g$ is almost cohomologically 
trivial (see Exercise \ref{kopell-liouv}).

\vspace{0.08cm}

\noindent{\underbar{Hint.}} Assuming the opposite, Shalom's super-rigidity theorem 
provides us with a unitary vector $K$ in $\mathcal{H}$ which is invariant either by 
$\Psi (f)$ or $\Psi (g)$. Show that the probability measure $\mu_K$ on $[0,1]$ 
defined by 
$$\mu_K(A) = \int_0^1 \int_A |K(x,y)|^2 \thinspace d \overline{Lv},$$
is invariant by either $f$ or $g$. 
Then using the fact that $f$ has no fixed point, show that $g$ preserves $\mu_K$,  
and that this measure is supported on the set of fixed points of $g$. Finally, 
use the fact that the latter set is countable to obtain a contradiction.

\vspace{0.08cm}

\noindent{\underbar{Remark.}} It would be interesting to give an 
explicit sequence of coboundaries converging to the cocycle $c$.
\label{liouv-red}
\end{ejer}

\begin{ejer} Given a product $G = G_1 \times G_2$ of compactly 
generated topological groups, let \esp $A = \Psi + c$ \esp be a 
representation of $G$ by isometries of a Hilbert space $\mathcal{H}$. 
Suppose that the isometric representation of $G_1$ obtained by restriction 
does not have invariant vectors. Show that the corresponding 
unitary representation of $G_2$ almost has invariant vectors.

\vspace{0.08cm}

\noindent{\underbar{Hint.}} Let $(g_n)$ be a sequence in $G_1$ such that $\| c(g_n) \|$ 
tends to infinity. Using the commutativity between $G_1$ and $G_2$, show that the 
sequence of unitary vectors $c(g_n) / \| c(g_n) \|$ is almost invariant by $\Psi(G_2)$.
\label{mon}
\end{ejer}

\begin{ejer}
Let $\Psi \!: \mathbb{Z}^2 \rightarrow U(\mathcal{H})$ be a unitary representation and 
$c\!: \mathbb{Z}^2 \rightarrow \mathcal{H}$ an associated cocycle. Show that, if 
both $\Psi \big( (1,0) \big)$ and $\Psi \big( (0,1) \big)$ do not almost have 
invariant vectors, then $c$ is cohomologically trivial.

\vspace{0.08cm}

\noindent{\underbar{Hint.}} Denote by $A_1=\Psi_1 + c_1$ and $A_2=\Psi_2 + c_2$
the isometries of $\mathcal{H}$ associated to the generators of $\mathbb{Z}^2$. 
Check that for all $i \! \in \! \{1,2\}$, all $K \in \mathcal{H}$, and all 
$n \in \mathbb{Z}$, 
$$(Id - \Psi_i) A_i^{n} (K) = \Psi^{n}_i (K) - \Psi^{1+n}_i (K) - \Psi^n_i (c_i) + c_i.$$
Using this relation and the commutativity between $A_1$ and $A_2$, find a uniform 
upper bound for \esp 
$\| (Id - \Psi_1) (Id - \Psi_2) A_1^{n_1} A_2^{n_2} (K) \|$. \esp 
From the fact that both $\Psi_1$ and $\Psi_2$ do not almost have invariant 
vectors, conclude that the orbits by the isometric action of $\mathbb{Z}^2$
\index{Tits!Center Lemma}
are bounded. Finally, apply the Tits Center Lemma 
({\em c.f.,} Exercise \ref{centro-tits}).
\label{fati}
\end{ejer}

\begin{ejer} Given $C \geq 0$, a subset $\mathrm{M}_0$ of a metric space 
M is {\bf{\em $C$-dense}} if for every $K \in \mathrm{M}$ there exists 
$K_0 \in \mathrm{M}_0$ such that $dist (K_0,K) \leq C$.

\vspace{0.08cm}

\noindent (i) Show that there is no isometry of a Hilbert space of dimension bigger 
than 1 having $C$-dense orbits (compare with the remark in Exercise \ref{rec-frappend}).

\noindent{\underbar{Hint.}} Following an argument due to Fathi, 
\index{Fathi} 
let $A \!=\! \Psi \!+\! c$ be an isometry of a Hilbert space $\mathcal{H}$ with 
a $C$-dense orbit $\{ A^n(K_0)\!: n \in \mathbb{Z} \}$. Using the equality  
\esp $A^n(K_0) = \Psi^n(K_0) + \sum_{i=0}^{n-1} \Psi^i(c)$, \esp show that 
the set \esp $\{ (Id \!-\! \Psi) A^n(K_0)\!: n \in \mathbb{Z} \}$ \esp 
is $2 C$-dense in the space $(Id \!-\! \Psi) \mathcal{H}$. 
Using the identity  
$$ (Id - \Psi) A^n (K_0) \esp 
= \esp \Psi^n (K_0) - \Psi^{1+n}(K_0) - \Psi^n(c) + c$$
conclude that the norms of the vectors in this space are bounded from above by 
\esp $2\big(\| K_0 \| + \| c \|\big)$, \esp and hence $\Psi = Id$. Thus, $A$ 
is a translation, which implies that $\mathcal{H}$ has dimension 0 or 1. 

\noindent (ii) As a generalization of (i), show that no action of $\mathbb{Z}^k$ 
by isometries of a finite dimensional Hilbert space admits $C$-dense orbits.

\noindent{\underbar{Hint.}} Let $A_i \!=\! \Psi_i + c_i$ be commuting isometries of 
an infinite dimensional Hilbert space $\mathcal{H}$, where $i \in \{ 1, \ldots, k \}$. 
Assuming that the orbit of $K_0 \in \mathcal{H}$ by the group generated by them is 
$C$-dense, show that the set 
$$\{ (Id - \Psi_1) \cdots (Id - \Psi_k)
A_k^{n_k} \cdots A_1^{n_1}(K_0) \! : n_{i_j} \in \mathbb{Z} \}$$
is $2^k C$-dense in the space $(Id-\Psi_1) \cdots (Id-\Psi_k) \mathcal{H}$.
Using a similar argument to that of Exercise \ref{fati}, show that this set 
is contained in the ball in $\mathcal{H}$ centered at the origin with  
radius \esp $2^k \big( \| K_0 \| + \sum_{i=1}^k \| c_i \| \big)$. \esp 
Deduce that for at least one of the $\Psi_i$'s, say $\Psi_k$, 
the set $\bar{\mathcal{H}}$ of invariant vectors 
is an infinite dimensional subspace of $\mathcal{H}$. 
Denoting by $\bar{c}_k$ the projection of $c_k$ into $\bar{\mathcal{H}}$, 
show that the orthogonal projection into the space 
$\mathcal{H}^* = \langle \bar{c}_k \rangle^{\bot} 
\cap \bar{\mathcal{H}}$ induces $k-1$ commuting isometries  
generating a group with $C$-dense orbits. Complete the proof 
by means of an inductive argument.

\noindent{\underbar{Remark.}} It is not difficult to construct isometries of 
an infinite dimensional Hilbert space generating a free group for which all 
the orbits are dense. The problem of determining what are the finitely generated 
groups which may act minimally on such a space seems to be interesting. According 
to (ii) above, such a group cannot be Abelian, and by the next exercise, it cannot   
be nilpotent either (compare \cite{con-chu}).
\label{fat}
\end{ejer}

\begin{ejer} Suppose that there exist finitely generated nilpotent  
\index{group!nilpotent} 
groups admitting isometric actions on infinitely dimensional Hilbert 
spaces with $C$-dense orbits for some $C>0$. Fix such a group $\Gamma$ 
with the smallest possible nilpotence degree, and consider the 
corresponding affine action $A \!=\! \Psi \!+\! c$. 

\noindent (i) Show that the restriction of the unitary action $\Psi$ 
to the center $H$ of $\Gamma$ is trivial. 
\index{center of a group}

\noindent{\underbar{Hint.}} Consider a point $K_0$ in the underlying Hilbert 
space $\mathcal{H}$ with a $C$-dense orbit, and for each $K \!\in\! \mathcal{H}$ 
choose $g \!\in\! \Gamma$ such that $\|A(g)(K_0) - K\| \leq C$. For  
all $h \in H$ one has 
\begin{eqnarray*}
\| A(h)(K) \!-\! K \| &\leq& \|A(h)(K) \!-\! A(hg)(K_0)\| \esp 
+ \esp \|A(hg)(K_0)\!-\!A(g)(K_0)\| \esp + \esp \|A(g)(K_0)\!-\!K\|\\
&\leq&  2C \esp + \esp \|A(gh)(K_0)\!-\!A(g)(K_0)\|\\
&=&  2C \esp + \esp \|A(h)(K_0) \!-\! K_0\|.
\end{eqnarray*}
In particular, if $K \neq 0$, then replacing $K$ by $\lambda K$ 
in the above inequality and letting $\lambda$ tend to infinity, 
conclude that $K$ is invariant by $H$.

\noindent (ii) By (i), for every $h \in H$ the isometry $A(h)$ is a translation, 
say by a vector $K_h$. Show that $K_h$ is invariant under $\Psi(\Gamma)$. 
Conclude that the subspace $\mathcal{H}_0$ formed by the $\Psi(\Gamma)$-invariant 
vectors is not reduced to $\{0\}$ (notice that the action of $H$ on 
$\mathcal{H}_0$ cannot be 
trivial). By projecting orthogonally into $\mathcal{H}_0$ and its orthogonal 
complement $\mathcal{H}_0^{\bot}$, we obtain isometric representations $A_0$ 
and $A_0^{\bot}$, both with $C$-dense orbits. Using the fact that the unitary part 
of $A_0$ is trivial (that is, $A_0$ is a representation by translations) and 
that $\Gamma$ is finitely generated, conclude that $A_0^{\bot}$ has infinite 
dimension. Obtain a contradiction by noticing that over $\mathcal{H}_0^{\bot}$
the affine action of $H$ is trivial, and hence $A_0^{\bot}$ induces an affine 
action of the quotient group $\Gamma / H$ whose nilpotence degree is 
smaller than that of $\Gamma$.
\label{ejer-gelander}
\end{ejer}

\begin{ejer} Give examples of non minimal, isometric actions on infinite 
dimensional Hilbert spaces of finitely generated solvable groups.
\index{group!solvable}
\end{ejer}
\end{small}


\subsection{Super-rigidity for actions on the circle}

\hspace{0.45cm} Let $\Gamma$ be a subgroup of $\mathrm{Diff}_+^{1+\tau}(\clo)$, 
where $\tau > 1/2$. According to \S \ref{pruebabacan} 
(see Remark \ref{pillo}), if the restriction of the Liouville cocycle to $\Gamma$ 
is cohomologically trivial, then $\Gamma$ is topologically conjugate to a subgroup 
of $\mathrm{PSL}(2,\mathbb{R})$. Using Delorme's lemma we will study the case 
where this cocycle is almost cohomologically trivial. The next lemma should 
be compared with Exercise \ref{sh-ejer}, item (iii).
\index{cocycle!Liouville}

\vspace{0.12cm}

\begin{lem} {\em Suppose that the Liouville cocycle restricted to $\Gamma$ is an almost 
coboundary which is not cohomologically trivial. Then $\Gamma$ preserves a probability 
measure on the circle.}
\label{inv-mes}
\end{lem}

\noindent{\bf Proof.} By Delorme's lemma, if $c$ is an almost coboundary which is 
not cohomologically trivial, then $\Psi$ almost has invariant vectors. Therefore, 
there exists a sequence $(K_n)$ of unitary vectors in $\mathcal{H} \!=\! 
\mathcal{L}^{2,\Delta}_{\mathbb{R}}(\clo\times\clo,Lv)$ such that, for every 
$g \!\in\! \Gamma$, the value of \esp $\|K_n - \Psi(g)K_n\|$ \esp converges 
to zero as $n$ goes to infinity. Let $\mu_n$ be the probability measure on 
$\clo$ defined by 
$$\mu_n(A) = \int_{\clo}\int_A K_n^2(x,y) d Lv.$$
For every continuous function 
$\varphi: \clo \rightarrow \mathbb{R}$ we have
\begin{eqnarray*}
\big| \mu_n(\varphi) - g(\mu_n)(\varphi) \big| \!
&\leq&  \! \|\varphi\|_{\mathcal{L}^{\infty}} \int_{\clo}\int_{\clo}
\left| K_n^2 - (\Psi(g)K_n)^2 \right| d Lv\\
\! &\leq& \! \|\varphi\|_{\mathcal{L}^{\infty}}
\|K_n + \Psi(g)K_n \|_{\mathcal{L}^2} \|K_n - \Psi(g)K_n \|_{\mathcal{L}^2}\\
\! &\leq& \! 2\|\varphi\|_{\mathcal{L}^{\infty}} \|K_n - \Psi(g)K_n \|_{\mathcal{L}^2}.
\end{eqnarray*}
From this it follows that \esp 
$|\mu_n(\varphi) - g(\mu_n)(\varphi)|$ \esp goes to zero as $n$ goes 
to infinity. Therefore, if $\mu$ is an accumulation point of $(\mu_n)$, then $\mu$ 
is a probability measure on $\clo$ which is invariant by $\Gamma$. $\hfill\square$

\vspace{0.5cm}

\noindent{\bf Proof of Theorem \ref{primero1}.} For each unitary vector $K\!\in\!\mathcal{H}\! 
= \!\mathcal{L}^{2,\Delta}_{\mathbb{R}}(\clo \times \clo,Lv)$, let $\mu_{K}$ be the probability 
measure of $\clo$ obtained by projecting on the first coordinate, that is,
$$\mu_{K}(A) = \int_{\clo}\int_{A} K^2(x,y) d Lv.$$
Let $prob$ be the map $prob (K) = \mu_K$ defined on the unit sphere of $\mathcal{H}$ and taking 
values in the space of probability measures of the circle. This map $prob$ is $\Gamma$-equivariant, 
in the sense that for all $g \in \Gamma$ and every unitary vector 
$K \in \mathcal{L}^{2,\Delta}_{\mathbb{R}}(\clo \times \clo,Lv)$,
\begin{equation}
prob(\Psi(g)K) = \Phi(g) (prob(K)).
\label{equivariant}
\end{equation}

Suppose that $\Phi(\Gamma)$ preserves no probability measure on $\clo$ and that 
$\Phi(\Gamma)$ is not conjugate to a subgroup of $\mathrm{PSL}(2,\mathbb{R})$. 
By Lemma \ref{inv-mes} and the discussion before it, Shalom's super-rigidity theorem 
provides us with a family $\{ \mathcal{H}_1, \ldots, \mathcal{H}_k \}$ of 
$\Psi(\Gamma)$-invariant subspaces of $\mathcal{H}$, 
as well as of cocycles $c_i \!: \Gamma \rightarrow \mathcal{H}_i$, such that 
at least one of them is not identically zero, and such that over each $\mathcal{H}_i$ 
the isometric action associated to $c_i$ extends continuously to $G$ and factors 
through $G_i$. Take an index $i \!\in\! \{1,\ldots,k\}$ so that $\mathcal{H}_i$ is nontrivial. 
We claim that the image of the unit sphere in $\mathcal{H}_i$ by the map $prob$ consists 
of at least two different measures. Indeed, if this image were made of a single measure 
$prob(K)$ then, due to (\ref{equivariant}) and to the fact that $\mathcal{H}_i$ is 
a $\Psi(\Gamma)$-invariant subspace, $prob(K)$ would be invariant by 
$\Gamma$, thus contradicting our hypothesis. Fix an orthonormal basis 
$\{K_1,K_2,\ldots\}$ of $\mathcal{H}_i$, and define 
$$\overline{K} = \sum_{n \geq 0} \frac{| K_n |}{2^n}, \qquad
K = \frac{\overline{K}}{ \| \overline{K} \|}.$$
The measure $\mu_{K}$ has ``maximal support'' among those obtained by projecting functions in 
$\mathcal{H}_i$. Moreover, $\mu_{K}$ has no atom, and it is absolutely continuous with respect 
to Lebesgue measure. We denote by $\Lambda$ the closure of the support of $\mu_K$, which is 
a compact set without isolated points. Since $\mathcal{H}_i$ is $\Psi(\Gamma)$-invariant, 
$\Lambda$ is $\Gamma$-invariant, and since $\Phi(\Gamma)$ has no invariant measure, 
$\Lambda$ is not reduced to the union of finitely many disjoint intervals.

If $\Lambda$ is not the whole circle, let us collapse to a point the closure of each 
connected component of $\clo \setminus \Lambda$. We thus obtain a topological circle 
$\clo_{\Lambda}$ over which the original action $\Phi$ induces an action by 
homeomorphisms $\Phi_{\Lambda}$. However, notice that the orbits of this induced 
action are not necessarily dense, since $\Lambda$ might be bigger than the 
(non-empty) minimal invariant closed set of the original action. If $\Lambda$ 
is the whole circle, we let $\clo_{\Lambda} = \clo$. In any case, $\clo_{\Lambda}$ 
is endowed with a natural metric, since it may parameterized by means of the 
mesure $\mu_{K}$.

Let $K'$ be a function in the unit sphere of $\mathcal{H}_i$ such that the measure 
$\mu_{K'}$ is different from $\mu_K$, and let $\Gamma_{\mu_K}$ (resp. $\Gamma_{\mu_{K'}}$) 
be the group of homeomorphisms of $\clo_{\Lambda}$ that preserve the measure (induced on 
$\clo_{\Lambda}$ by) $\mu_K$ (resp. $\mu_{K'}$). Notice that $\Gamma_{\mu_K}$ is topologically 
conjugate to the group of rotations. If $(g_n)$ is a sequence of elements in $\Gamma$ such 
that $\lim_{n \rightarrow +\infty} pr_i(g_n) \!=\! id_{G_i}$, then both 
$\| \Psi(g_n)K \!-\! K \|$ and $\|\Psi(g_n) K' \!-\! K' \|$ converge to 
zero as $n$ goes to infinity. An analogous argument to that of the proof 
of Lemma \ref{inv-mes} then shows that $(\Phi(g_n))_{*}(\mu_K)$ (resp. 
$(\Phi(g_n))_{*}(\mu_{K'})$) converges to $\mu_K$ (resp. $\mu_{K'}$) as $n$ 
goes to infinity. Using this, one easily concludes that $(\Phi_{\Lambda}(g_n))$ 
has accumulation points in $\mathrm{Homeo}_{+}(\clo_{\Lambda})$, all of them  
contained in $\Gamma_{\mu_{K}} \cap \Gamma_{\mu_{K'}}$. Since $\mu_{K'}$ is 
different from $\mu_K$ and its support is contained in that of $\mu_K$, the 
group $\Gamma_{\mu_{K}} \cap \Gamma_{\mu_{K'}}$ is strictly contained in 
$\Gamma_{K}$. Since $\Gamma_{\mu_{K}} \cap \Gamma_{\mu_{K'}}$ is closed in
$\mathrm{Homeo}_+(\clo_{\Lambda})$, and since every non-dense subgroup of 
the group of rotations is finite, we conclude that 
$\Gamma_{\mu_{K}} \cap \Gamma_{\mu_{K'}}$ must be finite.

Let $H$ be the set of elements $h \!\in\! \mathrm{Homeo}_+(\clo_{\Lambda})$ 
such that \esp $h \!=\! \lim_{n \rightarrow +\infty} \Phi_{\Lambda}(g_n)$ 
\esp for some sequence $(g_n)$ in $\Gamma$ satisfying \esp 
$\lim_{n \rightarrow +\infty} pr_i(g_n) \!=\! id_{G_i}$. \esp 
By definition, $H$ is a closed subgroup of $\mathrm{Homeo}_+(\clo_{\Lambda})$. 
Moreover, the argument above shows that $H$ is contained in 
$\Gamma_{\mu_{K}} \cap \Gamma_{\mu_{K'}}$. It must therefore be finite and 
cyclic, say of order $d$. In the case $d>1$, we choose a generator $h$ for $H$, 
and we notice that the rotation number $\rho(h)$ is non-zero. Let $(g_n)$ be a 
sequence in $\Gamma$ such that $\lim_{n \rightarrow +\infty} pr_i(g_n)=id_{G_i}$
and $h = \lim_{n \rightarrow +\infty} \Phi_{\Lambda}(g_n)$.

We now show that $H$ is contained in the centralizer of $\Phi_{\Lambda}(\Gamma)$ 
in $\mathrm{Homeo}_+(\clo_{\Lambda})$.
\index{centralizer} 
To do this, notice that for each $g \in \Gamma$ the sequence of maps $pr_i(g^{-1} g_n g)$ 
also tends to $id_{G_i}$ as $n$ goes to infinity. By definition, 
$(\Phi_{\Lambda}(g^{-1}g_ng))$ converges to some element $h^j \!\in\! H$, 
where $j \in \{1,\ldots,d\}$. From the equality \esp 
$\rho(\Phi_{\Lambda}(g^{-1}g_ng)) = \rho(g^{-1}g_ng) 
= \rho(g_n) = \rho(\Phi_{\Lambda}(g_n))$ \esp 
it easily follows that $j \!=\! 1$, which implies that 
$\Phi_{\Lambda}(g)$ commutes with $h$. Since $g \!\in\! \Gamma$ was an 
arbitrary element, the group $H$ centralizes $\Phi_{\Lambda}(\Gamma)$.

Let $\clo_{\Lambda}/\!\!\sim$ be the topological circle obtained by identifying the points 
in $\clo_{\Lambda}$ which are in the same orbit of $H$. This circle $\clo_{\Lambda}$ is a 
finite covering of degree $d$ of $\clo_{\Lambda}/\!\!\sim$. 
\index{degree!of a covering}
Moreover, \esp $\Phi_{\Lambda}\!: \Gamma \!\rightarrow \!\mathrm{Homeo}_{+}(\clo_{\Lambda})$
\esp naturally induces a representation \esp $\tilde{\Phi}\!: \Gamma \!\rightarrow\!
\mathrm{Homeo}_+ \big( \clo_{\Lambda}/\!\!\sim\!\big)$ \esp 
such that, if $(g_n)$ is a sequence in $\Gamma$ for which $pr_i(g_n)$ tends to $id_{G_i}$, 
then $\tilde{\Phi}(g_n)$ tends to the identity (on $\clo_{\Lambda}/\!\!\sim$). Lemma 
\ref{extender} then allows concluding that $\tilde{\Phi}$ extends to a representation 
$\hat{\Phi}\!: G \!\rightarrow\! \mathrm{Homeo}_+(\clo_{\Lambda}/\!\!\sim)$ which 
factors through $G_i$. This representation $\hat{\Phi}$ is the one that extends $\Phi$ 
up to a semiconjugacy and a finite covering, as we wanted to show. The proof of 
Theorem \ref{primero1} is thus concluded. $\hfill\square$

\vspace{0.5cm}

Now recall that the finite subgroups of $\mathrm{Homeo}_+(\clo)$ are topologically 
conjugate to groups of rotations, and thus to subgroups of $\mathrm{PSL}(2,\mathbb{R})$. 
Hence, to prove Corollary \ref{corolarito} it suffices to show that, if $\Phi(\Gamma)$ 
preserves a probability measure on the circle and 
$\mathrm{H}^1_{\mathbb{R}}(\Gamma_0) \! = \! \{ 0 \}$ for every finite index normal subgroup 
$\Gamma_0$ of $\Gamma$, then $\Phi(\Gamma)$ is finite. For this, notice that if the invariant 
measure has no atom, then $\Phi(\Gamma)$ is semiconjugate to a group of rotations; in  
the other case, $\Phi(\Gamma)$ has a finite orbit. We claim that this implies that 
$\Phi(\Gamma)$ must be finite. Indeed, in the case of a finite orbit, this follows 
from Thurston's stability theorem, whereas in the case of a semiconjugacy to a group 
of rotations, this follows by taking the rotation number homomorphism 
(recall that $\Gamma$ --and hence $\Gamma_0$-- is finitely generated). 
The proof of Corollary \ref{corolarito2} is thus concluded.
\index{Thurston!stability theorem}

\vspace{0.425cm}

\noindent{\bf Proof of Theorem \ref{ultimito}.} Since $G$ is a connected Lie group, its action 
$\hat{\Phi}$ on $\clo_{\Lambda}/\!\!\sim$ factors through homomorphisms into $(\mathbb{R},+)$,
$\mathrm{Aff}_+(\mathbb{R})$, $\mathrm{SO}(2,\mathbb{R})$, $\widetilde{\mathrm{PSL}}(2,\mathbb{R})$, 
or $\mathrm{PSL}_k(2,\mathbb{R})$ for some $k \!\geq\! 1$ (see \S \ref{Lie}). In what follows,  
suppose that $\Phi(\Gamma)$ is infinite, which according to the proof of Corollary 
\ref{corolarito2} is equivalent to that $\Phi(\Gamma)$ has no invariant 
probability measure on $\clo$.

Let us first consider the case of hypothesis (i). The circle $\clo_{\Lambda}$ 
then identifies with the original circle $\clo$. From the fact that the kernel of 
$\Phi$ is finite, one concludes that there exist sequences $(g_n)$ in $\Gamma$ 
such that $pr_i(g_n)$ converges to $id_{G_i}$ and the homeomorphisms 
$\hat{\Phi}(g_n)$ are two-by-two different. This implies that the Lie 
group $\hat{\Phi}(G_i)$ cannot be discrete. From the previously recalled 
classification, the connected component of the identity in this group 
$\hat{\Phi}(G_i)_0$ corresponds either to $\mathrm{SO}(2,\mathbb{R})$, 
to $\mathrm{PSL}_{k}(2,\mathbb{R})$ for some $k \!\geq\! 1$, or to a 
subgroup of a product of groups of translations, of affine groups, and/or of 
groups conjugate to $\widetilde{\mathrm{PSL}}(2,\mathbb{R})$, all of them acting 
on disjoint intervals. The first case cannot arise, since $\phi(\Gamma)$ does not 
fix any probability measure on $\clo$. The last case cannot arise either, because 
the orbits by $\Phi(\Gamma)$ are dense and $\hat{\Phi}(G_i)_0$ is normal in 
$\hat{\Phi}(G_i)$ (the family of the intervals fixed by $\hat{\Phi}(G_i)_0$ 
must be preserved by $\Gamma$). The group $\hat{\Phi}(G_i)_0$ is therefore 
conjugate to $\mathrm{PSL}_k(2,\mathbb{R})$ for some $k \!\geq\! 1$, and 
since $\mathrm{PSL}_k (2,\mathbb{R})$ coincides with its normalizer in 
$\mathrm{Homeo}_+(\clo)$, the same holds for $\hat{\Phi}(G_i)$.

Let us now consider assumption (ii), that is, let us suppose 
that $\Phi$ takes values in the group of real-analytic 
circle diffeomorphisms. We have already observed that the orbits of the action of $\Gamma$ on 
$\clo_{\Lambda}$ are not necessarily dense. Let $\dot{\Lambda}$ be the minimal invariant 
closed set of this action, and let $\dot{\Phi}\!: \Gamma
\!\rightarrow \!\mathrm{Homeo}(\clo_{\dot{\Lambda}})$ 
be the action induced on the topological circle $\clo_{\dot{\Lambda}}$ obtained after collapsing 
the closure of each connected component of $\clo_{\Lambda} \setminus \dot{\Lambda}$. The orbits 
by $\dot{\Phi}$ are dense. For applying the arguments used in case (i), it suffices to show that 
the kernel of $\dot{\Phi}$ is finite. But this is evident, since the fixed points of the 
nontrivial elements in $\Phi(\Gamma)$ are isolated (the kernel of the restriction of 
$\dot{\Phi}$ a $\Gamma$ coincides with the kernel of $\Phi$).

Finally, let us assume hypothesis (iii), that is, $\Gamma$ satisfies 
Margulis' normal subgroup theorem. Once again, we need to show that the 
kernel of $\dot{\Phi}$ is finite. Now if this were not the case, then this 
kernel would have finite index in $\Gamma$. This would imply that the orbits 
of points in $\dot{\Lambda}$ are finite. However, this is impossible, since 
all the orbits of $\dot{\Phi}$ are dense. The proof of Theorem \ref{ultimito} 
is thus concluded. $\hfill\square$

\index{super-rigidity theorem!for actions on the circle|)}

\newpage
\thispagestyle{empty}
${}$
\newpage

\appendix

\addcontentsline{toc}{chapter}{APPENDIX}

\chapter{Some Basic Concepts in Group Theory}

\hspace{0.45cm} Let $\Gamma_1$ and $\Gamma_2$ be subgroups of a group $\Gamma$. We 
denote by $[\Gamma_1,\Gamma_2]$ the group generated by the elements of the form 
$[f,g]=f^{-1}g^{-1}fg$, where $f \in \Gamma_1$ and $g \in \Gamma_2$. The group 
$\Gamma' = [\Gamma,\Gamma]$ is called the {\bf \em{derived}} (or {\bf \em{commutator}}) 
subgroup of $\Gamma$}. One says that $\Gamma$ is {\bf\em{metabelian}} if 
$[\Gamma,\Gamma]$ is Abelian, and {\bf \textit{perfect}}  
\index{group!perfect}\index{group!metabelian}if $[\Gamma,\Gamma] = \Gamma$.

Recall that a subgroup $\Gamma_{0}$ of $\Gamma$ is {\bf\textit{normal}} if for 
every $h \in \Gamma_0$ and all $g \in \Gamma$ one has $ghg^{-1} \in \Gamma_0$. 
The group $\Gamma$ is {\bf\textit{simple}} if the only normal subgroups are 
$\{ id \}$ and $\Gamma$ itself. 
\index{group!simple} 
Notice that the subgroup $[\Gamma,\Gamma]$ is normal in $\Gamma$. Hence, if 
$\Gamma$ is non-Abelian and simple, then it is perfect. The {\bf \em{center}} 
of a group is the subgroup formed by the elements which commute with all the 
elements in the group. More generally, given a subset $A$ of $\Gamma$, the 
{\bf\em{centralizer}} 
\index{centralizer} 
of $A$ (in $\Gamma$) is the subgroup of $\Gamma$ formed by the 
elements which commute with all the elements in $A$.\\

Given a group $\Gamma$, we define inductively the subgroups
$$\Gamma_0^{\mathrm{nil}} = \Gamma,
\qquad \Gamma_0^{\mathrm{sol}} = \Gamma,$$
$$\Gamma_{i+1}^{\mathrm{nil}} = [\Gamma,\Gamma^{\mathrm{nil}}_i],
\qquad \Gamma^{\mathrm{sol}}_{i+1} =
[\Gamma^{\mathrm{sol}}_{i},\Gamma^{\mathrm{sol}}_{i}].$$
The series of subgroups $\Gamma^{\mathrm{nil}}_i$ (resp. $\Gamma^{\mathrm{sol}}_{i}$) 
is called the {\bf\em{central series}} 
\index{central series} 
(resp {\bf\em{derived series}}) of $\Gamma$.  
The group is {\bf\textit{nilpotent}} (resp. {\bf \textit{solvable}}) if there 
exists $n \!\in\! \mathbb{N}$ such that $\Gamma^{\mathrm{nil}}_{n} \!=\! \{ id \}$ 
(resp. $\Gamma^{\mathrm{sol}}_{n} \!=\! \{ id \}$). The minimum integer $n$ for 
which this happens is called the {\bf \em{degree}} (also called {\bf\em{order}}) 
\index{degree!of nilpotence}\index{degree!of solvability}of 
nilpotence (resp. solvability) of the group. From the definitions it easily 
follows that every nilpotent group is solvable. A group is {\bf\textit{virtually 
nilpotent}} (resp. {\bf \textit{virtually solvable}}) if it contains a 
nilpotent (resp. solvable) subgroup of finite index.
\index{group!nilpotent} 

\vspace{0.05cm}

\begin{small} \begin{ejer}
Prove that the center of 
\index{center of a group}
every nontrivial nilpotent group is nontrivial. 
Show that each of the subgroups $\Gamma^{\mathrm{nil}}_{i}(\Gamma)$ and 
$\Gamma^{\mathrm{sol}}_{i}(\Gamma)$ is normal in $\Gamma$. Conclude that every 
nontrivial solvable group contains a nontrivial normal subgroup which is Abelian.
\end{ejer} \end{small}

\vspace{0.05cm}

Let $\mathrm{P}$ be some property for groups. For instance, $\mathrm{P}$ could be 
the property of being finite, Abelian, nilpotent, free, etc. One says that 
a group $\Gamma$ is {\bf\em{residually}} $\mathrm{P}$ if for every 
$g \!\in\! \Gamma$ different from the identity there exists a group $\Gamma_g$
satisfying $\mathrm{P}$ and a group homomorphism from $\Gamma$ into 
$\Gamma_g$ so that the image of $g$ is not trivial. One says that 
$\Gamma$ is {\bf\em{locally}} $\mathrm{P}$ if every finitely 
generated subgroup of $\Gamma$ has property P.

\vspace{0.05cm}

\begin{small} 
\begin{ejer} Show that a group $\Gamma$ is residually nilpotent (resp. residually solvable) 
if and only if \esp $\cap_{i \geq 0}  \Gamma^{\mathrm{nil}}_i = \{id\}$ \esp 
(resp. \esp $\cap_{i \geq 0} \Gamma^{\mathrm{sol}}_{i} = \{id\}$).
\label{se-uso}
\end{ejer}

\begin{ejer}
After reading \S \ref{sec-witte} and \S \ref{super-witte-morris}, show that every locally 
orderable (resp. locally bi-orderable, locally $\mathcal{C}$-orderable) group is orderable 
(resp. bi-orderable, $\mathcal{C}$-orderable).

\index{group!orderable}

\noindent{\underbar{Hint.}} For a general orderable group, endow $\mathcal{O}(\Gamma)$ 
with a natural topology and use Tychonov's theorem to prove compactness 
(see \cite{ordering2} in case of problems with this).
\label{nosecomo}
\end{ejer} \end{small}

A solvable group $\Gamma$ for which the subgroups $\Gamma^{\mathrm{sol}}_{i}$ are 
finitely generated is said to be {\bf \em{polycyclic}}. If these subgroups are 
moreover torsion-free, then $\Gamma$ is {\bf \em{strongly polycyclic}}. It is 
easy to check that every subgroup of a polycyclic group $\Gamma$ is polycyclic, 
and hence finitely generated. It follows that every family of subgroups 
of $\Gamma$ has a maximal element (with respect to the inclusion). The maximal 
element of the family of normal nilpotent subgroups is called the {\bf \em{nilradical}}
\index{group!polycyclic}
of $\Gamma$, and is commonly denoted by $N(\Gamma)$. Notice that $N(\Gamma)$ 
is not only normal in $\Gamma$ but also 
{\bf \em{characteristic}}, that is, stable under any automorphism of $\Gamma$.
\index{group!strongly polycyclic} 
It is possible to show that there exists a finite index subgroup $\Gamma_0$ 
\index{characteristic subgroup}
of $\Gamma$ such that $[\Gamma_0,\Gamma_0] \subset N(\Gamma)$ (see \cite[Chapter IV]{Ra}). 
In particular, the quotient $\Gamma_0 / N(\Gamma)$ is Abelian.
\index{nilradical}


\chapter{Invariant Measures and Amenable Groups\index{group!amenable|(}}


\hspace{0.45cm} Amenability is one of the most profound concepts in Group Theory, and 
it is certainly impossible to give a full treatement in just a few pages. We will 
therefore content ourselves by exploring a dynamical view of this notion which 
has been exploited throughout the text. For further information, we refer the 
reader to \cite{greenleaf,lubotsky,wagon,Zi}. 

We begin by recalling a classical theorem in Ergodic Theory due to Bogolioubov and 
Krylov: Every homeomorphism of a compact metric space preserves a probability measure. 
This result fails to hold for general group actions: consider for instance the action 
of a Schottky group on the circle (see \S \ref{sec-minimal}). 
\index{Bogolioubov-Krylov's theorem}\index{group!Schottky}

\vspace{0.1cm}

\begin{defn} A group $\Gamma$ is {\bf \em{amenable}} ({\bf \em{moyennable}}, in the 
French terminology) if every action of $\Gamma$ by homeomorphisms of a compact metric 
space admits an invariant probability measure.
\end{defn}

\vspace{0.1cm}

In order to get some insight into this definition (at least for the case of 
finitely generated groups), let us first recall the strategy of proof 
of Bogolioubov-Krylov's theorem. Given a homeomorphism $g$ of a compact 
metric space $\mathrm{M}$, we fix a probability measure $\mu$ on 
it, and we consider the sequence $(\nu_n)$ defined by 
$$\nu_n = \frac{1}{n} \big[ \mu + g(\mu) + g^2 (\mu)+ \ldots + g^{n-1} (\mu) \big].$$
Since the space of probability measures on $\mathrm{M}$ is compact when endowed with 
\index{weak-* topology} the weak-* topology, 
there exists a subsequence $(\nu_{n_k})$ of $(\nu_n)$ weakly 
converging to a probability measure $\nu$. We claim that this limit measure $\nu$ 
is invariant by $g$. Indeed, for every $k$ we have 
$$g(\nu_{n_k}) = \nu_{n_k} + \frac{1}{n_k} \left[ g^{n_k} (\mu) - \mu \right],$$
and this implies that 
\begin{equation}
g (\nu) \!=\! g \big( \lim\limits_{k \rightarrow \infty} \nu_{n_k} \big)
\!=\! \lim\limits_{k \rightarrow \infty} g (\nu_{n_k}) 
\!=\! \lim\limits_{k \rightarrow \infty} \nu_{n_k} +
\lim\limits_{k \rightarrow \infty} \frac{1}{n_k} \big[ g^{n_k} (\mu)- \mu \big]
\!=\! \lim\limits_{k \rightarrow \infty} \nu_{n_k} \!\!=\! \nu.
\label{prueba de b-k}
\end{equation}

Let us now try to repeat this argument for a group $\Gamma$ generated by a 
finite symmetric family of elements $\mathcal{G} \!=\! \{ g_1, \ldots , g_m\}$ 
(recall that \textit{symmetric} means that $g^{-1}$ belongs to $\mathcal{G}$ 
for every $g \!\in\!\mathcal{G}$). For each $g \in \Gamma$ define the 
{\bf \em{length}} of $g$ as the minimal number of (non-necessarily distinct) 
elements in $\mathcal{G}$ which are necessary for writing $g$ as a product. 
More precisely, 
$$\lh (g) = \mathrm{min} \{ n \in \mathbb{N} \!: 
\esp g = g_{i_n} g_{i_{n-1}}\!\! \cdots g_{i_1},
\hspace{0.35cm} g_{i_j} \in \mathcal{G} \}.$$
We call the {\bf \em{ball of radius $n$}} (with respect to $\mathcal{G}$) the 
set $B_{\mathcal{G}}(n)$ formed by the elements in $\Gamma$ having length smaller 
than or equal to $n$, and we denote by $L_{\mathcal{G}}(n)$ its cardinal.

Consider now an action of $\Gamma$ by homeomorphisms of a compact metric space 
$\mathrm{M}$, and let $\mu$ be a probability measure on $\mathrm{M}$. For each 
$n \!\in\! \mathbb{N}$ let us consider the probability measure 
$$\nu_n \esp = \esp \frac{1}{L_{\mathcal{G}}(n-1)}
\sum\limits_{g \in B_{\mathcal{G}}(n-1)} g \esp\! (\mu).$$
Passing to some subsequence $(\nu_{n_k})$, we have the convergence to some 
probability measure $\nu$. The problem now is that $\nu$ is not necessarily 
invariant. Indeed, if we try to repeat the arguments of proof of equality 
(\ref{prueba de b-k}), then we should need to estimate the value 
of an expression of the form
$$\frac{1}{L_{\mathcal{G}}(n_k)}
\sum\limits_{\small\begin{array}{c}\lh (g)=n_k,
\\ g_i \in \mathcal{G} \end{array}} g_i g \esp \! (\mu).$$
However, this expression does not necessarily converge 
to zero, since it may happen that the number of elements in 
$B_{\mathcal{G}}(n_k) \setminus B_{\mathcal{G}}(n_k-1)$ is not 
negligible with respect to $L_{\mathcal{G}}(n_k)$. To deal 
with this problem, the following definitions become natural.

\vspace{0.01cm}

\begin{defn} The {\bf \em{geometric boundary}} $\partial A$ of a 
non-empty subset $A \subset \Gamma$ is defined as
\index{geometric boundary}
$$\partial A = \bigcup_{g \in \mathcal{G}} \big( A \hspace{0.02cm}
\Delta \hspace{0.02cm} g(A) \big),$$
where $\Delta$ stands for the symmetric difference between sets.
\end{defn}

\begin{defn} A {\bf \em{F\o lner sequence}} for a group $\Gamma$ is a sequence 
\index{F\o lner sequence}
$(A_{n})$ of finite subsets of $\Gamma$ such that 
$$\lim\limits_{n \rightarrow \infty}
\frac{card \esp \! (\partial A_n)}{card \esp \!(A_n)}= 0.$$
\end{defn}

\vspace{0.01cm}

With this terminology, Bogolioubov-Krylov's argument shows that if $\Gamma$ 
is a (finitely generated) group having a F\o lner sequence, then every action of 
$\Gamma$ by homeomorphisms of a compact metric space admits an invariant probability 
measure. In other words, groups having F\o lner sequences are amenable. The converse 
to this is also true, according to a deep result due to  F\o lner.

\vspace{0.15cm}

\begin{thm} \textit{A finitely generated group is 
amenable if and only if it admits a F\o lner sequence.}
\end{thm}

\vspace{0.1cm}

It is important to notice that the preceding characterization of amenability is independent 
of the system of generators. Indeed, it is not difficult to check that the quotient of the 
length functions with respect to two different systems of generators are bounded (by a 
constant which is independent of the element in the group). Using this, one easily checks 
that each F\o lner sequence with respect to one system naturally induces a F\o lner 
sequence with respect to the other one.

\begin{ejer} Show that finite groups are amenable. 
Show that the same holds for Abelian groups. 

\noindent{\underbar{Hint.}} For finitely generated Abelian groups, find an 
explicit F\o lner sequence. Alternatively, notice that 
if $f$ and $g$ are commuting homeomorphisms of a compact 
metric space, then the argument leading to (\ref{prueba de b-k}) and applied 
to a probability measure $\mu$ which is already invariant by $f$ gives a 
probability measure $\nu$ which is invariant by both $f$ and $g$.
\end{ejer}

\begin{small} \begin{ejer} Show that if $\Gamma$ is a finitely generated group 
containing a free subgroup on two generators, then $\Gamma$ is not amenable.

\noindent{\underbar{Hint.}} Consider the natural action of $\Gamma$ on $\{0,1\}^{\Gamma}$.
\label{librenomoy}
\end{ejer}

\begin{ejer} Show that amenability is stable under the so-called 
{\bf \em{elementary operations}}. More precisely, show that:

\vspace{0.05cm}

\noindent{(i) subgroups of amenable groups are amenable;}

\vspace{0.05cm}

\noindent{(ii) every quotient of an amenable group is amenable;}

\vspace{0.05cm}

\noindent{(iii) if $\Gamma$ is a group containing an amenable normal subgroup $\Gamma_0$ 
so that the quotient $\Gamma / \Gamma_0$ is amenable, then $\Gamma$ itself is amenable;}

\vspace{0.05cm}

\noindent{(iv) if $\Gamma$ is the union of amenable subgroups $\Gamma_i$ such that for 
all indexes $i,i'$ there exists $j$ such that both $\Gamma_i$ and $\Gamma_{i'}$ are 
contained in $\Gamma_j$, then $\Gamma$ is amenable.}

\vspace{0.05cm}

\noindent As an application, conclude that every virtually solvable group is amenable.
\label{clasemoy}
\end{ejer}

\begin{ejer} A finitely generated group $\Gamma$ has {\bf \em{polynomial growth}} 
if there exists a real polynomial $Q$ such that $L_{\mathcal{G}}(n) \leq Q(n)$ 
for every $n \in \mathbb{N}$. Show that every group of polynomial growth is 
amenable.\\

\noindent{\underbar{Remark.}} A celebrated theorem by Gromov establishes that 
a \index{Gromov!theorem on polynomial growth groups} 
\index{growth!of groups}finitely generated group 
has polynomial growth if and only if it is virtually 
nilpotent \cite{Gr2}. Let us point out that the ``easy'' implication of 
this theorem ({\em i.e.}, the fact that the growth of nilpotent groups is 
polynomial) is prior to Gromov's work and independently due to Bass and Guivarch.
\index{Bass}\index{Guivarch}
\label{polinomial}
\end{ejer}

\begin{ejer} A finitely generated group $\Gamma$ has 
{\bf\em{sub-exponential growth}} if for every $C > 0$ one has  
$$\liminf_{n \rightarrow \infty} \frac{L_{\mathcal{G}}(n)}{\exp(Cn)} = 0.$$
Show that finitely generated groups of sub-exponential growth are amenable.

\noindent{\underbar{Hint.}} Show that, for groups of sub-exponential growth, 
the sequence of balls in the group contains a F\o lner sequence.
\label{sub-prom}
\end{ejer}
\end{small}
\index{group!amenable|)}

\newpage
\thispagestyle{empty}
${}$
\newpage


\printindex 

\end{document}